\documentclass[12pt,thmsb,titlepage]{article}%
\usepackage{graphicx}
\usepackage{amsmath}
\usepackage[english]{babel}
\usepackage{makeidx}
\usepackage{amsfonts}
\usepackage{amssymb}
\usepackage{numinsec}
\usepackage{symlist}
\usepackage{longtable}%
\newtheorem{theorem}{Theorem}[section]

\newtheorem{algorithm}[theorem]{Algorithm}

\newtheorem{conjecture}[theorem]{Conjecture}
\newtheorem{corollary}[theorem]{Corollary}

\newtheorem{definition}[theorem]{Definition}
\newtheorem{example}[theorem]{Example}

\newtheorem{lemma}[theorem]{Lemma}

\newtheorem{proposition}[theorem]{Proposition}
\newtheorem{remark}[theorem]{Remark}

\makeindex
\begin{document}

\title{Mirror symmetry and tropical geometry}
\author{Janko B\"{o}hm}
\date{31.08.2007}
\maketitle

\begin{abstract}
Using tropical geometry we propose a mirror construction for monomial
degenerations of Calabi-Yau varieties in toric Fano varieties. The
construction reproduces the mirror constructions by Batyrev for Calabi-Yau
hypersurfaces and by Batyrev and Borisov for Calabi-Yau complete
intersections. We apply the construction to Pfaffian examples and recover the
mirror given by R\o dland for the degree $14$ Calabi-Yau threefold in
$\mathbb{P}^{6}$ defined by the Pfaffians of a general linear $7\times7$
skew-symmetric matrix.

We provide the necessary background knowledge entering into the tropical
mirror construction such as toric geometry, Gr\"{o}bner bases, tropical
geometry, Hilbert schemes and deformations. The tropical approach yields an
algorithm which we illustrate in a series of explicit examples.

\end{abstract}
\tableofcontents
\listoffigures

\newpage

\listofsymbols\newpage

\setcounter{section}{-1}

\section{Introduction}

Mirror symmetry is a phenomenon postulated by theoretical physics in the
context of string theory. The goal of string theory is the unification of
general relativity, describing gravity, with the standard model, which
describes the electroweak and strong coupling. These theories model nature in
the large respectively in small scales in such astonishing precision that it
is very hard to obtain experimental data on this unification. String theory
follows the idea to replace point particles by extended objects like a
$1$-sphere and to replace the $4$-dimensional spacetime by a Riemannian
manifold of dimension $10$, which is locally the product of a $4$-dimensional
Minkowski space and a $6$-dimensional compact Riemannian manifold $X$, too
small to appear in measurements. For two out of five possible string theories
the manifold $X$ turns out to be a $3$-dimensional complex manifold with
trivial canonical sheaf. These kind of manifolds are called Calabi-Yau
manifolds and were studied in mathematics for a long time before their
appearance in theoretical physics.
\index{Hodge theory}%
Hodge theory associates to $X$ the Hodge numbers%
\[
h^{p,q}\left(  X\right)  =\dim H_{\bar{\partial}}^{p,q}\left(  X\right)  =\dim
H^{q}\left(  X,\Omega_{X}^{p}\right)  \text{, }p,q=0,...,3
\]
The general framework of string theory predicts that one type of string theory
obtained from a Calabi-Yau manifold $X$ is equivalent to the second type of
string theory on another Calabi-Yau manifold $X^{\circ}$ and the
\index{Hodge numbers}%
Hodge numbers of these are related by%
\[
h^{p,q}\left(  X\right)  =h^{3-p,q}\left(  X^{\circ}\right)  \text{ }%
\forall\text{ }0\leq p,q\leq3
\]
Such a pair $X$ and $X^{\circ}$ is called a mirror pair, and the question
arises how to get $X^{\circ}$ from $X$ and vice versa.

The first mirror construction was given by Greene and Plesser in\linebreak%
\cite{GrPl Duality in CalabiYau moduli space} for the general quintic
threefold $X\subset\mathbb{P}^{4}$. As for Calabi-Yau manifolds $T_{X^{\circ}%
}\cong\Omega_{X^{\circ}}^{2}$, the mirror $X^{\circ}$
\index{hypersurface}%
should
\index{Greene-Plesser}%
satisfy%
\[
\dim H^{1}\left(  X^{\circ},T_{X^{\circ}}\right)  =h^{2,1}\left(  X^{\circ
}\right)  =h^{1,1}\left(  X\right)  =1
\]
Greene and Plesser
\index{quintic threefold}%
construct the mirror as a general element in the $1$-parameter family of
quintics%
\[
X_{\lambda}=\left\{  x_{0}^{5}+x_{1}^{5}+x_{2}^{5}+x_{3}^{5}+x_{4}^{5}+\lambda
x_{0}x_{1}x_{2}x_{3}x_{4}=0\right\}
\]
with fibers in a $\mathbb{Z}_{5}^{3}$-quotient of $\mathbb{P}^{4}$. This
$1$-dimensional parameter space contains the degeneration point $\lambda
=\infty$ corresponding to the union of $5$ hyperplanes $\left\{  x_{0}%
x_{1}x_{2}x_{3}x_{4}=0\right\}  $. Indeed, degenerations of Calabi-Yau
manifolds to varieties given by monomial ideals appear naturally in the
context of various mirror constructions.

Generalizing the construction by Greene and Plesser, Batyrev considers in
\cite{Batyrev Dual polyhedra and mirror symmetry for CalabiYau hypersurfaces
in toric varieties} anticanonical hypersurfaces in Gorenstein toric Fano
varieties. There is a one-to-one correspondence of the Gorenstein toric Fano
varieties $\mathbb{P}\left(  \Delta\right)  $ of dimension $n$, polarized by
$-K_{\mathbb{P}\left(  \Delta\right)  }$ with the reflexive polytopes
$\Delta\subset M\otimes\mathbb{R}$, where $M=\mathbb{Z}^{n}$. Recall that a
polytope $\Delta$ is called reflexive if $\Delta$ and its dual polytope
$\Delta^{\ast}$ are integral and contain $0$. So duality of reflexive
polytopes is an involution of the set of Gorenstein toric Fano varieties.
Batyrev proves that general elements in $\left\vert -K_{\mathbb{P}\left(
\Delta\right)  }\right\vert $ and $\left\vert -K_{\mathbb{P}\left(
\Delta^{\ast}\right)  }\right\vert $ form a mirror pair in the sense of
mirrored Hodge numbers generalized to singular varieties. In the following we
associate to Batyrev%
\'{}%
s data a monomial degeneration. Denote by $\Sigma\subset N\otimes\mathbb{R}$,
where $N=\operatorname*{Hom}\left(  M,\mathbb{Z}\right)  $, the fan
representing $Y=\mathbb{P}\left(  \Delta\right)  $, i.e., the set of cones
over the faces of the dual polytope $\Delta^{\ast}$. Generalizing the
homogeneous coordinate ring of projective space, the
\index{Cox ring}%
Cox ring of $Y$ is the polynomial ring%
\[
S=\mathbb{C}\left[  y_{r}\mid r\in\Sigma\left(  1\right)  \right]
\]
with variables corresponding to the $1$-dimensional cones $\Sigma\left(
1\right)  $ in $\Sigma$ and graded by the Chow group of divisors
$A_{n-1}\left(  Y\right)  $ of $Y$ via the exact sequence%
\[
0\rightarrow M\overset{A}{\rightarrow}\mathbb{Z}^{\Sigma\left(  1\right)
}\overset{\deg}{\rightarrow}A_{n-1}\left(  Y\right)  \rightarrow0
\]
A reflexive polytope has $0$ as its unique interior lattice point, so a
generic toric hypersurface comes with a natural monomial degeneration
\[
\left\{  t\cdot\left(
{\textstyle\sum\limits_{m\in\partial\Delta}}
c_{m}\cdot A\left(  m\right)  \right)  +\prod_{r\in\Sigma\left(  1\right)
}y_{r}=0\right\}
\]
with generic coefficients $c_{m}$.

Note that toric varieties also appear in the context of monomial degenerations
in the sense that the special fiber is a union of toric varieties. Indeed
toric geometry plays an important role in the context of mirror symmetry as
toric varieties have non-trivial geometry and still can contain a reducible
special fiber of a degeneration in a natural description as a union of toric
strata.\smallskip

As general setup, we consider a one parameter degeneration $\mathfrak{X}$ of
Calabi-Yau varieties with fibers in a toric Fano variety $Y$ with Cox ring $S$
and with reduced monomial special fiber $X_{0}$. The toric Fano variety $Y$ is
given by the fan over the faces of a Fano polytope $P$, which is an integral
polytope in $N\otimes\mathbb{R}$ containing $0$ as the unique interior lattice
point. So we generalize Gorenstein toric Fano varieties to the Mori category
of $\mathbb{Q}$-Gorenstein toric Fano varieties. Let $I\subset\mathbb{C}%
\left[  t\right]  \otimes S$ be the ideal of the total space of the
degeneration and $I_{0}\subset S$ the ideal of the monomial special fiber.

Given a polytope $\Delta$ we denote by $\operatorname*{Poset}\left(
\Delta\right)  $ the complex of faces of $\Delta$, which is a partially
ordered set with respect to inclusion. For the polytope $\Delta=P^{\ast}$ the
complex $\operatorname*{Poset}\left(  \Delta\right)  $ is isomorphic to the
complex of the strata of $Y$, which we denote by $\operatorname*{Strata}%
\left(  Y\right)  $. So the complex of strata $\operatorname*{Strata}\left(
X_{0}\right)  $ of the monomial special fiber of $\mathfrak{X}$ can be
considered as a subcomplex of $\operatorname*{Poset}\left(  \Delta\right)  $.

Using Gr\"{o}bner basis techniques, we construct from the degeneration
$\mathfrak{X}$ a new polytope $\nabla$ with a new subcomplex $B\left(
I\right)  \subset\operatorname*{Poset}\left(  \nabla\right)  $. We begin by
associating to $\mathfrak{X}$ the Gr\"{o}bner cone $C_{I_{0}}\left(  I\right)
$ of weights on $\mathbb{C}\left[  t\right]  \otimes S$ selecting the monomial
special fiber ideal $I_{0}$ as lead ideal of $I$. For every face $F$ of
$C_{I_{0}}\left(  I\right)  $ we have an initial ideal $\operatorname*{in}%
_{F}\left(  I\right)  $ of $I$, which is no longer monomial for the proper
faces of $C_{I_{0}}\left(  I\right)  $. Denote by $\operatorname*{Poset}%
\left(  C_{I_{0}}\left(  I\right)  \right)  $ the fan of the faces
\newsym[$\operatorname*{Poset}\left(  C  \right)  $]{fan generated by the cone $C$}{}of
$C_{I_{0}}\left(  I\right)  $. We consider the subfan $BF_{I_{0}}\left(
I\right)  \subset\operatorname*{Poset}\left(  C_{I_{0}}\left(  I\right)
\right)  $ of those faces of $C_{I_{0}}\left(  I\right)  $, which have no
monomial in their initial ideal. This fan is the intersection of the fan
$\operatorname*{Poset}\left(  C_{I_{0}}\left(  I\right)  \right)  $ with the
Bergman fan $BF\left(  I\right)  $, introduced in \cite{Bergman The
logarithmic limitset of an algebraic variety}. Essentially equivalent to
Bergman%
\'{}%
s original definition, we define the Bergman fan $BF\left(  I\right)  $ as the
closure of the image of the vanishing locus of $I$ over the field
$L=\overline{\mathbb{C}\left(  s\right)  }$ of Puisseux series under the
valuation map%
\begin{gather*}
L^{\ast}\times\left(  L^{\ast}\right)  ^{n}\rightarrow\mathbb{R}^{n+1}\\
\left(  t,y_{1},...,y_{n}\right)  \mapsto\left(  val\left(  t\right)
,val\left(  y_{1}\right)  ,...,val\left(  y_{n}\right)  \right)
\end{gather*}
Here we consider the torus $\left(  L^{\ast}\right)  ^{n}\cong\left(  L^{\ast
}\right)  ^{\Sigma\left(  1\right)  }/\operatorname*{Hom}\nolimits_{\mathbb{Z}%
}\left(  A_{n-1}\left(  Y\right)  ,L^{\ast}\right)  $ and $val$ denotes the
valuation associating to a power series its vanishing order, i.e., the
exponent of its lowest order term.

The cone $C_{I_{0}}\left(  I\right)  $ is contained in the half-space of
$t$-local orderings. Hence, intersecting it transversally with the hyperplane
of $t$-weight $w_{t}=1$, i.e., identifying the parameters $t$ and $s$, we
obtain a convex polytope $\nabla$. The polytope $\nabla$ is naturally
contained in $N\otimes\mathbb{R}$ and it turns out that $\nabla^{\ast}$ is
again a Fano polytope. Corresponding to $BF_{I_{0}}\left(  I\right)
=BF\left(  I\right)  \cap\operatorname*{Poset}\left(  C_{I_{0}}\left(
I\right)  \right)  $ we also obtain a subcomplex $B\left(  I\right)  $ of the
complex $\operatorname*{Poset}\left(  \nabla\right)  $ of faces of $\nabla$.
Let $K=\overline{\mathbb{C}\left(  t\right)  }$ be the field of Puisseux
series in the parameter $t$.

If $w$ is a weight vector in a face of the Bergman complex $B\left(  I\right)
$, we can consider the power series solutions of $I$ lying over $w$ via the
valuation map%
\begin{gather*}
\left(  K^{\ast}\right)  ^{n}\rightarrow\mathbb{R}^{n}\\
\left(  y_{1},...,y_{n}\right)  \mapsto\left(  val\left(  y_{1}\right)
,...,val\left(  y_{n}\right)  \right)
\end{gather*}
Taking the limit $t\rightarrow0$ of these solutions induces an inclusion
reversing map%
\[%
\begin{tabular}
[c]{llll}%
$\lim:$ & $B\left(  I\right)  $ & $\rightarrow$ & $\operatorname*{Strata}%
\left(  X_{0}\right)  $%
\end{tabular}
\]
from the complex $B\left(  I\right)  $ to the complex of strata of the special
fiber $X_{0}$ of $\mathfrak{X}$, which is a subcomplex of faces of
$\Delta=P^{\ast}$. It turns out that the complex $B\left(  I\right)  $
essentially is dual to the complex of strata of $X_{0}$.

The complex $\operatorname*{Poset}\left(  \nabla^{\ast}\right)  $ describes
the initial ideals of $I$ at the faces of $\nabla$. Consider the reduced
standard basis of $I$ in $S\otimes\mathbb{C}\left[  t\right]  /\left\langle
t^{2}\right\rangle $ with respect to a monomial ordering in the interior of
$\nabla$. If $F$ is a face of $\nabla$, then all initial forms with respect to
$F$ of the standard basis elements involve a minimal generator of $I_{0}$.
Hence, dividing for all initial forms the non special fiber monomials by the
special fiber monomial, we obtain a set of degree $0$ Cox Laurent monomials.
These monomials correspond via the Chow presentation sequence%
\[
0\rightarrow M\overset{A}{\rightarrow}\mathbb{Z}^{\Sigma\left(  1\right)
}\overset{\deg}{\rightarrow}A_{n-1}\left(  Y\right)  \rightarrow0
\]
to the lattice points of $F^{\ast}$.

In the following we denote the first order deformations of $X_{0}$ which are
characters of the torus $\left(  \mathbb{C}^{\ast}\right)  ^{\Sigma\left(
1\right)  }$ as $\left(  \mathbb{C}^{\ast}\right)  ^{\Sigma\left(  1\right)
}$-deformations. Note that the vector space $\operatorname*{Hom}\left(
I_{0},S/I_{0}\right)  _{0}$ of degree $0$ deformations has a basis of $\left(
\mathbb{C}^{\ast}\right)  ^{\Sigma\left(  1\right)  }$-deformations.

The lattice points of the faces of $\nabla^{\ast}$ have a two-fold interpretation:

\begin{itemize}
\item Let $F\in\nabla$ be a face, $\delta\in F^{\ast}$ a lattice point and
$A\left(  \delta\right)  =\frac{m_{1}}{m_{0}}\in\mathbb{Z}^{\Sigma\left(
1\right)  }$ with relative prime $m_{0}$ and $m_{1}$ the corresponding degree
$0$ Cox Laurent monomial. Then $\delta$ can be considered as a $\left(
\mathbb{C}^{\ast}\right)  ^{\Sigma\left(  1\right)  }$-deformation of $X_{0}$
of degree $0$ by associating to it the element of $\operatorname*{Hom}\left(
I_{0},S/I_{0}\right)  _{0}$ defined for minimal generators $m\in I_{0}$ by
$\delta\left(  m\right)  =\frac{m_{1}}{m_{0}}\cdot m$ if $m_{0}\mid m$ and $0$
otherwise. Here we use $\delta\in F^{\ast}$.

It turns out that $B\left(  I\right)  ^{\ast}$ represents the tangent space of
the component of the Hilbert scheme of $I_{0}$ containing $\mathfrak{X}$,
assuming that we took the tangent vector of $\mathfrak{X}$ general enough.

\item The fan over the faces of $\nabla^{\ast}$ defines a toric Fano variety
$Y^{\circ}$, so the vertices of $\nabla^{\ast}$ are the variables of the Cox
ring of $Y^{\circ}$, i.e., the torus invariant divisors on $Y^{\circ}$. Hence
in particular the vertices of the faces of $B\left(  I\right)  ^{\ast}$ have
an interpretation as torus invariant divisors on $Y^{\circ}$. Passing from
vertices to lattice points amounts to a toric blowup.
\end{itemize}

Mirror symmetry of a pair $X\ $and $X^{\circ}$ identifies $h^{1,\dim\left(
X\right)  -1}\left(  X\right)  $ and $h^{1,1}\left(  X^{\circ}\right)  $ and
vice versa and more generally the complex moduli space of $X$ with the
K\"{a}hler moduli space of $X^{\circ}$ and vice versa. From this point of view
$Y^{\circ}$ is the toric Fano variety with sufficient divisors to represent
locally the component of the Hilbert scheme at $X_{0}$ containing
$\mathfrak{X}$. In the same way as $B\left(  I\right)  ^{\ast}$ describes the
tangent space at $X_{0}$ of the component of the Hilbert scheme along
$\mathfrak{X}$, we expect that the K\"{a}hler classes given by the lattice
points of $B\left(  I\right)  ^{\ast}$ suffice to represent the K\"{a}hler
moduli of the mirror.

The subcomplex $B\left(  I\right)  \subset\operatorname*{Poset}\left(
\nabla\right)  $ defines a monomial ideal $I_{0}^{\circ}$ in the Cox ring
$S^{\circ}$ of $Y^{\circ}$. The ideal $I_{0}^{\circ}$ is the intersection over
all facets (i.e., faces of maximal dimension) $F$ of $B\left(  I\right)  $ of
the ideals generated by the set of all facets of $\nabla$ containing $F$. This
generalizes the idea of Stanley-Reisner rings.

So we have constructed a toric Fano variety $Y^{\circ}$ and a monomial ideal
$I_{0}^{\circ}$, whose zero locus $X_{0}^{\circ}$ essentially is
combinatorially dual to the complex of strata of the special fiber $X_{0}$ of
$\mathfrak{X}$.

We know that the lattice points of $B\left(  I\right)  ^{\ast}\subset
\operatorname*{Poset}\left(  \nabla^{\ast}\right)  $ have an interpretation as
first order deformations of $X_{0}$ contributing to tangent vector of the
family $\mathfrak{X}$. Hence the first order deformations of the mirror
special fiber $X_{0}^{\circ}$ contributing to the tangent vector the mirror
degeneration $\mathfrak{X}^{\circ}$ should be given by the lattice points of
the dual $\left(  \lim\left(  B\left(  I\right)  \right)  \right)  ^{\ast
}\subset\operatorname*{Poset}\left(  \Delta^{\ast}\right)  $ of the image of
the limit map. Again the lattice points of $\left(  \lim\left(  B\left(
I\right)  \right)  \right)  ^{\ast}$ have the two fold interpretation as
deformations of $X_{0}^{\circ}$ and torus divisors on a blowup of $Y$.
Applying these deformations to $I_{0}^{\circ}$ we obtain the conjectural
mirror degeneration up to first order.

If the ideal $I_{0}^{\circ}$ obeys a structure theorem, e.g., the Koszul
resolution for complete intersections or the structure theorem of Buchsbaum
and Eisenbud for codimension $3$ subcanonical varieties, we can (in the case
of complete intersection trivially) extend the first order mirror family to a
degeneration over $\operatorname*{Spec}\mathbb{C}\left[  t\right]
$.\smallskip

The tropical mirror construction formalizes as follows:

\begin{itemize}
\item Let $N=\mathbb{Z}^{n}$, let $P\subset N\otimes\mathbb{R}$ be a Fano
polytope and $Y$ the corresponding toric Fano variety with Cox ring $S$. Let
$\mathfrak{X}$ be a one parameter monomial degeneration of Calabi-Yau
varieties with fibers in $Y$ and let $\mathfrak{X}$ be given by the ideal
$I\subset\mathbb{C}\left[  t\right]  \otimes S$. Suppose that the ideal
$I_{0}$ of the special fiber is a reduced monomial ideal.

\item Fix a monomial ordering $>$ on $\mathbb{C}\left[  t\right]  \otimes S$,
which is respecting the Chow grading on $S$ and which is local in $t$, and
denote by $>_{w}$ the weight ordering by $w$ refined by $>$. Then define
\[
C_{I_{0}}\left(  I\right)  =\left\{  -\left(  w_{t},w_{y}\right)
\in\mathbb{R}\oplus N_{\mathbb{R}}\mid L_{>_{\left(  w_{t},\varphi\left(
w_{y}\right)  \right)  }}\left(  I\right)  =I_{0}\right\}
\]
Note that we add the minus sign as $L$ is defined as selecting the monomial of
maximal weight.

\item Intersecting $C_{I_{0}}\left(  I\right)  $ with the hyperplane of
$t$-weight one, we obtain a polytope%
\[
\nabla=C_{I_{0}}\left(  I\right)  \cap\left\{  w_{t}=1\right\}  \subset
N_{\mathbb{R}}%
\]
and $\nabla^{\ast}$ is a Fano polytope, so gives a toric Fano variety
$Y^{\circ}$.

\item The complex of the faces of the polytope $\nabla$ has the subcomplex%
\begin{align*}
B\left(  I\right)   &  =\left(  BF\left(  I\right)  \cap\operatorname*{Poset}%
\left(  C_{I_{0}}\left(  I\right)  \right)  \right)  \cap\left\{
w_{t}=1\right\} \\
&  =\left\{  F\text{ face of }\nabla\mid\operatorname*{in}\nolimits_{F}\left(
I\right)  \text{ does not contain a monomial}\right\}
\end{align*}
the Bergman subcomplex or tropical subcomplex of $\nabla$. The intersection of
the fan $BF\left(  I\right)  \cap\operatorname*{Poset}\left(  C_{I_{0}}\left(
I\right)  \right)  $ with $\left\{  w_{t}=1\right\}  $ is defined as the
complex, whose faces are the intersections of the cones of $BF\left(
I\right)  \cap\operatorname*{Poset}\left(  C_{I_{0}}\left(  I\right)  \right)
$ with the hyperplane $\left\{  w_{t}=1\right\}  $.

\item The complex $B\left(  I\right)  $ is a subdivision of the dual of the
complex of strata $\operatorname*{Strata}\left(  X_{0}\right)  $ of the
special fiber $X_{0}$ of $\mathfrak{X}$ via the map%
\[%
\begin{tabular}
[c]{llll}%
$\lim:$ & $B\left(  I\right)  $ & $\rightarrow$ & $\operatorname*{Strata}%
\left(  X_{0}\right)  \subset\operatorname*{Strata}\left(  Y\right)  $\\
& \multicolumn{1}{c}{$F$} & $\mapsto$ & $\left\{  \lim_{t\rightarrow0}a\left(
t\right)  \mid a\in\operatorname*{val}\nolimits^{-1}\left(
\operatorname*{int}\left(  F\right)  \right)  \right\}  $%
\end{tabular}
\]

taking the limit of arc solutions of $I$. Here $\operatorname*{int}\left(
F\right)  $ denotes the relative interior of $F$.

\item Denote by $\Sigma^{\circ}$ the fan over the faces of $\nabla^{\ast}$
defining $Y^{\circ}$ and by%
\[
S^{\circ}=\mathbb{C}\left[  z_{r}\mid r\in\Sigma^{\circ}\left(  1\right)
\right]
\]
the Cox ring of $Y^{\circ}$ graded via%
\[%
\begin{tabular}
[c]{lllllll}%
$0\rightarrow$ & $N$ & $\overset{A^{\circ}}{\rightarrow}$ & $\mathbb{Z}%
^{\Sigma^{\circ}\left(  1\right)  }$ & $\overset{\deg}{\rightarrow}$ &
$A_{n-1}\left(  Y^{\circ}\right)  $ & $\rightarrow0$%
\end{tabular}
\]
Then the monomial ideal defining the special fiber $X_{0}^{\circ}\subset
Y^{\circ}$ of the mirror degeneration $\mathfrak{X}^{\circ}$ of $\mathfrak{X}$
is%
\begin{align*}
I_{0}^{\circ}  &  =\left\langle
{\displaystyle\prod\limits_{v\in J}}
z_{v}\mid J\subset\Sigma^{\circ}\left(  1\right)  \text{ with }%
\operatorname*{supp}\left(  B\left(  I\right)  \right)  \subset%
{\displaystyle\bigcup\limits_{r\in J}}
F_{r}\right\rangle \\
&  =%
{\textstyle\bigcap\nolimits_{F\in B\left(  I\right)  }}
\left\langle z_{G^{\ast}}\mid G\text{ a facet of }\nabla\text{ with }F\subset
G\right\rangle \subset S^{\circ}%
\end{align*}
where $F_{r}$ denotes the facet of $\nabla$ corresponding to the
$1$-dimensional cone $r$ of $\Sigma^{\circ}=\operatorname*{NF}\left(
\nabla\right)  $. Note that in the second description of $I_{0}$ it is
sufficient to take the intersection over the maximal faces of $B\left(
I\right)  $.

\item Let $M=\operatorname*{Hom}\left(  N,\mathbb{Z}\right)  $ and
$\Delta=P^{\ast}\subset M\otimes\mathbb{R}$. The image of $\lim$ naturally is
a subcomplex of the complex of faces of $\Delta$. Hence we obtain a subcomplex
$\left(  \lim\left(  B\left(  I\right)  \right)  \right)  ^{\ast}$ of the
complex of faces of $\Delta^{\ast}=P$, which describes the first order
deformations of the mirror degeneration at $X_{0}^{\circ}$ as degree zero Cox
Laurent monomials. So the conjectural mirror degeneration up to first order is
given by
\[
\left\langle m+t\cdot\sum_{\alpha\in\operatorname*{supp}\left(  \left(
\lim\left(  B\left(  I\right)  \right)  \right)  ^{\ast}\right)  \cap
N}c_{\alpha}\cdot\alpha\left(  m\right)  \mid m\in I_{0}^{\circ}\right\rangle
\subset\mathbb{C}\left[  t\right]  /\left\langle t^{2}\right\rangle \otimes
S^{\circ}%
\]
with generic coefficients $c_{\alpha}$.

Note that the description of the first order deformations as lattice points of
$\Delta^{\ast}$ is independent of the toric variety $Y^{\circ}$. This easily
allows to replace $Y^{\circ}$ by different birational models in the Mori
category.\smallskip
\end{itemize}

The tropical mirror construction reproduces known mirror constructions.
Batyrev and Borisov extend in \cite{Borisov Towards the mirror symmetry for
CalabiYau complete intersections in Gorenstein Fano toric varieties} and
\cite{BB On CalabiYau complete intersections in toric varieties in
HigherDimensional Complex Varieties Trento 1994} the mirror construction for
hypersurfaces in toric Fano varieties to complete intersections given by nef
partitions. In an analogous way as we obtained the degeneration associated to
an anticanonical hypersurface in a Gorenstein toric Fano variety, we also
obtain a monomial degeneration for a complete intersection. We show that, when
applied to this complete intersection degeneration, the tropical mirror
construction gives the degeneration associated to the Batyrev-Borisov mirror.
In particular, this also holds true in the case of Batyrev%
\'{}%
s mirror construction for hypersurfaces.

We introduce the notion of Fermat deformations in order to relate the mirror
degenerations to birational models with fibers in toric Fano varieties with
Chow group of rank $1$. Applying this, we connect the mirror degeneration
associated to the complete intersection of two general cubics in
$\mathbb{P}^{5}$ to a Greene Plesser type orbifolding mirror family given in
\cite{LT Lines on CalabiYau complete intersections mirror symmetry and
PicardFuchs equations}.

In the same way, applying the tropical mirror construction to a monomial
degeneration of non complete intersection Calabi-Yau threefolds of degree $14$
in $\mathbb{P}^{6}$ defined by the Pfaffians of a general linear skew
symmetric map $7\mathcal{O}\left(  -1\right)  \rightarrow7\mathcal{O}$, we
reproduce the orbifolding mirror given by R\o dland in \cite{Ro dland The
Pfaffian CalabiYau its Mirror and their link to the Grassmannian mathbbG27}.

We also apply the tropical mirror construction to a monomial degeneration of
non complete intersection Calabi-Yau threefolds of degree $13$ in
$\mathbb{P}^{6}$ defined by the Pfaffians of a general skew symmetric map
$\mathcal{O}\left(  -2\right)  \oplus4\mathcal{O}\left(  -1\right)
\rightarrow\mathcal{O}\left(  1\right)  \oplus4\mathcal{O}$. From the mirror
degeneration given by the tropical mirror construction we obtain, via the
concept of Fermat deformations, a flat degeneration with fibers in an orbifold
of $\mathbb{P}^{6}$, which again obeys the structure theorem of
Buchsbaum-Eisenbud.\medskip

In the following, we give a short overview of the individual sections.

\textbf{Section \ref{Sec prerequisites}.} This section provides an
introduction to various concepts used in the tropical mirror construction.

Section \ref{Sec calabi yau and mirror symmetry} recalls some facts on
Calabi-Yau manifolds and their relation to string theory and mirror symmetry.
A manifold $X$ of dimension $d$ is called a Calabi-Yau manifold if
$K_{X}=\mathcal{O}_{X}$ and $h^{i}\left(  X,\mathcal{O}_{X}\right)  =0$ for
$0<i<d$.

In Section
\ref{Mirror symmetry for singular Calabi-Yau varieties and stringy Hodge numbers}
we give a short introduction to the concept of stringy Hodge numbers
introduced by Batyrev to generalize Hodge numbers to singular varieties. Given
a normal projective variety $X$ with log-terminal singularities one associates
to $X$, via a resolution $f:Y\rightarrow X$ of singularities, a function
$E_{st}\left(  X;u,v\right)  $, which Batyrev proves to be independent of the
choice of the resolution. If $E_{st}$ is a polynomial, then stringy Hodge
numbers can be defined via the coefficients of $E_{st}$. In any case,
topological mirror symmetry of a pair of Calabi-Yau varieties $X$ and
$X^{\circ}$ of dimension $d$ can be defined via the stringy $E$-functions as
the relation $E_{st}\left(  X;u,v\right)  =\left(  -u\right)  ^{d}%
E_{st}\left(  X^{\circ};u^{-1},v\right)  $. If $X$ admits a crepant resolution
$f:Y\rightarrow X$ then
\[
E_{st}\left(  X;u,v\right)  =\sum_{0\leq p,q\leq d}\left(  -1\right)
^{p+q}h^{p,q}\left(  Y\right)  u^{p}v^{q}%
\]

In Section \ref{Sec facts from toric geometry} we continue with an overview of
toric geometry. The Sections \ref{Sec affine toric varieties}%
--\ref{Sec morphisms of toric varieties} give the standard description of
toric varieties and morphisms. If $N\cong\mathbb{Z}^{n}$, $N_{\mathbb{R}%
}=N\otimes\mathbb{R}$, $\sigma\subset N_{\mathbb{R}}$ is a rational convex
polyhedral cone, $M=\operatorname*{Hom}\left(  N,\mathbb{Z}\right)  $ and
\[
\check{\sigma}=\left\{  m\in M_{\mathbb{R}}\mid\left\langle m,w\right\rangle
\geq0\ \forall w\in\sigma\right\}
\]
is the dual cone, then $\check{\sigma}\cap M$ is a finitely generated
semigroup and defines an affine toric variety $U\left(  \sigma\right)
=\operatorname*{Spec}\left(  \mathbb{C}\left[  \check{\sigma}\cap M\right]
\right)  $. Given a fan $\Sigma$ in $N_{\mathbb{R}}$, i.e., a finite set of
strongly convex rational polyhedral cones such that every face of a cone in
$\Sigma$ is again a cone in $\Sigma$, the $U\left(  \sigma\right)  $,
$\sigma\in\Sigma$ glue to a toric variety $Y=X\left(  \Sigma\right)  $. The
torus $\operatorname*{Spec}\left(  \mathbb{C}\left[  \mathbb{Z}^{n}\right]
\right)  \hookrightarrow Y$ acts on $Y$. There is an inclusion reversing
bijection between the cones of $\Sigma$ and the torus orbit closures. Let
$\Sigma\left(  1\right)  $ be the set of rays in $\Sigma$, i.e., the set of
$1$-dimensional cones. We denote by $D_{r}$ the torus invariant divisor on $Y$
corresponding to the ray $r\in\Sigma\left(  1\right)  $.

As explained in Section \ref{Sec dualizing sheaf of a toric variety}, one can
describe the dualizing sheaf of a toric variety $X\left(  \Sigma\right)  $ as%
\[
\Omega_{X\left(  \Sigma\right)  }^{n}\cong\mathcal{O}_{X\left(  \Sigma\right)
}\left(  -\sum_{v\in\Sigma\left(  1\right)  }D_{v}\right)
\]

Section \ref{Divisors on toric varieties} shows how to represent Weil and
Cartier divisors, the Chow group $A_{n-1}\left(  Y\right)  $ of Weil divisors
modulo linear equivalence on a toric variety $Y=X\left(  \Sigma\right)  $ and
the Picard group $\operatorname*{Pic}\left(  Y\right)  $. Classes in
$A_{n-1}\left(  Y\right)  $ can be represented by torus invariant Weil
divisors via the exact sequence%
\[%
\begin{tabular}
[c]{lllllll}%
$0\rightarrow$ & $M$ & $\overset{A}{\rightarrow}$ & $\mathbb{Z}^{\Sigma\left(
1\right)  }$ & $\rightarrow$ & $A_{n-1}\left(  Y\right)  $ & $\rightarrow0$%
\end{tabular}
\]
where the rows of $A$ are formed by the minimal lattice generators of the rays.

In Section \ref{Sec projective toric varieties} we
\index{projective toric variety}%
describe the correspondence of integral polytopes in $M_{\mathbb{R}}$ and
projective toric varieties. To an integral polytope $\Delta\subset
M_{\mathbb{R}}$ one can associate the graded ring%
\[%
\begin{tabular}
[c]{ll}%
$S\left(  \Delta\right)  =\mathbb{C}\left[  t^{k}x^{m}\mid m\in k\Delta
\right]  $ & $\deg t^{k}x^{m}=k$%
\end{tabular}
\]
with $k\Delta=\left\{  km\mid m\in\Delta\right\}  $ and $t^{k}x^{m}\cdot
t^{l}x^{m^{\prime}}=t^{k+l}x^{m+m^{\prime}}$, and hence the projective toric
variety $\mathbb{P}\left(  \Delta\right)  =\operatorname*{Proj}\left(
S\left(  \Delta\right)  \right)  $. Consider for any face $F$ of $\Delta$ the
cone of linear forms $w\in N_{\mathbb{R}}$, which take their minimum on
$\Delta$ at the points of $F$. These cones form a fan, the normal fan
$\Sigma=\operatorname*{NF}\left(  \Delta\right)  $ of $\Delta$. If $0$ is in
the interior of $\Delta$, then $\operatorname*{NF}\left(  \Delta\right)  $ is
the fan formed by the cones over the faces of the dual polytope%
\[
\Delta^{\ast}=\left\{  n\in N_{\mathbb{R}}\mid\left\langle m,n\right\rangle
\geq-1\ \forall m\in\Delta\right\}
\]
of $\Delta$. Furthermore, $\Delta$ defines a divisor on $X\left(
\Sigma\right)  $
\[
D_{\Delta}=\sum_{r\in\Sigma\left(  1\right)  }-\min_{m\in\Delta}\left\langle
m,\hat{r}\right\rangle D_{r}%
\]
Then as a toric variety $\mathbb{P}\left(  \Delta\right)  \cong X\left(
\Sigma\right)  $ with choice of an ample line bundle $\mathcal{O}%
_{\mathbb{P}\left(  \Delta\right)  }\left(  1\right)  \cong\mathcal{O}%
_{X\left(  \Sigma\right)  }\left(  D_{\Delta}\right)  $.

The Cox ring of a toric variety is explained in Section
\ref{Sec Cox ring of a toric variety} and homogeneous coordinate presentations
of toric varieties in Section \ref{1homogeneouscoordinate}. The Cox ring of a
toric variety $Y=X\left(  \Sigma\right)  $ is the polynomial ring
$S=\mathbb{C}\left[  y_{r}\mid r\in\Sigma\left(  1\right)  \right]  $ graded
via the above presentation of the Chow group considering monomials in $S$ as
elements of $\mathbb{Z}^{\Sigma\left(  1\right)  }$. In an analogous way to
the representation of projective space as%
\[
\mathbb{P}^{n}=\left(  \mathbb{C}^{n+1}-V\left(  \left\langle y_{0}%
,...,y_{n}\right\rangle \right)  \right)  /\mathbb{C}^{\ast}%
\]
there is a similar description of toric varieties as a categorial quotient%
\[
X\left(  \Sigma\right)  =\left(  \mathbb{C}^{\Sigma\left(  1\right)
}-V\left(  B\left(  \Sigma\right)  \right)  \right)  //G\left(  \Sigma\right)
\]
with some irrelevant ideal $B\left(  \Sigma\right)  \subset S$ and the action
of%
\[
G\left(  \Sigma\right)  =\operatorname*{Hom}\nolimits_{\mathbb{Z}}\left(
A_{n-1}\left(  Y\right)  ,\mathbb{C}^{\ast}\right)
\]
induced by the above sequence.

The application of the Cox ring to represent subvarieties and sheaves is
treated in Sections \ref{Sec Global sections a cox monomials} and
\ref{Sec Homogeneous coordinate representation of subvarieties and sheaves}.
For example the vector space of global sections of the reflexive sheaf of
sections $\mathcal{O}_{X\left(  \Sigma\right)  }\left(  D\right)  $ of a Weil
divisor $D$ on $Y$ is isomorphic to the degree $\left[  D\right]  $-part of
the Cox ring.

Section \ref{Sec Kaehler cone Mori cone} gives an algorithm to compute the
Mori cone $\overline{NE}\left(  Y\right)  _{\mathbb{R}}\subset A_{1}\left(
Y\right)  \otimes\mathbb{R}$ of effective $1$-cycles for a simplicial toric
variety $Y$.

The one-to-one correspondence of Gorenstein toric Fano varieties
$\mathbb{P}\left(  \Delta\right)  $ of dimension $n$, polarized by
$-K_{\mathbb{P}\left(  \Delta\right)  }$ and reflexive polytopes
$\Delta\subset\mathbb{Z}^{n}\otimes\mathbb{R}$ is treated in Section
\ref{Sec toric fano varieties}. The involution of Gorenstein toric Fano
varieties induced by duality of reflexive polytopes is the foundation of
Batyrev%
\'{}%
s mirror construction for anticanonical hypersurfaces.

In order to understand, which torus invariant deformations represented by Cox
Laurent monomials are trivial, we have to describe the automorphism group. If
$Y=X\left(  \Sigma\right)  $ is simplicial, then the connected component of
the identity of $\operatorname*{Aut}\left(  Y\right)  $ is generated by
automorphisms induced by the torus in $Y$ and by the so called root
automorphisms. Represented as Cox Laurent monomials a root automorphism is a
degree $0$ Cox Laurent monomial in $\mathbb{Z}^{\Sigma\left(  1\right)  }$ of
the form%
\[
\frac{\prod_{r\in\Sigma\left(  1\right)  }y_{r}^{a_{r}}}{y_{v}}%
\]
with relative prime numerator and denominator, and the corresponding
$1$-parameter family of automorphisms is%
\[%
\begin{tabular}
[c]{cccc}%
$y_{v}\mapsto$ & $y_{v}+\lambda\prod_{s\in\Sigma\left(  1\right)  -\left\{
v\right\}  }y_{s}^{a_{r}}$ &  & \\
$y_{r}\mapsto$ & \multicolumn{1}{l}{$y_{r}$} & for & $r\in\Sigma\left(
1\right)  -\left\{  v\right\}  $%
\end{tabular}
\]

Toric Mori theory will be used to relate Calabi-Yau degenerations to
orbifolding mirror families by relating the polarizing toric Fano variety of
the degeneration to a different birational model. Section
\ref{Sec toric mori theory} gives an overview of Reid%
\'{}%
s toric interpretation of Mori theory, i.e., cone theorem, contraction
theorem, existence and termination of flips and the minimal model program.
Given a finite set $\mathcal{R}$ of $1$-dimensional rational cones of a
projective fan, the set of all closures of K\"{a}hler cones
$\operatorname{cpl}\left(  \Sigma\right)  $ of projective simplicial fans
$\Sigma$ with $\Sigma\left(  1\right)  \subset\mathcal{R}$ fit together as
$\left(  \left\vert \mathcal{R}\right\vert -n\right)  $-dimensional cones of a
fan in $A_{n-1}\left(  \mathcal{R}\right)  _{\mathbb{R}}\cong\mathbb{R}%
^{\mathcal{R}}/M_{\mathbb{R}}$. To justify the notation $A_{n-1}\left(
\mathcal{R}\right)  _{\mathbb{R}}$, observe that the presentation of the Chow
group of a toric variety $X\left(  \Sigma\right)  $ only depends on the
$1$-dimensional cones of the fan $\Sigma$. The fan generated by the maximal
cones $\operatorname{cpl}\left(  \Sigma\right)  $ is called the
Gelfand-Kapranov-Zelevinsky decomposition associated to $\mathcal{R}$ and can
be extended to a complete fan, called the secondary fan $\Sigma\left(
\mathcal{R}\right)  $. We explain an algorithm to compute the secondary fan
via triangulations of marked polytopes.

The next Section \ref{Sec Groebner basics} gives a short account of
Gr\"{o}bner bases, weight orderings and the Mora algorithm computing standard
bases in the non global setting. The concept of Gr\"{o}bner bases plays an
important role both for the theory of flat degenerations and for computing
tropical varieties, so also for the tropical mirror construction. With regard
to flat degenerations see also the remarks about Section
\ref{Groebnerbasesandflatness} below. Gr\"{o}bner basis theory is the
algorithmic object connecting tropical geometry to degenerations and mirror symmetry.

\textbf{Section \ref{Sec mirror constructions to generalize}.} In this section
we summarize the mirror constructions, which will be generalized in a common
framework by the tropical mirror construction.

We begin in Section \ref{Sec Batyrev} with a short overview of the mirror
construction given by Batyrev for anticanonical hypersurfaces in toric Fano
varieties. Reflexive polytopes $\Delta\subset M_{\mathbb{R}}$ correspond to
Gorenstein toric Fano varieties $Y=\mathbb{P}\left(  \Delta\right)  $
polarized by $-K_{\mathbb{P}\left(  \Delta\right)  }$. A general element of
$\left\vert -K_{\mathbb{P}\left(  \Delta\right)  }\right\vert $ is a
Calabi-Yau hypersurface in $Y$. Duality is an involution of the set of
reflexive polytopes. Batyrev proves that general elements of $\left\vert
-K_{\mathbb{P}\left(  \Delta\right)  }\right\vert $ and $\left\vert
-K_{\mathbb{P}\left(  \Delta^{\ast}\right)  }\right\vert $ form a mirror pair
in the sense of stringy Hodge numbers. In the original approach, Batyrev
constructs, via maximal projective subdivisions of the fan of $Y$, a partial
crepant resolution of the hypersurface. A maximal projective subdivision of
the normal fan $\Sigma$ of $\Delta$ is a simplicial refinement $\bar{\Sigma}$
of $\Sigma$ defining a projective toric variety $X\left(  \bar{\Sigma}\right)
$ with the property that the non-zero lattice points of $\Delta^{\ast}$ span
the $1$-dimensional cones of $\bar{\Sigma}$.

Batyrev%
\'{}%
s construction for hypersurfaces has a generalization to the case of complete
intersections given by nef partitions of reflexive polytopes. This mirror
construction was introduced by Borisov and is explained in Section
\ref{Sec Batyrev and Borisov mirror construction}. Let $\Delta\subset
M_{\mathbb{R}}$ be a reflexive polytope and $\Sigma=\operatorname*{NF}\left(
\Delta\right)  $ its normal fan. Let%

\[
\Sigma\left(  1\right)  =I_{1}\cup...\cup I_{c}%
\]
be a disjoint union and suppose that the corresponding divisors $E_{j}%
=\sum_{v\in I_{j}}D_{v}$ are Cartier, spanned by global sections, and let
$\Delta_{j}\subset M_{\mathbb{R}}$ be the polytope of sections of $E_{j}$.
Note that $\sum_{j=1}^{c}E_{j}=-K_{Y}$. With $\nabla_{j}%
=\operatorname*{convexhull}\left\{  \left\{  0\right\}  \cup I_{j}\right\}  $
the Minkowski sum%
\[
\nabla_{BB}=\nabla_{1}+...+\nabla_{c}%
\]
is again a reflexive polytope with $\nabla_{BB}^{\ast}%
=\operatorname*{convexhull}\left(  \Delta_{1}\cup...\cup\Delta_{c}\right)  $.
Let $\Sigma^{\circ}=\operatorname*{NF}\left(  \nabla_{BB}\right)  $, let%
\[
\Sigma^{\circ}\left(  1\right)  =J_{1}\cup...\cup J_{c}%
\]
be the disjoint union corresponding to the partition $\operatorname*{vertices}%
\left(  \nabla_{BB}^{\ast}\right)  \cap\Delta_{j}$ of the vertices of
$\nabla_{BB}^{\ast}$ and $E_{j}^{\circ}=\sum_{v\in J_{j}}D_{v}^{\circ}$. Then
$X$ in $Y=\mathbb{P}\left(  \Delta\right)  $ given by general sections of
$\mathcal{O}\left(  E_{1}\right)  ,...,\mathcal{O}\left(  E_{c}\right)  $ and
$X^{\circ}$ in $Y^{\circ}=\mathbb{P}\left(  \nabla_{BB}\right)  $ defined by
general sections of $\mathcal{O}\left(  E_{1}^{\circ}\right)  ,...,\mathcal{O}%
\left(  E_{c}^{\circ}\right)  $ form a mirror pair with respect to stringy
Hodge numbers.

Section \ref{Roedlandexample} introduces the Greene-Plesser orbifolding mirror
family given by R\o dland for the general Calabi-Yau threefold $X$ of degree
$14$ in $\mathbb{P}^{6}$ defined by the Pfaffians of a general linear skew
symmetric map%
\[
7\mathcal{O}\left(  -1\right)  \rightarrow7\mathcal{O}%
\]
The mirror is given as a general element of the $1=h^{1,1}\left(  X\right)
$-parameter family with fibers in a $\mathbb{Z}_{7}$-quotient of
$\mathbb{P}^{6}$ defined by the Pfaffians%
\[
\left(
\begin{array}
[c]{ccccccc}%
0 & tx_{1} & x_{2} & 0 & 0 & -x_{5} & -tx_{6}\\
-tx_{1} & 0 & tx_{3} & x_{4} & 0 & 0 & -x_{0}\\
-x_{2} & -tx_{3} & 0 & tx_{4} & x_{6} & 0 & 0\\
0 & -x_{4} & -tx_{4} & 0 & tx_{0} & x_{1} & 0\\
0 & 0 & -x_{6} & -tx_{0} & 0 & tx_{2} & x_{3}\\
x_{5} & 0 & 0 & -x_{1} & -tx_{2} & 0 & tx_{4}\\
tx_{6} & x_{0} & 0 & 0 & -x_{3} & -tx_{4} & 0
\end{array}
\right)
\]
in $\mathbb{C}\left[  t\right]  \otimes\mathbb{C}\left[  x_{0},...,x_{6}%
\right]  $, i.e., by the square roots of the $6\times6$ diagonal minors.

\textbf{Section \ref{Sec Degenerations and mirror symmetry}. }The next main
section introduces examples of monomial degenerations of Calabi-Yau varieties,
which will serve as an input for the tropical mirror construction. Section
\ref{Degenerationcompleteintersection} defines the natural monomial
degenerations associated to hypersurfaces given by reflexive polyhedra and to
complete intersections given by nef partitions. Let $\Delta\subset
M_{\mathbb{R}}$ be a reflexive polytope, $Y=\mathbb{P}\left(  \Delta\right)  $
a toric Fano variety with Cox ring $S$, $\Sigma=\operatorname*{NF}\left(
\Delta\right)  $ and $\Sigma\left(  1\right)  =I_{1}\cup...\cup I_{c}$ a nef
partition. Then we obtain a degeneration given by%
\begin{gather*}
I=\left\langle t\cdot g_{j}+m_{j}\mid j=1,...,c\right\rangle \subset
\mathbb{C}\left[  t\right]  \otimes S\\
\text{with }m_{j}=\prod_{v\in I_{j}}y_{v}%
\end{gather*}
and monomial special fiber%
\[
I_{0}=\left\langle m_{j}\mid j=1,...,c\right\rangle
\]
Section \ref{Sec degenerations of pfaffian calabi-yau varieties} gives
monomial degenerations of some non complete intersection Pfaffian Calabi-Yau
varieties. A monomial degeneration of a general Pfaffian elliptic curve in
$\mathbb{P}^{4}$ defined by the Pfaffians of a general linear skew symmetric
map $A:5\mathcal{O}_{\mathbb{P}^{4}}\left(  -1\right)  \rightarrow
5\mathcal{O}_{\mathbb{P}^{4}}$ is given by the Pfaffians of%
\[
t\cdot A+\left(
\begin{array}
[c]{ccccc}%
0 & 0 & x_{1} & -x_{4} & 0\\
0 & 0 & 0 & x_{2} & -x_{0}\\
-x_{1} & 0 & 0 & 0 & x_{3}\\
x_{4} & -x_{2} & 0 & 0 & 0\\
0 & x_{0} & -x_{3} & 0 & 0
\end{array}
\right)
\]
If $A:7\mathcal{O}_{\mathbb{P}^{6}}\left(  -1\right)  \rightarrow
7\mathcal{O}_{\mathbb{P}^{6}}$ is a general skew symmetric map then the
Pfaffians of%
\[
t\cdot A+\left(
\begin{array}
[c]{ccccccc}%
0 & 0 & x_{2} & 0 & 0 & -x_{5} & 0\\
0 & 0 & 0 & x_{4} & 0 & 0 & -x_{0}\\
-x_{2} & 0 & 0 & 0 & x_{6} & 0 & 0\\
0 & -x_{4} & 0 & 0 & 0 & x_{1} & 0\\
0 & 0 & -x_{6} & 0 & 0 & 0 & x_{3}\\
x_{5} & 0 & 0 & -x_{1} & 0 & 0 & 0\\
0 & x_{0} & 0 & 0 & -x_{3} & 0 & 0
\end{array}
\right)
\]
define a monomial degeneration of a general degree $14$ Pfaffian Calabi-Yau
threefold. In the same way there is a monomial degeneration of a general
Calabi-Yau threefold of degree $13$ in $\mathbb{P}^{6}$ defined by the
Pfaffians of a general skew symmetric map $\mathcal{O}\left(  -2\right)
\oplus4\mathcal{O}\left(  -1\right)  \rightarrow\mathcal{O}\left(  1\right)
\oplus4\mathcal{O}$.

\textbf{Section \ref{Sec tropical geometry ingredients}. }The next main
section introduces fundamental facts from tropical geometry used to formulate
the mirror construction. Section \ref{Sec amoebas} defines the amoeba of a
subvariety of a torus as its image under the map
\begin{align*}
\log_{t}  &  :\left(  \mathbb{C}^{\ast}\right)  ^{n}\rightarrow\mathbb{R}%
^{n}\\
\left(  z_{1},...,z_{n}\right)   &  \mapsto\left(  \log_{t}\left\vert
z_{1}\right\vert ,...,\log_{t}\left\vert z_{n}\right\vert \right)
\end{align*}
Let $K$ be the metric completion of the field of Puisseux series
$\overline{\mathbb{C}\left(  t\right)  }$ with respect to the norm $\left\Vert
f\right\Vert =e^{-val\left(  f\right)  }$, where $val\left(  f\right)  $
denotes the exponent of the lowest weight term of $f$, and let $I\subset
K\left[  x_{1},...,x_{n}\right]  $ be an ideal. Section
\ref{1NonArchimedianAmoebas} relates the limit for $t\rightarrow\infty$ of the
amoeba given by $I$ to the non-Archimedian amoeba. This is the image under
\begin{gather*}
\operatorname*{val}\nolimits_{-}:\left(  K^{\ast}\right)  ^{n+1}%
\rightarrow\mathbb{R}^{n+1}\\
\left(  f_{1},...,f_{n}\right)  \mapsto\left(  -val\left(  f_{1}\right)
,...,-val\left(  f_{n}\right)  \right)
\end{gather*}
of the vanishing locus $V_{K}\left(  I\right)  $ of $I$ over $K$. The
non-Archimedian amoeba of $I$ is also called the tropical variety
$\operatorname*{tropvar}\left(  I\right)  $ of $I$. Note that here we take the
negative of the vanishing order in the definition of the valuation map, as in
the context of tropical geometry one usually considers the point of view of
the $\left(  \max,+\right)  $ algebra.

Section \ref{Sec tropical varieties} lists the basic properties of tropical
varieties, in particular their characterization as the set of weight vectors
$w\in\mathbb{R}^{n}$ such that $\operatorname*{in}_{w}\left(  I\right)  $
contains no monomial. In Section \ref{Sec tropical prevarieties} we recall the
algebraic description of tropical varieties. Given a polynomial $f\in K\left[
x_{1},...,x_{n}\right]  $ we replace $+$ by $\max$, multiplication by $+$ and
the coefficients $c$ by $-val\left(  c\right)  $ hence associating to $f$ a
piecewise linear function $\operatorname*{trop}\left(  f\right)  $. Then
$\operatorname*{tropvar}\left(  \left\langle f\right\rangle \right)  $ is the
non-differentiability locus of $\operatorname*{trop}\left(  f\right)  $. In
the same way $\operatorname*{tropvar}\left(  I\right)  $ is the intersection
$T\left(  \operatorname*{trop}\left(  I\right)  \right)  $ of the
non-differentiability loci of all $\operatorname*{trop}\left(  f\right)  $ for
$f\in I$. In Section \ref{1tropicalvarietiesandtheBergmanfan} we relate the
tropical variety of $I\subset\mathbb{C}\left[  t,x_{1},...,x_{n}\right]  $ to
a complex $BC_{-}\left(  I\right)  $, defined via its underlying set, which is
the set of those points on the unit sphere that are the limit of projections
of points of $\log\left(  V\left(  I\right)  \right)  $ on an expanding sphere
$jS^{n}$ for $j\rightarrow\infty$. The fan $BF_{-}\left(  I\right)  $ is
defined as the fan over $BC_{-}\left(  I\right)  $. Note that the fan
$BF_{-}\left(  I\right)  $ is known in the literature as the Bergman fan,
which differs by reflection at the origin from the Bergman fan $BF\left(
I\right)  $ as we defined above. The relation between $\operatorname*{tropvar}%
\left(  I\right)  $ and $BC_{-}\left(  I\right)  $ is given by stereographic
projection $\pi_{-}$ of the lower half unit sphere from $0$ to the plane
$\left\{  w_{t}=-1\right\}  =\mathbb{R}^{n}$ of $t$-weight $-1$.

In the definition of the amoeba, of the non-Archimedian amoeba, of the
tropical variety, of the non-differentiability locus of $\operatorname*{trop}%
\left(  f\right)  $ and in the definition of $BC_{-}\left(  I\right)  $ and
$BF_{-}\left(  I\right)  $ we adopt the Gr\"{o}bner basis point of view,
looking at the maximal weight term and take $weight\left(  c\right)
=-val\left(  c\right)  $ for constants $c\in K$. From the point of view of
degenerations and local arc solutions of the total space of a degeneration at
the special fiber, it is more natural to consider the minimal weight term
combined with the definition $weight\left(  c\right)  =val\left(  c\right)  $
for constants $c\in K$. Summarizing, in our notation we have%
\begin{align*}
\operatorname*{val}\left(  V_{K}\left(  I\right)  \right)   &  =\pi\left(
BF\left(  I\right)  \cap S^{n}\cap\left\{  w_{t}>0\right\}  \right) \\
&  =-\lim_{t\rightarrow\infty}\left(  \log_{t}V\left(  I_{t}\right)  \right)
=-\operatorname*{val}\nolimits_{-}\left(  V_{K}\left(  I\right)  \right) \\
&  =-\operatorname*{tropvar}\left(  I\right)  =-T\left(  \operatorname*{trop}%
\left(  I\right)  \right) \\
&  =-\pi_{-}\left(  BF_{-}\left(  I\right)  \cap S^{n}\cap\left\{
w_{t}<0\right\}  \right)
\end{align*}
where $\pi$ is the stereographic projection of the upper half unit sphere from
$0$ to the plane $\left\{  w_{t}=1\right\}  =\mathbb{R}^{n}$ of $t$-weight $1$
and in the same way $\pi_{-}$ from the lower half unit sphere.

\textbf{Section \ref{Groebnerbasesandflatness}. }This section gives the
standard characterization of flatness via Gr\"{o}bner bases, e.g., a first
order degeneration $\mathfrak{X}$ defined by
\[
\left\langle f_{1}+t\cdot g_{1},...,f_{r}+t\cdot g_{r}\right\rangle \subset
R\otimes\mathbb{C}\left[  t\right]  /\left\langle t^{2}\right\rangle
\]
with special fiber given by $\left\langle f_{1},...,f_{r}\right\rangle $ is
flat if and only if any syzygy $\sum_{i}a_{i}f_{i}=0\in R$ lifts to a syzygy
between $f_{1}+tg_{1},...,f_{r}+tg_{r}$, i.e., there are $b_{i}\in R$ such
that%
\[
\sum_{i}\left(  a_{i}+tb_{i}\right)  \left(  f_{i}+tg_{i}\right)  =0\in
R\otimes k\left[  t\right]  /\left\langle t^{2}\right\rangle
\]

\textbf{Section \ref{Sec Computing the Bergman fan}.} We recall the definition
of the Gr\"{o}bner fan of an ideal introduced by Mora, its dual description
via state polytopes and the construction of multigraded Hilbert schemes.
Furthermore, we connect stability of the Hilbert point with state polytopes.
The existence of the multigraded Hilbert scheme shows that ideals in the Cox
ring provide the right framework to describe subvarieties in toric varieties.

Our main interest in the Gr\"{o}bner fan is the computation of tropical
varieties, so in Section \ref{sec computing the Bergman fan subsection} we
begin with a concept for computing the Bergman fan. Consider an ideal
$J\subset\mathbb{C}\left[  x_{1},...,x_{n}\right]  $ such that every weight
vector is equivalent to a non-negative weight vector, e.g., a homogeneous
ideal. Section \ref{Sec tropical representation of Groebner cones} introduces
the Gr\"{o}bner cone of weight vectors equivalent to a given global ordering
$>$ and the Gr\"{o}bner fan $GF\left(  J\right)  $. The maximal cones of the
fan $GF\left(  J\right)  $ correspond to the monomial initial ideals of $J$.
Section \ref{Sec computing the Groebner fan} gives a simple algorithm
terminating with the fan $GF\left(  J\right)  $. We take a cone $C$ in a non
complete subfan of $GF\left(  J\right)  $ and move into the complement of the
fan along an outer normal vector of a face, which appears in the fan only
once. Then we compute the corresponding Gr\"{o}bner cone. Note that this is
well suited for using Gr\"{o}bner walk algorithms.

The second part of Section \ref{Sec Computing the Bergman fan} deals with the
Hilbert scheme and stability. Generalizing step by step, we begin in Section
\ref{Sec Hilbert scheme and state polytope projective setup} with the setup of
homogeneous ideals $J$ with fixed Hilbert polynomial $P_{S/J}$ in
$S=\mathbb{C}\left[  x_{0},...,x_{n}\right]  $ with respect to the grading%
\[
0\rightarrow\mathbb{Z}^{n}\overset{A}{\rightarrow}\mathbb{Z}^{n+1}%
\overset{\deg}{\rightarrow}\mathbb{Z}\rightarrow0
\]
with%
\[
A=\left(
\begin{tabular}
[c]{lll}%
$1$ &  & \\
& $\ddots$ & \\
&  & $1$\\
$-1$ & $\cdots$ & $-1$%
\end{tabular}
\right)
\]
i.e., ideals in the Cox ring of $\mathbb{P}^{n}$. We recall the construction
of the Hilbert scheme, the state polytope and the characterization of
stability via the state polytope as given by Bayer and Morrison in \cite{BaMo
Initial Ideals and State Polytopes}.

After summarizing in Section \ref{Sec Linearizations} some facts on
$G$-linearizations of line bundles for an affine algebraic group $G$ acting
rationally on an algebraic variety, we generalize to the toric setting. In the
same way as for the case of subvarieties of $\mathbb{P}^{n}$, the key
ingredients are the Grassmann functor and the Hilbert functor explained in
Sections \ref{Sec Grassmann functor} and \ref{Sec Hilbert functor} as given by
Haiman and Sturmfels in \cite{HS Multigraded Hilbert schemes}. Let $k$ be a
commutative ring, $A\ $a set and $S$ a polynomial ring graded by a set $A$. If
$h:A\rightarrow\mathbb{N}$ is a function and $R$ is a $k$-algebra, then the
Hilbert functor $\mathbb{H}_{\left(  S,F\right)  }^{h}$ is defined via%
\[
\mathbb{H}_{\left(  S,F\right)  }^{h}\left(  R\right)  =\left\{  L\mid%
\begin{tabular}
[c]{l}%
$L\subset R\otimes S$ is an $F$-submodule $\text{with}$\\
$\left(  R\otimes S_{a}\right)  /L_{a}\text{ locally free of rank }h\left(
a\right)  $ $\forall a\in A$%
\end{tabular}
\ \right\}
\]
Here the notion of an $F$-submodule is defined via sets of operators
$F_{a,b}\subset\operatorname*{Hom}_{k}\left(  S_{a},S_{b}\right)  $. Under
appropriate conditions $\mathbb{H}_{\left(  S,F\right)  }^{h}$ is represented
by a closed subscheme of a projective Grassmann scheme. The key point is the
restriction to a finite set of degrees. In the case of the homogeneous
coordinate ring of $\mathbb{P}^{n}$ one can restrict to one degree.

Considering an example by Haiman and Sturmfels, Section
\ref{Sec multigraded hilbert scheme of admissable ideals} explains the
application of this construction to the Hilbert scheme of admissible ideals,
i.e., ideals with the property that $\left(  S/I\right)  _{a}=S_{a}/I_{a}$ is
a locally free $k$-module of finite rank for all $a\in A$. Note that this
setup is not directly applicable even to the homogeneous coordinate ring of
projective space. As a second example, Section
\ref{Sec classical hilbert functor} applies the above construction of
$\mathbb{H}_{\left(  S,F\right)  }^{h}$ to obtain the classical Hilbert scheme
via truncation $I_{\geq a}$ of ideals at an appropriate degree $a$.

The tangent space at $I\in\mathbb{H}_{\left(  S,F\right)  }^{h}\left(
k\right)  $ of the scheme representing $\mathbb{H}_{\left(  S,F\right)  }^{h}$
is described as $\operatorname*{Hom}_{S}\left(  I,S/I\right)  _{0}$ in Section
\ref{Sec tangent space and deformations}.

Sections \ref{Sec stanley decompositions} and \ref{Sec Multigraded regularity}
give an overview of Stanley filtrations and multigraded regularity as
introduced by Maclagan and Smith in \linebreak\cite{MaSm Multigraded
CastelnuovoMumford regularity} and \cite{MaSm Uniform bounds on multigraded
regularity}. If $Y$ is a smooth toric variety, $M$ is a finitely generated
$S$-module then for $m\in A_{n-1}\left(  Y\right)  $ the notion of
$m$-regularity is defined via local cohomology. The regularity of $M$ is the
set of all degrees $m\in A_{n-1}\left(  Y\right)  $ such that $M$ is
$m$-regular. In Section \ref{Sec Multigraded Hilbert schemes} the multigraded
Hilbert functor $R\mapsto\mathbb{H}_{Y}^{P}\left(  R\right)  $ associating to
a $\mathbb{C}$-algebra $R$ the set of ideal sheaves $\mathcal{J}$ of families
of subschemes $X\subset Y\times_{\mathbb{C}}\operatorname*{Spec}%
R\rightarrow\operatorname*{Spec}R$ with fixed multivariate Hilbert polynomial
$P$ is given. Using the above construction of $\mathbb{H}_{\left(  S,F\right)
}^{h}$, the functor $\mathbb{H}_{Y}^{P}$ is represented by a projective scheme
over $\mathbb{C}$. The finite set of degrees to represent $\mathbb{H}_{Y}^{P}$
as a subscheme of a Grassmann scheme can be computed algorithmically.

In Section \ref{Sec State polytope general toric} we introduce the state
polytope in the multigraded setting and characterize the sets of stable and
semistable points via the state polytope. If $I$ is an ideal in the Cox ring
of a smooth toric variety with Hilbert function $h$ and $\mathbb{H}$ the
corresponding Hilbert scheme, then the Hilbert point $H\left(  I\right)
\in\mathbb{H}$ is in the stable locus $\mathbb{H}^{s}$ if and only if $0$ is
in the interior of the state polytope $\operatorname*{State}\left(  I\right)
$ of $I$.

Finally, given a toric variety $Y$ defined by a fan $\Sigma\subset
N_{\mathbb{R}}$, we identify in Section \ref{1torichomogeneoussetting} the
weight vectors on the Cox ring $S$ of $Y$ with the vectors in $N_{\mathbb{R}}$
by dualizing the presentation of the Chow group. Hence, the Gr\"{o}bner fan of
an ideal in the Cox ring of $Y$ can be considered as a fan in $N_{\mathbb{R}}$.

\textbf{Section \ref{Sec Q-Gorenstein varieties fano polytopes}.} This section
considers toric Fano varieties $Y$ in the sense that some multiple of $-K_{Y}$
is an ample Cartier divisor. It explains how Fano polytopes $P\subset
N_{\mathbb{R}}$ represent $\mathbb{Q}$-Gorenstein toric Fano varieties defined
by the fan over the faces of $P$. This is the right category of toric Fano
varieties with respect to toric Mori theory.

\textbf{Section
\ref{Sec tropical mirror construction for complete intersections}.} Here we
formulate the tropical mirror construction for complete intersections. The
construction takes as an input the degenerations associated in Section
\ref{Degenerationcompleteintersection} to complete intersections defined by
nef partitions. So let $\Delta\subset M_{\mathbb{R}}$ be a reflexive polytope,
$Y=\mathbb{P}\left(  \Delta\right)  $ the corresponding toric Fano variety
with Cox ring $S$, presentation matrix $A$ of $A_{n-1}\left(  Y\right)  $ and
$\Sigma=\operatorname*{NF}\left(  \Delta\right)  \subset N_{\mathbb{R}}$. Let
$\Sigma\left(  1\right)  =I_{1}\cup...\cup I_{c}$ be a nef partition and
$\Delta_{j}\subset M_{\mathbb{R}}$ the polytopes of sections of the
corresponding divisors and $\nabla_{j}=\operatorname*{convexhull}\left\{
\left\{  0\right\}  \cup I_{j}\right\}  $. Let $\mathfrak{X}\subset
Y\times\operatorname*{Spec}\left(  \mathbb{C}\left[  t\right]  \right)  $ be
the associated degeneration as defined above by the ideal
\[
I=\left\langle f_{j}=t\cdot g_{j}+m_{j}\mid j=1,...,c\right\rangle
\subset\mathbb{C}\left[  t\right]  \otimes S
\]
and let $I_{0}\subset S$ be the ideal with minimal generators $m_{j}$. We
begin in Section \ref{degenerationtoriccompleteintersections} by exploring the
properties of these degenerations and describe in Section
\ref{Groebnerconeassocitatedtothespecialfiber ci} the special fiber
Gr\"{o}bner cone $C_{I_{0}}\left(  I\right)  \subset\mathbb{R}\oplus
N_{\mathbb{R}}$ and the special fiber polytope
\[
\nabla=C_{I_{0}}\left(  I\right)  \cap\left\{  w_{t}=1\right\}  \subset
N_{\mathbb{R}}%
\]
The special fiber cone $C_{I_{0}}\left(  I\right)  $ is cut out by the
half-space equations corresponding to first order deformations contributing to
the tangent vector of $\mathfrak{X}$ at the special fiber. We show that the
reflexive polytope $\nabla$ coincides with the Batyrev-Borisov mirror
polytope. Section \ref{Sec initial ideals of the faces ci} gives an explicit
description of the initial ideals of the faces of $\nabla$ and Section
\ref{Sec dual complex ci} introduces the map%
\[
\operatorname*{dual}:\operatorname*{Poset}\left(  \nabla\right)
\rightarrow\operatorname*{Poset}\left(  \nabla^{\ast}\right)
\]
between the complexes of faces of $\nabla$ and $\nabla^{\ast}$ associating to
a face $F$ of $\nabla$ the convex hull of all first order deformations
appearing the initial ideal of $I$ with respect to $F$. So, if%
\[
\operatorname*{in}\nolimits_{F}\left(  f_{j}\right)  =t\sum_{m\in G_{j}\left(
F\right)  }c_{m}m+m_{j}%
\]
then%
\[
\operatorname*{dual}\left(  F\right)  =\operatorname*{convexhull}\left(
A^{-1}\left(  \frac{m}{m_{j}}\right)  \mid m\in G_{j}\left(  F\right)  ,\text{
}j=1,...,c\right)  \subset M_{\mathbb{R}}%
\]
Note that first order deformations correspond to Cox Laurent monomials of
degree $0$ hence via the presentation matrix $A$ of the Chow group to elements
of $M$. Indeed we show that $\operatorname*{dual}\left(  F\right)  =F^{\ast}$
is the face of $\nabla^{\ast}$ dual to $F$.

By considering the faces of $\nabla$ which correspond to cones in the Bergman
fan we obtain in Section \ref{Sec Bergman subcomplex ci} the Bergman
subcomplex $B\left(  I\right)  \subset\operatorname*{Poset}\left(
\nabla\right)  $ of the poset of faces of $\nabla$. In Section
\ref{Sec mirror complex complete intersection} we give an inclusion reversing
map from $B\left(  I\right)  $ to the complex of faces of $\Delta$%
\[
\mu:B\left(  I\right)  \rightarrow\operatorname*{Poset}\left(  \Delta\right)
\]
by taking the Minkowski sum over the faces of $\nabla^{\ast}$ corresponding to
deformations of the individual equations, i.e.,%
\[
\mu\left(  F\right)  =\sum_{j=1}^{c}\operatorname*{convexhull}\left(
A^{-1}\left(  \frac{m}{m_{j}}\right)  \mid m\in G_{j}\left(  F\right)
\right)
\]

We relate $\mu\left(  F\right)  $ and $\operatorname*{dual}\left(  F\right)  $
in Section \ref{Sec dual complex BI ci} via%
\[
\mu\left(  F\right)  =\sum_{j=1}^{c}\operatorname*{dual}\left(  F\right)
\cap\Delta_{j}%
\]

\begin{figure}
[h]
\begin{center}
\includegraphics[
height=2.0816in,
width=4.644in
]%
{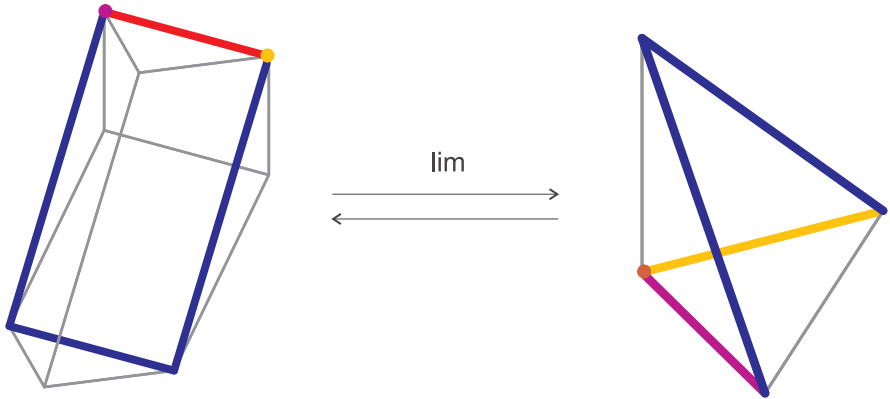}%
\caption{Limit correspondence for the monomial degeneration of the complete
intersection of two general quadrics in $\mathbb{P}^{3}$ and its mirror}%
\label{Fig lim 22}%
\end{center}
\end{figure}

In Section \ref{Sec intersection complex ci} we relate the maps $\lim$ and
$\mu$. If $F$ is a face of $B\left(  I\right)  $, then
\[
\lim\left(  F\right)  =V\left(  \left(  \mu\left(  F\right)  \right)  ^{\ast
}\right)
\]
is the toric stratum of $Y$ corresponding to the face $\mu\left(  F\right)  $
of $\Delta$. Figure \ref{Fig lim 22} shows the complexes $B\left(  I\right)  $
and $\lim\left(  B\left(  I\right)  \right)  $ and the polyhedra $\nabla$ and
$\Delta$ for the monomial degeneration of the complete intersection of two
general quadrics in $\mathbb{P}^{3}$ as given above.%

Section \ref{Sec covering complete intersection} gives a $c:1$ covering of the
complex $B\left(  I\right)  ^{\vee}$ by faces of $\operatorname*{dual}\left(
B\left(  I\right)  \right)  $. The covering is induced by associating to $F\in
B\left(  I\right)  $ the faces%
\[
\left\{  F^{\ast}\cap\Delta_{j}\mid j=1,...,c\right\}
\]
Note that this covering can have degenerate faces in the sense that
$\dim\left(  F^{\ast}\cap\Delta_{j}\right)  $ can be less than $\dim\left(
F^{\vee}\right)  =d-\dim\left(  F\right)  $, but the faces are always
non-empty. We give an algorithm computing this covering from the complex
$\operatorname*{dual}\left(  B\left(  I\right)  \right)  $.

Let $Y^{\circ}=\mathbb{P}\left(  \nabla\right)  $ and $A^{\circ}$ be the
presentation matrix of $A_{n-1}\left(  Y^{\circ}\right)  $. Applying in
Section \ref{Sec mirror covering ci} the covering algorithm on the mirror side
to the complex $\left(  \mu\left(  F\right)  \right)  ^{\ast}$, via the sheets
we obtain the ideal
\[
I^{\circ}=\left\langle t\cdot\sum_{\delta\in\left(  \mu\left(  B\left(
I\right)  \right)  \right)  ^{\ast}\cap\nabla_{j}}c_{\delta}\cdot\delta\left(
m_{j}^{\circ}\right)  +m_{j}^{\circ}\mid j=1,...,c\right\rangle \subset
S^{\circ}\otimes\mathbb{C}\left[  t\right]
\]
defining the mirror degeneration $\mathfrak{X}^{\circ}\subset\mathbb{P}\left(
\nabla\right)  \times\operatorname*{Spec}\mathbb{C}\left[  t\right]  $. Here
the monomials $m_{j}^{\circ}$ are the least common multiples of denominators
of the Cox Laurent monomials $A^{\circ}\left(  \delta\right)  $ for $\delta
\in\left(  \mu\left(  B\left(  I\right)  \right)  \right)  ^{\ast}\cap
\nabla_{j}$. Passing from $I^{\circ}$ to the ideal of the tropical mirror as
defined above by applying the deformations in $\left(  \mu\left(  B\left(
I\right)  \right)  \right)  ^{\ast}\cap N$ to the special fiber ideal%
\[%
{\textstyle\bigcap\nolimits_{F\in B\left(  I\right)  }}
\left\langle z_{G^{\ast}}\mid G\text{ a facet of }\nabla\text{ with }F\subset
G\right\rangle \subset S^{\circ}%
\]
is the toric analogue of saturation, also valid for non-simplicial toric
varieties. It does not change the geometry of the degeneration $\mathfrak{X}%
^{\circ}$ or the objects involved in the tropical mirror construction.

Along Section
\ref{Sec tropical mirror construction for complete intersections} we visualize
the objects introduced in the tropical mirror construction for the example of
the general complete intersection of two quadrics in $\mathbb{P}^{3}$.

In Section \ref{Sec Examples ci} we apply the tropical mirror construction to
some complete intersection examples. In particular, by considering a set of
Fermat deformations in order to relate $Y^{\circ}$ to a different birational
model, we obtain the Greene-Plesser type orbifolding mirror family of the
complete intersection of two cubics in $\mathbb{P}^{5}$ as given in \cite{LT
Lines on CalabiYau complete intersections mirror symmetry and PicardFuchs
equations}. Note that the text is computer generated by the implementation of
the tropical mirror construction in the Maple package \textsf{tropicalmirror}.

\textbf{Section \ref{Sec tropical mirror construction}.} This Section gives
the tropical mirror construction in its general form, as outlined above.

We begin in Section \ref{Sec concept of the tropical mirror construction} with
a summary and continue in the following sections by introducing fundamental
concepts used in the tropical mirror construction. Section
\ref{first order deformations and degree 0 Cox Laurent monomials} represents
torus invariant first order deformations of monomial ideals by lattice
monomials. Let $Y=X\left(  \Sigma\right)  $ be a toric variety given by the
fan $\Sigma$ in $N_{\mathbb{R}}$, let $S=\mathbb{C}\left[  y_{r}\mid
r\in\Sigma\left(  1\right)  \right]  $ be the Cox ring of $Y$ and
$I_{0}\subset S$ a monomial ideal. The space of degree $0$ first order
deformations $\operatorname*{Hom}\left(  I_{0},S/I_{0}\right)  _{0}$ has a
basis of $\left(  \mathbb{C}^{\ast}\right)  ^{\Sigma\left(  1\right)  }%
$-deformations. Any such homomorphism $\delta:I_{0}\rightarrow S/I_{0}$ is
representable by a degree $0$ Cox Laurent monomial $\frac{q_{1}}{q_{0}}$ with
relatively prime monomials $q_{0},q_{1}\in S$ via%
\[
\delta\left(  m\right)  =\left\{
\begin{tabular}
[c]{ll}%
$\frac{q_{1}}{q_{0}}\cdot m$ & if $q_{0}\mid m$\\
$0$ & otherwise
\end{tabular}
\right\}
\]
for minimal generators $m\in I_{0}$. Via the sequence%
\[
0\rightarrow M\overset{A}{\rightarrow}\mathbb{Z}^{\Sigma\left(  1\right)
}\overset{\deg}{\rightarrow}A_{n-1}\left(  X\left(  \Sigma\right)  \right)
\rightarrow0
\]
$\frac{q_{1}}{q_{0}}$ corresponds to a lattice monomial in
$M=\operatorname*{Hom}\left(  N,\mathbb{Z}\right)  $.

In Section
\ref{Monomial ideals in the Cox ring and the stratified toric primary decomposition}
we give a combinatorial description of the vanishing loci in $Y$ of reduced
monomial ideals $I_{0}\subset S$. Given a monomial $m\in I_{0}$, denote by
\[
\operatorname*{rays}\nolimits_{m}\left(  \Sigma\right)  =\left\{  r\in
\Sigma\left(  1\right)  \mid y_{r}\text{ divides }m\right\}
\]
the set of rays of $\Sigma$ such that $y_{r}$ divides $m$. We define the
stratified toric primary decomposition $SP\left(  I_{0}\right)  $ as the
complex, which has as faces of dimension $s$ the ideals $\left\langle
y_{r}\mid r\subset\sigma\right\rangle $ for all cones $\sigma\in\Sigma\ $of
dimension $n-s$ which contain a ray in $\operatorname*{rays}\nolimits_{m}%
\left(  \Sigma\right)  $ for all monomials $m\in I_{0}$.

Suppose $\Delta$ is a polytope with $\Sigma=\operatorname*{NF}\left(
\Delta\right)  $. Then $SP\left(  I_{0}\right)  $ is naturally isomorphic to
the complex $\operatorname*{Strata}\nolimits_{\Delta}\left(  I_{0}\right)  $
of strata of $I_{0}$. We define $\operatorname*{Strata}\nolimits_{\Delta
}\left(  I_{0}\right)  $ as the complex which has as faces of dimension $s$
those faces $F$ of $\Delta$ such that for all monomials $m\in I_{0}$ the set
\[
\left\{  G\mid G\text{ facet of }\Delta\text{ with }y_{G^{\ast}}\mid
m\right\}
\]
contains a facet $G$ with $F\subset G$. Suppose that the vanishing locus
$X_{0}$ of $I_{0}$ in $Y$ is equidimensional of dimension $d$. The complexes
$SP\left(  I_{0}\right)  \cong\operatorname*{Strata}\nolimits_{\Delta}\left(
I_{0}\right)  $ describe the vanishing locus of $I_{0}$ and they define the
ideal%
\begin{align*}
I_{0}^{\Sigma}  &  =%
{\textstyle\bigcap\nolimits_{J\in SP\left(  I_{0}\right)  _{d}}}
J\\
&  =%
{\textstyle\bigcap\nolimits_{F\in\operatorname*{Strata}\nolimits_{\Delta
}\left(  I_{0}\right)  _{d}}}
\left\langle y_{G^{\ast}}\mid G\text{ a facet of }\Delta\text{ with }F\subset
G\right\rangle
\end{align*}
naturally associated to $X_{0}$. Passing from $I_{0}$ to $I_{0}^{\Sigma}$ is
the toric analogue of saturation. Note that we do not assume $Y$ to be simplicial.

If $\Delta$ is a simplex and $I_{0}$ is a Stanley-Reisner ideal given by a
simplicial subcomplex $Z$ of the complex of cones of $\Sigma
=\operatorname*{NF}\left(  \Delta\right)  $, then we relate
$\operatorname*{Strata}\nolimits_{\Delta}\left(  I_{0}\right)  $ to $Z$ via
the map associating to a face $F\in\operatorname*{Strata}\nolimits_{\Delta
}\left(  I_{0}\right)  $ the hull of the rays of $\Sigma$ not contained in
$\operatorname*{hull}\left(  F^{\ast}\right)  $.

In Section \ref{Sec locally relevant deformations} we introduce the notion of
locally irrelevant deformations. Let $I_{0}\subset S$ be a reduced monomial
ideal defining $X_{0}\subset Y$, $X_{i}$ a stratum of $X_{0}$ and
$\mathfrak{X}$ a first order deformation of $X_{0}$. Then $\mathfrak{X}$ is
called locally irrelevant at $X_{i}$ if there is a formal analytic open
neighborhood $\tilde{U}\subset Y$ of $X_{i}$ and an isomorphism%
\[
\left(  \tilde{U}\cap X_{0}\right)  \times\operatorname*{Spec}\left(
\mathbb{C}\left[  t\right]  /\left\langle t^{2}\right\rangle \right)
\cong\mathfrak{X}\cap\left(  \tilde{U}\times\operatorname*{Spec}\left(
\mathbb{C}\left[  t\right]  /\left\langle t^{2}\right\rangle \right)  \right)
\]
extending $X_{i}\times\operatorname*{Spec}\left(  \mathbb{C}\left[  t\right]
/\left\langle t^{2}\right\rangle \right)  \subset\mathfrak{X}$.

In Section \ref{genericy condition} we give the setup for the tropical mirror
construction. Consider $N\cong\mathbb{Z}^{n}$, $M=\operatorname*{Hom}\left(
N,\mathbb{Z}\right)  $, $P$ a Fano polytope, $\Sigma=\Sigma\left(  P\right)  $
the fan over the faces of $P$ and $Y=X\left(  \Sigma\right)  $ with Cox ring
$S$. Let $I_{0}\subset S$ be a reduced monomial ideal with $I_{0}%
=I_{0}^{\Sigma}$ and equidimensional vanishing locus. Let $\mathfrak{X}\subset
Y\times\operatorname*{Spec}\mathbb{C}\left[  \left[  t\right]  \right]  $ be a
degeneration of Calabi-Yau varieties of codimension $c$, which is given by
$I\subset S\otimes\mathbb{C}\left[  t\right]  $ and with special fiber $X_{0}$
defined by $I_{0}$.

We give the conditions assumed to be satisfied for the input degeneration
$\mathfrak{X}$. Formulated in an explicit and testable form, these conditions are:

\begin{enumerate}
\item $C_{I_{0}}\left(  I\right)  \cap\left\{  w_{t}=0\right\}  =\left\{
0\right\}  $

\item $C_{I_{0}}\left(  I\right)  $ is the cone defined by the half-space
equations corresponding to the torus invariant first order deformations
appearing in the reduced standard basis of $I$ in $S\times\mathbb{C}\left[
t\right]  /\left\langle t^{2}\right\rangle $ with respect to a monomial
ordering in the interior of $C_{I_{0}}\left(  I\right)  $.

All lattice points of $F^{\ast}$ appear as deformations in $I$.

\item $\nabla^{\ast}\subset\Delta$, which is equivalent to the condition that
any first order deformation appearing in $I$ is also a deformation of the
anticanonical Calabi-Yau hypersurface in $Y$.

\item Any facet of $\operatorname*{Strata}\nolimits_{\Delta}\left(
I_{0}\right)  $ is contained in precisely $c$ facets of $\Delta$.

\item Any facet of $B\left(  I\right)  $ is contained in precisely $c$ facets
of $\nabla$.
\end{enumerate}

\noindent An interpretation of these conditions with respect to the geometry
of $\mathfrak{X}$ is given. We can satisfy requirement

\begin{enumerate}
\item via a condition on the position of the Hilbert point of $I_{0}$ with
respect to the state polytope of the general fiber,

\item via a genericity condition on the tangent vector with respect to the
tangent space of the component of the Hilbert scheme containing $\mathfrak{X}$,

\item via the condition that $\mathcal{O}_{X_{0}}$ has a resolution%
\[
0\rightarrow\mathcal{O}_{Y}\left(  -K_{Y}\right)  \rightarrow...\rightarrow
\mathcal{F}_{1}\rightarrow\mathcal{O}_{Y}\rightarrow\mathcal{O}_{X_{0}%
}\rightarrow0
\]
with direct sums $\mathcal{F}_{j}=%
{\textstyle\bigoplus\nolimits_{i}}
\mathcal{O}_{Y}\left(  D_{ji}\right)  $,

\item via the components of $X_{0}$ being given by $c$ linear equations,

\item via a condition on the locally relevant deformations of $\mathfrak{X}$
at the zero dimensional strata of $X_{0}$.
\end{enumerate}

In Sections \ref{Sec the groebner cone associated to the special fiber}%
-\ref{Sec first order mirror degeneration} we formulate the tropical mirror
construction in the general setting as already outlined above. Section
\ref{Sec the groebner cone associated to the special fiber} describes the
special fiber Gr\"{o}bner cone $C_{I_{0}}\left(  I\right)  \subset
\mathbb{R}\oplus N_{\mathbb{R}}$ and the special fiber polytope%
\[
\nabla=C_{I_{0}}\left(  I\right)  \cap\left\{  w_{t}=1\right\}  \subset
N_{\mathbb{R}}%
\]
The polytope $\nabla^{\ast}\subset M_{\mathbb{R}}$ is a Fano polytope, so the
Fan $\Sigma^{\circ}=\Sigma\left(  \nabla^{\ast}\right)  $ over the faces of
$\nabla^{\ast}$ defines a toric Fano variety $Y^{\circ}=X\left(  \Sigma
^{\circ}\right)  $ with Cox ring $S^{\circ}$.

Let%
\[
0\rightarrow M\overset{A}{\rightarrow}\mathbb{Z}^{\Sigma\left(  1\right)
}\overset{\deg}{\rightarrow}A_{n-1}\left(  X\left(  \Sigma\right)  \right)
\rightarrow0
\]
be the presentation of the Chow group of divisors of $X\left(  \Sigma\right)
$. Consider a face $F$ of $\nabla^{\ast}$. The initial forms of the reduced
standard basis of $I$ in $S\otimes\mathbb{C}\left[  t\right]  /\left\langle
t^{2}\right\rangle $ with respect to a monomial ordering in the interior of
$\nabla$ involve a minimal generator of $I_{0}$. Dividing the non special
fiber monomials by the special fiber monomial of the initial forms and
applying $A^{-1}$ we obtain a set lattice monomials. In analogy to complete
intersections associating to $F$ of $\nabla$ the convex hull of these lattice
monomials, we define in Section \ref{Sec dual complex general setting} the map%
\[
\operatorname*{dual}:\operatorname*{Poset}\left(  \nabla\right)
\rightarrow\operatorname*{Poset}\left(  \nabla^{\ast}\right)
\]
We observe that for $F\in\operatorname*{Poset}\left(  \nabla\right)  $%
\[
\operatorname*{dual}\left(  F\right)  =F^{\ast}%
\]
and all lattice points appear in the initial ideal. The non-special fiber
monomials of the initial forms decompose into characters of the torus $\left(
\mathbb{C}^{\ast}\right)  ^{\Sigma\left(  1\right)  }$, which are just the Cox
Laurent monomials associated to the lattice points of $F^{\ast}$. The
characters correspond to deformations of $X_{0}$ in $\operatorname*{Hom}%
\left(  I_{0},S/I_{0}\right)  _{0}$.

The complex $\operatorname*{dual}\left(  B\left(  I\right)  \right)  $ can be
seen as a polyhedral representation of the structure of the ideal $I$, as
described by structure theorems like the Koszul resolution for complete
intersections or the Buchsbaum-Eisenbud theorem for codimension $3$
subcanonical varieties.

We define in Section \ref{Sec Bergman subcomplex of Nabla general setting} the
special fiber Bergman complex%
\[
B\left(  I\right)  =\left(  BF\left(  I\right)  \cap\operatorname*{Poset}%
\left(  C_{I_{0}}\left(  I\right)  \right)  \right)  \cap\left\{
w_{t}=1\right\}  \subset\operatorname*{Poset}\left(  \nabla\right)
\]
of those faces $F$ of $\nabla$ such that $\operatorname*{in}_{F}\left(
I\right)  $ does not contain a monomial. Note that the Bergman fan $BF\left(
I\right)  $ contains more information than the combinatorial objects derived
from $B\left(  I\right)  $.

Section \ref{Sec covering structure in dual B} explores the covering structure
in $\operatorname*{dual}\left(  B\left(  I\right)  \right)  $ over $B\left(
I\right)  ^{\vee}$, generalizing the $c:1$ trivial covering in the case of
complete intersection.

In Section \ref{Sec limit map general setting} we describe the limit map
$\lim:B\left(  I\right)  \rightarrow\operatorname*{Strata}\left(
X_{0}\right)  $. Let $K$ be the metric completion of the ring of Puisseux
series as defined above. We introduce the notion of Cox arcs as elements of
\[
\left(  K^{\ast}\right)  ^{\Sigma\left(  1\right)  }/\operatorname*{Hom}%
\nolimits_{\mathbb{Z}}\left(  A_{n-1}\left(  Y\right)  ,K^{\ast}\right)
\cong\operatorname*{Hom}\nolimits_{\mathbb{Z}}\left(  M,K^{\ast}\right)
=\left(  K^{\ast}\right)  ^{n}%
\]
representing via the presentation of $A_{n-1}\left(  Y\right)  $ elements of
the torus $\left(  K^{\ast}\right)  ^{n}$ of dimension $n=\dim\left(
Y\right)  $. Then $V_{K}\left(  I\right)  \subset\left(  K^{\ast}\right)
^{n}$ is the image of the vanishing locus of $I\subset\mathbb{C}\left[
t\right]  \otimes S$ in $\left(  K^{\ast}\right)  ^{\Sigma\left(  1\right)
}/\operatorname*{Hom}\nolimits_{\mathbb{Z}}\left(  A_{n-1}\left(  X\left(
\Sigma\right)  \right)  ,K^{\ast}\right)  $. The limit map%
\[%
\begin{tabular}
[c]{llll}%
$\lim:$ & $B\left(  I\right)  $ & $\rightarrow$ & $\operatorname*{Strata}%
\left(  X_{0}\right)  \subset\operatorname*{Strata}\left(  Y\right)  $\\
& \multicolumn{1}{c}{$F$} & $\mapsto$ & $\left\{  \lim_{t\rightarrow0}a\left(
t\right)  \mid a\in\operatorname*{val}\nolimits^{-1}\left(
\operatorname*{int}\left(  F\right)  \right)  \right\}  $%
\end{tabular}
\]
associating to a face $F$ of $B\left(  I\right)  $ the stratum of $X_{0}$ of
limit points of arc solutions of $I$ over the interior of $F$. If $F$ is a
face of the special fiber Bergman complex $B\left(  I\right)  $, then there is
a unique cone $\tau$ of $\Sigma$ such that $\operatorname*{int}\left(
F\right)  \subset\operatorname*{int}\left(  \tau\right)  $ and%
\[
\lim\left(  F\right)  =V\left(  \tau\right)
\]
is the torus stratum of $Y$ corresponding to $\tau$.

In Sections \ref{Sec special fiber mirror degeneration general setup} and
\ref{Sec first order mirror degeneration} we define the mirror special fiber
and the first order conjectural mirror degeneration. This generalizes the case
of complete intersections. Denote by $d=\dim\left(  B\left(  I\right)
\right)  $ the fiber dimension of $\mathfrak{X}$. As noted above, the mirror
special fiber is the vanishing locus of the ideal%
\begin{align*}
I_{0}^{\circ}  &  =\left\langle
{\displaystyle\prod\limits_{v\in J}}
z_{v}\mid J\subset\Sigma^{\circ}\left(  1\right)  \text{ with }%
\operatorname*{supp}\left(  B\left(  I\right)  \right)  \subset%
{\displaystyle\bigcup\limits_{v\in J}}
F_{v}\right\rangle \\
&  =%
{\textstyle\bigcap\nolimits_{F\in B\left(  I\right)  _{d}}}
\left\langle z_{G^{\ast}}\mid G\text{ a facet of }\nabla\text{ with }F\subset
G\right\rangle \subset S^{\circ}%
\end{align*}
in $Y^{\circ}$ and the first order mirror degeneration is the vanishing locus
of%
\[
\left\langle m+t\cdot\sum_{\alpha\in\operatorname*{supp}\left(  \left(
\lim\left(  B\left(  I\right)  \right)  \right)  ^{\ast}\right)  \cap
N}a_{\alpha}\cdot\alpha\left(  m\right)  \mid m\in I_{0}^{\circ}\right\rangle
\subset\mathbb{C}\left[  t\right]  /\left\langle t^{2}\right\rangle \otimes
S^{\circ}%
\]

In Section \ref{sec orbifolding mirror families} we propose the notion of a
set of Fermat deformations associated to a monomial degeneration
$\mathfrak{X}$. The goal is to relate, if possible, the tropical mirror
degeneration $\mathfrak{X}^{\circ}$ to an orbifolding mirror family. For
simplicity we assume that the fibers of $\mathfrak{X}\subset\mathbb{P}%
^{n}\times\operatorname*{Spec}\mathbb{C}\left[  t\right]  $ are in projective
space. A set of Fermat deformations of $\mathfrak{X}$ is a set $\mathfrak{F}$
of non-trivial first order deformations of $X_{0}$ in $\mathfrak{X}$
corresponding to vertices of faces of $\operatorname*{dual}\left(  B\left(
I\right)  \right)  $. So the elements of $\mathfrak{F}$ have an interpretation
as Cox variables of $Y^{\circ}$. We require $\mathfrak{F}$ to satisfy the
following properties:

\begin{itemize}
\item For all zero dimensional strata $p$ of $\mathbb{P}^{n}$ precisely one of
the vertices of the faces $F\in\operatorname*{dual}\left(  B\left(  I\right)
\right)  $ with $\lim\left(  F\right)  =p$ is an element of $\mathfrak{F}$.

\item The convex hull of $\mathfrak{F}$ in $M_{\mathbb{R}}$ is a full
dimensional polytope containing $0$ in its interior, so $\mathfrak{F}$ spans a
fan $\hat{\Sigma}^{\circ}$ over a lattice simplex.

\item The elements of $\mathfrak{F}$ are incomparable with respect to the
preordering of Cox Laurent monomials given by divisibility of the denominators
(assumed to be relatively prime to the numerators).
\end{itemize}

The fan $\hat{\Sigma}^{\circ}$ defines a toric Fano variety $\hat{Y}^{\circ
}=X\left(  \hat{\Sigma}^{\circ}\right)  $ which is an orbifold of a weighted
projective space. We describe the special fiber of the monomial degeneration
$\mathfrak{\hat{X}}^{\circ}$ induced by $\mathfrak{X}^{\circ}$ via a
birational map $Y^{\circ}\rightarrow\hat{Y}^{\circ}$. The degeneration
$\mathfrak{\hat{X}}^{\circ}$ involves the first order deformations represented
by the degree $0$ Cox Laurent monomials%
\[
\left\{
{\textstyle\prod\nolimits_{r\in\hat{\Sigma}^{\circ}\left(  1\right)  }}
z_{r}^{\left\langle \hat{r},w\right\rangle }\mid w\in\left(  \lim\left(
B\left(  I\right)  \right)  \right)  ^{\ast}\cap N\right\}
\]
in the Cox ring of $\hat{Y}^{\circ}$.

\textbf{Section
\ref{Sec tropical mirror construction for Pfaffian Calabi-Yau varieties}.
}Here, we apply the tropical mirror construction to non complete intersection
Pfaffian examples.

We begin in Section \ref{1PfaffianCalabiYauThreefolds} by recalling the
structure theorem of Buchsbaum and Eisenbud for Pfaffian subschemes of
$\mathbb{P}^{n}$. Excluding special cases, locally Gorenstein subcanonical
schemes of codimension $3$ of $\mathbb{P}^{n}$ are locally given by the
Pfaffians of order $2k$ of a skew symmetric map $\varphi:\mathcal{E}\left(
-t\right)  \rightarrow\mathcal{E}^{\ast}$ for some vector bundle
$\mathcal{E}\rightarrow\mathbb{P}^{n}$ of rank $2k+1$. By the theorem of
Buchsbaum-Eisenbud, they have a locally free resolution of the form%
\[
0\rightarrow\mathcal{O}_{\mathbb{P}^{n}}\left(  -t-2s\right)  \rightarrow
\mathcal{E}\left(  -t-s\right)  \overset{\varphi}{\rightarrow}\mathcal{E}%
^{\ast}\left(  -s\right)  \rightarrow\mathcal{O}_{\mathbb{P}^{n}}%
\rightarrow\mathcal{O}_{X}\rightarrow0
\]
with $s=c_{1}\left(  \mathcal{E}\right)  +kt$. We recall the list of Pfaffian
Calabi-Yau threefolds given in \cite{Tonoli Canonical surfaces in mathbbP^5
and CalabiYau threefolds in mathbbP^6}.

In Section \ref{Sec deformations of Pfaffians} we make some remarks on the
deformation theory of Pfaffian ideals. The deformations of arithmetically
Gorenstein Pfaffian varieties are unobstructed and the base space is smooth
given by the independent entries of the skew symmetric syzygy matrix. These
observations allow to apply the tropical mirror construction to monomial
degenerations of Pfaffian ideals and to extend first order mirror families
$\mathfrak{X}^{\circ1}\subset Y^{\circ}\times\operatorname*{Spec}%
\mathbb{C}\left[  t\right]  /\left\langle t^{2}\right\rangle $, which obey the
structure theorem of Buchsbaum and Eisenbud, to flat families $\mathfrak{X}%
^{\circ}\subset Y\times\operatorname*{Spec}\mathbb{C}\left[  t\right]  $.

\begin{figure}
[h]
\begin{center}
\includegraphics[
height=2.9776in,
width=3.8683in
]%
{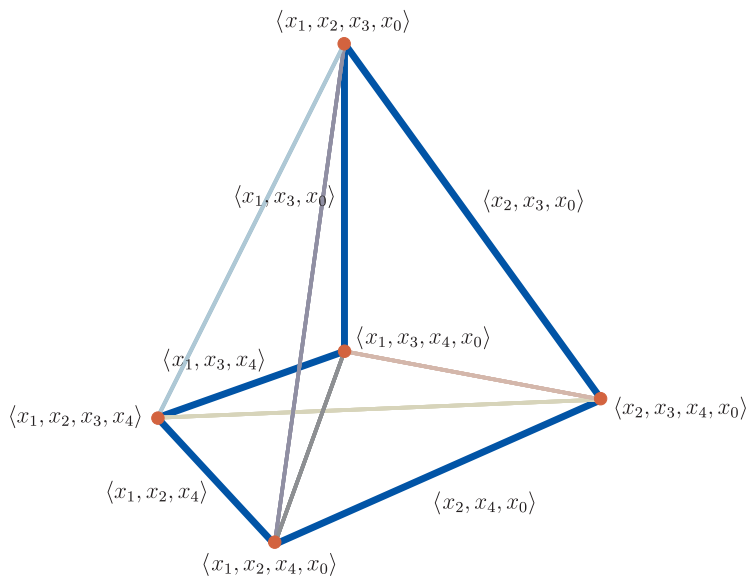}%
\caption{Projection of the complex $\lim\left(  B\left(  I\right)  \right)
\subset\operatorname*{Poset}\left(  \Delta\right)  $ for the monomial
degeneration of the general Pfaffian elliptic curve in $\mathbb{P}^{4}$}%
\label{Fig proj limB Pfaff ell}%
\end{center}
\end{figure}

Section \ref{Sec ex pfaffian elliptic curve} applies the tropical mirror
construction to the monomial degeneration $\mathfrak{X}$ of the general
Pfaffian elliptic curve in $\mathbb{P}^{4}$ as defined in Section
\ref{Sec degenerations of pfaffian calabi-yau varieties}. Note that the total
space of $\mathfrak{X}$ is a local complete intersection and the tropical
mirror construction treats this example much like a complete intersection. The
covering of $B\left(  I\right)  ^{\vee}$ in $\operatorname*{dual}\left(
B\left(  I\right)  \right)  $ is $c:1$ unbranched, but not a trivial $c:1$
covering like it would be for a codimension $c=3$ complete intersection.
Figure \ref{Fig proj limB Pfaff ell} shows a projection into $3$-space of the
complex $\lim\left(  B\left(  I\right)  \right)  \subset\operatorname*{Poset}%
\left(  \Delta\right)  $.%
\begin{figure}
[hh]
\begin{center}
\includegraphics[
height=2.6299in,
width=2.7121in
]%
{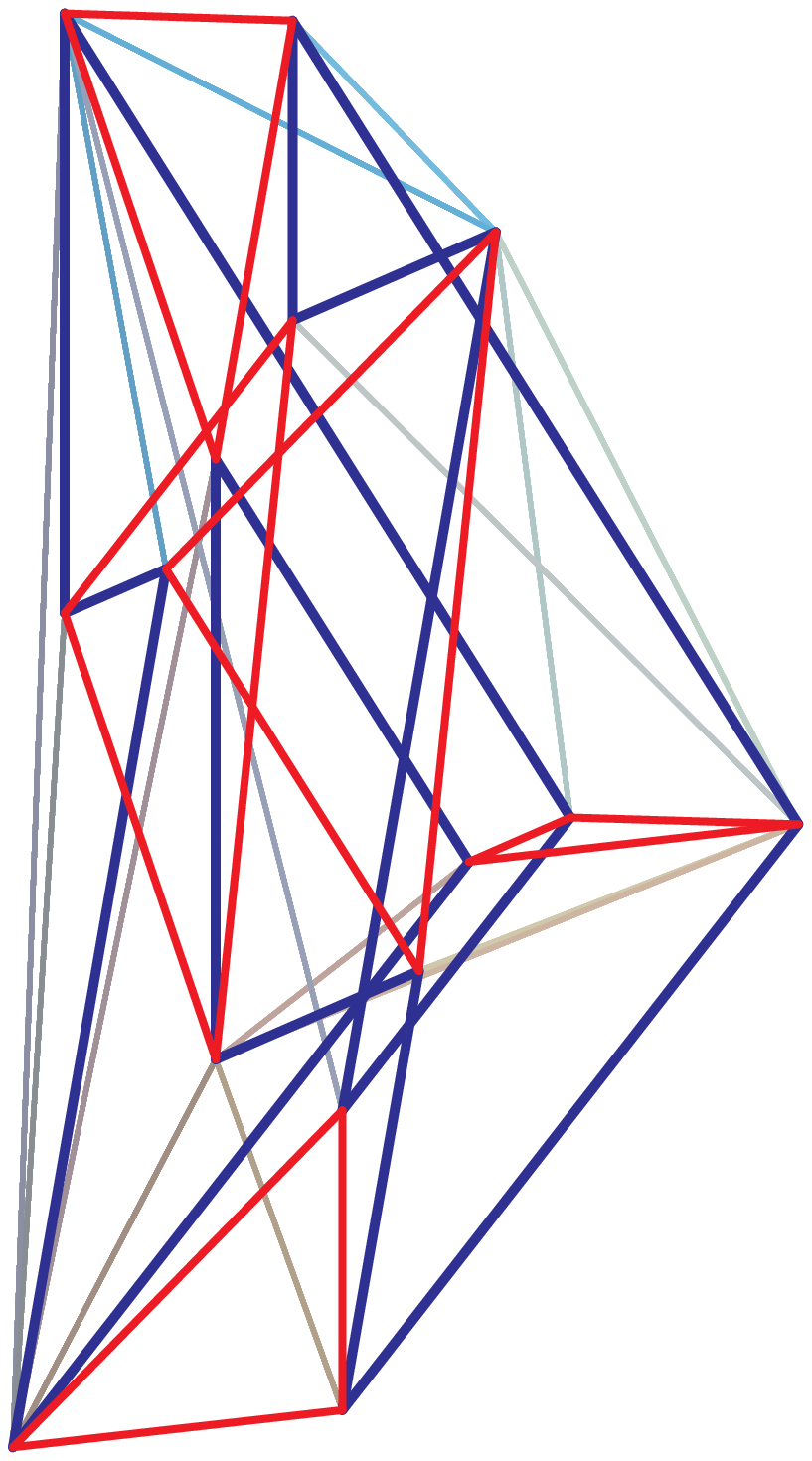}%
\caption{Projection of the complex $\operatorname*{dual}\left(  B\left(
I\right)  \right)  $ for a monomial degeneration of the generic Pfaffian
elliptic curve in $\mathbb{P}^{4}$}%
\label{Fig polyhedral plot dual B pfaff ell}%
\end{center}
\end{figure}
Figure \ref{Fig polyhedral plot dual B pfaff ell} visualizes a projection of
the facets of the polytope $\nabla^{\ast}$ and the subcomplex
$\operatorname*{dual}\left(  B\left(  I\right)  \right)  \subset
\operatorname*{Poset}\left(  \nabla^{\ast}\right)  $ of the boundary of
$\nabla^{\ast}$. Figure \ref{Fig dual Pfaff ell defs} shows the topology of
$\operatorname*{dual}\left(  B\left(  I\right)  \right)  $. The faces of the
complex are labeled by their image under $\lim$ and the lattice points of the
faces are labeled by the corresponding deformations of $I_{0}$. The complex
$\operatorname*{dual}\left(  B\left(  I\right)  \right)  $ has $5$ prisms as
facets of dimension $3$ and $5$ triangles as faces of dimension $2$. Every
prism intersects two of the other prisms along triangles. The vertices of the
triangles and the edges of the prisms connecting these vertices form the $3:1$
unbranched covering of $B\left(  I\right)  ^{\vee}$.%
\begin{figure}
[hhh]
\begin{center}
\includegraphics[
height=2.898in,
width=4.3137in
]%
{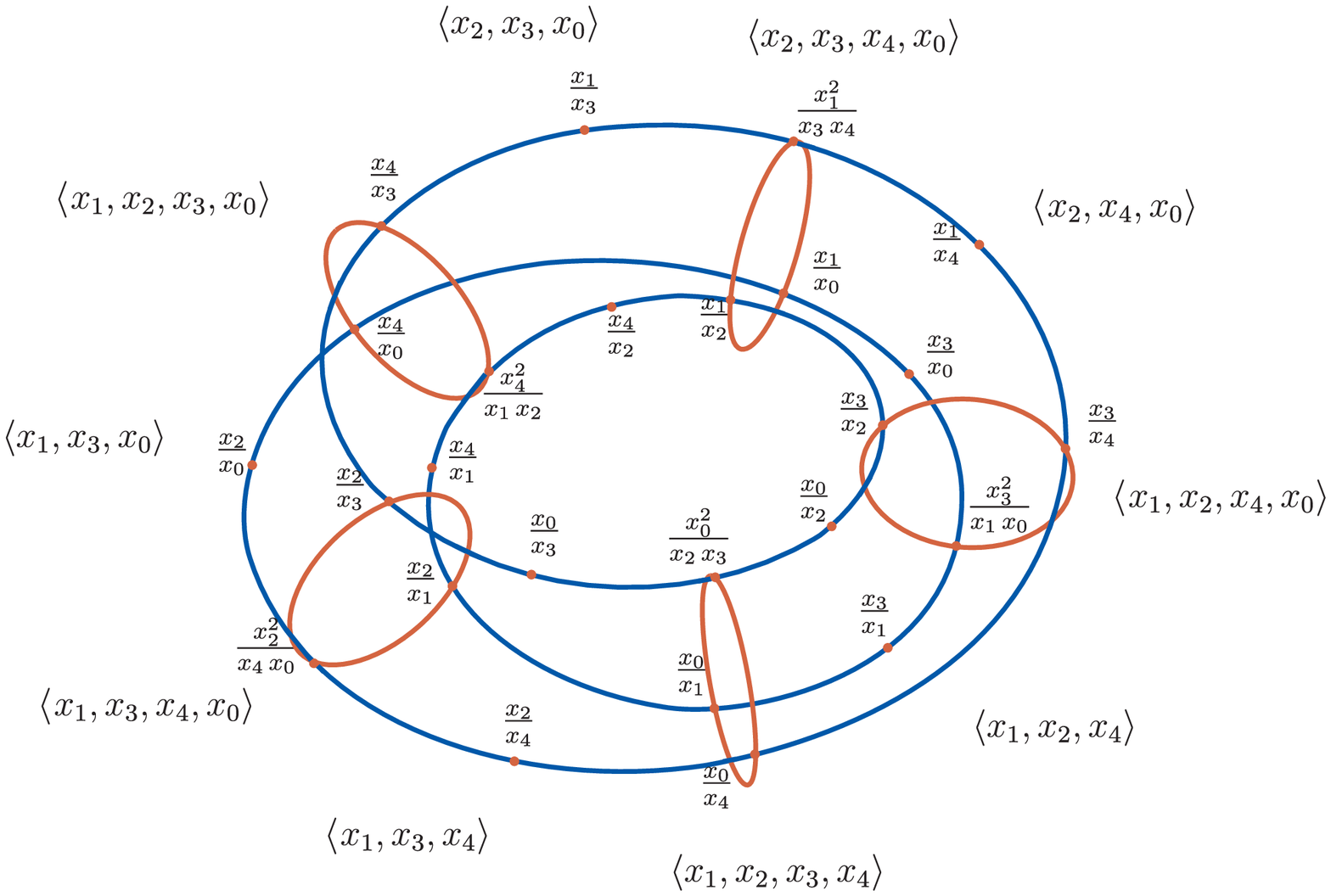}%
\caption{Topology of the complex $\operatorname*{dual}\left(  B\left(
I\right)  \right)  $ for a monomial degeneration of the generic Pfaffian
elliptic curve in $\mathbb{P}^{4}$}%
\label{Fig dual Pfaff ell defs}%
\end{center}
\end{figure}

In Section \ref{Sec ex pfaffian degree 14} we apply the tropical mirror
construction to the monomial degeneration given in Section
\ref{Sec degenerations of pfaffian calabi-yau varieties} for a general
Calabi-Yau threefold of degree $14$ in $\mathbb{P}^{6}$ defined by the
Pfaffians a general skew symmetric map $7\mathcal{O}\left(  -1\right)
\rightarrow7\mathcal{O}$. Using the concept of Fermat deformations from
Section \ref{sec orbifolding mirror families}, we relate the mirror
degeneration to the orbifolding mirror given by R\o dland.

Section \ref{Sec Ex 23333} applies the tropical mirror construction to a
monomial degeneration of the general Pfaffian Calabi-Yau threefold of degree
$13$ in $\mathbb{P}^{6}$ defined via a general skew symmetric map
$\mathcal{O}\left(  -2\right)  \oplus4\mathcal{O}\left(  -1\right)
\rightarrow\mathcal{O}\left(  1\right)  \oplus4\mathcal{O}$, as given in
Section \ref{Sec degenerations of pfaffian calabi-yau varieties}. Applying the
concept of Fermat deformations to switch to another birational model of
$Y^{\circ}$ with Chow group of rank $1$, we relate the mirror degeneration to
a Greene-Plesser orbifolding mirror family. This degeneration satisfies the
structure theorem of Buchsbaum-Eisenbud, in particular, allows extension of
the first order mirror degeneration.

Note again that the text of these examples is computer generated from the
output of the Maple package \textsf{tropicalmirror}, so all examples use the
same text fragments.

\textbf{Section \ref{Sec tropical stringy E}.} The next main section contains
some remarks on the tropical computation of stringy $E$-functions.

As this gives the general direction, we recall in Section
\ref{Sec Remarks on the Hodge numbers of Calabi-Yau} the relation of
$h^{d-1,1}\left(  X\right)  $, $h^{0}\left(  X,N_{X/\mathbb{P}^{n}}\right)  $
and $\operatorname*{Aut}\left(  \mathbb{P}^{n}\right)  $ for Calabi-Yau
manifolds of dimension $d$ in projective space $\mathbb{P}^{n}$.

Section \ref{Sec Batyrevs Hodge formula} explains Batyrev%
\'{}%
s original formulas for $h^{1,1}\left(  \bar{X}\right)  $ and $h^{d-1,1}%
\left(  \bar{X}\right)  $ via MPCP (maximal projective crepant partial)
desingularizations $\bar{X}\rightarrow X$.

In Section \ref{Sec First approximation of a tropical Hodge formula} we
explain a tropical method to compute $h^{1,\dim\left(  X\right)  -1}\left(
X\right)  $ for the general fiber $X$ of a Calabi-Yau monomial degeneration
with fibers in $Y=\mathbb{P}^{n}$, which is given by the ideal $I$. We
consider the lattice points of the dual complex $\operatorname*{dual}\left(
B\left(  I\right)  \right)  $, which do not correspond to roots of the toric
variety $Y$ (which are trivial deformations), and then divide out the torus.

In Section \ref{Sec string cohomology} we recall some known formulas for
stringy $E$-functions and give some ideas on a formula for the stringy
$E$-function of a Calabi-Yau variety computed from the tropical data we
associated to a monomial degeneration. We begin in Section
\ref{Sec Stringy Efunction of toric variety} with the example of the stringy
$E$-function of a toric variety. Section \ref{Sec the combinatorics of posets}
gives Batyrev%
\'{}%
s general concepts for the computation of stringy $E$-functions and we
consider the example of Calabi-Yau hypersurfaces in toric varieties in Section
\ref{Sec String theoretic Hodge formula for hypersurfaces}. Section
\ref{Sec string theoretic Hodge formula for complete intersections} recalls
Batyrev%
\'{}%
s and Borisov%
\'{}%
s computation of the stringy $E$-function of a complete intersection, which
works by relating the complete intersection to a hypersurface. Finally,
Section
\ref{Sec first approximation of a tropical computation of the stringy E function}
makes some remarks on a tropical formula for the stringy $E$-function.

\textbf{Section \ref{Sec implementation tropical mirror construction}.} This
section gives some remarks on computer algebra libraries which have been
written by the author in the context of the tropical mirror construction. See
also Section \ref{Sec toric mori theory} for some remarks on the
implementation within the Maple package \textsf{tropicalmirror} of Algorithm
\ref{Alg secondary fan GKZ} computing the secondary fan and the GKZ decomposition.

Section \ref{Sec moram2} explains the syntax of the Macaulay2 library
\textsf{mora.m2} which implements the standard basis algorithm. See also the
remarks in Section \ref{Sec Groebner basics}. The goal was to provide a
simple, transparent and flexible implementation capable of intermediate output
useful for testing and didactical purposes.

The Macaulay2 library \textsf{homology.m2}, described in Section
\ref{Sec homologym2}, computes the homology groups of a simplicial cell
complex and is useful in the context of Stanley filtrations.

Section \ref{Sec stanley filtration} gives an outline of the Macaulay2 library
\textsf{stanleyfiltration.m2} providing functions computing a Stanley
decomposition and Stanley filtration of an ideal and the set of monomial
ideals in a multigraded polynomial ring with given multigraded Hilbert polynomial.

Finally, Section \ref{Sec tropicalmirror implementation} explains the syntax
of the key functions of the Maple package \textsf{tropicalmirror} which gives
a full implementation of the tropical mirror construction. It takes as an
input a Fano polytope $P$ and a Calabi-Yau monomial degeneration
$\mathfrak{X}$ with fibers in the toric Fano variety $Y=X\left(
\Sigma\right)  $, where $\Sigma=\Sigma\left(  P\right)  $ is the fan over the
faces of $P$. The degeneration is defined by equations $m_{i}+tg_{i}\in
S\otimes\mathbb{C}\left[  t\right]  $, where the $m_{i}$ are the minimal
generators of a monomial ideal $I_{0}\ \subset S$ in the Cox ring $S$ of $Y$
defining the special fiber $X_{0}$ of $\mathfrak{X}$. The library outputs the
mirror Fano polytope $P^{\circ}$ and the first order mirror degeneration
$\mathfrak{X}^{\circ}$ with fibers in $X\left(  \Sigma^{\circ}\right)  $,
$\Sigma^{\circ}=\Sigma\left(  P^{\circ}\right)  $ specified in an analogous
way to $\mathfrak{X}$. The package \textsf{tropicalmirror} provides functions
to compute the various intermediate objects introduced in the tropical mirror
construction. The library also contains a function to find sets of Fermat
deformations and the corresponding contracted degeneration $\mathfrak{\hat{X}%
}^{\circ}$. We provide a function, which tests whether a first order
degeneration satisfies the Koszul or Buchsbaum-Eisenbud structure theorem and
extends degenerations to a flat family over $\operatorname*{Spec}%
\mathbb{C}\left[  t\right]  $, provided the arithmetically Gorenstein
Buchsbaum-Eisenbud structure theorem applies.

\textbf{Section \ref{Sec perspectives}.} The last main section explores the
perspectives of the tropical mirror construction, the underlying concepts and
technical formalisms.

In Section
\ref{Sec tropical computation of the stringy E function and string cohomology}
we note that the natural next step is to compute from the tropical objects the
stringy $E$-functions and more generally string cohomology in order to deal
with the singular general fibers appearing in the tropical mirror
construction. The tropical formula should generalize Batyrev%
\'{}%
s formula for the stringy $E$-function of anticanonical hypersurfaces in
Gorenstein toric Fano varieties.

Section \ref{Sec Hilbert schemes and moduli spaces} raises the question of the
computation of the local Hilbert scheme and moduli stack. The concepts of
Section \ref{Sec Hilbert scheme state polytope Cox homogeneou} introduced by
Haiman and Sturmfels allow algorithmic computation of the local equations of
the Hilbert scheme for ideals in the Cox ring of a smooth toric variety. Using
the ideas noted in Sections \ref{Sec Multigraded Hilbert schemes} and
\ref{Monomial ideals in the Cox ring and the stratified toric primary decomposition}%
, one should be able to generalize the multigraded regularity and Hilbert
scheme to the setting of simplicial and further to non simplicial toric varieties.

Section \ref{Sec integrally affine structures} raises the question of relating
the tropical mirror construction to the mirror construction by Gross and
Siebert via integrally affine structures.

As noted in Section \ref{Sec torus fibrations}, the non-Archimedian amoeba map
gives a degenerate torus fibration of the special fiber, hence we ask the
question to obtain from this, via the amoeba map, a torus fibration of the
general fiber.

Section \ref{1FurtherStanleyReisnerexamples} suggests to apply the tropical
mirror construction to further Calabi-Yau degenerations with fibers in
projective space. Altmann and Christophersen compute the first order
deformations and obstructions of Stanley-Reisner rings. Applying these
algorithms in the case of triangulations of spheres one can obtain the
necessary data to apply the tropical mirror construction for smoothable
examples. We also ask in Section
\ref{Sec deformations and obstructions of a toric generalization of stanley reisner rings}
how to generalize the work of Altmann and Christophersen to the non-simplicial
toric setting.

Section \ref{Sec mirror pfaffian resolution in Cox ring} suggests to extend
the structure theorem of Buchsbaum-Eisenbud to describe codimension $3$
Calabi-Yau ideals in the Cox ring of a toric variety.

Tropical geometry is known to reflect the $p$-adic geometry, so Section
\ref{Sec tropical geometry finite fields} raises the relation of the tropical
mirror construction to mirror symmetry over finite fields and $\zeta
$-functions as explored by Candelas et al.

A central question is the extension of topological mirror symmetry to the
stronger condition of mathematical mirror symmetry via Frobenius manifolds. So
Section \ref{Sec instanton numbers} makes some remarks on the question of
computation of instanton numbers and the $A$-model correlation functions.
Section \ref{Sec GKZ} raises the question of describing quantum cohomology
rings from the tropical data via $GKZ$ hypergeometric differential equations
associated to toric data.\medskip

\textbf{Acknowledgements.}

My hearty thanks go to Frank-Olaf Schreyer for his great support and the
privilige to be his assistant at the University of Saarland. He always had an
open ear to my questions and supported my mathematical development. He made
many helpful suggestions and introduced me to the mathematical foundations
necessary for the tropical mirror construction. Also many thanks to Wolfram
Decker and all the members of the algebraic geometry group in Saarbr\"{u}cken
for the great working environment. I wish to heartily thank Stavros Papadakis
always being willing to answer so many of my questions about algebraic
geometry. Many thanks to Christine Wilk-Pitz for being a great colleague. My
thanks also go to all members of the algebraic geometry group in Bayreuth.

I would like to thank the University of Saarland, the Forschungsschwerpunkt
"Globale Methoden in der komplexen Geometrie" and the University of Oslo for
their support.

I am grateful to Jan Arthur Christophersen for a lot of helpful suggestions
and fruitful discussions during my stay at the University of Oslo in spring
2006 and also all other members of the algebraic geometry group in Oslo. I
wish to thank Bernd Siebert and Mark Gross for their great hospitality and the
helpful discussions during my visits in 2004 and 2006 at
Albert-Ludwigs-Universit\"{a}t in Freiburg. Also I would like to thank Klaus
Altmann for the visit at the Freie Universit\"{a}t Berlin in summer 2006 and
Duco van Straten for his hospitality at Johannes Gutenberg-Universit\"{a}t in Mainz.

Many thanks to the organizers David Eisenbud, Joe Harris and Frank-Olaf
Schreyer of the Oberwolfach Classical Algebraic Geometry meeting in 2004, and
Gregory Mikhalkin whose talk on the tropical curve count on Fano surfaces
triggered my interest in considering tropical geometry in the context of
mirror symmetry. I would like to thank the organizers and lecturers of the
Nordfjordeid conference on tropical geometry and algebraic statistics in fall
2006 and especially Bernd Sturmfels from UC Berkeley for our great discussions
concerning the tropical mirror construction.\newpage

\section{Prerequisites\label{Sec prerequisites}}

\subsection{Calabi-Yau varieties and mirror
symmetry\label{Sec calabi yau and mirror symmetry}}

\begin{definition}
A normal projective $d$-dimensional algebraic variety $X$ is called a
\index{Calabi Yau variety|textbf}%
\textbf{Calabi-Yau variety} if it has at worst Gorenstein canonical
singularities, $K_{X}=\mathcal{O}_{X}$ and $h^{i}\left(  X,\mathcal{O}%
_{X}\right)  =0$ for $0<i<d$.
\end{definition}

\begin{remark}
The
\index{Hodge diamond|textbf}%
\textbf{Hodge diamond }of a Calabi-Yau $d$-fold $X$, formed by the
\index{Hodge numbers|textbf}%
Hodge numbers $h^{p,q}\left(  X\right)  =\dim H_{\bar{\partial}}^{p,q}\left(
X\right)  =\dim H^{q}\left(  X,\Omega_{X}^{p}\right)  $, has horizontal and
vertical symmetry by
\index{Serre duality}%
Serre duality and
\index{Hodge duality}%
Hodge duality, and $\Omega_{X}^{d}=K_{X}=\mathcal{O}_{X}$, hence
\[
H^{0,i}\left(  X\right)  \cong H^{i}\left(  X,\mathcal{O}_{X}\right)  \cong
H^{i}\left(  X,\Omega_{X}^{d}\right)  \cong H^{d,i}\left(  X\right)
\]
e.g., for $d=3$
\[%
\begin{tabular}
[c]{ccccccc}
&  &  & $1$ &  &  & \\
&  & $0$ &  & $0$ &  & \\
& $0$ &  & $h^{2,2}\left(  X\right)  $ &  & $0$ & \\
$1$ &  & $h^{2,1}\left(  X\right)  $ &  & $h^{1,2}\left(  X\right)  $ &  &
$1$\\
& $0$ &  & $h^{1,1}\left(  X\right)  $ &  & $0$ & \\
&  & $0$ &  & $0$ &  & \\
&  &  & $1$ &  &  &
\end{tabular}
\]

\end{remark}

\begin{definition}
A pair of smooth Calabi-Yau $d$-folds $X$ and $X^{\circ}$ is called a
\index{topological mirror pair|textbf}%
\textbf{topological mirror pair} if their Hodge numbers satisfy%
\begin{equation}
h^{p,q}\left(  X\right)  =h^{d-p,q}\left(  X^{\circ}\right)  \text{ }%
\forall0\leq p,q\leq d \label{4topologmirrorsym}%
\end{equation}
i.e., the Hodge diamond is mirrored at the diagonal.
\end{definition}

This definition can be extended via the stringy $E$-functions for varieties
with log-terminal singularities. For the precise definition see Section
\ref{Stingy E function}.

\begin{definition}
Calabi-Yau varieties $X$ and $X^{\circ}$ of dimension $d$ are called a
\index{stringy topological mirror pair|textbf}%
\textbf{stringy topological mirror pair} if the stringy $E$-functions satisfy%
\[
E_{st}\left(  X;u,v\right)  =u^{d}E_{st}\left(  X^{\circ};u^{-1},v\right)
\]
For Gorenstein varieties \newsym{Stringy $E$-function}{$E_{st}$} is the
\index{generating function}%
generating function for the
\index{stringy Hodge numbers}%
stringy Hodge numbers, which coincide with the Hodge numbers of a
\index{crepant resolution}%
crepant resolution if such one exists.
\end{definition}

\begin{remark}
To a Calabi-Yau manifold we can associate two
\index{Frobenius manifold}%
Frobenius manifolds called
\index{A-model}%
$A$- and
\index{B-model}%
$B$-model (see e.g. \cite{Manin Frobenius Manifolds Quantum Cohomology and
Moduli Spaces} and \cite{CK Mirror Symmetry and Algebraic Geometry}). $X$ and
$X^{\circ}$ are called a
\index{mathematical mirror pair|textbf}%
\textbf{mathematical mirror pair} if the $A$-model of $X$ is isomorphic to the
$B$-model of $X^{\circ}$ and vice versa.
\end{remark}

\begin{remark}
String theory
\index{String theory}%
replaces particles by extended objects, e.g., by an $S^{1}$ or an interval.
Whereas a point sweeps out a real $1$-dimensional object in spacetime, a
propagating string gives a surface, called its
\index{worldsheet}%
worldsheet.

There are $5$ possible
\index{superstring theory}%
superstring theories, which are defined on a real $10$-dimensional
\index{Riemannian manifold}%
Riemannian manifold. One assumes that this manifold is locally the product of
a real $4$-dimensional Riemannian manifold $M_{4}$ and a $6$-dimensional
compact Riemannian manifold $X$ too small to appear in measurements.

In the case of type $IIA$ and $IIB$ superstring theory one concludes that
$M_{4}$ is a
\index{Minkowski space}%
Minkowski space, that $X$ has a complex structure $J$, that there is a
K\"{a}hler metric $g$ on $\left(  X,J\right)  $ and that $g$ has holonomy
group $\operatorname*{Hol}\left(  g\right)  \subset SU\left(  d\right)  $. As
explained in Remark \ref{CalabiYau}, these conditions are satisfied by a
Calabi-Yau manifold.

The
\index{worldsheet}%
worldsheets project to algebraic curves on this threefold. The
\index{Hodge numbers}%
Hodge numbers of $X$ are important characteristica of the physical theory,
e.g.,
\[
\frac{1}{2}\left\vert \chi\left(  X\right)  \right\vert =\left\vert
h^{1,1}\left(  X\right)  -h^{2,1}\left(  X\right)  \right\vert
\]
is the number of fermion generations. So this number is identical for two
manifolds, which form a mirror pair. Experiments indicate that real world has
$3$ fermion generations.

From the point of view of physics,
\index{mirror symmetry}%
mirror symmetry of two Calabi-Yau threefolds $X$ and $X^{\circ}$ is the
\index{duality}%
duality of two of these types of compactified string theories, which is again
a stronger condition than $X$ and $X^{\circ}$ forming a mathematical mirror pair.
\end{remark}

\begin{remark}
The most simple case of duality in physics is found in
\index{Maxwell equations}%
Maxwell%
\'{}%
s equations, describing the
\index{electromagnetic interaction}%
electromagnetic interaction:%
\[%
\begin{tabular}
[c]{lll}%
$\partial_{\nu}F^{\mu\nu}=-j^{\mu}$ &  & $\partial_{\nu}\tilde{F}^{\mu\nu
}=-k^{\mu}$%
\end{tabular}
\]
with%
\begin{align*}
F^{\mu\nu}  &  =\left(
\begin{array}
[c]{cccc}%
0 & -E_{x} & -E_{y} & -E_{z}\\
E_{x} & 0 & -B_{z} & B_{y}\\
E_{y} & B_{z} & 0 & -B_{x}\\
E_{z} & -B_{y} & B_{x} & 0
\end{array}
\right) \\
\tilde{F}^{\mu\nu}  &  =\frac{1}{2}\varepsilon^{\mu\nu\alpha\beta}%
F_{\alpha\beta}=\left(
\begin{array}
[c]{cccc}%
0 & -B_{x} & -B_{y} & -B_{z}\\
B_{x} & 0 & E_{z} & -E_{y}\\
B_{y} & -E_{z} & 0 & E_{x}\\
B_{z} & E_{y} & -E_{x} & 0
\end{array}
\right)
\end{align*}
$E$ is the electric and $B$ the magnetic field, $j=\left(  \rho,j_{x}%
,j_{y},j_{z}\right)  $ the electric four-current with charge density $\rho$
and electric three-current $\left(  j_{x},j_{y},j_{z}\right)  $ and $k=\left(
\sigma,k_{x},k_{y},k_{z}\right)  $ is the magnetic four-current introduced by
Dirac, and index manipulations are done with respect to the flat
\index{Minkowski space}%
Minkowski metric with signature $\left(  +,-,-,-\right)  $. These equations
are invariant under an $\operatorname*{SO}\left(  2\right)  $ rotating $E$ and
$B$, in particular under%
\[%
\begin{tabular}
[c]{lllllll}%
$E\mapsto B$ &  & $B\mapsto-E$ &  & $j\mapsto k$ &  & $k\mapsto-j$%
\end{tabular}
\]
One can deduce that electrostatic theory for high interaction energies, which
is difficult to solve, is equivalent to magnetic theory for low interaction
energies, which is easy to solve. In the case of
\index{mirror symmetry}%
mirror symmetry,
\index{duality}%
duality allows for example the treatment of enumerative problems in algebraic geometry.
\end{remark}

\begin{remark}
(see also \cite{GHJ CalabiYau Manifolds and Related Geometries}%
)\label{CalabiYau} Let $X$ be a manifold of dimension $d$, $E$ a vector bundle
on $X$ and
\[
\nabla:\mathcal{A}\left(  E\right)  \rightarrow\mathcal{A}\left(  E\otimes
T^{\ast}\left(  X\right)  \right)
\]
a connection on $E$, where $\mathcal{A}\left(  E\right)  $ is the sheaf of
smooth sections of $E$. Then for any smooth curve $\gamma:\left[  0,1\right]
\rightarrow X$ with $\gamma\left(  0\right)  =x$ and $\gamma\left(  1\right)
=y$ with $x,y\in X$ and any $v\in E_{x}$ there is a unique smooth section
$\sigma\in\gamma^{\ast}\left(  E\right)  $ with $\nabla_{\gamma\left(
t\right)  }\sigma\left(  t\right)  =0$ for all $t\in\left[  0,1\right]  $ and
$\sigma\left(  0\right)  =v$. So one can associate the \textbf{parallel
transport map}%
\[%
\begin{tabular}
[c]{llll}%
$P_{\gamma}:$ & $E_{x}$ & $\rightarrow$ & $E_{y}$\\
& $v$ & $\mapsto$ & $\sigma\left(  1\right)  $%
\end{tabular}
\]

The \textbf{holonomy group} $\operatorname*{Hol}_{x}\left(  \nabla\right)  $
of $\nabla$ based at $x\in X$ is%
\[
\operatorname*{Hol}\nolimits_{x}\left(  \nabla\right)  =\left\{  P_{\gamma
}\mid\gamma\text{ a loop based at }x\right\}  \subset\operatorname*{GL}\left(
E_{x}\right)
\]

If $g$ is a Riemannian metric on $X$, there is a unique torsion free
connection $\nabla$ on $X$ with $\nabla g=0$, which is called the
\textbf{Levi-Civita connection}.

If $x\in X$, then $\operatorname*{Hol}\nolimits_{x}\left(  g\right)  $ is the
holonomy group of the Levi-Civita connection on the Riemannian manifold
$\left(  X,g\right)  $. As $\nabla g=0$, it follows that $g$ is invariant
under the natural action of the holonomy group, so $\operatorname*{Hol}%
\nolimits_{x}\left(  g\right)  $ is up to conjugation a subgroup of
$\operatorname*{O}\left(  d\right)  $ and is denoted by
$\newsym[\operatorname*{Hol}]{Holonomy group}{\operatorname*{Hol}\left(  g\right)}$%
.

If $\left(  X,g\right)  $ is a Riemannian manifold of dimension $r$ such that
$X$ is simply-connected, $g$ is irreducible (i.e., $\left(  X,g\right)  $ is
not locally isometric to a Riemannian product) and $g$ is non-symmetric, then
Berger%
\'{}%
s classification (see \cite{Berger Sur les groupes dholonomie des varietes a
connexion affine et des varietes Riemanniennes}) shows that

\begin{enumerate}
\item $\operatorname*{Hol}\left(  g\right)  =\operatorname*{SO}\left(
r\right)  $. This is the case of the generic Riemannian metric.

\item $\operatorname*{Hol}\left(  g\right)  =\operatorname*{U}\left(
d\right)  $ with $r=2d$, $d\geq2$. This is the case of a generic K\"{a}hler
manifold, in particular $X$ is a complex manifold.

\item $\operatorname*{Hol}\left(  g\right)  =\operatorname*{SU}\left(
d\right)  \subset\operatorname*{SO}\left(  r\right)  $ with $r=2d$, $d\geq2$.
Then $X$ is a Ricci-flat K\"{a}hler manifold. $X$ is a Calabi-Yau manifold
(omitting the condition algebraic) if $X$ is compact .

\item $\operatorname*{Hol}\left(  g\right)  =\operatorname*{Sp}\left(
a\right)  \subset\operatorname*{SO}\left(  r\right)  $ with $r=4a$, $a\geq2$.
Then $X$ is a Ricci-flat K\"{a}hler manifold of complex dimension $2a$. If $X$
is compact, then $X$ is called compact \textbf{hyperk\"{a}hler manifold} (it
admits many K\"{a}hler metrics).

\item $\operatorname*{Hol}\left(  g\right)  =\operatorname*{Sp}\left(
a\right)  \operatorname*{Sp}\left(  1\right)  \subset\operatorname*{SO}\left(
r\right)  $ with $d=4a$, $a\geq2$. In this case $X$ is called
\textbf{quaternionic K\"{a}hler} (note that $\operatorname*{Sp}\left(
a\right)  $ and $\operatorname*{Sp}\left(  a\right)  \operatorname*{Sp}\left(
1\right)  $ are groups of automorphism of $\mathbb{H}^{d}$, denoting by
$\mathbb{H}$ the quaternions), it is not K\"{a}hler.

\item $\operatorname*{Hol}\left(  g\right)  =G_{2}\subset\operatorname*{SO}%
\left(  7\right)  $ and $r=7$, a so called exceptional case.

\item $\operatorname*{Hol}\left(  g\right)  =\operatorname*{Spin}\left(
7\right)  \subset\operatorname*{SO}\left(  8\right)  $ and $r=8$, the other
exceptional case.
\end{enumerate}

A Riemannian manifold $\left(  X,g\right)  $ of dimension $r=2d$ is K\"{a}hler
if and only if $\operatorname*{Hol}\left(  g\right)  \subset\operatorname*{U}%
\left(  d\right)  \subset\operatorname*{O}\left(  r\right)  $. Then $X$ has a
complex structure $J$.

If $\left(  X,J,g\right)  $ is a K\"{a}hler manifold and $\rho$ its Ricci
form, then%
\[
\left[  \rho\right]  =2\pi c_{1}\left(  X\right)
\]
in $H^{2}\left(  X,\mathbb{R}\right)  $.

Let $\left(  X,J,g\right)  $ be a K\"{a}hler manifold of dimension $d$. If
$\operatorname*{Hol}\left(  g\right)  \subset\operatorname*{SU}\left(
d\right)  $, then $g$ is Ricci-flat. If $g$ is Ricci-flat and $K_{X}%
=\mathcal{O}_{X}$, then $\operatorname*{Hol}\left(  g\right)  \subset
\operatorname*{SU}\left(  d\right)  $. If $X$ is Ricci-flat and simply
connected, then $K_{X}=\mathcal{O}_{X}$.

Yau%
\'{}%
s proof of the Calabi conjecture implies: If $\left(  X,J\right)  $ is a
compact complex manifold, admitting K\"{a}hler metrics, and $c_{1}\left(
X\right)  =0$, then in each K\"{a}hler class, there is a unique Ricci-flat
K\"{a}hler metric. The Ricci-flat K\"{a}hler metrics on $X$ form a smooth
family of dimension $h^{1,1}\left(  X\right)  $, which is isomorphic to the
K\"{a}hler cone of $X$.

Any compact Ricci-flat K\"{a}hler manifold $\left(  X,J,g\right)  $ is up to a
finite cover isometric to the product of

\begin{itemize}
\item a flat K\"{a}hler torus

\item a compact simply-connected Riemannian manifold $N$.
\end{itemize}

$N$ is a Riemannian product of non-symmetric Ricci-flat irreducible Riemannian
K\"{a}hler manifolds $X_{j}$ of real dimension $r_{j}$ with
$\operatorname*{Hol}\left(  g_{j}\right)  \subset\operatorname*{SU}\left(
d_{j}\right)  $ and $r_{j}=2d_{j}$, $d_{j}\geq2$, which are

\begin{itemize}
\item a Calabi-Yau manifold (omitting the condition algebraic), i.e.,
$\operatorname*{Hol}\left(  g\right)  =\operatorname*{SU}\left(  d_{j}\right)
$, or

\item a Hyperk\"{a}hler manifold, i.e., $\operatorname*{Hol}\left(  g\right)
=\operatorname*{Sp}\left(  a_{j}\right)  $ with $r_{j}=2d_{j}=4a_{j}$ and
$a_{j}\geq2$.
\end{itemize}

Let $\left(  X,J,g\right)  $ be a compact K\"{a}hler manifold of dimension
$r=2d$ with $d\geq2$ and $\operatorname*{Hol}\left(  g\right)
=\operatorname*{SU}\left(  d\right)  $, then $X$ has finite fundamental group,
$h^{0,0}\left(  X\right)  =h^{d,0}\left(  X\right)  =1$ and $h^{i,0}\left(
X\right)  =0$ for $0<i<d$. If $d$ is even, then $X$ is simply connected.

If $d\geq3$, then $\left(  X,J\right)  $ is isomorphic to a complex
submanifold of some $\mathbb{P}_{\mathbb{C}}^{N}$ and is algebraic.

If $d=2$, then $\operatorname*{SU}\left(  2\right)  =\operatorname*{Sp}\left(
1\right)  $. The moduli space of Calabi-Yau twofolds (omitting the condition
algebraic), i.e., \textbf{K3-surfaces}, is a connected complex space of
dimension $20$. All Calabi-Yau twofolds are diffeomorphic. The algebraic $K3$
surfaces form a countable dense union of subvarieties of dimension $19$ inside
the moduli space.
\end{remark}

The following example was the first known
\index{mirror construction}%
mirror construction for Calabi-Yau varieties and was given by
\index{Greene-Plesser}%
Greene and Plesser, see \cite{GrPl Duality in CalabiYau moduli space} and
\cite{COGP A pair of CalabiYau manifolds as an exactly soluble superconformal
field theory}.

\begin{example}
\label{GreenePlesserQuinticExample}For a general
\index{quintic threefold}%
quintic threefold $X\subset\mathbb{P}^{4}$, by $T_{X^{\circ}}\cong
\Omega_{X^{\circ}}^{2}$ the mirror $X^{\circ}$
\index{hypersurface}%
should
\index{Greene-Plesser}%
satisfy%
\[
\dim H^{1}\left(  X^{\circ},T_{X^{\circ}}\right)  =h^{2,1}\left(  X^{\circ
}\right)  =h^{1,1}\left(  X\right)  =1
\]
hence in order to construct the
\index{mirror construction}%
mirror one looks for a $1$-parameter family. It turns out that the right
choice is
\begin{equation}
X_{\lambda}=\left\{  x_{0}^{5}+x_{1}^{5}+x_{2}^{5}+x_{3}^{5}+x_{4}^{5}+\lambda
x_{0}x_{1}x_{2}x_{3}x_{4}=0\right\}  \label{4GPFermatQuinticfamily}%
\end{equation}
divided out by the action of
\[
\frac{\left\{  \left(  a_{0},...,a_{4}\right)  \in\mathbb{Z}_{5}^{5}\mid
\sum_{i=0}^{4}a_{i}\equiv0\operatorname{mod}5\right\}  }{\mathbb{Z}_{5}\left(
1,...,1\right)  }%
\]
via%
\[
\left(  a_{0},...,a_{4}\right)  \left(  x_{0}:...:x_{4}\right)  =\left(
\mu^{a_{0}}x_{0}:...:\mu^{a_{4}}x_{4}\right)
\]
where $\mu$ is a $5$th
\index{root of unity}%
root of unity. Resolving the singularities of this singular quotient without
destroying the Calabi-Yau property gives the mirror of $X$.
\end{example}

\begin{remark}
These kind of orbifolding constructions were generalized for some
\index{hypersurface}%
hypersurfaces in weighted projective space, for
\index{complete intersection}%
complete intersections in $\mathbb{P}^{n}$ and complete intersections in
products of weighted projective spaces. See, e.g., \cite{CLS CalabiYau
Manifolds in Weighted mathbbP^4}, \cite{BH A Generalized Construction of
Mirror Manifolds}, \cite{CDLS Complete intersection CalabiYau manifolds},
\cite{LT Lines on CalabiYau complete intersections mirror symmetry and
PicardFuchs equations},\linebreak\cite{KS LandauGinzburg String Vacua}.

As some of the examples did not have a mirror in their classes, these
approaches were unified by
\index{Batyrev}%
Batyrev for
\index{hypersurface}%
hypersurfaces in toric varieties and by Batyrev and Borisov for
\index{complete intersection}%
complete intersections in toric varieties, see \cite{Batyrev Dual polyhedra
and mirror symmetry for CalabiYau hypersurfaces in toric varieties},
\cite{Borisov Towards the mirror symmetry for CalabiYau complete intersections
in Gorenstein Fano toric varieties}, \cite{BB On CalabiYau complete
intersections in toric varieties in HigherDimensional Complex Varieties Trento
1994}.
\end{remark}

\subsection{Mirror symmetry for singular Calabi-Yau varieties and stringy
Hodge
numbers\label{Mirror symmetry for singular Calabi-Yau varieties and stringy Hodge numbers}%
}

The following considerations allow to introduce a well-defined notion of
mirror symmetry for a certain class of singular varieties. This justifies to
give
\index{mirror construction}%
mirror constructions for and leading to singular Calabi-Yau varieties.

\subsubsection{Setup}

In constructing mirror pairs we encounter several problems: Even if we start
with a manifold, we encounter singular varieties. See, e.g., the quintic in
$\mathbb{P}^{4}$ in Example \ref{GreenePlesserQuinticExample}. First of all we
know that we can resolve the singularities by a sequence of
\index{blowup}%
blowups:

\begin{theorem}
[Hironaka]\cite{Hironaka Resolution of singularities of an algebraic variety
over a field of characteristic zero} Let $X$ be
\index{Hironaka theorem}%
a normal projective variety over an algebraically closed field of
characteristic $0$. For any proper subvariety $D\subset X$ there exists a
smooth projective variety $Y$ and a birational morphism $f:Y\rightarrow X$
such that $f^{-1}\left(  D\right)  $ is a divisor with only
\index{simple normal crossings}%
simple normal crossings (and $f$ is a composition of
\index{blowup}%
blowups in smooth closed centers).
\end{theorem}

For a proof and an algorithmic implementation of Hironaka%
\'{}%
s theorem, see for example \cite{Villamayor Constructiveness of Hironakas
resolution}, \cite{EnHa Strong resolution of singularities in characteristic
zero}, \cite{Hauser The Hironaka Theorem on resolution of singularities Or: A
proof we always wanted to understand} and \cite{FrPf Algorithmic Resolution of
Singularities in C. Lossen G. Pfister: Singularities and Computer Algebra}. Of
course we want the resolved variety to be still a Calabi-Yau:

\begin{definition}
A birational projective morphism $f:Y\rightarrow X$ with $Y$ smooth and $X$ at
worst Gorenstein canonical singularities is called
\index{crepant resolution|textbf}%
\textbf{crepant }(or non discrepant)\textbf{ desingularization} of $X$ if
$f^{\ast}K_{X}=K_{Y}$ (i.e., if the
\index{discrepancy|textbf}%
\textbf{discrepancy} $K_{Y}-f^{\ast}K_{X}$ is zero).
\end{definition}

If the crepant desingularizations of $Y\rightarrow X$ resp. $Y^{\circ
}\rightarrow X^{\circ}$ exist, we can define a topological mirror pair by%
\[
h^{p,q}\left(  Y\right)  =h^{d-p,q}\left(  Y^{\circ}\right)  \text{ }%
\forall0\leq p,q\leq d
\]
However it is not obvious that this is well defined: If a crepant
desingularization exists, it is not necessarily unique. In particular, given
two crepant resolutions $Y_{1}\rightarrow X$ and $Y_{2}\rightarrow X$ it is
not clear a priori that the
\index{Hodge numbers}%
Hodge numbers of $Y_{1}$ and $Y_{2}$ are equal. We will see in Theorem
\ref{TheoremCrepant} that they indeed are.

\begin{example}
\label{ExampleEfunction1}Let $X_{0}$ be a smooth
\index{Fano}%
Fano embedded by a very ample line bundle $L$ with $L^{k}=K_{X_{0}}^{-l}$
($k,l\in\mathbb{N}$), let $E=\mathcal{O}_{X_{0}}\oplus L$ and consider the
map
\[%
\begin{tabular}
[c]{llll}%
$\pi:$ & $Y=\mathbb{P}\left(  E\right)  $ & $\rightarrow$ & $X\subset
\overset{\cong H^{0}\left(  \mathbb{P}\left(  E\right)  ,\mathcal{O}%
_{\mathbb{P}\left(  E\right)  }\left(  1\right)  \right)  }{\mathbb{P}\left(
H^{0}\left(  X_{0},\mathcal{O}_{X_{0}}\oplus L\right)  \right)  }$\\
& $\downarrow\uparrow\sigma$ &  & \\
& $X_{0}$ &  &
\end{tabular}
\]
which is the contraction of $\sigma\left(  X_{0}\right)  \cong X_{0}$ to $p\in
X$ where $\sigma:X_{0}\rightarrow\mathbb{P}\left(  E\right)  $ is the section
of the $\mathbb{P}^{1}$-bundle $\mathbb{P}\left(  E\right)  $ corresponding to
the natural embedding $\mathcal{O}_{X_{0}}\hookrightarrow\mathcal{O}_{X_{0}%
}\oplus L$. Hence $X=C\left(  X_{0}\right)  $ is a cone over $X_{0}$.

We now calculate the discrepancy: $\pi$ is the
\index{blowup}%
blowup of $X$ in the singular point of $X$ and with exceptional locus
$D=\sigma\left(  X_{0}\right)  \cong X_{0}$. So
\[
\mathcal{O}_{Y}\left(  D\right)  \mid_{D}=\mathcal{N}_{D/Y}=L^{-1}%
\]
Write%
\[
K_{Y}=\pi^{\ast}K_{X}\otimes\mathcal{O}_{Y}\left(  D\right)  ^{a}%
\]
and restrict to $D$%
\[
K_{Y}\mid_{D}=\mathcal{O}_{Y}\left(  D\right)  ^{a}\mid_{D}=L^{-a}%
\]
The
\index{adjunction formula}%
adjunction formula yields%
\[
L^{-\frac{k}{l}}=K_{D}=\left(  K_{Y}\otimes\mathcal{O}_{Y}\left(  D\right)
\right)  \mid_{D}=L^{-a}\otimes L^{-1}=L^{-\left(  a+1\right)  }%
\]
so $a=\frac{k}{l}-1$.

Now consider the case of a smooth quadric $X_{0}\cong\mathbb{P}^{1}%
\times\mathbb{P}^{1}\subset\mathbb{P}^{3}$. Then we can write $X=S\left(
1,1,0\right)  $ as%
\[
X=\left\{  \det\left(
\begin{array}
[c]{cc}%
y_{0} & y_{2}\\
y_{1} & y_{3}%
\end{array}
\right)  =0\right\}  \subset\mathbb{P}^{4}%
\]
so $P=\left(  0:0:0:1\right)  $ is the singular point of $X$. The
\index{discrepancy}%
discrepancy is%
\[
K_{Y}-\pi^{\ast}K_{X}=D
\]
We now calculate a small and hence crepant resolution of $X=C\left(
X_{0}\right)  $. Let%
\begin{align*}
E_{1}  &  :=\mathcal{O}_{\mathbb{P}^{1}}\left(  2\right)  \oplus
\mathcal{O}_{\mathbb{P}^{1}}\left(  2\right)  \oplus\mathcal{O}_{\mathbb{P}%
^{1}}\left(  1\right) \\
E_{2}  &  :=\mathcal{O}_{\mathbb{P}^{1}}\left(  1\right)  \oplus
\mathcal{O}_{\mathbb{P}^{1}}\left(  1\right)  \oplus\mathcal{O}_{\mathbb{P}%
^{1}}%
\end{align*}
The maps from $\mathbb{P}\left(  E_{1}\right)  =\mathbb{P}\left(
E_{2}\right)  $ to $\mathbb{P}\left(  H^{0}\left(  \mathbb{P}\left(
E_{i}\right)  ,\mathcal{O}_{\mathbb{P}\left(  E_{i}\right)  }\left(  1\right)
\right)  \right)  =\mathbb{P}\left(  H^{0}\left(  \mathbb{P}^{1},E_{i}\right)
\right)  $ give rise to a diagram
\[%
\begin{tabular}
[c]{lll}%
$\mathbb{P}\left(  \mathcal{O}_{\mathbb{P}^{1}}\left(  2\right)
\oplus\mathcal{O}_{\mathbb{P}^{1}}\left(  2\right)  \oplus\mathcal{O}%
_{\mathbb{P}^{1}}\left(  1\right)  \right)  $ & $\overset{\sim}{\rightarrow}$
& $S\left(  2,2,1\right)  =:Y_{small}$\\
\multicolumn{1}{c}{$\shortparallel$} &  & \\
$\mathbb{P}\left(  \mathcal{O}_{\mathbb{P}^{1}}\left(  1\right)
\oplus\mathcal{O}_{\mathbb{P}^{1}}\left(  1\right)  \oplus\mathcal{O}%
_{\mathbb{P}^{1}}\right)  $ & $\rightarrow$ & $S\left(  1,1,0\right)  =X$%
\end{tabular}
\]
and hence to a morphism $Y_{small}\rightarrow X$. With%
\[
S\left(  2,2,1\right)  =\left\{  \operatorname*{minors}\left(  2,\left(
\begin{array}
[c]{ccccc}%
x_{0} & x_{1} & x_{3} & x_{4} & x_{6}\\
x_{1} & x_{2} & x_{4} & x_{5} & x_{7}%
\end{array}
\right)  \right)  =0\right\}  \subset\mathbb{P}^{7}%
\]
a morphism $g:Y_{small}\rightarrow X$ is given by
\[
g\left(  x_{0}:...:x_{7}\right)  =\left(  x_{0}:x_{1}:x_{3}:x_{4}%
:x_{6}\right)
\]
and the exceptional locus is $\mathbb{P}^{1}$.
\end{example}

Finally there are also Calabi-Yau varieties, which do not have
\index{crepant resolution}%
crepant desingularizations. Nevertheless we want to have a notion of
\index{mirror symmetry}%
mirror symmetry also for them.

To resolve these issues, the idea is to define so called
\index{stringy Hodge numbers}%
stringy Hodge numbers $h_{st}^{p,q}\left(  X\right)  $ for singular varieties.
They should coincide with the usual
\index{Hodge numbers}%
Hodge numbers for smooth varieties, if there is a
\index{crepant resolution}%
crepant desingularization $Y\rightarrow X$ they should coincide with the
\index{Hodge numbers}%
Hodge numbers of $Y$, and even if there is no
\index{crepant resolution}%
crepant desingularization there should still be a notion of mirror symmetry.

One can even consider an enlarged class of varieties for which there is in
general no notion of
\index{stringy Hodge numbers}%
stringy
\index{Hodge numbers}%
Hodge numbers, but as there is a (not necessarily polynomial)
\index{generating function}%
generating function encoding equivalent information, there is still a notion
of
\index{mirror symmetry}%
mirror symmetry.

In the following let $X$ be an irreducible normal algebraic variety of
dimension $d$ over $\mathbb{C}$.

\subsubsection{The Hodge weight filtration and the $E$-polynomial}

The cohomology groups $H^{k}\left(  X,\mathbb{Q}\right)  $ of a complex
algebraic variety $X$ carry a natural
\index{mixed Hodge structure|textbf}%
mixed Hodge structure \cite{Deligne Theorie de Hodge II}, \cite{Deligne
Theorie de Hodge III}, which is given by the following data:

An increasing filtration
\[
0=W_{-1}\subset W_{0}\subset...\subset W_{2k}=H^{k}\left(  X,\mathbb{Q}%
\right)
\]
on $H^{k}\left(  X,\mathbb{Q}\right)  $ called
\index{weight filtration|textbf}%
weight filtration, and a decreasing filtration
\[
H^{k}\left(  X,\mathbb{C}\right)  =F^{0}\supset F^{1}\supset...\supset
F^{k}\supset F^{k+1}=0
\]
on $H^{k}\left(  X,\mathbb{C}\right)  =H^{k}\left(  X,\mathbb{Q}\right)
\otimes\mathbb{C}$ called
\index{Hodge filtration|textbf}%
Hodge filtration. We then have%
\[
H^{p,q}\left(  H^{k}\left(  X,\mathbb{C}\right)  \right)  =F^{p}Gr_{p+q}%
H^{k}\left(  X,\mathbb{C}\right)  \cap\overline{F^{q}Gr_{p+q}H^{k}\left(
X,\mathbb{C}\right)  }%
\]
where%
\begin{align*}
Gr_{l}H^{k}\left(  X,\mathbb{Q}\right)   &  :=\left(  W_{l}/W_{l-1}\right) \\
F^{p}Gr_{l}H^{k}\left(  X,\mathbb{C}\right)   &  :=\operatorname{Im}\left(
F^{p}\cap\left(  W_{l}\otimes\mathbb{C}\right)  \rightarrow Gr_{l}H^{k}\left(
X,\mathbb{Q}\right)  \otimes\mathbb{C}\right)
\end{align*}
and the filtrations have the property that $F^{p}Gr_{l}H^{k}\left(
X,\mathbb{C}\right)  $ gives a (pure) Hodge structure of weight $l$ on
$Gr_{l}H^{k}\left(  X,\mathbb{Q}\right)  $.

We therefore have a decomposition%
\[
H^{k}\left(  X,\mathbb{C}\right)  =\bigoplus_{p,q}H^{p,q}\left(  H^{k}\left(
X,\mathbb{C}\right)  \right)
\]

In \cite{DK Newton polyhedra and an algorithm for computing HodgeDeligne
numbers} one can find a proof that also the cohomology with compact support
$H_{c}^{i}\left(  X,\mathbb{Q}\right)  $ admits a mixed Hodge structure.

\begin{definition}
The $E$\textbf{-polynomial} $E\left(  X;u,v\right)  \in\mathbb{Q}\left[
u,v\right]  $ (coefficients in $\mathbb{Z}$)
\index{E-polynomial|textbf}%
of a complex normal algebraic variety $X$ of dimension $d$ is then defined as%
\[
E\left(  X;u,v\right)  :=\sum_{0\leq p,q\leq d}\sum_{0\leq i\leq2d}\left(
-1\right)  ^{i}h^{p,q}\left(  H_{c}^{i}\left(  X\right)  \right)  u^{p}v^{q}%
\]
So we have a map from the category of normal algebraic varieties
$\mathcal{V}_{\mathbb{C}}$ to $\mathbb{Q}\left[  u,v\right]  $ by%
\[
E:ob\mathcal{V}_{\mathbb{C}}\rightarrow\mathbb{Q}\left[  u,v\right]  ,\text{
}X\mapsto E\left(  X;u,v\right)
\]
associating to each $X$ its $E$ polynomial.
\end{definition}

Important properties of the $E$-polynomial:

\begin{proposition}
Let $X$ and $X_{i}$ be complex normal algebraic variety.

\begin{enumerate}
\item If $X=\bigcup_{i}X_{i}$ is
\index{stratification}%
stratified by a disjoint union of locally closed subvarieties then%
\[
E\left(  X\right)  =\sum_{i}E\left(  X_{i}\right)
\]

\item
\[
E\left(  X_{1}\times X_{2}\right)  =E\left(  X_{1}\right)  \cdot E\left(
X_{2}\right)
\]

\item If $X\rightarrow B$ is a
\index{locally trivial fibration}%
locally trivial fibration and $F$ the fiber over the closed point then%
\[
E\left(  X\right)  =E\left(  F\right)  \cdot E\left(  B\right)
\]

\end{enumerate}
\end{proposition}

A proof can be found in the previously mentioned paper\linebreak\cite{DK
Newton polyhedra and an algorithm for computing HodgeDeligne numbers}. Note
that the number of $\mathbb{F}_{q}$-points of a variety has similar properties
as $E$.

\begin{remark}
For smooth compact $X$ of dimension $d$%
\[
E\left(  X;u,v\right)  =\sum_{0\leq p,q\leq d}\left(  -1\right)  ^{p+q}%
h^{p,q}\left(  X\right)  u^{p}v^{q}%
\]
with $h^{p,q}\left(  X\right)  =\dim H_{\bar{\partial}}^{p,q}\left(  X\right)
=\dim H^{q}\left(  X,\Omega_{X}^{p}\right)  $

\begin{itemize}
\item Hodge
\index{Hodge duality}%
duality for $X$ is equivalent to
\[
E\left(  X;u,v\right)  =E\left(  X;v,u\right)
\]

\item Poincar\'{e}
\index{Poincare duality}%
duality for $X$ is equivalent to
\[
E\left(  X;u,v\right)  =\left(  uv\right)  ^{d}E\left(  X;u^{-1}%
,v^{-1}\right)
\]

\item Topological
\index{topological mirror symmetry}%
mirror symmetry for a pair of varieties $X$ and $X^{\circ}$ is equivalent to
\[
E\left(  X;u,v\right)  =\left(  -u\right)  ^{d}E\left(  X^{\circ}%
;u^{-1},v\right)
\]

\end{itemize}
\end{remark}

\begin{remark}
Consider a
\index{stratification}%
stratification $X=U\cup C$ with $X$ and $C$ compact. The long exact sequence
for cohomology with compact support reads as%
\[
...\rightarrow H_{c}^{k}\left(  U\right)  \overset{\varphi_{k}}{\rightarrow
}H^{k}\left(  X\right)  \overset{\psi_{k}}{\rightarrow}H^{k}\left(  C\right)
\overset{\delta_{k}}{\rightarrow}H_{c}^{k+1}\left(  U\right)  \rightarrow...
\]
where $\varphi_{k}$ is given by continuation by $0$, $\psi_{k}$ is given by
restriction and the boundary map $\delta_{k}$ is given by $\varpi\mapsto
d\left(  \beta\cdot r^{\ast}\varpi\right)  $ where $r$ is the retract of a
tubular neighborhood of $C$ and $\beta$ is a bump function on this neighborhood.
\end{remark}

We illustrate this with the following two examples:

\begin{example}
For $X=\mathbb{P}^{1}$, $U=\mathbb{C}$ and $C=\left\{  pt\right\}  $%
\[%
\begin{tabular}
[c]{l|lllll}%
$k$ & $H_{c}^{k}\left(  \mathbb{C}\right)  $ & $\rightarrow$ & $H^{k}\left(
\mathbb{P}^{1}\right)  $ & $\rightarrow$ & $H^{k}\left(  pt\right)  $\\\hline
$2$ & $uv$ &  & $uv$ &  & $0$\\
$1$ & $0$ &  & $0$ &  & $0$\\
$0$ & $0$ &  & $1$ &  & $1$%
\end{tabular}
\]
where we denote the
\index{Hodge filtration}%
Hodge filtration by the corresponding $E$ monomials. The $E$-polynomials
\index{E-polynomial}%
are%
\[%
\begin{tabular}
[c]{lll}%
$E\left(  \mathbb{C}\right)  =uv$ & $E\left(  \mathbb{P}^{1}\right)  =1+uv$ &
$E\left(  pt\right)  =1$%
\end{tabular}
\]
Remark: The long exact sequence decomposes in short ones if all varieties have
only even cohomology.
\end{example}

\begin{example}
\label{elliptic curve E polynomial}For $X=\mathbb{P}^{3}$, $C$ an elliptic
curve and $U=\mathbb{P}^{3}-C$ we have%
\[%
\begin{tabular}
[c]{l|lllll}%
$k$ & $H_{c}^{k}\left(  U\right)  $ & $\rightarrow$ & $H^{k}\left(
\mathbb{P}^{3}\right)  $ & $\rightarrow$ & $H^{k}\left(  C\right)  $\\\hline
$6$ & $\left(  uv\right)  ^{3}$ &  & $\left(  uv\right)  ^{3}$ &  & $0$\\
$5$ & $0$ &  & $0$ &  & $0$\\
$4$ & $\left(  uv\right)  ^{2}$ &  & $\left(  uv\right)  ^{2}$ &  & $0$\\
$3$ & $0$ &  & $0$ &  & $0$\\
$2$ & $u+v$ &  & $uv$ &  & $uv$\\
$1$ & $0$ &  & $0$ &  & $-u-v$\\
$0$ & $0$ &  & $1$ &  & $1$%
\end{tabular}
\]
The $E$-polynomials
\index{E-polynomial}%
are%
\begin{align*}
E\left(  U\right)   &  =u+v+\left(  uv\right)  ^{2}+\left(  uv\right)  ^{3}\\
E\left(  \mathbb{P}^{3}\right)   &  =1+uv+\left(  uv\right)  ^{2}+\left(
uv\right)  ^{3}\\
E\left(  C\right)   &  =1-u-v+uv
\end{align*}
So a shift in the cohomological weight occurs (Notice also the sign of the
$u+v$ term).
\end{example}

\begin{example}
Continuing Example \ref{elliptic curve E polynomial} the corresponding
\index{Hodge filtration}%
Hodge filtration for the cohomology of $U$:%
\begin{align*}
&
%
$}%
\end{tabular}
\\
0  &  =W_{-1}=...=W_{3}\subset W_{4}=...=W_{8}=H^{4}\left(  U,\mathbb{Q}%
\right)
\end{align*}
and similar for $k=6$.
\end{example}

\begin{example}
\label{ExampleEfunction2}We continue Example \ref{ExampleEfunction1} of the
cone over the quadric calculating the $E$-polynomials:

The
\index{cohomology ring}%
cohomology ring $H^{\ast}\left(  X_{0}\right)  $ of $X_{0}$ is generated by
$h_{i}=pr_{i}^{\ast}c_{1}\left(  \mathbb{P}^{1}\right)  $, $i=1,2$ and hence
$1,h_{1},h_{2},h_{1}h_{2}$ form a basis as a vector space, so%
\[
E\left(  X_{0}\right)  =1+2uv+\left(  uv\right)  ^{2}%
\]
which agrees with the product formula $E\left(  X_{0}\right)  =E\left(
\mathbb{P}^{1}\right)  ^{2}=\left(  1+uv\right)  ^{2}$.

$H^{\ast}\left(  Y\right)  $ is a free module over $H^{\ast}\left(
X_{0}\right)  $ with basis $1,c=c_{1}\left(  \mathcal{O}_{Y}\left(  1\right)
\right)  $ and hence $1,c,h_{1},h_{2},ch_{1},ch_{2},h_{1}h_{2},ch_{1}h_{2}$ is
a vector space basis (where $h_{i}$ is short for $\pi^{\ast}h_{i}$), so%
\[
E\left(  Y\right)  =1+3uv+3\left(  uv\right)  ^{2}+\left(  uv\right)  ^{3}%
\]

$H^{\ast}\left(  Y_{small}\right)  $ is a free module over $H^{\ast}\left(
\mathbb{P}^{1}\right)  $ with basis $1,c,c^{2}$ with $c=c_{1}\left(
\mathcal{O}_{Y_{small}}\left(  1\right)  \right)  $ and hence $1,h,c,c^{2}%
,ch,hc^{2}$ is a vector space basis ($h=\pi^{\ast}c_{1}\left(  \mathbb{P}%
^{1}\right)  $), so%
\[
E\left(  Y_{small}\right)  =1+2uv+2\left(  uv\right)  ^{2}+\left(  uv\right)
^{3}%
\]

So
\index{E-polynomial}%
the $E$ polynomials%
\begin{align*}
E\left(  Y\backslash X_{0}\right)   &  =E\left(  Y\right)  -E\left(
X_{0}\right)  =\left(  1+3uv+3\left(  uv\right)  ^{2}+\left(  uv\right)
^{3}\right)  -\left(  1+2uv+\left(  uv\right)  ^{2}\right) \\
&  =uv+2\left(  uv\right)  ^{2}+\left(  uv\right)  ^{3}\\
E\left(  Y_{small}\backslash\mathbb{P}^{1}\right)   &  =E\left(
Y_{small}\right)  -E\left(  \mathbb{P}^{1}\right)  =\left(  1+2uv+2\left(
uv\right)  ^{2}+\left(  uv\right)  ^{3}\right)  -\left(  1+uv\right) \\
&  =uv+2\left(  uv\right)  ^{2}+\left(  uv\right)  ^{3}%
\end{align*}
agree as expected because of $Y\backslash X_{0}\cong X\backslash P\cong
Y_{small}\backslash\mathbb{P}^{1}$. Using this we can also calculate%
\[
E\left(  X\right)  =1+uv+2\left(  uv\right)  ^{2}+\left(  uv\right)  ^{3}%
\]

\end{example}

\subsubsection{Varieties with canonical singularities}

\begin{definition}
\label{CanonicalLogTerminalSingularities}Let $X$ be a normal projective
variety $X$, which
\index{Q-Gorenstein|textbf}%
is $\mathbb{Q}$-Gorenstein, i.e., $K_{X}\in Div\left(  X\right)
\otimes\mathbb{Q}$, and let $f:Y\rightarrow X$ be a resolution of
singularities such that the exceptional locus of $f$ is a divisor $E$, whose
irreducible components $D_{1},...,D_{r}$ are smooth divisors with only
\index{simple normal crossings}%
simple normal crossings, and let $K_{Y}=f^{\ast}K_{X}+\sum_{i=1}^{r}a_{i}%
D_{i}$. Then $X$ is said to have

\begin{itemize}
\item \textbf{terminal}
\index{terminal singularities|textbf}%
singularities if $a_{i}>0$ for all $i$

\item \textbf{canonical }%
\index{canonical singularities|textbf}%
singularities if $a_{i}\geq0$ for all $i$

\item \textbf{log-terminal }%
\index{log terminal singularities|textbf}%
singularities if $a_{i}>-1$ for all $i$

\item \textbf{log-canonical }%
\index{log canonical singularities|textbf}%
singularities if $a_{i}\geq-1$ for all $i$.
\end{itemize}
\end{definition}

\subsubsection{The stringy $E$-function\label{Stingy E function}}

In the following, we consider a normal projective $d$-dimensional variety $X$
with log-terminal singularities, let $f:Y\rightarrow X$ be a resolution of
singularities and $D_{1},...,D_{r}$ the smooth components of the exceptional
locus with only
\index{simple normal crossings}%
simple normal crossings, and $K_{Y}=f^{\ast}K_{X}+\sum_{i=1}^{r}a_{i}D_{i}$.

Let $I=\left\{  1,...,r\right\}  $ and set for any $J\subset I$%
\begin{align*}
D_{J}  &  =Y\cap\bigcap_{j\in J}D_{j}\\
D_{J}^{\circ}  &  =D_{J}\backslash\bigcup_{i\in I\backslash J}D_{i}%
\end{align*}
This gives a
\index{stratification}%
stratification $D_{J}=\bigcup_{J^{\prime},J\subset J^{\prime}}D_{J^{\prime}%
}^{\circ}$.

\begin{definition}
Define the
\index{stringy E-function|textbf}%
\textbf{stringy }$E$\textbf{-function }$E_{st}$ of $X$ as%
\[
E_{st}\left(  X;u,v\right)  :=\sum_{J\subset I}E\left(  D_{J}^{\circ
};u,v\right)  \prod_{j\in J}\frac{uv-1}{\left(  uv\right)  ^{a_{j}+1}-1}%
\]

\end{definition}

\begin{remark}
If $X$ is
\index{Gorenstein}%
Gorenstein, then the $a_{j}\in\mathbb{Z}_{\geq0}$ and hence $E_{st}\left(
X;u,v\right)  \in\mathbb{Z}\left[  \left[  u,v\right]  \right]  \cap
\mathbb{Q}\left(  u,v\right)  $. $E_{st}\left(  X;u,v\right)  $ is not a
rational function in general.
\end{remark}

The following key theorem by
\index{Batyrev}%
Batyrev \cite{Batyrev Stringy Hodge numbers of varieties with Gorenstein
canonical singularities. Integrable systems and algebraic geometry KobeKyoto
1997}, using ideas by Kontsevich and Denef and Loeser, assures that
$E_{st}\left(  X;u,v\right)  $ is well defined. See also \cite{DL Germs of
arcs on singular algebraic varieties and motivic integration}.

\begin{theorem}
\label{Est_independent}$E_{st}\left(  X;u,v\right)  $ does not depend on the
resolution $f:Y\rightarrow X$, in particular, $E_{st}\left(  X;u,v\right)  $
is well defined.
\end{theorem}

This is also true in the case of
\index{log terminal singularities}%
log terminal singularities. As direct corollary, we have:

\begin{corollary}
If $X$ is smooth, then $E_{st}\left(  X;u,v\right)  =E\left(  X;u,v\right)  $.
\end{corollary}

\begin{remark}
First make an easy but important observation: $E_{st}$ is not affected by the
\index{blowup}%
blowup $f:Y\rightarrow X$ of a point $P$ in smooth $X$: The exceptional locus
of $f$ is $D=\mathbb{P}^{d-1}$ and the
\index{discrepancy}%
discrepancy is%
\[
K_{Y}-f^{\ast}K_{X}=\left(  d-1\right)  D
\]%
\begin{align*}
E_{st}\left(  X\right)   &  =E\left(  Y\backslash D\right)  +E\left(
D\right)  \frac{uv-1}{\left(  uv\right)  ^{a_{1}+1}-1}=E\left(  Y\backslash
D\right)  +E\left(  \mathbb{P}^{d-1}\right)  \frac{uv-1}{\left(  uv\right)
^{d}-1}\\
&  =E\left(  Y\backslash D\right)  +\left(  1+uv+...+\left(  uv\right)
^{d-1}\right)  \frac{uv-1}{\left(  uv\right)  ^{d}-1}\\
&  =E\left(  Y\backslash D\right)  +1=E\left(  X\right)
\end{align*}

\end{remark}

\begin{remark}
The idea of the proof of Theorem \ref{Est_independent}, considering for
simplicity Gorenstein canonical singularities, is the following (see
\cite{Batyrev Stringy Hodge numbers of varieties with Gorenstein canonical
singularities. Integrable systems and algebraic geometry KobeKyoto 1997},
\cite{DL Germs of arcs on singular algebraic varieties and motivic
integration} and reviews in \cite{Blicke Jet schemes and motivic integration}
and \cite{Craw An introduction to motivic integration}):

Consider the
\index{Grothendieck ring|textbf}%
\textbf{Grothendieck ring of complex algebraic varieties} $\mathcal{M}$, which
is the free abelian group of isomorphism classes of complex algebraic
varieties modulo the subgroup generated by $\left[  X\right]  -\left[
V\right]  -\left[  X-V\right]  $ for closed subsets $V\subset X$, with a ring
structure given by%
\[
\left[  X\right]  \cdot\left[  X^{\prime}\right]  =\left[  X\times X^{\prime
}\right]
\]
Call the neutral element $\left[  point\right]  =:1$ and $\left[
\mathbb{C}\right]  =:\mathbb{L}$.

The map $\left[  -\right]  :ob\mathcal{V}_{\mathbb{C}}\rightarrow\mathcal{M}$
is the universal map being additive on disjoint unions of constructible sets
(i.e., a finite disjoint union of locally closed subvarieties with respect to
the Zariski topology) and multiplicative on products, so any other map
$E:ob\mathcal{V}_{\mathbb{C}}\rightarrow\mathbb{Q}\left[  u,v\right]  $ with
the same properties factors through $\left[  -\right]  $. So the universality
of $\left[  -\right]  $ gives a factorization of $E:\mathcal{V}_{\mathbb{C}%
}\rightarrow\mathbb{Q}\left[  u,v\right]  $ through the Grothendieck ring%
\[%
\begin{tabular}
[c]{lcl}%
$ob\mathcal{V}_{\mathbb{C}}$ & $\overset{E}{\longrightarrow}$ & $\mathbb{Q}%
\left[  u,v\right]  $\\
\multicolumn{1}{r}{$\left[  -\right]  \searrow$} & \multicolumn{1}{l}{} &
$\nearrow E$\\
& $\mathcal{M}$ &
\end{tabular}
\]

The goal is to write%
\[
E_{st}\left(  X;u,v\right)  =\sum_{J\subset I}E\left(  D_{J}^{\circ
};u,v\right)  \prod_{j\in J}\frac{uv-1}{\left(  uv\right)  ^{a_{j}+1}%
-1}=E\left(  \int_{J_{\infty}\left(  Y\right)  }F_{D}d\mu\mathbb{L}%
^{d}\right)
\]
for a suitable function $F_{D}$ associated to the
\index{discrepancy}%
discrepancy divisor ($J_{\infty}\left(  Y\right)  $ is the bundle of formal
arcs on $Y$), after extending $E$ to $K_{0}\left(  \mathcal{V}_{\mathbb{C}%
}\right)  \left[  \mathbb{L}^{-1}\right]  $ and then to an appropriate
completion.
\index{transformation rule of motivic integration}%
The transformation rule for motivic integration implies that the motivic
integral does not depend on the resolution.
\end{remark}

\subsubsection{Stringy Hodge numbers}

\begin{theorem}
[Poincar\'{e} Duality]\cite{Batyrev Stringy Hodge numbers of varieties with
Gorenstein canonical singularities. Integrable systems and algebraic geometry
KobeKyoto 1997} $E_{st}\left(  X;u,v\right)  $ has the following properties:%
\begin{align*}
E_{st}\left(  X;u,v\right)   &  =\left(  uv\right)  ^{d}E_{st}\left(
X;u^{-1},v^{-1}\right) \\
E_{st}\left(  X;0,0\right)   &  =1
\end{align*}

\end{theorem}

This is also true in the case of log terminal singularities.

\begin{corollary}
If $E_{st}$ is a polynomial, then $\deg\left(  E\right)  =2d$.
\end{corollary}

\begin{definition}
If $E_{st}$ is a polynomial, the
\index{stringy Hodge numbers|textbf}%
\textbf{stringy Hodge numbers} of $X$ are defined as%
\[
h_{st}^{p,q}\left(  X\right)  =\left(  -1\right)  ^{p+q}\operatorname*{coeff}%
\left(  E_{st},u^{p}v^{q}\right)
\]

\end{definition}

So $h_{st}^{p,q}\left(  X\right)  =0$ outside the
\index{Hodge diamond}%
Hodge diamond and $h_{st}^{0,0}\left(  X\right)  =h_{st}^{d,d}\left(
X\right)  =1$.

\subsubsection{Crepant resolutions and mirror symmetry}

\begin{theorem}
\label{TheoremCrepant}\cite{Batyrev Stringy Hodge numbers of varieties with
Gorenstein canonical singularities. Integrable systems and algebraic geometry
KobeKyoto 1997} If $X$ is $\mathbb{Q}$-Gorenstein with at worst log-terminal
singularities and $f:Y\rightarrow X$ is a projective birational morphism with
$K_{Y}=f^{\ast}K_{X}$ then $E_{st}\left(  X;u,v\right)  =E_{st}\left(
Y;u,v\right)  $. So if $X$ admits a
\index{crepant resolution}%
crepant resolution $f:Y\rightarrow X$ then $E_{st}\left(  X;u,v\right)
=E\left(  Y;u,v\right)  $.
\end{theorem}

\begin{remark}
In particular if $X$ admits a crepant resolution, then $E_{st}\left(
X;u,v\right)  $ is polynomial and hence the
\index{stringy Hodge numbers}%
stringy Hodge numbers of $X$ exist. If $E_{st}\left(  X;u,v\right)  $ is not
polynomial, then $X$ admits no
\index{crepant resolution}%
crepant resolution.
\end{remark}

\begin{definition}
Two Calabi-Yau varieties $X$ and $X^{\circ}$ are called
\index{stringy topological mirror pair|textbf}%
\textbf{stringy topological mirror pair} if their
\index{stringy E-function}%
stringy $E$-functions satisfy%
\[
E_{st}\left(  X;u,v\right)  =\left(  -u\right)  ^{d}E_{st}\left(  X^{\circ
};u^{-1},v\right)
\]

\end{definition}

This is well defined even in the case when $E_{st}$ is not polynomial.

\begin{example}
Now we return to the Example \ref{ExampleEfunction1} and
\ref{ExampleEfunction2}: Let $X_{0}\subset\mathbb{P}^{d}$ be a smooth quadric
($k=d-1,$ $l=1$) and
\[
\mathbb{P}\left(  \mathcal{O}_{X_{0}}\left(  1\right)  \oplus\mathcal{O}%
_{X_{0}}\right)  =Y\rightarrow X=C\left(  X_{0}\right)
\]
with discrepancy divisor $\left(  d-2\right)  D$ with $D\cong X_{0}$.

For $d=3$ we had $X=S\left(  1,1,0\right)  $, we computed a small resolution
\[
S\left(  2,2,1\right)  =Y_{small}\rightarrow X
\]
and calculated%
\begin{align*}
E\left(  Y_{small}\right)   &  =1+2uv+2\left(  uv\right)  ^{2}+\left(
uv\right)  ^{3}\\
E\left(  Y\right)   &  =1+3uv+3\left(  uv\right)  ^{2}+\left(  uv\right)
^{3}\\
E\left(  D\right)   &  =1+2uv+\left(  uv\right)  ^{2}\\
E\left(  Y\backslash D\right)   &  =uv+2\left(  uv\right)  ^{2}+\left(
uv\right)  ^{3}%
\end{align*}
So the
\index{stringy E-function}%
stringy $E$ function $E_{st}$ is%
\begin{align*}
E_{st}\left(  X\right)   &  =E\left(  D_{\emptyset}^{\circ}\right)  +E\left(
D_{\left\{  1\right\}  }^{\circ}\right)  \frac{uv-1}{\left(  uv\right)
^{2}-1}=E\left(  Y\backslash D\right)  +E\left(  D\right)  \frac{1}{uv+1}\\
&  =\left(  uv+2\left(  uv\right)  ^{2}+\left(  uv\right)  ^{3}\right)
+\left(  1+uv\right)  ^{2}\frac{1}{uv+1}\\
&  =1+2uv+2\left(  uv\right)  ^{2}+\left(  uv\right)  ^{3}=E\left(
Y_{small}\right)
\end{align*}
and, as predicted by Theorem \ref{TheoremCrepant}, the
\index{stringy Hodge numbers}%
stringy Hodge numbers of $X$ indeed coincide with the Hodge numbers of the
\index{small resolution}%
small resolution.

$E_{st}\left(  X\right)  $ is not a polynomial for $d>3$, in particular $X$
does not admit a
\index{crepant resolution}%
crepant resolution: As $X\backslash p\cong Y\backslash D$ is isomorphic to the
totalspace of $L\rightarrow X_{0}$, we have%
\[
E\left(  X\backslash p\right)  =E\left(  Y\backslash D\right)  =\left(
uv\right)  E\left(  X_{0}\right)
\]
and
\[
E_{st}\left(  X\right)  =\left(  uv\right)  E\left(  D\right)  +E\left(
D\right)  \frac{uv-1}{\left(  uv\right)  ^{d-1}-1}=E\left(  D\right)
\frac{\left(  uv\right)  ^{d}-1}{\left(  uv\right)  ^{d-1}-1}%
\]
As%
\[
E\left(  D\right)  =\left\{
\begin{tabular}
[c]{ll}%
$\frac{\left(  \left(  uv\right)  ^{\frac{d-1}{2}}+1\right)  \left(  \left(
uv\right)  ^{\frac{d+1}{2}}-1\right)  }{uv-1}$ & for $d$ odd\\
$\frac{\left(  uv\right)  ^{d}-1}{uv-1}$ & for $d$ even
\end{tabular}
\right\}
\]
the
\index{stringy E-function}%
stringy $E$ function $E_{st}\left(  X\right)  $ is not a polynomial.
\end{example}

\subsubsection{Birational Calabi-Yau manifolds
\index{Birational Calabi-Yau manifolds}%
}

\begin{theorem}
\cite{Batyrev Birational CalabiYau nfolds have equal Betti numbers},
\cite{Batyrev Stringy Hodge numbers of varieties with Gorenstein canonical
singularities. Integrable systems and algebraic geometry KobeKyoto 1997}
Birational Calabi-Yau manifolds have equal
\index{Hodge numbers}%
Hodge numbers.
\end{theorem}

Actually one proves that $\left[  X_{1}\right]  =\left[  X_{2}\right]  $,
i.e., $X_{1}$ and $X_{2}$ represent the same class in the
\index{Grothendieck ring}%
Grothendieck ring $\mathcal{M}$.

\subsection{Some facts and notations from toric
geometry\label{Sec facts from toric geometry}}

The key example of a toric variety is the projective space $\mathbb{P}^{n}$.
Let $\left(  y_{0}:...:y_{n}\right)  $ be the homogeneous coordinates. On the
open set $U_{i}=\left\{  y\in\mathbb{P}^{n}\mid y_{i}\neq0\right\}  $ the
functions
\[
x_{k}^{i}=\frac{y_{k}}{y_{i}}%
\]
give an isomorphism
\begin{align*}
U_{i}  &  \rightarrow\mathbb{A}^{n}\\
\left(  y_{0}:...:y_{n}\right)   &  \mapsto\left(  x_{0}^{i},...,\widehat
{x_{i}^{i}},...,x_{n}^{i}\right)
\end{align*}
Considering another chart $U_{j}\rightarrow\mathbb{A}^{n}$, on $U_{i}\cap
U_{j}$%
\[
x_{k}^{j}=\frac{y_{k}}{y_{j}}=\frac{y_{k}}{y_{i}}\frac{y_{i}}{y_{j}}=x_{k}%
^{i}\left(  x_{j}^{i}\right)  ^{-1}%
\]
i.e., the coordinate functions in one chart are given as
\index{Laurent monomial|textbf}%
Laurent monomials (i.e., monomials which also can have negative exponents) in
the coordinates of the other chart, a key property of toric varieties.

\subsubsection{Affine toric varieties\label{Sec affine toric varieties}}

If $S\subset M=\mathbb{Z}^{n}$ is \newsym[$M$]{lattice of torus characters}{}a
finitely generated commutative semigroup with $0$, we can associate to $S$
its
\index{semigroup algebra}%
\textbf{semigroup algebra} $\mathbb{C}\left[  S\right]  $, consisting of all
finite formal sums $\sum_{m\in S}a_{m}x^{m}$, $a_{m}\in\mathbb{C}$ with
multiplication $x^{m}\cdot x^{m^{\prime}}=x^{m+m^{\prime}}$.

\begin{example}
\label{examplecone1}The semigroup algebra of $S=\left\langle \left(
1,0\right)  ,\left(  1,1\right)  ,\left(  1,2\right)  \right\rangle
\subset\mathbb{Z}^{2}$ is $\mathbb{C}\left[  S\right]  =\mathbb{C}\left[
x,xy,xy^{2}\right]  $.
\end{example}

To the semigroup algebra we can
\newsym[$U\left(  \sigma\right)  $]{affine toric variety associated to $\sigma$}{}associate
an
\index{affine toric variety|textbf}%
\textbf{affine toric variety} $Y=\operatorname*{Spec}\mathbb{C}\left[
S\right]  $. Considering $\mathbb{C=C}^{\ast}\cup\left\{  0\right\}  $ as a
semigroup with respect to multiplication, the maximal points of $Y$ are the
semigroup homomorphisms $\operatorname*{Hom}_{sg}\left(  S,\mathbb{C}\right)
$. If $y\in\operatorname*{Hom}_{sg}\left(  S,\mathbb{C}\right)  $ and
$x^{m}\in S$, then $x^{m}\left(  y\right)  =y\left(  m\right)  $.

For generators $m_{1},...,m_{r}$ of $S$, the toric ideal
\index{toric ideal|textbf}%
of $Y$ is the kernel $I_{S}$ of%
\begin{align*}
\mathbb{C}\left[  y_{1},...,y_{r}\right]   &  \rightarrow\mathbb{C}\left[
S\right] \\
y_{i}  &  \mapsto x^{m_{i}}%
\end{align*}
It is given by the
\index{binomial ideal}%
binomial ideal%
\[
I_{S}=\left\langle y^{u^{+}}-y^{u^{-}}\mid u\in\ker\left(  m_{1}%
,...,m_{r}\right)  \right\rangle
\]
where $u=u^{+}-u^{-}$ with $u^{+},u^{-}$ with non-negative entries and
disjoint support (see \cite{Sturmfels Equations defining toric varieties}).

\begin{example}
\label{examplecone2}For $S=\left\langle \left(  1,0\right)  ,\left(
1,1\right)  ,\left(  1,2\right)  \right\rangle $ as in Example
\ref{examplecone1} we have
\[
\ker\left(
\begin{array}
[c]{ccc}%
1 & 1 & 1\\
0 & 1 & 2
\end{array}
\right)  =\left(  1,-2,1\right)  ^{t}%
\]
hence
\[
\mathbb{C}\left[  S\right]  \cong\mathbb{C}\left[  y_{1},y_{2},y_{3}\right]
/\left\langle y_{2}^{2}-y_{1}y_{3}\right\rangle
\]
so $Y=\left\{  y_{2}^{2}-y_{1}y_{3}=0\right\}  $ is a quadric cone.
\end{example}

The inclusion $S\subset M$ gives an embedding of the
\index{torus|textbf}%
torus%
\[
T=\operatorname*{Hom}\nolimits_{\mathbb{Z}}\left(  M,\mathbb{C}^{\ast}\right)
=\left(  \mathbb{C}^{\ast}\right)  ^{n}=\operatorname*{Spec}\left(
\mathbb{C}\left[  \mathbb{Z}^{n}\right]  \right)  \hookrightarrow
\operatorname*{Spec}\left(  \mathbb{C}\left[  S\right]  \right)  =Y
\]
If $t\in T=\operatorname*{Hom}\nolimits_{\mathbb{Z}}\left(  M,\mathbb{C}%
^{\ast}\right)  $ considered as a group homomorphism $t:M\rightarrow
\mathbb{C}^{\ast}$, and $y\in Y$ considered as a semigroup homomorphism
$y:S\rightarrow\mathbb{C}$, then $T$ acts
\index{torus action|textbf}%
on $X$ by%
\[%
\begin{tabular}
[c]{lll}%
$T\times Y$ & $\rightarrow$ & $Y$\\
$\left(  t,y\right)  $ & $\longmapsto$ &
\begin{tabular}
[c]{llll}%
$ty:$ & $S$ & $\rightarrow$ & $\mathbb{C}$\\
& $u$ & $\mapsto$ & $t\left(  u\right)  y\left(  u\right)  $%
\end{tabular}
\end{tabular}
\]
respectively%
\begin{align*}
T\times\mathbb{C}\left[  S\right]   &  \rightarrow\mathbb{C}\left[  S\right]
\\
\left(  t,x^{m}\right)   &  \mapsto t\left(  m\right)  x^{m}%
\end{align*}
for $m\in S$.

\begin{figure}
[h]
\begin{center}
\includegraphics[
height=1.8758in,
width=1.7184in
]%
{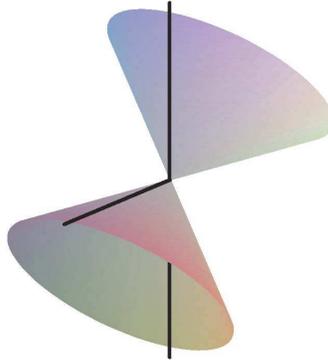}%
\caption{Torus orbits of the quadric cone in $\mathbb{A}^{3}$}%
\label{Fquadriccone}%
\end{center}
\end{figure}

\begin{example}
\label{examplecone3}In Example \ref{examplecone2}, the torus is $\left\{
y_{1}\neq0,y_{2}\neq0,y_{3}\neq0\right\}  \subset X$, i.e., the complement of
the two lines as shown in Figure \ref{Fquadriccone}.
\end{example}

%

We recall some standard facts and notations from polyhedral geometry:

\begin{definition}
A finite intersection of closed half-spaces
\newsym[$M_{\mathbb{R}}$]{$M\otimes_{\mathbb{Z}}\mathbb{R}$}{}in
$M_{\mathbb{R}}=M\otimes_{\mathbb{Z}}\mathbb{R}$ is called a
\index{polyhedron|textbf}%
\textbf{polyhedron}.

A subset $\sigma\subset M_{\mathbb{R}}$ is called a \textbf{polyhedral cone}
if there are $u_{1},...,u_{s}\in M_{\mathbb{R}}$ such that%
\[
\sigma=\left\{  \lambda_{1}u_{1}+...+\lambda_{s}u_{s}\mid\lambda
_{1},...,\lambda_{s}\in\mathbb{R}_{\geq0}\right\}
\]
so any polyhedral cone is a polyhedron.

A polyhedral cone $\sigma\subset M_{\mathbb{R}}$ is
\newsym[$\sigma$]{convex polyhedral cone}{}called
\index{rational polyhedral cone|textbf}%
\textbf{rational polyhedral cone }if there are $u_{1},...,u_{s}\in M$ with
$\sigma=\left\{  \lambda_{1}u_{1}+...+\lambda_{s}u_{s}\mid\lambda
_{1},...,\lambda_{s}\in\mathbb{R}_{\geq0}\right\}  $.

It is called
\index{strongly convex}%
\textbf{strongly convex} if $-\sigma\cap\sigma=\left\{  0\right\}  $ and the
dimension of $\sigma$ is the dimension of the subspace of $M_{\mathbb{R}}$,
spanned by the elements of $\sigma$.

The
\index{convex hull|textbf}%
\textbf{convex hull}
\newsym[$\operatorname*{convexhull}$]{convex hull of points}{}of a subset
$V\subset M_{\mathbb{R}}\cong\mathbb{R}^{n}$ is the intersection of all convex
sets containing $V$. It is denoted by $\operatorname*{convexhull}\left(
V\right)  $ and%
\[
\operatorname*{convexhull}\left(  V\right)  =\left\{  \sum_{i=1}^{s}%
\lambda_{i}v_{i}\mid\lambda_{1},...,\lambda_{s}\geq0\text{ with }\sum
_{i=1}^{s}\lambda_{i}=1\text{ and }v_{1},...,v_{s}\in V\right\}
\]

\end{definition}

\begin{theorem}
[Carath\'{e}odory]If $V\subset\mathbb{R}^{n}$ then any point of
$\operatorname*{convexhull}\left(  V\right)  $ is a convex combination of at
most $n+1$ points of $V$.
\end{theorem}

\begin{definition}
A
\index{polytope|textbf}%
\textbf{polytope }$\Delta\subset M_{\mathbb{R}}$ is the convex hull of a
finite set of points. The dimension of $\Delta$ is the dimension of the
subspace spanned by the points $m-m^{\prime}$ with $m,m^{\prime}\in\Delta$.

A polytope $\Delta$ is called
\index{integral polytope}%
\textbf{integral}
\index{lattice polytope}%
or \textbf{lattice polytope} if it is the convex hull of a finite set of
points in $M$.
\end{definition}

\begin{theorem}
Any bounded polyhedron is a polytope and vice versa.
\end{theorem}

\begin{definition}
A
\index{face|textbf}%
\textbf{face} of a polyhedron $\Delta$ is either $\Delta$ or a subset
$\Delta\cap h$ of $\Delta$, where $h$ is a hyperplane such that $\Delta$ is
contained in one of the closed halfspaces given by $h$. A
\index{facet|textbf}%
\textbf{facet} of $\Delta$ is a codimension one face. Any face of a polyhedron
$\Delta$ is a polyhedron.
\end{definition}

\begin{proposition}
\cite[Sec. 5.3]{GKZ Discriminants Resultants and Multidimensional
Determinants} Let $\tau=\operatorname*{hull}\left(  S\right)  $ be the
rational polyhedral cone in $M_{\mathbb{R}}$ defined as the hull of the
semigroup $S$ introduced above. There is a bijective inclusion respecting map
\[
\left\{  \text{faces of }\tau\right\}  \overset{1:1}{\longrightarrow}\left\{
\text{torus orbit closures in }X\right\}
\]
If $\sigma$ is a face of $\tau$, then the corresponding
\index{torus orbit closure|textbf}%
torus
\index{torus action}%
orbit
\index{torus orbit|textbf}%
is given by $x^{m}=0$ $\forall m\notin S\cap\tau$ and $x^{m}\neq0$ $\forall
m\in S\cap\tau$. The closure of the torus orbit is isomorphic to
$\operatorname*{Spec}\left(  \mathbb{C}\left[  S\cap\tau\right]  \right)  $.
\end{proposition}

\begin{example}
\label{examplecone4}For Example \ref{examplecone2}, the torus orbits in $X$
are given by
\begin{align*}
&  \left\{  y_{1}\neq0,y_{2}\neq0,y_{3}\neq0\right\} \\
&  \left\{  y_{1}\neq0,y_{2}=0,y_{3}=0\right\} \\
&  \left\{  y_{1}=0,y_{2}=0,y_{3}\neq0\right\} \\
&  \left\{  \left(  0,0,0\right)  \right\}
\end{align*}

\end{example}

\subsubsection{Toric varieties from fans\label{Sec Toric varieties from fans}}

Let $N\cong\mathbb{Z}^{n}$
\newsym[$N$]{lattice of $1$-parameter subgroups}{}and let
$M=\operatorname*{Hom}\left(  N,\mathbb{Z}\right)  $
\newsym[$N_{\mathbb{R}}$]{$N\otimes_{\mathbb{Z}}\mathbb{R}$}{}be the
\index{dual lattice|textbf}%
dual lattice of $N$, and denote by
\[
\left\langle -,-\right\rangle :M\times N\rightarrow\mathbb{Z}%
\]
the canonical bilinear pairing. Given a rational
\newsym[$\check{\sigma}$]{dual cone of $\sigma$}{}convex polyhedral cone
$\sigma$ in $N_{\mathbb{R}}$ consider the
\index{dual cone|textbf}%
dual cone
\[
\check{\sigma}=\left\{  m\in M_{\mathbb{R}}\mid\left\langle m,w\right\rangle
\geq0\ \forall w\in\sigma\right\}
\]
of non-negative linear forms on $\sigma$.

\begin{proposition}
[Gordan%
\'{}%
s Lemma]\cite[Sec. 1.1]{Oda Convex Bodies and Algebraic Geometry} If
$\sigma\subset N_{\mathbb{R}}$
\index{Gordan lemma}%
is a rational convex polyhedral cone, then $\check{\sigma}\cap M$ is a
finitely generated semigroup.
\end{proposition}

Given a strongly convex rational polyhedral cone $\sigma$, i.e., $\sigma
\cap\left(  -\sigma\right)  =\left\{  0\right\}  $, we get a finitely
generated semigroup $\check{\sigma}\cap M$ generating $M$ as a group, i.e.,
$\check{\sigma}\cap M+\left(  -\check{\sigma}\cap M\right)  =M$, and hence
\newsym[$U\left(  \sigma\right)  $]{affine toric variety associated to $\sigma$}{}an
\index{affine toric variety}%
affine toric variety $U\left(  \sigma\right)  =\operatorname*{Spec}\left(
\mathbb{C}\left[  \check{\sigma}\cap M\right]  \right)  $.

\begin{proposition}
\cite[Sec. 2.1]{Fulton Introduction to Toric Varieties} $\mathbb{C}\left[
\check{\sigma}\cap M\right]  $ is integrally closed, i.e., $U\left(
\sigma\right)  $ is
\index{normal toric variety}%
normal.
\end{proposition}

The semigroup $\check{\sigma}\cap M$ is
\index{saturated semigroup|textbf}%
\textbf{saturated}, which means that if $a\cdot m\in\check{\sigma}\cap M$ for
$a\in\mathbb{Z}_{>0}$ and $m\in M$, then $m\in\check{\sigma}\cap M$. Indeed a
semigroup $S\subset M$ is saturated if and only if $\mathbb{C}\left[
S\right]  $ is integrally closed. The saturation
\[
\left\{  m\in M\mid a\cdot m\in S\text{ for some }a\in\mathbb{Z}_{>0}\right\}
\]
gives the normalization
\index{normalization of a toric variety}%
of $\operatorname*{Spec}\mathbb{C}\left[  S\right]  $.

The following proposition gives a characterization of the semigroups given by
the duals of strongly convex rational polyhedral cones.

\begin{proposition}
\cite[Sec. 1.1]{Oda Convex Bodies and Algebraic Geometry} Let $S\ $be an
additive subsemigroup of $M$. Then there is a unique strongly convex rational
polyhedral cone $\sigma$ in $N_{\mathbb{R}}$ with $S=\check{\sigma}\cap M$ if
and only if the following conditions are satisfied:

\begin{enumerate}
\item $S$ contains $0\in M$.

\item $S$ is finitely generated as an additive semigroup, i.e., there are
$m_{1},...,m_{r}\in S$ with%
\[
S=\left\{  a_{1}m_{1}+...+a_{r}m_{r}\mid m_{1},...,m_{r}\in\mathbb{Z}_{\geq
0}\right\}
\]

\item $S$ generates $M$ as a group, i.e., $S+\left(  -S\right)  =M$.

\item $S$ is saturated.
\end{enumerate}
\end{proposition}

\begin{lemma}
If $\sigma$ is a rational strongly convex polyhedral cone of dimension $n$ in
$N_{\mathbb{R}}$, then the dual cone $\check{\sigma}$ is also a rational
strongly convex polyhedral cone of dimension $n$, and there is a canonical
bijective inclusion reversing correspondence between the faces of $\sigma$ and
$\check{\sigma}$, given by%
\[
F\mapsto F^{\vee}=\left\{  m\in\check{\sigma}\mid\left\langle m,w\right\rangle
=0\ \forall w\in F\right\}
\]
if $F$ is a face of $\sigma$.
\end{lemma}

\begin{proposition}
\cite{Danilov The geometry of toric varieties} Any toric $U\left(
\sigma\right)  $ is
\index{Cohen-Macaulay}%
Cohen-Macaulay.
\end{proposition}

\begin{proposition}
\cite[2.1]{Fulton Introduction to Toric Varieties} $U\left(  \sigma\right)  $
is
\index{nonsingular toric variety}%
nonsingular if and only if $\sigma$ is generated by a subset of a basis of
$N$. Then%
\[
U\left(  \sigma\right)  =\mathbb{C}^{\dim\sigma}\times\left(  \mathbb{C}%
^{\ast}\right)  ^{n-\dim\sigma}%
\]

\end{proposition}

\begin{example}
For $\sigma=\operatorname*{hull}\left\{  \left(  0,1\right)  ,\left(
2,-1\right)  \right\}  $ we get the semigroup $S=\check{\sigma}\cap
M=\left\langle \left(  1,0\right)  ,\left(  1,1\right)  ,\left(  1,2\right)
\right\rangle \subset\mathbb{Z}^{2}$ given in Example \ref{examplecone1}.
\end{example}

Given two such cones $\sigma_{1}$ and $\sigma_{2}$ intersecting along a face
$\tau$ of both cones, the inclusions of $\tau\subset\sigma_{1},\sigma_{2}$
give inclusions of $U\left(  \tau\right)  \subset U\left(  \sigma_{1}\right)
,U\left(  \sigma_{2}\right)  $, hence we can glue the corresponding
\index{affine toric variety}%
affine toric varieties $U\left(  \sigma_{1}\right)  $ and $U\left(  \sigma
_{2}\right)  $ along $U\left(  \tau\right)  $.

A finite set $\Sigma$ of
\index{stongly convex}%
strongly convex rational polyhedral cones in $N_{\mathbb{R}}$ with the
property that every face of a cone in $\Sigma$ is again a cone in $\Sigma$,
\newsym[$\Sigma$]{fan}{} is called a
\index{fan|textbf}%
\textbf{fan}. Given a fan $\Sigma$ we can glue all $U\left(  \sigma\right)  $,
$\sigma\in\Sigma$, and
\newsym[$X\left(  \Sigma\right)  $]{toric variety associated to the fan $\Sigma$}{}get
a
\index{toric variety|textbf}%
\textbf{toric variety} $X\left(  \Sigma\right)  $.

Denote by
\[
\operatorname*{supp}\left(  \Sigma\right)  =\bigcup_{\sigma\in\Sigma}\sigma
\]
the
\index{support of a fan|textbf}%
\textbf{support }of a
\newsym[$\operatorname*{supp}\left(  \Sigma\right)  $]{support of the fan $\Sigma$}{}fan
$\Sigma$ and by $\Sigma\left(  m\right)  $ the set of $m$-dimensional cones of
$\Sigma$. The elements of $\Sigma\left(  1\right)  $ are
\newsym[$\Sigma\left(  1\right)  $]{rays of the fan $\Sigma$}{}called
\index{ray|textbf}%
\textbf{rays}, and for each $r\in\Sigma\left(  1\right)  $ let $\hat{r}$ be
the
\index{minimal lattice generator|textbf}%
minimal \newsym[$\hat{r}$]{minimal lattice generator of the ray $r$}{}lattice
generator of $r$, i.e., the unique generator of the semigroup $r\cap N$.

A cone is called
\index{simplicial|textbf}%
\textbf{simplicial}
\index{simplicial cone|textbf}%
if it is generated by linearly independent generators. A fan $\Sigma$ and the
toric
\index{simplicial toric variety|textbf}%
variety $X\left(  \Sigma\right)  $ are called
\index{simplicial fan|textbf}%
simplicial if all cones of $\Sigma$ are simplicial.

\begin{figure}
[h]
\begin{center}
\includegraphics[
trim=0.226009in 0.226201in 0.000000in 0.000000in,
height=1.4408in,
width=1.452in
]%
{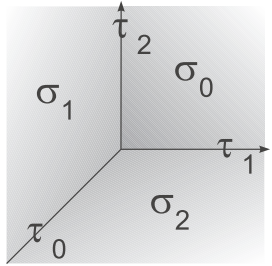}%
\caption{Fan representing $\mathbb{P}^{2}$ as toric variety}%
\label{FigP2}%
\end{center}
\end{figure}

\begin{figure}
[hh]
\begin{center}
\includegraphics[
height=1.343in,
width=5.8306in
]%
{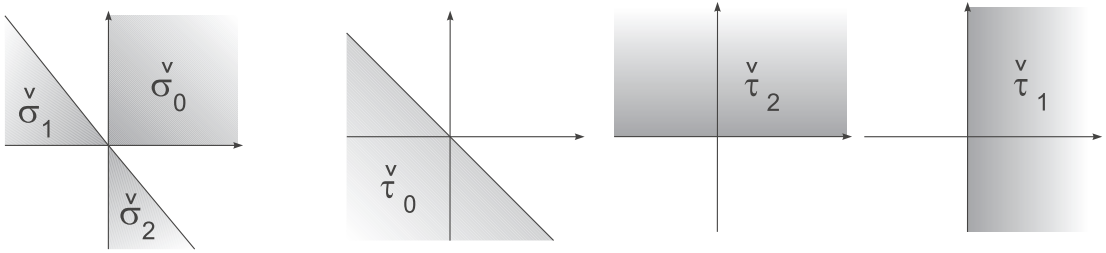}%
\caption{Duals of the cones of fan representing $\mathbb{P}^{2}$}%
\label{FigP2dualcones}%
\end{center}
\end{figure}

\begin{example}
\label{2exampleP2fan}The fan $\Sigma$, as depicted in Figure \ref{FigP2},
formed by the cones $\left\{  0\right\}  ,\tau_{0},\tau_{1},\tau_{2}%
,\sigma_{0},\sigma_{1},\sigma_{2}$ where%
\[%
\begin{tabular}
[c]{ll}%
$\tau_{0}=\operatorname*{hull}\left\{  \left(  -1,-1\right)  \right\}  $ &
$\sigma_{0}=\operatorname*{hull}\left\{  \left(  1,0\right)  ,\left(
0,1\right)  \right\}  $\\
$\tau_{1}=\operatorname*{hull}\left\{  \left(  1,0\right)  \right\}  $ &
$\sigma_{1}=\operatorname*{hull}\left\{  \left(  0,1\right)  ,\left(
-1,-1\right)  \right\}  $\\
$\tau_{2}=\operatorname*{hull}\left\{  \left(  0,1\right)  \right\}  $ &
$\sigma_{2}=\operatorname*{hull}\left\{  \left(  1,0\right)  ,\left(
-1,-1\right)  \right\}  $%
\end{tabular}
\]
gives $X\left(  \Sigma\right)  =\mathbb{P}^{2}$. The
\index{dual cone}%
dual cones of $\sigma_{i}$ and $\tau_{i}$ are shown in Figure
\ref{FigP2dualcones}$\ $and%
\[%
\begin{tabular}
[c]{lll}%
$U\left(  \sigma_{i}\right)  \cong\mathbb{A}^{2}$ & $U\left(  \tau_{i}\right)
\cong\mathbb{C\times C}^{\ast}$ & $U\left(  0\right)  \cong\left(
\mathbb{C}^{\ast}\right)  ^{2}$%
\end{tabular}
\]

\end{example}%

The
\index{torus action}%
torus actions on the
\index{affine toric variety}%
affine toric varieties $U\left(  \sigma\right)  $ give an action of the torus
on $X\left(  \Sigma\right)  $ extending the product in the
\index{torus}%
torus (see \cite[Sec. 1.4]{Fulton Introduction to Toric Varieties}).

\begin{proposition}
\cite[Sec. 3.1]{Fulton Introduction to Toric Varieties} The
\index{torus}%
torus
\index{torus action}%
acts on $X\left(  \Sigma\right)  $ and we get an inclusion reversing bijection
between the cones $\tau$ of $\Sigma$ and the
\index{torus orbit closure|textbf}%
closures $V\left(  \tau\right)  $
\newsym[$O\left(  \tau\right)  $]{torus orbit}{}of
\newsym[$V\left(  \tau\right)  $]{torus orbit closure}{}the
\index{torus orbit|textbf}%
torus orbits $O\left(  \tau\right)  $
\[%
\begin{tabular}
[c]{lll}%
$\Sigma$ & $\overset{1:1}{\longrightarrow}$ & $\left\{  \text{torus orbit
closures in }X\left(  \Sigma\right)  \right\}  $\\
$\tau$ & $\mapsto$ & $V\left(  \tau\right)  =\overline{O\left(  \tau\right)
}$%
\end{tabular}
\]

\end{proposition}

For any rational convex polyhedral cone $\sigma$ in $N_{\mathbb{R}}$,
$U\left(  \sigma\right)  $ contains
\newsym[$x_{\sigma}$]{distinguished point}{}a
\index{distinguished point|textbf}%
\textbf{distinguished point }$x_{\sigma}$ (see \cite[Sec. 2.1]{Fulton
Introduction to Toric Varieties}) given by%
\[%
\begin{tabular}
[c]{llll}%
$x_{\sigma}:$ & $\check{\sigma}\cap M$ & $\rightarrow$ & $\mathbb{C}$\\
& $m$ & $\mapsto$ & $\left\{
\begin{tabular}
[c]{ll}%
$1$ & if $m\in\sigma^{\perp}$\\
$0$ & otherwise
\end{tabular}
\right\}  $%
\end{tabular}
\]
with $\sigma^{\perp}=\left\{  m\in M_{\mathbb{R}}\mid\left\langle
m,w\right\rangle =0\text{ }\forall w\in\sigma\right\}  $, which is well
defined as $\sigma^{\perp}\cap\check{\sigma}$ is a face of $\check{\sigma}$.

If $\sigma$ spans $N_{\mathbb{R}}$, then $x_{\sigma}$ is the unique
\index{fixed point of torus action}%
fixed point of the
\index{torus action}%
torus action on $U\left(  \sigma\right)  $.

For the multiplicative group $\mathbb{C}^{\ast}$, we have $\operatorname*{Hom}%
\left(  \mathbb{C}^{\ast},\mathbb{C}^{\ast}\right)  =\mathbb{Z}$, hence
\index{torus}%
for $T=\operatorname*{Hom}\left(  M,\mathbb{C}^{\ast}\right)  $, there is a
one-to-one correspondence between lattice points $w\in N$ and $1$-parameter
subgroups $\lambda_{w}$ of $T$:%
\[%
\begin{tabular}
[c]{lll}%
$N\cong\operatorname*{Hom}\left(  \mathbb{Z},N\right)  $ & $\cong$ &
$\operatorname*{Hom}\left(  \mathbb{C}^{\ast},T\right)  $\\
$w$ & $\mapsto$ & $\left(
\begin{tabular}
[c]{llll}%
$\lambda_{w}:$ & $\mathbb{C}^{\ast}$ & $\rightarrow$ & $\operatorname*{Hom}%
\left(  M,\mathbb{C}^{\ast}\right)  $\\
& \multicolumn{1}{c}{$t$} & $\mapsto$ & $\left(
\begin{tabular}
[c]{llll}%
$\lambda_{w}\left(  t\right)  :$ & $M$ & $\rightarrow$ & $\mathbb{C}^{\ast}$\\
& $m$ & $\mapsto$ & $t^{\left\langle m,w\right\rangle }$%
\end{tabular}
\right)  $%
\end{tabular}
\right)  $%
\end{tabular}
\]

Note that by $M=\operatorname*{Hom}\left(  N,\mathbb{Z}\right)
=\operatorname*{Hom}\left(  T,\mathbb{C}^{\ast}\right)  $ there is also a
\newsym[$\lambda_{w}$]{one parameter subgroup}{}one-to-one correspondence
between the elements of $M$ and the elements of the
\index{character group}%
character group $\widehat{T}$ of the
\index{torus}%
torus.

\begin{proposition}
\cite[Sec. 2.3]{Fulton Introduction to Toric Varieties}%
\label{Prop Toric limit points} If $\tau$ is a cone of $\Sigma$ and
$w\in\operatorname*{int}\left(  \tau\right)  $ in the relative interior,
\index{limit}%
then%
\[
\lim_{t\rightarrow0}\lambda_{w}\left(  t\right)  =x_{\tau}%
\]
If $v$ is not in any cone of $\Sigma$, then the limit does not exist in
$X\left(  \Sigma\right)  $.
\end{proposition}

Note that this characterizes $\sigma\cap N$ as the set%
\[
\sigma\cap N=\left\{  w\in N\mid\lim_{t\rightarrow0}\lambda_{w}\left(
t\right)  \text{ exists in }U\left(  \sigma\right)  \right\}
\]
hence allows to recover the fan from the
\index{torus action}%
torus action.

In the above one-to-one correspondence between cones $\sigma$ of $\Sigma$ and
torus orbits, $O\left(  \sigma\right)  $ is the unique
\index{torus orbit}%
torus orbit containing $x_{\sigma}$. As
\[
V\left(  \sigma\right)  =\bigcup\nolimits_{\substack{\tau\in\Sigma
\\\tau\subset\sigma}}O\left(  \tau\right)
\]
$V\left(  \sigma\right)  $ contains precisely the
\index{distinguished point}%
distinguished points $x_{\tau}$ for $\tau\subset\sigma$.

\begin{figure}
[h]
\begin{center}
\includegraphics[
height=1.9363in,
width=1.9363in
]%
{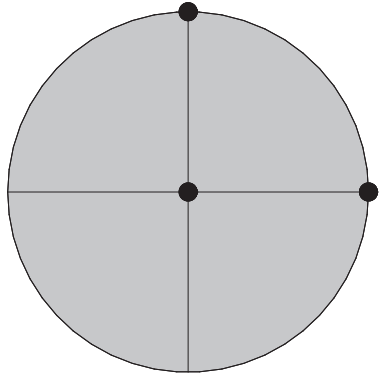}%
\caption{Torus orbits of $\mathbb{P}^{2}$}%
\label{FigP2Strata}%
\end{center}
\end{figure}

\begin{example}
In Example \ref{2exampleP2fan} the
\index{torus orbit}%
torus orbits and their
\index{torus orbit closure}%
closures are%
\[%
\begin{tabular}
[c]{l|l|l|l}%
$\sigma$ & $O\left(  \sigma\right)  $ & $V\left(  \sigma\right)  $ &
$x_{\sigma}$\\\hline
$\sigma_{k}$ & $\left\{  X_{i}=0,X_{j}=0,X_{k}\neq0\right\}  $ & $\left\{
X_{i}=0,X_{j}=0\right\}  $ &
\begin{tabular}
[c]{l}%
$x_{\sigma_{0}}=\left(  0:0:1\right)  $\\
$x_{\sigma_{1}}=\left(  1:0:0\right)  $\\
$x_{\sigma_{2}}=\left(  0:1:0\right)  $%
\end{tabular}
\\
$\tau_{k}$ & $\left\{  X_{i}=0,X_{j}\neq0,X_{k}\neq0\right\}  $ & $\left\{
X_{i}=0\right\}  $ &
\begin{tabular}
[c]{l}%
$x_{\tau_{0}}=\left(  0:1:1\right)  $\\
$x_{\tau_{1}}=\left(  1:0:1\right)  $\\
$x_{\tau_{2}}=\left(  1:1:0\right)  $%
\end{tabular}
\\
$0$ & $\left\{  X_{i}\neq0,X_{j}\neq0,X_{k}\neq0\right\}  $ & $\mathbb{P}%
^{2}=X\left(  \Sigma\right)  $ &
\begin{tabular}
[c]{l}%
$x_{0}=\left(  1:1:1\right)  $%
\end{tabular}
\end{tabular}
\]
Figure \ref{FigP2Strata} shows the real picture of the torus orbits,
identifying opposite points of the outer circle.
\end{example}

%

\subsubsection{Morphisms of toric
varieties\label{Sec morphisms of toric varieties}}

Suppose $\sigma^{\prime}\subset N^{\prime}$ is a
\index{morphisms of toric varieties}%
strongly convex rational polyhedral cone and
\[
\varphi:N^{\prime}\rightarrow N
\]
is a homomorphism of lattices such that $\varphi_{\mathbb{R}}:N_{\mathbb{R}%
}^{\prime}\rightarrow N_{\mathbb{R}}$ is mapping $\sigma^{\prime}$ into a
strongly convex rational polyhedral cone $\sigma\subset N_{\mathbb{R}}$. Hence
the dual $\varphi^{\ast}:M\rightarrow M^{\prime}$ maps $\check{\sigma}\cap M$
to $\check{\sigma}^{\prime}\cap M^{\prime}$ and gives a homomorphism%
\[
\mathbb{C}\left[  \check{\sigma}\cap M\right]  \rightarrow\mathbb{C}\left[
\check{\sigma}^{\prime}\cap M^{\prime}\right]
\]
hence a morphism%
\[
U\left(  \sigma^{\prime}\right)  \rightarrow U\left(  \sigma\right)
\]

\begin{proposition}
\cite[Sec. 1.4]{Fulton Introduction to Toric Varieties} Suppose $\Sigma$ is a
fan in $N$ and $\Sigma^{\prime}$ is a fan in $N^{\prime}$ and $\varphi
:N^{\prime}\rightarrow N$ is a homomorphism of lattices. If for each cone
$\sigma^{\prime}$ in $\Sigma^{\prime}$ there is some cone $\sigma$ in $\Sigma$
such that $\varphi\left(  \sigma^{\prime}\right)  \subset\sigma$, then there
is a morphism $U\left(  \sigma^{\prime}\right)  \rightarrow U\left(
\sigma\right)  \subset X\left(  \Sigma\right)  $, and the
\index{morphisms of toric varieties}%
morphism $U\left(  \sigma^{\prime}\right)  \rightarrow X\left(  \Sigma\right)
$ is independent of the choice of $\sigma$. These morphisms patch together to
a
\index{toric blowup|textbf}%
morphism%
\[
\varphi_{\ast}:X\left(  \Sigma^{\prime}\right)  \rightarrow X\left(
\Sigma\right)
\]

\end{proposition}

If $X\left(  \Sigma\right)  $ is compact, then $\operatorname*{supp}\left(
\Sigma\right)  =N_{\mathbb{R}}$. Otherwise, there would be a $w\in\left(
N_{\mathbb{R}}-\operatorname*{supp}\left(  \Sigma\right)  \right)  \cap N$ and
$\lim_{t\rightarrow0}\lambda_{w}\left(  t\right)  $ would not exist in
$X\left(  \Sigma\right)  $. The converse is given by:

\begin{proposition}
Let $\Sigma$ be a fan in $N$ and $\Sigma^{\prime}$ a fan in $N^{\prime}$ and
$\varphi:N^{\prime}\rightarrow N$ a homomorphism of lattices inducing a
morphism $\varphi_{\ast}:X\left(  \Sigma^{\prime}\right)  \rightarrow X\left(
\Sigma\right)  $. The morphism $\varphi_{\ast}$ is
\index{proper morphism}%
proper if and only if $\varphi^{-1}\left(  \operatorname*{supp}\left(
\Sigma\right)  \right)  =\operatorname*{supp}\left(  \Sigma^{\prime}\right)  $.
\end{proposition}

\begin{corollary}
The toric variety $X\left(  \Sigma\right)  $ is
\index{complete toric variety}%
complete if and only if $\Sigma$
\index{complete fan}%
is
\index{support of a fan}%
complete, i.e., $\operatorname*{supp}\left(  \Sigma\right)  =N_{\mathbb{R}}$.
\end{corollary}

\begin{example}
Suppose $v_{1},...,v_{n}$ are a basis of $N$ generating the cone
\[
\sigma=\operatorname*{hull}\left\{  v_{1},...,v_{n}\right\}
\]
and $\Sigma$ is the fan generated by the cone $\sigma$ (i.e., the fan
consisting of all faces of $\sigma$), so $x_{\sigma}=\left(  0,...,0\right)
\in\mathbb{C}^{n}=U\left(  \sigma\right)  =X\left(  \Sigma\right)  $. Write
$x_{i}=x^{e_{i}^{\ast}}$, $i=1,...,n$. Set $v_{0}=v_{1}+...+v_{n}$ and
consider the subdivision of $\sigma$ with respect to $v_{0}$, i.e., the fan
$\Sigma^{\prime}$ generated by the cones
\[
\sigma_{i}=\operatorname*{hull}\left\{  v_{0},v_{1},...,v_{i-1},v_{i+1}%
,...,v_{n}\right\}
\]
for $i=1,...,n$. Then $X\left(  \Sigma^{\prime}\right)  $ is the
\index{blowup}%
blowup of $X\left(  \Sigma\right)  $ at $x_{\sigma}$:

To describe $X\left(  \Sigma^{\prime}\right)  $ note that%
\[
\check{\sigma}_{i}=\operatorname*{hull}\left\{  v_{i}^{\ast},v_{1}^{\ast
}-v_{i}^{\ast},...,v_{i-1}^{\ast}-v_{i}^{\ast},v_{i+1}^{\ast}-v_{i}^{\ast
},...,v_{n}^{\ast}-v_{i}^{\ast}\right\}
\]
hence%
\[
\mathbb{C}\left[  \check{\sigma}_{i}\cap M\right]  =\mathbb{C}\left[
x_{i},x_{1}x_{i}^{-1},...,x_{i-1}x_{i}^{-1},x_{i+1}x_{i}^{-1},...,x_{n}%
x_{i}^{-1}\right]
\]
so $U\left(  \sigma_{i}\right)  =\mathbb{C}^{n}$.

The
\index{blowup}%
blowup
\index{toric blowup}%
of $U\left(  \sigma\right)  $ at $x_{\sigma}$ is $\left\{  x_{i}y_{j}%
-x_{j}y_{i}\mid i,j=1,...,n\right\}  \subset\mathbb{C}^{n}\times
\mathbb{P}^{n-1}$, where $y_{1},...,y_{n}$ are homogeneous coordinates on
$\mathbb{P}^{n-1}$, and it is covered by the open sets $U_{i}=\left\{
y_{i}\neq0\right\}  =\mathbb{C}^{n}$ for $i=1,...,n$, which by $x_{j}%
=x_{i}\frac{y_{j}}{y_{i}}$ and $\frac{y_{j}}{y_{i}}=\frac{x_{j}}{x_{i}}$ have
coordinates $x_{i},x_{1}x_{i}^{-1},...,x_{i-1}x_{i}^{-1},x_{i+1}x_{i}%
^{-1},...,x_{n}x_{i}^{-1}$.
\end{example}

Any cone of a given fan $\Sigma$ can be subdivided such that it becomes
\index{simplicial}%
simplicial. Given a simplicial cone $\sigma$ of dimension $d$ with
\index{minimal lattice generator}%
minimal lattice generators $v_{1},...,v_{d}$ of the rays of $\sigma$ and the
lattice $N_{\sigma}=\left\langle \sigma\cap N\right\rangle $
\newsym[$\operatorname*{mult}\left(  \sigma\right)  $]{multiplicity of a cone}{}generated
by $\sigma$, the \textbf{multiplicity}
\index{multiplicity of a cone|textbf}%
of $\sigma$ is defined as the index of $\mathbb{Z}v_{1}+...+\mathbb{Z}v_{d}$
in $N_{\sigma}$%
\[
\operatorname*{mult}\left(  \sigma\right)  =\left[  N_{\sigma}:\mathbb{Z}%
v_{1}+...+\mathbb{Z}v_{d}\right]
\]
Then $U\left(  \sigma\right)  $ is
\index{nonsingular toric variety}%
nonsingular if and only if $\operatorname*{mult}\left(  \sigma\right)  =1$.

\begin{figure}
[h]
\begin{center}
\includegraphics[
height=1.2903in,
width=2.079in
]%
{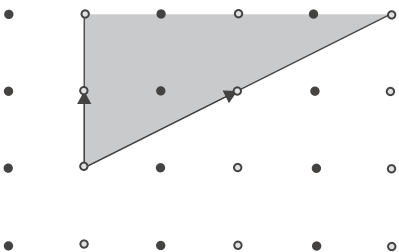}%
\caption{Cone of multiplicity $2$}%
\label{FigMultCone}%
\end{center}
\end{figure}

\begin{example}
If $N=\mathbb{Z}e_{1}+\mathbb{Z}e_{2}$ and $\sigma=\left\langle v_{1}%
,v_{2}\right\rangle $ with $v_{2}=e_{2}$ and $v_{1}=2e_{1}+e_{2}$, then
$N_{\sigma}=N$ and $\mathbb{Z}v_{1}+\mathbb{Z}v_{2}=\mathbb{Z}\left(
2e_{1}\right)  +\mathbb{Z}e_{2}$. Figure \ref{FigMultCone} shows the cone
$\sigma$ and the groups $N_{\sigma}$ and $\mathbb{Z}v_{1}+\mathbb{Z}v_{2}$.
\end{example}

%

Any simplicial
\index{simplicial}%
cone $\sigma$ can be subdivided until it has multiplicity $1$: If
$\operatorname*{mult}\left(  \sigma\right)  >1$ there is a $v\in\sigma\cap N$
such that $v=\sum_{j=1}^{d}a_{j}v_{j}$ with $0\leq a_{j}<1$. Subdividing
$\sigma$ with respect to $\operatorname*{hull}\left\{  v\right\}  $, we obtain
the cones $\sigma_{i}=\operatorname*{hull}\left\{  v,v_{1},...,v_{i-1}%
,v_{i+1},...,v_{d}\right\}  $ for all $i$ with $a_{i}\neq0$, and if $v$ is
the
\index{minimal lattice generator}%
minimal lattice generator of $\operatorname*{hull}\left\{  v\right\}  $, then
\[
\operatorname*{mult}\left(  \sigma_{i}\right)  =a_{i}\operatorname*{mult}%
\left(  \sigma\right)
\]
Hence:

\begin{proposition}
\cite[Sec. 2.6]{Fulton Introduction to Toric Varieties} Given a fan $\Sigma$
in $N$ and a fan $\Sigma^{\prime}$ in $N$ refining $\Sigma$, the identity
$id:N\rightarrow N$ on the lattice induces a proper birational morphism
$id_{\ast}:X\left(  \Sigma^{\prime}\right)  \rightarrow X\left(
\Sigma\right)  $.

There is a refinement $\Sigma^{\prime}$ of $\Sigma$ inducing a
\index{resolution of singularities}%
resolution of singularities $X\left(  \Sigma^{\prime}\right)  \rightarrow
X\left(  \Sigma\right)  $.
\end{proposition}

\subsubsection{Divisors on toric varieties\label{Divisors on toric varieties}}

Let $X\left(  \Sigma\right)  $ be a toric variety of dimension $n$. The
\index{T-invariant}%
$T$-invariant prime
\index{Weil divisor}%
Weil divisors on $X\left(  \Sigma\right)  $ are the components of dimension
$n-1$ of the complement $X\left(  \Sigma\right)  -O\left(  \left\{  0\right\}
\right)  $ of the
\index{torus orbit}%
torus orbit $O\left(  \left\{  0\right\}  \right)  $, i.e., they are the
closures of the torus orbits of dimension $n-1$. We denote by $D_{r}$ the
\index{T-invariant}%
$T$-invariant prime Weil divisor corresponding to the ray $r\in\Sigma\left(
1\right)  $. Denote by $\operatorname*{WDiv}_{T}\left(  X\left(
\Sigma\right)  \right)  $ the
\newsym[$\operatorname*{WDiv}_{T}\left(  X\left(  \Sigma\right)  \right)  $]{group of $T$-invariant Weil divisors of $X\left(  \Sigma\right)$}{}group
of
\index{T-invariant}%
$T$-invariant Weil divisors on $X\left(  \Sigma\right)  $, which is isomorphic
to $\mathbb{Z}^{\Sigma\left(  1\right)  }$ by%
\[%
\begin{tabular}
[c]{lll}%
$\mathbb{Z}^{\Sigma\left(  1\right)  }$ & $\rightarrow$ &
$\operatorname*{WDiv}_{T}\left(  X\left(  \Sigma\right)  \right)  $\\
$\left(  a_{r}\right)  _{r}$ & $\mapsto$ & $\sum_{r}a_{r}D_{r}$%
\end{tabular}
\]%
\index{T-invariant}%
$T$-invariant divisors
\index{T-divisor|textbf}%
are
\index{T-Weil divisor|textbf}%
also
\index{T-Cartier divisor|textbf}%
called $T$\textbf{-divisors}.

By a
\index{T-invariant}%
$T$-invariant
\index{Cartier divisor}%
Cartier divisor $D$ on $X\left(  \Sigma\right)  $ a collection of rational
functions $\varphi_{D,\sigma}$, $\sigma\in\Sigma$ is given such that
$\varphi_{D,\sigma}$ defines $D$ on $U\left(  \sigma\right)  $, $\varphi
_{D,\sigma}$ is invariant under the
\index{torus action}%
torus action up to multiplication by a non zero constant, $\varphi_{D,\sigma}$
is unique up to multiplication with an invertible function on $U\left(
\sigma\right)  $, and $\frac{\varphi_{D,\sigma_{1}}}{\varphi_{D,\sigma2}}$ is
invertible on $U\left(  \sigma_{1}\right)  \cap U\left(  \sigma_{2}\right)  $.
As $\varphi_{D,\sigma}$ is an eigenvector of the action of the
\index{torus}%
torus, we can write%
\[
\varphi_{D,\sigma}=x^{-m\left(  D,\sigma\right)  }%
\]
with $m\left(  D,\sigma\right)  \in M$.

As $\varphi_{D,\sigma}$ is unique up to multiplication with an invertible
function on $U\left(  \sigma\right)  $, the lattice point $m\left(
D,\sigma\right)  \in M$ is unique modulo the
\newsym[$M_{\sigma}$]{orthogonal sublattice}{}sublattice%
\[
M_{\sigma}=\left\langle \sigma\cap N\right\rangle ^{\perp}=\left\{  m\in
M\mid\left\langle m,w\right\rangle =0\text{ }\forall w\in\sigma\right\}
\]
of $M$ orthogonal to the sublattice $\left\langle \sigma\cap N\right\rangle $
of $N$ generated by $\sigma$. Hence giving $D\mid_{U\left(  \sigma\right)  }$
is equivalent to specifying the function $\left\langle m\left(  D,\sigma
\right)  ,-\right\rangle $.

Invertibility of
\[
\frac{\varphi_{D,\sigma_{1}}}{\varphi_{D,\sigma2}}=x^{-m\left(  D,\sigma
_{1}\right)  +m\left(  D,\sigma_{2}\right)  }%
\]
on $U\left(  \sigma_{1}\right)  \cap U\left(  \sigma_{2}\right)  =U\left(
\tau\right)  $ with $\tau=\sigma_{1}\cap\sigma_{2}$ is equivalent to the
condition that $m\left(  D,\sigma_{1}\right)  -m\left(  D,\sigma_{2}\right)
\in M_{\tau}$ or, equivalently, that the functions $\left\langle m\left(
D,\sigma_{1}\right)  ,-\right\rangle $ and $\left\langle m\left(  D,\sigma
_{2}\right)  ,-\right\rangle $ agree on $\tau$, i.e.,%
\[
\left\langle m\left(  D,\sigma_{1}\right)  ,-\right\rangle \mid_{\tau
}=\left\langle m\left(  D,\sigma_{2}\right)  ,-\right\rangle \mid_{\tau}%
\]

Hence associated to $D$ there is a well defined
\index{piecewise linear}%
piecewise linear continuous function
\begin{gather*}
\Phi_{D}:\operatorname*{supp}\left(  \Sigma\right)  \rightarrow\mathbb{R}\\
\Phi_{D}\left(  w\right)  =\left\langle m\left(  D,\sigma\right)
,w\right\rangle \text{ for }w\in\sigma
\end{gather*}
on the
\index{support of a fan}%
support\textbf{ }of the fan $\Sigma$. The function $\Phi_{D}$ is called
\newsym[$\Phi_{D}$]{support function}{}the
\index{support function|textbf}%
\textbf{support function} of $D$. A
\index{piecewise linear}%
piecewise linear continuous function on $\operatorname*{supp}\left(
\Sigma\right)  $, which is given on $\sigma$ by $\left\langle m\left(
\sigma\right)  ,-\right\rangle $ with $m\left(  \sigma\right)  \in M$ is
called
\index{integral function|textbf}%
\textbf{integral}.

On the other hand, any
\index{piecewise linear}%
piecewise linear continuous integral function $\Phi:\operatorname*{supp}%
\left(  \Sigma\right)  \rightarrow\mathbb{R}$ is the support function of a
unique Cartier divisor which, written as a
\index{T-Weil divisor}%
Weil divisor, is given by%
\[
D=\sum_{r\in\Sigma\left(  1\right)  }-\Phi\left(  \hat{r}\right)  D_{r}%
\]

A $T$-Weil
\index{T-Weil divisor}%
divisor $D=\sum_{r\in\Sigma\left(  1\right)  }a_{r}D_{r}$ is
\index{T-Cartier divisor}%
Cartier if and only if for all $\sigma\in\Sigma$ there is an $m\left(
D,\sigma\right)  \in M$ such that%
\[
-a_{r}=\left\langle m\left(  D,\sigma\right)  ,\hat{r}\right\rangle \text{ for
all }r\in\Sigma\left(  1\right)  \text{ with }\hat{r}\in\sigma
\]

\begin{proposition}
\cite[Sec. 4.2]{Voisin Mirror Symmetry} By associating to
\index{T-Cartier divisor}%
a $T$-Cartier divisor $D$ the function $\Phi_{D}$, we get a one-to-one
correspondence between $T$-Cartier divisors on $X\left(  \Sigma\right)  $ and
\index{piecewise linear}%
piecewise linear continuous integral functions on $\operatorname*{supp}\left(
\Sigma\right)  $.
\end{proposition}

For any $m\in M$ the
\index{Laurent monomial}%
Laurent monomial $x^{m}$ is a holomorphic function on the
\index{torus}%
torus $T$, hence a rational function on $X\left(  \Sigma\right)  $ defining a
\index{T-invariant}%
$T$-invariant
\newsym[$\operatorname*{div}\left(  x^{m}\right)  $]{principal Cartier divisor}{}principal
Cartier divisor%
\[
\operatorname*{div}\left(  x^{m}\right)  =\sum_{r\in\Sigma\left(  1\right)
}\left\langle m,\hat{r}\right\rangle D_{r}%
\]
The principal
\index{T-Cartier divisor}%
Cartier divisor $\operatorname*{div}\left(  x^{m}\right)  $ corresponds to the
support function $\Phi_{\operatorname*{div}\left(  x^{m}\right)
}=-\left\langle m,-\right\rangle $ defined globally by an element of $M$.

\begin{proposition}
\cite[Sec. 3.1]{CK Mirror Symmetry and Algebraic Geometry}, \cite[Sec.
3.4]{Fulton Introduction to Toric Varieties} Classes in the
\newsym[$\operatorname*{Pic}\left(  X\left(  \Sigma\right)  \right)  $]{Picard group of $X\left(  \Sigma\right)$}{}Picard
group
\index{Picard group|textbf}%
$\operatorname*{Pic}\left(  X\left(  \Sigma\right)  \right)  $ of line bundles
on $X\left(  \Sigma\right)  $ modulo isomorphism and the
\newsym[$A_{n-1}\left(  X\left(  \Sigma\right)  \right)  $]{Chow group of divisors of $X\left(  \Sigma\right)$}{}Chow
group
\index{Chow group|textbf}%
$A_{n-1}\left(  X\left(  \Sigma\right)  \right)  $ of Weil divisors on
$X\left(  \Sigma\right)  $ modulo linear equivalence can be represented by
\index{T-invariant}%
$T$-invariant
\index{T-Cartier divisor}%
Cartier respectively
\index{T-Weil divisor}%
Weil divisors via
\index{Weil divisor}%
the
\index{Cartier divisor}%
exact
\index{ray}%
sequences%
\[%
\begin{tabular}
[c]{lllllll}
& $m$ & $\mapsto$ & $\operatorname*{div}x^{m}$ &  &  & \\
$0\rightarrow$ & $M$ & $\rightarrow$ & $\operatorname*{Div}_{T}\left(
X\left(  \Sigma\right)  \right)  $ & $\rightarrow$ & $\operatorname*{Pic}%
\left(  X\left(  \Sigma\right)  \right)  $ & $\rightarrow0$\\
& \multicolumn{1}{c}{$\parallel$} &  & $\cap$ &  & $\cap$ & \\
$0\rightarrow$ & $M$ & $\overset{A}{\rightarrow}$ & $\mathbb{Z}^{\Sigma\left(
1\right)  }$ & $\rightarrow$ & $A_{n-1}\left(  X\left(  \Sigma\right)
\right)  $ & $\rightarrow0$\\
& $m$ & $\mapsto$ & $\left(  \left\langle m,\hat{r}\right\rangle \right)
_{r\in\Sigma\left(  1\right)  }$ &  &  & \\
&  &  & $\left(  a_{r}\right)  _{r}$ & $\mapsto$ & $\sum_{r}a_{r}D_{r}$ &
\end{tabular}
\]
$\operatorname*{Pic}\left(  X\left(  \Sigma\right)  \right)  $ is torsion free.
\end{proposition}

\begin{example}
For the fan $\Sigma$ of $\mathbb{P}^{2}$ given in Example \ref{2exampleP2fan}
we
\index{hypersurface}%
get
\[
A_{1}\left(  X\left(  \Sigma\right)  \right)  =\operatorname*{coker}\left(
\begin{array}
[c]{cc}%
1 & 0\\
0 & 1\\
-1 & -1
\end{array}
\right)  \cong\mathbb{Z}%
\]

\end{example}

\begin{example}
\label{2exquinticchow}Considering
\index{quintic threefold}%
the fan $\Sigma$ over the
\index{hypersurface}%
faces of the degree $5$
\index{Veronese}%
Veronese polytope of $\mathbb{P}^{4}$
\[
\operatorname*{convexhull}\left(  \left(  4,-1,-1,-1\right)  ,...,\left(
-1,-1,-1,4\right)  ,\left(  -1,-1,-1,-1\right)  \right)
\]
we get%
\[
A_{3}\left(  X\left(  \Sigma\right)  \right)  =\operatorname*{coker}\left(
\begin{array}
[c]{cccc}%
4 & -1 & -1 & -1\\
-1 & 4 & -1 & -1\\
-1 & -1 & 4 & -1\\
-1 & -1 & -1 & 4\\
-1 & -1 & -1 & -1
\end{array}
\right)  \cong\mathbb{Z\times}H
\]
with%
\[
H=\frac{\left\{  \left(  a_{0},...,a_{4}\right)  \in\mathbb{Z}_{5}^{5}\mid
\sum_{i=0}^{4}a_{i}=0\operatorname{mod}5\right\}  }{\mathbb{Z}_{5}\left(
1,1,1,1,1\right)  }\cong\mathbb{Z}_{5}^{3}%
\]

\end{example}

In the following suppose that all
\index{maximal cones}%
maximal dimensional cones of the fan $\Sigma$ have dimension $n$, and denote
by $\sigma_{1},...,\sigma_{s}$ the maximal dimensional cones of $\Sigma$. As
shown above a $T$ -Cartier divisor is given by a collection $m_{i}\in
M/M_{\sigma_{i}}=M$ for all $i$ such that $m_{i}=m\left(  D,\sigma_{i}\right)
$ maps to $m\left(  D,\sigma_{i}\cap\sigma_{j}\right)  $ under the canonical
map $M/M_{\sigma_{i}}\rightarrow M/M_{\sigma_{i}\cap\sigma_{j}}$, hence:

\begin{lemma}
The group of
\index{T-Cartier divisor}%
$T$-Cartier divisors $\operatorname*{Div}_{T}\left(  X\left(  \Sigma\right)
\right)  $ on $X\left(  \Sigma\right)  $ is given by the kernel of the map%
\[%
\begin{tabular}
[c]{ccc}%
$\bigoplus_{i=1}^{s}M/M_{\sigma_{i}}$ & $\rightarrow$ & $\bigoplus
_{i<j}M/M_{\sigma_{i}\cap\sigma_{j}}$\\
$\left(  m_{i}\right)  _{i}$ & $\mapsto$ & $\left(  m_{i}-m_{j}\right)
_{i<j}$%
\end{tabular}
\]

\end{lemma}

\begin{lemma}
\cite[Sec. 3.2]{Fulton Introduction to Toric Varieties} $H^{2}\left(  X\left(
\Sigma\right)  ,\mathbb{Z}\right)  $ is given by the kernel of the map%
\[
\bigoplus_{i<j}M_{\sigma_{i}\cap\sigma_{j}}\rightarrow\bigoplus_{i<j<l}%
M_{\sigma_{i}\cap\sigma_{j}\cap\sigma_{l}}%
\]

\end{lemma}

\begin{corollary}
The map%
\[%
\begin{tabular}
[c]{ccc}%
$\ker\left(  \bigoplus_{i=1}^{s}M/M_{\sigma_{i}}\rightarrow\bigoplus
_{i<j}M/M_{\sigma_{i}\cap\sigma_{j}}\right)  $ & $\rightarrow$ & $\ker\left(
\bigoplus_{i<j}M_{\sigma_{i}\cap\sigma_{j}}\rightarrow\bigoplus_{i<j<l}%
M_{\sigma_{i}\cap\sigma_{j}\cap\sigma_{l}}\right)  $\\
$\left(  m_{i}\right)  _{i}$ & $\mapsto$ & $\left(  m_{i}-m_{j}\right)
_{i<j}$%
\end{tabular}
\]
induces an
\index{Picard group}%
isomorphism%
\[%
\begin{tabular}
[c]{ccc}%
$\operatorname*{Pic}\left(  X\left(  \Sigma\right)  \right)  $ & $\cong$ &
$H^{2}\left(  X\left(  \Sigma\right)  ,\mathbb{Z}\right)  $%
\end{tabular}
\]

\end{corollary}

\begin{proposition}
\cite[Sec. 3.4]{Fulton Introduction to Toric Varieties} The following
conditions are equivalent:

\begin{enumerate}
\item $X\left(  \Sigma\right)  $ is simplicial.

\item All Weil divisors on $X\left(  \Sigma\right)  $ are $\mathbb{Q}$-Cartier.

\item $\operatorname*{Pic}\left(  X\left(  \Sigma\right)  \right)
\otimes\mathbb{Q\rightarrow}A_{n-1}\left(  X\left(  \Sigma\right)  \right)
\otimes\mathbb{Q}$ is an isomorphism.

\item $\operatorname*{rank}\left(  \operatorname*{Pic}\left(  X\left(
\Sigma\right)  \right)  \right)  =\left\vert \Sigma\left(  1\right)
\right\vert -n$.
\end{enumerate}
\end{proposition}

\begin{proposition}
\cite[Sec. 3.2]{CK Mirror Symmetry and Algebraic Geometry}, \cite[Sec.
3.4]{Fulton Introduction to Toric Varieties} Suppose $X\left(  \Sigma\right)
$ is complete. For any divisor $D=\sum_{r\in\Sigma\left(  1\right)  }%
a_{r}D_{r}$ the
\index{global sections}%
global sections of the
\index{reflexive sheaf}%
reflexive sheaf $\mathcal{O}_{X\left(  \Sigma\right)  }\left(  D\right)  $
correspond to the \newsym[$\Delta_{D}$]{polytope of global sections}{}lattice
points of the polytope%
\[
\Delta_{D}=\left\{  m\in M_{\mathbb{R}}\mid\left\langle m,\hat{r}\right\rangle
\geq-a_{r}\forall r\in\Sigma\left(  1\right)  \right\}
\]
i.e.,%
\[
H^{0}\left(  X\left(  \Sigma\right)  ,\mathcal{O}_{X\left(  \Sigma\right)
}\left(  D\right)  \right)  \cong\bigoplus_{m\in\Delta_{D}\cap M}%
\mathbb{C}x^{m}%
\]

\end{proposition}

\begin{remark}
If $D$ is
\index{Cartier divisor}%
Cartier, then%
\[
\Delta_{D}=\left\{  m\in M_{\mathbb{R}}\mid\left\langle m,-\right\rangle
\geq\Phi_{D}\text{ on }N_{\mathbb{R}}\right\}
\]

\end{remark}

\begin{lemma}
\cite[Sec. 3.2]{CK Mirror Symmetry and Algebraic Geometry}, \cite[Sec.
3.4]{Fulton Introduction to Toric Varieties} A
\index{Cartier divisor}%
Cartier divisor $D=\sum_{r\in\Sigma\left(  1\right)  }a_{r}D_{r}$ on a
\index{complete toric variety}%
complete toric variety $X\left(  \Sigma\right)  $ is
\index{generated by global sections}%
generated by global sections if and only if for all $\sigma\in\Sigma$%
\[
\left\langle m\left(  D,\sigma\right)  ,\hat{r}\right\rangle \geq-a_{r}\text{
}\forall r\in\Sigma\left(  1\right)  \text{ with }\hat{r}\notin\sigma
\]
Note that by definition%
\[
\left\langle m\left(  D,\sigma\right)  ,\hat{r}\right\rangle =-a_{r}\text{
}\forall r\in\Sigma\left(  1\right)  \text{ with }\hat{r}\in\sigma
\]

Hence by $\Phi_{D}\left(  \hat{r}\right)  =a_{r}$, it follows that $D$ is
\index{generated by global sections}%
generated by its global sections if and only if the graph of $\Phi_{D}$ lies
below the graphs of the functions $\left\langle m\left(  D,\sigma\right)
,-\right\rangle $ for all $\sigma\in\Sigma$, i.e., $\Phi_{D}$ is
\index{upper convex|textbf}%
\textbf{upper convex}.

Reformulating via the polytope of sections: $D$ is
\index{generated by global sections}%
generated by global sections if and only if%
\[
\Delta_{D}=\operatorname*{convexhull}\left\{  m\left(  D,\sigma\right)
\mid\sigma\in\Sigma\left(  n\right)  \right\}
\]

So in particular $\Delta_{D}$ is a
\index{lattice polytope}%
lattice polytope and%
\[
\Phi_{D}\left(  w\right)  =\min_{\sigma\in\Sigma\left(  n\right)
}\left\langle m\left(  D,\sigma\right)  ,w\right\rangle =\min_{m\in\Delta
_{D}\cap M}\left\langle m,w\right\rangle
\]
hence $\Phi_{D}$ or equivalently
\index{T-Cartier divisor}%
the $T$-Cartier divisor $D$ can be reconstructed from $\Delta_{D}$.
\end{lemma}

\begin{lemma}
\cite[Sec. 3.2]{CK Mirror Symmetry and Algebraic Geometry}, \cite[Sec.
3.4]{Fulton Introduction to Toric Varieties} A Cartier divisor $D=\sum
_{r\in\Sigma\left(  1\right)  }a_{r}D_{r}$ on a
\index{complete toric variety}%
complete toric variety $X\left(  \Sigma\right)  $ is
\index{ample}%
ample if and only if for all $\sigma\in\Sigma$%
\[
\left\langle m\left(  D,\sigma\right)  ,\hat{r}\right\rangle >-a_{r}\text{
}\forall r\in\Sigma\left(  1\right)  \text{ with }\hat{r}\notin\sigma
\]
i.e., for all $\sigma\in\Sigma$ the graph of $\Phi_{D}$ on the complement of
$\sigma$ lies strictly below the graph of $\left\langle m\left(
D,\sigma\right)  ,-\right\rangle $, i.e., $\Phi_{D}$ is
\index{strictly upper convex}%
\textbf{strictly upper convex}.

Reformulating via the polytope of sections: $D$ is
\index{ample}%
ample if and only if $\Delta_{D}$ is a polytope of dimension $n$ with%
\[
\operatorname*{vertices}\left(  \Delta_{D}\right)  =\left\{  m\left(
D,\sigma\right)  \mid\sigma\in\Sigma\left(  n\right)  \right\}
\]
and all $m\left(  D,\sigma\right)  $, $\sigma\in\Sigma\left(  n\right)  $ are
pairwise different.
\end{lemma}

\begin{lemma}
\cite[Sec. 3.4]{Fulton Introduction to Toric Varieties} Any
\index{ample}%
ample
\index{Cartier divisor}%
Cartier divisor on a
\index{complete toric variety}%
complete toric variety $X\left(  \Sigma\right)  $ is
\index{generated by global sections}%
generated by its global sections.
\end{lemma}

\begin{lemma}
\cite[Sec. 3.2]{CK Mirror Symmetry and Algebraic Geometry}, \cite[Sec.
3.4]{Fulton Introduction to Toric Varieties} A
\index{T-Cartier divisor}%
Cartier divisor $D=\sum_{r\in\Sigma\left(  1\right)  }a_{r}D_{r}$ on a
\index{complete toric variety}%
complete toric variety $X\left(  \Sigma\right)  $ is
\index{very ample}%
very ample if and only if $\Phi_{D}$ is
\index{strictly upper convex}%
strictly upper convex and for all $\sigma\in\Sigma\left(  n\right)  $%
\[
\check{\sigma}\cap M=\left\langle m-m\left(  D,\sigma\right)  \mid m\in
\Delta_{D}\cap M\right\rangle
\]

\end{lemma}

\begin{lemma}
\cite[Sec. 3.4]{Fulton Introduction to Toric Varieties} If $X\left(
\Sigma\right)  $ is complete and non-singular, then a $T$-divisor is
\index{ample}%
ample if and only if it is
\index{very ample}%
very ample.
\end{lemma}

\subsubsection{Dualizing sheaf of a toric
variety\label{Sec dualizing sheaf of a toric variety}}

Suppose $X\left(  \Sigma\right)  $ is a
\index{dualizing sheaf}%
nonsingular
\index{anticanonical class}%
toric variety, $e_{1},...,e_{n}$ form a basis of $N$ and $x_{i}=x^{e_{i}%
^{\ast}}$, $i=1,...,n$ are the corresponding coordinates, then the divisor of
the rational section
\[
\omega=\frac{dx_{1}}{x_{1}}\wedge...\wedge\frac{dx_{n}}{x_{n}}%
\]
of $\Omega_{X\left(  \Sigma\right)  }^{n}$ is $-\sum_{v\in\Sigma\left(
1\right)  }D_{v}$, hence:

\begin{proposition}
\cite[Sec. 4.3]{Fulton Introduction to Toric Varieties} If $X\left(
\Sigma\right)  $ is a nonsingular toric variety, then%
\[
\Omega_{X\left(  \Sigma\right)  }^{n}\cong\mathcal{O}_{X\left(  \Sigma\right)
}\left(  -\sum_{v\in\Sigma\left(  1\right)  }D_{v}\right)
\]

\end{proposition}

Suppose $X\left(  \Sigma\right)  $ is any toric variety of dimension $n$, then
$\sum_{v\in\Sigma\left(  1\right)  }D_{v}$ is not Cartier or $\mathbb{Q}%
$-Cartier in general, but still the coherent sheaf
\[
\hat{\Omega}_{X\left(  \Sigma\right)  }^{n}=\mathcal{O}_{X\left(
\Sigma\right)  }\left(  -\sum_{v\in\Sigma\left(  1\right)  }D_{v}\right)
\]
gives the
\newsym[$\hat{\Omega}_{X\left(  \Sigma\right)  }^{n}$]{dualizing sheaf}{}dualizing
sheaf.

\begin{proposition}
\cite[Sec. 4.4]{Fulton Introduction to Toric Varieties} Suppose $X\left(
\Sigma\right)  $ is a toric variety given by the fan $\Sigma$.

If $\Sigma^{\prime}$ is a refinement of $\Sigma$ inducing a
\index{resolution of singularities}%
resolution of singularities%
\[
f:X\left(  \Sigma^{\prime}\right)  \rightarrow X\left(  \Sigma\right)
\]
then%
\[
f_{\ast}\left(  \Omega_{X\left(  \Sigma^{\prime}\right)  }^{n}\right)
=\hat{\Omega}_{X\left(  \Sigma\right)  }^{n}%
\]
and $R^{i}f_{\ast}\left(  \Omega_{X\left(  \Sigma^{\prime}\right)  }%
^{n}\right)  =0$ $\forall i>0$.

If $X\left(  \Sigma\right)  $ is complete, then for any line bundle $L$ on
$X\left(  \Sigma\right)  $%
\[
H^{n-i}\left(  X\left(  \Sigma\right)  ,L^{\ast}\otimes\hat{\Omega}_{X\left(
\Sigma\right)  }^{n}\right)  \cong H^{i}\left(  X\left(  \Sigma\right)
,L\right)  ^{\ast}%
\]

\end{proposition}

\subsubsection{Projective toric
varieties\label{Sec projective toric varieties}}

\paragraph{The normal fan}

If $P\subset M_{\mathbb{R}}$ is a polyhedron and $w\in N_{\mathbb{R}}$, then
define%
\[
\operatorname*{face}\nolimits_{w}\left(  P\right)  =\left\{  m^{\prime}\in
P\mid\left\langle m^{\prime},w\right\rangle \leq\left\langle m,w\right\rangle
\text{ for all }m\in P\right\}
\]
With respect to
\newsym[$\operatorname*{face}\nolimits_{w}\left(  P\right)  $]{face of a polytope}{}Minkowski
sums, $\operatorname*{face}\nolimits_{w}$ has the property that
\[
\operatorname*{face}\nolimits_{w}\left(  P+P^{\prime}\right)
=\operatorname*{face}\nolimits_{w}\left(  P\right)  +\operatorname*{face}%
\nolimits_{w}\left(  P^{\prime}\right)
\]

If $F$ is a face of $P$ define
\newsym[$\sigma_{P}\left(  F\right)  $]{normal cone of the face $F$ of the polytope $P$}{}the
\index{normal cone|textbf}%
\textbf{normal cone} of $F$ as%
\[
\sigma_{P}\left(  F\right)  =\left\{  w\in N_{\mathbb{R}}\mid
\operatorname*{face}\nolimits_{w}\left(  P\right)  =F\right\}
\]
The normal cone has dimension $\dim\left(  \sigma_{P}\left(  F\right)
\right)  =n-\dim\left(  F\right)  $. $F^{\prime}$ is a face of $F$ if and only
if $\sigma_{P}\left(  F\right)  $ is a face of $\sigma_{P}\left(  F^{\prime
}\right)  $. Hence the set of normal cones $\sigma_{P}\left(  F\right)  $ for
all faces $F$ of $P$ forms a fan,
\newsym[$NF\left(  P\right)  $]{normal fan of the polytope $P$}{}the
\index{normal fan|textbf}%
\textbf{normal fan} $\operatorname*{NF}\left(  P\right)  $ of $P$.

Any polyhedron $P$ may be written as%
\[
P=\Delta+C=\left\{  m+m^{\prime}\mid m\in\Delta\text{ and }m^{\prime}\in
C\right\}
\]
with a polytope $\Delta$ and a cone $C$. The cone $C$ is unique and $C^{\ast}$
is the support $\operatorname*{supp}\left(  \operatorname*{NF}\left(
P\right)  \right)  $ of the normal fan of $P$. It is the set of linear forms
on $P$, which have a bounded minimum on $P$. If $P$ is a polytope, then
$\operatorname*{NF}\left(  P\right)  $ is complete.

The normal fan $\Sigma=\operatorname*{NF}\left(  \Delta\right)  $ of the
polytope $\Delta$ consists of all duals%
\[
\sigma_{P}\left(  F\right)  =\left\{  w\in N_{\mathbb{R}}\mid\left\langle
m^{\prime},w\right\rangle \leq\left\langle m,w\right\rangle \text{ for all
}m\in\Delta\text{ and }m^{\prime}\in F\right\}
\]
of the cones%
\[
\left\{  \lambda\left(  m-m^{\prime}\right)  \in M_{\mathbb{R}}\mid m\in
\Delta,\text{ }m^{\prime}\in F,\text{ }\lambda\geq0\right\}
\]
for all nonempty faces $F$ of $\Delta$.

If $0\in\operatorname*{int}\left(  \Delta\right)  $, its normal fan
$\operatorname*{NF}\left(  \Delta\right)  $ is the fan over
\newsym[$\Delta^{\ast}$]{dual polytope of $\Delta$}{}the
\index{dual polytope|textbf}%
\textbf{dual polytope}
\[
\Delta^{\ast}=\left\{  n\in N_{\mathbb{R}}\mid\left\langle m,n\right\rangle
\geq-1\ \forall m\in\Delta\right\}
\]

\paragraph{The projective toric variety associated to an integral polytope}

Given an
\index{integral polytope}%
integral
\index{lattice polytope}%
polytope $\Delta\subset M$, we can
\newsym[$S\left(  \Delta\right)  $]{polytope ring}{}associate to it the
\index{polytope ring|textbf}%
\textbf{polytope ring}%
\[%
\begin{tabular}
[c]{ll}%
$S\left(  \Delta\right)  =\mathbb{C}\left[  t^{k}x^{m}\mid m\in k\Delta
\right]  $ & $\deg t^{k}x^{m}=k$%
\end{tabular}
\]
with
\[
k\Delta=\left\{  km\mid m\in\Delta\right\}  =\overset{k}{\overbrace
{\Delta+...+\Delta}}%
\]
and multiplication $t^{k}x^{m}\cdot t^{l}x^{m^{\prime}}=t^{k+l}x^{m+m^{\prime
}}$, and hence
\newsym[$\mathbb{P}\left(  \Delta\right)  $]{projective toric variety}{}define
a
\index{projective toric variety|textbf}%
\textbf{projective toric variety} $\mathbb{P}\left(  \Delta\right)
=\operatorname*{Proj}S\left(  \Delta\right)  $.

On the other hand we can associate to $\Delta$ its
\index{normal fan}%
normal fan $\Sigma=\operatorname*{NF}\left(  \Delta\right)  $ and a
\index{piecewise linear}%
piecewise linear continuous integral convex function%
\[%
\begin{tabular}
[c]{l}%
$\Phi:N_{\mathbb{R}}\rightarrow\mathbb{R}$\\
$\Phi\left(  w\right)  =\min\limits_{m\in\Delta}\left\langle m,w\right\rangle
=\min\limits_{m\in\Delta\cap M}\left\langle m,w\right\rangle =\min
\limits_{m\in\operatorname*{vertices}\left(  \Delta\right)  }\left\langle
m,w\right\rangle $%
\end{tabular}
\]
giving a
\index{T-Cartier divisor}%
Cartier divisor%
\[
D_{\Delta}=\sum_{r\in\Sigma\left(  1\right)  }-\min_{m\in\Delta}\left\langle
m,\hat{r}\right\rangle D_{r}%
\]
which satisfies $\Delta_{D_{\Delta}}=\Delta$ and is
\index{ample}%
ample.

\begin{theorem}
\cite{Batyrev Dual polyhedra and mirror symmetry for CalabiYau hypersurfaces
in toric varieties}, \cite[Sec. 3.2]{CK Mirror Symmetry and Algebraic
Geometry}, \cite[Sec. 3.4]{Fulton Introduction to Toric Varieties} With this
notation
\begin{align*}
\mathbb{P}\left(  \Delta\right)   &  \cong X\left(  \operatorname*{NF}\left(
\Delta\right)  \right) \\
\mathcal{O}_{\mathbb{P}\left(  \Delta\right)  }\left(  1\right)   &
\cong\mathcal{O}_{\mathbb{P}\left(  \Delta\right)  }\left(  D_{\Delta}\right)
\end{align*}

If $X\left(  \Sigma\right)  $ is complete and $D$ is an
\index{ample}%
ample
\index{T-Cartier divisor}%
$T$-Cartier divisor on $X\left(  \Sigma\right)  $,
\index{normal fan}%
then $\operatorname*{NF}\left(  \Delta_{D}\right)  =\Sigma$.
\end{theorem}

\begin{remark}
Choosing a basis of $M$ gives coordinates $t_{1},...,t_{n}$ on the
\index{torus}%
torus $T=\operatorname*{Hom}\nolimits_{\mathbb{Z}}\left(  M,\mathbb{C}^{\ast
}\right)  $. Writing $m=\left(  a_{1},...,a_{n}\right)  $ we have $x^{m}%
=\prod_{i=1}^{n}t_{i}^{a_{i}}=:t^{m}$. Given $\Delta$ choose $k$ such that
$kD_{\Delta}$ is
\index{very ample}%
very ample on $\mathbb{P}\left(  \Delta\right)  $. The lattice points
$k\Delta\cap M=\left\{  m_{0},...,m_{r}\right\}  $ of $k\Delta$ correspond to
monomials $t^{m_{0}},...,t^{m_{r}}$. $\mathbb{P}\left(  \Delta\right)  $ is
the closure of the image of the map%
\begin{align*}
T  &  \rightarrow\mathbb{P}^{r}\\
t  &  \mapsto\left(  t^{m_{0}},...,t^{m_{r}}\right)
\end{align*}

\end{remark}

\begin{example}
For%
\[
\Delta=\operatorname*{convexhull}\left(  \left(  -1,-1\right)  ,\left(
2,-1\right)  ,\left(  -1,2\right)  \right)
\]
the
\index{normal fan}%
normal fan
\index{hypersurface}%
is the fan of $\mathbb{P}^{2}$ given in Example \ref{2exampleP2fan} and
$\mathbb{P}\left(  \Delta\right)  $ is the closure of the image of the
\index{torus}%
torus under the degree $3$ monomials in $3$ variables, hence it is the degree
$3$
\index{Veronese}%
Veronese embedding of $\mathbb{P}^{2}$.
\end{example}

\subsubsection{The Cox ring of a toric
variety\label{Sec Cox ring of a toric variety}}

In \cite{Cox The homogeneous coordinate ring of a toric variety} the
representation of $\mathbb{P}^{n}$ as
\[
\mathbb{P}^{n}=\left(  \mathbb{C}^{n+1}-V\left(  \left\langle y_{0}%
,...,y_{n}\right\rangle \right)  \right)  /\mathbb{C}^{\ast}%
\]
was generalized to arbitrary toric varieties $X\left(  \Sigma\right)  $. In
order to do so, introduce the homogeneous coordinate ring of a toric variety
$X\left(  \Sigma\right)  $:

\begin{definition}
The
\index{homogeneous coordinate ring|textbf}%
\textbf{homogeneous coordinate ring}
\newsym[$S\left(  X\left(  \Sigma\right)  \right)  $]{homogeneous coordinate ring of the toric variety $X\left(  \Sigma\right)  $}{}or
\index{Cox ring|textbf}%
\textbf{Cox ring}
\index{ray}%
of $X\left(  \Sigma\right)  $ is%
\[
S=S\left(  X\left(  \Sigma\right)  \right)  =\mathbb{C}\left[  y_{r}\mid
r\in\Sigma\left(  1\right)  \right]
\]
with
\index{Cox degree|textbf}%
the
\index{grading}%
grading%
\[
\deg\left(  \prod\nolimits_{r}y_{r}^{a_{r}}\right)  =\left[  \sum
\nolimits_{r}a_{r}D_{r}\right]  \in A_{n-1}\left(  X\left(  \Sigma\right)
\right)
\]

\end{definition}

If $D=\sum\nolimits_{r}a_{r}D_{r}$ write $y^{D}=\prod\nolimits_{r}y_{r}%
^{a_{r}}$, so $\deg\left(  y^{D}\right)  =\left[  D\right]  $. The homogeneous
coordinate ring is the direct
\newsym[$S_{\alpha}$]{degree $\alpha$ subspace of $S$}{}sum
\[
S=\bigoplus\nolimits_{\alpha\in A_{n-1}\left(  X\left(  \Sigma\right)
\right)  }S_{\alpha}%
\]
with%
\[
S_{\alpha}=\bigoplus\limits_{\left[  D\right]  =\alpha}\mathbb{C\cdot}y^{D}%
\]
and it holds $S_{\alpha}\cdot S_{\beta}\subset S_{\alpha+\beta}$.

\subsubsection{Global sections as Cox
monomials\label{Sec Global sections a cox monomials}}

\begin{proposition}
\cite[Sec. 3.2]{CK Mirror Symmetry and Algebraic Geometry} The
\index{global sections}%
global sections of the
\index{reflexive sheaf}%
reflexive sheaf of sections $\mathcal{O}_{X\left(  \Sigma\right)  }\left(
D\right)  $ of a
\index{Weil divisor}%
Weil divisor $D$ is isomorphic to the degree $\left[  D\right]  $-part of the
\index{Cox ring}%
Cox ring%
\begin{align*}
H^{0}\left(  X\left(  \Sigma\right)  ,\mathcal{O}_{X\left(  \Sigma\right)
}\left(  D\right)  \right)   &  \longrightarrow S_{\left[  D\right]  }\\
x^{m}  &  \mapsto\prod_{r}y_{r}^{\left\langle m,\hat{r}\right\rangle +b_{r}}%
\end{align*}
where $D=\sum b_{r}D_{r}$.
\end{proposition}

In particular $S_{\left[  D\right]  }$ is finite dimensional of dimension
$\dim_{\mathbb{C}}\left(  S_{\left[  D\right]  }\right)  =\left\vert
\Delta_{D}\cap M\right\vert $.

\begin{remark}
The
\index{homogeneous coordinate ring}%
homogeneous coordinate ring contains all possible
\index{polytope ring}%
polytope rings associated to
\index{ample}%
ample divisors of a
\index{projective toric variety}%
projective toric variety. Indeed if $D$ is ample on $Y$, then%
\[
S\left(  \Delta_{D}\right)  \cong\bigoplus_{d=0}^{\infty}S_{k\left[  D\right]
}%
\]

\end{remark}

Given $D=\sum b_{r}D_{r}$, we can describe $S_{\left[  D\right]  }$
explicitly: In the presentation of the
\index{Chow group}%
Chow group of divisors of $X\left(  \Sigma\right)  $%
\[%
\begin{tabular}
[c]{lllllll}%
$0\rightarrow$ & $M$ & $\overset{A}{\rightarrow}$ & $\mathbb{Z}^{\Sigma\left(
1\right)  }$ & $\overset{\deg}{\rightarrow}$ & $A_{n-1}\left(  X\left(
\Sigma\right)  \right)  $ & $\rightarrow0$\\
& $m$ & $\mapsto$ & $\left(  \left\langle m,\hat{r}\right\rangle \right)
_{r\in\Sigma\left(  1\right)  }$ &  &  & \\
&  &  & $\left(  a_{r}\right)  _{r\in\Sigma\left(  1\right)  }$ & $\mapsto$ &
$\sum_{r\in\Sigma\left(  1\right)  }a_{r}D_{r}$ &
\end{tabular}
\
\]
with the rows of $A$ being the minimal lattice generators of the elements
\index{ray}%
of $\Sigma\left(  1\right)  $, we have $\operatorname*{image}\left(  A\right)
=\ker\left(  \deg\right)  $. Hence the Cox monomials of the same
\index{Cox degree}%
Cox degree as $D$, i.e., giving divisors linearly equivalent to $D$ (so these
torus invariant elements form a vector space basis of the space of
\index{global sections}%
global sections of $D$), are
\[
\left\{  y^{a}\mid a\in\left(  b_{r}\right)  +\operatorname*{image}\left(
A\right)  \text{, }a\in\mathbb{Z}_{\geq0}^{\Sigma\left(  1\right)  }\right\}
\]

\begin{example}
\label{2quinticCoxmon}For $X\left(  \Sigma\right)  $ given by the fan over the
faces of%
\[
\operatorname*{convexhull}\left(  \left(  4,-1,-1,-1\right)  ,...,\left(
-1,-1,-1,4\right)  ,\left(  -1,-1,-1,-1\right)  \right)
\]
as in
\index{quintic threefold}%
Example \ref{2exquinticchow} the Cox monomials in $S_{\left[  -K_{X\left(
\Sigma\right)  }\right]  }$, i.e., the monomials of the same
\index{hypersurface}%
degree as the
\index{anticanonical class}%
anticanonical class, are%
\[
\left(  \left(
\begin{array}
[c]{c}%
1\\
1\\
1\\
1\\
1
\end{array}
\right)  +\operatorname*{image}\left(
\begin{array}
[c]{cccc}%
4 & -1 & -1 & -1\\
-1 & 4 & -1 & -1\\
-1 & -1 & 4 & -1\\
-1 & -1 & -1 & 4\\
-1 & -1 & -1 & -1
\end{array}
\right)  \right)  \cap\mathbb{Z}_{\geq0}^{5}=\left\{  \left(
\begin{array}
[c]{c}%
1\\
1\\
1\\
1\\
1
\end{array}
\right)  ,\left(
\begin{array}
[c]{c}%
5\\
0\\
0\\
0\\
0
\end{array}
\right)  ,...,\left(
\begin{array}
[c]{c}%
0\\
0\\
0\\
0\\
5
\end{array}
\right)  \right\}
\]
i.e., the monomials $y_{1}^{5},...,y_{5}^{5},y_{1}\cdot...\cdot y_{5}$.

For $\mathbb{P}^{4}$, which
\index{Greene-Plesser}%
is given by the fan over the faces of the polytope%
\[
\operatorname*{convexhull}\left(  \left(  1,0,0,0\right)  ,...,\left(
0,0,0,1\right)  ,\left(  -1,-1,-1,-1\right)  \right)
\]
the Cox monomials in $S_{\left[  -K_{\mathbb{P}^{4}}\right]  }$ are%
\[
\left(  \left(
\begin{array}
[c]{c}%
1\\
1\\
1\\
1\\
1
\end{array}
\right)  +\operatorname*{image}\left(
\begin{array}
[c]{cccc}%
1 & 0 & 0 & 0\\
0 & 1 & 0 & 0\\
0 & 0 & 1 & 0\\
0 & 0 & 0 & 1\\
-1 & -1 & -1 & -1
\end{array}
\right)  \right)  \cap\mathbb{Z}_{\geq0}^{5}%
\]
yielding all monomials of homogeneous degree $5$ in $5$ variables.
\end{example}

\begin{algorithm}
In order to compute $S_{\left[  D\right]  }$ for given $D\in\mathbb{Z}%
^{\Sigma\left(  1\right)  }$ consider a basis $v_{1},...,v_{s}$ of
$\ker\left(  A^{t}\right)  $ and consider the cone $C\subset\left(
\mathbb{R}^{\Sigma\left(  1\right)  }\right)  ^{\ast}\oplus\mathbb{R}$ with
the rays
\[
\operatorname*{hull}\left\{  \left(  e_{r},0\right)  \right\}  \text{ for
}r\in\Sigma\left(  1\right)
\]
and lineality space (i.e., largest linear space contained in $C$) spanned by
\[
\left\{  \left(  v_{i},-v_{i}\cdot D\right)  \mid i=1,...s\right\}
\]
Then the intersection of $C^{\ast}$ with the hyperplane defined by setting the
last coordinate equal to $1$ is the polytope $P=\left(
D+\operatorname*{image}\left(  A\right)  \right)  \cap\mathbb{R}_{\geq
0}^{\Sigma\left(  1\right)  }$. Consider the preimage in $M_{\mathbb{R}}$ of
this polytope under the map $m\mapsto A\cdot m+D$. The lattice points of this
polytope map via $m\mapsto A\cdot m+D$ to a basis of $S_{\left[  D\right]  }$.
\end{algorithm}

For $y^{D},y^{E}\in S$ define%
\[
y^{D}<y^{E}\Leftrightarrow\exists y^{F}\in S\text{ such that }\left[
F\right]  =\left[  E\right]  \text{, }y^{D}\mid y^{F}\text{ and }y^{D}\neq
y^{F}%
\]
Under this condition $\left[  E\right]  -\left[  D\right]  =\left[  F\right]
-\left[  D\right]  =\left[  F-D\right]  $ is the class of an effective divisor.

\begin{lemma}
\cite{Cox The homogeneous coordinate ring of a toric variety} If $X\left(
\Sigma\right)  $ is complete, then $>$ is a transitive, antisymmetric,
multiplicative ordering on the monomials of $S$.
\end{lemma}

\subsubsection{Homogeneous coordinate presentation of toric
varieties\label{1homogeneouscoordinate}}

\paragraph{Quotient presentations}

Let $q:Y^{\prime}\rightarrow Y$ be a surjective morphism of toric varieties
with tori $T\subset Y$ and $T^{\prime}\subset Y^{\prime}$, and denote by
$q^{\ast}:\operatorname*{Div}_{T}\left(  Y\right)  \rightarrow
\operatorname*{Div}_{T}\left(  Y^{\prime}\right)  $ the pullback. With
$U=T\cup\bigcup_{r\in\Sigma\left(  1\right)  }O\left(  r\right)  $ there is
the strict transform $q^{\#}$%
\[%
\begin{tabular}
[c]{ccc}%
$\operatorname*{Div}_{T}\left(  U\right)  $ & $\overset{q^{\ast}}{\rightarrow
}$ & $\operatorname*{Div}_{T^{\prime}}\left(  U^{\prime}\right)  $\\
&  & $\cap$\\
$\parallel$ &  & $\operatorname*{WDiv}_{T^{\prime}}\left(  U^{\prime}\right)
$\\
&  & $\cap$\\
$\operatorname*{WDiv}_{T}\left(  Y\right)  $ & $\overset{q^{\#}}%
{\hookrightarrow}$ & $\operatorname*{WDiv}_{T^{\prime}}\left(  Y^{\prime
}\right)  $%
\end{tabular}
\]
If $Y$ is a toric variety, then a
\index{quotient presentation|textbf}%
\textbf{quotient presentation} of $Y$ is a quasiaffine toric variety
$Y^{\prime}$ and a surjective, affine toric morphism $q:Y^{\prime}\rightarrow
Y$ such that $q^{\#}$ is bijective. This can be tested locally for all
invariant affine open $U\subset Y$.

\begin{theorem}
\cite{ACHS Homogeneous coordinates and quotient presentations for toric
varieties} Suppose $Y=X\left(  \Sigma\right)  $ and $Y^{\prime}=X\left(
\Sigma^{\prime}\right)  $ are toric varieties given by fans $\Sigma\subset
N_{\mathbb{R}}$ and $\Sigma^{\prime}\subset N_{\mathbb{R}}^{\prime}$ and
$q:Y^{\prime}\rightarrow Y$ is a toric morphism given by a homomorphism of
lattices $\varphi:N^{\prime}\rightarrow N$ as described in Section
\ref{Sec morphisms of toric varieties}. Then $q$ is a quotient presentation if
and only if the following conditions are satisfied:

\begin{itemize}
\item $\operatorname*{coker}\left(  \varphi\right)  $ is finite.

\item there is a strongly convex rational polyhedral cone $\overline{\sigma
}\subset N_{\mathbb{R}}^{\prime}$ such that $\Sigma^{\prime}$ is a subfan of a
fan $\Sigma^{\prime\prime}$ spanned by $\overline{\sigma}$ (so $Y^{\prime
}\subset U\left(  \overline{\sigma}\right)  $).

\item the map $\sigma\mapsto\varphi_{\mathbb{R}}\left(  \sigma\right)  $ is a
bijection $\Sigma^{\prime\max}\rightarrow\Sigma^{\max}$ and $\Sigma^{\prime
}\left(  1\right)  \rightarrow\Sigma\left(  1\right)  $.

\item for all rays $r^{\prime}\in\Sigma^{\prime}\left(  1\right)  $ the image
$\varphi\left(  \hat{r}^{\prime}\right)  $ of a
\index{minimal lattice generator}%
minimal lattice generator $\hat{r}^{\prime}$ is a primitive lattice element of
$N$.
\end{itemize}
\end{theorem}

If $q:Y^{\prime}\rightarrow Y$ is a quotient presentation of $Y$, then via the
isomorphism $q^{\#}$ we get a commutative diagram%
\[%
\begin{tabular}
[c]{lcll}
& $0$ &  & \\
& $\downarrow$ &  & \\
$0\longrightarrow$ & \multicolumn{1}{l}{$M$} & $\overset{\operatorname*{div}%
}{\longrightarrow}$ & $\operatorname*{WDiv}\nolimits_{T}\left(  Y\right)  $\\
& $\downarrow$ & $\nearrow$ & \\
& \multicolumn{1}{l}{$M^{\prime}$} &  &
\end{tabular}
\]

\begin{definition}
A triangle is a lattice $M^{\prime}$ and a commutative diagram%
\[%
\begin{tabular}
[c]{lll}%
$M$ & $\overset{\operatorname*{div}}{\longrightarrow}$ & $\operatorname*{WDiv}%
\nolimits_{T}\left(  Y\right)  $\\
\multicolumn{1}{c}{$\downarrow$} & $\nearrow$ & \\
$M^{\prime}$ &  &
\end{tabular}
\]
such that $M\longrightarrow M^{\prime}$ is injective and for all $T$-invariant
open $U\subset Y$ there is an $m^{\prime}\in M^{\prime}$ such that the
associated divisor on $Y$ is effective with support $Y\backslash U$.{}
\end{definition}

\begin{theorem}
\cite{ACHS Homogeneous coordinates and quotient presentations for toric
varieties} Let $Y$ be a toric variety. Above commutative diagram associated to
a quotient presentation is a triangle. Up to isomorphism, this assignment is a
bijection between quotient presentations and triangles.
\end{theorem}

\begin{example}
The triangle given by $M^{\prime}=\operatorname*{WDiv}\nolimits_{T}\left(
Y\right)  \overset{id}{\rightarrow}\operatorname*{WDiv}\nolimits_{T}\left(
Y\right)  $ defines the Cox quotient presentation explored in detail in the
following Section \ref{Sec Cox quotient presentation}.
\end{example}

\begin{example}
If $D$ is an ample Cartier divisor on $Y$, then%
\[
M\rightarrow M\oplus\mathbb{Z}D\rightarrow\operatorname*{WDiv}\nolimits_{T}%
\left(  Y\right)
\]
is a triangle and the corresponding quotient presentation is the associated
$\mathbb{C}^{\ast}$-bundle of $\mathcal{O}_{Y}\left(  D\right)  $.
\end{example}

Suppose $q:Y^{\prime}\rightarrow Y$ is a quotient presentation given by the
triangle $M\rightarrow M^{\prime}\rightarrow\operatorname*{WDiv}%
\nolimits_{T}\left(  Y\right)  $. Denote by $T$ and $T^{\prime}$ the tori of
$Y$ and $Y^{\prime}$ and let%
\[
G=\ker\left(  T^{\prime}\rightarrow T\right)
\]
With $A=M^{\prime}/M$ we have $G=\operatorname*{Spec}\left(  \mathbb{C}\left[
A\right]  \right)  $ and $\widehat{G}=A$. Then $q_{\ast}\mathcal{O}%
_{Y^{\prime}}$ is graded by $A$ with $\mathcal{O}_{Y}$-modules $R_{a}$%
\[
q_{\ast}\mathcal{O}_{Y^{\prime}}=\bigoplus_{a\in A}R_{a}%
\]

The group $G$ acts on $Y^{\prime}$ and the morphism $q$ is a good quotient if
\[
\left(  q_{\ast}\mathcal{O}_{Y^{\prime}}\right)  ^{G}=\mathcal{O}_{Y}%
\]
One can test this condition locally, so assume that $q$ is given by an
inclusion $\mathbb{C}\left[  \sigma^{\vee}\cap M\right]  \subset
\mathbb{C}\left[  \sigma^{\prime\vee}\cap M^{\prime}\right]  $. One can show
\[
\mathbb{C}\left[  \sigma^{\prime\vee}\cap M^{\prime}\right]  ^{G}%
=\mathbb{C}\left[  \sigma^{\vee}\cap M\right]
\]
hence:

\begin{proposition}
\cite{ACHS Homogeneous coordinates and quotient presentations for toric
varieties} Any quotient presentation of a toric variety is a good quotient.
\end{proposition}

The morphism $q$ is a categorial quotient and for all closed invariant $W_{i}$
it holds $q\left(  \bigcap_{i}W_{i}\right)  =\bigcap_{i}q\left(  W_{i}\right)
$.

If $Y$ is simplicial, then $q$ is a geometric quotient. Any quotient
presentation is geometric in codimension $2$.

\paragraph{Cox quotient presentation of toric
varieties\label{Sec Cox quotient presentation}}

Suppose $\Sigma\left(  1\right)  $ spans $N_{\mathbb{R}}$. Applying
$\operatorname*{Hom}\nolimits_{\mathbb{Z}}\left(  -,\mathbb{C}^{\ast}\right)
$ to the presentation%
\[%
\begin{tabular}
[c]{lllllll}%
$0\rightarrow$ & $M$ & $\overset{A}{\rightarrow}$ & $\overset{\mathbb{Z}%
^{\Sigma\left(  1\right)  }}{\overset{\cong}{\operatorname*{WDiv}_{T}}}\left(
X\left(  \Sigma\right)  \right)  $ & $\overset{\deg}{\rightarrow}$ &
$A_{n-1}\left(  X\left(  \Sigma\right)  \right)  $ & $\rightarrow0$%
\end{tabular}
\
\]
of
\index{Chow group}%
$A_{n-1}\left(  X\left(  \Sigma\right)  \right)  $, we get an exact sequence%
\[%
\begin{tabular}
[c]{ccccccc}%
$1\rightarrow$ & $G\left(  \Sigma\right)  $ & $\rightarrow$ &
$\operatorname*{Hom}\nolimits_{\mathbb{Z}}\left(  \operatorname*{WDiv}%
_{T}\left(  X\left(  \Sigma\right)  \right)  ,\mathbb{C}^{\ast}\right)  $ &
$\rightarrow$ & $\operatorname*{Hom}\nolimits_{\mathbb{Z}}\left(
M,\mathbb{C}^{\ast}\right)  $ & $\rightarrow1$\\
&  &  & $\shortparallel$ &  & $\shortparallel$ & \\
&  &  & $\left(  \mathbb{C}^{\ast}\right)  ^{\Sigma\left(  1\right)  }$ &  &
$T$ &
\end{tabular}
\
\]
with the kernel
\[
G\left(  \Sigma\right)  =\operatorname*{Hom}\nolimits_{\mathbb{Z}}\left(
A_{n-1}\left(  X\left(  \Sigma\right)  \right)  ,\mathbb{C}^{\ast}\right)
\]
of the map of tori, hence the inclusion
\newsym[$G\left(  \Sigma\right)  $]{group acting on the Cox quotient presentation space}{}of
$G\left(  \Sigma\right)  $ in $\operatorname*{Hom}\nolimits_{\mathbb{Z}%
}\left(  \mathbb{Z}^{\Sigma\left(  1\right)  },\mathbb{C}^{\ast}\right)  $
gives an action
\begin{align*}
G\left(  \Sigma\right)  \times\operatorname*{Hom}\nolimits_{\mathbb{Z}}\left(
\operatorname*{WDiv}\nolimits_{T}\left(  X\left(  \Sigma\right)  \right)
,\mathbb{C}^{\ast}\right)   &  \rightarrow\operatorname*{Hom}%
\nolimits_{\mathbb{Z}}\left(  \operatorname*{WDiv}\nolimits_{T}\left(
X\left(  \Sigma\right)  \right)  ,\mathbb{C}^{\ast}\right) \\
\left(  g,a\right)   &  \mapsto%
\begin{tabular}
[c]{llll}%
$ga:$ & $\operatorname*{WDiv}\nolimits_{T}\left(  X\left(  \Sigma\right)
\right)  $ & $\rightarrow$ & $\mathbb{C}^{\ast}$\\
& $D_{r}$ & $\mapsto$ & $g\left(  \left[  D_{r}\right]  \right)  a\left(
D_{r}\right)  $%
\end{tabular}
\end{align*}
which induces an action of $G\left(  \Sigma\right)  $ on
\[
\operatorname*{Hom}\nolimits_{sg}\left(  \operatorname*{WDiv}\nolimits_{T}%
\left(  X\left(  \Sigma\right)  \right)  ,\mathbb{C}\right)
=\operatorname*{Specm}\left(  S\right)  =\mathbb{C}^{\Sigma\left(  1\right)  }%
\]
considering $\mathbb{C=C}^{\ast}\cup\left\{  0\right\}  $ as a semigroup with
respect to multiplication. This action is given by%
\begin{align*}
G\left(  \Sigma\right)  \times\operatorname*{Hom}\nolimits_{sg}\left(
\operatorname*{WDiv}\nolimits_{T}\left(  X\left(  \Sigma\right)  \right)
,\mathbb{C}\right)   &  \rightarrow\operatorname*{Hom}\nolimits_{sg}\left(
\operatorname*{WDiv}\nolimits_{T}\left(  X\left(  \Sigma\right)  \right)
,\mathbb{C}\right) \\
\left(  g,a\right)   &  \mapsto%
\begin{tabular}
[c]{llll}%
$ga:$ & $\operatorname*{WDiv}\nolimits_{T}\left(  X\left(  \Sigma\right)
\right)  $ & $\rightarrow$ & $\mathbb{C}$\\
& $D_{r}$ & $\mapsto$ & $g\left(  \left[  D_{r}\right]  \right)  a\left(
D_{r}\right)  $%
\end{tabular}
\end{align*}

The group $G\left(  \Sigma\right)  $ is isomorphic to the product of a torus
$\left(  \mathbb{C}^{\ast}\right)  ^{\operatorname*{rank}\left(
A_{n-1}\left(  X\left(  \Sigma\right)  \right)  \right)  }$ and the finite
group $\operatorname*{Hom}\nolimits_{\mathbb{Z}}\left(  A_{n-1}\left(
X\left(  \Sigma\right)  \right)  _{tor},\mathbb{Q}/\mathbb{Z}\right)  $.

\begin{definition}
If $\sigma\in\Sigma$ is a cone define the
\newsym[$D_{\widehat{\sigma}}$]{irrelevant divisor of $\sigma$}{}divisor%
\[
D_{\widehat{\sigma}}=\sum_{r\in\Sigma\left(  1\right)  ,\text{ }%
r\not \subset \sigma}D_{r}%
\]
and the
\index{irrelevant ideal|textbf}%
\textbf{irrelevant ideal} of $X\left(  \Sigma\right)  $
\newsym[$B\left(  \Sigma\right)  $]{irrelevant ideal}{}by%
\[
B\left(  \Sigma\right)  =\left\langle y^{D_{\widehat{\sigma}}}\mid\sigma
\in\Sigma\right\rangle =\left\langle \prod_{r\in\Sigma\left(  1\right)
,\text{ }r\not \subset \sigma}y_{r}\mid\sigma\in\Sigma\right\rangle \subset S
\]

\end{definition}

If $\sigma\in\Sigma$ is a cone, then%
\[
U_{\sigma}=\mathbb{C}^{\Sigma\left(  1\right)  }-V\left(  y^{D_{\widehat
{\sigma}}}\right)
\]
is invariant \newsym[$U_{\sigma}$]{principal open set}{}under the action of
$G\left(  \Sigma\right)  $. So%
\[
\mathbb{C}^{\Sigma\left(  1\right)  }-V\left(  B\left(  \Sigma\right)
\right)  =\bigcup\limits_{\sigma\in\Sigma}U_{\sigma}%
\]
is invariant under $G\left(  \Sigma\right)  $. The localization $S_{\sigma
}=S_{y^{\widehat{\sigma}}}$ is the coordinate ring of the affine variety
$U_{\sigma}$ and the invariants under the action of $G\left(  \Sigma\right)  $
are%
\[
\left(  S_{\sigma}\right)  ^{G\left(  \Sigma\right)  }=\left(  S_{\sigma
}\right)  _{0}\cong\mathbb{C}\left[  \check{\sigma}\cap M\right]
\]
so%
\[
U_{\sigma}/G\left(  \Sigma\right)  =\operatorname*{Spec}\left(  \left(
S_{\sigma}\right)  ^{G\left(  \Sigma\right)  }\right)  =\operatorname*{Spec}%
\mathbb{C}\left[  \check{\sigma}\cap M\right]  =U\left(  \sigma\right)
\]
is the affine toric variety $U\left(  \sigma\right)  \subset X\left(
\Sigma\right)  $.

\begin{theorem}
\cite[Sec. 3.2]{CK Mirror Symmetry and Algebraic Geometry} Suppose
$\Sigma\left(  1\right)  $ spans $N_{\mathbb{R}}$. Then with the above action
of
\[
G\left(  \Sigma\right)  =\operatorname*{Hom}\nolimits_{\mathbb{Z}}\left(
A_{n-1}\left(  X\left(  \Sigma\right)  \right)  ,\mathbb{C}^{\ast}\right)
\]
on $\mathbb{C}^{\Sigma\left(  1\right)  }$ and above irrelevant ideal
$B\left(  \Sigma\right)  $ it holds%
\begin{equation}
X\left(  \Sigma\right)  =\left(  \mathbb{C}^{\Sigma\left(  1\right)
}-V\left(  B\left(  \Sigma\right)  \right)  \right)  //G\left(  \Sigma\right)
\label{4homogquotienttoric}%
\end{equation}

The quotient is geometric if and only if $\Sigma$ is
\index{simplicial}%
simplicial.
\end{theorem}

\begin{example}
\label{2exquinticquotientdescription}If $\Sigma$ is the fan over the degree
$5$
\index{Veronese}%
Veronese polytope of $\mathbb{P}^{4}$ as considered
\index{Chow group}%
Example \ref{2exquinticchow}
\begin{gather*}
A_{3}\left(  X\left(  \Sigma\right)  \right)  \cong\mathbb{Z\times}H\\
\text{with }H=\frac{\left\{  \left(  a_{0},...,a_{4}\right)  \in\mathbb{Z}%
_{5}^{5}\mid\sum_{i=0}^{4}a_{i}=0\operatorname{mod}5\right\}  }{\mathbb{Z}%
_{5}\left(  1,1,1,1,1\right)  }\cong\mathbb{Z}_{5}^{3}%
\end{gather*}
and $\operatorname*{Hom}\nolimits_{\mathbb{Z}}\left(  A_{3}\left(  X\left(
\Sigma\right)  \right)  ,\mathbb{C}^{\ast}\right)  =\mathbb{C}^{\ast}%
\times\mathbb{Z}_{5}^{3}$ acts on $\mathbb{C}^{5}$ by
\[
\left(  \lambda,\left(  \mu^{a_{0}},...,\mu^{a_{4}}\right)  \right)
\cdot\left(  y_{0},...,y_{4}\right)  =\left(  \lambda\mu^{a_{0}}%
y_{0},...,\lambda\mu^{a_{4}}y_{4}\right)
\]
where $\mu$ is a $5$th
\index{root of unity}%
root of unity. Furthermore,%
\[
B\left(  \Sigma\right)  =\left\langle y_{0},y_{1},y_{2},y_{3},y_{4}%
\right\rangle \subset\mathbb{C}\left[  y_{0},y_{1},y_{2},y_{3},y_{4}\right]
\]
hence $X\left(  \Sigma\right)  =\mathbb{P}^{4}/\mathbb{Z}_{5}^{3}$, which is
precisely the quotient of $\mathbb{P}^{4}$, the
\index{quintic threefold}%
Greene-Plesser
\index{Greene-Plesser}%
mirror of the generic quintic sits inside (see Example
\ref{GreenePlesserQuinticExample}).
\end{example}

\begin{remark}
For the practical representation of the action of $G\left(  \Sigma\right)  $
on $\mathbb{C}^{\Sigma\left(  1\right)  }$ we proceed as follows: Choose a
numbering of rays of $\Sigma$, let $r=\left\vert \Sigma\left(  1\right)
\right\vert $ be the number of rays and denote by $A$ the presentation matrix
of $A_{n-1}\left(  Y\right)  $. By Smith normal form we obtain $W\in
\operatorname*{GL}\left(  n,\mathbb{Z}\right)  $ and $U\in\operatorname*{GL}%
\left(  r,\mathbb{Z}\right)  $ and a commutative diagram%
\[%
\begin{tabular}
[c]{lllllll}%
$0\rightarrow$ & $\mathbb{Z}^{n}$ & $\overset{A}{\rightarrow}$ &
$\mathbb{Z}^{r}$ & $\rightarrow$ & $A_{n-1}\left(  Y\right)  $ &
$\rightarrow0$\\
& $\downarrow W$ &  & $\downarrow U$ &  & $\downarrow\cong$ & \\
$0\rightarrow$ & $\mathbb{Z}^{n}$ & $\overset{A^{\prime}}{\rightarrow}$ &
$\mathbb{Z}^{r}$ & $\rightarrow$ & $H$ & $\rightarrow0$%
\end{tabular}
\]
such that $A^{\prime}$ is a matrix with non zero entries only on the diagonal.
Then
\[
G\left(  \Sigma\right)  ^{\prime}=\operatorname*{Hom}\nolimits_{\mathbb{Z}%
}\left(  H,\mathbb{C}^{\ast}\right)
\]
acts by%
\begin{align*}
G\left(  \Sigma\right)  ^{\prime}\times\mathbb{C}^{r}  &  \rightarrow
\mathbb{C}^{r}\\
\left(  \left(  t_{j}\right)  ,\left(  a_{j}\right)  \right)   &
\mapsto\left(
{\textstyle\prod\nolimits_{i=1}^{r}}
t_{i}^{u_{ij}}a_{j}\right)  _{j=1,...,r}%
\end{align*}
where $U=\left(  u_{ij}\right)  $, and it holds%
\[
X\left(  \Sigma\right)  =\left(  \mathbb{C}^{\Sigma\left(  1\right)
}-V\left(  B\left(  \Sigma\right)  \right)  \right)  //G\left(  \Sigma\right)
^{\prime}%
\]

\end{remark}

\subsubsection{Homogeneous coordinate representations of subvarieties and
sheaves\label{Sec Homogeneous coordinate representation of subvarieties and sheaves}%
}

Let $Y=X\left(  \Sigma\right)  $ be simplicial and $I\subset S$ a graded
ideal. Then $V\left(  I\right)  -V\left(  B\left(  \Sigma\right)  \right)
\subset\mathbb{C}^{\Sigma\left(  1\right)  }-V\left(  B\left(  \Sigma\right)
\right)  $ is $G\left(  \Sigma\right)  $-invariant. As
\[
X\left(  \Sigma\right)  =\left(  \mathbb{C}^{\Sigma\left(  1\right)
}-V\left(  B\left(  \Sigma\right)  \right)  \right)  /G\left(  \Sigma\right)
\]
is a geometric quotient, the $G\left(  \Sigma\right)  $-invariant Zariski
closed subsets of $\mathbb{C}^{\Sigma\left(  1\right)  }-V\left(  B\left(
\Sigma\right)  \right)  $ are in one-to-one correspondence to the Zariski
closed subsets of $X\left(  \Sigma\right)  $. Denote by $V_{Y}\left(
I\right)  $ the
\newsym[$V_{Y}\left(  I\right)  $]{Zariski closed subset of toric variety}{}Zariski
closed subset of $Y$ corresponding to $V\left(  I\right)  -V\left(  B\left(
\Sigma\right)  \right)  $.

So $V_{Y}\left(  I\right)  =\varnothing$ if and only if $V\left(  I\right)
\subset V\left(  B\left(  \Sigma\right)  \right)  $, which is equivalent to
the existence of an $m$ with $B\left(  \Sigma\right)  ^{m}\subset I$ by the Nullstellensatz.

\begin{proposition}
\cite{Cox The homogeneous coordinate ring of a toric variety} Let $Y=X\left(
\Sigma\right)  $ be a simplicial toric variety. Then

\begin{enumerate}
\item For any graded ideal $I\subset S$%
\[
V_{Y}\left(  I\right)  =\varnothing\Leftrightarrow\exists m:B\left(
\Sigma\right)  ^{m}\subset I
\]

\item There is a one-to-one correspondence%
\[%
\begin{tabular}
[c]{ccc}%
$\left\{  \text{graded radical ideals }I\subset S\text{ with }B\left(
\Sigma\right)  \subset I\right\}  $ & $\rightarrow$ & $\left\{  \text{Zariski
closed subsets of }Y\right\}  $\\
$I$ & $\mapsto$ & $V_{Y}\left(  I\right)  $%
\end{tabular}
\]

\end{enumerate}
\end{proposition}

A graded $S$-module $F$ has a decomposition into a direct sum%
\[
F=\bigoplus\nolimits_{\alpha\in A_{n-1}\left(  X\left(  \Sigma\right)
\right)  }F_{\alpha}%
\]
with $S_{\alpha}F_{\beta}\subset F_{\alpha+\beta}$.

Let $\sigma\in\Sigma$ be a cone. The degree $0$ part $\left(  F_{\sigma
}\right)  _{0}$ of the graded $S_{\sigma}$-module $F_{\sigma}=F\otimes
_{S}S_{\sigma}$ is an $\left(  S_{\sigma}\right)  _{0}$-module, which defines
a quasi-coherent sheaf $\widetilde{\left(  F_{\sigma}\right)  }_{0}$ on the
affine toric variety $U\left(  \sigma\right)  =\operatorname*{Spec}\left(
\left(  S_{\sigma}\right)  _{0}\right)  \subset X\left(  \Sigma\right)  $.
According to the fan the sheaves $\widetilde{\left(  F_{\sigma}\right)  }_{0}$
patch to a quasi-coherent sheaf $\widetilde{F}$.

\begin{theorem}
\cite{Cox The homogeneous coordinate ring of a toric variety} Let $Y=X\left(
\Sigma\right)  $ with Cox ring $S$. The map $F\mapsto\widetilde{F}$ is an
exact functor from the graded $S$-modules to quasi-coherent $O_{Y}$-modules.
It has the following properties:

\begin{itemize}
\item If $Y$ is simplicial, then every quasi-coherent sheaf $\mathcal{F}$
arises in this way as $\mathcal{F}\cong\widetilde{F}$ with%
\[
F=\bigoplus\limits_{\alpha\in A_{n-1}\left(  X\left(  \Sigma\right)  \right)
}H^{0}\left(  Y,\mathcal{F}\otimes_{\mathcal{O}_{Y}}\widetilde{S\left(
\alpha\right)  }\right)
\]
where $S\left(  \alpha\right)  _{\beta}=S_{\alpha+\beta}$.

\item If $F$ is finitely generated, then $\widetilde{F}$ is coherent.

\item If $Y$ is simplicial, then every coherent sheaf on $Y$ is of the form
$\widetilde{F}$ with $F$ finitely generated.

$\widetilde{F}=0$ if and only if there is some $k>0$ such that $B\left(
\Sigma\right)  ^{k}F_{\alpha}=\left\{  0\right\}  $ for all $\alpha
\in\operatorname*{Pic}\left(  Y\right)  $.

\item If $Y$ is smooth, then $\widetilde{F}=0$ if and only if there is some
$k>0$ such that $B\left(  \Sigma\right)  ^{k}F=\left\{  0\right\}  $.
\end{itemize}
\end{theorem}

\begin{theorem}
\cite{Cox The homogeneous coordinate ring of a toric variety} Let $Y=X\left(
\Sigma\right)  $ with Cox ring $S$.

\begin{enumerate}
\item If $Y$ is simplicial, then any closed subscheme of $Y$ is given by a
graded ideal $I\subset S$, and graded ideals $I,J\subset S$ correspond to the
same closed subscheme of $Y$ if and only if $\left(  I:B\left(  \Sigma\right)
^{\infty}\right)  _{\alpha}=\left(  J:B\left(  \Sigma\right)  ^{\infty
}\right)  _{\alpha}$ for all $\alpha\in\operatorname*{Pic}\left(  Y\right)  $.

\item If $Y$ is smooth, then graded ideals $I,J\subset S$ correspond to the
same closed subscheme of $Y$ if and only if $\left(  I:B\left(  \Sigma\right)
^{\infty}\right)  =\left(  J:B\left(  \Sigma\right)  ^{\infty}\right)  $, so
there is a one-to-one correspondence between the graded ideals $I\subset S$
which are saturated in $B\left(  \Sigma\right)  $ and the closed subschemes of
$Y$.
\end{enumerate}
\end{theorem}

\subsubsection{K\"{a}hler cone and Mori cone\label{Sec Kaehler cone Mori cone}%
}

Suppose $Y=X\left(  \Sigma\right)  $ is a
\index{simplicial}%
simplicial projective toric variety of dimension $n$ given by the fan
$\Sigma\subset N_{\mathbb{R}}$. Then%
\[
A_{n-1}\left(  Y\right)  \otimes\mathbb{R}\cong H^{2}\left(  Y,\mathbb{R}%
\right)
\]

The
\index{K\"{a}hler cone|textbf}%
\textbf{K\"{a}hler cone }$K\left(  Y\right)  $
\newsym[$K\left(  Y\right)  $]{K\"{a}hler cone of $Y$}{}of $Y$ is the cone of
all
\index{K\"{a}hler classes}%
K\"{a}hler classes on $Y$ considered as a subset in $A_{n-1}\left(  Y\right)
\otimes\mathbb{R}$ or $H^{2}\left(  Y,\mathbb{R}\right)  $.

The
\newsym[$A_{n-1}^{+}\left(  Y\right)  \otimes\mathbb{R}$]{cone generated by prime $T$-divisor classes}{}cone
$A_{n-1}^{+}\left(  Y\right)  \otimes\mathbb{R}$ is defined as the cone
generated by the divisor classes $\left[  D_{r}\right]  \in A_{n-1}\left(
Y\right)  $ for $r\in\Sigma\left(  1\right)  $.

\begin{proposition}
If $a=\sum_{r\in\Sigma\left(  1\right)  }a_{r}\left[  D_{r}\right]  \in
A_{n-1}^{+}\left(  Y\right)  \otimes\mathbb{R}$, for any $\sigma\in\Sigma$
there is an $m_{\sigma}\in M_{\mathbb{R}}$ such that $\left\langle m_{\sigma
},\hat{r}\right\rangle =-a_{r}$ for all rays $r\subset\sigma$. If
$\left\langle m_{\sigma},\hat{r}\right\rangle \geq-a_{r}$ for all
$r\not \subset \sigma$, then $a$ is
\newsym[$\operatorname*{cpl}\left(  \Sigma\right)  $]{cone of convex classes}{}called
\textbf{convex}. The set $\operatorname*{cpl}\left(  \Sigma\right)  $ of all
convex $a\in A_{n-1}^{+}\left(  Y\right)  \otimes\mathbb{R}$ is a $\left\vert
\Sigma\left(  1\right)  \right\vert -n$ dimensional convex cone.
\end{proposition}

$a$ is in the interior of $\operatorname*{cpl}\left(  \Sigma\right)  $ if and
only if $\left\langle m_{\sigma},\hat{r}\right\rangle >-a_{r}$ for all maximal
dimensional cones $\sigma\in\Sigma$ and all $r\not \subset \sigma$.

\begin{proposition}
\cite[Sec. 3.3.]{CK Mirror Symmetry and Algebraic Geometry} The K\"{a}hler
cone of $Y$ is the interior of $\operatorname*{cpl}\left(  \Sigma\right)  $.
\end{proposition}

\begin{corollary}
The
\index{Mori cone|textbf}%
\textbf{Mori cone }$\overline{NE}\left(  Y\right)  _{\mathbb{R}}$
\newsym[$\overline{NE}\left(  Y\right)  _{\mathbb{R}}$]{Mori cone}{}of
effective $1$-cycles in $A_{1}\left(  Y\right)  \otimes\mathbb{R}\cong
H_{2}\left(  Y,\mathbb{R}\right)  $ is dual to $\operatorname*{cpl}\left(
\Sigma\right)  $.
\end{corollary}

\begin{proposition}
\cite{Reid Decomposition of toric morphisms in Arithmetic and Geometry Vol.
II} $\overline{NE}\left(  Y\right)  _{\mathbb{R}}$ is generated by the
\index{torus orbit closure}%
torus orbit closures $V\left(  \sigma\right)  $ where $\sigma\in\Sigma$ is a
cone of dimension $n-1$.
\end{proposition}

Suppose $\sigma\in\Sigma$ is a cone of dimension $n-1$ generated by
$v_{1},...,v_{n-1}\in N$. The cone $\sigma$ is contained in exactly two $n$
dimensional cones $C_{1}$ and $C_{2}$. There are $v_{n},v_{n+1}\in N$ such
that%
\begin{align*}
C_{1}  &  =\operatorname*{hull}\left(  v_{1},...,v_{n-1},v_{n}\right) \\
C_{2}  &  =\operatorname*{hull}\left(  v_{1},...,v_{n-1},v_{n+1}\right)
\end{align*}
There are relatively prime integers $\lambda_{1},...,\lambda_{n+1}%
\in\mathbb{Z}$ with $\lambda_{n},\lambda_{n+1}>0$ such that $\sum_{i=1}%
^{n+1}\lambda_{i}v_{i}=0$. Denote the relation $\left(  \lambda_{i}\right)  $
by $\lambda_{\sigma}$.

Consider%
\[
\Lambda_{\mathbb{Q}}=\left\{  \left(  \lambda_{v}\right)  \in\mathbb{Q}%
^{\Sigma\left(  1\right)  }\mid%
{\textstyle\sum\nolimits_{i=1}^{n+1}}
\lambda_{i}v_{i}=0\right\}
\]
Applying $\operatorname*{Hom}\nolimits_{\mathbb{Z}}\left(  \mathbb{-}%
,\mathbb{Z}\right)  $ to the sequence%
\[
0\rightarrow M\rightarrow\mathbb{Z}^{\Sigma\left(  1\right)  }\rightarrow
A_{n-1}\left(  X\left(  \Sigma\right)  \right)  \rightarrow0\text{ }%
\]
and tensoring with $\mathbb{Q}$ we get a natural isomorphism%
\[
A_{1}\left(  Y\right)  \otimes\mathbb{Q}\cong\Lambda_{\mathbb{Q}}%
\]
and $V\left(  \sigma\right)  $ is mapped to a multiple $c_{\sigma}%
\lambda_{\sigma}\in\Lambda_{\mathbb{Q}}$ of the relation $\lambda_{\sigma}%
\in\Lambda_{\mathbb{Q}}$ given by $%
{\textstyle\sum\nolimits_{i=1}^{n+1}}
\lambda_{i}v_{i}=0$.

\begin{proposition}
\cite[Sec. 3.3]{CK Mirror Symmetry and Algebraic Geometry} If $\sigma\in
\Sigma\left(  n-1\right)  $ and $v_{1},...,v_{n+1}\in N$ with
\[
\sigma=\operatorname*{hull}\left(  v_{1},...,v_{n-1}\right)
=\operatorname*{hull}\left(  v_{1},...,v_{n-1},v_{n}\right)  \cap
\operatorname*{hull}\left(  v_{1},...,v_{n-1},v_{n+1}\right)
\]
and $\sum_{i=1}^{n+1}\lambda_{i}v_{i}=0$ with $\lambda_{1},...,\lambda
_{n+1}\in\mathbb{Z}$ and $\lambda_{n},\lambda_{n+1}>0$, then there is a
$c_{\sigma}>0$%
\[%
\begin{tabular}
[c]{ccc}%
$A_{1}\left(  Y\right)  \otimes\mathbb{Q}$ & $\overset{\cong}{\rightarrow}$ &
$\Lambda_{\mathbb{Q}}\subset\mathbb{Q}^{\Sigma\left(  1\right)  }$\\
$V\left(  \sigma\right)  $ & $\mapsto$ & $\lambda_{\sigma}=c_{\sigma}%
\cdot\left(  \lambda_{i}\right)  $%
\end{tabular}
\]
A Cartier divisor $D$ is ample if and only if $D.V\left(  \sigma\right)  >0$
for all $\sigma\in\Sigma\left(  n-1\right)  $, so:

A divisor $D=\sum_{r\in\Sigma\left(  1\right)  }a_{r}D_{r}$ is ample if and
only if it is Cartier and $\left(  a_{r}\right)  \cdot\lambda_{\sigma}>0$ for
all $\sigma\in\Sigma\left(  n-1\right)  $.
\end{proposition}

Note that by the sequence%
\[%
\begin{tabular}
[c]{lllllll}%
$0\rightarrow$ & $\operatorname*{Hom}\nolimits_{\mathbb{R}}\left(
A_{n-1}\left(  Y\right)  \otimes\mathbb{R},\mathbb{R}\right)  $ &
$\rightarrow$ & $\operatorname*{Hom}\nolimits_{\mathbb{R}}\left(
\mathbb{R}^{\Sigma\left(  1\right)  },\mathbb{R}\right)  $ & $\overset{\_\circ
A}{\rightarrow}$ & $N_{\mathbb{R}}$ & $\rightarrow0$%
\end{tabular}
\]
the Mori cone of $Y$ is%
\begin{align*}
\overline{NE}\left(  Y\right)  _{\mathbb{R}}  &  =\operatorname*{hull}\left\{
\lambda_{\sigma}\mid\sigma\in\Sigma\text{, }\dim\left(  \sigma\right)
=n-1\right\} \\
&  \subset\ker\left(  A^{t}\right)  ^{\ast}=\operatorname*{Hom}%
\nolimits_{\mathbb{R}}\left(  A_{n-1}\left(  Y\right)  \otimes\mathbb{R}%
,\mathbb{R}\right)  \cong A_{1}\left(  Y\right)  \otimes\mathbb{R}%
\end{align*}

\subsubsection{Toric Fano varieties\label{Sec toric fano varieties}}

Recall
\index{Fano}%
that
\index{toric Fano}%
any
\index{Cohen-Macaulay}%
Cohen-Macaulay variety $Y$ of dimension $n$ has a
\index{dualizing sheaf}%
dualizing sheaf $\hat{\Omega}_{Y}^{n}$ and that this is a line bundle if and
only if $Y$ is
\index{Gorenstein}%
Gorenstein.

\begin{definition}
\label{Def Gorenstein Fano}A complete Gorenstein variety $Y$ is
\index{toric Fano|textbf}%
called
\index{Fano|textbf}%
\textbf{Fano} if the dual of $\hat{\Omega}_{Y}^{n}$ is ample.
\end{definition}

For any toric variety $Y$%
\[
\hat{\Omega}_{Y}^{n}=\mathcal{O}_{Y}\left(  -\sum_{v\in\Sigma\left(  1\right)
}D_{v}\right)
\]
so a toric variety $Y$ is Gorenstein if and only if $-K_{Y}=\sum_{v\in
\Sigma\left(  1\right)  }D_{v}$ is
\index{Cartier}%
Cartier, hence:

\begin{lemma}
If $Y$ is a
\index{complete toric variety}%
complete toric variety $Y$, then it is Fano if and only if $\sum_{v\in
\Sigma\left(  1\right)  }D_{v}$
\index{T-Cartier divisor}%
is Cartier and
\index{ample}%
ample.
\end{lemma}

\begin{definition}
A polytope $\Delta\subset M_{\mathbb{R}}\cong\mathbb{R}^{n}$ of dimension $n$
is called
\index{reflexive|textbf}%
\textbf{reflexive} if $\Delta$ and its
\index{dual polytope}%
dual $\Delta^{\ast}$ are
\index{integral polytope}%
integral
\index{lattice polytope}%
and contain $0$.
\end{definition}

If $\Delta\subset M_{\mathbb{R}}$ is reflexive, then the vertices of
$\Delta^{\ast}$ are in the lattice $M$, hence, for each facet $F$ of $\Delta$
there is an $m_{F}\in M$ with
\[
F=\Delta\cap\left\{  w\in N_{\mathbb{R}}\mid\left\langle m_{F},w\right\rangle
=-1\right\}
\]
and $\Delta$ is cut out by the inequalities $\left\langle m_{F},w\right\rangle
\geq-1$ for all $F$, so for any lattice point in the interior of $\Delta$ we
have $\left\langle m_{F},w\right\rangle >-1$ and $\left\langle m_{F}%
,w\right\rangle \in\mathbb{Z}$ for all facets $F$ of $\Delta$, hence, $0$ is
the unique interior lattice point of $\Delta$.

\begin{lemma}
If $\Delta\subset M_{\mathbb{R}}$ is reflexive, then $0$ is the unique
interior lattice point of $\Delta$.
\end{lemma}

\begin{theorem}
\cite[Sec. 3.5]{CK Mirror Symmetry and Algebraic Geometry}, \cite[Sec.
4.4]{Voisin Mirror Symmetry} The Gorenstein toric Fano varieties
$\mathbb{P}\left(  \Delta\right)  $ of dimension $n$, polarized by
$-K_{\mathbb{P}\left(  \Delta\right)  }$ are in one-to-one correspondence with
the reflexive polytopes $\Delta\subset M_{\mathbb{R}}$, $\operatorname*{rank}%
M=n$. Hence
\index{dual polytope}%
duality of
\index{reflexive}%
reflexive polytopes is an involution of the set of Gorenstein toric Fano varieties.
\end{theorem}

This involution is used by Batyrev in his
\index{mirror construction}%
mirror construction for
\index{hypersurface}%
hypersurfaces in toric varieties.

\begin{figure}
[h]
\begin{center}
\includegraphics[
height=1.8299in,
width=1.8299in
]%
{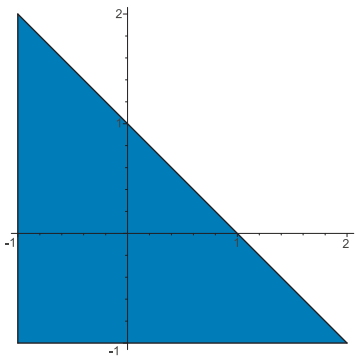}%
\caption{Polytope representing the degree $3$ Veronese of $\mathbb{P}^{2}$}%
\label{FigP2polytope}%
\end{center}
\end{figure}

\begin{figure}
[hh]
\begin{center}
\includegraphics[
height=1.3301in,
width=1.3301in
]%
{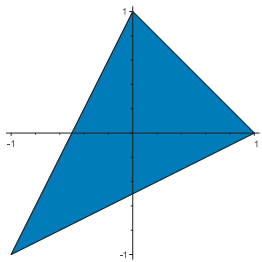}%
\caption{Dual polytope of the degree $3$ Veronese polytope of $\mathbb{P}^{2}%
$}%
\label{FigP2dualpolytope}%
\end{center}
\end{figure}

\begin{example}
The polytope
\[
\Delta=\operatorname*{convexhull}\left(  \left(  2,-1\right)  ,\left(
-1,2\right)  ,\left(  -1,-1\right)  \right)
\]
which is shown in Figure \ref{FigP2polytope}, giving the degree $3$
\index{Veronese}%
Veronese of $\mathbb{P}^{2}$ is reflexive with
\index{dual polytope}%
dual
\[
\Delta^{\ast}=\operatorname*{convexhull}\left(  \left(  1,0\right)  ,\left(
0,1\right)  ,\left(  -1,-1\right)  \right)
\]
shown in Figure \ref{FigP2dualpolytope}.
\end{example}

%


\subsubsection{The automorphism group of a toric
variety\label{Sec automorphism group of a toric variety}}

Suppose $X\left(  \Sigma\right)  $ is a
\index{automorphism}%
complete toric variety
\index{complete toric variety}%
given by a
\index{simplicial}%
simplicial fan $\Sigma$ and
\begin{align*}
X\left(  \Sigma\right)   &  =\left(  \mathbb{C}^{\Sigma\left(  1\right)
}-V\left(  B\left(  \Sigma\right)  \right)  \right)  /G\left(  \Sigma\right)
\\
G\left(  \Sigma\right)   &  =\operatorname*{Hom}\nolimits_{\mathbb{Z}}\left(
A_{n-1}\left(  X\left(  \Sigma\right)  \right)  ,\mathbb{C}^{\ast}\right)
\end{align*}
the homogeneous coordinate representation. The following three possible types
of
\index{automorphism}%
automorphisms of $X\left(  \Sigma\right)  $ are given as automorphisms of
$\mathbb{C}^{\Sigma\left(  1\right)  }-V\left(  B\left(  \Sigma\right)
\right)  $ commuting with the action of $G\left(  \Sigma\right)  $ on
$\mathbb{C}^{\Sigma\left(  1\right)  }-V\left(  B\left(  \Sigma\right)
\right)  $:

\begin{enumerate}
\item By exactness of%
\[
1\rightarrow G\left(  \Sigma\right)  \rightarrow\left(  \mathbb{C}^{\ast
}\right)  ^{\Sigma\left(  1\right)  }\rightarrow T\rightarrow1
\]
the elements of $\left(  \mathbb{C}^{\ast}\right)  ^{\Sigma\left(  1\right)
}$ induce the
\index{automorphism}%
automorphisms of $X\left(  \Sigma\right)  $, which are in the
\index{torus}%
torus $T\subset\operatorname*{Aut}\left(  X\left(  \Sigma\right)  \right)  $.

\item A
\index{root|textbf}%
\textbf{root} of $X\left(  \Sigma\right)  $ is a pair $\left(  y_{v}%
,\prod_{r\in\Sigma\left(  1\right)  }y_{r}^{a_{r}}\right)  $ of a Cox variable
$y_{v}$ and a Cox monomial $\prod_{r\in\Sigma\left(  1\right)  }y_{r}^{a_{r}}%
$, which are not equal, but have the same
\index{Cox degree}%
Cox degree, i.e.,%
\[
\left[  \sum_{r\in\Sigma\left(  1\right)  }a_{r}D_{r}\right]  =\left[
D_{v}\right]  \in A_{n-1}\left(  X\left(  \Sigma\right)  \right)
\]

The Cox monomial $\prod_{r}y_{r}^{a_{r}}$ is not divisible by $y_{v}$, as
otherwise the quotient would be a nontrivial degree $0$ Cox monomial.

Any root $\left(  y_{v},\prod_{r}y_{r}^{a_{r}}\right)  $ induces a
$1$-parameter family of
\index{automorphism}%
automorphisms of $\mathbb{C}^{\Sigma\left(  1\right)  }-V\left(  B\left(
\Sigma\right)  \right)  $ commuting with $G\left(  \Sigma\right)  $
\[%
\begin{tabular}
[c]{cccc}%
$y_{v}\mapsto$ & $y_{v}+\lambda\prod_{s\in\Sigma\left(  1\right)  -\left\{
v\right\}  }y_{s}^{a_{r}}$ &  & \\
$y_{r}\mapsto$ & \multicolumn{1}{l}{$y_{r}$} & for & $r\in\Sigma\left(
1\right)  -\left\{  v\right\}  $%
\end{tabular}
\]
where $y_{v}$ denote coordinates on $\mathbb{C}^{\Sigma\left(  1\right)  }$.

Denote
\newsym[$\operatorname*{Roots}\left(  X\left(  \Sigma\right)  \right)  $]{roots of $X\left(  \Sigma\right)  $}{}by
$\operatorname*{Roots}\left(  X\left(  \Sigma\right)  \right)  $ the set of
roots of $X\left(  \Sigma\right)  $.

\item Any
\index{automorphism}%
automorphism of $N$, which permutes the cones of the fan, gives a permutation
of the rays of $\Sigma$, i.e., of the Cox variables.
\end{enumerate}

\begin{theorem}
\cite{Cox The homogeneous coordinate ring of a toric variety}, \cite[Sec.
3.6]{CK Mirror Symmetry and Algebraic Geometry} If $X\left(  \Sigma\right)  $
is
\index{simplicial}%
simplicial, then torus,
\index{root}%
root and fan
\index{automorphism}%
automorphisms generate $\operatorname*{Aut}\left(  X\left(  \Sigma\right)
\right)  $. Torus and
\index{root}%
root
\index{automorphism}%
automorphisms
\newsym[$\operatorname*{Aut}\left(  X\left(  \Sigma\right)  \right)  $]{automorphism group}{}generate
the connected component of the identity of $\operatorname*{Aut}\left(
X\left(  \Sigma\right)  \right)  $ and%
\begin{equation}
\dim\operatorname*{Aut}\left(  X\left(  \Sigma\right)  \right)  =\dim\left(
T\right)  +\left\vert \operatorname*{Roots}\left(  X\left(  \Sigma\right)
\right)  \right\vert \label{4dimAutoftoric}%
\end{equation}

\end{theorem}

Note that
\[
\left\vert \operatorname*{Roots}\left(  X\left(  \Sigma\right)  \right)
\right\vert =\sum_{r\in\Sigma\left(  1\right)  }\left(  \dim\left(  S_{\left[
D_{r}\right]  }\right)  -1\right)
\]
see also Section \ref{Sec Global sections a cox monomials}.

If $X\left(  \Sigma\right)  $ is
\index{simplicial}%
simplicial and
\index{Gorenstein}%
Gorenstein, then there is a one-to-one correspondence between the lattice
points in the relative interior of the facets (i.e., codimension one faces) of
the polytope
\[
\Delta_{-K_{X\left(  \Sigma\right)  }}=\left\{  m\in M_{\mathbb{R}}%
\mid\left\langle m,\hat{r}\right\rangle \geq-1\forall r\in\Sigma\left(
1\right)  \right\}
\]
and the
\index{root}%
roots of $X\left(  \Sigma\right)  $:

\begin{itemize}
\item If $\rho=\left(  y_{v},\prod_{r\in\Sigma\left(  1\right)  -\left\{
v\right\}  }y_{r}^{a_{r}}\right)  $ is a
\index{root}%
root of $X\left(  \Sigma\right)  $, then $\deg\left(  \frac{\prod_{r\in
\Sigma\left(  1\right)  -\left\{  v\right\}  }y_{r}^{a_{r}}}{y_{v}}\right)
=0$, so $\left(  b_{r}\right)  \in\mathbb{Z}^{\Sigma\left(  1\right)  }%
\ $with
\[
b_{r}=\left\{
\begin{tabular}
[c]{cl}%
$a_{r}$ & if $r\in\Sigma\left(  1\right)  -\left\{  v\right\}  $\\
$-1$ & if $r=v$%
\end{tabular}
\right\}
\]
is in the image of $A$ in%
\[%
\begin{tabular}
[c]{lllllll}%
$0\rightarrow$ & $M$ & $\overset{A}{\rightarrow}$ & $\mathbb{Z}^{\Sigma\left(
1\right)  }$ & $\rightarrow$ & $A_{n-1}\left(  X\left(  \Sigma\right)
\right)  $ & $\rightarrow0$\\
& $m$ & $\mapsto$ & $\left(  \left\langle m,\hat{r}\right\rangle \right)
_{r\in\Sigma\left(  1\right)  }$ &  &  &
\end{tabular}
\]
i.e., there is a unique $m_{\rho}\in M$ such that%
\[
\left\langle m_{\rho},\hat{v}\right\rangle =-1
\]
and%
\[
\left\langle m_{\rho},\hat{r}\right\rangle =a_{r}\text{ for all }r\neq v
\]
By%
\[
\Delta_{-K_{X\left(  \Sigma\right)  }}=\left\{  m\in M_{\mathbb{R}}%
\mid\left\langle m_{\rho},\hat{r}\right\rangle \geq-1\forall r\in\Sigma\left(
1\right)  \right\}
\]
and $\left\langle m_{\rho},\hat{r}\right\rangle =a_{r}\geq0>-1$ for $r\neq v$
and $\left\langle m_{\rho},\hat{v}\right\rangle =-1$ we conclude that
$m_{\rho}$ is in the interior of the facet of $\Delta_{-K_{X\left(
\Sigma\right)  }}$ given by $\left\langle m,\hat{v}\right\rangle =-1$, i.e.,
\[
m_{\rho}\in\operatorname*{int}\left(  \Delta_{-K_{X\left(  \Sigma\right)  }%
}\cap\left\{  \left\langle m,\hat{v}\right\rangle =-1\right\}  \right)  \cap
M
\]

\item If $m_{\rho}\in\operatorname*{int}\left(  \Delta_{-K_{X\left(
\Sigma\right)  }}\cap\left\{  \left\langle m,\hat{v}\right\rangle =-1\right\}
\right)  \cap M$ is a lattice point in the relative interior of a facet of
$\Delta_{-K_{X\left(  \Sigma\right)  }}$, then%
\begin{align*}
\left\langle m_{\rho},\hat{v}\right\rangle  &  =-1\\
\left\langle m_{\rho},\hat{r}\right\rangle  &  >-1\text{ }\forall r\in
\Sigma\left(  1\right)  \text{ with }r\neq v
\end{align*}
so with%
\[
a_{r}=\left\langle m_{\rho},\hat{r}\right\rangle \in\mathbb{Z}_{\geq-1}%
\]
for $r\neq v$ we have%
\[
\frac{\prod_{r\in\Sigma\left(  1\right)  -\left\{  v\right\}  }y_{r}^{a_{r}}%
}{y_{v}}=A\left(  m_{\rho}\right)
\]
i.e.,%
\[
\deg\left(  \frac{\prod_{r\in\Sigma\left(  1\right)  -\left\{  v\right\}
}y_{r}^{a_{r}}}{y_{v}}\right)  =0
\]
so $\rho=\left(  y_{v},\prod_{r\in\Sigma\left(  1\right)  -\left\{  v\right\}
}y_{r}^{a_{r}}\right)  $ is a
\index{root}%
root of $X\left(  \Sigma\right)  $. Summarizing:
\end{itemize}

\begin{proposition}
\cite{AGM The MonomialDivisor Mirror Map} If $X\left(  \Sigma\right)  $ is
\index{simplicial}%
simplicial and Gorenstein, then the
\index{root}%
roots of $X\left(  \Sigma\right)  $ are in one-to-one correspondence with the
lattice points in the relative interior of the facets of $\Delta_{-K_{X\left(
\Sigma\right)  }}\subset M_{\mathbb{R}}$.
\end{proposition}

The polytope $\Delta_{-K_{X\left(  \Sigma\right)  }}$ is not a lattice
polytope in general.

\begin{corollary}
If $X\left(  \Sigma\right)  $ is
\index{simplicial}%
simplicial and
\index{Gorenstein}%
Gorenstein, then%
\[
\dim\left(  \operatorname*{Aut}\left(  X\left(  \Sigma\right)  \right)
\right)  =\dim\left(  T\right)  +\sum_{Q\text{ facet of }\Delta_{-K_{X\left(
\Sigma\right)  }}}\left\vert \operatorname*{int}\nolimits_{M}\left(  Q\right)
\right\vert
\]
with $\operatorname*{int}_{M}$ denoting the set of lattice points in the
relative interior of $Q$.
\end{corollary}

\begin{example}
For $X\left(  \Sigma\right)  =\mathbb{P}^{n}$ the
\index{root}%
roots are the pairs $\left(  x_{i},x_{j}\right)  $ for $i\neq j$ , hence%
\[
\dim\left(  \operatorname*{Aut}\left(  \mathbb{P}^{n}\right)  \right)
=n+\left(  n+1\right)  ^{2}-\left(  n+1\right)  =\left(  n+1\right)  ^{2}-1
\]
i.e., the dimension
\index{PGL}%
of $\operatorname*{PGL}\left(  n+1,\mathbb{C}\right)  $.
\end{example}

\subsubsection{Toric Mori theory\label{Sec toric mori theory}}

Recall that for normal varieties $X$ and $Y$ a proper birational morphism
$f:X\rightarrow Y$ is called
\index{small morphism|textbf}%
\textbf{small} if it is an isomorphism in codimension one. A normal variety
$X$ is
\index{Q-factorial|textbf}%
called $\mathbb{Q}$-factorial if all prime divisors on $X$ are $\mathbb{Q}$-Cartier.

\begin{lemma}
\cite{Reid Decomposition of toric morphisms in Arithmetic and Geometry Vol.
II} A toric variety $Y$ is $\mathbb{Q}$-factorial if and only if $Y$ is simplicial.
\end{lemma}

Recall that for any toric variety $X$ there is a small projective toric
morphism $X^{\prime}\rightarrow X$ such that $X^{\prime}$ is $\mathbb{Q}$-factorial.

Let $X$ and $Y$ be normal varieties and $f:X\rightarrow Y$ be a proper
morphism. A
\index{1-cycle|textbf}%
$1$\textbf{-cycle} of $X/Y$ is a formal sum $\sum a_{i}C_{i}$ with $a_{i}%
\in\mathbb{Z}$ of complete curves $C_{i}$ in the fibers of $f$. Denote by
\begin{align*}
Z_{1}\left(  X/Y\right)   &  =\left\{  1\text{-cycles of }X/Y\right\} \\
Z_{1}\left(  X/Y\right)  _{\mathbb{Q}}  &  =Z_{1}\left(  X/Y\right)
\otimes\mathbb{Q}%
\end{align*}
There is a bilinear pairing%
\[%
\begin{tabular}
[c]{lll}%
$\operatorname*{Pic}\left(  X\right)  \times Z_{1}\left(  X/Y\right)
_{\mathbb{Q}}$ & $\rightarrow$ & $\mathbb{Q}$\\
\multicolumn{1}{c}{$\left(  \mathcal{L},C\right)  $} & $\mapsto$ & $\deg
_{C}\left(  \mathcal{L}\right)  $%
\end{tabular}
\]
Consider two line bundles respectively $1$-cycles numerically equivalent
$\equiv$ if they induce the same linear form on $Z_{1}\left(  X/Y\right)
_{\mathbb{Q}}$ respectively $\operatorname*{Pic}\left(  X\right)  $. So we get
with
\begin{align*}
N^{1}\left(  X/Y\right)   &  =\left(  \operatorname*{Pic}\left(  X\right)
\otimes_{\mathbb{Z}}\mathbb{Q}\right)  /\equiv\\
N_{1}\left(  X/Y\right)   &  =\left(  Z_{1}\left(  X/Y\right)  _{\mathbb{Q}%
}\right)  /\equiv
\end{align*}
the induced perfect pairing
\[
N^{1}\left(  X/Y\right)  \times N_{1}\left(  X/Y\right)  \rightarrow\mathbb{Q}%
\]
Consider the cone of effective $1$-cycles%
\[
NE\left(  X/Y\right)  =\left\{  C\in N_{1}\left(  X/Y\right)  \mid C=\sum
a_{i}C_{i}\text{ with }a_{i}\geq0\right\}
\]

\begin{definition}
A subcone $W$ of a cone $V$ is called
\index{extremal cone|textbf}%
\textbf{extremal}, if for all $u,v\in V$ with $u+v\in W$ it holds $u,v\in W$.
An
\index{extremal ray|textbf}%
\textbf{extremal ray} is an extremal cone of dimension $1$.
\end{definition}

So for a strongly convex cone $V$ a subcone $W$ is extremal if there is a
linear form $l$ such that
\[
W=\left\{  v\in\partial V\mid l\left(  v\right)  =0\right\}
\]

The
\index{relative Picard number|textbf}%
\textbf{relative Picard number}
\newsym[$\rho\left(  X/Y\right)$]{relative Picard number}{}of $X/Y$ is%
\[
\rho\left(  X/Y\right)  =\dim_{\mathbb{Q}}\left(  N^{1}\left(  X/Y\right)
\right)
\]
$D\in N^{1}\left(  X/Y\right)  $
\index{f-nef|textbf}%
is $f$\textbf{-nef} if $D\geq0$ on $NE\left(  X/Y\right)  $.

\begin{theorem}
[Cone Theorem]\label{thm cone theorem}\cite{FuSa Introduction to the toric
Mori theory} If $f:X\rightarrow Y$ is a proper toric morphism, then $NE\left(
X/Y\right)  $ is a convex polyhedral cone. If $f$ is projective, then
$NE\left(  X/Y\right)  $ is strongly convex.
\end{theorem}

\begin{theorem}
[Contraction Theorem]\label{thm contraction theorem}\cite{FuSa Introduction to
the toric Mori theory} Let $f:X\rightarrow Y$ be a projective toric morphism
and let $R$ be an extremal face of $NE\left(  X/Y\right)  $. There is a
projective surjective toric morphism $g:X\rightarrow W$ over $Y$ such that

\begin{itemize}
\item $Z$ is a toric variety which is projective over $Y$,

\item $g$ has connected fibers,

\item if $C$ is a curve in a fiber of $f$, then $\left[  C\right]  \in R$ if
and only if $g\left(  C\right)  $ is a point.
\end{itemize}

If $R$ is an extremal ray and $X$ is $\mathbb{Q}$-factorial, then also $W$ is
$\mathbb{Q}$-factorial and if $g$ is not small, then $\rho\left(  W/Y\right)
=\rho\left(  X/Y\right)  -1$.
\end{theorem}

\begin{theorem}
[Existence of flips]\label{thm existence of flips}\cite{FuSa Introduction to
the toric Mori theory} Suppose $f:X\rightarrow Z$ is a small toric morphism,
$D$ is a torus invariant $\mathbb{Q}$-Cartier divisor on $X$ and $-D$ is
$f$-ample and $r$ is an integer with $rD$ Cartier. Then there is a small
projective toric morphism%
\[
h:X^{+}=\operatorname*{Proj}\nolimits_{Z}\left(
{\textstyle\bigoplus\nolimits_{m\geq0}}
f_{\ast}\mathcal{O}_{X}\left(  m\cdot r\cdot D\right)  \right)  \rightarrow Z
\]
such that the proper transform $D^{+}$ of $D$ on $X^{+}$ is an $h$-ample
$\mathbb{Q}$-Cartier divisor. Then the birational map%
\[%
\begin{tabular}
[c]{lll}%
$X\quad$ & $\longrightarrow$ & $\quad X^{+}$\\
\multicolumn{1}{r}{${\small f}\searrow$} &  & $\swarrow{\small h}$\\
& $Z$ &
\end{tabular}
\]
is called the elementary transformation with respect to $D$. If $X$ is
$\mathbb{Q}$-factorial and $\rho\left(  X/Z\right)  =1$, then $X^{+}$ is
$\mathbb{Q}$-factorial and $\rho\left(  X^{+}/Z\right)  =1$.
\end{theorem}

As the $1$-skeleton of the fan is not changed by elementary transformations we have:

\begin{theorem}
[Termination of flips]\label{thm termination of flips}\cite{FuSa Introduction
to the toric Mori theory} Let%
\[%
\begin{tabular}
[c]{llclcll}%
$X_{0}\quad$ & $\longrightarrow$ & $X_{1}$ & $\longrightarrow$ & $X_{2}$ &
$\longrightarrow$ & $\cdots$\\
\multicolumn{1}{r}{$\searrow$} &  & \multicolumn{1}{l}{$\swarrow\quad\searrow
$} &  & \multicolumn{1}{l}{$\swarrow\quad\searrow$} &  & \\
& $Z_{0}$ & \multicolumn{1}{l}{} & $Z_{1}$ & \multicolumn{1}{l}{} &  &
\end{tabular}
\]
be a sequence of elementary transformations with respect to the divisors
$D_{i}$, $i=0,1,...$ where $D_{i+1}$ is the proper transform of $D_{i}$. Then
this sequence terminates after finitely many steps.
\end{theorem}

\begin{algorithm}
[Toric minimal model program]\cite{Reid Decomposition of toric morphisms in
Arithmetic and Geometry Vol. II} Given a $\mathbb{Q}$-factorial toric variety,
a projective toric morphism $f:X\rightarrow Y$ and a $\mathbb{Q}$-divisor $D$
on $X$ we have two possibilities:

\begin{enumerate}
\item $D$ is $f$-nef, i.e., $D.C\geq0$ for all curves $C$ contracted by $f$.

In this case the process stops, and we call $X$
\index{minimal model|textbf}%
a
\index{relative minimal model|textbf}%
\textbf{relative }$D$\textbf{-minimal model} over $Y$.

\item $D$ is not $f$-nef:

In this case the Cone Theorem \ref{thm cone theorem} gives the existence of an
extremal ray $R$ in $NE\left(  X/Y\right)  _{D<0}$ and the Contraction Theorem
\ref{thm contraction theorem} yields the associated extremal contraction
$g:X\rightarrow W$ over $Y$ and we have

\begin{enumerate}
\item if $\dim W<\dim X$ then the process stops with a Mori fiber space.

\item if $g$ is birational and contracts a divisor, then $\rho\left(
W/Y\right)  =\rho\left(  X/Y\right)  -1$ and $g$ is called a
\index{divisorial contraction|textbf}%
\textbf{divisorial contraction}.

Continue with $W\rightarrow Y$ and the divisor $g_{\ast}D$.

\item if $g$ is small, then by Theorem \ref{thm existence of flips} there an
elementary transformation $h:X\rightarrow X^{+}$ with respect to $D$. The
birational map $h$ is also called a log-flip .

Continue with $X^{+}\rightarrow Y$ and the divisor $h_{\ast}D$.
\end{enumerate}

This process stops, as $\rho\left(  X/Y\right)  $ drops by divisorial
contractions and any sequence of log-flips terminates by Theorem
\ref{thm termination of flips}.
\end{enumerate}
\end{algorithm}

In the standard Mori theory for toric varieties $f$ is birational, $Y$ is
projective and $D=K_{X}$. Note that $K_{X}$ is $f$-nef if and only if
$K_{X}=f^{\ast}K_{Y}$ in the sense of $\mathbb{Q}$-Cartier divisors.

\begin{theorem}
\cite{Fujino Notes on toric varieties from the Mori theoretic viewpoint} Let
$X$ be a $\mathbb{Q}$-factorial toric variety, $Y$ a complete toric variety
and $f:X\rightarrow Y$ a birational toric morphism. If $E$ is a subset of the
exceptional divisors of $f$, then $f$ factors%
\[%
\begin{tabular}
[c]{lll}%
$X\quad$ & $\overset{f}{\longrightarrow}$ & $\quad Y$\\
\multicolumn{1}{r}{${\small g}\searrow$} &  & $\nearrow{\small h}$\\
& $Y^{\prime}$ &
\end{tabular}
\]
such that the birational map $g:X\longrightarrow Y^{\prime}$

\begin{enumerate}
\item contracts all divisors in $E$

\item is a local isomorphism at every generic point of the divisors not in $E$

\item $g^{-1}:Y^{\prime}\rightarrow X$ contracts no divisor

\item $Y^{\prime}$ is projective over $Y$ and $\mathbb{Q}$-factorial.
\end{enumerate}

So if $E$ is the set of $f$-exceptional divisors then $h$ is a small
projective $\mathbb{Q}$-factorialization.
\end{theorem}

\begin{proposition}
\cite{Fujino Notes on toric varieties from the Mori theoretic viewpoint} If
$Y$ is a complete toric variety and $f_{i}:Y_{i}\rightarrow Y$, $i=1,2$ are
small projective $\mathbb{Q}$-factorializations, then there is a finite
composition of elementary transformations $Y_{1}\rightarrow Y_{2}$.
\end{proposition}

\begin{remark}
We illustrate in the following example the contraction process, the
corresponding Mori cones and the linear forms $K_{X_{i}}$:

\begin{itemize}
\item Let $f_{1}:X_{1}\rightarrow X_{0}=\mathbb{P}^{2}$ be the blowup of
$X_{0}$ in a point $p_{1}$ and $E_{1}$ the exceptional. $X_{1}=\mathbb{P}%
\left(  \mathcal{O}_{\mathbb{P}^{1}}\oplus\mathcal{O}_{\mathbb{P}^{1}}\left(
1\right)  \right)  $ has a fibration over $\mathbb{P}^{1}$ by the $0$-curves.
Denote by $H$ a line in $X_{0}$ not meeting $p_{1}$. The map $f_{1}$ has a
toric representation as a map of fans $\Sigma_{1}\rightarrow\Sigma$.%
\[%
{\includegraphics[
trim=0.329150in 0.329670in 0.000000in 0.000000in,
height=1.3344in,
width=3.2854in
]%
{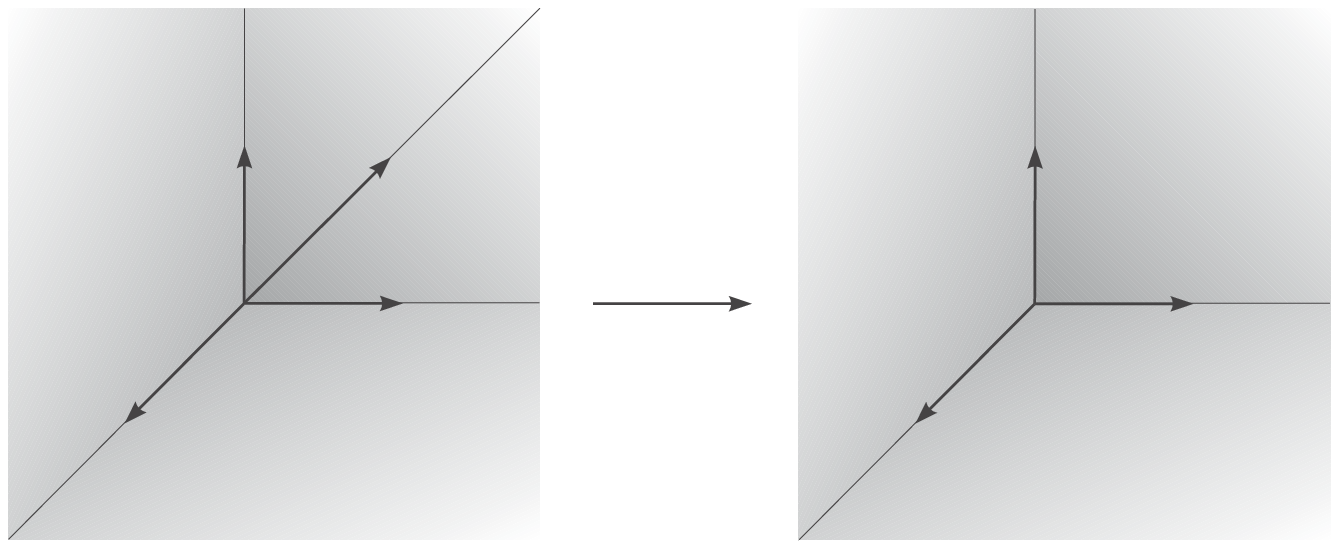}%
}%
\]

\item Further consider the blowup $f_{2}:X_{2}\rightarrow X_{1}$ of $X_{1}$ in
a point $p_{2}\in E_{1}$ with exceptional $E_{2}$ and denote $\widetilde
{E_{1}}$ the strict transform of $E_{1}$. Denote by $L$ the line through
$p_{1}$ such that $\widetilde{L}$ meets $E_{1}$ in $p_{2}$.%
\[%
{\includegraphics[
height=1.7789in,
width=5.2537in
]%
{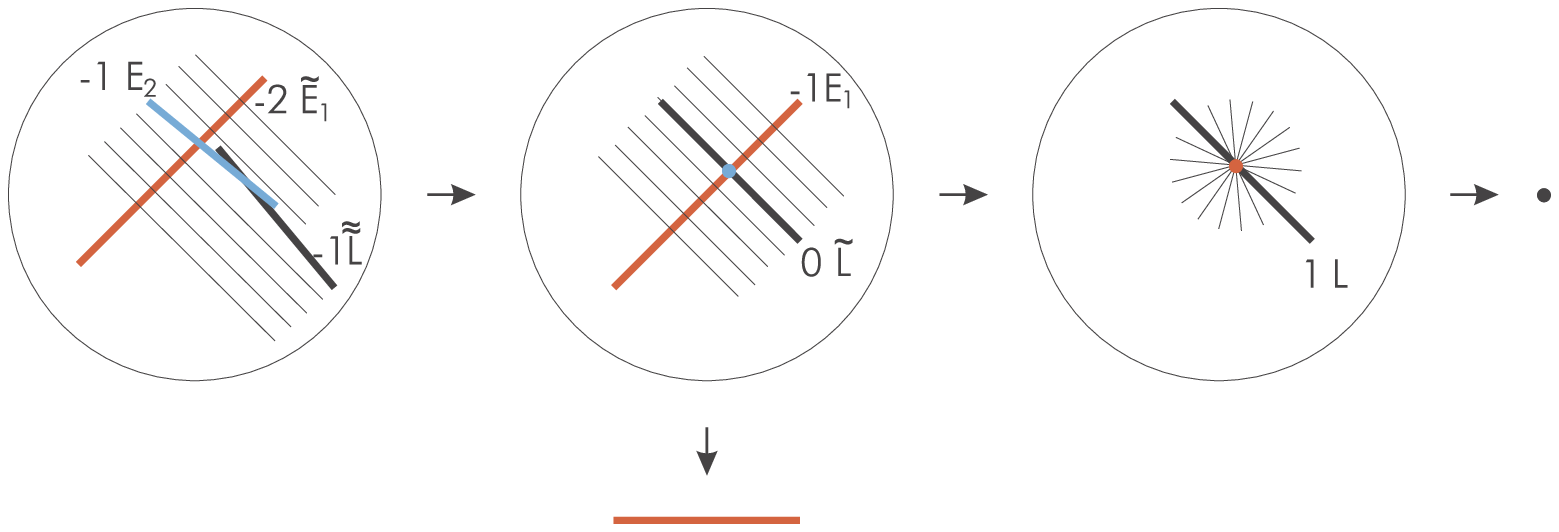}%
}%
\]
Any of the depicted maps corresponds to the contraction of a curve, whose
class $\left[  E\right]  $ generates an extremal ray of $\overline{NE}\left(
X_{i}\right)  _{\mathbb{R}}$. To see this, we calculate the Mori cones:
\end{itemize}

\begin{enumerate}
\item $\overline{NE}\left(  \mathbb{P}^{2}\right)  _{\mathbb{R}}$:

$N_{1}\left(  \mathbb{P}^{2}\right)  =\mathbb{R}$, $\overline{NE}\left(
\mathbb{P}^{2}\right)  _{\mathbb{R}}=\mathbb{R}_{\geq0}$ and $K_{\mathbb{P}%
^{2}}=-3H$%
\[%
{\includegraphics[
height=0.5647in,
width=2.7985in
]%
{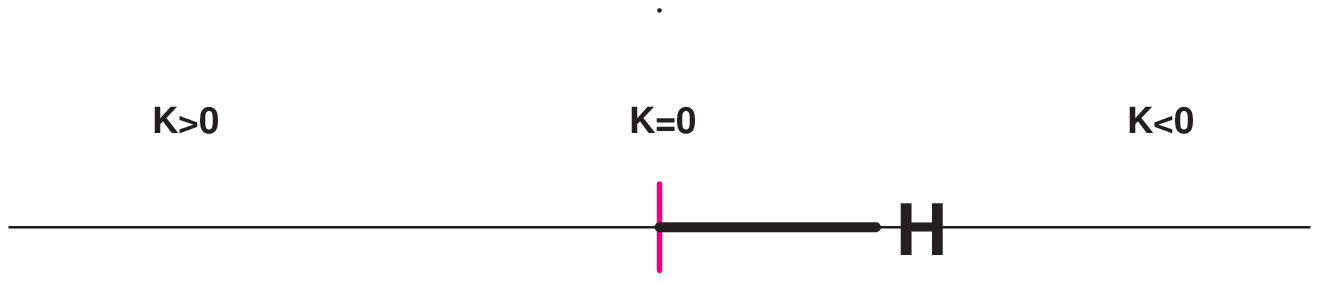}%
}%
\]
From the toric point of view%
\[
\overline{NE}\left(  \mathbb{P}^{2}\right)  _{\mathbb{R}}=\operatorname*{hull}%
\left\{  \left(  1,1,1\right)  \right\}  \subset\left\langle \left(
1,1,1\right)  \right\rangle \subset\mathbb{R}^{\Sigma\left(  1\right)  }%
\]

\item $\overline{NE}\left(  X_{1}\right)  _{\mathbb{R}}$:

\begin{itemize}
\item $\left[  \widetilde{H}\right]  $ and $\left[  E_{1}\right]  $ form a
basis of $N_{1}\left(  X_{1}\right)  \cong\mathbb{R}^{2}$, and we choose
\[%
\begin{tabular}
[c]{ll}%
$\left[  \widetilde{H}\right]  =\left(
\begin{array}
[c]{c}%
1\\
0
\end{array}
\right)  $ & $\left[  E_{1}\right]  =\left(
\begin{array}
[c]{c}%
0\\
1
\end{array}
\right)  $%
\end{tabular}
\]

\item $\widetilde{H},E_{1}\in N_{1}\left(  X_{1}\right)  $ are determined by
$\widetilde{H}^{2}=1$, $\widetilde{H}.E_{1}=0$ and $E_{1}^{2}=-1$.

\item From $E_{1}^{2}=-1$ and $\widetilde{L}^{2}=0$ we know that $\left[
E_{1}\right]  $ generates an extremal ray of $\overline{NE}\left(
X_{1}\right)  _{\mathbb{R}}$ and $\left[  \widetilde{L}\right]  $ is on the
boundary of $\overline{NE\left(  X_{1}\right)  }$, so we can conclude that
they span $\overline{NE}\left(  X_{1}\right)  _{\mathbb{R}}\subset
\mathbb{R}^{2}$ and $\left[  \widetilde{L}\right]  $ generates an extremal
ray. By $\widetilde{L}.\widetilde{H}=L.H-E_{1}.\widetilde{H}=1$,
$\widetilde{L}.E_{1}=1$ we calculate
\[
\left[  \widetilde{L}\right]  =\left(
\begin{array}
[c]{c}%
1\\
-1
\end{array}
\right)
\]

\item $K_{X_{1}}=-3\widetilde{H}+E_{1}$ hence
\[
K_{X_{1}}.\left(  x_{1}\left[  \widetilde{H}\right]  +x_{2}\left[
E_{1}\right]  \right)  =-3x_{1}-x_{2}%
\]%
\[%
{\includegraphics[
height=1.9182in,
width=1.9484in
]%
{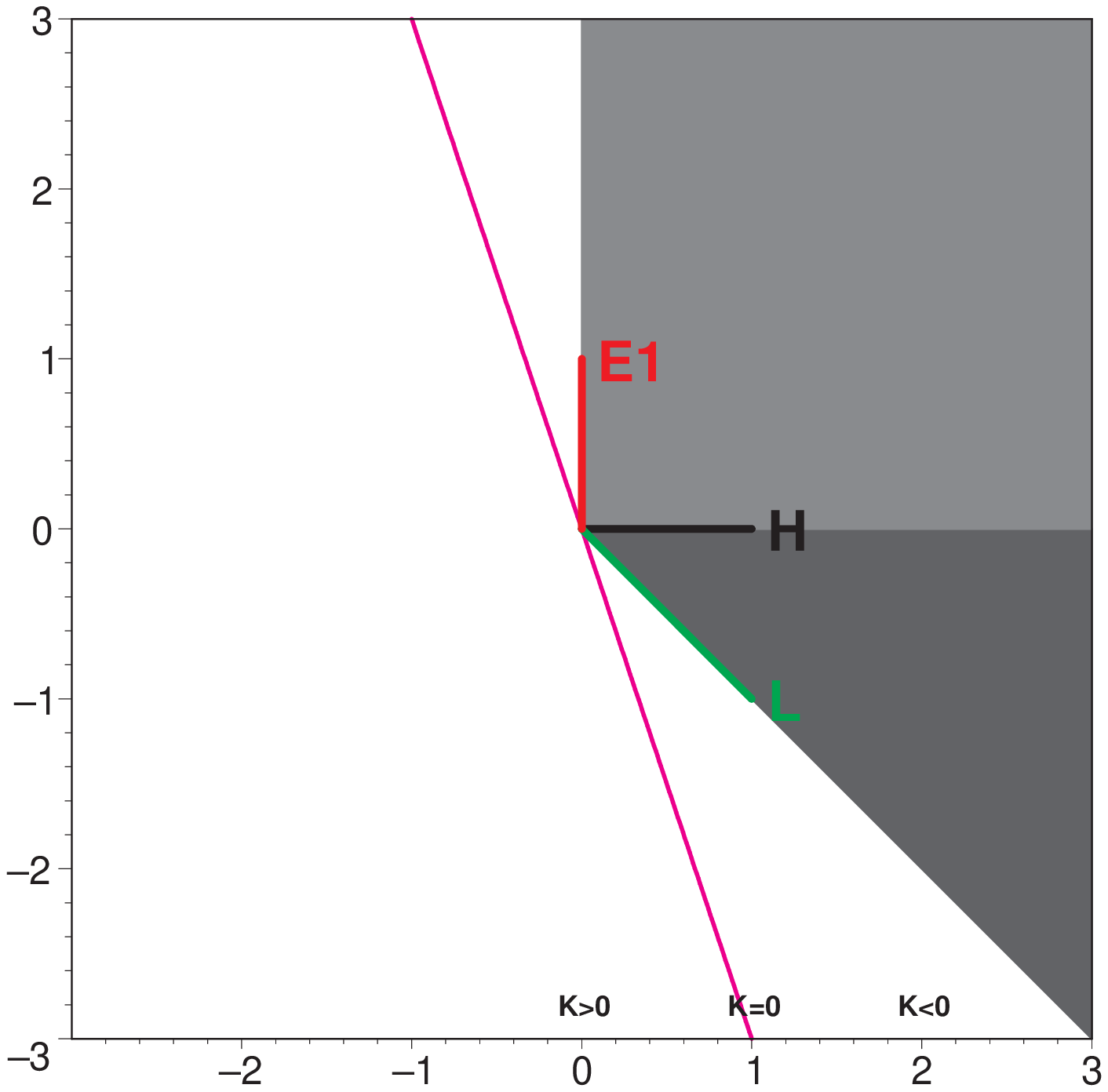}%
}%
\]
We observe that $X_{1}$ is Fano.

From the toric point of view $\overline{NE}\left(  X_{1}\right)  _{\mathbb{R}%
}$ is a cone in the subvectorspace%
\[
\ker\left(

\]

\item Finally we compute the coordinates of the classes $\left[  \widetilde
{C}\right]  $ for all reduced irreducible curves $C\subset\mathbb{P}^{2}$:

We have $C\sim dH$ with $d=C.H$, and $\widetilde{C}=f_{1}^{\ast}C-m_{p_{1}%
}\left(  C\right)  E_{1}$with the multiplicity of $C$ in $p_{1}$, hence
\[
\widetilde{H}.\widetilde{C}=H.C=d\text{ and }E_{1}.\widetilde{C}=f_{1\ast
}\left(  E_{1}\right)  .C-m_{p_{1}}\left(  C\right)  E_{1}^{2}=m_{p_{1}%
}\left(  C\right)
\]
i.e.,
\[
\left[  \widetilde{C}\right]  =\left(
\begin{array}
[c]{c}%
d\\
-m_{p_{1}}\left(  C\right)
\end{array}
\right)
\]
in particular the classes from $\mathbb{P}^{2}$ all lie inside the cone
generated by $\widetilde{H}$ and $L$, i.e., the dark grey area.
\end{itemize}

\item $\overline{NE}\left(  X_{2}\right)  _{\mathbb{R}}$:

\begin{itemize}
\item $\left[  \widetilde{\widetilde{H}}\right]  $, $\left[  \widetilde{E_{1}%
}\right]  $ and $\left[  E_{2}\right]  $ form a basis of $N_{1}\left(
X_{2}\right)  =\mathbb{R}^{3}$, and we choose coordinates
\[%
\begin{tabular}
[c]{lll}%
$\left[  \widetilde{H}\right]  =\left(
\begin{array}
[c]{c}%
1\\
0\\
0
\end{array}
\right)  $ & $\left[  \widetilde{E_{1}}\right]  =\left(
\begin{array}
[c]{c}%
0\\
1\\
0
\end{array}
\right)  $ & $\left[  E_{2}\right]  =\left(
\begin{array}
[c]{c}%
0\\
0\\
1
\end{array}
\right)  $%
\end{tabular}
\]

\item $\widetilde{\widetilde{H}},\widetilde{E_{1}},E_{2}\in N_{1}\left(
X_{2}\right)  $ are determined by $\widetilde{\widetilde{H}}^{2}=1$,
$\widetilde{E_{1}}^{2}=-2$, $E_{2}^{2}=-1$, $\widetilde{\widetilde{H}%
}.\widetilde{E_{1}}=0$, $\widetilde{\widetilde{H}}.E_{2}=0$ and $\widetilde
{E_{1}}.E_{2}=1$.

\item From $E_{2}^{2}=-1$, $\widetilde{\widetilde{L}}^{2}=-1$ and
$\widetilde{E_{1}}^{2}=-2$ we know that $\left[  E_{2}\right]  ,\left[
\widetilde{E_{1}}\right]  $ and $\left[  \widetilde{\widetilde{L}}\right]  $
generate extremal rays of $\overline{NE}\left(  X_{2}\right)  _{\mathbb{R}}$.
From%
\[
\widetilde{\widetilde{L}}.\widetilde{\widetilde{H}}=\left(  f_{2}^{\ast
}\widetilde{L}-E_{2}\right)  .\widetilde{\widetilde{H}}=\widetilde
{L}.\widetilde{H}=1
\]
$\widetilde{\widetilde{L}}.\widetilde{E_{1}}=0$, $\widetilde{\widetilde{L}%
}.E_{2}=1$ we get for the coordinates of $\widetilde{\widetilde{L}}%
=x_{1}\widetilde{\widetilde{H}}+x_{2}\widetilde{E_{1}}+x_{3}E_{2}$ that
$x_{1}=1$, $-2x_{2}+x_{3}=0$, $x_{2}-x_{3}=1$, hence
\[
\left[  \widetilde{\widetilde{L}}\right]  =\left(
\begin{array}
[c]{c}%
1\\
-1\\
-2
\end{array}
\right)
\]

To check that $\overline{NE}\left(  X_{2}\right)  _{\mathbb{R}}$ is really the
cone generated by $\left[  \widetilde{E_{1}}\right]  ,\left[  E_{2}\right]  $
and $\left[  \widetilde{\widetilde{L}}\right]  $, we check that it contains
all classes $\left[  \widetilde{C}\right]  $ for reduced irreducible curves
$C\subset\mathbb{P}^{2}$:

Let $d=C.H$, $m_{p_{1}}\left(  C\right)  $ be the multiplicity of $C$ in
$p_{1}$ and $m_{p_{2}}\left(  \widetilde{C}\right)  $ the tangency of $C$ to
the line $L$.
\begin{align*}
\widetilde{\widetilde{H}}.\widetilde{\widetilde{C}}  &  =\widetilde
{H}.\widetilde{C}=H.C=d\text{ }\\
\widetilde{E_{1}}.\widetilde{\widetilde{C}}  &  =\left(  f_{2}^{\ast}%
E_{1}-E_{2}\right)  .\left(  f_{2}^{\ast}\widetilde{C}-m_{p_{2}}\left(
\widetilde{C}\right)  E_{2}\right)  =E_{1}.\widetilde{C}-m_{p_{2}}\left(
\widetilde{C}\right)  E_{2}^{2}\\
&  =m_{p_{1}}\left(  C\right)  -m_{p_{2}}\left(  \widetilde{C}\right) \\
E_{2}.\widetilde{\widetilde{C}}  &  =m_{p_{2}}\left(  \widetilde{C}\right)
\end{align*}
we get for the coordinates of $\widetilde{\widetilde{C}}=x_{1}\widetilde
{\widetilde{H}}+x_{2}\widetilde{E_{1}}+x_{3}E_{2}$ that $x_{1}=d$,
$-2x_{2}+x_{3}=m_{p_{1}}\left(  C\right)  -m_{p_{2}}\left(  \widetilde
{C}\right)  $, $x_{2}-x_{3}=m_{p_{2}}\left(  \widetilde{C}\right)  $, hence
\[
\left[  \widetilde{C}\right]  =\left(
\begin{array}
[c]{c}%
d\\
-m_{p_{1}}\left(  C\right) \\
-m_{p_{1}}\left(  C\right)  -m_{p_{2}}\left(  \widetilde{C}\right)
\end{array}
\right)
\]
All the classes from $\mathbb{P}^{2}$ lie inside the cone spanned by
\[
\widetilde{\widetilde{L}}=\left(
\begin{array}
[c]{c}%
1\\
-1\\
-2
\end{array}
\right)  ,\left(
\begin{array}
[c]{c}%
1\\
-1\\
0
\end{array}
\right)  ,\left(
\begin{array}
[c]{c}%
1\\
0\\
-2
\end{array}
\right)  ,\widetilde{\widetilde{H}}=\left(
\begin{array}
[c]{c}%
1\\
0\\
0
\end{array}
\right)
\]
i.e., the part of $\overline{NE}\left(  X_{2}\right)  _{\mathbb{R}}$ with non
positive $\widetilde{E_{1}}$ and $E_{2}$ coordinate%
\[%
{\includegraphics[
height=2.3627in,
width=2.0755in
]%
{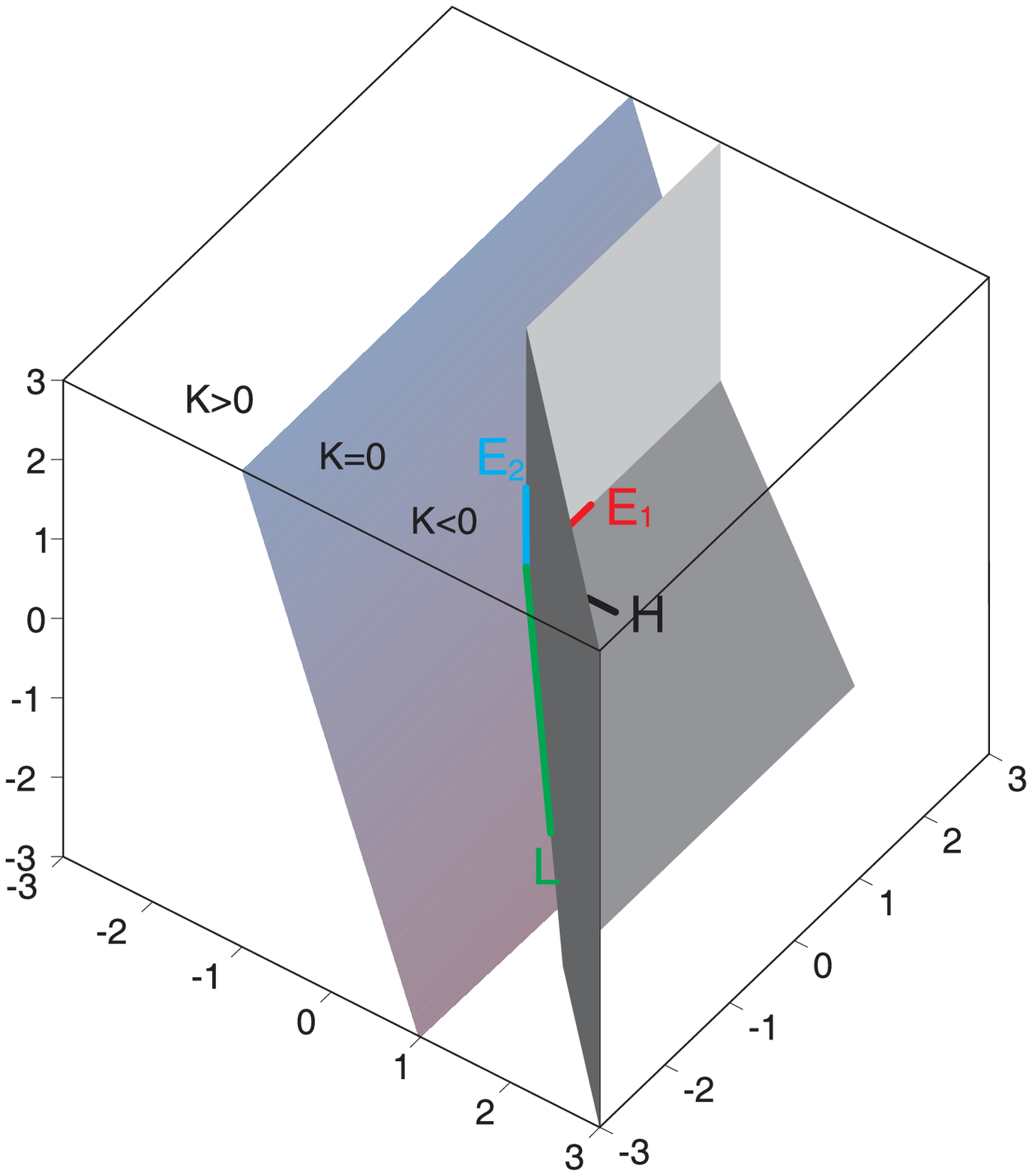}%
}%
\]
We see that $X_{2}$ is no longer Fano.
\end{itemize}
\end{enumerate}
\end{remark}

In the following for a given finite set $\mathcal{R}$ of $1$-dimensional
rational cones we describe the set of all K\"{a}hler cones of projective
simplicial fans $\Sigma$ with $\Sigma\left(  1\right)  \subset\mathcal{R}$.
These K\"{a}hler cones fit together as the maximal cones of a fan. We will
relate this fan to the birational geometry of the toric varieties given by
projective fans $\Sigma$ with $\Sigma\left(  1\right)  \subset\mathcal{R}$.

Given a toric variety $X\left(  \Sigma\right)  $, the Chow group of divisors
$A_{n-1}\left(  X\left(  \Sigma\right)  \right)  $ has the presentation%
\[%
\begin{tabular}
[c]{lllllll}%
$0\rightarrow$ & $M$ & $\rightarrow$ & $\mathbb{Z}^{\Sigma\left(  1\right)  }$
& $\rightarrow$ & $A_{n-1}\left(  X\left(  \Sigma\right)  \right)  $ &
$\rightarrow0$\\
& $m$ & $\mapsto$ & $\left(  \left\langle m,\hat{r}\right\rangle \right)
_{r\in\Sigma\left(  1\right)  }$ &  &  &
\end{tabular}
\]
hence depends only on the $1$-skeleton $\mathcal{R}=\Sigma\left(  1\right)  $
of the fan $\Sigma$. So denote $A_{n-1}\left(  X\left(  \Sigma\right)
\right)  $ by $A_{n-1}\left(  \mathcal{R}\right)  $ and $A_{n-1}\left(
X\left(  \Sigma\right)  \right)  \otimes\mathbb{R}$ by $A_{n-1}\left(
\mathcal{R}\right)  _{\mathbb{R}}$. If $\Sigma$ is a projective simplicial fan
with $\Sigma\left(  1\right)  =\mathcal{R}$ then the K\"{a}hler cone $K\left(
X\left(  \Sigma\right)  \right)  $ canonically lies in $A_{n-1}\left(
\mathcal{R}\right)  _{\mathbb{R}}\cong\mathbb{R}^{\mathcal{R}}/M_{\mathbb{R}}$.

So in the following let $\mathcal{R}$ be a finite set of $1$-dimensional
rational cones in $N_{\mathbb{R}}$, which is the set of rays of a complete fan
in $N_{\mathbb{R}}$. For any projective simplicial fan $\Sigma$ with
$\Sigma\left(  1\right)  =\mathcal{R}$ we get the cone $\operatorname{cpl}%
\left(  \Sigma\right)  $ of dimension $\left\vert \mathcal{R}\right\vert -n$
which is the closure of the K\"{a}hler cone $K\left(  X\left(  \Sigma\right)
\right)  $ and lies in the cone of effective divisor classes $A_{n-1}%
^{+}\left(  \mathcal{R}\right)  _{\mathbb{R}}$. In the same way:

\begin{lemma}
\cite{OP Linear Gale transforms and GelfandKapranovZelevinsky decompositions}
Let $\Sigma$ be a projective fan with $\Sigma\left(  1\right)  \subset
\mathcal{R}$. If $a=\sum_{r\in\Sigma\left(  1\right)  }a_{r}\left[
D_{r}\right]  \in A_{n-1}^{+}\left(  \mathcal{R}\right)  _{\mathbb{R}}$, then
for any $\sigma\in\Sigma$ there is an $m_{\sigma}\in M_{\mathbb{R}}$ such that
$\left\langle m_{\sigma},\hat{r}\right\rangle =-a_{r}$ for all rays
$r\subset\sigma$. If $\left\langle m_{\sigma},\hat{r}\right\rangle \geq-a_{r}$
for all $r\in\mathcal{R}$, $r\not \subset \sigma$, then $a$ is
\index{Sigma-convex|textbf}%
called $\Sigma$-\textbf{convex}. The set of all $\Sigma$-convex $a\in
A_{n-1}^{+}\left(  \mathcal{R}\right)  _{\mathbb{R}}$ is an $\left\vert
\mathcal{R}\right\vert -n$ dimensional convex cone, which we also denote by
$\operatorname{cpl}\left(  \Sigma\right)  $.
\end{lemma}

\begin{theorem}
\cite{OP Linear Gale transforms and GelfandKapranovZelevinsky decompositions}
The set of all $\operatorname{cpl}\left(  \Sigma\right)  $ for projective fans
$\Sigma$ with $\Sigma\left(  1\right)  \subset\mathcal{R}$ form the set of the
$\left\vert \mathcal{R}\right\vert -n$ dimensional cones of a fan with support
$A_{n-1}^{+}\left(  \mathcal{R}\right)  _{\mathbb{R}}$. It is called the
\index{GKZ decomposition|textbf}%
Gelfand-Kapranov-Zelevinsky
\newsym[$GKZ\left(  \mathcal{R}\right)  $]{GKZ decomposition}{}decomposition
$GKZ\left(  \mathcal{R}\right)  $ associated to $\mathcal{R}$.
\end{theorem}

Non simplicial fans $\Sigma$ with $\Sigma\left(  1\right)  \subset\mathcal{R}$
correspond to cones in $GKZ\left(  \mathcal{R}\right)  $ of dimension less
than $\left\vert \mathcal{R}\right\vert -n$. Note that the converse is not true.

\begin{proposition}
\cite{OP Linear Gale transforms and GelfandKapranovZelevinsky decompositions},
\cite{CK Mirror Symmetry and Algebraic Geometry} Two cones $\operatorname{cpl}%
\left(  \Sigma\right)  $ and $\operatorname{cpl}\left(  \Sigma^{\prime
}\right)  $ of dimension $\left\vert \mathcal{R}\right\vert -n$ of $GKZ\left(
\mathcal{R}\right)  $ have a common facet if and only if the toric varieties
$X\left(  \Sigma\right)  $ and $X\left(  \Sigma^{\prime}\right)  $ are related
by a birational extremal contraction.
\end{proposition}

The fan $GKZ\left(  \mathcal{R}\right)  $ may be extended to a complete fan
$\Sigma\left(  \mathcal{R}\right)  $ in $A_{n-1}\left(  \mathcal{R}\right)
_{\mathbb{R}}$:

\begin{definition}
A
\index{marked polytope|textbf}%
\textbf{marked polytope} is a pair $\left(  P,M\right)  $ where $P\subset
\mathbb{R}^{n}$ is a convex polytope and $M\subset P$ is a finite subset with
$\operatorname*{vertices}\left(  P\right)  \subset M$.
\end{definition}

So we may view a marked polytope as just a finite set $M$ of points in
$\mathbb{R}^{n}$, and $P=\operatorname*{convexhull}\left(  M\right)  $.

\begin{definition}
A
\index{polyhedral subdivision|textbf}%
\textbf{polyhedral subdivision} of a marked polytope $\left(  P,M\right)  $ in
$\mathbb{R}^{n}$ is a set of marked polytopes $\left(  P_{i},M_{i}\right)  $
with $\dim\left(  P_{i}\right)  =\dim\left(  P\right)  $ such that%
\[%
{\textstyle\bigcup\nolimits_{i}}
P_{i}=P
\]
and for all $i,j$ the intersection $F=P_{i}\cap P_{j}$ is a face of $P_{i}$
and $P_{j}$ (which may be empty) and%
\[
M_{i}\cap F=M_{j}\cap F
\]
i.e., $M_{i}\cap\operatorname*{convexhull}\left(  M_{j}\right)  =M_{j}%
\cap\operatorname*{convexhull}\left(  M_{i}\right)  $.

A polyhedral subdivision is called
\index{triangulation|textbf}%
\textbf{triangulation}, if all $P_{i}$ are simplices and the $M_{i}$ is the
set of vertices of $P_{i}$.

If $\left\{  \left(  P_{i},M_{i}\right)  \right\}  $ and $\left\{  \left(
P_{j}^{\prime},M_{j}^{\prime}\right)  \right\}  $ are polyhedral subdivisions
of $\left(  P,M\right)  $,
\index{refining triangulations|textbf}%
then $\left\{  \left(  P_{i},M_{i}\right)  \right\}  $ \textbf{refines}
$\left\{  \left(  P_{j}^{\prime},M_{j}^{\prime}\right)  \right\}  $, if for
all $j$
\[
\left\{  \left(  P_{i},M_{i}\right)  \mid P_{i}\subset P_{j}^{\prime}\right\}
\]
is a polyhedral subdivision of $\left(  P_{j}^{\prime},M_{j}^{\prime}\right)
$.
\end{definition}

Hence the set of polyhedral subdivisions of $\left(  P,M\right)  $ form a
poset and the triangulations are the minimal elements.

\begin{lemma}
\cite[Sec. 7.2.]{GKZ Discriminants Resultants and Multidimensional
Determinants} Let $\left(  P,M\right)  $ be a marked polytope in
$\mathbb{R}^{n}$. If $f:M\rightarrow\mathbb{R}^{n}$, i.e., $f\in\mathbb{R}%
^{M}$, is a function let%
\[
G_{f}=\operatorname*{convexhull}\left\{  \left(  x,y\right)  \in\mathbb{R}%
^{n}\oplus\mathbb{R}\mid x\in M\text{, }y\in\mathbb{R}\text{ with }y\leq
f\left(  x\right)  \right\}
\]
Then%
\[%
\begin{tabular}
[c]{l}%
$g_{f}:P\rightarrow\mathbb{R}$\\
$g_{f}\left(  x\right)  =\max\left\{  y\in\mathbb{R\mid}\left(  x,y\right)
\in G_{f}\right\}  $%
\end{tabular}
\ \ \
\]
is a piecewise linear function on $P$. Denote by $P_{i}$ the domains on
linearity of $g_{f}$ and let
\begin{align*}
M_{i}  &  =\left\{  x\in M\cap P_{i}\mid g_{f}\left(  x\right)  =f\left(
x\right)  \right\} \\
&  =\left\{  x\in M\cap P_{i}\mid\left(  x,f\left(  x\right)  \right)
\in\partial G_{f}\right\}
\end{align*}
Then $\left\{  \left(  P_{i},M_{i}\right)  \right\}  $ is a polyhedral
subdivision of $\left(  P,M\right)  $ denoted as $\mathcal{S}\left(  f\right)
$.
\end{lemma}

\begin{definition}
A polyhedral subdivision of a marked polytope $\left(  P,M\right)  $ is
called
\index{coherent|textbf}%
\textbf{coherent}, if it is of the form $S\left(  f\right)  $ for some
$f\in\mathbb{R}^{M}$.
\end{definition}

Let $\mathcal{R}$ be a set of $1$-dimensional rational cones in $N_{\mathbb{R}%
}$, denote by $\hat{r}$, $r\in\mathcal{R}$ the primitive lattice generators of
the elements of $\mathcal{R}$. With
\[
\mathcal{R}^{\prime}=\left\{  \left(  \hat{r},1\right)  \mid r\in
\mathcal{R}\right\}  \cup\left\{  \left(  0,1\right)  \right\}  \subset
N_{\mathbb{R}}\oplus\mathbb{R}%
\]
we have an exact sequence%
\[
0\rightarrow M_{\mathbb{R}}\oplus\mathbb{R}\overset{A^{\prime}}{\rightarrow
}\mathbb{R}^{\mathcal{R}^{\prime}}\overset{\pi}{\rightarrow}A_{n-1}\left(
\mathcal{R}\right)  _{\mathbb{R}}\rightarrow0
\]
with the elements of $\mathcal{R}^{\prime}$ forming the rows of $A^{\prime}$.
Consider the marked polytope $\left(  P,M\right)  =\left(
\operatorname*{convexhull}\left(  \mathcal{R}^{\prime}\right)  ,\mathcal{R}%
^{\prime}\right)  $ in $N_{\mathbb{R}}\oplus\mathbb{R}$.

\begin{definition}
If $\mathcal{S}=\left\{  \left(  P_{i},M_{i}\right)  \right\}  $ is a coherent
polyhedral subdivision of the marked polytope $\left(  P,M\right)  $, then
let
\[
C\left(  \mathcal{S}\right)  =\left\{  \pi\left(  f\right)  \mid
f\in\mathbb{R}^{\mathcal{R}^{\prime}}\text{, }\mathcal{S}\text{ is a
subdivision of }\mathcal{S}\left(  f\right)  \right\}
\]
be the image under $\pi$ of the cone of those functions $f\in\mathbb{R}%
^{\mathcal{R}^{\prime}}$ such that $\mathcal{S}$ is a subdivision of
$\mathcal{S}\left(  f\right)  $.
\end{definition}

\begin{proposition}
\cite[Sec. 7.2.]{GKZ Discriminants Resultants and Multidimensional
Determinants},\newline\cite[Sec. 3.4.]{CK Mirror Symmetry and Algebraic
Geometry} The cones $C\left(  \mathcal{S}\right)  $ form a complete fan in
$A_{n-1}\left(  \mathcal{R}\right)  _{\mathbb{R}}$. This fan is
\newsym[$\Sigma\left(  \mathcal{R}\right)  $]{secondary fan}{}called the
\index{secondary fan|textbf}%
\textbf{secondary fan }$\Sigma\left(  \mathcal{R}\right)  $ of $\mathcal{R}$.
\end{proposition}

\begin{lemma}
\cite[Sec. 3.4.]{CK Mirror Symmetry and Algebraic Geometry} Let $\mathcal{S}%
=\left\{  \left(  P_{i},M_{i}\right)  \right\}  $ be a coherent polyhedral
subdivision of $\left(  P,M\right)  $. Then%
\[
C\left(  \mathcal{S}\right)  =%
{\textstyle\bigcap\nolimits_{i}}
\operatorname*{hull}\left\{  \pi\left(  e_{r^{\prime}}\right)  \mid r^{\prime
}\notin M_{i}\right\}
\]
where $e_{r^{\prime}}\in\mathbb{R}^{\mathcal{R}^{\prime}}$ denotes the
standard basis vector corresponding to $r^{\prime}$.
\end{lemma}

If $C\left(  \mathcal{S}\right)  $ is a coherent polyhedral subdivision
involving $\left(  0,1\right)  $, then the cones over the polytopes of
$\mathcal{S}$ form a complete fan in $N_{\mathbb{R}}\cong N_{\mathbb{R}}%
\times\left\{  1\right\}  $, and any complete fan in $N_{\mathbb{R}}$ arises
this way.

\begin{lemma}
\cite[Sec. 3.4.]{CK Mirror Symmetry and Algebraic Geometry} The secondary fan
$\Sigma\left(  \mathcal{R}\right)  $ contains $GKZ\left(  \mathcal{R}\right)
$ as a subfan. The cones of $GKZ\left(  \mathcal{R}\right)  $ are those
corresponding to coherent polyhedral subdivisions of $\left(
\operatorname*{convexhull}\left(  \mathcal{R}^{\prime}\right)  ,\mathcal{R}%
^{\prime}\right)  $ involving $\left(  0,1\right)  $.
\end{lemma}

\begin{figure}
[h]
\begin{center}
\includegraphics[
height=3.2327in,
width=3.0856in
]%
{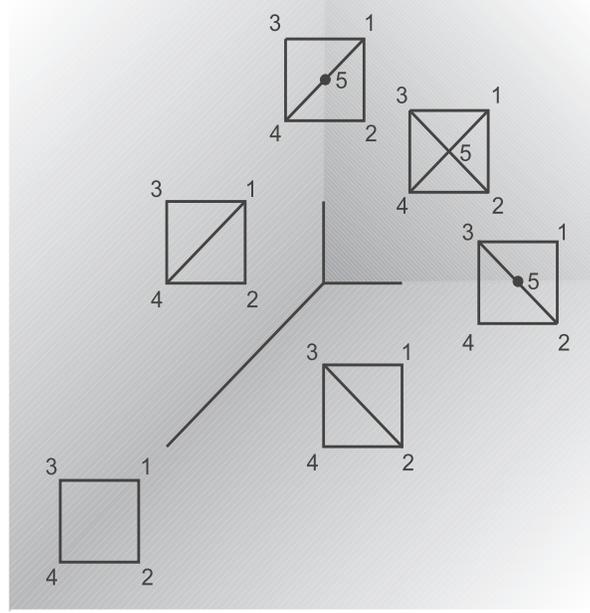}%
\caption{Secondary fan of $\mathbb{P}_{1}\times\mathbb{P}_{1}$}%
\label{Fig secondary fan}%
\end{center}
\end{figure}

\begin{example}
Consider the fan $\Sigma$ given by the rays spanned by
\[
r_{1}=\left(  1,1\right)  ,\text{ }r_{2}=\left(  -1,1\right)  ,\text{ }%
r_{3}=\left(  1,-1\right)  ,\text{ }r_{4}=\left(  -1,-1\right)  \in
N=\mathbb{Z}^{2}%
\]
in $N_{\mathbb{R}}=\mathbb{R}^{2}$ and denote $r_{5}=\left(  0,0\right)  $. We
choose a basis of $\ker\left(  \_\circ A\right)  =\left\langle \left(
1,0,0,1,-2\right)  ,\left(  0,1,1,0,-1\right)  \right\rangle $, so we have the
sequence%
\[%
\begin{tabular}
[c]{cccccccc}%
$0$ & $\rightarrow$ & $M_{\mathbb{R}}$ & $\overset{A}{\rightarrow}$ &
$\mathbb{R}^{\Sigma\left(  1\right)  }$ & $\rightarrow$ & $A_{n-1}\left(
\mathcal{R}\right)  _{\mathbb{R}}$ & $\rightarrow0$\\
&  &  &  &  &  & $\cong$ & \\
$0$ & $\rightarrow$ & $M_{\mathbb{R}}\oplus\mathbb{R}$ & $\overset{A^{\prime}%
}{\rightarrow}$ & $\mathbb{R}^{\Sigma\left(  1\right)  }\oplus\mathbb{R}$ &
$\overset{B}{\rightarrow}$ & $\mathbb{R}^{2}$ & $\rightarrow0$%
\end{tabular}
\]
with%
\[%
\begin{tabular}
[c]{lll}%
$A^{\prime}=\left(
\begin{array}
[c]{ccc}%
1 & 1 & 1\\
-1 & 1 & 1\\
1 & -1 & 1\\
-1 & 1 & 1\\
0 & 0 & 1
\end{array}
\right)  $ &  & $B=\left(
\begin{tabular}
[c]{lllll}%
$1$ & $1$ & $0$ & $1$ & $-2$\\
$0$ & $0$ & $1$ & $0$ & $-2$%
\end{tabular}
\right)  $%
\end{tabular}
\]
Hence considering the secondary fan as a subfan of $\mathbb{R}^{2}\cong
A_{n-1}\left(  \mathcal{R}\right)  _{\mathbb{R}}$ the cones corresponding to
triangulations are%
\begin{align*}
C%
\raisebox{-0.1531in}{\includegraphics[
height=0.3745in,
width=0.3304in
]%
{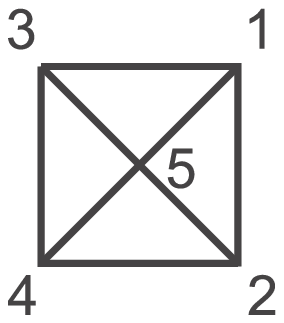}%
}%
&  =\operatorname*{hull}\left(  \left(
\begin{array}
[c]{c}%
1\\
0
\end{array}
\right)  ,\left(
\begin{array}
[c]{c}%
0\\
1
\end{array}
\right)  \right) \\
C%
\raisebox{-0.1531in}{\includegraphics[
height=0.3736in,
width=0.3338in
]%
{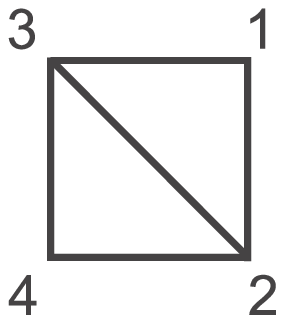}%
}%
&  =\operatorname*{hull}\left(  \left(
\begin{array}
[c]{c}%
1\\
0
\end{array}
\right)  ,\left(
\begin{array}
[c]{c}%
-2\\
-2
\end{array}
\right)  \right) \\
C%
\raisebox{-0.1323in}{\includegraphics[
height=0.3745in,
width=0.3338in
]%
{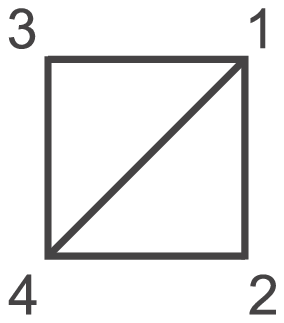}%
}%
&  =\operatorname*{hull}\left(  \left(
\begin{array}
[c]{c}%
0\\
1
\end{array}
\right)  ,\left(
\begin{array}
[c]{c}%
-2\\
-2
\end{array}
\right)  \right)
\end{align*}
The secondary fan is shown in Figure \ref{Fig secondary fan}.
\end{example}

%

\begin{algorithm}
\label{Alg secondary fan GKZ}If $\mathcal{R}=\left\{  r_{1},...,r_{s}\right\}
\subset N_{\mathbb{R}}=\mathbb{Z}^{n}\otimes\mathbb{R}$ is a set of
$1$-dimensional rational cones, which are the rays of a projective fan, the
following algorithm computes the secondary fan $\Sigma\left(  \mathcal{R}%
\right)  $ of $\mathcal{R}$:

\begin{enumerate}
\item Let $\mathcal{\hat{R}}=\left\{  \hat{r}_{1},...,\hat{r}_{s},0\right\}  $
be the set of minimal lattice generators of the rays $r_{i}$ together with $0$.

\item Let $\Sigma\left(  \mathcal{R}\right)  =\{\}$.

\item Choose a random $f\in\mathbb{Z}_{\geq0}^{\mathcal{\hat{R}}}$. Let%
\[
G_{f}=\operatorname*{convexhull}\left\{  \left(  \hat{r}_{j},f_{\hat{r}_{j}%
}\right)  ,\left(  \hat{r}_{j},-1\right)  \mid j=1,...,s\right\}  \subset
N_{\mathbb{R}}\oplus\mathbb{R}%
\]
and compute the set $\mathcal{S}^{\prime}$ of all faces of $G_{f}$ which do
not involve one of the vertices $\left(  \hat{r}_{j},-1\right)  $.

\item Compute the set $\mathcal{S}$ of projections of the faces of
$\mathcal{S}^{\prime}$ under $N_{\mathbb{R}}\oplus\mathbb{R\rightarrow
}N_{\mathbb{R}}$.

\item Let%
\[
\mathcal{R}^{\prime}=\left\{  \left(  \hat{r},1\right)  \mid r\in
\mathcal{R}\right\}  \cup\left\{  \left(  0,1\right)  \right\}  \subset
N_{\mathbb{R}}\oplus\mathbb{R}%
\]
and $A^{\prime}$ be the linear map with rows $\left(  \hat{r},1\right)
,\left(  0,1\right)  $ and $\pi$ the map
\[
0\rightarrow M_{\mathbb{R}}\oplus\mathbb{R}\overset{A^{\prime}}{\rightarrow
}\mathbb{R}^{\mathcal{R}^{\prime}}\overset{\pi}{\rightarrow}A_{n-1}\left(
\mathcal{R}\right)  _{\mathbb{R}}\rightarrow0
\]
Compute%
\[
C\left(  \mathcal{S}\right)  =%
{\textstyle\bigcap\nolimits_{F\in\mathcal{S}^{\prime}}}
\operatorname*{hull}\left\{  \pi\left(  e_{r^{\prime}}\right)  \mid r^{\prime
}\notin F\right\}
\]

\item If $\dim\left(  C\left(  \mathcal{S}\right)  \right)  <s-n$ then we
found a non maximal cone of the secondary fan (which we may remember), and we
go back to $3$.

Otherwise we set $\Sigma\left(  \mathcal{R}\right)  =\Sigma\left(
\mathcal{R}\right)  \cup\left\{  C\left(  \mathcal{S}\right)  \right\}  $.

If $\Sigma\left(  \mathcal{R}\right)  $ is not complete, then we go back to
$3$.

In order to get the fans associated to the cones in the $GKZ$ decomposition,
we may also remember for each cone $C\left(  \mathcal{S}\right)  $ the
corresponding triangulation $\mathcal{S}$.
\end{enumerate}
\end{algorithm}

\begin{remark}
The Maple
\index{tropical mirror}%
package \textsf{tropicalmirror} (see also Section
\ref{Sec tropicalmirror implementation}) provides an implementation of this
algorithm. Given a set $\mathcal{R}=\left\{  r_{1},...,r_{s}\right\}  $ of
lattice vectors in $\mathbb{Z}^{n}$, which are the primitive lattice
generators of the rays of a projective fan, the function
$\mathsf{Triangulations}$ takes $\mathcal{R}$ as an argument and computes all
triangulations of the marked polytope $\left(  \operatorname*{convexhull}%
\left(  \mathcal{R}\right)  ,\mathcal{R}\right)  $. Let $B$ be a matrix such
that the sequence%
\[
0\rightarrow\mathbb{R}^{n+1}\overset{A^{\prime}}{\rightarrow}\mathbb{R}%
^{s+1}\overset{B}{\rightarrow}\mathbb{R}^{s-n}\rightarrow0
\]
with%
\[
A^{\prime}=\left(
\begin{array}
[c]{cc}%
r_{1} & 1\\
\vdots & \vdots\\
r_{s} & 1\\
0 & 1
\end{array}
\right)
\]
is exact, so choosing an isomorphism $A_{n-1}\left(  \mathcal{R}\right)
_{\mathbb{R}}\cong\mathbb{R}^{s-n}$. The function $\mathsf{SecondaryFan}$
takes as argument the list $\left(  r_{1},...,r_{s}\right)  $ and $B$ and
computes the secondary fan as a fan in $\mathbb{R}^{s-n}$. In the same way the
function $\mathsf{GKZFan}$ takes the argument $\left(  r_{1},...,r_{s}\right)
$ and $B$ and computes the $GKZ$ decomposition of $\mathcal{R}$.
\end{remark}

\subsection{Gr\"{o}bner basics\label{Sec Groebner basics}}

To compute combinatorial objects in tropical geometry, we will use
\index{Gr\"{o}bner bases}%
Gr\"{o}bner basis techniques, so we recall some notation and make some remarks
about their implementation in the Macaulay $2$ library
\index{mora.m2}%
\textsf{mora.m2} which is part of the computer algebra implementation of the
\index{mirror construction}%
mirror construction given here.

\subsubsection{Semigroup orderings\label{Sec monomial orderings}}

\begin{definition}
[Semigroup ordering]A
\index{semigroup ordering|textbf}%
\textbf{semigroup ordering}
\index{monomial ordering|textbf}%
(monomial ordering) on the semigroup of monomials in the variables
$x_{1},...,x_{n}$ is an ordering $>$ of the
\newsym[$>$]{semigroup ordering}{}monomials with the following properties

\begin{enumerate}
\item $>$ is a
\index{total ordering}%
total ordering

\item $>$ is compatible with multiplication, i.e., $x^{\alpha}>x^{\beta
}\Rightarrow x^{\alpha}x^{\gamma}>x^{\beta}x^{\gamma}$.
\end{enumerate}
\end{definition}

\begin{definition}
[global ordering]A
\index{global ordering|textbf}%
\textbf{global ordering} $>$ is a semigroup ordering with the following
equivalent properties

\begin{enumerate}
\item $x_{i}>1\forall i$

\item $x^{\alpha}>1$ for all $\alpha\neq0$

\item $>$ is a
\index{well ordering}%
well ordering.

\item $\alpha\geq\beta$ and $\alpha\neq\beta$ $\Rightarrow x^{\alpha}%
>x^{\beta}$
\end{enumerate}
\end{definition}

\begin{definition}
[local ordering]A
\index{local ordering|textbf}%
\textbf{local ordering} $>$ is a semigroup ordering with%
\[
x_{i}<1\text{ }\forall i
\]

\end{definition}

Local orderings are not
\index{well ordering}%
well orderings. This leads to problems with the termination of normal form algorithms.

\begin{remark}
In one variable all global (resp. local) orderings are equivalent.
\end{remark}

\begin{definition}
[weighted degree ordering]A monomial ordering $>$ is called a
\index{weighted degree ordering|textbf}%
\textbf{weighted degree ordering} if there is some $w\in\mathbb{R}^{n}$ with
non zero entries such that%
\[
w\alpha>w\beta\Rightarrow x^{\alpha}>x^{\beta}%
\]

\end{definition}

\begin{example}
If $>$ is any
\index{monomial ordering}%
monomial ordering and $w\in\mathbb{R}^{n}$, then $>_{w}$ given by%
\[
x^{\alpha}>_{w}x^{\beta}\Leftrightarrow w\alpha>w\beta\text{ or }\left(
w\alpha=w\beta\text{ and }x^{\alpha}>x^{\beta}\right)
\]
is a
\index{monomial ordering}%
monomial ordering. It is a weighted degree ordering, it is
\index{global ordering}%
global if $w_{i}>0$ for all $i$ and it is
\index{local ordering}%
local if $w_{i}<0$ for all $i$.
\end{example}

\begin{proposition}
\cite[Sec. 1.2]{GP A Singular Introduction to Commutative Algebra}%
\label{1anymonordbyweightord} Given any finite set of monomials $M$ and any
\index{semigroup ordering}%
semigroup ordering $>$, there is some $w\in\mathbb{Z}^{n}$ such that
\[
x^{\alpha}>x^{\beta}\Leftrightarrow\sum w_{i}\alpha_{i}>\sum w_{i}\beta_{i}%
\]
for all $x^{\alpha},x^{\beta}\in M$.

\noindent$w$ can be choosen such that $w_{i}>0$ if $x_{i}>1$ and $w_{i}<0$ if
$x_{i}<1$.

\noindent$w$ is called a
\index{weight vector|textbf}%
\textbf{weight vector} inducing $>$ on $M$.
\end{proposition}

\begin{example}
The following orderings are
\index{semigroup ordering}%
semigroup orderings:

\begin{itemize}
\item \textbf{lexicographical}
\index{lexicographic ordering|textbf}%
ordering \newsym[$lp$]{lexicographical ordering}{}%
\index{lp|textbf}%
$\mathsf{lp}$:

$x^{\alpha}<x^{\beta}\Leftrightarrow\exists1\leq i\leq n:\alpha_{1}=\beta
_{1},...,\alpha_{i-1}=\beta_{i-1},\alpha_{i}<\beta_{i}$.

\item \textbf{reverse lexicographical}
\index{reverse lexicographic ordering|textbf}%
ordering \newsym[$rp$]{reverse lexicographical ordering}{}%
\index{rp|textbf}%
$\mathsf{rp}$:

$x^{\alpha}<x^{\beta}\Leftrightarrow\exists1\leq i\leq n:\alpha_{n}=\beta
_{n},...,\alpha_{i+1}=\beta_{i+1},\alpha_{i}>\beta_{i}$.

\item \textbf{weighted reverse lexicographical}
\index{weighted reverse lexicographic ordering|textbf}%
ordering \newsym[$lp$]{weighted reverse lexicographical ordering}{}%
\index{wp|textbf}%
$\mathsf{wp}\left(  w\right)  $ for $w\in\mathbb{R}^{n}$:

$x^{\alpha}<x^{\beta}\Leftrightarrow\sum w_{i}\alpha_{i}<\sum w_{i}\beta_{i}$
or $\sum w_{i}\alpha_{i}=\sum w_{i}\beta_{i}$ and $\exists1\leq i\leq
n:\alpha_{n}=\beta_{n},...,\alpha_{i+1}=\beta_{i+1},\alpha_{i}>\beta_{i}$.

\item \textbf{weighted lexicographical}
\index{weighted lexicographic ordering|textbf}%
ordering \newsym[$Wp$]{weighted lexicographical ordering}{}%
\index{Wp|textbf}%
$\mathsf{Wp}\left(  w\right)  $ for $w\in\mathbb{R}^{n}$:

$x^{\alpha}<x^{\beta}\Leftrightarrow\sum w_{i}\alpha_{i}<\sum w_{i}\beta_{i}$
or $\sum w_{i}\alpha_{i}=\sum w_{i}\beta_{i}$ and $\exists1\leq i\leq
n:\alpha_{1}=\beta_{1},...,\alpha_{i-1}=\beta_{i-1},\alpha_{i}<\beta_{i}$.

\item \textbf{degree reverse lexicographical}
\index{degree reverse lexicographic ordering|textbf}%
ordering \newsym[$rp$]{degree reverse lexicographical ordering}{}%
\index{dp|textbf}%
$\mathsf{dp}=\mathsf{wp}\left(  1,...,1\right)  $.

\item \textbf{negative lexicographical}
\index{negative lexicographic ordering|textbf}%
ordering \newsym[$ls$]{negative lexicographical ordering}{}%
\index{ls|textbf}%
$\mathsf{ls}$:

$x^{\alpha}<x^{\beta}\Leftrightarrow\exists1\leq i\leq n:\alpha_{1}=\beta
_{1},...,\alpha_{i-1}=\beta_{i-1},\alpha_{i}>\beta_{i}$.

\item \textbf{matrix
\index{matrix ordering|textbf}%
ordering} associated to
\[
A=\left(
\begin{array}
[c]{c}%
a_{1}\\
\vdots\\
a_{m}%
\end{array}
\right)  \in Mat\left(  m\times n,\mathbb{R}\right)
\]
with $\operatorname*{rank}\left(  A\right)  =n$, is given by:

$x^{\alpha}<x^{\beta}\Leftrightarrow\exists1\leq i\leq m:a_{1}\alpha
=a_{1}\beta,...,a_{i-1}\alpha=a_{i-1}\beta,a_{i}\alpha<a_{i}\beta
\Leftrightarrow A\alpha<_{lex}A\beta$.

Note:

\item $\mathsf{lp}=\mathsf{Wp}\left(  0\right)  $ and $\mathsf{rp}%
=\mathsf{wp}\left(  0\right)  $.

\item If all weights are non negative, then $\mathsf{wp}$ and $\mathsf{Wp}$
are
\index{global ordering}%
global orderings.

\item $\mathsf{ls}$ is a local ordering.
\end{itemize}
\end{example}

\begin{remark}
$\mathsf{ws}\left(  w\right)  =\mathsf{wp}\left(  -w\right)  $
\index{ws|textbf}%
is \newsym[$ws$]{local weighted reverse lexicographical ordering}{}denoted as
\index{local weighted reverse lexicographical|textbf}%
\textbf{local weighted reverse lexicographical} ordering, $\mathsf{Ws}\left(
w\right)  =\mathsf{Wp}\left(  -w\right)  $
\index{Ws|textbf}%
is denoted \newsym[$Ws$]{local weighted lexicographical ordering}{}as
\index{local weighted lexicographical|textbf}%
\textbf{local weighted lexicographical} ordering.

The \textbf{local degree reverse lexicographical}
\index{local degree reverse lexicographic ordering|textbf}%
ordering is \newsym[$ds$]{local degree reverse lexicographical ordering}{}%
\index{ds|textbf}%
$\mathsf{ds}=\mathsf{ws}\left(  1,...,1\right)  $.
\end{remark}

\begin{example}
On a finite set of monomials the ordering $\mathsf{lp}$ can be represented by
the
\index{weight vector}%
weight vector $w=\left(  v^{n-1},...,v,1\right)  $ if all monomials are
contained in a cube of side length $\leq v$.
\end{example}

\begin{remark}
The above
\index{monomial ordering}%
monomial orderings are implemented in the \linebreak Macaulay $2$ package $%
\index{mora.m2}%
$\textsf{mora.m2}. They are selected by the value of the global method
$\mathsf{monord}$, which can be given the values $\mathsf{lp}$, $\mathsf{dp}$,
$\mathsf{wp}$, $\mathsf{ls}$, $\mathsf{Ws}$, $\mathsf{ws}$, $\mathsf{Wp}$ and
$\mathsf{Mat}$ for
\index{matrix ordering}%
matrix orderings. The
\index{weight vector}%
weight vector, if needed, is represented by the global list $\mathsf{ww}$ and
the matrix inducing above
\index{matrix ordering}%
matrix ordering by the global Macaulay~$2$ type matrix $\mathsf{mm}$.
\end{remark}

\begin{remark}
Any
\index{matrix ordering}%
matrix ordering can be represented by a matrix in $Gl\left(  n,\mathbb{R}%
\right)  $. Note that one can add multiples of $a_{i}$ to any lower row
$a_{j}$ with $j>i$ without changing the monomial order.
\end{remark}

\begin{example}
The
\index{weight ordering}%
weight ordering $Wp\left(  w_{1},...,w_{n}\right)  $ can be represented by the
matrix ordering given by%
\[
\left(
\begin{array}
[c]{cccc}%
w_{1} & \cdots & \cdots & w_{n}\\
1 &  &  & \\
& \ddots &  & \\
&  & \ddots & \\
&  &  & 1
\end{array}
\right)
\]
If $w_{j}\neq0$ and $w_{j+1}=0,...,w_{n}=0$, then this ordering is equivalent
to%
\[
\left(
\begin{array}
[c]{ccccccc}%
w_{1} & \cdots & w_{j-1} & w_{j} & 0 & \cdots & 0\\
1 &  &  &  &  &  & \\
& \ddots &  &  &  &  & \\
&  & 1 &  &  &  & \\
w_{1} &  & w_{j-1} & 0 & 0 & \cdots & 0\\
&  &  &  & 1 &  & \\
&  &  &  &  & \ddots & \\
&  &  &  &  &  & 1
\end{array}
\right)  \text{ }%
\]
hence to%
\[
\left(
\begin{array}
[c]{ccccccc}%
w_{1} & \cdots & w_{j-1} & w_{j} & 0 & \cdots & 0\\
1 &  &  &  &  &  & \\
& \ddots &  &  &  &  & \\
&  & 1 &  &  &  & \\
0 &  & 0 & 0 & 0 & \cdots & 0\\
&  &  &  & 1 &  & \\
&  &  &  &  & \ddots & \\
&  &  &  &  &  & 1
\end{array}
\right)
\]
i.e., to the matrix ordering given by%
\[
\left(
\begin{array}
[c]{ccccccc}%
w_{1} & \cdots & w_{j-1} & w_{j} & 0 & \cdots & 0\\
1 &  &  &  &  &  & \\
& \ddots &  &  &  &  & \\
&  & 1 &  &  &  & \\
&  &  &  & 1 &  & \\
&  &  &  &  & \ddots & \\
&  &  &  &  &  & 1
\end{array}
\right)
\]

\end{example}

\begin{proposition}
\cite{RS Subalgebra Bases} Every
\index{semigroup ordering}%
semigroup ordering is representable by a
\index{matrix ordering}%
matrix ordering.
\end{proposition}

\begin{definition}
Let $>$ be a
\index{monomial ordering}%
monomial ordering on the monomials of $K\left[  x_{1},...,x_{n}\right]  $. For
$f\in K\left[  x_{1},...,x_{n}\right]  $, denote by $L\left(  f\right)  $ the
\index{lead monomial|textbf}%
\textbf{lead monomial}, i.e.,
\newsym[$L\left(  f\right)  $]{lead monomial of $f$ with respect to fixed monomial ordering}{}the
largest monomial with respect to $>$ appearing in $f$, by $LC\left(  f\right)
$
\newsym[$LC\left(  f\right)  $]{lead coefficient of $f$ with respect to fixed monomial ordering}{}the
\index{lead coefficient|textbf}%
\textbf{lead coefficient}, i.e., the coefficient of $L\left(  f\right)  $ in
$f$, and by $LT\left(  f\right)  =LC\left(  f\right)  L\left(  f\right)  $
\newsym[$LT\left(  f\right)  $]{lead term of $f$ with respect to fixed monomial ordering}{}the
\index{lead term|textbf}%
\textbf{lead term} of $f$.
\end{definition}

\subsubsection{Localizations}

Let $K$ be a field.

\begin{remark}
The rings
\[
K\left[  x_{1},...,x_{n}\right]  _{\left\langle x_{1},...,x_{n}\right\rangle
}\subset K\left\{  \left\{  x_{1},...,x_{n}\right\}  \right\}  \subset
K\left[  \left[  x_{1},...,x_{n}\right]  \right]
\]
correspond to looking at increasingly smaller neighborhoods of the origin:

\begin{enumerate}
\item Elements of $K\left[  x_{1},...,x_{n}\right]  _{\left\langle
x_{1},...,x_{n}\right\rangle }$ are defined in the complement of an algebraic
set, i.e., in a Zariski open neighborhood of the origin, e.g., $\frac{f}{g}$
is defined in the complement of $V\left(  g\right)  $.

\item Elements
\newsym[$K\left\{\left\{x_{1},...,x_{n}\right\}\right\}$]{convergent power series ring}{}of
$K\left\{  \left\{  x_{1},...,x_{n}\right\}  \right\}  $ are defined in a
neighborhood of the origin in the analytic topology, which can be much
smaller, e.g., the geometric series $\sum_{k=0}^{\infty}x^{k}$ is defined for
$\left\vert x\right\vert <1$.

\item Elements of $K\left[  \left[  x_{1},...,x_{n}\right]  \right]  $ are
defined just at the origin in general.
\end{enumerate}

Nevertheless, they all share the property of being local rings (for $K\left[
\left[  x_{1},...,x_{n}\right]  \right]  $ this is shown by using the
geometric series).
\end{remark}

\begin{remark}
For any
\index{semigroup ordering}%
semigroup ordering $>$ on $K\left[  x_{1},...,x_{n}\right]  $
\begin{align*}
L\left(  gf\right)   &  =L\left(  g\right)  L\left(  f\right) \\
L\left(  g+f\right)   &  \leq\max\left\{  L\left(  g\right)  ,L\left(
f\right)  \right\}
\end{align*}
hence,%
\[
S_{>}=\left\{  u\in K\left[  x_{1},...,x_{n}\right]  \backslash\left\{
0\right\}  \mid L\left(  u\right)  =1\right\}
\]
is multiplicatively closed.

$S_{>}=K^{\ast}\Leftrightarrow$ $>$ is
\index{global ordering}%
global.

$S_{>}=K\left[  x_{1},...,x_{n}\right]  \backslash\left\langle x_{1}%
,...,x_{n}\right\rangle \Leftrightarrow$ $>$ is
\index{local ordering}%
local.
\end{remark}

\begin{definition}
Let $>$ be a
\index{semigroup ordering}%
semigroup ordering on $K\left[  x_{1},...,x_{n}\right]  $. The
\index{localization|textbf}%
\textbf{localization of }$K\left[  x_{1},...,x_{n}\right]  $\textbf{
associated to} $>$ is
\[
K\left[  x_{1},...,x_{n}\right]  _{>}=S_{>}^{-1}K\left[  x_{1},...,x_{n}%
\right]  =\left\{  \frac{f}{u}\mid f,u\in K\left[  x_{1},...,x_{n}\right]
,L\left(  u\right)  =1\right\}
\]

\end{definition}

\begin{lemma}
\cite{DS Varieties Groebner Bases and Algebraic Curves}, \cite[Sec. 1.5]{GP A
Singular Introduction to Commutative Algebra} Given a
\index{semigroup ordering}%
semigroup ordering $>$ on the monomials of $K\left[  x_{1},...,x_{n}\right]
$, there is a natural extension of the leading data to $K\left[
x_{1},...,x_{n}\right]  _{>}$: If $f\in K\left[  x_{1},...,x_{n}\right]  _{>}%
$, then there is a $u\in K\left[  x_{1},...,x_{n}\right]  $ with $LT\left(
u\right)  =1$ and $uf\in K\left[  x_{1},...,x_{n}\right]  $. The element
$L\left(  uf\right)  $ is independent of the choice of $u$ and is
\index{L}%
called $L\left(  f\right)  $, in the same
\index{LT}%
way define $LT\left(  f\right)  :=LT\left(  uf\right)  $
\index{LC}%
and $LC\left(  f\right)  :=LC\left(  uf\right)  $.

$LT\left(  f\right)  $
\index{LT}%
corresponds to a unique term in the power series expansion of $f$ and
subtracting this term gives the
\index{tail|textbf}%
\textbf{tail} of $f$ \newsym[$tail\left(  f\right)  $]{tail of $f$}{}denoted
by $tail\left(  f\right)  $.
\end{lemma}

\begin{definition}
Given a
\index{semigroup ordering}%
semigroup ordering $>$ for any subset $G\subset K\left[  x_{1},...,x_{n}%
\right]  _{>}$ define the
\index{lead ideal|textbf}%
\textbf{lead ideal} of $G$ as%
\[
L\left(  G\right)  =\text{ }_{K\left[  x_{1},...,x_{n}\right]  }\left\langle
L\left(  g\right)  \mid g\in G\backslash\left\{  0\right\}  \right\rangle
\]

\end{definition}

\subsubsection{Normal forms}

Fix a
\index{semigroup ordering}%
semigroup ordering $>$ on $K\left[  x_{1},...,x_{n}\right]  $ and let
$R=K\left[  x_{1},...,x_{n}\right]  _{>}$.

\begin{definition}
Let $\mathcal{G}$ be the set \newsym[$NF$]{normal form}{}of all finite lists
of elements in $R$. A map%
\[
NF:R\times\mathcal{G}\rightarrow R
\]
is called a
\index{weak normal form|textbf}%
\textbf{weak normal form} on $R$ if

\begin{enumerate}
\item $NF\left(  0,G\right)  =0$ $\forall G\in\mathcal{G}$

\item For all $G\in\mathcal{G}$ and $f\in R$%
\[
NF\left(  f,G\right)  \neq0\Rightarrow L\left(  NF\left(  f,G\right)  \right)
\notin L\left(  G\right)
\]

\item For all $G=\left\{  g_{1},...,g_{r}\right\}  \in\mathcal{G}$ and $f\in
R$ there is a unit $u\in R^{\ast}$ with either

\begin{itemize}
\item $uf=NF\left(  f,G\right)  $ or

\item $uf-NF\left(  f,G\right)  =\sum_{i=1}^{r}a_{i}g_{i}$ with $a_{i}\in R$
and for all $i$ with $a_{i}g_{i}\neq0$%
\[
L\left(  f\right)  \geq L\left(  a_{i}g_{i}\right)
\]

\end{itemize}
\end{enumerate}

Furthermore:

\begin{itemize}
\item $NF$ is called a
\index{normal form|textbf}%
\textbf{normal form} if one can always take $u=1$.

\item $NF$ is called
\index{polynomial normal form|textbf}%
\textbf{polynomial} if $f$ and $G$ are in $K\left[  x_{1},...,x_{n}\right]  $,
then also $u$ and $a_{i}$ can be taken in $K\left[  x_{1},...,x_{n}\right]  $.
A normal form is called
\index{reduced normal form|textbf}%
\textbf{reduced} if no monomial of $NF\left(  f,G\right)  $ is divisible by
some $L\left(  g_{i}\right)  $.

\item If the above properties are satisfied for some fixed $G\in\mathcal{G}$
we call $NF\left(  -,G\right)  :R\rightarrow R$ a (weak, polynomial, reduced)
\textbf{normal form with respect to }$G$.
\end{itemize}
\end{definition}

\begin{remark}
If $NF$ is polynomial, then $u\in R^{\ast}\cap K\left[  x_{1},...,x_{n}%
\right]  =S_{>}$.

From any weak normal form $NF$, we can build a normal form by dividing by $u$,
but the result will no longer be a polynomial.

Weak
\index{weak normal form}%
normal forms are introduced because they allow finite algorithmic computations
in $R=K\left[  x_{1},...,x_{n}\right]  _{>}$ and in $R$ a weak division
expression $uf-NF\left(  f,G\right)  =\sum_{i=1}^{r}a_{i}g_{i}$ is as good as
a division expression given by a normal form.
\end{remark}

\begin{definition}
If $f,g\in R\backslash\left\{  0\right\}  $, then
\newsym[$\operatorname*{SPolynomial}$]{$S$-polynomial}{}their
\index{S-polynomial|textbf}%
$S$\textbf{-polynomial} is%
\[
\operatorname*{SPolynomial}\left(  f,g\right)  =\frac{\operatorname{lcm}%
\left(  L\left(  f\right)  ,L\left(  g\right)  \right)  }{L\left(  f\right)
}f-\frac{LC\left(  f\right)  }{LC\left(  g\right)  }\frac{\operatorname{lcm}%
\left(  L\left(  f\right)  ,L\left(  g\right)  \right)  }{L\left(  g\right)
}g
\]

\end{definition}

\begin{remark}
This is implemented in
\index{mora.m2}%
\textsf{\textrm{mora.m2}} in the function $\mathsf{SPolynomial}$.
\end{remark}

\begin{algorithm}
[Gr\"{o}bner normal form]\label{Alg Groebner normal form}\cite{DS Varieties
Groebner Bases and Algebraic Curves},\newline\cite[Sec. 1.6]{GP A Singular
Introduction to Commutative Algebra} Let $>$ be a
\index{semigroup ordering}%
semigroup ordering
\index{Gr\"{o}bner normal form|textbf}%
and $G\in\mathcal{G}$.

Let $\operatorname*{divmon}:=\left(  h,G\right)  \mapsto\left(  g\in G\mid
L\left(  g\right)  \text{ divides }L\left(  h\right)  \right)  ;$

For any ordering of the sequences
\newsym[$NFG$]{Gr\"{o}bner normal form}{}produced by $\operatorname*{divmon}$,
the following algorithm is a
\index{NFG}%
normal form $f\mapsto NFG\left(  f,G\right)  :=h$ with respect to $G$.

$h:=f;$

$while$ $(h\neq0$ and $\operatorname*{divmon}\left(  h,G\right)  \neq
\emptyset)$ $\operatorname*{do}$ $($

\quad$g:=\operatorname*{divmon}\left(  h,G\right)  \#0;$

\quad$h:=\operatorname*{SPolynomial}\left(  h,g\right)  ;$

$);$

$h;$

This algorithm terminates if $>$ is a
\index{well ordering}%
well ordering. Otherwise
\index{NFG}%
$NFG$ may compute a power series convergent in the $\left\langle
x_{1},...,x_{n}\right\rangle $-adic topology.
\end{algorithm}

\begin{remark}
This is implemented in $%
\index{mora.m2}%
$\textsf{mora.m2} in the function $NFG$.
\end{remark}

\begin{algorithm}
[Gr\"{o}bner reduced normal form]\label{Alg Groebner reduced normal form}%
\cite{DS Varieties Groebner Bases and Algebraic Curves}, \cite[Sec. 1.6]{GP A
Singular Introduction to Commutative Algebra} Let $>$ be a
\index{semigroup ordering}%
semigroup ordering \newsym[$redNFG$]{Gr\"{o}bner reduced normal form}{}and
$G\in\mathcal{G}$. The
\index{Gr\"{o}bner reduced normal form|textbf}%
following
\index{redNFG}%
algorithm is a reduced normal form $f\mapsto redNFG\left(  f,G\right)
:=1/LC\left(  h\right)  \cdot h$ with respect to $G$:

$h:=0;g:=f;$

$\operatorname*{while}$ $g\neq0$ $\operatorname*{do}$ $($

\quad$g:=NFG\left(  g,G\right)  ;$

\quad$\operatorname*{if}$ $g\neq0$ $\operatorname*{then}$ $($

\quad\quad$h:=h+LT\left(  g\right)  ;$

\quad\quad$g:=g-LT\left(  g\right)  ;$

\quad$);$

$);$

$1/LC\left(  h\right)  \cdot h;$

This algorithm terminates if $>$ is a
\index{well ordering}%
well ordering.
\end{algorithm}

\begin{remark}
This is implemented in $%
\index{mora.m2}%
$\textsf{mora.m2} in the function $redNFG$.
\end{remark}

\begin{remark}
If we apply $NFG$ for the anti-degree order on $K\left[  x\right]  $, then
dividing $x$ by $x-x^{2}$, we get in $K\left[  \left[  x\right]  \right]  $
\[
x=\left(  \sum_{k=0}^{\infty}x^{k}\right)  \left(  x-x^{2}\right)  +0
\]
If we use a
\index{weak normal form}%
weak normal form, we can write%
\[
\left(  1-x\right)  \cdot x=1\cdot\left(  x-x^{2}\right)  +0
\]

\end{remark}

\begin{algorithm}
[Mora weak normal form]\label{Alg Mora normal form}\cite{Mora An Algorithm to
Compute the Equations of Tangent Cones},\newline\cite{DS Varieties Groebner
Bases and Algebraic Curves}, \cite[Sec. 1.6]{GP A Singular Introduction to
Commutative Algebra} Let $>$ be a
\index{semigroup ordering}%
semigroup ordering
\index{Mora normal form|textbf}%
and $G\in\mathcal{G}$. For polynomial input and any ordering of the sequences
produced by $\operatorname*{mecart}$ the following
\index{NFM}%
algorithm $f\mapsto NFM\left(  f,G\right)  :=h$ is a polynomial
\index{weak normal form}%
weak normal form with respect to $G$:

Let

$\operatorname*{ecart}:=f\mapsto\deg\left(  f\right)  -\deg LM\left(
f\right)  ;$

$\operatorname*{mecart}:=L\mapsto$the element of the sequence $L$ with minimal
ecart and minimal index;

$h:=f;$

$T:=G;$

$\operatorname*{while}$ $(h\neq0$ and $\operatorname*{divmon}\left(
h,T\right)  \neq\emptyset)$ $\operatorname*{do}$ $($

\quad$g:=\operatorname*{mecart}\left(  \operatorname*{divmon}\left(
h,T\right)  \right)  ;$

\quad$\operatorname*{if}$ $\operatorname*{ecart}\left(  g\right)
>\operatorname*{ecart}\left(  h\right)  $ $\operatorname*{then}$
$T:=\operatorname*{append}(T,h);$

\quad$h:=\operatorname*{SPolynomial}(h,g);$

$);$

$h;$

\noindent The algorithm terminates.
\end{algorithm}

\begin{remark}
This algorithm is implemented in $%
\index{mora.m2}%
$\textsf{mora.m2} in the function
\index{NF}%
$NF$.
\end{remark}

\begin{remark}
The Mora algorithm allows reductions also by the results of previous
reductions, in the above example%
\[
x=1\cdot\left(  x-x^{2}\right)  +x^{2}%
\]
so we also allow reduction by $x^{2}$, i.e.,%
\[
x=1\cdot\left(  x-x^{2}\right)  +1\cdot x^{2}+0
\]
which, as desired, also can be written%
\[
\left(  1-x\right)  \cdot x=1\cdot\left(  x-x^{2}\right)  +0
\]

\end{remark}

\begin{remark}
For homogeneous input $\operatorname*{ecart}$ is $0$, hence $NF$ and $NFG$
agree (for the same choice of the ordering of the list produced by
$\operatorname*{divmon}$).

If $>$ is a
\index{well ordering}%
well ordering, then any element $h$ appended to $T$ will not be used in
further steps: If it would be used, then $L\left(  h\right)  \mid L\left(
h_{new}\right)  $, hence $L\left(  h\right)  <L\left(  h_{new}\right)  $ or
$L\left(  h\right)  =L\left(  h_{new}\right)  $, as $>$ is a
\index{well ordering}%
well ordering. On the other hand the lead term of $h$ was canceled in a
previous step, so $L\left(  h\right)  >L\left(  h_{new}\right)  $.
\end{remark}

\subsubsection{Standard bases}

Let $R=K\left[  x_{1},...,x_{n}\right]  _{>}$ and fix a
\index{semigroup ordering}%
semigroup ordering $>$.

\begin{definition}
Let $I\subset R$ be an ideal. A finite subset $G\subset I$ is called a
\index{standard basis|textbf}%
\textbf{standard basis} (or
\index{Gr\"{o}bner basis|textbf}%
\textbf{Gr\"{o}bner basis} if $>$ is global) of $I$ if $L\left(  I\right)
=L\left(  G\right)  $, equivalently, if for any $f\in I\backslash\left\{
0\right\}  $ there is some $g\in G$ with $L\left(  g\right)  \mid L\left(
f\right)  $.
\end{definition}

\begin{proposition}
\cite{DS Varieties Groebner Bases and Algebraic Curves}, \cite[Sec. 1.6]{GP A
Singular Introduction to Commutative Algebra} Let $G\subset I\subset R$ be a
standard basis of the ideal $I$ and $NF\left(  -,G\right)  $ a
\index{weak normal form}%
weak normal form with respect to $G$, then:

\begin{enumerate}
\item For all $f\in R$ it holds%
\[
f\in I\Leftrightarrow NF\left(  f,G\right)  =0
\]

\item If $NF\left(  -,G\right)  $ is a
\index{reduced normal form}%
reduced normal form, then it is uniquely determined by $>$ and $I$ and denoted
by $NF\left(  -,I\right)  $.
\end{enumerate}
\end{proposition}

\begin{proposition}
\cite{DS Varieties Groebner Bases and Algebraic Curves}, \cite[Sec. 1.6]{GP A
Singular Introduction to Commutative Algebra}\label{IdealLeadideal} If
$G\subset I\subset R$ is a standard basis of the ideal $I$ and $NF\left(
-,G\right)  $ a
\index{weak normal form}%
weak normal form with respect to $G$, then it holds:

\begin{enumerate}
\item If $J\subset R$ is an ideal with $I\subset J$ and $L\left(  I\right)
=L\left(  J\right)  $, then $I=J$.

\item $I=\left\langle G\right\rangle $.
\end{enumerate}
\end{proposition}

\begin{theorem}
[Buchberger test]\cite{DS Varieties Groebner Bases and Algebraic Curves},
\cite[Sec. 1.6]{GP A Singular Introduction to Commutative Algebra} Let $NF$ be
a
\index{weak normal form}%
weak normal form, $G\in\mathcal{G}$ and $I\subset R$ an ideal. Then the
following properties are equivalent:

\begin{enumerate}
\item $G$ is a standard basis of $I$.

\item $I=\left\langle G\right\rangle $ and $NF\left(
\operatorname*{SPolynomial}\left(  g_{i},g_{j}\right)  ,G\right)  =0$ for all
$i,j$.

\item $NF\left(  f,G\right)  =0$ for all $f\in I$.
\end{enumerate}
\end{theorem}

This leads to the following algorithm:

\begin{algorithm}
[Gr\"{o}bner basis, Standard basis]\label{Alg standard basis}\cite{Mora An
Algorithm to Compute the Equations of Tangent Cones}, \newline\cite{DS
Varieties Groebner Bases and Algebraic Curves}, \cite[Sec. 1.6]{GP A Singular
Introduction to Commutative Algebra} Let $NF$ be a
\index{weak normal form}%
weak normal form. Given $G\in\mathcal{G}$, the following
\index{standard basis algorithm|textbf}%
algorithm computes a
\index{standard basis}%
standard basis $S$ of $\left\langle G\right\rangle \subset R$:

$S:=G;$

$P:=\left\{  \left(  f,g\right)  \mid f,g\in S,f\neq g\right\}  ;$

$\operatorname*{while}$ $P\neq\varnothing$ $\operatorname*{do}$ $($

\quad choose $\left(  f,g\right)  \in P;$

\quad$h:=NF\left(  \operatorname*{SPolynomial}\left(  f,g\right)  ,S\right)
;$

\quad$\operatorname*{if}$ $h\neq0$ $\operatorname*{then}$ $($

\quad\quad$P:=P\cup\left\{  \left(  h,f\right)  \mid f\in S\right\}  ;$

\quad\quad$S:=S\cup\left\{  h\right\}  ;$

\quad$);$

$);$

$S;$

This algorithm terminates.
\end{algorithm}

\begin{remark}
This algorithm is implemented in $%
\index{mora.m2}%
$\textsf{mora.m2} in the function
\index{std}%
$\mathsf{Std}$.
\end{remark}

\begin{remark}
Note that termination is only up to termination of $NF$. For a well-ordering
$NFG$ terminates, otherwise $NFG$ may compute a power series convergent in the
$\left\langle x_{1},...,x_{n}\right\rangle $-adic topology. In this case we
can use the
\index{Mora normal form}%
Mora normal form instead, which for polynomial input will terminate with
polynomial output, hence also the
\index{standard basis}%
standard basis algorithm will.
\end{remark}

Given any
\index{semigroup ordering}%
semigroup ordering $>$ on the monomials of $K\left[  x_{1},...,x_{n}\right]
$, the following ordering, introducing one additional variable $s$ to
homogenize the equations, can be used to compute standard bases via the
Gr\"{o}bner normal form.

\begin{definition}
For $f\in K\left[  x_{1},...,x_{n}\right]  $ of degree $d$ define
\[
f^{h}=s^{d}f\left(  \frac{x_{1}}{s},...,\frac{x_{n}}{s}\right)  \in K\left[
s,x_{1},...,x_{n}\right]
\]
to be \newsym[$f^{h}$]{homogenization}{}its
\index{homogenization|textbf}%
\textbf{homogenization}.
\end{definition}

\begin{definition}
\label{Lazard order}Let $A\in GL\left(  n,\mathbb{R}\right)  $ be the matrix
associated to $>$ and the
\index{semigroup ordering}%
semigroup ordering on the monomials of $K\left[  s,x_{1},...,x_{n}\right]  $
given by the matrix%
\[
\left(
\begin{tabular}
[c]{llll}%
$1$ & $1$ & $\cdots$ & $1$\\\cline{2-4}%
$0$ & \multicolumn{1}{|l}{} &  & \multicolumn{1}{l|}{}\\
$\vdots$ & \multicolumn{1}{|l}{} & $A$ & \multicolumn{1}{l|}{}\\
$0$ & \multicolumn{1}{|l}{} &  & \multicolumn{1}{l|}{}\\\cline{2-4}%
\end{tabular}
\ \right)
\]
i.e.,%
\[%
\begin{tabular}
[c]{lll}%
$s^{a}x^{\alpha}>s^{b}x^{\beta}$ & $\Leftrightarrow$ & $a+\left\vert
\alpha\right\vert >b+\left\vert \beta\right\vert $\\
&  & or\\
&  & $a+\left\vert \alpha\right\vert =b+\left\vert \beta\right\vert $ and
$x^{\alpha}>x^{\beta}$%
\end{tabular}
\ \text{ }%
\]

\end{definition}

\begin{algorithm}
[Lazard method]\cite[Sec. 1.7]{GP A Singular Introduction to Commutative
Algebra} Given polynomial $G=\left\{  g_{1},...,g_{r}\right\}  \in
\mathcal{G}$
\index{Lazard method|textbf}%
the following algorithm computes a
\index{standard basis}%
standard basis $S$ of $\left\langle G\right\rangle \subset K\left[
x_{1},...,x_{n}\right]  $ (note that we do not need $R=K\left[  x_{1}%
,...,x_{n}\right]  _{>}$ coefficients):

Using the
\index{Gr\"{o}bner normal form}%
Gr\"{o}bner normal form $NFG$, apply the
\index{standard basis}%
standard basis algorithm to $\left\{  g_{1}^{h},...,g_{r}^{h}\right\}  $ with
the induced monomial order from Definition \ref{Lazard order} to compute $S$
and put $s=1$.
\end{algorithm}

\begin{remark}
This algorithm is implemented in $%
\index{mora.m2}%
$\textsf{mora.m2} in the function
\index{Lstd}%
$\mathsf{LStd}$.
\end{remark}

Being a
\index{standard basis}%
standard basis depends only on finitely many monomials.

\begin{theorem}
\cite[Sec. 1.7]{GP A Singular Introduction to Commutative Algebra}%
\label{1finitedetermstd} For any ideal $I\subset K\left[  x_{1},...,x_{n}%
\right]  $ and
\index{standard basis}%
standard basis $S$ of $I$ with respect to $>$, there is a finite set of
monomials $F$ (i.e., all monomials appearing in the Buchberger test
computations) with the following property:

For all monomial orders $>_{1}$ identical to $>$ on $F$

\begin{enumerate}
\item $L_{>}\left(  g\right)  =L_{>_{1}}\left(  g\right)  $ $\forall g\in G$.

\item $G$ is also a
\index{standard basis}%
standard basis with respect to $>_{1}$.
\end{enumerate}
\end{theorem}

Hence for computing standard bases any monomial order can be represented by an
appropriate
\index{weight vector}%
weight vector.

Now consider the question of uniqueness:

\begin{definition}
A finite subset $G\subset R$ is called

\begin{itemize}
\item \textbf{interreduced}
\index{interreduced|textbf}%
(or
\index{minimal|textbf}%
\textbf{minimal}) if $0\notin G$ and $L\left(  f\right)  \nmid L\left(
g\right)  $ for all $f\neq g$.

\item \textbf{reduced}
\index{reduced|textbf}%
if it is interreduced and for all $f,g$ no term of
\index{tail}%
$tail\left(  g\right)  \in K\left[  \left[  x_{1},...,x_{n}\right]  \right]  $
is divisible by some $L\left(  f\right)  $.
\end{itemize}
\end{definition}

\begin{remark}
If $>$ is
\index{global ordering}%
global no term
\index{tail}%
of $tail\left(  g\right)  $ is divisible by $L\left(  g\right)  $, hence $G$
is reduced if for all $f\neq g$ no term of $g$ is divisible by $L\left(
f\right)  $.
\end{remark}

By Proposition \ref{IdealLeadideal} the following algorithm computes an
interreduced
\index{standard basis}%
standard basis of $I$:

\begin{algorithm}
Let $G$ be a
\index{standard basis}%
standard basis of $I$. Deleting successively all elements $g$ with $L\left(
f\right)  \mid L\left(  g\right)  $ for some $f\in G$, $f\neq g$ leads to an
interreduced
\index{standard basis}%
standard basis of $I$.
\end{algorithm}

\begin{remark}
This algorithm is implemented in $%
\index{mora.m2}%
$\textsf{mora.m2} in the function
\index{minimizeStd}%
$\mathsf{MinimizeStd}$ and
\index{mStd}%
$\mathsf{MStd}$ computes an interreduced
\index{standard basis}%
standard basis using
\index{Mora normal form}%
Mora normal form and applying $\mathsf{MinimizeStd}$.
\end{remark}

\begin{algorithm}
\label{Alg reduced std}If the input generators of $I$ for the standard basis
algorithm were reduced, and the
\index{standard basis}%
standard basis algorithm used a reduced normal form $NF$, then after
minimalization the output is also reduced. If the input was not reduced and
$G=\left\{  g_{1},...,g_{n}\right\}  $ is the interreduced output of the
\index{standard basis}%
standard basis algorithm using a reduced normal form $NF$, then%
\[
M=\left\{  NF_{>}\left(  g_{i},\left\{  g_{1},...,g_{i-1},g_{i+1}%
,...,g_{n}\right\}  \right)  \mid i=1,...,n\right\}
\]
is the
\index{reduced standard basis|textbf}%
reduced standard basis of $I$.
\end{algorithm}

\begin{remark}
Then $\left\{  L\left(  g_{i}\right)  \mid i=1,...,n\right\}  $ is the minimal
generating set of $L\left(  I\right)  $ and%
\[
M=\left\{  NF_{>}\left(  L\left(  g_{i}\right)  ,G\right)  \mid
i=1,...,n\right\}
\]
Hence $M$ is uniquely determined by $I$ and $>$, as $NF_{>}\left(  -,G\right)
$ is.
\end{remark}

\begin{algorithm}
\label{Alg reduced std mora}As reduced normal form $NF$ one can use the
reduced Gr\"{o}bner normal form $redNFG_{>}$. If $>$ is
\index{global ordering}%
global, then $redNFG_{>}$ terminates with an element in $K\left[
x_{1},...,x_{n}\right]  $, otherwise $redNFG_{>}$ computes an element in
$K\left[  \left[  x_{1},...,x_{n}\right]  \right]  $ in general.

If $G=\left\{  g_{1},...,g_{n}\right\}  $ is an
\index{interreduced}%
interreduced
\index{standard basis}%
standard basis of $I$ with respect to $>$, computed using any weak normal form
(e.g., Mora normal form), and $NF$ is a reduced normal form (e.g., Gr\"{o}bner
normal form), then%
\[
\left\{  L\left(  g_{i}\right)  +NF_{>}\left(  tail\left(  g_{i}\right)
,G\right)  \mid i=1,...,n\right\}  \subset K\left[  \left[  x_{1}%
,...,x_{n}\right]  \right]
\]
is the unique reduced standard basis of $I$ with respect $>$.
\end{algorithm}

\begin{remark}
This is implemented in $%
\index{mora.m2}%
$\textsf{mora.m2} in the function
\index{reduceGb}%
$\mathsf{ReduceGb}$. For the non global case the number iterations can be
limited by the global variable $\mathsf{iterlimit}$.
\end{remark}

\subsubsection{Localization in prime ideals}

\begin{proposition}
\cite[Sec. 1.5]{GP A Singular Introduction to Commutative Algebra} Let $K$ be
a field, $>$ a local ordering on the polynomial ring $K\left[  x_{1}%
,...,x_{n}\right]  $. Then the localization of the polynomial ring $K\left[
x_{1},...,x_{n},y_{1},...,y_{m}\right]  $ at the prime ideal $\left\langle
x_{1},...,x_{n}\right\rangle $ is the localization of $K\left(  y_{1}%
,...,y_{m}\right)  \left[  x_{1},...,x_{n}\right]  $ with respect to $>$%
\[
K\left(  y_{1},...,y_{m}\right)  \left[  x_{1},...,x_{n}\right]  _{>}=K\left[
x_{1},...,x_{n},y_{1},...,y_{m}\right]  _{\left\langle x_{1},...,x_{n}%
\right\rangle }%
\]

\end{proposition}

Recall that%
\begin{align*}
K\left(  y_{1},...,y_{m}\right)  \left[  x_{1},...,x_{n}\right]  _{>}  &
=U_{>}^{-1}\left(  K\left(  y_{1},...,y_{m}\right)  \left[  x_{1}%
,...,x_{n}\right]  \right) \\
&  =\left\{  \frac{f}{u}\mid f,u\in K\left(  y_{1},...,y_{m}\right)  \left[
x_{1},...,x_{n}\right]  ,L\left(  u\right)  =1\right\}
\end{align*}
with%
\[
U_{>}=\left\{  u\in K\left(  y_{1},...,y_{m}\right)  \left[  x_{1}%
,...,x_{n}\right]  \backslash\left\{  0\right\}  \mid L\left(  u\right)
=1\right\}
\]

This allows to do Gr\"{o}bner computations in localizations at prime ideals
$\left\langle x_{1},...,x_{n}\right\rangle $, e.g. at the ideals of the strata
of a toric variety in the Cox ring.

\section{Mirror constructions
\index{mirror construction}%
to generalize\label{Sec mirror constructions to generalize}}

\subsection{Batyrev%
\'{}%
s construction for hypersurfaces in toric varieties\label{Sec Batyrev}}

Let $Y=\mathbb{P}\left(  \Delta\right)  $
\index{Batyrev}%
be a toric
\index{Fano}%
Fano
\index{toric Fano}%
variety of $\dim Y=n$ represented by the
\index{reflexive}%
reflexive polytope $\Delta\subset M_{\mathbb{R}}$ and let
$N=\operatorname*{Hom}\left(  M,\mathbb{Z}\right)  $.

\begin{proposition}
\cite{Batyrev Dual polyhedra and mirror symmetry for CalabiYau hypersurfaces
in toric varieties}, \cite[Sec. 4.1.1]{CK Mirror Symmetry and Algebraic
Geometry}, \cite{Reid Canonical 3folds} A general element in $\left|
-K_{\mathbb{P}\left(  \Delta\right)  }\right|  $ is a
\index{Calabi Yau variety}%
Calabi Yau variety
\index{hypersurface}%
of dimension $n-1$.
\end{proposition}

\begin{theorem}
\cite{Batyrev Dual polyhedra and mirror symmetry for CalabiYau hypersurfaces
in toric varieties}, \cite{BB Mirror duality and stringtheoretic Hodge
numbers} For any
\index{reflexive}%
reflexive $\Delta$ general
\newsym[$\left\vert -K_{\mathbb{P}\left(  \Delta\right)  }\right\vert $]{anticanonical linear system}{}elements
$X$ of $\left\vert -K_{\mathbb{P}\left(  \Delta\right)  }\right\vert $ and
$X^{\circ}$ of $\left\vert -K_{\mathbb{P}\left(  \Delta^{\ast}\right)
}\right\vert $ are
\index{stringy topological mirror pair}%
stringy topological mirror pairs (indeed
\index{mathematical mirror pair}%
mathematical mirror pairs), and
\index{mirror construction}%
there are
\index{stringy Hodge numbers}%
explicit formulas computing $h_{st}^{d-1,1}\left(  X\right)  $ and
$h_{st}^{1,1}\left(  X\right)  $ from the polytope:%
\begin{align}
h_{st}^{d-1,1}\left(  X\right)   &  =\left\vert \Delta\cap M\right\vert
-n-1-\sum_{Q\text{ facet of }\Delta}\left\vert \operatorname*{int}%
\nolimits_{M}\left(  Q\right)  \right\vert \label{4BatyrevHodgeformula}\\
&  +\sum_{\substack{Q\text{ face of }\Delta\\\operatorname*{codim}%
Q=2}}\left\vert \operatorname*{int}\nolimits_{M}\left(  Q\right)  \right\vert
\cdot\left\vert \operatorname*{int}\nolimits_{N}\left(  Q^{\ast}\right)
\right\vert \nonumber\\
h_{st}^{1,1}\left(  X\right)   &  =\left\vert \Delta^{\ast}\cap M\right\vert
-n-1-\sum_{Q^{\ast}\text{ facet of }\Delta^{\ast}}\left\vert
\operatorname*{int}\nolimits_{N}\left(  Q^{\ast}\right)  \right\vert
\nonumber\\
&  +\sum_{\substack{Q^{\ast}\text{ face of }\Delta^{\ast}%
\\\operatorname*{codim}Q^{\ast}=2}}\left\vert \operatorname*{int}%
\nolimits_{N}\left(  Q^{\ast}\right)  \right\vert \cdot\left\vert
\operatorname*{int}\nolimits_{M}\left(  Q\right)  \right\vert \nonumber
\end{align}
Here,
\index{int}%
$\operatorname*{int}_{M}\left(  Q\right)  $ denotes the set of lattice points
in the relative interior of the face $Q$ with respect to the lattice $M$.
\end{theorem}

\begin{example}
As
\index{hypersurface}%
discussed
\index{quintic threefold}%
in Examples \ref{2exquinticchow} and \ref{2quinticCoxmon} the
\index{reflexive}%
reflexive degree $5$
\index{Veronese}%
Veronese polytope $\Delta$ of $\mathbb{P}^{4}$ and its
\index{dual polytope}%
dual yield $\mathbb{P}\left(  \Delta\right)  =\mathbb{P}^{4}$ and
$\mathbb{P}\left(  \Delta^{\ast}\right)  =\mathbb{P}^{4}/\mathbb{Z}_{5}^{3}$.
General anticanonical
\index{hypersurface}%
hypersurfaces $X$ and $X^{\circ}$ inside satisfy%
\begin{align*}
h^{1,1}\left(  X\right)   &  =h_{st}^{2,1}\left(  X^{\circ}\right)  =\left(
6-1\right)  -4=1\\
h^{2,1}\left(  X\right)   &  =h_{st}^{1,1}\left(  X^{\circ}\right)  =\left(
126-1\right)  -24=101
\end{align*}
As noticed in Example \ref{2exquinticquotientdescription}, $X$ is given by a
general element in $S_{\left[  -K_{\mathbb{P}\left(  \Delta\right)  }\right]
}$, which is a general degree $5$ polynomial in $\mathbb{C}\left[
x_{1},...,x_{5}\right]  $ and $X^{\circ}$ by a general element in $S_{\left[
-K_{\mathbb{P}\left(  \Delta^{\ast}\right)  }\right]  }$, which is
\[
\sum_{i=1}^{5}c_{i}y_{i}^{5}+c_{0}y_{1}\cdot...\cdot y_{5}%
\]
Modulo
\index{automorphism}%
automorphisms of $\mathbb{P}\left(  \Delta^{\ast}\right)  $, this is the one
parameter family obtained from the
\index{Greene-Plesser}%
Greene-Plesser construction, see Example \ref{GreenePlesserQuinticExample}.
\end{example}

The singularities may be resolved
\index{crepant resolution}%
crepantly via
\index{maximal projective subdivision}%
maximal projective subdivisions $\bar{\Sigma}$ of the fan $\Sigma$.

\begin{definition}
Given a
\index{reflexive}%
reflexive polytope $\Delta\subset M_{\mathbb{R}}$, a fan $\bar{\Sigma}$ in
$N_{\mathbb{R}}$ is called a
\index{projective subdivision|textbf}%
\textbf{projective subdivision} of the
\index{normal fan}%
normal fan $\Sigma$ of $\Delta$ if

\begin{enumerate}
\item $\bar{\Sigma}$ refines $\Sigma$.

\item $\bar{\Sigma}\left(  1\right)  \subset\Delta^{\ast}\cap N-\left\{
0\right\}  $

\item $X\left(  \bar{\Sigma}\right)  $ is projective and
\index{simplicial}%
simplicial.
\end{enumerate}

$\bar{\Sigma}$ is called
\index{maximal projective subdivision|textbf}%
\textbf{maximal} if $\bar{\Sigma}\left(  1\right)  =\Delta^{\ast}\cap
N-\left\{  0\right\}  $.
\end{definition}

\begin{proposition}
\cite{OP Linear Gale transforms and GelfandKapranovZelevinsky decompositions},
\cite[4.1.1]{CK Mirror Symmetry and Algebraic Geometry}, \cite[3.4]{CK Mirror
Symmetry and Algebraic Geometry} For any
\index{reflexive}%
reflexive $\Delta$, there exists a maximal projective subdivision of the
\index{normal fan}%
normal fan $\Sigma$ of $\Delta$.
\end{proposition}

\begin{proposition}
\cite{Batyrev Dual polyhedra and mirror symmetry for CalabiYau hypersurfaces
in toric varieties}, \cite[4.1.1]{CK Mirror Symmetry and Algebraic Geometry}
For any
\index{reflexive}%
reflexive $\Delta$, any projective subdivision $\bar{\Sigma}$ of the
\index{normal fan}%
normal fan $\Sigma$ of $\Delta$ gives a birational morphism%
\[
f:X\left(  \bar{\Sigma}\right)  \rightarrow\mathbb{P}\left(  \Delta\right)
\]
and

\begin{enumerate}
\item $X\left(  \bar{\Sigma}\right)  $ is a
\index{Gorenstein orbifold}%
Gorenstein orbifold.

\item $f$ is
\index{crepant resolution}%
crepant, i.e., $f^{\ast}K_{\mathbb{P}\left(  \Delta\right)  }=K_{X\left(
\bar{\Sigma}\right)  }$.

\item If the subdivision $\bar{\Sigma}$ is maximal, then $X\left(  \bar
{\Sigma}\right)  $ has
\index{terminal singularities}%
terminal singularities.
\end{enumerate}

Furthermore, for a general element $\bar{X}$ in $\left\vert -K_{X\left(
\bar{\Sigma}\right)  }\right\vert $ it holds:

\begin{enumerate}
\item $\bar{X}$ is a Calabi-Yau orbifold.

\item If the subdivision $\bar{\Sigma}$ is
\index{maximal subdivision}%
maximal, then $\bar{X}$ is called a
\index{maximal projective crepant partial desingularization|textbf}%
\textbf{maximal projective crepant partial
\index{MPCP desingularization|textbf}%
(MPCP) desingularization} of $X$ and has the following properties:

\begin{enumerate}
\item $\bar{X}$ is a
\index{minimal Calabi-Yau orbifold}%
minimal Calabi-Yau orbifold in the sense of Mori theory.

\item $\bar{X}$ is the
\index{proper transform}%
proper transform by $f$ of a general element $X$ in $\left\vert -K_{\mathbb{P}%
\left(  \Delta\right)  }\right\vert $. The induced map $f:\bar{X}\rightarrow
X$ is
\index{crepant resolution}%
crepant.

\item If $\dim X=3$, then $\bar{X}$ is smooth, as
\index{Gorenstein orbifold terminal singularities}%
Gorenstein orbifold terminal singularities are smooth in dimension $3$.
\end{enumerate}
\end{enumerate}
\end{proposition}

\begin{example}
Consider the
\index{weighted projective space}%
weighted projective space $\mathbb{P}\left(  1,1,2,2,2\right)  $ given by the
fan $\Sigma$ over
\[
\Delta^{\ast}=\operatorname*{convexhull}\left(  \left(  -1,-2,-2,-2\right)
,\left(  1,0,0,0\right)  ,\left(  0,1,0,0\right)  ,\left(  0,0,1,0\right)
,\left(  0,0,0,1\right)  \right)
\]
hence, via the
\index{Cox ring}%
Cox ring it has the description%
\[
\mathbb{P}\left(  1,1,2,2,2\right)  \cong\frac{\mathbb{C}^{5}-V\left(
\left\langle y_{0},...,y_{4}\right\rangle \right)  }{\mathbb{C}^{\ast}}%
\]
with $\mathbb{C}^{\ast}$-action%
\[
\mathbb{\lambda}\left(  y_{0},...,y_{4}\right)  =\left(  \mathbb{\lambda}%
y_{0},\mathbb{\lambda}y_{1},\mathbb{\lambda}^{2}y_{2},\mathbb{\lambda}%
^{2}y_{3},\mathbb{\lambda}^{2}y_{4}\right)
\]
A toric variety given by a fan $\Sigma$ is smooth if and only if for every
cone in $\Sigma$ the
\index{minimal lattice generator}%
minimal lattice generators are a subset of a $\mathbb{Z}$-basis of the lattice
$N$. So $\operatorname*{hull}\left(  \left(  -1,-2,-2,-2\right)  ,\left(
1,0,0,0\right)  \right)  $ is singular. The only lattice point of
$\Delta^{\ast}$, which is not a vertex is
\[
\frac{1}{2}\left(  \left(  -1,-2,-2,-2\right)  +\left(  1,0,0,0\right)
\right)  =\left(  0,-1,-1,-1\right)
\]
Hence, considering the fan $\bar{\Sigma}$ obtained by splitting all maximal
dimensional cones of $\Sigma$ into two via the new ray, we obtain a MPCP
desingularization and indeed a smooth toric variety $X\left(  \bar{\Sigma
}\right)  $, which is a blowup
\index{blowup}
of $\mathbb{P}\left(  1,1,2,2,2\right)  $. Via its
\index{Cox ring}%
Cox ring, $X\left(  \bar{\Sigma}\right)  $ is given by%
\[
X\left(  \bar{\Sigma}\right)  \cong\frac{\mathbb{C}^{6}-V\left(  \left\langle
y_{0},y_{1}\right\rangle \cap\left\langle y_{2}...,y_{4},y_{5}\right\rangle
\right)  }{\left(  \mathbb{C}^{\ast}\right)  ^{2}}%
\]
with $\left(  \mathbb{C}^{\ast}\right)  ^{2}$-action%
\[
\left(  \mathbb{\lambda},\mu\right)  \left(  y_{0},...,y_{5}\right)  =\left(
\mathbb{\lambda}y_{0},\mathbb{\lambda}y_{1},\mu y_{2},\mu y_{3},\mu
y_{4},\frac{\mu}{\mathbb{\lambda}^{2}}y_{5}\right)
\]

\end{example}

\subsection{Batyrev%
\'{}%
s and Borisov%
\'{}%
s construction for complete intersections in toric
varieties\label{Sec Batyrev and Borisov mirror construction}}

Consider
\index{mirror construction}%
a
\index{Batyrev}%
toric
\index{Fano}%
Fano
\index{toric Fano}%
variety $Y=\mathbb{P}\left(  \Delta\right)  $ represented by the
\index{reflexive}%
reflexive polytope $\Delta\subset M_{\mathbb{R}}$ with
\index{normal fan}%
normal fan $\Sigma\subset N_{\mathbb{R}}$. A disjoint union%
\[
\Sigma\left(  1\right)  =I_{1}\cup...\cup I_{c}%
\]
is
\index{ray}%
called a
\index{nef partition|textbf}%
\textbf{nef partition} if all $E_{j}=\sum_{v\in I_{j}}D_{v}$ are
\index{T-Cartier divisor}%
Cartier,
\index{spanned by global sections}%
spanned by global sections. By $\sum_{j=1}^{c}E_{j}=-K_{Y}$, general sections
of $\mathcal{O}\left(  E_{1}\right)  ,...,\mathcal{O}\left(  E_{c}\right)  $
give a Calabi-Yau
\index{complete intersection}%
complete intersection $X\subset Y$.

\begin{proposition}
\cite{BB On CalabiYau complete intersections in toric varieties in
HigherDimensional Complex Varieties Trento 1994} The polytopes $\Delta
_{j}=\Delta_{E_{j}}$ of sections of $E_{j}$ are
\index{lattice polytope}%
lattice polytopes (see Section \ref{Divisors on toric varieties}), and it
holds
\[
\Delta=\Delta_{1}+...+\Delta_{c}%
\]

\end{proposition}

\begin{figure}
[h]
\begin{center}
\includegraphics[
height=2.1577in,
width=1.5774in
]%
{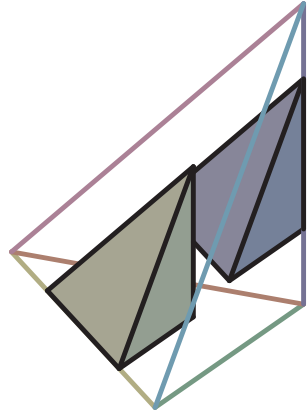}%
\caption{Batyrev-Borisov polytopes $\Delta_{1}$, $\Delta_{2}$ and their
Minkowski sum $\Delta$ for the $\left(  2,2\right)  $ complete intersection
elliptic curve in $\mathbb{P}^{3}$}%
\label{Fig22DeltaDelta1Delta2}%
\end{center}
\end{figure}

\begin{example}
\label{22example1}Let
\[
\Delta=\operatorname*{convexhull}\left(  \left(  -1,-1,-1\right)  ,\left(
3,-1,-1\right)  ,\left(  -1,3,-1\right)  ,\left(  -1,-1,3\right)  \right)
\]
be the
\index{reflexive}%
reflexive degree $4$
\index{Veronese}%
Veronese polytope of $\mathbb{P}^{3}$. By the partition of the $4$ vertices of
$\Delta^{\ast}$ into $I_{1}$ and $I_{2}$ each with $2$ elements%
\begin{align*}
I_{1}  &  =\left\{  \left(  -1,-1,-1\right)  ,\left(  0,0,1\right)  \right\}
\\
I_{2}  &  =\left\{  \left(  1,0,0\right)  ,\left(  0,1,0\right)  \right\}
\end{align*}
a general $\left(  2,2\right)  $
\index{elliptic curve!complete intersection}%
complete
\index{complete intersection}%
intersection elliptic curve in $\mathbb{P}^{3}$ is given and
\begin{align*}
\Delta_{1}  &  =\operatorname*{convexhull}\left(  \left(  1,-1,0\right)
,\left(  -1,-1,0\right)  ,\left(  -1,-1,2\right)  ,\left(  -1,1,0\right)
\right) \\
\Delta_{2}  &  =\operatorname*{convexhull}\left(  \left(  0,0,-1\right)
,\left(  0,0,1\right)  ,\left(  2,0,-1\right)  ,\left(  0,2,-1\right)
\right)
\end{align*}
are degree $2$
\index{Veronese}%
Veronese polytopes, which add up to $\Delta=\Delta_{1}+\Delta_{2}$, as shown
in Figure \ref{Fig22DeltaDelta1Delta2}.
\end{example}

%

Define the
\index{lattice polytope}%
lattice polytopes
\[
\nabla_{j}=\operatorname*{convexhull}\left\{  \left\{  0\right\}  \cup
I_{j}\right\}
\]
and define $\nabla_{BB}$ by%
\[
\nabla_{BB}^{\ast}=\operatorname*{convexhull}\left(  \Delta_{1}\cup
...\cup\Delta_{c}\right)
\]

\begin{proposition}
\cite{BB On CalabiYau complete intersections in toric varieties in
HigherDimensional Complex Varieties Trento 1994} $\nabla_{BB}=\nabla
_{1}+...+\nabla_{c}$.
\end{proposition}

In particular $\nabla_{BB}$ is a
\index{lattice polytope}%
lattice polytope containing $0$, hence:

\begin{corollary}
$\nabla_{BB}$ is
\index{reflexive}%
reflexive.
\end{corollary}

\begin{figure}
[h]
\begin{center}
\includegraphics[
trim=-0.010218in 0.000000in 0.010218in 0.000000in,
height=2.3981in,
width=1.593in
]%
{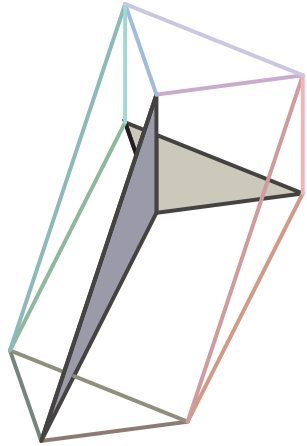}%
\caption{Batyrev-Borisov polytopes $\nabla_{1}$, $\nabla_{2}$ and their
Minkowski sum $\nabla$ for the mirror of the $\left(  2,2\right)  $ complete
intersection elliptic curve in $\mathbb{P}^{3}$}%
\label{FigNablaNabla1Nabla2}%
\end{center}
\end{figure}

\begin{example}
\label{22example2}In the above
\index{elliptic curve!complete intersection}%
Example \ref{22example1}%
\begin{align*}
\nabla_{1}  &  =\operatorname*{convexhull}\left\{  \left(  0,0,0\right)
,\left(  -1,-1,-1\right)  ,\left(  0,0,1\right)  \right\} \\
\nabla_{2}  &  =\operatorname*{convexhull}\left\{  \left(  0,0,0\right)
,\left(  1,0,0\right)  ,\left(  0,1,0\right)  \right\}
\end{align*}
Figure \ref{FigNablaNabla1Nabla2} shows the polytopes $\nabla_{1}$,
$\nabla_{1}$ and
\index{complete intersection}%
their
\index{Minkowski sum}%
Minkowski sum $\nabla_{BB}$.
\end{example}

%

Let $Y^{\circ}=\mathbb{P}\left(  \nabla_{BB}\right)  $ be the toric
\index{Fano}%
Fano
\index{toric Fano}%
variety associated to $\nabla_{BB}$. Then%
\[
\sum_{j=1}^{c}D_{\nabla_{j}}=-K_{Y^{\ast}}%
\]
is a
\index{nef partition}%
nef partition, and $X^{\circ}$ given by general sections of $\mathcal{O}%
\left(  D_{\nabla_{1}}\right)  ,...,\mathcal{O}\left(  D_{\nabla_{c}}\right)
$ is a Calabi-Yau
\index{complete intersection}%
complete intersection in $Y^{\circ}$.

\begin{theorem}
\cite{BB Mirror duality and stringtheoretic Hodge numbers} $X$ and $X^{\circ}$
form a stringy topological mirror pair.
\end{theorem}

A maximal projective subdivision $\bar{\Sigma}$ of $\Sigma=\operatorname*{NF}%
\left(  \Delta\right)  $ gives a maximal
\index{maximal projective subdivision}%
projective
\index{partial crepant resolution}%
partial crepant desingularization%
\[
f:X\left(  \bar{\Sigma}\right)  \rightarrow\mathbb{P}\left(  \Delta\right)
\]
such that the $T$-divisors of the projective toric variety $X\left(
\bar{\Sigma}\right)  $ correspond to the lattice points of the boundary of
$\Delta^{\ast}$. Then $f$ induces a resolution $\bar{X}\rightarrow X$ of the
complete intersection $X\subset\mathbb{P}\left(  \Delta\right)  $ such that
$\bar{X}$ is a complete intersection, has at most Gorenstein terminal abelian
quotient singularities and $K_{\bar{X}}=\mathcal{O}_{\bar{X}}$, for a
reference see \cite{BB On CalabiYau complete intersections in toric varieties
in HigherDimensional Complex Varieties Trento 1994}. In particular, if
$\dim\left(  \bar{X}\right)  \leq3$, then $\bar{X}$ is smooth.

\subsection{R\o dland%
\'{}%
s orbifolding
\index{mirror construction}%
mirror construction for the degree $14$ Pfaffian Calabi-Yau threefold in
$\mathbb{P}^{6}$\label{Roedlandexample}}

Consider
\index{mirror construction}%
a $7$-dimensional $\mathbb{C}$-vector space $V$ and the trivial vector bundle
$\mathcal{V}$ with fiber $V$ on $\mathbb{P}\left(  \bigwedge^{2}V\right)  $.
Define $M$ as the
\index{degeneracy locus}%
degeneracy locus of the
\index{universal skew symmetric linear map}%
universal skew symmetric linear map%
\[
\alpha:\mathcal{V}^{\ast}\left(  -1\right)  \rightarrow\mathcal{V}%
\]
i.e., as the locus $\operatorname*{rank}\alpha\leq4$.

$M$ is given by the $6\times6$
\index{Pfaffian}%
Pfaffians of $\alpha$, is locally
\index{Gorenstein}%
Gorenstein of codimension $3$ in $\mathbb{P}\left(  \bigwedge^{2}V\right)  $
and has $K_{M}=\mathcal{O}_{M}\left(  -14\right)  $. Its singular locus is
given by $\operatorname*{rank}\alpha\leq2$ and has codimension $7$ in $M$.

Intersecting $M$ with a general $\mathbb{P}^{d+3}\subset\mathbb{P}\left(
\bigwedge^{2}V\right)  $ gives $X^{d}=M\cap\mathbb{P}^{3+d}$ of dimension $d$
and $K_{X^{d}}=\mathcal{O}_{X^{d}}\left(  3-d\right)  $ and $X^{d}$ is smooth
for $d\leq6$.

$X=X^{3}$ is a
\index{locally complete intersection}%
local complete intersection
\index{Calabi Yau variety}%
Calabi-Yau threefold with $h^{1,1}\left(  X\right)  =1$ and $h^{1,2}\left(
X\right)  =50$:

\begin{remark}
By%
\[
0\rightarrow\mathcal{O}_{\mathbb{P}^{6}}\left(  -7\right)  \rightarrow
\mathcal{O}_{\mathbb{P}^{6}}\left(  -4\right)  ^{\oplus7}\rightarrow
\mathcal{O}_{\mathbb{P}^{6}}\left(  -3\right)  ^{\oplus7}\rightarrow
\mathcal{O}_{\mathbb{P}^{6}}\rightarrow\mathcal{O}_{X}%
\]
we get%
\[
H^{i}\left(  X,\mathcal{O}_{X}\right)  \cong H^{i+3}\left(  \mathbb{P}%
^{6},\mathcal{O}_{\mathbb{P}^{6}}\left(  -7\right)  \right)  \cong\left\{
\begin{array}
[c]{cc}%
0 & i=1,2\\
\mathbb{C} & i=3
\end{array}
\right\}
\]
and using the resolution of $\mathcal{J}_{X_{I}}^{2}$, the Euler sequence, the
definition of the normal sheaf and the conormal sequence%
\begin{align*}
0  &  \rightarrow\mathcal{O}_{\mathbb{P}^{6}}\left(  -8\right)  ^{\oplus
21}\rightarrow\mathcal{O}_{\mathbb{P}^{6}}\left(  -7\right)  ^{\oplus
48}\rightarrow\mathcal{O}_{\mathbb{P}^{6}}\left(  -6\right)  ^{\oplus
28}\rightarrow\mathcal{J}_{X}^{2}\rightarrow0\\
0  &  \rightarrow\Omega_{\mathbb{P}^{6}}\mid_{X}\rightarrow\mathcal{O}%
_{X}\left(  -1\right)  ^{\oplus7}\rightarrow\mathcal{O}_{X}\rightarrow0\\
0  &  \rightarrow\mathcal{J}_{X}^{2}\rightarrow\mathcal{J}_{X}\rightarrow
\mathcal{N}_{X/\mathbb{P}^{6}}^{\vee}\rightarrow0\\
0  &  \rightarrow\mathcal{N}_{X}^{\vee}\rightarrow\Omega_{\mathbb{P}^{6}}%
\mid_{X}\rightarrow\Omega_{X}\rightarrow0
\end{align*}
one
\index{Hodge numbers}%
computes%
\begin{align*}
h^{1,1}\left(  X\right)   &  =1\\
h^{1,2}\left(  X\right)   &  =50
\end{align*}

\end{remark}

The mirror is constructed via Greene-Plesser orbifolding in an analogous way
to \cite{GrPl Mirror Manifolds: A Brief Review and Progress Report} and
Example \ref{GreenePlesserQuinticExample}. As one expects for the mirror to
hold $h^{1,2}\left(  X^{\circ}\right)  =h^{1,1}\left(  X\right)  =1$, to
apply
\index{Greene-Plesser orbifolding}%
Greene-Plesser orbifolding, one looks for a $1$-parameter
\index{subfamily}%
subfamily, i.e., a $1$-parameter family of $\mathbb{P}^{6}\subset
\mathbb{P}\left(  \bigwedge^{2}V\right)  $. Choosing a basis $\left(
e_{i}\right)  $ of $V$ R\o dland considers the
\index{group action}%
action of the group $G=\left\langle \sigma,\tau\right\rangle \subset
\operatorname*{Aut}\mathbb{P}\left(  \bigwedge^{2}V\right)  $ of order $49$
generated by
\begin{align*}
\sigma &  =\left(
\begin{array}
[c]{ccccccc}%
0 & 0 & 0 & 0 & 0 & 0 & 1\\
1 & 0 & 0 & 0 & 0 & 0 & 0\\
0 & 1 & 0 & 0 & 0 & 0 & 0\\
0 & 0 & 1 & 0 & 0 & 0 & 0\\
0 & 0 & 0 & 1 & 0 & 0 & 0\\
0 & 0 & 0 & 0 & 1 & 0 & 0\\
0 & 0 & 0 & 0 & 0 & 1 & 0
\end{array}
\right) \\
\tau &  =\operatorname*{diag}\left(  \zeta^{i}\right)  _{i=0,...,6}%
\end{align*}
on $\mathbb{P}\left(  \bigwedge^{2}V\right)  $. Taking coordinates on
$\mathbb{P}^{6}$, the
\index{subfamily}%
subfamily $X_{y}$ invariant under the action of $G$ is given by the
\index{skew symmetric}%
skew symmetric matrix%
\begin{equation}
A_{y}=\left(
\begin{array}
[c]{ccccccc}%
0 & y_{1}x_{1} & y_{2}x_{2} & y_{3}x_{3} & -y_{3}x_{4} & -y_{2}x_{5} &
-y_{1}x_{6}\\
-y_{1}x_{1} & 0 & y_{1}x_{3} & y_{2}x_{4} & y_{3}x_{5} & -y_{3}x_{6} &
-y_{2}x_{0}\\
-y_{2}x_{2} & -y_{1}x_{3} & 0 & y_{1}x_{4} & y_{2}x_{6} & y_{3}x_{0} &
-y_{3}x_{1}\\
-y_{3}x_{3} & -y_{2}x_{4} & -y_{1}x_{4} & 0 & y_{1}x_{0} & y_{2}x_{1} &
y_{3}x_{2}\\
y_{3}x_{4} & -y_{3}x_{5} & -y_{2}x_{6} & -y_{1}x_{0} & 0 & y_{1}x_{2} &
y_{2}x_{3}\\
y_{2}x_{5} & y_{3}x_{6} & -y_{3}x_{0} & -y_{2}x_{1} & -y_{1}x_{2} & 0 &
y_{1}x_{4}\\
y_{1}x_{6} & y_{2}x_{0} & y_{3}x_{1} & -y_{3}x_{2} & -y_{2}x_{3} & -y_{1}x_{4}
& 0
\end{array}
\right)  \label{4P2familydeg14}%
\end{equation}
with $\left(  y_{1}:y_{2}:y_{2}\right)  \in\mathbb{P}^{2}$, and its general
element has the $49$
\index{double point}%
double points
\[
G\cdot\left(  0:y_{1}:y_{2}:y_{3}:-y_{3}:-y_{2}:-y_{1}\right)
\]

The induced
\index{group action}%
action on $\mathbb{P}^{6}$ is given by%
\[%
\begin{tabular}
[c]{lll}%
$\sigma\left(  x_{i}\right)  =x_{\left(  i+2\right)  \operatorname{mod}7}$ &
& $\tau\left(  x_{i}\right)  =\zeta^{2i}x_{i}$%
\end{tabular}
\]
Let $H=\left\langle \tau\right\rangle $ and consider the $\mathbb{P}^{1}$-%
\index{subfamily}%
subfamily $X_{s}=X_{\left(  0:1:s\right)  }$. The general element of $X_{s}$
has $56$
\index{double point}%
double points and $X_{s}$
\index{degeneration}%
degenerates for $s=0,\infty$ into a configuration of $14$ $\mathbb{P}^{3}$.

\begin{theorem}
\cite{Ro dland The Pfaffian CalabiYau its Mirror and their link to the
Grassmannian mathbbG27} The
\index{quotient}%
quotient of the general $X_{s}$ by $H$ has a
\index{crepant resolution}%
crepant resolution $\widetilde{X_{s}/H}$, and the Hodge numbers of
$\widetilde{X_{s}/H}$ coincide with the mirrored
\index{Hodge numbers}%
Hodge numbers of the general $X$.
\end{theorem}

R\o dland conjectured \cite{Ro dland The Pfaffian CalabiYau its Mirror and
their link to the Grassmannian mathbbG27} and Tj\o tta \cite{Tjo tta Quantum
cohomology of a Pfaffian CalabiYau variety: Verifying mirror symmetry
predictions} proved that the
\index{Picard-Fuchs equation}%
Picard-Fuchs equation of $\widetilde{X_{s}/H}$ at $s=\infty$ coincides with
the
\index{A-model}%
$A$-model of the general degree $14$ Pfaffian $X\subset\mathbb{P}^{6}$.

\section{Degenerations and mirror
symmetry\label{Sec Degenerations and mirror symmetry}}

Degenerations
\index{degeneration}%
to
\index{monomial degeneration}%
monomial ideals
\index{monomial ideal}%
in toric varieties play an important role in almost all known
\index{mirror construction}%
mirror constructions. For the concept of flat families%
\index{flat family}
see Section \ref{Groebnerbasesandflatness}.

\subsection{Degenerations associated to complete intersections in toric
varieties\label{Degenerationcompleteintersection}}

We want to associate to any
\index{complete intersection}%
complete intersection inside a toric variety $\mathbb{P}\left(  \Delta\right)
$, given by a
\index{nef partition}%
nef partition and represented by an ideal in the
\index{Cox ring}%
Cox ring $S$ of $\mathbb{P}\left(  \Delta\right)  $, a
\index{monomial degeneration}%
monomial
\index{degeneration}%
degeneration:

Suppose
\[
\Sigma\left(  1\right)  =I_{1}\cup...\cup I_{c}%
\]
i.e.,%
\[
-K_{Y}=\sum_{j=1}^{c}\overset{E_{j}}{\overbrace{\sum_{v\in I_{j}}D_{v}}}%
\]
is a
\index{nef partition}%
nef partition, i.e., all $E_{j}$ are
\index{T-Cartier divisor}%
Cartier,
\index{spanned by global sections}%
spanned by global sections.

\begin{example}
Consider $I\subset\mathbb{C}\left[  t\right]  \otimes S$ defined as%
\begin{align}
m_{j}  &  =\prod_{v\in I_{j}}y_{v}\label{4degencompleteintersection}\\
I_{0}  &  =\left\langle m_{j}\mid j=1,...,c\right\rangle \nonumber\\
I  &  =\left\langle f_{j}=t\cdot g_{j}+m_{j}\mid j=1,...,c\right\rangle
\subset\mathbb{C}\left[  t\right]  \otimes S\nonumber
\end{align}
where $g_{j}\in S_{\left[  E_{j}\right]  }$ corresponds to a general section
of $\mathcal{O}\left(  E_{j}\right)  $, i.e., a general linear combination of
the lattice points of $\Delta_{E_{j}}$ for $j=1,...,c$. Then $I$ defines a
\index{flat family}%
flat
\index{degeneration}%
degeneration $\mathfrak{X}\subset Y\times\operatorname*{Spec}\mathbb{C}\left[
\left[  t\right]  \right]  $ of a Calabi-Yau
\index{complete intersection}%
complete intersection in $Y=\mathbb{P}\left(  \Delta\right)  $, given by
general sections of $\mathcal{O}\left(  E_{1}\right)  ,...,\mathcal{O}\left(
E_{c}\right)  $, to the
\index{monomial degeneration}%
monomial
\index{special fiber}%
special fiber given by $I_{0}$.
\end{example}

We may assume that the $f_{j}$ are
\index{reduced standard basis}%
reduced with respect $I_{0}$ in the sense of Gr\"{o}bner bases. Flatness of
this family will be discussed in Section
\ref{degenerationtoriccompleteintersections}.

\begin{example}
In particular, a
\index{degeneration}%
degeneration $\mathfrak{X}$ of a Calabi-Yau
\index{hypersurface}%
hypersurface $X$ in $Y=\mathbb{P}\left(  \Delta\right)  $, defined by a
general section of $-K_{Y}$, to the monomial special fiber defined by
$\left\langle \prod_{v\in\Sigma\left(  1\right)  }y_{v}\right\rangle $ is
given by%
\[
I=\left\langle t\cdot\left\langle Am\mid m\in\partial\Delta\right\rangle
+\prod_{v\in\Sigma\left(  1\right)  }y_{v}\mid j=1,...,c\right\rangle
\subset\mathbb{C}\left[  t\right]  \otimes S
\]

\end{example}

\begin{example}
\label{22degen}The partitions for the above Example \ref{22example1} induce
\index{degeneration}%
degenerations given by the
\index{complete intersection}%
following ideals:

\begin{enumerate}
\item With
\index{elliptic curve!complete intersection}%
variables $x_{0},...,x_{3}$ of the
\index{Cox ring}%
Cox ring $S$ of $\mathbb{P}\left(  \Delta\right)  =\mathbb{P}^{3}$
corresponding to the vertices of $\Delta^{\ast}$ consider the ideal%
\[
I=\left\langle t\cdot g_{1}+x_{1}x_{2},\ t\cdot g_{2}+x_{0}x_{3}\right\rangle
\subset\mathbb{C}\left[  t\right]  \otimes S
\]

where $g_{1},g_{2}\in\mathbb{C}\left[  x_{1},...,x_{4}\right]  _{2}$ are
general not involving monomials in $I_{0}=\left\langle x_{1}x_{2},x_{0}%
x_{3}\right\rangle $. The ideal $I$ defines a flat degeneration $\mathfrak{X}%
\subset Y\times\operatorname*{Spec}\mathbb{C}\left[  \left[  t\right]
\right]  $ over $\operatorname*{Spec}\mathbb{C}\left[  \left[  t\right]
\right]  $ of an elliptic curve $X$ given as the complete intersection of two
quadrics in $\mathbb{P}^{3}$ to the monomial special fiber defined by $I_{0}$.

\item With variables $y_{1},...,y_{8}$ of the
\index{Cox ring}%
Cox ring $S^{\circ}$ of $Y^{\circ}=\mathbb{P}\left(  \nabla\right)  $
corresponding to the vertices of $\nabla^{\ast}$ consider the ideal%
\[%
\begin{tabular}
[c]{lll}%
$I^{\circ}=$ & $\langle t\cdot\left(  a_{1}\cdot y_{4}^{2}y_{8}^{2}+a_{2}\cdot
y_{3}^{2}y_{6}^{2}\right)  +y_{1}y_{2}y_{3}y_{4},$ & \\
& $t\cdot\left(  a_{3}\cdot y_{1}^{2}y_{5}^{2}+a_{4}\cdot y_{2}^{2}y_{7}%
^{2}\right)  +y_{5}y_{6}y_{7}y_{8}\rangle$ & $\subset\mathbb{C}\left[
t\right]  \otimes S^{\circ}$%
\end{tabular}
\
\]
with general coefficients $a_{i}$. The ideal $I^{\circ}$ defines a flat
degeneration $\mathfrak{X}\subset Y^{\circ}\times\operatorname*{Spec}%
\mathbb{C}\left[  \left[  t\right]  \right]  $ of the mirror $X^{\circ}$ of
$X$ to the monomial ideal
\[
I_{0}^{\circ}=\left\langle y_{1}y_{2}y_{3}y_{4},y_{5}y_{6}y_{7}y_{8}%
\right\rangle
\]
Note that the subvariety of $Y^{\circ}$ defined by the ideal $I_{0}^{\circ}$
decomposes into $4$ one-dimensional toric strata intersecting in $4$
zero-dimensional strata.

The stratification of the vanishing locus of reduced monomial ideals in the
Cox ring of a toric variety is explored in detail in Section
\ref{Monomial ideals in the Cox ring and the stratified toric primary decomposition}%
.
\end{enumerate}
\end{example}

\subsection{Degenerations of Pfaffian Calabi-Yau
manifolds\label{Sec degenerations of pfaffian calabi-yau varieties}}

Flatness of the following
\index{Pfaffian}%
Pfaffian
\index{degeneration}%
degenerations, which is obtained from the structure theorem of Buchsbaum and
Eisenbud \cite{BE Algebra structures for finite free resolutions and some
structure theorems for ideals of codimension 3}, is explored in Section
\ref{1PfaffianCalabiYauThreefolds}.

\begin{example}
\label{Ex Degeneration general Pfaffian elliptic curve}Let $S$ be the Cox ring
of $\mathbb{P}^{4}$, i.e., the homogeneous coordinate ring of $\mathbb{P}^{4}%
$. By the $4\times4$ Pfaffians in $\mathbb{C}\left[  t\right]  \otimes S$ of
the matrix
\[
A_{t}=t\cdot A+A_{0}%
\]
where%
\[
A_{0}=\left(
\begin{array}
[c]{ccccc}%
0 & 0 & x_{1} & -x_{4} & 0\\
0 & 0 & 0 & x_{2} & -x_{0}\\
-x_{1} & 0 & 0 & 0 & x_{3}\\
x_{4} & -x_{2} & 0 & 0 & 0\\
0 & x_{0} & -x_{3} & 0 & 0
\end{array}
\right)
\]
and $A$ is a general
\index{skew symmetric}%
skew symmetric $5\times5$ matrix linear in $x_{0},...,x_{4}$, we obtain a flat
degeneration $\mathfrak{X}\subset Y\times\operatorname*{Spec}\mathbb{C}\left[
\left[  t\right]  \right]  $ over $\operatorname*{Spec}\mathbb{C}\left[
\left[  t\right]  \right]  $ of a generic
\index{Pfaffian}%
Pfaffian
\index{elliptic curve!Pfaffian}%
elliptic curve in $\mathbb{P}^{4}$ to the monomial special fiber given by the
$4\times4$ Pfaffians of $A_{0}$.
\end{example}

Recall that for a skew symmetric matrix $A$ the determinants of the matrices
$A_{j}$ obtained by deleting the $j$-th row and column are squares, and the
$\sqrt{\det A_{j}}$ are called the
\index{Pfaffian}%
Pfaffians of $A$, for details see Section \ref{1PfaffianCalabiYauThreefolds}.

\begin{example}
\label{Ex degeneration 1 parameter Pfaffian degree 14}Let $H$ be the group
given in Section \ref{Roedlandexample}. The Cox ring of the quotient of
$\mathbb{P}^{6}$ by $H$ is a polynomial ring $S=\mathbb{C}\left[
x_{0},...,x_{6}\right]  $. By the $6\times6$
\index{Pfaffian}%
Pfaffians in $\mathbb{C}\left[  t\right]  \otimes S$ of
\[
A_{t}=\left(
\begin{array}
[c]{ccccccc}%
0 & tx_{1} & x_{2} & 0 & 0 & -x_{5} & -tx_{6}\\
-tx_{1} & 0 & tx_{3} & x_{4} & 0 & 0 & -x_{0}\\
-x_{2} & -tx_{3} & 0 & tx_{4} & x_{6} & 0 & 0\\
0 & -x_{4} & -tx_{4} & 0 & tx_{0} & x_{1} & 0\\
0 & 0 & -x_{6} & -tx_{0} & 0 & tx_{2} & x_{3}\\
x_{5} & 0 & 0 & -x_{1} & -tx_{2} & 0 & tx_{4}\\
tx_{6} & x_{0} & 0 & 0 & -x_{3} & -tx_{4} & 0
\end{array}
\right)
\]
a flat degeneration $\mathfrak{X}\subset\left(  \mathbb{P}^{6}/H\right)
\times\operatorname*{Spec}\mathbb{C}\left[  \left[  t\right]  \right]  $ over
$\operatorname*{Spec}\mathbb{C}\left[  \left[  t\right]  \right]  $ with
monomial special fiber is given.

This is the one
\index{Greene-Plesser}%
parameter family used to construct the mirror of a general degree $14$
Pfaffian Calabi-Yau threefold in $\mathbb{P}^{6}$ via
\index{Greene-Plesser orbifolding}%
Greene-Plesser orbifolding by $H$.
\end{example}

\begin{example}
\label{Ex Degeneration general Pfaffian degree 14}Let $S$ be the homogeneous
coordinate ring of $\mathbb{P}^{6}$. By the $6\times6$ Pfaffians in
$\mathbb{C}\left[  t\right]  \otimes S$ of
\[
A_{t}=t\cdot A+A_{0}%
\]
where%
\[
A_{0}=\left(
\begin{array}
[c]{ccccccc}%
0 & 0 & x_{2} & 0 & 0 & -x_{5} & 0\\
0 & 0 & 0 & x_{4} & 0 & 0 & -x_{0}\\
-x_{2} & 0 & 0 & 0 & x_{6} & 0 & 0\\
0 & -x_{4} & 0 & 0 & 0 & x_{1} & 0\\
0 & 0 & -x_{6} & 0 & 0 & 0 & x_{3}\\
x_{5} & 0 & 0 & -x_{1} & 0 & 0 & 0\\
0 & x_{0} & 0 & 0 & -x_{3} & 0 & 0
\end{array}
\right)
\]
and $A$ is a general
\index{skew symmetric}%
skew symmetric $7\times7$ matrix linear in $x_{0},...,x_{6}$, one obtains a
flat degeneration $\mathfrak{X}\subset Y\times\operatorname*{Spec}%
\mathbb{C}\left[  \left[  t\right]  \right]  $ over $\operatorname*{Spec}%
\mathbb{C}\left[  \left[  t\right]  \right]  $ of a general degree $14$
\index{Pfaffian}%
Pfaffian Calabi-Yau threefold in $\mathbb{P}^{6}$ to the
\index{monomial degeneration}%
monomial special fiber, given by the $6\times6$ Pfaffians of $A_{0}$.
\end{example}

\begin{example}
\label{Ex Degeneration general Pfaffian degree 13}Let $S$ be the homogeneous
coordinate ring of $\mathbb{P}^{6}$. The $5\times5$ Pfaffians in
$\mathbb{C}\left[  t\right]  \otimes S$ of $A_{t}=t\cdot A+A_{0}$, where%
\[
A_{0}=\left(
\begin{array}
[c]{ccccc}%
0 & 0 & x_{3}x_{4} & -x_{1}x_{2} & 0\\
0 & 0 & 0 & x_{7} & x_{6}\\
-x_{3}x_{4} & 0 & 0 & 0 & -x_{5}\\
x_{1}x_{2} & -x_{7} & 0 & 0 & 0\\
0 & -x_{6} & x_{5} & 0 & 0
\end{array}
\right)
\]
and $A$ is a general
\index{skew symmetric}%
skew symmetric map $\mathcal{E}^{\ast}\left(  -1\right)  \rightarrow
\mathcal{E}$ with
\[
\mathcal{E}=\mathcal{O}_{\mathbb{P}^{6}}\left(  1\right)  \oplus
\mathcal{O}_{\mathbb{P}^{6}}^{4}%
\]
one obtains a flat degeneration $\mathfrak{X}\subset Y\times
\operatorname*{Spec}\mathbb{C}\left[  \left[  t\right]  \right]  $ over
$\operatorname*{Spec}\mathbb{C}\left[  \left[  t\right]  \right]  $ of a
general degree $13$
\index{Pfaffian}%
Pfaffian
\index{elliptic curve!Pfaffian}%
Calabi-Yau threefold in $\mathbb{P}^{6}$ with
\index{monomial degeneration}%
monomial special fiber given by the $5\times5$ Pfaffians of $A_{0}$.

The
\index{special fiber}%
special fiber $X_{0}$ is obtained from a
\index{simplicial}%
simplicial $4$-polytope with $7$ vertices given in \cite{GS An Enumeration of
Simplicial 4Polytopes with 8 Vertices}.
\end{example}

\noindent For more monomial Calabi-Yau ideals obtained in this way see Section
\ref{1FurtherStanleyReisnerexamples}.

\section{Tropical geometry
ingredients\label{Sec tropical geometry ingredients}}

Tropical
\index{tropical geometry}%
geometry will be interpreted as a tool to explore
\index{one parameter}%
one parameter
\index{degeneration}%
degenerations inside toric varieties, as it associates to such a
\index{degeneration}%
degeneration a combinatorial object.

\subsection{Amoebas\label{Sec amoebas}}

\begin{definition}
Let $Y$ be a toric variety with
\index{torus}%
torus $\left(  \mathbb{C}^{\ast}\right)  ^{n}$ and $V\subset Y$ a subvariety.
The
\index{amoeba|textbf}%
\textbf{amoeba} of $V$ is \newsym[$\log$]{amoeba map}{}given as the image of
$V$ under
\begin{align*}
\log &  :\left(  \mathbb{C}^{\ast}\right)  ^{n}\rightarrow\mathbb{R}^{n}\\
\left(  z_{1},...,z_{n}\right)   &  \mapsto\left(  \log\left\vert
z_{1}\right\vert ,...,\log\left\vert z_{n}\right\vert \right)
\end{align*}

\end{definition}

\begin{remark}
The amoeba can be considered as a subset of a lower half sphere via%
\[%
\begin{tabular}
[c]{ccl}%
$\mathbb{R}^{n}$ & $\rightarrow$ & $S^{n}\cap\left\{  w_{t}\leq0\right\}  $\\
&  & $=\left\{  \left(  w_{t},w_{1},...,w_{n}\right)  \in\mathbb{R}^{n+1}\mid
w_{t}^{2}+w_{1}^{2}+...+w_{n}^{2}=1,\text{ }w_{t}\leq0\right\}  $\\
$\left(  w_{1},...,w_{n}\right)  $ & $\mapsto$ & $\frac{1}{\left\Vert \left(
-1,w_{1},...,w_{n}\right)  \right\Vert }\left(  -1,w_{1},...,w_{n}\right)  $%
\end{tabular}
\]

We refer to the points on the equator of the sphere, i.e., the points with
$w_{t}=0$, as the
\index{points at infinity|textbf}%
\textbf{points at infinity} of the amoeba.
\end{remark}

\begin{figure}
[h]
\begin{center}
\includegraphics[
height=2.0764in,
width=2.0764in
]%
{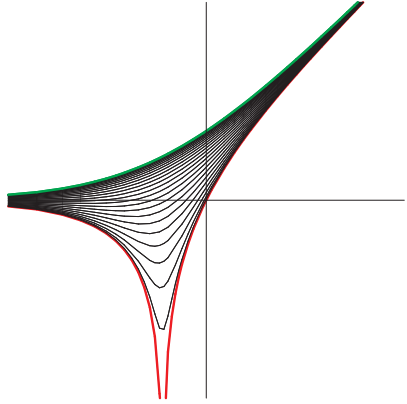}%
\caption{Amoeba of a line}%
\label{FigAmoebaLine}%
\end{center}
\end{figure}

\begin{figure}
[hh]
\begin{center}
\includegraphics[
height=2.0695in,
width=2.0695in
]%
{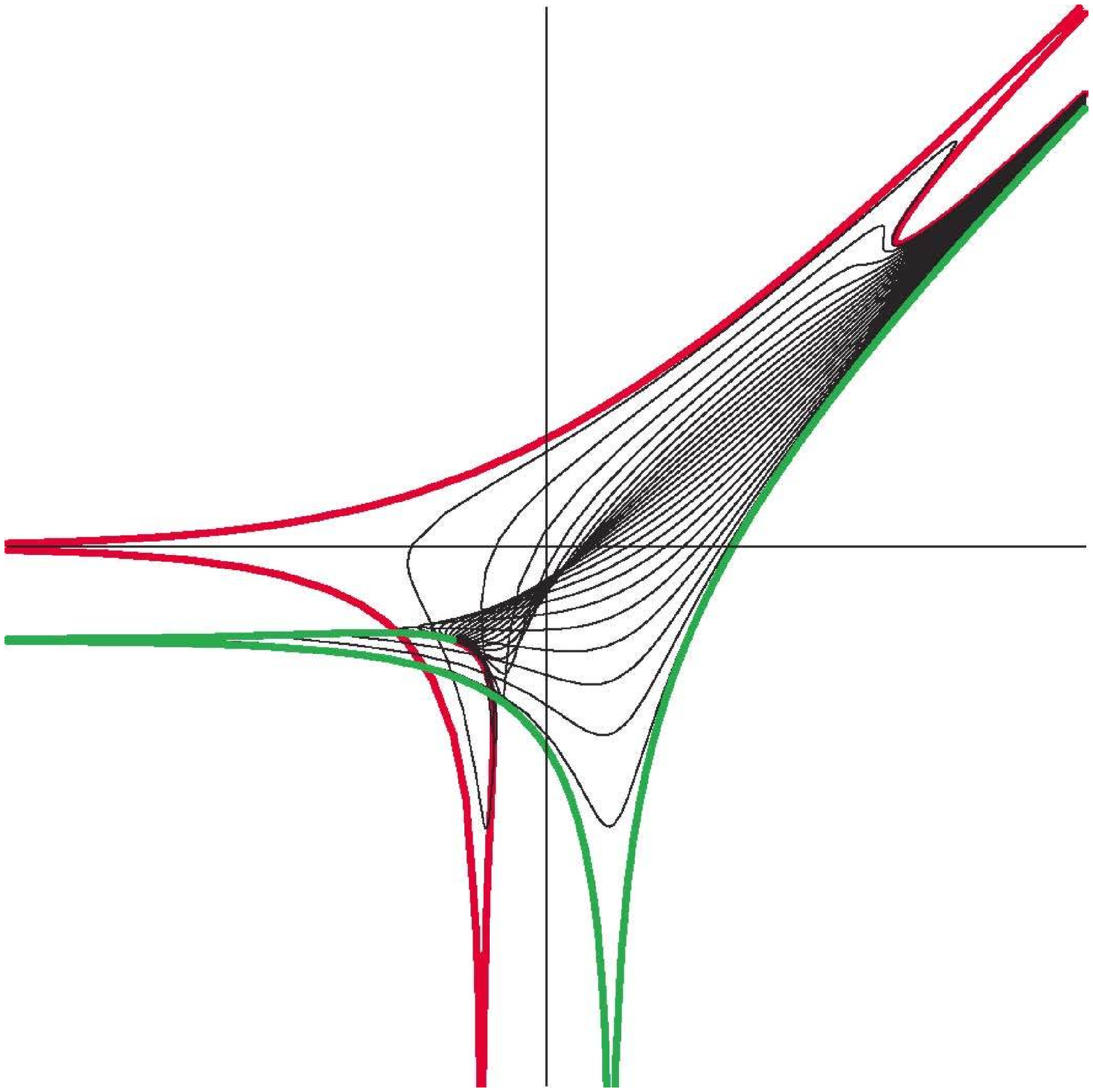}%
\caption{Amoeba of a conic}%
\label{FigAmoebaConic}%
\end{center}
\end{figure}

\begin{example}
\label{2amoebaslineconic}The amoeba of the line $L=\left\{  2x+y+1\right\}  $,
shown in Figure \ref{FigAmoebaLine}, is the image of%
\begin{align*}
&
\begin{tabular}
[c]{lllllll}%
$p:$ & $\mathbb{R}_{\geq0}\mathbb{\times}\left[  0,2\pi\right[  $ &
$\rightarrow$ & $\mathbb{C}$ & $\rightarrow$ & $\mathbb{C}^{2}$ &
$\rightarrow$\\
& $\left(  r,\varphi\right)  $ & $\mapsto$ & $re^{i\varphi}$ & $\mapsto$ &
$\left(  re^{i\varphi},-1-2re^{i\varphi}\right)  $ & $\mapsto$%
\end{tabular}
\medskip\\
&  \qquad\qquad\qquad%
\begin{tabular}
[c]{llll}%
$\rightarrow$ & $\mathbb{R}^{2}$ & $\rightarrow$ & $S^{2}\cap\left\{
w_{t}\leq0\right\}  $\\
$\mapsto$ & $\left(  \log r,\log\left\vert 1+2re^{i\varphi}\right\vert
\right)  $ & $\mapsto$ & $\frac{\left(  -1,\log r,\log\left\vert
1+2re^{i\varphi}\right\vert \right)  }{\left\Vert \left(  -1,\log
r,\log\left\vert 1+2re^{i\varphi}\right\vert \right)  \right\Vert }$%
\end{tabular}
\end{align*}
Considered as a subset of a lower half sphere via the last map the points at
infinity are%
\begin{align*}
\lim_{r\rightarrow0}p\left(  r,\varphi\right)   &  =\left(  0,1,0\right) \\
\lim_{r\rightarrow\infty}p\left(  r,\varphi\right)   &  =\left(  0,\frac
{1}{\sqrt{2}},\frac{1}{\sqrt{2}}\right) \\
\lim_{\left(  r,\varphi\right)  \rightarrow\left(  \frac{1}{2},\pi\right)
}p\left(  r,\varphi\right)   &  =\left(  0,1,0\right)
\end{align*}
The amoeba of the conic $\left\{  x^{2}+2y^{2}-3xy+x+y-1=0\right\}  $ is shown
in Figure \ref{FigAmoebaConic}.
\end{example}

%


\begin{figure}
[h]
\begin{center}
\includegraphics[
height=1.9726in,
width=1.9726in
]%
{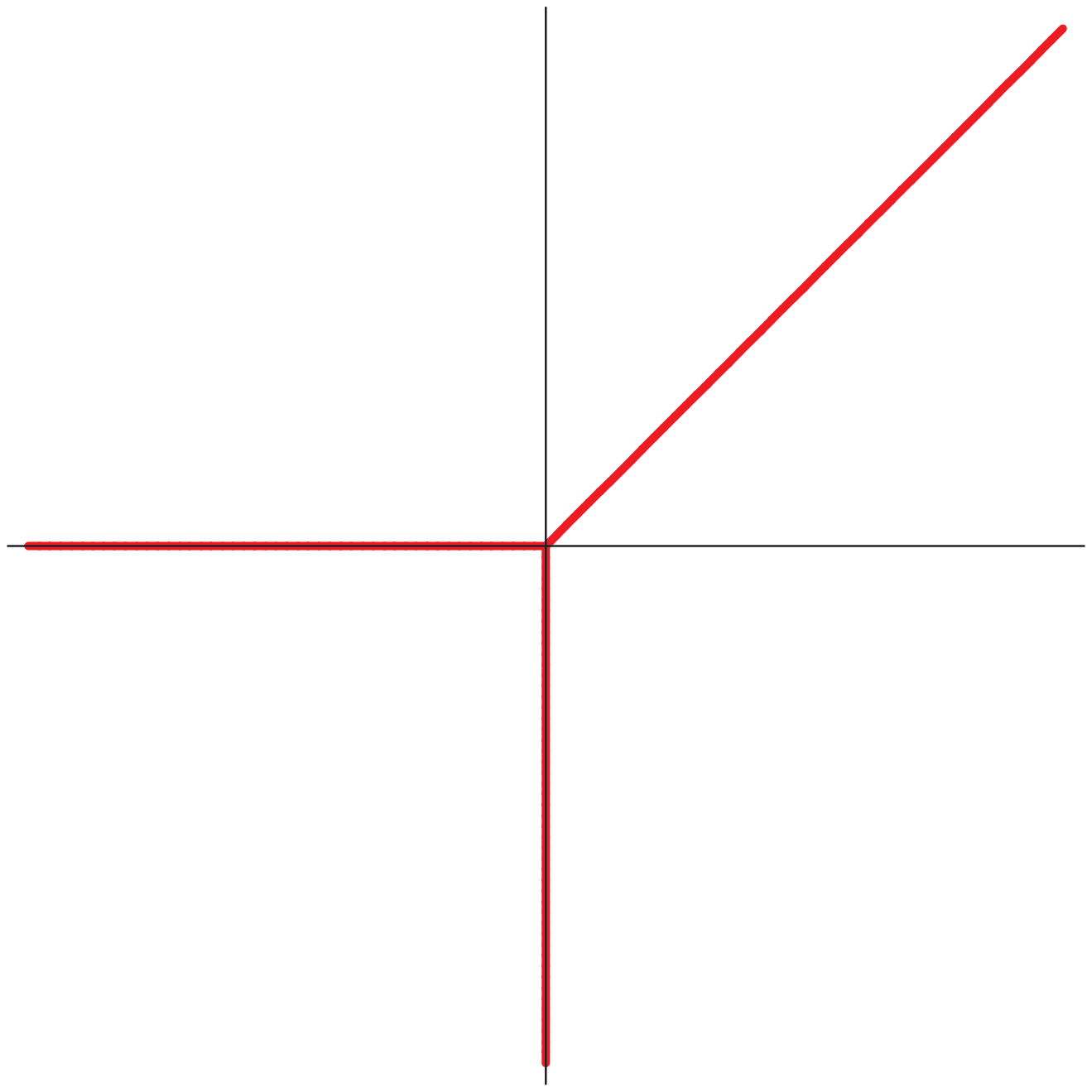}%
\caption{Limit amoeba}%
\label{FigLimitAmoeba}%
\end{center}
\end{figure}

\begin{example}
Replacing $\log$ by $\log_{t}$ the
\newsym[$\log_{t}$]{amoeba map to the base $t$}{}amoeba is rescaled, Figure
\ref{FigLimitAmoeba} shows the limit $t\rightarrow\infty$ of both amoebas
given in Example \ref{2amoebaslineconic}. For the conic one has to assign
multiplicity $2$ to each leg.
\end{example}

%

This limit process can be formalized in the following way:

\subsection{Non-Archimedian amoebas\label{1NonArchimedianAmoebas}}

Consider the field of Puisseux series $\overline{\mathbb{C}\left(  t\right)
}$, which is equipped with a \newsym[$val$]{valuation}{}valuation%
\begin{gather*}
val:\overline{\mathbb{C}\left(  t\right)  }\rightarrow\mathbb{Q\cup}\left\{
\infty\right\} \\
\sum_{j\in J}\alpha_{j}t^{j}\mapsto\min J
\end{gather*}
satisfying $val\left(  f+g\right)  \geq\min\left\{  val\left(  f\right)
,val\left(  g\right)  \right\}  $, and with a norm $\left\Vert f\right\Vert
=e^{-val\left(  f\right)  }$. Consider further the metric completion $K$ of
$\overline{\mathbb{C}\left(  t\right)  }$ containing those elements
$\sum_{j\in J}\alpha_{j}t^{j}$, which satisfy the condition that any subset of
$J$ has a minimum. Denote the corresponding valuation and norm on $K$ again by
$val$ and $\left\Vert -\right\Vert $. $K$ is a complete algebraically closed
non-Archimedian field with surjective
\index{valuation}%
valuation
\[
val:K\rightarrow\mathbb{R}\cup\left\{  \infty\right\}
\]
The term non-Archimedean means that the norm on $K$ satisfies the inequality%
\[
\left\Vert f+g\right\Vert \leq\max\left\{  \left\Vert f\right\Vert ,\left\Vert
g\right\Vert \right\}
\]
for all $f,g\in K$. This in particular implies that the Archimedian axiom is
not satisfied. If $f,g\in K$ with $\left\Vert f\right\Vert <\left\Vert
g\right\Vert $, then for all natural number $n$ we have $\left\Vert n\cdot
f\right\Vert \leq\left\Vert f\right\Vert <\left\Vert g\right\Vert $, indeed
$\left\Vert n\cdot f\right\Vert =\left\Vert f\right\Vert $.

Let $I$ be an ideal in $K\left[  x_{1},...,x_{n}\right]  $ and $V_{K}\left(
I\right)  $ be the
\newsym[$V_{K}\left(  I\right)  $]{algebraic variety in $\left(  K^{\ast}\right)  ^{n}$ given by $I$}{}algebraic
variety given by $I$ in $\left(  K^{\ast}\right)  ^{n}$.

As the norm is given by $\left\Vert -\right\Vert =e^{-val\left(  -\right)  }$,
the corresponding
\newsym[$\operatorname*{val}$]{non-Archimedian amoeba map}{}amoeba map
$\log\left\Vert -\right\Vert $ over $K$ is given by the valuations
\begin{gather*}
\operatorname*{val}\nolimits_{-}=\log\left\Vert -\right\Vert :\left(  K^{\ast
}\right)  ^{n}\rightarrow\mathbb{R}^{n}\\
\left(  z_{1},...,z_{n}\right)  \mapsto\left(  -val\left(  z_{1}\right)
,...,-val\left(  z_{n}\right)  \right)
\end{gather*}

\begin{definition}
The
\index{non-Archimedian amoeba|textbf}%
\textbf{non-Archimedian amoeba} of $V_{K}\left(  I\right)  $ is
$\operatorname*{val}\nolimits_{-}\left(  V_{K}\left(  I\right)  \right)  $.
\end{definition}

A proof of the following theorem for hypersurfaces can be found in\linebreak%
\cite[Sec. 6.1]{GKZ Discriminants Resultants and Multidimensional
Determinants}, the general statement in terms of the
\index{Bergman fan}%
Bergman fan (see Section \ref{1tropicalvarietiesandtheBergmanfan}) in
\cite[Sec. 9.4]{Sturmfels Solving Systems of Polynomial Equations}.

\begin{theorem}
The limit $\lim_{t\rightarrow\infty}\log_{t}V\left(  I_{t}\right)  $ exists as
the limit in the Hausdorff metric on compacts, and%
\[
\operatorname*{val}\nolimits_{-}\left(  V_{K}\left(  I\right)  \right)
=\lim_{t\rightarrow\infty}\log_{t}V\left(  I_{t}\right)
\]

\end{theorem}

Recall that the distance of two closed subsets of a metric space in the
\index{Hausdorff metric|textbf}%
Hausdorff metric is given by%
\[
d\left(  A,B\right)  =\max\left\{  \sup_{a\in A}d\left(  a,B\right)
,\sup_{b\in B}d\left(  A,b\right)  \right\}
\]
so above convergence means that for any compact $D\subset\mathbb{R}^{n}$
\[
\lim_{t\rightarrow\infty}d\left(  D\cap\log_{t}V\left(  I_{t}\right)
,D\cap\operatorname*{val}V_{K}\left(  I\right)  \right)  =0
\]
From the point of view of degenerations it will turn out to be more natural to
consider the image
\[
\operatorname*{val}\left(  V_{K}\left(  I\right)  \right)
=-\operatorname*{val}\nolimits_{-}\left(  V_{K}\left(  I\right)  \right)
\]
of $V_{K}\left(  I\right)  $ under the map%
\begin{gather*}
\operatorname*{val}:\left(  K^{\ast}\right)  ^{n}\rightarrow\mathbb{R}^{n}\\
\left(  z_{1},...,z_{n}\right)  \mapsto\left(  val\left(  z_{1}\right)
,...,val\left(  z_{n}\right)  \right)
\end{gather*}
associating to each component the minimal weight term.

\begin{figure}
[h]
\begin{center}
\includegraphics[
trim=0.000000in -0.015178in 0.000000in 0.015179in,
height=2.4163in,
width=2.4163in
]%
{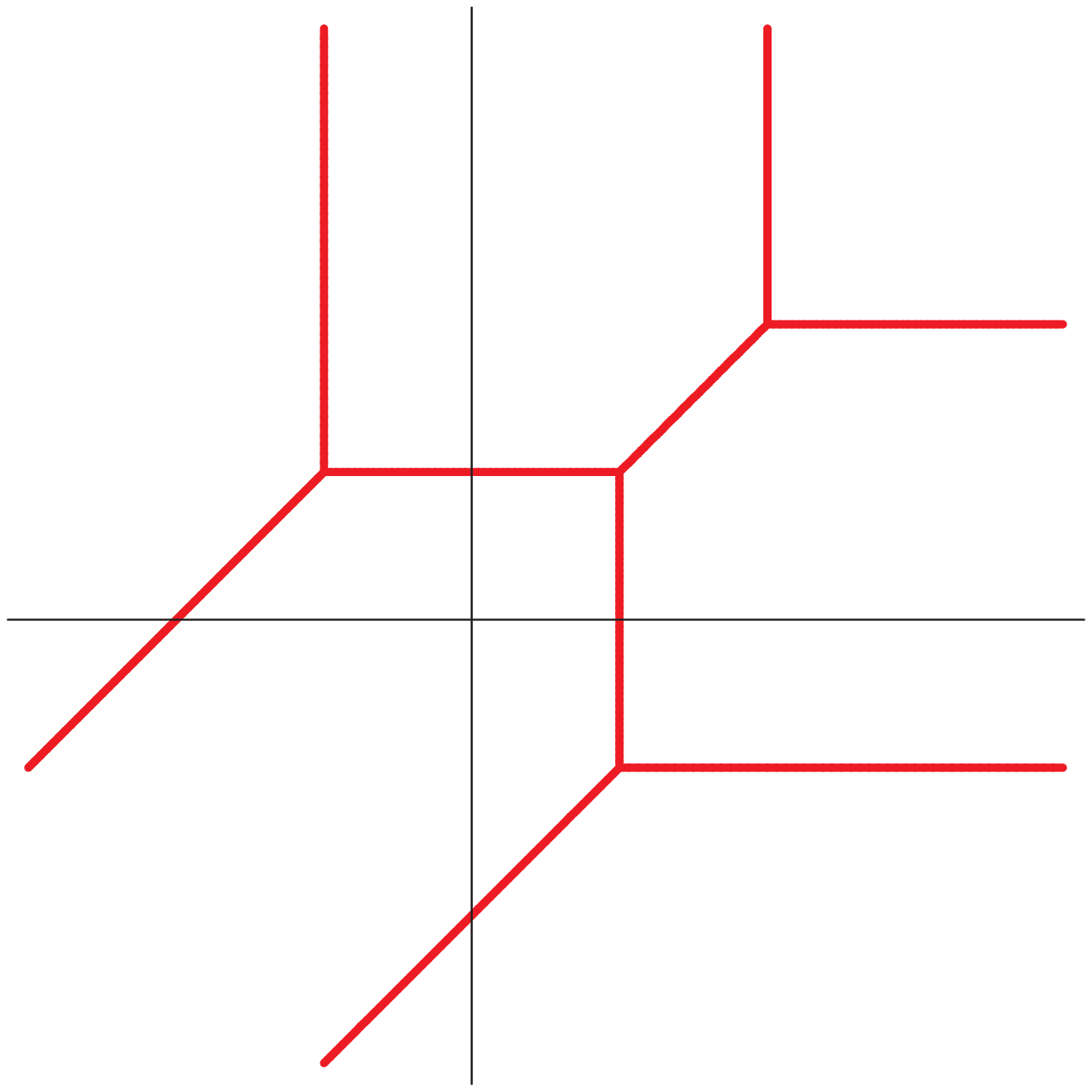}%
\caption{$\operatorname*{val}\left(  V_{K}\left(  I\right)  \right)  $ for the
ideal of a plane quadric with coefficients in $K$}%
\label{Figtropconic}%
\end{center}
\end{figure}

\begin{example}
$\operatorname*{val}\left(  V_{K}\left(  I\right)  \right)  $ for the ideal of
a plane quadric with coefficients in $K$ is shown in Figure \ref{Figtropconic}.
\end{example}

%

\begin{figure}
[h]
\begin{center}
\includegraphics[
height=1.9311in,
width=1.9311in
]%
{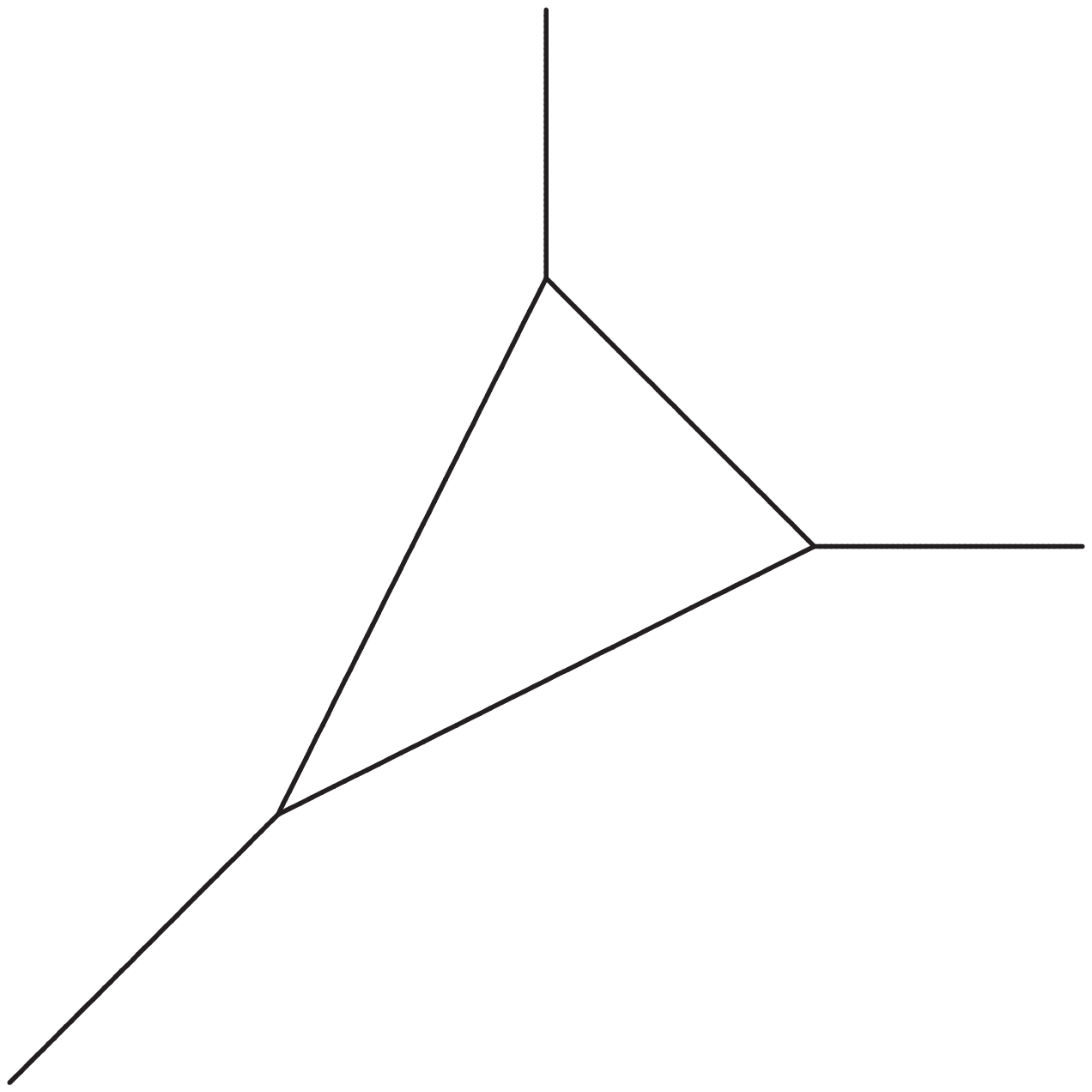}%
\caption{$\operatorname*{val}\left(  V_{K}\left(  I\right)  \right)  $ for the
degeneration of plane cubics}%
\label{FigtropCubicPlaneDegen}%
\end{center}
\end{figure}

\begin{example}
\label{1planecubic}$\operatorname*{val}\left(  V_{K}\left(  I\right)  \right)
$ for the
\index{monomial degeneration}%
degeneration
\index{hypersurface}%
of
\index{degeneration}%
general
\index{plane cubic}%
plane cubics%
\[
\left\{  x_{0}x_{1}x_{2}+tf=0\right\}
\]
with $f$ a general element in $\mathbb{C}\left[  x_{0},x_{1},x_{2}\right]
_{3}$ is shown in Figure \ref{FigtropCubicPlaneDegen}. Note that for an ideal
$I$, homogeneous with respect to the
\index{grading}%
grading $\deg x_{i}=1$ $\forall i$ on $K\left[  x_{1},...,x_{n}\right]  $, one
can consider $\operatorname*{val}\left(  V_{K}\left(  I\right)  \right)  \ $as
a subset of $\frac{\mathbb{R}^{3}}{\mathbb{R}\left(  1,1,1\right)  }%
\cong\mathbb{R}^{2}$, this will be explored in detail in Section
\ref{1torichomogeneoussetting}.
\end{example}

Having made the geometric connection between%
degenerations and tropical geometry, we consider the algebraic connection:

\subsection{Tropical varieties\label{Sec tropical varieties}}

\begin{definition}
A
\index{tropical variety|textbf}%
\textbf{tropical variety} is a
\newsym[$\operatorname*{tropvar}\left(  I\right)  $]{tropical variety}{}subset
\[
\operatorname*{tropvar}\left(  I\right)  =\operatorname*{val}\nolimits_{-}%
\left(  V_{K}\left(  I\right)  \right)  \subset\mathbb{R}^{n}%
\]
where $I$ is an ideal in $K\left[  x_{1},...,x_{n}\right]  $.
\end{definition}

For $w\in\mathbb{R}^{n}$ the
\index{initial form|textbf}%
\textbf{initial form} $in_{w}\left(  f\right)  $ of $f\in K\left[
x_{1},...,x_{n}\right]  $ is
\newsym[$in_{w}\left(  f\right)  $]{initial form}{}the sum of the terms of
maximal weight with respect to $w$ and $weight\left(  c\right)  =-val\left(
c\right)  $ for $c\in K$. For any ideal $J\subset K\left[  x_{1}%
,...,x_{n}\right]  $ \newsym[$in_{w}\left(  J\right)  $]{initial ideal}{}its
\index{initial ideal|textbf}%
\textbf{initial ideal} is%
\[
in_{w}\left(  J\right)  =\left\langle in_{w}\left(  f\right)  \mid f\in
J\right\rangle
\]

\begin{theorem}
\label{thm tropical variety properties}\cite[Sec. 2]{RGST First Steps in
Tropical Geometry}, \cite[Sec. 9.2]{Sturmfels Solving Systems of Polynomial
Equations}, \cite[Sec. 2]{SpS The tropical Grassmannian} Every ideal $I\subset
K\left[  x_{1},...,x_{n}\right]  $ has a finite subset $\mathcal{G}$, called
a
\index{tropical basis|textbf}%
\textbf{tropical basis} of $I$, such that

\begin{enumerate}
\item For all $w\in\operatorname*{tropvar}\left(  I\right)  $ the set
$\left\{  in_{w}\left(  g\right)  \mid g\in\mathcal{G}\right\}  $ generates
$in_{w}\left(  I\right)  $.

\item For all $w\notin\operatorname*{tropvar}\left(  I\right)  $ the set
$\left\{  in_{w}\left(  g\right)  \mid g\in\mathcal{G}\right\}  $ contains a monomial.
\end{enumerate}

$\operatorname*{tropvar}\left(  I\right)  $ is a finite intersection of the
tropical hypersurfaces
\index{tropvar}%
$\operatorname*{tropvar}\left\langle \left(  g\right)  \right\rangle $ for
$g\in\mathcal{G}$, it is a
\index{polyhedral cell complex}%
polyhedral cell complex, its dimension is the
\index{Krull dimension}%
Krull dimension of $\frac{K\left[  x_{1},...,x_{n}\right]  }{I}$, it is
equidimensional if $V_{K}\left(  I\right)  $ is, and%
\[
\operatorname*{tropvar}\left(  I\right)  =\left\{  w\in\mathbb{R}^{n}\mid
in_{w}\left(  I\right)  \text{ contains no monomial}\right\}
\]

\end{theorem}

Selecting the maximal weight term and defining the weight of a constant $c\in
K$ as $weight\left(  c\right)  =-val\left(  c\right)  $ is the Gr\"{o}bner
basis point of view. With respect to degenerations it is more natural to look
at the minimal weight term and take $weight\left(  c\right)  =val\left(
c\right)  $ for $c\in K$, i.e., to consider $\operatorname*{val}\left(
V_{K}\left(  I\right)  \right)  =-\operatorname*{val}\nolimits_{-}\left(
V_{K}\left(  I\right)  \right)  $.

\begin{figure}
[h]
\begin{center}
\includegraphics[
height=2.514in,
width=2.5054in
]%
{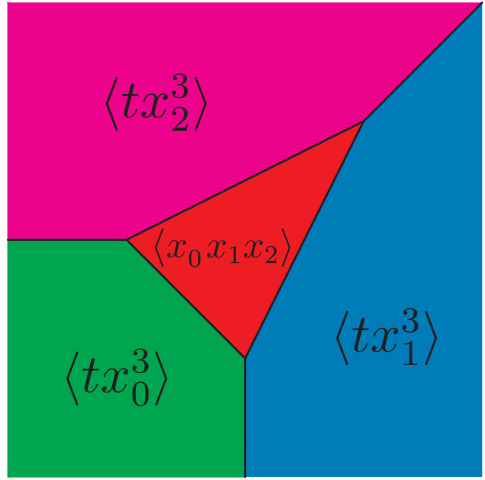}%
\caption{Tropical variety and initial ideals for the degeneration of plane
cubics}%
\label{FigInitialIdealsCubicDegen}%
\end{center}
\end{figure}

\begin{example}
The
\index{hypersurface}%
monomial
\index{initial ideal}%
initial ideals and the
\index{elliptic curve!plane}%
sets of
\index{weights vector}%
weight vectors leading to them for the plane cubic case as in Example
\ref{1planecubic} are depicted in Figure \ref{FigInitialIdealsCubicDegen}. For
$w\in\operatorname*{tropvar}\left(  I\right)  $ the
\index{initial ideal}%
initial ideal is generated by a sum of the initial terms appearing in a
neighborhood of $w$.
\end{example}

%

\subsection{Tropical prevarieties\label{Sec tropical prevarieties}}

\begin{definition}
The
\index{tropical semiring|textbf}%
\textbf{tropical semiring} is $\mathbb{R\cup}\left\{  -\infty\right\}  $ with
\index{tropical addition|textbf}%
\textbf{tropical addition} and
\index{tropical multiplication|textbf}%
\textbf{multiplication}%
\begin{align*}
a\oplus b  &  =\max\left(  a,b\right) \\
a\odot b  &  =a+b
\end{align*}

\end{definition}

The tropical semiring satisfies $\left(  a\oplus b\right)  \odot c=a\odot
c\oplus b\odot c$, the additive unit is $-\infty$, the multiplicative unit is
$0$. In general there is no additive inverse in the tropical semiring.

\begin{definition}
A
\index{tropical polynomial|textbf}%
\textbf{tropical polynomial} is a polynomial formed with $\oplus$ and $\odot$,
i.e., a
\index{piecewise linear}%
piecewise linear function%
\[%
\begin{tabular}
[c]{lll}%
$F:\mathbb{R}^{n}\rightarrow\mathbb{R}$ &  & $F\left(  x_{1},...,x_{n}\right)
=\max\left\{  a_{1j}x_{1}+...+a_{nj}x_{n}+c_{j}\mid j\right\}  $%
\end{tabular}
\]

\end{definition}

\begin{definition}
The
\index{tropical prevariety|textbf}%
\textbf{tropical prevariety} $T\left(  F\right)  $ of $F$ is the set where the
maximum is attained at least twice, and $T\left(  \mathcal{G}\right)
=\bigcap_{g\in\mathcal{G}}T\left(  g\right)  $ for any
\newsym[$T\left(  \mathcal{G}\right)  $]{tropical prevariety}{}set of tropical
polynomials $\mathcal{G}$.
\end{definition}

\begin{definition}
For any polynomial $f\in K\left[  x_{1},...,x_{n}\right]  $%
\[
f=\sum\nolimits_{a}b_{a}\left(  t\right)  \cdot x^{a}%
\]
define
\newsym[$\operatorname*{trop}\left(  f\right)  $]{tropicalization}{}its
\index{tropicalization of a polynomial|textbf}%
\textbf{tropicalization} as%
\[
\operatorname*{trop}\left(  f\right)  =\bigoplus\nolimits_{a}-val\left(
b_{a}\left(  t\right)  \right)  \odot x^{\odot a}%
\]

\end{definition}

So, consistent with the amoeba, the non-Archimedian amoeba and the tropical
variety, we again adopt the Gr\"{o}bner basis point of view, looking at the
maximal weight term and take $weight\left(  c\right)  =-val\left(  c\right)  $
for $c\in K$.

\begin{theorem}
\cite[Sec. 2]{RGST First Steps in Tropical Geometry}, \cite[Sec.
9.2]{Sturmfels Solving Systems of Polynomial Equations}, \cite[Sec. 2]{SpS The
tropical Grassmannian} Any tropical variety $\operatorname*{tropvar}\left(
I\right)  $, $I\subset K\left[  x_{1},...,x_{n}\right]  $ is a tropical
prevariety. For any ideal $I\subset K\left[  x_{1},...,x_{n}\right]  $
\[
\operatorname*{tropvar}\left(  I\right)  =T\left(  \operatorname*{trop}\left(
I\right)  \right)
\]

\end{theorem}

\begin{figure}
[h]
\begin{center}
\includegraphics[
height=1.6561in,
width=2.4578in
]%
{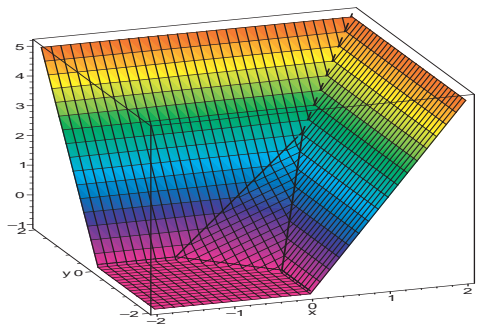}%
\caption{Piecewise linear function associated to the degeneration of plane
cubics}%
\label{FigPWLcubic}%
\end{center}
\end{figure}

\begin{example}
For
\index{hypersurface}%
the general plane
\index{elliptic curve!plane}%
elliptic curve in Example \ref{1planecubic} the tropical variety
$\operatorname*{tropvar}\left\langle f\right\rangle $ is the non
differentiability locus $T\left(  F\right)  $ of the piecewise linear function%
\[
F=\max\left\{  3x_{1}-1,2x_{1}+x_{2}-1,x_{1}+2x_{2}-1,3x_{2}-1,...,-1,x_{1}%
+x_{2}\right\}
\]
Figure \ref{FigPWLcubic} shows the graph of $F$.
\end{example}

%

\begin{figure}
[h]
\begin{center}
\includegraphics[
height=2.0159in,
width=2.3549in
]%
{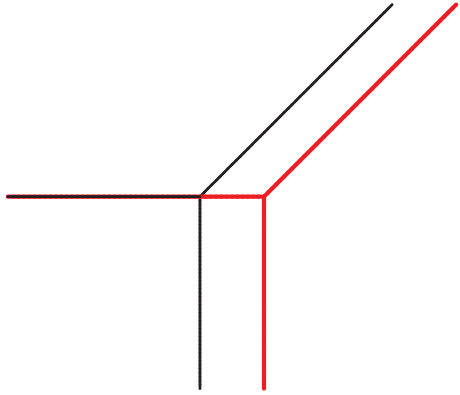}%
\caption{Intersecting two tropical lines}%
\label{FigIntersectionTropicalLines}%
\end{center}
\end{figure}

Not every tropical prevariety is a tropical variety:

\begin{example}
The intersection of the
\index{tropical line}%
tropical lines $L_{1}=T\left(  \operatorname*{trop}\left(  x+y+1\right)
\right)  $ and $L_{2}=T\left(  \operatorname*{trop}\left(  tx+y+1\right)
\right)  $, as depicted in Figure \ref{FigIntersectionTropicalLines}, is a
tropical prevariety, but not a
\index{tropical variety}%
tropical variety.
\end{example}

%

\subsection{Tropical varieties and the Bergman
fan\label{1tropicalvarietiesandtheBergmanfan}}

Let $I$ be an ideal in $\mathbb{C}\left[  t,x_{1},...,x_{n}\right]  $.

\begin{theorem}
\cite{Bergman The logarithmic limitset of an algebraic variety}, \cite{BG The
geometry of the set of characters induced by valuations}, \cite[Sec.
9.3]{Sturmfels Solving Systems of Polynomial Equations} Define
\[
suppBC_{-}\left(  I\right)  =\left\{
\begin{array}
[c]{c}%
p\in S^{n}\mid\exists\text{ sequence }\left(  p_{j}\right)  _{j\in\mathbb{N}%
}\text{ with }p_{j}\in\log\left(  V\left(  I\right)  \right)  \cap
jS^{n}\subset\mathbb{R}^{n+1}\\
\text{and }\lim_{j\rightarrow\infty}\frac{1}{j}p_{j}=p
\end{array}
\right\}
\]
and%
\[
suppBF_{-}\left(  I\right)  =\left\{  p\in\mathbb{R}^{n+1}\backslash\left\{
0\right\}  \mid\frac{p}{\left\Vert p\right\Vert }\in suppBC_{-}\left(
I\right)  \right\}  \cup\left\{  0\right\}
\]
If $V\left(  I\right)  $ is an
\newsym[$suppBC\left(  I\right)  $]{support of the Bergman complex}{}irreducible
\newsym[$suppBF\left(  I\right)  $]{support of the Bergman fan}{}subvariety
\index{suppBC|textbf}%
of $\left(  \mathbb{C}^{\ast}\right)  ^{n+1}$ of dimension $d+1$, then
\index{suppBF|textbf}%
$suppBF_{-}\left(  I\right)  $ is a finite union $d+1$-dimensional convex
polyhedral cones. The intersection of any two is a common face.

Denote by
\index{BF}%
$BF_{-}\left(  I\right)  $ the corresponding fan and by
\index{BC}%
$BC_{-}\left(  I\right)  $ the corresponding complex of dimension $d$.
\end{theorem}

Note:

\begin{itemize}
\item $V\left(  I\right)  \subset\left(  \mathbb{C}^{\ast}\right)  ^{n+1}$

\item These definitions are consistent with the Gr\"{o}bner basis $\left(
\max,+\right)  $ point of view looking at the maximum weight terms.

\item The definitions of $BF_{-}\left(  I\right)  $ and $BC_{-}\left(
I\right)  $ are symmetric in all variables $t,x_{1},...,x_{n}$.

\item The complex $BC_{-}\left(  I\right)  $ and the fan $BF_{-}\left(
I\right)  $ are known in the literature as Bergman complex and Bergman fan
respectively. However with degenerations in mind, i.e., the power series point
of view, it is more natural to consider the reflection of these objects at the
origin. So we use the following definition:
\end{itemize}

\begin{definition}
Analogous to $K$ denote by $L$ the metric completion of the field
$\overline{\mathbb{C}\left(  s\right)  }$ of Puisseux series in a new variable
$s$. If $I$ is an ideal in $\mathbb{C}\left[  t,x_{1},...,x_{n}\right]  $,
then
\index{Bergman fan|textbf}%
\textbf{Bergman fan }$BF\left(  I\right)  $ of $I$ is the image of the
vanishing locus of $V_{L}\left(  I\right)  $ of $I$ over $L$ under the map
\begin{gather*}
\left(  L^{\ast}\right)  ^{n+1}\rightarrow\mathbb{R}^{n+1}\\
\left(  t,x_{1},...,x_{n}\right)  \mapsto\left(  val\left(  t\right)
,val\left(  x_{1}\right)  ,...,val\left(  x_{n}\right)  \right)
\end{gather*}
The intersection of $BF\left(  I\right)  $ with the unit sphere $S^{n}$ is
called the
\index{Bergman complex|textbf}%
\textbf{Bergman complex} $BC\left(  I\right)  $ of $I$.

Note that this non-Archimedian type definition has the advantage that it
avoids problems with limit processes.
\end{definition}

If you prefer the $\left(  \max,+\right)  $ point of view you may replace in
this definiton $val$ by $-val$.

\begin{figure}
[h]
\begin{center}
\includegraphics[
height=1.4174in,
width=1.5446in
]%
{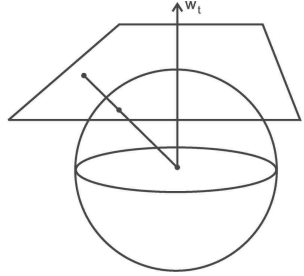}%
\caption{Stereographic projection relating the Bergman complex and the
tropical variety}%
\label{FigStereographicProjection}%
\end{center}
\end{figure}

Consider the
\index{stereographic projection}%
stereographic projection $\pi$, visualized in Figure
\ref{FigStereographicProjection}, of the upper half sphere
\[
S^{n}\cap\left\{  w_{t}>0\right\}  =\left\{  \left(  w_{t},w_{x_{1}%
},...,w_{x_{n}}\right)  \in\mathbb{R}^{n+1}\mid w_{t}^{2}+w_{x_{1}}%
^{2}+...+w_{x_{n}}^{2}=1,\text{ }w_{t}>0\right\}
\]
from $0$ to $\mathbb{R}^{n}=\left\{  w_{t}=1\right\}  $.%
Here we denote the coordinates of $\mathbb{R}^{n+1}$ corresponding to the
variables of $\mathbb{C}\left[  t,x_{1},...,x_{n}\right]  $ by $w_{t}%
,w_{x_{1}},...,w_{x_{n}}$, as they are weights on the variables.

In the same way denote by $\pi_{-}$ the stereographic projection of the lower
half sphere from $0$ to $\mathbb{R}^{n}=\left\{  w_{t}=-1\right\}  $.

Connecting the
\index{Bergman complex}%
Bergman complex to the
\index{tropical variety}%
tropical variety via $\pi$ (see \cite[Sec. 9.4]{Sturmfels Solving Systems of
Polynomial Equations}) and summarizing:

\begin{theorem}
For any ideal $I$ in $\mathbb{C}\left[  t,x_{1},...,x_{n}\right]  $ it holds%
\begin{align*}
\lim_{t\rightarrow\infty}\left(  \log_{t}V\left(  I_{t}\right)  \right)   &
=\operatorname*{val}\nolimits_{-}\left(  V_{K}\left(  I\right)  \right)
=\operatorname*{tropvar}\left(  I\right)  =T\left(  \operatorname*{trop}%
\left(  I\right)  \right) \\
&  =\pi_{-}\left(  BC_{-}\left(  I\right)  \cap\left\{  w_{t}<0\right\}
\right)  \subset\mathbb{R}^{n}%
\end{align*}
If $I\subset\mathbb{C}\left[  x_{1},...,x_{n}\right]  $, then $BF_{-}\left(
I\right)  \subset\mathbb{R}^{n}$ coincides with the above when considering $I$
as an ideal in $\mathbb{C}\left[  t,x_{1},...,x_{n}\right]  $.
\end{theorem}

\begin{remark}
Reflecting at the origin, we have
\begin{align*}
\operatorname*{val}\left(  V_{K}\left(  I\right)  \right)   &  =\pi\left(
BC\left(  I\right)  \cap\left\{  w_{t}>0\right\}  \right) \\
&  =-\lim_{t\rightarrow\infty}\left(  \log_{t}V\left(  I_{t}\right)  \right)
=-\operatorname*{tropvar}\left(  I\right)  =-T\left(  \operatorname*{trop}%
\left(  I\right)  \right)
\end{align*}
Our non-Archimedian definition of the Bergman fan relates to the limit
definition by%
\begin{align*}
BC\left(  I\right)   &  =-BC_{-}\left(  I\right) \\
BF\left(  I\right)   &  =-BF_{-}\left(  I\right)
\end{align*}
Going from the Bergman complex $BC\left(  I\right)  $ to $\operatorname*{val}%
\left(  V_{K}\left(  I\right)  \right)  $, i.e., intersecting with the plane
$\left\{  w_{t}=1\right\}  $, amounts to the identification of the parameters
$s$ and $t$.
\end{remark}

For the subset of $BC\left(  I\right)  $ lying inside the equator $\left\{
w_{t}=0\right\}  $ of the sphere, we introduce the notation:

\begin{definition}
$BC\left(  I\right)  \cap\left\{  w_{t}=0\right\}  $ is called the
\index{tropical variety at infinity|textbf}%
\textbf{tropical variety at infinity}.
\end{definition}

%

\begin{figure}
[h]
\begin{center}
\includegraphics[
height=1.6034in,
width=2.5019in
]%
{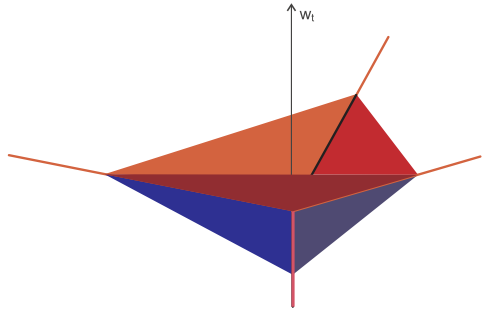}%
\caption{Bergman fan of the degeneration of plane elliptic curves}%
\label{FigBFplaneelliptic}%
\end{center}
\end{figure}
\begin{figure}
[hptb]
\begin{center}
\includegraphics[
height=2.1897in,
width=2.1577in
]%
{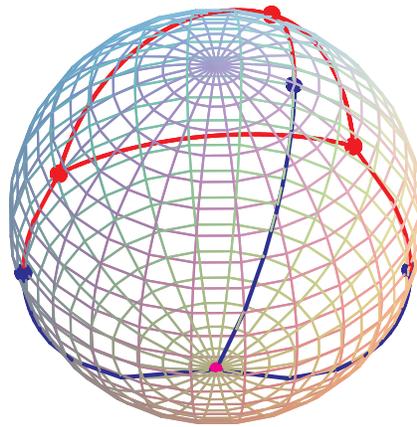}%
\caption{Bergman complex of the degeneration of plane elliptic curves}%
\label{FigBCplaneEllipticDegen}%
\end{center}
\end{figure}

\begin{example}
For the
\index{elliptic curve!plane}%
plane elliptic curve in
\index{hypersurface}%
above Example \ref{1planecubic} the Bergman fan $BF\left(  I\right)  $ is
shown in Figure \ref{FigBFplaneelliptic} (extending the depicted faces to
infinity). Applying $\pi$ to the $w_{t}>0$ part of Figure
\ref{FigBCplaneEllipticDegen}, which is visualizing the Bergman complex
$BC\left(  I\right)  $, gives $\operatorname*{val}V_{K}\left(  I\right)  $.
\end{example}

\section{Flatness, Gr\"{o}bner bases and the normal
sheaf\label{Groebnerbasesandflatness}}

\subsection{Flatness}

\begin{definition}
Let $A$ be a ring. An $A$-module $M$ is called
\index{flat module|textbf}%
\textbf{flat} over $A$ if for every injective homomorphism $N\rightarrow L$
the induced map $N\otimes_{A}M\rightarrow L\otimes_{A}M$ is injective.
\end{definition}

\begin{proposition}
\cite[Ch. 9]{Hartshorne Algebraic Geometry} Let $A$ be a ring and $M$ an
$A$-module. $M$ is flat over $A$ if and only if for all finitely generated
ideals $a\subset A$ the map%
\[
a\otimes M\rightarrow M
\]
is injective.
\end{proposition}

\begin{definition}
Given a morphism of schemes $f:Z\rightarrow Y$, an $\mathcal{O}_{Z}$-module
$\mathcal{F}$ is called flat over $Y$ at $z\in Z$ if $\mathcal{F}_{z}$ is flat
over $\mathcal{O}_{f\left(  z\right)  ,Y}$, which is considered as an
$\mathcal{O}_{f\left(  z\right)  ,Y}$-module via the natural map
$\mathcal{O}_{f\left(  z\right)  ,Y}\rightarrow\mathcal{O}_{z,Z}$.
$\mathcal{F}$ is called flat over $Y$ if it is flat over $Y$ for all $z\in Z$.

$Z$ is called a
\index{flat family|textbf}%
\textbf{flat family} over $Y$ if $\mathcal{O}_{Z}$ is flat over $Y$.
\end{definition}

\begin{proposition}
\cite[Ch. 9]{Hartshorne Algebraic Geometry} Let $A_{1}\rightarrow A_{2}$ be a
ring homomorphism and%
\[
f:\operatorname*{Spec}A_{2}\rightarrow\operatorname*{Spec}A_{1}%
\]
the corresponding morphism of affine schemes. If $M$ is an $A_{2}$-module,
then $\tilde{M}$ is flat over $\operatorname*{Spec}A_{1}$ if and only if $M$
is flat over $A_{1}$.
\end{proposition}

\subsection{First order deformations and the normal sheaf}

\begin{definition}
If $X_{0}\subset Y$ is a closed subscheme of a scheme $Y$ over $k$, a
\index{first order deformation|textbf}%
\textbf{first order
\index{deformation}%
deformation} of $X_{0}$ in $Y$ is a flat family $\mathfrak{X}\subset
Y\times_{k}\operatorname*{Spec}k\left[  t\right]  /\left\langle t^{2}%
\right\rangle $ over $\operatorname*{Spec}k\left[  t\right]  /\left\langle
t^{2}\right\rangle $ such that the fiber over $\operatorname*{Spec}%
k\subset\operatorname*{Spec}k\left[  t\right]  /\left\langle t^{2}%
\right\rangle $ is $X_{0}$.
\end{definition}

The tangent space of the
\index{Hilbert scheme}%
Hilbert scheme $\mathbb{H}_{Y}^{P}$ of subschemes with Hilbert polynomial $P$
of the projective scheme $Y$ at the point $X_{0}$ is the space of first order
\index{deformation}%
deformations of $X_{0}$ in $Y$.

We show that if $X_{0}\subset Y$ is a closed subscheme of a scheme $Y$ over
$k$, the space of first order
\index{deformation}%
deformations of $X_{0}$ in $Y$ coincides with the space of global sections of
$N_{X_{0}/Y}$:

Suppose $\mathfrak{X}\subset Y\times\operatorname*{Spec}k\left[  t\right]
/\left\langle t^{2}\right\rangle $ is a subscheme such that $X_{0}$ is
isomorphic to the fiber product%
\[%
\begin{tabular}
[c]{ccc}%
$\mathfrak{X}\times_{k\left[  t\right]  /\left\langle t^{2}\right\rangle
}\operatorname*{Spec}k$ & $\overset{\pi_{1}}{\rightarrow}$ & $\mathfrak{X}$\\
${\small \pi}_{{\small 2}}\downarrow$ &  & $\downarrow$\\
$\operatorname*{Spec}k$ & $\hookrightarrow$ & $\operatorname*{Spec}k\left[
t\right]  /\left\langle t^{2}\right\rangle $%
\end{tabular}
\]
and fix an isomorphism of $X_{0}$ and $\mathfrak{X}\times_{k\left[  t\right]
/\left\langle t^{2}\right\rangle }\operatorname*{Spec}k$. Consider an affine
open set $U\subset Y$:

Let $R=\mathcal{O}_{Y}\left(  U\right)  $ be the coordinate ring of $U$ and
$I\left(  X_{0}\cap U\right)  \subset R$ the ideal of $X_{0}\cap U$. Then
$N_{X_{0}/Y}\mid_{X_{0}\cap U}$ is the sheaf associated to
\[
\operatorname*{Hom}\nolimits_{R}\left(  I\left(  X_{0}\cap U\right)
,R/I\left(  X_{0}\cap U\right)  \right)
\]
The coordinate ring of $U\times\operatorname*{Spec}k\left[  t\right]
/\left\langle t^{2}\right\rangle $ is $R\otimes k\left[  t\right]
/\left\langle t^{2}\right\rangle $, so write the ideal of the intersection
$\mathfrak{U}$ of $\left(  U\times\operatorname*{Spec}k\left[  t\right]
/\left\langle t^{2}\right\rangle \right)  $ and $\mathfrak{X}$ as%
\[
I\left(  \mathfrak{U}\right)  =\left\langle f_{1}+t\cdot g_{1},...,f_{r}%
+t\cdot g_{r}\right\rangle
\]
with $I\left(  X_{0}\cap U\right)  =\left\langle f_{1},...,f_{r}\right\rangle
$ and $g_{i}\in R$.

We give different characterizations of flatness of $\mathfrak{U}$ over
$\operatorname*{Spec}k\left[  t\right]  /\left\langle t^{2}\right\rangle $:

$\mathfrak{U}$ is flat over $\operatorname*{Spec}k\left[  t\right]
/\left\langle t^{2}\right\rangle $ if and only if%
\[
\left\langle t\right\rangle \otimes\mathcal{O}_{\mathfrak{X}}\left(
\mathfrak{U}\right)  \rightarrow\mathcal{O}_{\mathfrak{X}}\left(
\mathfrak{U}\right)
\]
is injective, i.e., if and only if for all $f\in R$ it holds%
\[
tf\in I\left(  \mathfrak{U}\right)  \Rightarrow f\in I\left(  X_{0}\cap
U\right)
\]

\begin{itemize}
\item This is equivalent to the existence of a $\varphi_{U}\in
\operatorname*{Hom}\nolimits_{R}\left(  I\left(  X_{0}\cap U\right)
,R/I\left(  X_{0}\cap U\right)  \right)  $ with $\varphi_{U}\left(
f_{i}\right)  =g_{i}$:

First note that if $f\in R$ and $tf\in I\left(  \mathfrak{U}\right)  \subset
R\otimes k\left[  t\right]  /\left\langle t^{2}\right\rangle $, then there are
$a_{i},b_{i}\in R$ such that%
\[
tf=\sum_{i}\left(  a_{i}+tb_{i}\right)  \left(  f_{i}+tg_{i}\right)  =\sum
_{i}a_{i}f_{i}+t\sum_{i}\left(  a_{i}g_{i}+b_{i}f_{i}\right)
\]
hence $\sum_{i}a_{i}f_{i}=0$.

$\Longleftarrow:$ If there is a
\[
\varphi_{U}:I\left(  X_{0}\cap U\right)  \rightarrow R/I\left(  X_{0}\cap
U\right)
\]
with $\varphi_{U}\left(  f_{i}\right)  =g_{i}$, then
\begin{align*}
\sum_{i}a_{i}g_{i}+I\left(  X_{0}\cap U\right)   &  =\sum_{i}a_{i}\varphi
_{U}\left(  f_{i}\right)  =\varphi_{U}\left(  \sum_{i}a_{i}f_{i}\right) \\
&  =\varphi_{U}\left(  0\right)  =0\in R/I\left(  X_{0}\cap U\right)
\end{align*}
hence%
\[
f=\sum_{i}a_{i}g_{i}+\sum_{i}b_{i}f_{i}\in I\left(  X_{0}\cap U\right)
\]
i.e., for all $f\in R$ with $tf\in I\left(  \mathfrak{U}\right)  $ we have
$f\in I\left(  X_{0}\cap U\right)  $.

$\Longrightarrow:$ On the other hand if for all $f\in R$ with $tf\in I\left(
\mathfrak{U}\right)  $ it holds $f\in I\left(  X_{0}\cap U\right)  $, then the
homomorphism
\[
\varphi_{U}:I\left(  X_{0}\cap U\right)  \rightarrow R/I\left(  X_{0}\cap
U\right)
\]
is given in the following way: If $f=\sum_{i}a_{i}f_{i}\in I\left(  X_{0}\cap
U\right)  $, define $\varphi_{U}$ by
\[
\varphi_{U}\left(  f\right)  =\sum_{i}a_{i}g_{i}+I\left(  X_{0}\cap U\right)
\in R/I\left(  X_{0}\cap U\right)
\]
This is well defined: If%
\[
\sum_{i}a_{i}f_{i}=0
\]
then%
\[
t\sum_{i}a_{i}g_{i}=\sum_{i}a_{i}\left(  f_{i}+tg_{i}\right)  \in I\left(
\mathfrak{U}\right)
\]
hence $\sum_{i}a_{i}g_{i}\in I\left(  X_{0}\cap U\right)  $.

\item Existence of
\[
\varphi_{U}\in\operatorname*{Hom}\nolimits_{R}\left(  I\left(  X_{0}\cap
U\right)  ,R/I\left(  X_{0}\cap U\right)  \right)
\]
with $\varphi_{U}\left(  f_{i}\right)  =g_{i}$ is equivalent to the condition
that any
\index{syzygy}%
syzygy between $f_{1},...,f_{r}$ can be lifted to a
\index{syzygy}%
syzygy between
\[
f_{1}+tg_{1},...,f_{r}+tg_{r}%
\]
$\Longrightarrow:$ Suppose $\sum_{i}a_{i}f_{i}=0\in R$ and there is
$\varphi_{U}$ with $\varphi_{U}\left(  f_{i}\right)  =g_{i}$, then as above
\[
\sum_{i}a_{i}g_{i}+I\left(  X_{0}\cap U\right)  =\varphi_{U}\left(  \sum
_{i}a_{i}f_{i}\right)  =0\in R/I\left(  X_{0}\cap U\right)
\]
i.e., $\sum_{i}a_{i}g_{i}\in I\left(  X_{0}\cap U\right)  $, hence there are
$b_{i}\in R$ such that $-\sum_{i}a_{i}g_{i}=\sum_{i}b_{i}f_{i}$. So%
\[
\sum_{i}a_{i}\left(  f_{i}+tg_{i}\right)  =-t\sum_{i}b_{i}f_{i}%
\]
hence%
\[
\sum_{i}\left(  a_{i}+tb_{i}\right)  \left(  f_{i}+tg_{i}\right)  =0\in
R\otimes k\left[  t\right]  /\left\langle t^{2}\right\rangle
\]
$\Longleftarrow:$ On the other hand if $f=\sum_{i}a_{i}f_{i}\in I\left(
X_{0}\cap U\right)  $ and any
\index{syzygy}%
syzygy lifts, define as above $\varphi_{U}$ by
\[
\varphi_{U}\left(  f\right)  =\sum_{i}a_{i}g_{i}+I\left(  X_{0}\cap U\right)
\in R/I\left(  X_{0}\cap U\right)
\]

\end{itemize}

This ist well defined: If $\sum_{i}a_{i}f_{i}=0$, then there are $b_{i}\in R$
such that%
\[
t\left(  \sum_{i}a_{i}g_{i}+\sum_{i}b_{i}f_{i}\right)  =\sum_{i}\left(
a_{i}+tb_{i}\right)  \left(  f_{i}+tg_{i}\right)  =0\in R\otimes k\left[
t\right]  /\left\langle t^{2}\right\rangle
\]
hence $\sum_{i}a_{i}g_{i}\in I\left(  X_{0}\cap U\right)  $.

Summarizing these statements:

\begin{proposition}
\cite{EiHa Schemes The Language of Modern Algebraic Geometry}%
\label{Prop flatness over t mod t2} Let $X_{0}\subset Y$ be a closed subscheme
of a scheme $Y$ over $k$ and $\mathfrak{X}\subset Y\times\operatorname*{Spec}%
k\left[  t\right]  /\left\langle t^{2}\right\rangle $ a subscheme such that
$X_{0}\cong\mathfrak{X}\times_{k\left[  t\right]  /\left\langle t^{2}%
\right\rangle }\operatorname*{Spec}k$. Consider an affine open $U\subset Y$,
set $R=\mathcal{O}_{Y}\left(  U\right)  $ and write the ideal of the
intersection $\mathfrak{U}$ of $\mathfrak{X}$ and $\left(  U\times
\operatorname*{Spec}k\left[  t\right]  /\left\langle t^{2}\right\rangle
\right)  $ as%
\[
I\left(  \mathfrak{U}\right)  =\left\langle f_{1}+t\cdot g_{1},...,f_{r}%
+t\cdot g_{r}\right\rangle
\]
with $I\left(  X_{0}\cap U\right)  =\left\langle f_{1},...,f_{r}\right\rangle
$ and $g_{i}\in R$. Then the following statements are
\index{flat family}%
equivalent:

\begin{enumerate}
\item $\mathfrak{U}\rightarrow\operatorname*{Spec}k\left[  t\right]
/\left\langle t^{2}\right\rangle $ is flat

\item $\left\langle t\right\rangle \otimes\mathcal{O}_{\mathfrak{X}}\left(
\mathfrak{U}\right)  \rightarrow\mathcal{O}_{\mathfrak{X}}\left(
\mathfrak{U}\right)  $ is injective

\item For all $f\in R$ it holds%
\[
tf\in I\left(  \mathfrak{U}\right)  \Rightarrow f\in I\left(  X_{0}\cap
U\right)
\]

\item There is a unique
\[
\varphi_{U}\in\operatorname*{Hom}\nolimits_{R}\left(  I\left(  X_{0}\cap
U\right)  ,R/I\left(  X_{0}\cap U\right)  \right)
\]
with
\[
\varphi_{U}\left(  f_{i}\right)  =g_{i}%
\]

\item Any
\index{syzygy}%
syzygy between $f_{1},...,f_{r}$ lifts to a syzygy between $f_{1}%
+tg_{1},...,f_{r}+tg_{r}$, i.e., if
\[
\sum_{i}a_{i}f_{i}=0\in R
\]
with $a_{i}\in R$, there are $b_{i}\in R$ such that%
\[
\sum_{i}\left(  a_{i}+tb_{i}\right)  \left(  f_{i}+tg_{i}\right)  =0\in
R\otimes k\left[  t\right]  /\left\langle t^{2}\right\rangle
\]

\end{enumerate}
\end{proposition}

So if $\mathfrak{X}\rightarrow\operatorname*{Spec}k\left[  t\right]
/\left\langle t^{2}\right\rangle $ is
\index{flat family}%
flat, then for any affine open set $U$ there is a unique $\varphi_{U}$, and
the $\varphi_{U}$ patch together to a section of $N_{X_{0}/Y}$. On the other
hand, if $\varphi$ is a global section of $N_{X_{0}/Y}$, then define the
associated
\index{first order deformation}%
first order
\index{deformation}%
deformation $\mathfrak{X}$ by the local equations%
\[
\left\{  f+t\cdot\varphi\left(  f\right)  =0\mid f\in I\left(  X_{0}\cap
U\right)  \right\}
\]
on $U\times\operatorname*{Spec}k\left[  t\right]  /\left\langle t^{2}%
\right\rangle $ hence:

\begin{theorem}
If $X_{0}\subset Y$ is a closed subscheme of a scheme $Y$ over $k$, the space
of
\index{first order deformation}%
first order
\index{deformation}%
deformations of $X_{0}$ in $Y$ equals the space of
\index{global sections}%
global sections of $N_{X_{0}/Y}$.
\end{theorem}

\subsection{Flatness over $k\left[  \left[  t\right]  \right]  $}

The statement analogous to Proposition \ref{Prop flatness over t mod t2} for
base $\operatorname*{Spec}k\left[  \left[  t\right]  \right]  $ is given in
the following.

\begin{proposition}
Let $X_{0}\subset Y$ be a closed subscheme of a scheme $Y$ over $k$ and
$\mathfrak{X}\subset Y\times\operatorname*{Spec}k\left[  \left[  t\right]
\right]  $ a subscheme such that $X_{0}\cong\mathfrak{X}\times_{k\left[
\left[  t\right]  \right]  }\operatorname*{Spec}k$. Consider an affine open
$U\subset Y$, let $R=\mathcal{O}_{Y}\left(  U\right)  $ and write the ideal of
the intersection $\mathfrak{U}$ of $\mathfrak{X}$ and $\left(  U\times
\operatorname*{Spec}k\left[  \left[  t\right]  \right]  \right)  $ as%
\[
I\left(  \mathfrak{U}\right)  =\left\langle f_{1}+t\cdot g_{1},...,f_{r}%
+t\cdot g_{r}\right\rangle
\]
with $I\left(  X_{0}\cap U\right)  =\left\langle f_{1},...,f_{r}\right\rangle
$ and $g_{i}\in R\otimes k\left[  \left[  t\right]  \right]  $. Then the
following statements are
\index{flat family}%
equivalent:

\begin{enumerate}
\item $\mathfrak{U}\rightarrow\operatorname*{Spec}k\left[  \left[  t\right]
\right]  $ is flat

\item $\left\langle t\right\rangle \otimes\mathcal{O}_{\mathfrak{X}}\left(
\mathfrak{U}\right)  \rightarrow\mathcal{O}_{\mathfrak{X}}\left(
\mathfrak{U}\right)  $ is injective

\item For all $f\in R$ it holds%
\[
tf\in I\left(  \mathfrak{U}\right)  \Rightarrow f\in I\left(  X_{0}\cap
U\right)
\]

\item Any
\index{syzygy}%
syzygy between $f_{1},...,f_{r}$ lifts to a syzygy between $f_{1}%
+tg_{1},...,f_{r}+tg_{r}$, i.e., if
\[
\sum_{i}a_{i}f_{i}=0\in R
\]
with $a_{i}\in R$, there are $b_{i}\in R\otimes k\left[  \left[  t\right]
\right]  $ such that%
\[
\sum_{i}\left(  a_{i}+tb_{i}\right)  \left(  f_{i}+tg_{i}\right)  =0\in
R\otimes k\left[  \left[  t\right]  \right]
\]

\end{enumerate}
\end{proposition}

\section{Gr\"{o}bner fan, state polytope, Hilbert scheme and
stability\label{Sec Computing the Bergman fan}}

\subsection{Concept for computing the Bergman
fan\label{sec computing the Bergman fan subsection}}

Let $I$ be
\index{Bergman fan}%
an ideal in $\mathbb{C}\left[  x_{0},...,x_{n}\right]  $ and $w$ a
\index{global ordering}%
global
\index{weight vector}%
weight vector on the variables of $\mathbb{C}\left[  x_{0},...,x_{n}\right]
$, i.e., $w_{i}\geq0$ for all $i$. Given a
\index{monomial ordering}%
monomial ordering $>$ we have $L_{>}\left(  in_{w}\left(  g\right)  \right)
=L_{>_{w}}\left(  g\right)  $ for every $g\in\mathbb{C}\left[  x_{0}%
,...,x_{n}\right]  $, so the subsets
\index{initial form}%
of
\index{initial ideal}%
monomials in $L_{>}\left(  in_{w}\left(  I\right)  \right)  $ and $L_{>_{w}%
}\left(  I\right)  $ are equal, hence:

\begin{proposition}
If $>$ is any
\index{monomial ordering}%
monomial ordering,
\index{initial ideal}%
then%
\[
L_{>}\left(  in_{w}\left(  I\right)  \right)  =L_{>_{w}}\left(  I\right)
\]

\end{proposition}

\begin{corollary}
\label{1reducedgbofinitialideal}If $g=\left(  g_{1},...,g_{r}\right)  $ is
the
\index{reduced standard basis}%
reduced
\index{Gr\"{o}bner basis}%
Gr\"{o}bner basis of $I$ with
\index{initial form}%
respect to $>_{w}$, then%
\[
\left(  in_{w}\left(  g_{i}\right)  \mid i=1,...,r\right)
\]
is the
\index{reduced standard basis}%
reduced
\index{Gr\"{o}bner basis}%
Gr\"{o}bner basis of $in_{w}\left(  I\right)  $ with respect to $>$.
\end{corollary}

\begin{proposition}
If $g=\left(  g_{1},...,g_{r}\right)  $ is the
\index{reduced standard basis}%
reduced
\index{Gr\"{o}bner basis}%
Gr\"{o}bner basis of $I$ with respect to $>_{w}$, then $\operatorname*{in}%
_{w}\left(  I\right)  $ contains a
\index{initial form}%
monomial if and only if%
\[
\left(  \left\langle in_{w}\left(  g_{i}\right)  \mid i=1,...,r\right\rangle
:\left\langle x_{0}\cdot...\cdot x_{n}\right\rangle ^{\infty}\right)
=\left\langle 1\right\rangle
\]

\end{proposition}

\begin{remark}
To speed up computations, one first checks if any of
\index{initial form}%
the $in_{w}\left(  g_{i}\right)  $, $i=1,...,r$ is a monomial, before doing
the saturation.
\end{remark}

Let $I$ be an ideal in $\mathbb{C}\left[  x_{1},...,x_{n}\right]  $ and
$J\subset\mathbb{C}\left[  x_{0},x_{1},...,x_{n}\right]  $ be the ideal of
the
\index{projective closure}%
projective closure of $V\left(  I\right)  $.

\begin{proposition}
\cite[Sec. 9.2]{Sturmfels Solving Systems of Polynomial Equations}%
\[
BF\left(  I\right)  =\left\{  w\in\mathbb{R}^{n}\mid in_{\left(  0,-w\right)
}\left(  J\right)  \text{ does not contain a monomial}\right\}
\]

\end{proposition}

This allows to
\index{initial ideal}%
homogenize, so any
\index{monomial ordering}%
monomial ordering will be equivalent to a
\index{global ordering}%
global ordering, hence the
\index{reduced standard basis}%
reductions in Gr\"{o}bner computations stay finite.

To compute $BF\left(  I\right)  $ we have to understand, which
\index{initial ideal}%
initial ideals can occur. This is described by the Gr\"{o}bner fan.

\subsection{Tropical representation of Gr\"{o}bner
cones\label{Sec tropical representation of Groebner cones}}

Although there are infinitely many
\index{global ordering}%
global
\index{semigroup ordering}%
semigroup orderings on the monomials of $\mathbb{C}\left[  x_{1}%
,...,x_{n}\right]  $, if we fix an ideal $J\subset\mathbb{C}\left[
x_{1},...,x_{n}\right]  $ and consider $>_{1}$ and $>_{2}$ equivalent if
$L_{>_{1}}\left(  J\right)  =L_{>_{2}}\left(  J\right)  $, there are only
finitely many equivalence classes of
\index{monomial ordering}%
monomial orderings.

\begin{definition}
Given a
\index{global ordering}%
global ordering $>$ on the
\newsym[$C_{>}\left(  J\right)  $]{Gr\"{o}bner cone with respect to monomial ordering}{}monomials
of $\mathbb{C}\left[  x_{1},...,x_{n}\right]  $ and an ideal $J\subset
\mathbb{C}\left[  x_{0},...,x_{n}\right]  $
\index{Gr\"{o}bner cone|textbf}%
define%
\[
C_{>}\left(  J\right)  =\left\{  w^{\prime}\in\mathbb{R}^{n}\mid
in_{w^{\prime}}\left(  J\right)  =L_{>}\left(  J\right)  \right\}
\]
If $J_{0}=L_{>}\left(  J\right)  $ for
\index{lead ideal}%
some
\index{initial ideal}%
ordering $>$, then
\newsym[$C_{J_{0}}\left(  J\right)$]{closed Gr\"{o}bner cone with respect to lead ideal}{}define
$C_{J_{0}}\left(  J\right)  =\overline{C_{>}\left(  J\right)  }$.
\end{definition}

By Proposition \ref{1anymonordbyweightord} and Proposition
\ref{1finitedetermstd} we have:

\begin{lemma}
Given a
\index{global ordering}%
global ordering $>$ on the monomials of $\mathbb{C}\left[  x_{1}%
,...,x_{n}\right]  $ and an ideal $J\subset\mathbb{C}\left[  x_{0}%
,...,x_{n}\right]  $, there is some $w\in\mathbb{Z}^{n}$ with positive
entries, such
\index{lead ideal}%
that $in_{w}\left(  J\right)  =L_{>}\left(  J\right)  $,
\index{initial ideal}%
hence, $C_{>}\left(  J\right)  $ is non empty.
\end{lemma}

\begin{definition}
Given $w\in\mathbb{Z}^{n}$ with positive
\newsym[$C_{w}\left(  J\right)  $]{Gr\"{o}bner cone with respect to weight vector $w$}{}entries
and $J\subset\mathbb{C}\left[  x_{0},...,x_{n}\right]  $
\index{Gr\"{o}bner cone|textbf}%
define%
\[
C_{w}\left(  J\right)  =\left\{  w^{\prime}\in\mathbb{R}^{n}\mid
in_{w^{\prime}}\left(  J\right)  =in_{w}\left(  J\right)  \right\}
\]

\end{definition}

Consider some
\index{tie break ordering}%
tie break ordering $>$ and the unique
\index{reduced standard basis}%
reduced
\index{Gr\"{o}bner basis}%
Gr\"{o}bner basis $\mathcal{G}=\left(  g_{1},...,g_{r}\right)  $ of $J$ with
respect to $>_{w}$,
\index{initial form}%
then%
\begin{equation}
C_{w}\left(  J\right)  =\left\{  w^{\prime}\in\mathbb{R}^{n}\mid
in_{w^{\prime}}\left(  g_{i}\right)  =in_{w}\left(  g_{i}\right)  \forall
i=1,...,r\right\}  \label{4equCwJreducedgroebnerbasis}%
\end{equation}

To see this, suppose $w^{\prime}\in\mathbb{R}^{n}$ with $in_{w^{\prime}%
}\left(  g_{i}\right)  =in_{w}\left(  g_{i}\right)  $,
\index{initial form}%
then%
\[
in_{w}\left(  J\right)  =\left\langle in_{w^{\prime}}\left(  g_{i}\right)
\mid i=1,...,r\right\rangle \subset in_{w^{\prime}}\left(  J\right)
\]
so%
\[
L_{>_{w}}\left(  J\right)  =L_{>}\left(  in_{w}\left(  J\right)  \right)
\subset L_{>}\left(  in_{w^{\prime}}\left(  J\right)  \right)
=L_{>_{w^{\prime}}}\left(  J\right)
\]
therefore
\index{initial ideal}%
the
\index{lead ideal}%
lead ideals $L_{>_{w}}\left(  J\right)  $ and $L_{>_{w^{\prime}}}\left(
J\right)  $ are equal, hence by Proposition \ref{IdealLeadideal}%
\[
in_{w}\left(  J\right)  =in_{w^{\prime}}\left(  J\right)
\]
i.e., $w^{\prime}\in C_{w}\left(  J\right)  $.

On the other hand if $w^{\prime}\in\mathbb{R}^{n}$
\index{initial ideal}%
with%
\[
in_{w}\left(  J\right)  =in_{w^{\prime}}\left(  J\right)
\]
then by Corollary \ref{1reducedgbofinitialideal} the
\index{reduced standard basis}%
reduced
\index{Gr\"{o}bner basis}%
Gr\"{o}bner basis
\index{initial ideal}%
of $in_{w^{\prime}}\left(  J\right)  $ with respect to $>$ is given
\index{initial form}%
by $G=\left(  in_{w}\left(  g_{1}\right)  ,...,in_{w}\left(  g_{r}\right)
\right)  $. So for all $i=1,...,r$ we have $NF_{>}\left(  in_{w^{\prime}%
}\left(  g_{i}\right)  ,G\right)  =0$, hence $m_{i}=L_{>_{w}}\left(
g_{i}\right)  $ is a monomial of $in_{w^{\prime}}\left(  g_{i}\right)  $, as
by
\index{reduced standard basis}%
reducedness, $m_{i}$ is the only monomial of $g_{i}$ in $L_{>_{w}}\left(
J\right)  $. Write%
\begin{align*}
in_{w}\left(  g_{i}\right)   &  =m_{i}+h_{i}\\
in_{w^{\prime}}\left(  g_{i}\right)   &  =m_{i}+h_{i}^{\prime}%
\end{align*}
then $h_{i}$ and $h_{i}^{\prime}$ do not involve monomial of $L_{>_{w}}\left(
J\right)  $. The first step of the division with remainder,
\index{initial ideal}%
calculating $NF_{>}\left(  in_{w^{\prime}}\left(  g_{i}\right)  ,G\right)  $,
gives $h_{i}^{\prime}-h_{i}\in in_{w^{\prime}}\left(  J\right)  =in_{w}\left(
J\right)  $. On the other hand, no term of $h_{i}^{\prime}-h_{i}$ is
\index{lead ideal}%
in $L_{>_{w}}\left(  J\right)  =L_{>}\left(  in_{w}\left(  J\right)  \right)
$, hence, $h_{i}^{\prime}=h_{i}$, i.e., $in_{w^{\prime}}\left(  g_{i}\right)
=in_{w}\left(  g_{i}\right)  $.

Suppose now $w\in\mathbb{R}^{n}$ is
\index{lead term}%
representing the
\index{monomial ordering}%
monomial ordering $>$, then define
\[
m_{i}=c_{i}x^{a_{i}}=LT_{>}\left(  g_{i}\right)  =in_{w}\left(  g_{i}\right)
\]
and write $g_{i}=m_{i}+h_{i}$ with the
\index{tail}%
tail $h_{i}$ of $g_{i}$. By the description of $C_{w}\left(  J\right)  $ via
Equation \ref{4equCwJreducedgroebnerbasis}%
\begin{align*}
C_{w}\left(  J\right)   &  =\left\{  w^{\prime}\in\mathbb{R}^{n}\mid
w^{\prime}b_{i}<w^{\prime}a_{i}\text{ }\forall\text{ monomials }x^{b_{i}%
}\text{ of the tail }h_{i}\text{ and }\forall\text{ }i=1,...,r\right\} \\
&  =\left\{  w^{\prime}\in\mathbb{R}^{n}\mid\operatorname*{trop}\left(
g-in_{w}\left(  g\right)  \right)  \left(  w^{\prime}\right)
<\operatorname*{trop}\left(  in_{w}\left(  g\right)  \right)  \forall
g\in\mathcal{G}\right\}
\end{align*}
so summarizing:

\begin{theorem}
If $w\in\mathbb{Z}^{n}$ with positive entries, $J\subset\mathbb{C}\left[
x_{1},...,x_{n}\right]  $ and $>$ some
\index{global ordering}%
global ordering on the monomials of $\mathbb{C}\left[  x_{1},...,x_{n}\right]
$ and $\mathcal{G}=\left(  g_{1},...,g_{r}\right)  $ is the unique
\index{reduced standard basis}%
reduced
\index{Gr\"{o}bner basis}%
Gr\"{o}bner basis of $J$ with respect to $>_{w}$, then%
\[
C_{w}\left(  J\right)  =\left\{  w^{\prime}\in\mathbb{R}^{n}\mid
in_{w^{\prime}}\left(  g_{i}\right)  =in_{w}\left(  g_{i}\right)  \forall
i=1,...,r\right\}
\]
in particular $C_{w}\left(  J\right)  $ is a relatively open convex
\index{Gr\"{o}bner cone}%
polyhedral
\index{initial form}%
cone.

If $in_{w}\left(  J\right)  =L_{>}\left(  J\right)  $, then%
\begin{equation}
C_{w}\left(  J\right)  =\left\{  w^{\prime}\in\mathbb{R}^{n}\mid
\operatorname*{trop}\left(  g-in_{w}\left(  g\right)  \right)  \left(
w^{\prime}\right)  <\operatorname*{trop}\left(  in_{w}\left(  g\right)
\right)  \forall g\in\mathcal{G}\right\}  \label{tropical equations for cone}%
\end{equation}

\end{theorem}

\begin{remark}
If $in_{w}\left(  J\right)  $ is not
\index{Gr\"{o}bner cone}%
monomial, then $C_{w}\left(  J\right)  $ is given by these inequalities
together with the equalities coming from the condition that for each
$g\in\mathcal{G}$ all monomials of the
\index{initial form}%
initial form $in_{w}\left(  g\right)  $ have the same weight.
\end{remark}

\begin{definition}
Let $J\subset\mathbb{C}\left[  x_{1},...,x_{n}\right]  $ be an ideal. The
Gr\"{o}bner region of $J$ is%
\[
GR\left(  J\right)  =\left\{  w\in\mathbb{R}^{n}\mid\exists w^{\prime}%
\in\mathbb{R}_{\geq0}^{n}\text{ with }in_{w}\left(  J\right)  =in_{w^{\prime}%
}\left(  J\right)  \right\}
\]

\end{definition}

\begin{lemma}
\cite[Ch. 1]{Sturmfels Groebner Bases and Convex Polytopes} If $J\subset
\mathbb{C}\left[  x_{1},...,x_{n}\right]  $ is homogeneous, then $GR\left(
J\right)  =\mathbb{R}^{n}$.
\end{lemma}

\begin{definition}
Let $J\subset\mathbb{C}\left[  x_{1},...,x_{n}\right]  $ and
\newsym[$GF\left(  J\right)  $]{Gr\"{o}bner fan}{}suppose $GR\left(  J\right)
=\mathbb{R}^{n}$. The
\index{Gr\"{o}bner fan|textbf}%
\textbf{Gr\"{o}bner fan} $GF\left(  J\right)  $ is the set of all closures
$\overline{C_{w}\left(  J\right)  }$ of cones $C_{w}\left(  J\right)  $ for
all $w\in\mathbb{R}^{n}$.
\end{definition}

If $f=\sum_{\alpha}c_{\alpha}x^{\alpha}$ is a Laurent polynomial in the
\newsym[$N\left(  f\right)  $]{Newton polytope}{}variables $x_{1},...,x_{n}$,
then its
\index{Newton polytope|textbf}%
\textbf{Newton polytope} is%
\[
N\left(  f\right)  =\operatorname*{convexhull}\left\{  \alpha\mid c_{\alpha
}\neq0\right\}  \subset\mathbb{R}^{n+1}%
\]

\begin{lemma}
\cite[Section 6.1]{GKZ Discriminants Resultants and Multidimensional
Determinants} The Newton polytope of $f$ lies in the hyperplane $\left\{
\alpha\in\mathbb{R}^{n+1}\mid w\cdot\alpha=a\right\}  $ for some
$a\in\mathbb{Z}$ and $w\in\mathbb{Z}^{n+1}$ if and only if $f$ is
$w$-homogeneous, i.e $f\left(  t^{w_{0}}x_{1},...,t^{w_{n}}x_{n}\right)
=t^{a}f\left(  x_{1},...,x_{n}\right)  $.

If $f,g$ are Laurent polynomials, then $N\left(  f\cdot g\right)  =N\left(
f\right)  +N\left(  g\right)  $.
\end{lemma}

\begin{lemma}
If $w\in\mathbb{Z}^{n}$ with positive entries, $I\subset\mathbb{C}\left[
x_{1},...,x_{n}\right]  $, $>$ some
\index{global ordering}%
global ordering on the monomials of $\mathbb{C}\left[  x_{0},...,x_{n}\right]
$ and $\mathcal{G}$ the unique
\index{reduced standard basis}%
reduced
\index{Gr\"{o}bner basis}%
Gr\"{o}bner basis of $J$ with respect to $>_{w}$, then%
\[
C_{w}\left(  J\right)  =\sigma_{Q}\left(  \operatorname*{face}\nolimits_{w}%
\left(  Q\right)  \right)
\]
is the normal cone of the face $\operatorname*{face}\nolimits_{w}\left(
Q\right)  $ of $Q$ with%
\[
Q=N\left(  \prod\nolimits_{g\in\mathcal{G}}g\right)  =\sum\nolimits_{g\in
\mathcal{G}}N\left(  g\right)
\]

\end{lemma}

Using this representation of $C_{w}\left(  J\right)  $ one can conclude:

\begin{proposition}
\cite[Ch. 2]{Sturmfels Groebner Bases and Convex Polytopes} $GF\left(
J\right)  $ is a fan.
\end{proposition}

\subsection{Computing the Gr\"{o}bner
fan\label{Sec computing the Groebner fan}}

\begin{algorithm}
Given an ideal $J\subset k\left[  x_{1},...,x_{n}\right]  $ with $GR\left(
J\right)  =\mathbb{R}^{n}$ and a subfan $F\subset\mathbb{R}^{n}$ of the
Gr\"{o}bner fan of $I$ the following
\index{findRandomCone}%
algorithm $\operatorname*{findRandomCone}$ computes some cone of the
Gr\"{o}bner fan which is not in $F$:

Choose some random $w\in\mathbb{R}^{n}-\operatorname*{supp}\left(  F\right)
;$

Let $>_{w}$ be the corresponding
\index{weight ordering}%
weight ordering$;$

$g:=\operatorname*{redStd}_{Wp\left(  w\right)  }\left(  J\right)  $, i.e.,
the
\index{reduced standard basis}%
reduced
\index{Gr\"{o}bner basis}%
Gr\"{o}bner basis of $J$ with respect to $Wp\left(  w\right)  ;$

$\operatorname*{if}$ $\operatorname*{in}_{w}\left(  J\right)  =\left\langle
in_{w}\left(  g_{i}\right)  \mid i=1,...,r\right\rangle $ contains a monomial
repeat with different $w;$

Compute $\overline{C_{w}\left(  J\right)  }$ from $g$ via
\ref{tropical equations for cone}$;$

$\operatorname*{return}\left(  \overline{C_{w}\left(  J\right)  }\right)  ;$
\end{algorithm}

The following randomized algorithm computes the
\index{Gr\"{o}bner fan}%
Gr\"{o}bner fan:

\begin{algorithm}
Given an ideal $J\subset k\left[  x_{1},...,x_{n}\right]  $ with $GR\left(
J\right)  =\mathbb{R}^{n}$ the following algorithm computes the
\index{Gr\"{o}bner fan}%
Gr\"{o}bner fan of $J$:

Let $F$ be the empty fan in $\mathbb{R}^{n}$.

$\operatorname*{while}$ $\operatorname*{isComplete}\left(  F\right)  =false$
$\operatorname*{do}$

\quad$F:=$ the fan generated by the cones of $F$
\index{findRandomCone}%
and $\operatorname*{findRandomCone}\left(  J,F\right)  ;$

$\operatorname*{od};$
\end{algorithm}

The following algorithm avoids searching a
\index{weight vector}%
weight vector in the complement of the support of a non
\index{complete fan}%
complete fan and it integrates the test for completeness:

\begin{algorithm}
Given an ideal $J\subset k\left[  x_{1},...,x_{n}\right]  $ with $GR\left(
J\right)  =\mathbb{R}^{n}$ the following algorithm computes the
\index{Gr\"{o}bner fan}%
Gr\"{o}bner fan:

$F:=$ the fan generated
\index{findRandomCone}%
by $\operatorname*{findRandomCone}\left(  J,F\right)  ;$

$remainingfacets:=facets\left(  cones\left(  F\right)  \left[  1\right]
\right)  ;$

$\operatorname*{while}$ $remainingfacets<>\left\{  {}\right\}  $
$\operatorname*{do}$

\quad$fc:=remainingfacets\left[  1\right]  ;$

\quad$outernormal:=-\operatorname*{rays}\left(  \operatorname*{dual}\left(
fc\right)  \right)  \left[  1\right]  ;$

\quad$internal:=\operatorname*{sum}\left(  \operatorname*{rays}\left(
fc\right)  \right)  ;$

\quad$s:=1;$

\quad$w:=s\cdot internal+outernormal;$

\quad$\operatorname*{while}$ $w\in\operatorname*{support}\left(  F\right)  $
$\operatorname*{do}$

\quad\quad$s:=10\cdot s;$

\quad\quad$w:=s\cdot internal+outernormal;$

\quad$\operatorname*{od};$

\quad$F:=$ the fan given by the cones of $F$ and all faces of $C_{>_{w}%
}\left(  J\right)  ;$

\quad$\operatorname*{for}$ $\operatorname*{all}$ $fct\in facets\left(
C_{>_{w}}\left(  J\right)  \right)  $ $\operatorname*{do}$

\quad\quad$\operatorname*{if}$ $fct\in remainingfacets$ $\operatorname*{then}$

\quad\quad\quad$remainingfacets:=remainingfacets-\left\{  fct\right\}  ;$

\quad\quad$else$

\quad\quad\quad$remainingfacets:=remainingfacets\cup\left\{  fct\right\}  ;$

\quad\quad$\operatorname*{fi};$

\quad$\operatorname*{od};$

$\operatorname*{od};$

Note that this algorithm necessarily stops with a
\index{complete fan}%
complete fan, as the set $remainingfacets$ is empty if and only if all facets
of cones have appeared twice.
\end{algorithm}

\begin{remark}
To compute the
\index{Bergman fan}%
Bergman fan out of the
\index{Gr\"{o}bner fan}%
Gr\"{o}bner fan, it is computationally important to note that if a cone $F$ of
the
\index{Gr\"{o}bner fan}%
Gr\"{o}bner fan is not contained in the
\index{Bergman fan}%
Bergman fan, then also all higher dimensional cones containing $F$ are not in
the
\index{Bergman fan}%
Bergman fan.
\end{remark}

\subsection{Hilbert scheme and state polytope: Projective
setup\label{Sec Hilbert scheme and state polytope projective setup}}

\subsubsection{Gr\"{o}bner fan and state polytopes}

Let $I\subset S=\mathbb{C}\left[  x_{0},...,x_{n}\right]  $ be a homogeneous
ideal. For $d\geq1$ define%
\[
P_{d}\left(  I\right)  =\operatorname*{convexhull}\left\{  \prod_{x^{\alpha
}\in M_{d}}x^{\alpha}\mid M=in_{>}\left(  I\right)  \text{,}>\text{ a monomial
ordering}\right\}  \subset\mathbb{R}^{n+1}%
\]
If $d_{0}$ is the maximum degree appearing in a minimal universal Gr\"{o}bner
basis of $I$, define%
\[
P\left(  I\right)  =\sum_{d=1}^{d_{0}}P_{d}\left(  I\right)
\]

\begin{definition}
A
\index{state polytope|textbf}%
\textbf{state polytope for} $I$ is a
\newsym[$P\left(  I\right)  $]{state polytope}{}polytope $P\subset
M_{\mathbb{R}}$ with $GF\left(  I\right)  =\operatorname*{NF}\left(  P\right)
$.
\end{definition}

\begin{proposition}
\cite{Sturmfels Groebner Bases and Convex Polytopes} The Gr\"{o}bner fan
$GF\left(  I\right)  $ of $I$ is the normal fan of $P\left(  I\right)  $%
\[
GF\left(  I\right)  =\operatorname*{NF}\left(  P\left(  I\right)  \right)
\]
so $P\left(  I\right)  $ is a state polytope for $I$.

If $w=\left(  w_{0},...,w_{n}\right)  \in\mathbb{R}^{n+1}$, then%
\[
\operatorname*{face}\nolimits_{w}\left(  P_{d}\left(  I\right)  \right)
=P_{d}\left(  in_{w}\left(  I\right)  \right)
\]
and%
\[
\operatorname*{face}\nolimits_{w}\left(  P\left(  I\right)  \right)
=\sum_{d=1}^{d_{0}}P_{d}\left(  in_{w}\left(  I\right)  \right)
\]

If $>$ and $>^{\prime}$ are monomial orderings, then $\prod_{x^{\alpha}\in
in_{>}\left(  I\right)  _{d}}x^{\alpha}=\prod_{x^{\alpha}\in in_{>^{\prime}%
}\left(  I\right)  _{d}}x^{\alpha}$ if and only if $in_{>}\left(  I\right)
_{d}=in_{>^{\prime}}\left(  I\right)  _{d}$.
\end{proposition}

\begin{proposition}
\cite{Sturmfels Groebner Bases and Convex Polytopes} If $\mathcal{G}$ is a
universal Gr\"{o}bner basis of $I$ which is reduced with respect to any
monomial ordering, then%
\[
\sum_{g\in\mathcal{G}}N\left(  g\right)
\]
is a state polytope for $I$.
\end{proposition}

\subsubsection{State polytope and the Hilbert scheme}

Suppose $V=$ $_{\mathbb{C}}\left\langle x_{0},...,x_{n}\right\rangle
=\mathbb{C}^{n}$, $S=\operatorname*{Sym}\left(  V\right)  \cong\mathbb{C}%
\left[  x_{0},...,x_{n}\right]  $ and $I\subset S$ is a homogeneous ideal such
that $S/I$ has Hilbert polynomial $P=P_{S/I}$.

\begin{lemma}
\cite{Bayer The division algorithm and the Hilbert scheme}
\label{Lem hilbert point degree d}There is a degree $d_{0}$ such that for all
$d\geq d_{0}$ and all homogeneous saturated ideals $J\subset S$ with Hilbert
polynomial $P_{S/J}=P$

\begin{itemize}
\item $J$ is determined by the degree $d$ part $J_{d}$ of $J$, i.e.,
$J=\left(  \left\langle J_{d}\right\rangle :\left\langle x_{0},...,x_{n}%
\right\rangle ^{\infty}\right)  $

\item $\dim_{\mathbb{C}}\left(  S_{d}/J_{d}\right)  =P\left(  d\right)  $

\item For all semigroup orderings $>$
\[
in_{>}\left(  J\right)  =\left\langle in_{>}\left(  f\right)  \mid f\in
J\text{ with }\deg\left(  f\right)  \leq d\right\rangle
\]

\end{itemize}
\end{lemma}

\begin{definition}
With above notation $I_{d}$ is a point in the
\newsym[$\mathbb{G}\left(  r,V\right)  $]{Grassmannian}{}Grassmannian
$\mathbb{G}\left(  P\left(  d\right)  ,S_{d}\right)  $ of subspaces with
codimension $P\left(  d\right)  $. This point is denoted as the $d$-th
\index{Hilbert Point|textbf}%
\textbf{Hilbert point} $H\left(  I\right)  $
\newsym[$H\left(  I\right)  $]{Hilbert point}{}of $I$. The Hilbert point
$H\left(  I\right)  $ determines $\left(  I:\left\langle x_{0},...,x_{n}%
\right\rangle ^{\infty}\right)  $.

The set of all $d$-th Hilbert points $H\left(  J\right)  $ of homogeneous
ideals $J\subset S$ with $P_{S/J}=P$ is a closed
\newsym[$\mathbb{H}_{n}^{P}$]{classical Hilbert scheme}{}subscheme
$\mathbb{H}_{n}^{P}\subset\mathbb{G}\left(  P\left(  d\right)  ,S_{d}\right)
$, the $P$-th
\index{Hilbert scheme|textbf}%
\textbf{Hilbert scheme}.
\end{definition}

\begin{remark}
Let $>$ be a
\index{total ordering}%
total ordering of the monomials of degree $d$ of $S$ and $x^{\alpha_{1}%
},...,x^{\alpha_{s}}$ with $s=\binom{n+d-1}{d}$ a monomial basis of $S_{d}$
ordered with respect to $>$. If $B=\left(  f_{1},...,f_{r}\right)  $ with
$r=s-P\left(  d\right)  $ is a basis of $I_{d}$, then writing%
\[
f_{j}=\sum_{\left\vert \alpha\right\vert =d}a_{j,\alpha}x^{\alpha}%
\]
we obtain
\index{Hilbert matrix|textbf}%
the $>$-\textbf{Hilbert matrix} of $I$%
\[
A=\left(  a_{j,\alpha_{i}}\right)  _{j=1,...r,\text{ }i=1,...,s}\in k^{r\times
s}%
\]
representing $H\left(  I\right)  $ with respect to above bases%
\[%
\begin{tabular}
[c]{lll}%
$I_{d}$ & $\hookrightarrow$ & $S_{d}$\\
$\uparrow\cong$ &  & $\uparrow\cong$\\
$\mathbb{C}^{r}$ & $\overset{A^{t}}{\hookrightarrow}$ & $\mathbb{C}^{s}$%
\end{tabular}
\]
If $B^{\prime}=\left(  f_{1}^{\prime},...,f_{r}^{\prime}\right)  $ is another
basis of $I_{d}$ and $A^{\prime}$ the corresponding $>$-Hilbert matrix, then
there is a $Q\in\operatorname*{GL}\left(  r,k\right)  $ with $A^{\prime}=QA$.
\end{remark}

Let $s=\dim S_{d}$ and $r=s-P\left(  d\right)  $. The
\index{Pl\"{u}cker embedding|textbf}%
\textbf{Pl\"{u}cker embedding}%
\[%
\begin{tabular}
[c]{cccc}%
$\mathfrak{p}:$ & $\mathbb{G}\left(  P\left(  d\right)  ,S_{d}\right)  $ &
$\mathbb{\rightarrow}$ & $\mathbb{P}\left(  \bigwedge\nolimits^{r}%
S_{d}\right)  $\\
& $\left\langle f_{1},...,f_{r}\right\rangle $ & $\mapsto$ & $f_{1}%
\wedge...\wedge f_{r}$%
\end{tabular}
\]
of $\mathbb{G}\left(  P\left(  d\right)  ,S_{d}\right)  $ in the
\newsym[$\mathfrak{p}$]{Pl\"{u}cker embedding}{}projective space
$\mathbb{P}\left(  W\right)  $ with
\[
W=\bigwedge\nolimits^{r}S_{d}%
\]
is given by the positive line bundle $L=\det U^{\ast}$, where $U$ is the
universal subbundle $U\rightarrow\mathbb{G}\left(  P\left(  d\right)
,S_{d}\right)  $ of $\mathbb{C}^{s}\times\mathbb{G}\left(  P\left(  d\right)
,S_{d}\right)  \rightarrow\mathbb{G}\left(  P\left(  d\right)  ,S_{d}\right)
$ with fiber over a point of $\mathbb{G}\left(  P\left(  d\right)
,S_{d}\right)  $ the corresponding subspace of $\mathbb{C}^{s}$.

\begin{remark}
With respect to the basis
\[
x_{B}=x^{\alpha_{b_{1}}}\wedge...\wedge x^{\alpha_{b_{r}}}%
\]
with $B=\left\{  b_{1},...,b_{r}\right\}  \subset\left\{  1,...,s\right\}  $,
$\left\vert B\right\vert =r$ of $\bigwedge\nolimits^{r}S_{d}$ the Pl\"{u}cker
embedding is given by the $r\times r$ minors of the matrix representative%
\begin{align*}
A  &  =\left(  a_{j,\alpha_{i}}\right)  _{j=1,...r,\text{ }i=1,...,s}%
\in\mathbb{C}^{r\times s}\\
m_{j}  &  =\sum_{\left\vert \alpha\right\vert =d}a_{j,\alpha}x^{\alpha}%
\end{align*}
Denoting by $A_{B}$ the matrix formed by the columns of $A$ with indices
$b_{1},...,b_{r}$
\[%
\begin{tabular}
[c]{cccc}%
$\mathfrak{p}:$ & $\mathbb{G}\left(  P\left(  d\right)  ,S_{d}\right)  $ &
$\mathbb{\rightarrow}$ & $\mathbb{P}^{\binom{s}{r}-1}$\\
& $\left\langle m_{1},...,m_{r}\right\rangle $ & $\mapsto$ & $\left(  \det
A_{B}\mid\left\vert B\right\vert =r\right)  $%
\end{tabular}
\
\]
Note that if $A^{\prime}=UA$ is another matrix representative, then $\det
A_{B}^{\prime}=\det U\det A_{B}$ hence the homogeneous coordinates are well defined.
\end{remark}

The action of $\operatorname*{SL}\left(  V\right)  $ on $V$ gives a
representation of $\operatorname*{SL}\left(  V\right)  $ on $S_{d}%
=\operatorname*{Sym}_{d}\left(  V\right)  $ and on $W=\bigwedge\nolimits^{r}%
S_{d}$ so inducing an action%
\[
\operatorname*{SL}\left(  V\right)  \times\mathbb{P}\left(  W\right)
\rightarrow\mathbb{P}\left(  W\right)
\]
The
\index{Grassmannian}%
Grassmannian $\mathbb{G}\left(  P\left(  d\right)  ,S_{d}\right)
\hookrightarrow\mathbb{P}\left(  W\right)  $ and the
\index{Hilbert scheme}%
Hilbert scheme $\mathbb{H}_{n}^{P}\subset\mathbb{G}\left(  P\left(  d\right)
,S_{d}\right)  $ are invariant under this action.

Let $T\subset\operatorname*{SL}\left(  V\right)  $ be a maximal torus. For
$\chi\in\widehat{T}$ define the
\newsym[$\operatorname*{State}\left(  W\right)  $]{states}{}subspace%
\[
W_{\chi}=\left\{  v\in W\mid\Lambda v=\chi\left(  \Lambda\right)  v\text{
}\forall\Lambda\in T\right\}
\]
and%
\[
\operatorname*{State}\left(  W\right)  =\left\{  \chi\in\widehat{T}\mid
W_{\chi}\neq\left\{  0\right\}  \right\}
\]
so%
\[
W=\bigoplus\nolimits_{\chi\in\operatorname*{State}\left(  W\right)  }W_{\chi}%
\]
If $H\left(  I\right)  $ is the $d$-th Hilbert point of $I$ and $h^{\ast}\in
W$ is a representative of $\mathfrak{p}\left(  H\left(  I\right)  \right)  $,
then we get the corresponding decomposition%
\[
h^{\ast}=\sum\nolimits_{\chi\in\operatorname*{State}\left(  W\right)  }%
h_{\chi}\left(  I\right)
\]
with $h_{\chi}\left(  I\right)  \in W_{\chi}$. The statements $h_{\chi}\left(
I\right)  =0$ and $h_{\chi}\left(  I\right)  \neq0$ are independent of the
choice of $h^{\ast}$, as different representatives are $\mathbb{C}^{\ast}$
multiples of each other.

\begin{definition}
The
\index{state polytope|textbf}%
$d$\textbf{-th} \textbf{state polytope}
\newsym[$\operatorname*{State}\left(  I\right)  $]{state polytope}{}of $I$ is%
\[
\operatorname*{State}\left(  I\right)  =\operatorname*{convexhull}\left\{
\chi\in\operatorname*{State}\left(  W\right)  \mid h_{\chi}\left(  I\right)
\neq0\right\}  \subset\widehat{T}\otimes_{\mathbb{Z}}\mathbb{R}%
\]

\end{definition}

If the elements of $T\subset\operatorname*{SL}\left(  V\right)  $ are diagonal
with respect to the basis $x_{0},...,x_{n}$, then any
\index{one parameter subgroup}%
one parameter subgroup $\lambda:\mathbb{C}^{\ast}\rightarrow T$ is of the form%
\[
\lambda\left(  t\right)  =\operatorname*{diag}\left(  t^{w_{0}},...,t^{w_{n}%
}\right)
\]
with
\index{weight vector}%
weight vector $w=\left(  w_{0},...,w_{n}\right)  \in\mathbb{Z}^{n+1}$ and
$\sum_{i=0}^{n}w_{i}=0$. By the action of $\operatorname*{SL}\left(  V\right)
$ a one parameter subgroup $\lambda:\mathbb{C}^{\ast}\rightarrow
T\subset\operatorname*{SL}\left(  V\right)  $ assigns a weight to any monomial
of $S$ and to any
\index{Pl\"{u}cker embedding}%
Pl\"{u}cker coordinate.

\begin{definition}
A one parameter subgroup $\lambda:\mathbb{C}^{\ast}\rightarrow T$ is called
\index{m-generic}%
$d$-generic if it induces a total ordering on the monomials of $S$ of degree
less or equal to $d$.
\end{definition}

\begin{remark}
As seen in Section \ref{Sec Groebner basics} for any $d$ and any
\index{semigroup ordering}%
semigroup ordering $>$ there is a $d$-generic one parameter subgroup
$\lambda:\mathbb{C}^{\ast}\rightarrow T$ representing $>$ on the monomials of
degree at most $d$.
\end{remark}

Suppose $\lambda:\mathbb{C}^{\ast}\rightarrow T$, $\lambda\left(  t\right)
=\operatorname*{diag}\left(  t^{w_{0}},...,t^{w_{n}}\right)  $ is representing
$>$ on the monomials of degree at most $d$ and%
\[
A=\left(  a_{j,\alpha_{i}}\right)  _{j=1,...r,\text{ }i=1,...,s}\in k^{r\times
s}%
\]
is
\index{Hilbert matrix}%
the $>$-Hilbert matrix representing the $d$-th Hilbert point $H\left(
I\right)  $ of $I$ with respect to the basis%
\[
f_{j}=\sum_{i=1}^{s}a_{j,\alpha_{i}}x^{\alpha_{i}},\text{ }j=1,...,r
\]
of $I_{d}$ and the $>$-ordered basis $x^{\alpha_{1}},...,x^{\alpha_{s}}$ of
$S_{d}$. Then%
\[
in_{>}\left(  f_{j}\right)  =a_{j,\alpha_{b_{j}}}x^{\alpha_{b_{j}}}%
\]
if and only if $a_{j,\alpha_{l}}=0$ $\forall l=1,...,b_{j}-1$ and
$a_{j,\alpha_{b_{j}}}\neq0$.

Suppose that $f_{j},$ $j=1,...,r$ is a basis of $I_{d}$ such that $A$ is in
row echelon form%
\begin{equation}%
\begin{tabular}
[c]{l|lllllll}
& $x^{\alpha_{1}}$ & $\cdots$ &  &  &  & $\cdots$ & $x^{\alpha_{s}}$\\\hline
$f_{1}$ &  & $a_{1,\alpha_{b_{1}}}$ & $\cdots$ &  &  &  & \\
&  &  &  &  &  &  & \\
$\vdots$ &  &  & $a_{j,\alpha_{b_{j}}}$ & $\cdots$ &  &  & \\
&  &  &  &  &  &  & \\
$f_{r}$ &  &  &  &  & $a_{r,\alpha_{b_{r}}}$ & $\cdots$ &
\end{tabular}
\label{Equ Hilbert matrix row echelon form}%
\end{equation}
Then the
\index{Pl\"{u}cker embedding}%
Pl\"{u}cker coordinate $\mathfrak{p}\left(  H\left(  I\right)  \right)
_{B}\neq0$ for%
\[
B=\left(  b_{1},...,b_{r}\right)
\]
Note that the
\index{Pl\"{u}cker embedding}%
Pl\"{u}cker coordinates are independent of the choice of the basis of $I_{d}$.

If $B^{\prime}=\left(  b_{1}^{\prime},...,b_{r}^{\prime}\right)  \neq B$ is
some other
\index{Pl\"{u}cker embedding}%
Pl\"{u}cker coordinate with $x^{\alpha_{b_{j}^{\prime}}}>x^{\alpha_{b_{j}}}$
for some $j$, then $\mathfrak{p}\left(  H\left(  I\right)  \right)
_{B^{\prime}}=0$, hence:

\begin{lemma}
\cite{BaMo Initial Ideals and State Polytopes} Let $H\left(  I\right)
\in\mathbb{H}_{n}^{P}\subset\mathbb{G}\left(  P\left(  d\right)
,S_{d}\right)  $ be the $d$-th
\index{Hilbert point}%
Hilbert point of $I$. Fix a basis $V=$ $_{k}\left\langle x_{0},...,x_{n}%
\right\rangle $ and let $\lambda:\mathbb{C}^{\ast}\rightarrow T$,
$\lambda\left(  t\right)  =\operatorname*{diag}\left(  t^{w_{0}},...,t^{w_{n}%
}\right)  $, $w=\left(  w_{0},...,w_{n}\right)  $ be a $d$-generic one
parameter subgroup of the torus $T$.

Then there is a unique Pl\"{u}cker coordinate $x_{B}$ with $\mathfrak{p}%
\left(  H\left(  I\right)  \right)  _{B}\neq0$ such that all Pl\"{u}cker
coordinates $B^{\prime}\neq B$ with $\mathfrak{p}\left(  H\left(  I\right)
\right)  _{B^{\prime}}\neq0$ have smaller weight $w\left(  x_{B^{\prime}%
}\right)  <w\left(  x_{B}\right)  $.

If $x^{\alpha_{1}},...,x^{\alpha_{s}}$ is a $\lambda$-ordered basis of $S_{d}$
and $f_{1},...,f_{r}$ is a basis of $I_{d}$ such that the corresponding
\index{Hilbert matrix}%
Hilbert matrix of $I$ is in row echelon form as in Equation
\ref{Equ Hilbert matrix row echelon form}, then $f_{1},...,f_{r}$ form a
standard basis of $I_{\geq d}$ and%
\[
in_{>}\left(  I_{\geq d}\right)  =\left\langle x^{\alpha_{b_{1}}%
},...,x^{\alpha_{b_{r}}}\right\rangle
\]

\end{lemma}

If $\lambda:\mathbb{C}^{\ast}\rightarrow T$ is a $1$-parameter subgroup, then
$\mathbb{C}^{\ast}$ acts on $\mathbb{G}\left(  P\left(  d\right)
,S_{d}\right)  $ by%
\[%
\begin{tabular}
[c]{ccc}%
$\mathbb{C}^{\ast}\times\mathbb{G}\left(  P\left(  d\right)  ,S_{d}\right)  $
& $\rightarrow$ & $\mathbb{G}\left(  P\left(  d\right)  ,S_{d}\right)  $\\
$\left(  t,z\right)  $ & $\mapsto$ & $\lambda\left(  t\right)  z$%
\end{tabular}
\
\]
Consider $\mathbb{C}^{\ast}\hookrightarrow\mathbb{A}^{1}$ via $\mathbb{C}%
\left[  t\right]  \hookrightarrow\mathbb{C}\left[  t,t^{-1}\right]  $, so
$\mathbb{C}^{\ast}=\mathbb{A}^{1}-\left\{  0\right\}  $. If $z\in
\mathbb{G}\left(  P\left(  d\right)  ,S_{d}\right)  $ and $\mathbb{C}^{\ast
}\rightarrow\mathbb{G}\left(  P\left(  d\right)  ,S_{d}\right)  $,
$t\mapsto\lambda\left(  t\right)  z$ extends to a morphism $\mathbb{A}%
^{1}\rightarrow\mathbb{G}\left(  P\left(  d\right)  ,S_{d}\right)  $, then
call the image of $0\in\mathbb{A}^{1}$ the limit of $z$ under $\lambda$,
written $\lim_{t\rightarrow0}\lambda\left(  t\right)  z$.

\begin{lemma}
\label{Lem path groebner degeneration}Let $H\left(  I\right)  \in
\mathbb{H}_{n}^{P}\subset\mathbb{G}\left(  P\left(  d\right)  ,S_{d}\right)  $
be the $d$-th
\index{Hilbert point}%
Hilbert point of $I$ and let $\lambda:\mathbb{C}^{\ast}\rightarrow T$ be a
$d$-generic one parameter subgroup. With the action of $\operatorname*{SL}%
\left(  V\right)  $ on $\mathbb{G}\left(  P\left(  d\right)  ,S_{d}\right)  $%
\[
\lim_{t\rightarrow0}\lambda\left(  t\right)  H\left(  I\right)  =H^{\prime}%
\in\mathbb{H}_{n}^{P}\subset\mathbb{G}\left(  P\left(  d\right)
,S_{d}\right)
\]
as $\mathbb{H}_{n}^{P}$ is projective, so there is a homogeneous ideal
$I^{\prime}\subset S$ with $P_{S/I^{\prime}}=P$ such that%
\[
H^{\prime}=H\left(  I^{\prime}\right)
\]

\end{lemma}

\begin{lemma}
\cite{BaMo Initial Ideals and State Polytopes} With the setup of Lemma
\ref{Lem path groebner degeneration}%
\[
in_{>}\left(  I\right)  =I^{\prime}%
\]

\end{lemma}

\begin{proposition}
\cite{BaMo Initial Ideals and State Polytopes} The monomial initial ideals
$in_{>}\left(  I\right)  $ for all semigroup orderings $>$ correspond to the
vertices of the $d$-th state polytope $\operatorname*{State}\left(  I\right)
$.
\end{proposition}

\begin{proposition}
The Gr\"{o}bner fan $GF\left(  I\right)  $ considered as a fan in
\[
N_{\mathbb{R}}=\frac{\mathbb{R}^{n+1}}{\mathbb{R}\left(  1,...,1\right)  }%
\]
is the normal fan of the $d$-th state polytope $\operatorname*{State}\left(
I\right)  $%
\[
GF\left(  I\right)  =\operatorname*{NF}\left(  \operatorname*{State}\left(
I\right)  \right)
\]

\end{proposition}

\subsubsection{State polytope and stability}

\begin{definition}
Suppose $d\geq d_{0}$ as in Lemma \ref{Lem hilbert point degree d}, let
$H\left(  I\right)  =I_{d}\in\mathbb{G}\left(  P\left(  d\right)
,S_{d}\right)  $ be the $d$-th
\index{Hilbert point}%
Hilbert point of $I$ and let $h^{\ast}$ be a lift of $\mathfrak{p}\left(
H\left(  I\right)  \right)  $ to%
\[
W=\bigoplus\nolimits_{\chi\in\operatorname*{State}\left(  W\right)  }W_{\chi}%
\]
The ideal $I$ is called
\index{semi-stable|textbf}%
\textbf{semi-stable} if $0\notin\overline{\operatorname*{SL}\left(  V\right)
h^{\ast}}$, otherwise it is called
\index{unstable|textbf}%
\textbf{unstable}.
\end{definition}

\begin{theorem}
\cite{BaMo Initial Ideals and State Polytopes} With the setup of the previous
definition, the following conditions are equivalent:

\begin{enumerate}
\item $I$ is semi-stable.

\item For any choice of a basis $V=$ $_{k}\left\langle x_{0},...,x_{n}%
\right\rangle $ and any $1$-parameter subgroup $\lambda:\mathbb{C}^{\ast
}\rightarrow D\subset\operatorname*{SL}\left(  V\right)  $, $\lambda\left(
t\right)  =\operatorname*{diag}\left(  t^{w_{0}},...,t^{w_{n}}\right)  $,
$w=\left(  w_{0},...,w_{n}\right)  $ with $\sum_{i=0}^{n}w_{i}=0$ there are
\index{Pl\"{u}cker embedding}%
Pl\"{u}cker coordinates $x_{B}$ and $x_{B^{\prime}}$ such that $\mathfrak{p}%
\left(  H\left(  I\right)  \right)  _{B}\neq0$ and $\mathfrak{p}\left(
H\left(  I\right)  \right)  _{B^{\prime}}\neq0$ and for the corresponding
weights it holds%
\[
w\left(  x_{B}\right)  \leq0\leq w\left(  x_{B^{\prime}}\right)
\]

\item For any choice of a basis $V=$ $_{k}\left\langle x_{0},...,x_{n}%
\right\rangle $ the
\index{state polytope}%
state polytope $\operatorname*{State}\left(  I\right)  $ contains the origin.
\end{enumerate}
\end{theorem}

\subsection{Hilbert scheme and state polytope: Polarized toric setup}

\subsubsection{Linearizations\label{Sec Linearizations}}

Before reformulating and generalizing this setup, we need some general facts
about linearizations of group actions on line bundles.

If $G$ is an affine algebraic group over $K$ acting rationally on an algebraic
variety $Y$ via%
\[
\sigma:G\times Y\rightarrow Y
\]
then the pair $\left(  Y,\sigma\right)  $ is called a
\index{G-variety|textbf}%
$G$\textbf{-variety}.

\begin{definition}
If $Y$ is a $G$-variety by $\sigma:G\times Y\rightarrow Y$ and $L$ is a line
bundle on $Y$, then
\index{G-linearization}%
a
\index{linearization|textbf}%
$G$-\textbf{linearization} of $L$ is an action $\overline{\sigma}:G\times
L\rightarrow L$ such that the diagram%
\[%
\begin{tabular}
[c]{rclll}%
$G$ & $\times$ & $L$ & $\overset{\overline{\sigma}}{\rightarrow}$ & $L$\\
$id\times\pi$ & $\downarrow$ &  &  & $\downarrow\pi$\\
$G$ & $\times$ & $Y$ & $\overset{\sigma}{\rightarrow}$ & $Y$%
\end{tabular}
\]
is commutative and the action is linear on the fibers, i.e., for all $y\in Y$
the maps $\overline{\sigma}_{y}\left(  g\right)  :L_{y}\rightarrow L_{g\cdot
y}$ are linear.

The pair $\left(  L,\overline{\sigma}\right)  $ is called a
\index{G-linearized line bundle|textbf}%
$G$-\textbf{linearized line bundle}.
\end{definition}

Here $\pi:L\rightarrow Y$ denotes the projection of the total space of $L$ to
$Y$. If $g\in G$, then the group action induces an isomorphism%
\[%
\begin{tabular}
[c]{llll}%
$\sigma\left(  g\right)  :$ & $X$ & $\rightarrow$ & $X$\\
& $x$ & $\mapsto$ & $g\cdot x$%
\end{tabular}
\]
and for $y\in Y$ the maps $\overline{\sigma}_{y}\left(  g\right)
:L_{y}\rightarrow L_{g\cdot y}$ are isomorphisms of vector spaces giving an
isomorphism of line bundles%
\[
\overline{\sigma}\left(  g\right)  :L\rightarrow g^{\ast}L
\]
With%
\[
pr_{2}:G\times Y\rightarrow Y
\]
the isomorphisms of vector bundles $\overline{\sigma}\left(  g\right)  $ for
$g\in G$ form an isomorphism of vector bundles%
\[
\Phi:pr_{2}^{\ast}\left(  L\right)  \rightarrow\sigma^{\ast}\left(  L\right)
\]
Indeed also the converse is true:

\begin{lemma}
\cite[Knop, Kraft, Luna and Vust, Sec. 4]{KSS Algebraische
Transformationsgruppen und Invariantentheorie} If $G$ is a connected affine
algebraic group, $Y$ is a $G$-variety and $L$ is a line bundle on $Y$, then
$L$ has a $G$-linearization if and only if there is an isomorphism of line
bundles%
\[
\Phi:pr_{2}^{\ast}\left(  L\right)  \rightarrow\sigma^{\ast}\left(  L\right)
\]

\end{lemma}

The set of $G$-bundles on the $G$-variety $X$ carries the structure
\newsym[$\operatorname*{Pic}^{G}\left(  Y\right)  $]{group of $G$-bundles}{}of
an abelian group, the \textbf{group of }$G$\textbf{-bundles}
\index{group of G-bundles|textbf}%
by $\operatorname*{Pic}^{G}\left(  Y\right)  $: If $L$ and $L^{\prime}$ are
$G$-bundles with linearizations given by the isomorphisms $\Phi:pr_{2}^{\ast
}\left(  L\right)  \rightarrow\sigma^{\ast}\left(  L\right)  $ and
$\Phi^{\prime}:pr_{2}^{\ast}\left(  L^{\prime}\right)  \rightarrow\sigma
^{\ast}\left(  L^{\prime}\right)  $, then on $L\otimes L^{\prime}$ a
$G$-linearization is given by the isomorphism%
\[%
\begin{tabular}
[c]{lccc}%
$\Phi\otimes\Phi^{\prime}:$ & $pr_{2}^{\ast}\left(  L\otimes L^{\prime
}\right)  $ & $\rightarrow$ & $\sigma^{\ast}\left(  L\otimes L^{\prime
}\right)  $\\
& $\shortparallel$ &  & $\shortparallel$\\
& $pr_{2}^{\ast}\left(  L\right)  \otimes pr_{2}^{\ast}\left(  L^{\prime
}\right)  $ & $\rightarrow$ & $\sigma^{\ast}\left(  L\right)  \otimes
\sigma^{\ast}\left(  L^{\prime}\right)  $%
\end{tabular}
\]
The neutral element of $\operatorname*{Pic}^{G}\left(  Y\right)  $ is the line
bundle $Y\times K\rightarrow Y$ with the $G$-linearization%
\[
\sigma\times id:G\times Y\times K\rightarrow Y\times K
\]
If $L$ is a $G$-bundle with linearization given by $\Phi:pr_{2}^{\ast}\left(
L\right)  \rightarrow\sigma^{\ast}\left(  L\right)  $, then its inverse is
$L^{\ast}$ with the linearization%
\[
\left(  \Phi^{\ast}\right)  ^{-1}:pr_{2}^{\ast}\left(  L^{\ast}\right)
\rightarrow\sigma^{\ast}\left(  L^{\ast}\right)
\]

The map%
\[
\alpha:\operatorname*{Pic}\nolimits^{G}\left(  Y\right)  \rightarrow
\operatorname*{Pic}\left(  Y\right)
\]
forgetting the linearization is a homomorphism.

\begin{proposition}
\cite[Sec. 7]{Dolgachev Lectures on Invariant Theory} If $Y$ is connected and
proper over $\mathbb{C}$, then%
\[
\ker\left(  \alpha\right)  \cong\chi\left(  G\right)
\]

\end{proposition}

\begin{lemma}
\cite[Knop, Kraft, Luna and Vust, Sec. 4]{KSS Algebraische
Transformationsgruppen und Invariantentheorie} If $G$ is a connected affine
algebraic group, $Y$ is a normal $G$-variety and $E$ is a line bundle on
$G\times Y$, then for all $y_{0}\in Y$%
\[
E\cong pr_{1}^{\ast}\left(  L\mid_{G\times\left\{  y_{0}\right\}  }\right)
\otimes pr_{2}^{\ast}\left(  L\mid_{\left\{  e\right\}  \times Y}\right)
\]

\end{lemma}

So if $y_{0}\in Y$, define the homomorphism%
\[%
\begin{tabular}
[c]{llll}%
$\delta:$ & $\operatorname*{Pic}\left(  Y\right)  $ & $\rightarrow$ &
$\operatorname*{Pic}\left(  G\right)  $\\
& $L$ & $\mapsto$ & $pr_{2}^{\ast}\left(  L\right)  \otimes\sigma^{\ast
}\left(  L^{\ast}\right)  \mid_{G\times y_{0}}$%
\end{tabular}
\
\]
which has $\ker\left(  \delta\right)  =\operatorname*{image}\left(
\alpha\right)  $.

\begin{theorem}
\cite[Sec. 7]{Dolgachev Lectures on Invariant Theory} If $G$ is a connected
affine algebraic group and $Y$ is a normal $G$-variety, then the sequence%
\[
0\rightarrow\ker\left(  \alpha\right)  \rightarrow\operatorname*{Pic}%
\nolimits^{G}\left(  Y\right)  \rightarrow\operatorname*{Pic}\left(  Y\right)
\rightarrow\operatorname*{Pic}\left(  G\right)
\]
is exact.
\end{theorem}

If $G$ is a connected affine algebraic group, then $\operatorname*{Pic}\left(
G\right)  $ is finite, see\linebreak\cite[Knop, Kraft, Luna and Vust, Prop.
4.5]{KSS Algebraische Transformationsgruppen und Invariantentheorie}, so:

\begin{remark}
$\operatorname*{Pic}\nolimits^{G}\left(  Y\right)  $ has finite index in
$\operatorname*{Pic}\left(  Y\right)  $, hence for all $L\in
\operatorname*{Pic}\left(  Y\right)  $ there is an $m$ such that $L^{\otimes
m}$ is a $G$-bundle.

If $G$ is $\operatorname*{GL}\left(  n,\mathbb{C}\right)  $ or a torus
$\left(  \mathbb{C}^{\ast}\right)  ^{n}$ or $\operatorname*{SL}\left(
n,\mathbb{C}\right)  $, then $\operatorname*{Pic}\left(  G\right)  =0$.

Hence if $T=\left(  \mathbb{C}^{\ast}\right)  ^{n}$ is a torus and $Y$ is a
$T$-variety, then we have an exact sequence%
\[
0\rightarrow\widehat{T}\rightarrow\operatorname*{Pic}\nolimits^{T}\left(
Y\right)  \rightarrow\operatorname*{Pic}\left(  Y\right)  \rightarrow0
\]
so any line bundle $L$ on $Y$ has a $T$-linearization and any two
linearizations differ by a translation in the lattice $\widehat{T}%
\cong\mathbb{Z}^{n}$.
\end{remark}

\begin{remark}
If $Y$ is a toric variety with torus $T$, then the sheaf of Zariski
differential forms $\Omega_{Y}^{p}$ has a canonical linearization given by the
pullback of differential forms with respect to the isomorphism%
\[%
\begin{tabular}
[c]{llll}%
$\sigma\left(  g\right)  :$ & $X$ & $\rightarrow$ & $X$\\
& $x$ & $\mapsto$ & $g\cdot x$%
\end{tabular}
\]

\end{remark}

\subsubsection{Setup for subvarieties of a projective toric variety}

Let $Y=X\left(  \Sigma\right)  $ be a simplicial toric variety of dimension
$n$ given by the fan $\Sigma\subset N_{\mathbb{R}}=N\otimes_{\mathbb{Z}%
}\mathbb{R}$ with $N\cong\mathbb{Z}^{n}$ and let $L=\mathcal{O}_{Y}\left(
D\right)  $ be a very ample line bundle on $Y$. The lattice
\[
M=\operatorname*{Hom}\left(  N,\mathbb{Z}\right)  =\operatorname*{Hom}\left(
T,\mathbb{C}^{\ast}\right)  =\widehat{T}%
\]
is the character group of the torus $T\subset Y$. Let $S=\mathbb{C}\left[
y_{v}\mid v\in\Sigma\left(  1\right)  \right]  $ be the Cox ring of $Y$,%
\[
0\rightarrow M\overset{A}{\rightarrow}\mathbb{Z}^{\Sigma\left(  1\right)
}\rightarrow A_{n-1}\left(  Y\right)  \rightarrow0
\]
the presentation of the Chow group and
\begin{gather*}
1\rightarrow G\left(  \Sigma\right)  \rightarrow\left(  \mathbb{C}^{\ast
}\right)  ^{\Sigma\left(  1\right)  }\rightarrow T\rightarrow1\\
G\left(  \Sigma\right)  =\operatorname*{Hom}\nolimits_{\mathbb{Z}}\left(
A_{n-1}\left(  X\left(  \Sigma\right)  \right)  ,\mathbb{C}^{\ast}\right)
\end{gather*}
the corresponding sequence involving the tori $\left(  \mathbb{C}^{\ast
}\right)  ^{\Sigma\left(  1\right)  }\ $and $T$.

By Section \ref{Sec Linearizations} there is a linearization of the action of
$T$ on $Y$ on the line bundle $L$%
\[%
\begin{tabular}
[c]{rclll}%
$T$ & $\times$ & $L$ & $\overset{\overline{\sigma}}{\rightarrow}$ & $L$\\
$id\times\pi$ & $\downarrow$ &  &  & $\downarrow\pi$\\
$T$ & $\times$ & $Y$ & $\overset{\sigma}{\rightarrow}$ & $Y$%
\end{tabular}
\]

With%
\[
V=H^{0}\left(  Y,L\right)  =S_{\left[  D\right]  }%
\]
the line bundle $L$ defines an embedding
\begin{gather*}
\phi_{V}:Y\rightarrow\mathbb{P}\left(  V^{\ast}\right) \\
\phi_{V}\left(  y\right)  =\left\{  s\in V\mid s\left(  y\right)  =0\right\}
\end{gather*}
identifying elements of $\mathbb{P}\left(  V^{\ast}\right)  $ with hyperplanes
in $V$. The map $\phi_{V}$ is $T$-equivariant with respect to the action%
\begin{gather*}
T\times\mathbb{P}\left(  V^{\ast}\right)  \rightarrow\mathbb{P}\left(
V^{\ast}\right) \\
g\cdot H=g^{-1}\left(  H\right)
\end{gather*}

The toric variety $Y$ embedded by $\phi_{V}$ is isomorphic to the projective
toric variety%

\[
\mathbb{P}\left(  \Delta_{D}\right)  =\operatorname*{Proj}S\left(  \Delta
_{D}\right)
\]
and the polytope ring $S\left(  \Delta_{D}\right)  $ is isomorphic to%
\[
S\left(  \Delta_{D}\right)  \cong\bigoplus_{d=0}^{\infty}S_{d\left[  D\right]
}\subset S
\]
which is $\mathbb{Z}_{\geq0}$-graded by $d$.

If $I\subset S$ is a homogeneous ideal, then it corresponds under the
embedding $\phi_{V}$ of $Y$ via $L=\mathcal{O}_{Y}\left(  D\right)  $ to the
ideal%
\[
I_{\mathbb{N}\left[  D\right]  }=\bigoplus_{d=0}^{\infty}I_{d\left[  D\right]
}\subset S\left(  \Delta_{D}\right)
\]
The Hilbert function%
\[
h_{S/I}\left(  k\right)  =\dim_{\mathbb{C}}\left(  S_{k\left[  D\right]
}/I_{k\left[  D\right]  }\right)
\]
agrees for large $k$ with a polynomial $P$, the Hilbert polynomial of $S/I$
under the embedding of $Y$ given by $L$.

Furthermore via this embedding there is a $d_{0}$ such that for any
homogeneous ideal $J\subset S\left(  \Delta_{D}\right)  $ with Hilbert
polynomial $P$ under the embedding of $Y$ given by $L$
\[
\left(  J_{\mathbb{N}\left[  D\right]  }:B\left(  \Sigma\right)
_{\mathbb{N}\left[  D\right]  }^{\infty}\right)  _{\geq d_{0}}%
\]
is generated by the Hilbert point%
\[
J_{d_{0}\left[  D\right]  }\in\mathbb{G=G}\left(  P\left(  d_{0}\right)
,S_{d_{0}\left[  D\right]  }\right)
\]
and these points form the $P$-th Hilbert scheme%
\[
\mathbb{H}_{L}^{P}\subset\mathbb{G}%
\]
in the embedding of $Y$ via $L$.

The action of $T$ on $V=H^{0}\left(  Y,L\right)  =S_{\left[  D\right]  }$
induces an action of $T$ on%
\[
W=\bigwedge\nolimits^{r}S_{d_{0}\left[  D\right]  }%
\]
with $r=\dim S_{d_{0}\left[  D\right]  }-P\left(  d_{0}\right)  $, and the
Pl\"{u}cker embedding
\[
\mathfrak{p}:\mathbb{G\rightarrow P}\left(  W\right)
\]
is $T$-equivariant. With
\[
W_{\chi}=\left\{  v\in W\mid\Lambda v=\chi\left(  \Lambda\right)  v\text{
}\forall\Lambda\in T\right\}
\]
for $\chi\in\widehat{T}=M$ and%
\[
\operatorname*{State}\nolimits_{L}\left(  W\right)  =\left\{  \chi\in M\mid
W_{\chi}\neq\left\{  0\right\}  \right\}
\]
we get a decomposition%
\[
W=\bigoplus\nolimits_{\chi\in\operatorname*{State}\nolimits_{L}\left(
W\right)  }W_{\chi}%
\]
If $h^{\ast}\in W$ is a representative of $\mathfrak{p}\left(  I_{d_{0}\left[
D\right]  }\right)  $, there is the corresponding decomposition%
\[
h^{\ast}=\sum\nolimits_{\chi\in\operatorname*{State}\nolimits_{L}\left(
W\right)  }h_{\chi}%
\]
with $h_{\chi}\in W_{\chi}$. With%
\[
\operatorname*{State}\nolimits_{L}\left(  h\right)  =\left\{  \chi\in M\mid
h_{\chi}\neq0\right\}
\]
define the state polytope of $I$ with respect to $L$ as the convex hull
\[
\operatorname*{State}\nolimits_{L}\left(  I\right)
=\operatorname*{convexhull}\left(  \operatorname*{State}\nolimits_{L}\left(
h\right)  \right)  \subset\widehat{T}\otimes_{\mathbb{Z}}\mathbb{R}%
=M_{\mathbb{R}}%
\]

\subsubsection{Hilbert-Mumford stability}

Suppose $G$ is a reductive group and $Y$ is an irreducible $G$-variety.

\begin{definition}
Let $L$ be a $G$-bundle on $Y$ and $y\in Y$.

\begin{enumerate}
\item $y$ is called
\index{semi-stable|textbf}%
\textbf{semi-stable with respect to }$L$ if there is an $a>0$ and
\newsym[$Y^{ss}\left(  L\right)  $]{semistable points}{}an $\alpha\in
H^{0}\left(  Y,L^{a}\right)  ^{G}$ such that
\[
Y_{\alpha}=\left\{  y\in Y\mid\alpha\left(  y\right)  \neq0\right\}
\]
is affine and $y\in Y_{\alpha}$.

\item $y$ is called
\index{unstable|textbf}%
\textbf{unstable with respect to }$L$ if
\newsym[$Y^{us}\left(  L\right)  $]{unstable points}{}it is not semi-stable
with respect to $L$.

\item $y$ is called
\index{stable|textbf}%
\textbf{stable with respect to }$L$ if the
\newsym[$Y^{s}\left(  L\right)  $]{stable points}{}isotropy group $G_{y}$ is
finite and the $G$-orbits in $Y_{\alpha}$ are closed.
\end{enumerate}

Denote by $Y^{ss}\left(  L\right)  ,Y^{us}\left(  L\right)  $ and
$Y^{s}\left(  L\right)  $ the set of semi-stable, unstable and stable points
of $Y$, respectively.
\end{definition}

\begin{lemma}
\cite[Sec. 8]{Dolgachev Lectures on Invariant Theory} With the notation from above:

The sets $Y^{ss}\left(  L\right)  ,Y^{us}\left(  L\right)  $ and $Y^{s}\left(
L\right)  $ do not change when replacing $L$ by $L^{a}$ for $a\in
\mathbb{Z}_{>0}$.

If $L$ is ample and $Y$ is projective, then $Y_{\alpha}$ is always affine.
\end{lemma}

\subsubsection{Stability on a variety with a torus action}

Suppose $X$ is a projective variety, the torus $T$ acts on $X$ and $E$ is a
very ample $T$-linearized line bundle on $X$. Let $W=H^{0}\left(  X,E\right)
$, $s=\dim_{\mathbb{C}}W$ and let $\phi_{W}:X\rightarrow\mathbb{P}\left(
W\right)  $ be the corresponding embedding. So $T$ acts on $X$ via a linear
representation%
\[
T\rightarrow\operatorname*{GL}\left(  W^{\ast}\right)
\]

If $x\in X$ and $x^{\ast}$ is a representative of $\phi_{W}\left(  x\right)
$, then
\[
x\in X^{us}\left(  E\right)  \Leftrightarrow0\in\overline{T\cdot x^{\ast}}%
\]
So if $0\in\overline{\lambda\left(  \mathbb{C}^{\ast}\right)  \cdot x^{\ast}}$
for some one parameter subgroup $\lambda\in\widehat{T}^{\ast}$ of $T$
\[
\lambda:\mathbb{C}^{\ast}\rightarrow T
\]
then $x$ is unstable.

Choosing a basis of $W$ such that $T$ acts via diagonal matrices write%
\[
x^{\ast}=\left(  x_{1},...,x_{s}\right)
\]
with respect to this basis. Then%
\[
\lambda\left(  t\right)  x^{\ast}=\left(  t^{\beta_{1}}x_{1},...,t^{\beta_{s}%
}x_{s}\right)
\]
with some $\beta_{i}$.

\begin{itemize}
\item If $\beta_{i}>0$ for all $i$ with $x_{i}\neq0$, then
\[%
\begin{tabular}
[c]{llll}%
$\lambda_{x^{\ast}}:$ & $\mathbb{A}^{1}\backslash\left\{  0\right\}  $ &
$\rightarrow$ & $\mathbb{A}^{s}$\\
& \multicolumn{1}{c}{$t$} & $\mapsto$ & $\lambda\left(  t\right)  x^{\ast}$%
\end{tabular}
\]
extends to a regular map%
\[%
\begin{tabular}
[c]{lll}%
$\mathbb{A}^{1}$ & $\rightarrow$ & $\mathbb{A}^{s}$\\
\multicolumn{1}{c}{$t$} & $\mapsto$ & $\lambda\left(  t\right)  x^{\ast}$ for
$t\neq0$\\
\multicolumn{1}{c}{$0$} & $\mapsto$ & $0$%
\end{tabular}
\]
so $0\in\overline{\lambda\left(  \mathbb{C}^{\ast}\right)  \cdot x^{\ast}}$,
hence $x$ is unstable.

\item If $\beta_{i}<0$ for all $i$ with $x_{i}\neq0$, then above argument
applied to $\lambda^{-1}$ shows that $x$ is unstable.
\end{itemize}

Define%
\[
\mu^{E}\left(  x,\lambda\right)  =\min\left\{  \beta_{i}\mid x_{i}%
\neq0\right\}
\]
so if $\mu^{E}\left(  x,\lambda\right)  >0$, then $x\in X^{us}\left(
E\right)  $, hence%
\[
x\in X^{ss}\left(  E\right)  \Rightarrow\mu^{E}\left(  x,\lambda\right)
\leq0\text{ }\forall\lambda\in\widehat{T}^{\ast}%
\]

On the other hand if $\mu^{E}\left(  x,\lambda\right)  \leq0$ $\forall
\lambda\in\widehat{T}^{\ast}$ and there is a $\lambda\in\widehat{T}^{\ast}$
with $\mu^{E}\left(  x,\lambda\right)  =0$, then $y^{\ast}=\left(
y_{i}\right)  $ with%
\[
y_{i}=\left\{
\begin{tabular}
[c]{ll}%
$0$ & $\text{ if }x_{i}\neq0\text{ and }\beta_{i}>0$\\
$x_{i}$ & $\text{otherwise}$%
\end{tabular}
\ \right\}
\]
is in the closure of $\lambda\left(  \mathbb{C}^{\ast}\right)  \cdot x^{\ast}%
$, i.e.,%
\[
y^{\ast}\in\overline{\lambda\left(  \mathbb{C}^{\ast}\right)  \cdot x^{\ast}}%
\]
If $x$ would be stable, then it would have to hold that $y^{\ast}\in
\lambda\left(  \mathbb{C}^{\ast}\right)  \cdot x^{\ast}$, but this is not
possible as%
\[
\lambda\left(  t\right)  \cdot y^{\ast}=y^{\ast}\text{ for all }t\in
\mathbb{C}^{\ast}%
\]
hence%
\[
x\in X^{s}\left(  E\right)  \Rightarrow\mu^{E}\left(  x,\lambda\right)
<0\text{ }\forall\lambda\in\widehat{T}^{\ast}%
\]
indeed both statements are characterizations of the semi stable and stable points:

\begin{theorem}
\cite{Dolgachev Lectures on Invariant Theory}
\label{Thm Hilbert-Mumford numerical stability}With the setup from above%
\[%
\begin{tabular}
[c]{lll}%
$x\in X^{ss}\left(  E\right)  $ & $\Leftrightarrow$ & $\mu^{E}\left(
x,\lambda\right)  \leq0\text{ }\forall\lambda\in\widehat{T}^{\ast}$\\
$x\in X^{s}\left(  E\right)  $ & $\Leftrightarrow$ & $\mu^{E}\left(
x,\lambda\right)  <0\text{ }\forall\lambda\in\widehat{T}^{\ast}$%
\end{tabular}
\]

\end{theorem}

\subsubsection{State polytope and Stability}

The bilinear pairing between characters and one parameter subgroups of $T$%
\[%
\begin{tabular}
[c]{lll}%
$\widehat{T}\times\widehat{T}^{\ast}$ & $\rightarrow$ & $\widehat
{\mathbb{C}^{\ast}}=\mathbb{Z}$\\
$\left(  \chi,\lambda\right)  $ & $\mapsto$ & $\left\langle \chi
,\lambda\right\rangle =\chi\circ\lambda$%
\end{tabular}
\]
corresponds via the identification $\widehat{T}=M$ and $\widehat{T}^{\ast}=N$
to the canonical bilinear pairing
\[
\left\langle -,-\right\rangle :M\times N\rightarrow\mathbb{Z}%
\]

Fix a $T$-invariant basis $x_{0},...,x_{n}$ of $V$ and let%
\[
x_{B}=x_{b_{1}}\wedge...\wedge x_{b_{r}}%
\]
with $B=\left\{  b_{1},...,b_{r}\right\}  \subset\left\{  0,...,n\right\}  $,
$\left\vert B\right\vert =r$ be the corresponding $T$-invariant basis of $W$,
which is compatible with the decomposition%
\[
W=\bigoplus\nolimits_{\chi\in\operatorname*{State}_{L}\left(  W\right)
}W_{\chi}%
\]
With respect to the basis $\left(  x_{B}\right)  $ the representation%
\[
\rho:T\rightarrow\operatorname*{GL}\left(  W\right)
\]
given by the action $T\times W\rightarrow W$ is of the form%
\[
\rho\left(  x\right)  =\operatorname*{diag}\left(  x^{m_{1}},...,x^{m_{\dim
W}}\right)
\]
with $m_{i}\in M$.

If%
\begin{align*}
\lambda &  :\mathbb{C}^{\ast}\rightarrow T\\
\lambda\left(  t\right)   &  =\operatorname*{diag}\left(  t^{w_{1}%
},...,t^{w_{n}}\right)
\end{align*}
is a one parameter subgroup of $T$, then the composition is%
\[%
\begin{tabular}
[c]{llll}%
$\rho\circ\lambda:$ & $\mathbb{C}^{\ast}$ & $\rightarrow$ &
$\operatorname*{GL}\left(  W\right)  $\\
& \multicolumn{1}{c}{$t$} & $\mapsto$ & $\operatorname*{diag}\left(
t^{\left\langle w,m_{1}\right\rangle },...,t^{\left\langle w,m_{\dim
W}\right\rangle }\right)  $%
\end{tabular}
\]

Now suppose $h^{\ast}\in W$ is a representative of the image of the Hilbert
point $I_{d_{0}\left[  D\right]  }$ under the Pl\"{u}cker embedding
$\mathfrak{p}$, and write%
\[
h^{\ast}=\sum\nolimits_{\chi\in\operatorname*{State}_{L}\left(  W\right)
}h_{\chi}%
\]
with $h_{\chi}\in W_{\chi}$. Write%
\[
h^{\ast}=\left(  \alpha_{1},...,\alpha_{\dim W}\right)
\]
with respect to the basis $\left(  x_{B}\right)  $. So%
\[
\lambda\left(  t\right)  \cdot h^{\ast}=\operatorname*{diag}\left(
t^{\left\langle w,m_{1}\right\rangle }\alpha_{1},...,t^{\left\langle w,m_{\dim
W}\right\rangle }\alpha_{\dim W}\right)
\]
hence with the line bundle%
\[
E=\overline{\mathfrak{p}}^{\ast}\left(  \mathcal{O}_{\mathbb{P}\left(
W\right)  }\left(  1\right)  \right)
\]
where $\overline{\mathfrak{p}}:\mathbb{H}_{L}^{P}\mathbb{\rightarrow P}\left(
W\right)  $ is the embedding of the Hilbert scheme induced by the Pl\"{u}cker
embedding, we have%
\[
\mu^{E}\left(  h,\lambda\right)  =\min\left\{  \left\langle w,m_{i}%
\right\rangle \mid\alpha_{i}\neq0\right\}  =\min_{\chi\in\operatorname*{State}%
_{L}\left(  h\right)  }\left\langle \chi,\lambda\right\rangle
\]
so by Theorem \ref{Thm Hilbert-Mumford numerical stability} we obtain:

\begin{theorem}
Suppose $Y=X\left(  \Sigma\right)  $ is a toric variety given by the fan
$\Sigma\subset N_{\mathbb{R}}$, $L=\mathcal{O}_{Y}\left(  D\right)  $ is a
very ample $T$-line bundle on $Y$ and $S$ is the Cox ring of $Y$. If $I\subset
S$ is homogeneous, then stability and semi-stability of the Hilbert point
\[
H_{L}\left(  I\right)  \in\mathbb{H=H}_{L}^{P}%
\]
are characterized via the state polytope $\operatorname*{State}\nolimits_{L}%
\left(  I\right)  \subset M_{\mathbb{R}}$ as follows:
\[%
\begin{tabular}
[c]{lll}%
$H_{L}\left(  I\right)  \in\mathbb{H}^{ss}\left(  E\right)  $ &
$\Leftrightarrow$ & $0\in\operatorname*{State}\nolimits_{L}\left(  I\right)
$\\
$H_{L}\left(  I\right)  \in\mathbb{H}^{s}\left(  E\right)  $ &
$\Leftrightarrow$ & $0\in\operatorname*{int}\left(  \operatorname*{State}%
\nolimits_{L}\left(  I\right)  \right)  $%
\end{tabular}
\]

\end{theorem}

\subsection{Hilbert scheme and state polytope: Cox homogeneous
setup\label{Sec Hilbert scheme state polytope Cox homogeneou}}

\subsubsection{Grassmann functor\label{Sec Grassmann functor}}

Let $k$ be a commutative ring.

\begin{definition}
Let $N$ be a finitely generated $k$-module.
\newsym[$\mathbb{G}_{N}^{r}$]{Grassmann functor}{}The
\index{Grassmann functor|textbf}%
\textbf{Grassmann functor} $\mathbb{G}_{N}^{r}:$\underline{$k-Alg$%
}$\rightarrow$\underline{$Set$} is defined as%
\[
\mathbb{G}_{N}^{r}\left(  R\right)  =\left\{  L\mid L\subset R\otimes N\text{
submodule with }\left(  R\otimes N\right)  /L\text{ locally free of rank
}r\right\}
\]

\end{definition}

An $R$-module $W$ is locally free of rank $r$ if there are $f_{1},...,f_{k}\in
R$ with $\left\langle f_{1},...,f_{k}\right\rangle =\left\langle
1\right\rangle \subset R$ such that $W_{f_{j}}\cong R_{f_{j}}^{r}$ for all
$j=1,...,k$.

The Grassmann functor $\mathbb{G}_{N}^{r}$ is represented by the
\index{Grassmann scheme|textbf}%
\textbf{Grassmann scheme} $\mathbb{G}_{N}^{r}$ described in coordinates as follows:

\begin{itemize}
\item If $N=k^{m}$:

Let $v_{1},...,v_{m}$ be a basis of $N$ and let $B=\left\{  v_{i_{1}%
},...,v_{i_{r}}\right\}  $. The subfunctor $\mathbb{G}_{k^{m}\backslash B}%
^{r}$ of $\mathbb{G}_{k^{m}}^{r}$ is defined as%
\[
\mathbb{G}_{k^{m}\backslash B}^{r}\left(  R\right)  =\left\{  L\mid L\subset
R^{m}\text{ submodule with }R^{m}/L\text{ free with basis }B\right\}
\]
It is represented by the affine space $\mathbb{A}^{r\left(  m-r\right)  }$
associating $L\in\mathbb{G}_{k^{m}\backslash B}^{r}\left(  R\right)  $ to
$\left(  \lambda_{j}^{i}\right)  $ with%
\[
R^{m}/L\ni\overline{v_{i}}=\sum_{j=1}^{r}\lambda_{j}^{i}\overline{v_{i_{j}}%
}\text{ for }i\notin\left\{  i_{1},...,i_{r}\right\}
\]
Via the Pl\"{u}cker embedding the Grassmann functor $\mathbb{G}_{k^{m}}^{r}$
is represented by a projective scheme covered by affine open subsets
representing $\mathbb{G}_{k^{m}\backslash B}^{r}$.

\item If $N=k^{m}/J$ is a finitely generated $k$-module:

Then $R\otimes N\cong R^{m}/RJ$ for any $k$-algebra $R$ and
\[
\mathbb{G}_{N}^{r}\left(  R\right)  =\left\{  L\in\mathbb{G}_{k^{m}}%
^{r}\left(  R\right)  \mid RJ\subset L\right\}
\]
If $v_{1},...,v_{m}$ is a basis of $k^{m}$ and $B=\left\{  v_{i_{1}%
},...,v_{i_{r}}\right\}  $, then the subfunctor $\mathbb{G}_{k^{m}\backslash
B}^{r}\cap\mathbb{G}_{N}^{r}$ of $\mathbb{G}_{k^{m}}^{r}$ is represented by
the subscheme%
\[
\left\{  \left(  \lambda_{j}^{i}\right)  \mid a_{i_{j}}^{u}+\sum
_{i\notin\left\{  i_{1},...,i_{r}\right\}  }a_{i}^{u}\lambda_{j}^{i}=0\text{
}\forall u\in J\text{ }\forall j=1,...,r\right\}  \subset\mathbb{A}^{r\left(
m-r\right)  }%
\]
where for $u\in J$ the $a_{i}^{u}\in k$ are defined by%
\[
u=\sum_{i=1}^{m}a_{i}^{u}v_{i}%
\]

\end{itemize}

\begin{proposition}
\cite{HS Multigraded Hilbert schemes} If $N$ is a finitely generated
$k$-module, then the functor $\mathbb{G}_{N}^{r}$ is represented by a closed
subscheme of the scheme representing $\mathbb{G}_{k^{m}}^{r}$.
\end{proposition}

Let $N$ be a finitely generated $k$-module and $M\subset N$ a submodule and
consider the subfunctor $\mathbb{G}_{N\backslash M}^{r}\subset\mathbb{G}%
_{N}^{r}$%
\begin{align*}
\mathbb{G}_{N\backslash M}^{r}\left(  R\right)   &  =\left\{  L\in
\mathbb{G}_{N}^{r}\left(  R\right)  \mid\left(  R\otimes N\right)  /L\text{
locally free with bases in }M\right\} \\
&  =\left\{  L\in\mathbb{G}_{N}^{r}\left(  R\right)  \mid%
\begin{tabular}
[c]{l}%
$\exists f_{1},...,f_{k}\in R$ with $\left\langle f_{1},...,f_{k}\right\rangle
=\left\langle 1\right\rangle $\\
such that $\left(  \left(  R\otimes N\right)  /L\right)  _{f_{j}}\text{ has a
basis in }M$%
\end{tabular}
\right\} \\
&  =\left\{  L\in\mathbb{G}_{N}^{r}\left(  R\right)  \mid M\text{ generates
}\left(  R\otimes N\right)  /L\right\}
\end{align*}

\begin{proposition}
\cite{HS Multigraded Hilbert schemes} $\mathbb{G}_{N\backslash M}^{r}$ is
represented by an open subscheme of the scheme representing $\mathbb{G}%
_{N}^{r}$, so by a
\newsym[$\mathbb{G}_{N\backslash M}^{r}$]{relative Grassmann functor}{}quasiprojective
scheme over $k$. It is called the
\index{relative Grassmann functor|textbf}%
\textbf{relative Grassmann functor of }$M\subset N$.
\end{proposition}

If $A$ is a \newsym[$\mathbb{G}_{N}^{h}$]{graded Grassmann functor}{}finite
set and $N=\bigoplus_{a\in A}N_{a}$ is a finitely generated graded $k$-module
and $h:A\rightarrow\mathbb{N}$ is some function, then the
\index{graded Grassmann functor}%
\textbf{graded Grassmann functor} $\mathbb{G}_{N}^{h}$ is defined as%
\[
\mathbb{G}_{N}^{h}\left(  R\right)  =\left\{  L\mid%
\begin{tabular}
[c]{l}%
$L\subset R\otimes N$ homogeneous submodule $\text{with}$\\
$\left(  R\otimes N_{a}\right)  /L_{a}\text{ locally free of rank }h\left(
a\right)  $ $\forall a\in A$%
\end{tabular}
\right\}
\]
and $\mathbb{G}_{N}^{h}$ is naturally isomorphic to $%
{\textstyle\prod\nolimits_{a\in A}}
\mathbb{G}_{N_{a}}^{h\left(  a\right)  }$ hence is projective.

If $M\subset N$ a homogeneous
\newsym[$\mathbb{G}_{N\backslash M}^{h}$]{relative graded Grassmann functor}{}submodule
the
\index{relative graded Grassmann functor|textbf}%
\textbf{relative graded Grassmann functor of }$M\subset N$ is defined by
\[
\mathbb{G}_{N\backslash M}^{h}\left(  R\right)  =\left\{  L\in\mathbb{G}%
_{N}^{h}\left(  R\right)  \mid\left(  R\otimes N_{a}\right)  /L_{a}\text{
locally free with bases in }M\text{ }\forall a\in A\right\}
\]
and is represented by a quasiprojective scheme over $k$.

$\mathbb{G}_{N}^{h}$ and $\mathbb{G}_{N\backslash M}^{h}$ are subfunctors of
$\mathbb{G}_{N}^{r}$ respectively $\mathbb{G}_{N\backslash M}^{r}$ with
$r=\sum_{a\in A}h\left(  a\right)  $ and the corresponding morphisms of
schemes are closed embeddings.

\subsubsection{Hilbert functor\label{Sec Hilbert functor}}

Let $k$ be a commutative ring, $A\ $a set and let
\[
S=\bigoplus_{a\in A}S_{a}%
\]
be a graded $k$-module. For all $a,b\in A$ let $F_{a,b}\subset
\operatorname*{Hom}_{k}\left(  S_{a},S_{b}\right)  $ a subset such that
$F_{bc}\circ F_{a,b}\subset F_{a,c}$ $\forall a,b\in A$ and $id_{S_{a}}\in
F_{a,a}$ $\forall a\in A$ and call $F=\bigcup_{a,b\in A}F_{a,b}$ a set of
\index{operator|textbf}%
operators on $S$. So $\left(  S,F\right)  $ is a small category of $k$-modules.

If $R$ is a $k$-algebra, then%
\[
R\otimes S=\bigoplus_{a\in A}R\otimes S_{a}%
\]
is a graded $R$-module with operators%
\[
F_{a,b}^{R}=\left(  1_{R}\otimes_{k}\_\right)  \left(  F_{a,b}\right)
=\left\{  1_{R}\otimes_{k}f\mid f\in F_{a,b}\right\}
\]

A homogeneous submodule $L=\bigoplus_{a\in A}L_{a}\subset R\otimes S$ is
called an
\index{F-submodule|textbf}%
$F$\textbf{-submodule} if $F_{a,b}^{R}\left(  L_{a}\right)  \subset L_{b}$ for
all $a,b\in A$.

\begin{definition}
If $h:A\rightarrow\mathbb{N}$ is a function and $R$ is a $k$-algebra, then
define%
\[
\mathbb{H}_{\left(  S,F\right)  }^{h}\left(  R\right)  =\left\{  L\mid%
\begin{tabular}
[c]{l}%
$L\subset R\otimes S$ is an $F$-submodule $\text{with}$\\
$\left(  R\otimes S_{a}\right)  /L_{a}\text{ locally free of rank }h\left(
a\right)  $ $\forall a\in A$%
\end{tabular}
\right\}
\]
If $\phi:R\rightarrow R^{\prime}$ is a homomorphism of commutative rings and
$L\in\mathbb{H}_{\left(  S,F\right)  }^{h}\left(  R\right)  $, then
$L^{\prime}=R^{\prime}\otimes_{R}L$ is an $F$-submodule of $R^{\prime}\otimes
S$ and $\left(  R^{\prime}\otimes_{k}S_{a}\right)  /L_{a}^{\prime}$ is locally
free of rank $h\left(  a\right)  $ for all $a\in A$, so define $\mathbb{H}%
_{\left(  S,F\right)  }^{h}\left(  \phi\right)  :\mathbb{H}_{\left(
S,F\right)  }^{h}\left(  R\right)  \rightarrow\mathbb{H}_{\left(  S,F\right)
}^{h}\left(  R^{\prime}\right)  $, $L\mapsto L^{\prime}$. These
\newsym[$\mathbb{H}_{\left(  S,F\right)  }^{h}$]{Hilbert functor}{}assignments
make $\mathbb{H}_{\left(  S,F\right)  }^{h}$ into a functor \underline
{$k-Alg$}$\rightarrow$\underline{$Set$}, the
\index{Hilbert functor|textbf}%
\textbf{Hilbert functor}.
\end{definition}

If $D\subset A$ is a subset the restriction $\left(  S_{D},F_{D}\right)  $ of
$\left(  S,F\right)  $ to degree $D$ is defined by%
\[%
\begin{tabular}
[c]{lll}%
$S_{D}=\bigoplus_{a\in D}S_{a}$ &  & $F_{D}=\bigcup_{a,b\in D}F_{a,b}$%
\end{tabular}
\]
It is a full subcategory of $\left(  S,F\right)  $ and there is the natural
restriction map%
\[%
\begin{tabular}
[c]{ccc}%
$\mathbb{H}_{\left(  S,F\right)  }^{h}$ & $\rightarrow$ & $\mathbb{H}_{\left(
S_{D},F_{D}\right)  }^{h}$\\
$L$ & $\mapsto$ & $L_{D}=\bigoplus_{a\in D}L_{a}$%
\end{tabular}
\]

\begin{lemma}
If $L^{\prime}\subset R\otimes S_{D}$ is an $F_{D}$-submodule and $L\subset
R\otimes S$ is the $F$-submodule generated by $L^{\prime}$, then $L_{b}%
=\sum_{a\in D}F_{a,b}\left(  L_{a}^{\prime}\right)  $ for all $b\in A$, so
$L_{D}=L^{\prime}$.
\end{lemma}

\begin{theorem}
\label{thm hilbert scheme for finite degrees}\cite{HS Multigraded Hilbert
schemes} Let $k$ be a commutative ring, $A$ a set, $S$ an $A$-graded
$k$-module with operators $F$ and $h:A\rightarrow\mathbb{N}$ a function. If
there are homogeneous $k$-submodules $M\subset N\subset S$ such that

\begin{enumerate}
\item $N$ is a finitely generated $k$-module,

\item $N$ generates $S$ as an $F$-module,

\item for all fields $K\in$\underline{$k-Alg$} and for all $L\in
\mathbb{H}_{\left(  S,F\right)  }^{h}\left(  K\right)  $ the submodule
$M\subset S$ spanns $\left(  K\otimes S\right)  /L$,

\item there is a subset $G\subset F$ which generates $F$ as a category such
that $GM\subset N$,
\end{enumerate}

then

\begin{itemize}
\item $N$ spans $\left(  K\otimes S\right)  /L$ so $\infty>\dim_{K}\left(
\left(  K\otimes S\right)  /L\right)  =\sum_{a\in A}h\left(  a\right)  $ hence
$h$ has finite support,

\item $\mathbb{H}_{\left(  S,F\right)  }^{h}$ is represented by a
quasiprojective closed subscheme of $\mathbb{G}_{N\backslash M}^{h}$ over $k$,
the
\index{Hilbert scheme}%
Hilbert scheme.
\end{itemize}
\end{theorem}

\begin{corollary}
\cite{HS Multigraded Hilbert schemes} If $A$ is finite and $S_{a}$ is a
finitely generated $k$-module for all $a\in A$, then in above theorem one can
choose $M=N=S$ and $G=F$, so $\mathbb{H}_{\left(  S,F\right)  }^{h}$ is
represented by a closed subscheme of the projective Grassmann scheme
$\mathbb{G}_{N\backslash M}^{h}=\mathbb{G}_{N}^{h}$, hence is projective.
\end{corollary}

\begin{theorem}
\label{thm Hilbert scheme reduction to finite degree}\cite{HS Multigraded
Hilbert schemes} Let $k$ be a commutative ring, $A$ a set, $S$ an $A$-graded
$k$-module with operators $F$ and $h:A\rightarrow\mathbb{N}$ a function.
Suppose $D\subset A$ such that

\begin{enumerate}
\item $\mathbb{H}_{\left(  S_{D},F_{D}\right)  }^{h}$ is represented by a
scheme over $k$,

\item for all $a\in A$ there is a finite set of operators $E\subset
\bigcup_{b\in D}F_{b,a}$ such that $S_{a}/\sum_{b\in D}E_{b,a}\left(
S_{b}\right)  $ is a finitely generated $k$-module,

\item for all fields $K\in$\underline{$k-Alg$} and for all $L^{\prime}%
\in\mathbb{H}_{\left(  S_{D},F_{D}\right)  }^{h}\left(  K\right)  $%
\[
\dim\left(  \left(  K\otimes S_{a}\right)  /L_{a}\right)  \leq h\left(
a\right)
\]
for all $a\in A$, where $L\subset K\otimes S$ is the $F$-submodule generated
by $L^{\prime}$.
\end{enumerate}

Then $\mathbb{H}_{\left(  S,F\right)  }^{h}$ is a subfunctor of $\mathbb{H}%
_{\left(  S_{D},F_{D}\right)  }^{h}$ via the natural restriction map $L\mapsto
L_{D}$ and is represented by a closed subscheme of the Hilbert scheme
representing $\mathbb{H}_{\left(  S_{D},F_{D}\right)  }^{h}$.

If $D$ is finite, then $\mathbb{H}_{\left(  S_{D},F_{D}\right)  }^{h}$ is
projective, hence $\mathbb{H}_{\left(  S,F\right)  }^{h}$ is projective.
\end{theorem}

\subsubsection{Example: Multigraded Hilbert schemes of admissible
ideals\label{Sec multigraded hilbert scheme of admissable ideals}}

Let $k$ be a commutative ring, $A$ an abelian group and $S=k\left[
x_{1},...,x_{r}\right]  $ a polynomial ring graded by a homomorphism of
semigroups $\deg:\mathbb{N}^{n}\rightarrow A$ via $\deg x^{u}=\deg u$, so%
\[
S=\bigoplus_{a\in A}S_{a}%
\]
and
\[
S_{a}\cdot S_{b}\subset S_{a+b}%
\]

The ring $S$ comes with operators $F=\bigcup_{a,b\in A}F_{a,b}$ where%
\[
F_{a,b}=\left\{  \left\{
\begin{tabular}
[c]{lll}%
$S_{a}$ & $\rightarrow$ & $S_{b}$\\
$f$ & $\mapsto$ & $m\cdot f$%
\end{tabular}
\ \right\}  \in\operatorname*{Hom}\nolimits_{k}\left(  S_{a},S_{b}\right)
\mid%
\begin{tabular}
[c]{l}%
$m\in S\text{ a monomial with}$\\
$\deg m=b-a$%
\end{tabular}
\right\}
\]

If $L\subset R\otimes S=\bigoplus_{a\in A}R\otimes S_{a}$ is an $F$-submodule,
then $L$ is a homogeneous ideal with respect to the grading of $R\otimes S$ by
$A$.

A homogeneous ideal $I\subset S$ is called
\index{admissable|textbf}%
\textbf{admissible} if $\left(  S/I\right)  _{a}=S_{a}/I_{a}$ is a locally
free $k$-module of finite rank for all $a\in A$. Denote by
\[%
\begin{tabular}
[c]{llll}%
$h_{S/I}:$ & $A$ & $\rightarrow$ & $\mathbb{N}$\\
& $a$ & $\mapsto$ & $\operatorname*{rank}_{k}\left(  \left(  S/I\right)
_{a}\right)  $%
\end{tabular}
\]
the Hilbert function of $S/I$. Denote by $A_{+}=\left\langle a_{1}%
,...,a_{r}\right\rangle \subset A$ the subgroup generated by $a_{i}=\deg
x_{i}$. The support of $h_{S/I}$ is contained in $A_{+}$.

If $h:A\rightarrow\mathbb{N}$ is a function with support on $A_{+}$, then for
any $R\in$\underline{$k-Alg$}%
\begin{align*}
\mathbb{H}_{\left(  S,F\right)  }^{h}\left(  R\right)   &  =\left\{  I\mid%
\begin{tabular}
[c]{l}%
$I\subset R\otimes S$ homogeneous ideal such that\\
$\left(  R\otimes S_{a}\right)  /I_{a}\text{ locally free of rank }h\left(
a\right)  $ $\forall a\in A$%
\end{tabular}
\right\} \\
&  =\left\{  I\mid I\subset R\otimes S\text{ admissible ideal with }%
h_{S/I}=h\right\}
\end{align*}
consists of the admissible ideals in $R\otimes S$ with Hilbert function $h$.

An antichain of monomial ideals in $S$ is a set $C$ of monomial ideals such
that for all $I_{1},I_{2}\in C$ it holds $I_{1}\not \subset I_{2}$.

\begin{lemma}
\cite{HS Multigraded Hilbert schemes} If $C$ is an antichain in $S$ then $C$
is finite.
\end{lemma}

So if $C$ is the set of all monomial ideals in $S$ with Hilbert function $h$
then $C$ is finite.

\begin{corollary}
\label{Cor existence finite set degrees}\cite{HS Multigraded Hilbert schemes}
If $h:A\rightarrow\mathbb{N}$ is a function with support on $A_{+}$, then
there is a finite set $D\subset A$ such that

\begin{enumerate}
\item any monomial ideal $I\subset S$ with $h_{S/I}=h$ is generated by
monomials in degrees $D$,

\item any monomial ideal $I\subset S$ generated in degrees $D$ satisfies: If
$h_{S/I}\left(  a\right)  =h\left(  a\right)  $ for all $a\in D$, then
$h_{S/I}\left(  a\right)  \leq h\left(  a\right)  $ for all $a\in A$.
\end{enumerate}
\end{corollary}

For any finite $D\subset A$ the assumptions of Theorem
\ref{thm hilbert scheme for finite degrees} hold for $\left(  S_{D}%
,F_{D}\right)  $, hence $\mathbb{H}_{\left(  S_{D},F_{D}\right)  }^{h}$ is
represented by a quasiprojective scheme.

For $D$ as given by Corollary \ref{Cor existence finite set degrees} the
assumptions of Theorem \ref{thm Hilbert scheme reduction to finite degree} are
satisfied, hence:

\begin{theorem}
\cite{HS Multigraded Hilbert schemes} If $h:A\rightarrow\mathbb{N}$ is a
function with support on $A_{+}$, then $\mathbb{H}_{\left(  S,F\right)  }^{h}$
is represented by a quasiprojective scheme.
\end{theorem}

Note that this setup is not directly applicable to the Cox ring of a toric
variety or the homogeneous coordinate ring of projective space.

\subsubsection{Example: Classical Hilbert
functor\label{Sec classical hilbert functor}}

The Grothendieck
\index{Hilbert scheme}%
Hilbert scheme
\newsym[$\mathbb{H}_{n}^{P}$]{classical Hilbert functor}{}represents the
functor $\mathbb{H}_{n}^{P}$ with%
\[
\mathbb{H}_{n}^{P}\left(  R\right)  =\left\{  X\mid X\subset\mathbb{P}%
^{n}\left(  R\right)  \text{ flat family with Hilbert polynomial }P\right\}
\]
for $R\in$\underline{$k-Alg$}. These $X$ correspond to saturated homogeneous
ideals $I\subset R\left[  x_{0},...,x_{n}\right]  $ with Hilbert polynomial
$P$.

Given $P$ and $n$ there is a degree $d_{0}$, the maximum of the
Castelnuovo-Mumford regularities of all saturated monomial ideals in $R\left[
x_{0},...,x_{n}\right]  $ with Hilbert polynomial $P$ such that for all
saturated homogeneous ideals $I\subset R\left[  x_{0},...,x_{n}\right]  $ with
Hilbert polynomial $P$%
\[
h_{S/I}\left(  a\right)  =P\left(  a\right)  \text{ for all }a\geq d_{0}%
\]

\begin{proposition}
\cite{HS Multigraded Hilbert schemes} Consider $S=k\left[  x_{0}%
,...,x_{n}\right]  $, let $F$ be the multiplication by monomials, $P$ some
Hilbert polynomial and%
\[
h\left(  a\right)  =\left\{
\begin{tabular}
[c]{ll}%
$\binom{n+a-1}{a}$ & for $a<d_{0}$\\
$P\left(  a\right)  $ & for $a\geq d_{0}$%
\end{tabular}
\right\}
\]
The Grothendieck
\index{Hilbert scheme}%
Hilbert scheme representing $\mathbb{H}_{n}^{P}$ is isomorphic to the Hilbert
scheme representing $\mathbb{H}_{\left(  S,F\right)  }^{h}$ via the bijection%
\[%
\begin{tabular}
[c]{ccc}%
$\mathbb{H}_{n}^{P}\left(  R\right)  $ & $\rightleftarrows$ & $\mathbb{H}%
_{\left(  S,F\right)  }^{h}\left(  R\right)  $\\
$I_{\geq a_{0}}$ & $\leftarrow$ & $I$\\
$J$ & $\mapsto$ & $\left(  J:\left\langle x_{0},...,x_{n}\right\rangle
^{\infty}\right)  $%
\end{tabular}
\]

\end{proposition}

\subsubsection{Tangent space and
deformations\label{Sec tangent space and deformations}}

Let $k$ be a field, $A$ an abelian group, $S=k\left[  x_{1},...,x_{r}\right]
$ graded by $\deg:\mathbb{N}^{n}\rightarrow A$ and $F$ the multiplication by
monomials. Let $h:A\rightarrow\mathbb{N}$ be a function with support on
$A_{+}$ and let $I\in\mathbb{H}_{\left(  S,F\right)  }^{h}\left(  k\right)  $.
The $S$-module $\operatorname*{Hom}_{S}\left(  I,S/I\right)  $ is graded by
$A$ and $\operatorname*{Hom}_{S}\left(  I,S/I\right)  _{a}$ is a finite
dimensional $k$-vector space for all $a\in A$.

Let $R=k\left[  t\right]  /\left\langle t^{2}\right\rangle $ and
$\phi:R\rightarrow k$, $s\mapsto s/\left\langle t\right\rangle $ so the map%
\[
\mathbb{H}_{\left(  S,F\right)  }^{h}\left(  \phi\right)  :\mathbb{H}_{\left(
S,F\right)  }^{h}\left(  R\right)  \rightarrow\mathbb{H}_{\left(  S,F\right)
}^{h}\left(  k\right)
\]
is given by $J\mapsto J/\left\langle t\right\rangle $. The Zariski tangent
space of the scheme representing $\mathbb{H}_{\left(  S,F\right)  }^{h}$ at
$I\in\mathbb{H}_{\left(  S,F\right)  }^{h}\left(  k\right)  $ is%
\begin{align*}
&  \left\{  J\in\mathbb{H}_{\left(  S,F\right)  }^{h}\left(  R\right)
\mid\mathbb{H}_{\left(  S,F\right)  }^{h}\left(  \phi\right)  \left(
J\right)  =I\right\} \\
&  =\left\{  J\mid%
\begin{tabular}
[c]{l}%
$J\subset R\otimes S$ an $A$-homogeneous ideal with $J/\left\langle
t\right\rangle =I$\\
such that $R\left[  x_{1},...,x_{r}\right]  /J$ is a free $R$-module
\end{tabular}
\ \right\}
\end{align*}
and is isomorphic to $\operatorname*{Hom}_{S}\left(  I,S/I\right)  _{0}$ by
associating to $J$ the homomorphism%
\[
S\overset{t\cdot}{\rightarrow}t\cdot R\left[  x_{1},...,x_{r}\right]
\rightarrow t\cdot R\left[  x_{1},...,x_{r}\right]  /\left(  J\cap t\cdot
R\left[  x_{1},...,x_{r}\right]  \right)  \cong S/I
\]

\begin{proposition}
\cite{HS Multigraded Hilbert schemes} The Zariski tangent space of the scheme
representing $\mathbb{H}_{\left(  S,F\right)  }^{h}$ at $I\in\mathbb{H}%
_{\left(  S,F\right)  }^{h}\left(  k\right)  $ is canonically isomorphic to
$\operatorname*{Hom}_{S}\left(  I,S/I\right)  _{0}$.
\end{proposition}

\subsubsection{Stanley decompositions\label{Sec stanley decompositions}}

\paragraph{Setup}

Let $Y=X\left(  \Sigma\right)  $ be a smooth complete toric variety of
dimension $n$ given by the fan $\Sigma\subset N_{\mathbb{R}}=N\otimes
_{\mathbb{Z}}\mathbb{R}$ with $N\cong\mathbb{Z}^{n}$. Let $S=\mathbb{C}\left[
y_{r}\mid r\in\Sigma\left(  1\right)  \right]  $ be the Cox ring of $Y$ graded
by $A_{n-1}\left(  Y\right)  $, and
\begin{gather*}
B\left(  \Sigma\right)  =\left\langle y^{D_{\widehat{\sigma}}}\mid\sigma
\in\Sigma\right\rangle \subset S\\
\text{ with }D_{\widehat{\sigma}}=\sum_{r\in\Sigma\left(  1\right)  ,\text{
}r\not \subset \sigma}D_{r}%
\end{gather*}
the
\index{irrelevant ideal}%
irrelevant ideal of $Y$. Write
\[
0\rightarrow M\overset{A}{\rightarrow}\operatorname*{WDiv}\nolimits_{T}\left(
Y\right)  \overset{\deg}{\rightarrow}A_{n-1}\left(  Y\right)  \rightarrow0
\]
for the presentation of the
\index{Chow group}%
Chow group of $Y$ and set $a_{i}=\deg D_{i}$. Denote by%
\[
Y^{\prime}=\mathbb{A}^{\Sigma\left(  1\right)  }-V\left(  B\left(
\Sigma\right)  \right)  =X\left(  \Sigma^{\prime}\right)  \rightarrow Y
\]
the Cox quotient presentation of $Y$ as defined in Section
\ref{1homogeneouscoordinate} and set $Y^{\prime\prime}=\mathbb{A}%
^{\Sigma\left(  1\right)  }=X\left(  \Sigma^{\prime\prime}\right)  $ with the
fan $\Sigma^{\prime\prime}\subset\mathbb{Z}^{\Sigma\left(  1\right)  }$ over
the standard simplex. For
\[
D\in\operatorname*{WDiv}\nolimits_{T}\left(  Y^{\prime\prime}\right)
\cong\operatorname*{WDiv}\nolimits_{T}\left(  Y^{\prime}\right)
\cong\operatorname*{WDiv}\nolimits_{T}\left(  Y\right)  \cong\mathbb{Z}%
^{\Sigma\left(  1\right)  }%
\]
denote by $x^{D}$ the corresponding (Laurent-) monomial in the Cox ring $S$.

Denote by $\mathcal{K}$ the set of integral points in the closure of the
\index{K\"{a}hler cone}%
K\"{a}hler cone%
\[
\operatorname*{cpl}\left(  \Sigma\right)  \subset A_{n-1}^{+}\left(  Y\right)
\otimes\mathbb{R\subset}A_{n-1}\left(  Y\right)  \otimes\mathbb{R}\cong
H^{2}\left(  Y,\mathbb{R}\right)
\]
as described in Section \ref{Sec Kaehler cone Mori cone}.

\paragraph{Primary decompositions and Stanley decompositions of monomial
ideals}

Consider first the vanishing locus of a monomial ideal in the affine space
$Y^{\prime\prime}$.

\begin{definition}
If $I\subset S$ is a \newsym[$\mathcal{S}$]{Stanley decomposition}{}monomial
ideal, then a
\index{Stanley decomposition|textbf}%
\textbf{Stanley decomposition} of $I$ is a subset%
\[
\mathcal{S}\subset\left\{  \left(  D,\sigma\right)  \mid D\in
\operatorname*{WDiv}\nolimits_{T}\left(  Y^{\prime\prime}\right)  \text{,
}D\text{ effective, }\sigma\in\Sigma^{\prime\prime}\right\}
\]
such that%
\[
S/I\cong%
{\displaystyle\bigoplus\limits_{\left(  D,\sigma\right)  \in\mathcal{S}}}
S_{\sigma}\left(  -\left[  D\right]  \right)
\]
where $S_{\sigma}=\mathbb{C}\left[  y_{r}\mid r\notin\sigma\right]  \cong
S/I\left(  V_{Y^{\prime\prime}}\left(  \sigma\right)  \right)  $ is the
\newsym[$S_{\sigma}$]{Cox ring with variables not in $\sigma$}{}Cox ring of
$U_{Y^{\prime\prime}}\left(  \sigma\right)  $. Here $V_{Y^{\prime\prime}%
}\left(  \sigma\right)  \subset Y^{\prime\prime}$ is the torus orbit closure
associated to $\sigma\in\Sigma^{\prime\prime}$ and $U_{Y^{\prime\prime}%
}\left(  \sigma\right)  =\operatorname*{Spec}\left(  \mathbb{C}\left[
\check{\sigma}\cap M\right]  \right)  \subset Y^{\prime\prime}$.
\end{definition}

\begin{remark}
Note that for the Cox quotient representation $Y^{\prime\prime}\supset
Y^{\prime}\rightarrow Y$ it holds%
\[
\operatorname*{WDiv}\nolimits_{T}\left(  Y^{\prime\prime}\right)
\cong\operatorname*{WDiv}\nolimits_{T}\left(  Y^{\prime}\right)
\cong\operatorname*{WDiv}\nolimits_{T}\left(  Y\right)
\]
and%
\[
\Sigma^{\prime\prime}\supset\Sigma^{\prime}\supset\Sigma^{\prime}\left(
1\right)  \overset{1:1}{\rightleftarrows}\Sigma\left(  1\right)
\]
so $\Sigma$ can be considered as a subfan of $\Sigma^{\prime\prime}$. If
$\sigma\in\Sigma$, then
\[%
\begin{tabular}
[c]{lll}%
$Y^{\prime\prime}$ & $\supset$ & $V_{Y^{\prime\prime}}\left(  \sigma\right)
=\left\{  y\in Y^{\prime\prime}\mid y_{r}=0\text{ }\forall r\in\Sigma\left(
1\right)  \text{, }r\subset\sigma\right\}  $\\
$\cup$ &  & $\cup$\\
$Y^{\prime}=Y^{\prime\prime}-V\left(  B\left(  \Sigma\right)  \right)  $ &
$\supset$ & $V_{Y^{\prime}}\left(  \sigma\right)  =\left\{  y\in Y^{\prime
}\mid y_{r}=0\text{ }\forall r\in\Sigma\left(  1\right)  \text{, }%
r\subset\sigma\right\}  $\\
$\downarrow$ &  & $\downarrow$\\
$Y$ & $\supset$ & $V_{Y}\left(  \sigma\right)  $%
\end{tabular}
\]
so the prime ideal
\[
\left\langle y_{r}\mid r\in\Sigma\left(  1\right)  \text{, }r\subset
\sigma\right\rangle \subset S
\]
corresponds to the torus orbit closure $V_{Y}\left(  \sigma\right)  \subset Y$.

Recall also that with%
\begin{align*}
D_{\widehat{\sigma}}  &  =\sum_{r\in\Sigma\left(  1\right)  ,\text{
}r\not \subset \sigma}D_{r}\\
U_{Y^{\prime}}\left(  \sigma\right)   &  =Y^{\prime\prime}-V\left(
y^{D_{\widehat{\sigma}}}\right)
\end{align*}
we have%
\[
U_{Y^{\prime}}\left(  \sigma\right)  /G\left(  \Sigma\right)  =U_{Y}\left(
\sigma\right)
\]

\end{remark}

Any associated prime of $I$ is of the form $\left\langle y_{r}\mid r\in
\Sigma\left(  1\right)  \text{, }r\subset\sigma\right\rangle $ for some
$\sigma\in\Sigma^{\prime\prime}$.

\begin{lemma}
\cite{MaSm Uniform bounds on multigraded regularity} Let $I\subset S$ be a
monomial ideal. Then $I$ is $B\left(  \Sigma\right)  $-saturated if and only
if all associated primes of $I$ are of the form
\[
\left\langle y_{r}\mid r\in\Sigma\left(  1\right)  \text{, }r\subset
\sigma\right\rangle
\]
for $\sigma\in\Sigma$.
\end{lemma}

A pair $\left(  D,\sigma\right)  $ with $D\in\operatorname*{WDiv}%
\nolimits_{T}\left(  Y^{\prime\prime}\right)  $, $D$ effective and $\sigma
\in\Sigma^{\prime\prime}$ is called
\index{admissable pair|textbf}%
\textbf{admissible} if $\operatorname*{supp}\left(  D\right)  \cap
\operatorname*{supp}\left(  D_{\widehat{\sigma}}\right)  =\varnothing$, i.e.,
if $D\subset U_{Y^{\prime\prime}}\left(  \sigma\right)  $.

A partial order on the set of admissible pairs is given by%
\[%
\begin{tabular}
[c]{lll}%
$\left(  D_{1},\sigma_{1}\right)  \leq\left(  D_{2},\sigma_{2}\right)  $ &
$\Leftrightarrow$ &
\begin{tabular}
[c]{l}%
$D_{2}-D_{1}\geq0\text{ and }$\\
$\operatorname*{supp}\left(  \left(  D_{2}-D_{1}\right)  +D_{\widehat
{\sigma_{2}}}\right)  \subset\operatorname*{supp}\left(  D_{\widehat
{\sigma_{1}}}\right)  $%
\end{tabular}
\\
& $\Leftrightarrow$ &
\begin{tabular}
[c]{l}%
$D_{2}-D_{1}\geq0\text{ and }$\\
$U_{Y^{\prime\prime}}\left(  \sigma_{1}\right)  \subset U_{Y^{\prime\prime}%
}\left(  \sigma_{2}\right)  \cap\left(  Y^{\prime\prime}-\operatorname*{supp}%
\left(  D_{2}-D_{1}\right)  \right)  $%
\end{tabular}
\\
& $\Leftrightarrow$ & $y^{D_{2}}S_{\sigma_{2}}\subset y^{D_{1}}S_{\sigma_{1}}$%
\end{tabular}
\]

An admissible pair $\left(  D,\sigma\right)  $ is called
\index{standard|textbf}%
\textbf{standard} with respect to $I$ if $\left(  D,\sigma\right)  $ is
minimal with respect to $\leq$ with the property $y^{D}S_{\sigma}\cap
I=\left\{  0\right\}  $.

\begin{lemma}
If $\mathcal{S}$ gives a Stanley decomposition of the monomial ideal $I\subset
S$, then%
\[
I=%
{\displaystyle\bigcap\limits_{\left(  D,\sigma\right)  \in\mathcal{S}}}
\left\langle y_{r}^{u_{r}+1}\mid r\in\Sigma\left(  1\right)  \text{, }%
r\subset\sigma\text{, }D=%
{\textstyle\sum\nolimits_{r\in\Sigma\left(  1\right)  }}
u_{r}D_{r}\right\rangle
\]
If $\left(  D,\sigma\right)  \in\mathcal{S}$ is a standard pair, then
$\left\langle y_{r}\mid r\in\Sigma\left(  1\right)  \text{, }r\subset
\sigma\right\rangle =I\left(  V_{Y^{\prime\prime}}\left(  \sigma\right)
\right)  $ is an associated prime of $I$.
\end{lemma}

\begin{algorithm}
\label{alg stanley decomposition}\cite{MaSm Uniform bounds on multigraded
regularity} The following algorithm computes a Stanley decomposition of the
monomial ideal $I\subset S$:

\begin{itemize}
\item If $I$ is a prime ideal and $I=\left\langle y_{r}\mid r\in\Sigma\left(
1\right)  \text{, }r\subset\sigma\right\rangle $ with $\sigma\in\Sigma
^{\prime\prime}$, then return $\left\{  \left(  0,\sigma\right)  \right\}  $.

\item Otherwise, let $r\in\Sigma\left(  1\right)  $ such that $I\neq\left(
I:\left\langle y_{r}\right\rangle \right)  \neq\left\langle 1\right\rangle $ .

Compute Stanley decompositions $\mathcal{S}_{1}$ of $S/\left(  I+\left\langle
y_{r}\right\rangle \right)  $ and $\mathcal{S}_{2}$ of $S/\left(
I:\left\langle y_{r}\right\rangle \right)  $.

Return
\[
\mathcal{S}=\left\{  \left(  D_{1},\sigma_{1}\right)  \mid\left(  D_{1}%
,\sigma_{1}\right)  \in\mathcal{S}_{1}\right\}  \cup\left\{  \left(
D_{2}+D_{r},\sigma_{2}\right)  \mid\left(  D_{2},\sigma_{2}\right)
\in\mathcal{S}_{2}\right\}
\]

\end{itemize}
\end{algorithm}

\begin{example}
A Stanley decomposition of the reduced ideal
\[
I=\left\langle y_{1}y_{2}y_{3}\right\rangle \subset S=\mathbb{C}\left[
y_{0},y_{1},y_{2}\right]
\]
is given by
\begin{align*}
S/I  &  =1\cdot\mathbb{C\oplus}\\
&  y_{0}\cdot\mathbb{C}\left[  y_{0}\right]  \mathbb{\oplus}y_{1}%
\cdot\mathbb{C}\left[  y_{1}\right]  \mathbb{\oplus}y_{2}\cdot\mathbb{C}%
\left[  y_{2}\right]  \mathbb{\oplus}\\
&  y_{0}y_{1}\cdot\mathbb{C}\left[  y_{0},y_{1}\right]  \mathbb{\oplus}%
y_{1}y_{2}\cdot\mathbb{C}\left[  y_{1},y_{2}\right]  \mathbb{\oplus}y_{0}%
y_{2}\cdot\mathbb{C}\left[  y_{0},y_{2}\right]
\end{align*}
writing $I$ as the intersection%
\[
I=\underset{\text{irrelevant}}{\underbrace{\left\langle y_{0},y_{1}%
,y_{2}\right\rangle }}\cap\left\langle y_{1},y_{2}\right\rangle \cap
\left\langle y_{0},y_{2}\right\rangle \cap\left\langle y_{0},y_{1}%
\right\rangle \cap\underset{\text{associated primes}}{\underbrace{\left\langle
y_{2}\right\rangle \cap\left\langle y_{0}\right\rangle \cap\left\langle
y_{1}\right\rangle }}%
\]
corresponding to the toric stratification of $I$. The ideal $\left\langle
y_{0},y_{1},y_{2}\right\rangle $ defines $0\in Y^{\prime\prime}$ but does not
correspond to a subvariety of $Y^{\prime}$ and $Y$, hence it is irrelevant.
The corresponding maximal cone $\sigma$ of the fan $\Sigma^{\prime\prime}$ is
not a cone of $\Sigma^{\prime}$.

The Stanley decomposition given by above algorithm is%
\begin{align*}
S/I  &  =1\cdot\mathbb{C}\left[  y_{1},y_{2}\right]  \mathbb{\oplus}y_{0}%
\cdot\mathbb{C}\left[  y_{0},y_{2}\right]  \mathbb{\oplus}y_{0}y_{1}%
\cdot\mathbb{C}\left[  y_{0},y_{1}\right] \\
I  &  =\left\langle y_{0}\right\rangle \cap\left\langle y_{1}\right\rangle
\cap\left\langle y_{2}\right\rangle
\end{align*}

\end{example}

Note that there are ideals which do not admit a Stanley decomposition, where
every factor corresponds to an associated prime, e.g.,
\[
I=\left\langle y_{1},y_{2}\right\rangle \cap\left\langle y_{0},y_{3}%
\right\rangle \subset\mathbb{C}\left[  y_{0},...,y_{3}\right]
\]

\begin{example}
Consider the ideal $I=\left\langle y_{1}y_{2},y_{0}y_{3}\right\rangle \subset
S=\mathbb{C}\left[  y_{0},...,y_{3}\right]  $. Then
\[%
\begin{tabular}
[c]{ll}%
$S/I=$ & $1\cdot\mathbb{C\oplus}$\\
& $y_{0}\cdot\mathbb{C}\left[  y_{0}\right]  \mathbb{\oplus}y_{1}%
\cdot\mathbb{C}\left[  y_{1}\right]  \mathbb{\oplus}y_{2}\cdot\mathbb{C}%
\left[  y_{2}\right]  \mathbb{\oplus}y_{3}\cdot\mathbb{C}\left[  y_{3}\right]
\mathbb{\oplus}$\\
& $y_{0}y_{1}\cdot\mathbb{C}\left[  y_{0},y_{1}\right]  \mathbb{\oplus}%
y_{0}y_{2}\cdot\mathbb{C}\left[  y_{0},y_{2}\right]  \mathbb{\oplus}y_{1}%
y_{3}\cdot\mathbb{C}\left[  y_{1},y_{3}\right]  \mathbb{\oplus}y_{2}y_{3}%
\cdot\mathbb{C}\left[  y_{2},y_{3}\right]  $%
\end{tabular}
\]
is a Stanley decomposition of $I$ representing the ideal as%
\begin{align*}
I  &  =\overset{\text{irrelevant}}{\overbrace{\left\langle y_{0},y_{1}%
,y_{2},y_{3}\right\rangle }}\cap\\
&  \left\langle y_{1},y_{2},y_{3}\right\rangle \cap\left\langle y_{0}%
,y_{2},y_{3}\right\rangle \cap\left\langle y_{0},y_{1},y_{3}\right\rangle
\cap\left\langle y_{0},y_{1},y_{2}\right\rangle \cap\\
&  \underset{\text{associated primes}}{\underbrace{\left\langle y_{0}%
,y_{1}\right\rangle \cap\left\langle y_{0},y_{2}\right\rangle \cap\left\langle
y_{1},y_{3}\right\rangle \cap\left\langle y_{2},y_{3}\right\rangle }}%
\end{align*}
The Stanley decomposition given by above algorithm is%
\begin{align*}
S/I  &  =1\cdot\mathbb{C}\left[  y_{2},y_{3}\right]  \mathbb{\oplus}y_{1}%
\cdot\mathbb{C}\left[  y_{1},y_{3}\right]  \mathbb{\oplus}y_{0}\cdot
\mathbb{C}\left[  y_{0},y_{2}\right]  \mathbb{\oplus}y_{0}y_{1}\cdot
\mathbb{C}\left[  y_{0},y_{1}\right] \\
I  &  =\left\langle y_{0},y_{1}\right\rangle \cap\left\langle y_{0}%
,y_{2}\right\rangle \cap\left\langle y_{1},y_{3}\right\rangle \cap\left\langle
y_{2},y_{3}\right\rangle
\end{align*}
The first decomposition is obtained from the second by further subdivision:
\begin{align*}
1\cdot\mathbb{C}\left[  y_{2},y_{3}\right]   &  =1\cdot\mathbb{C\oplus}%
y_{2}\cdot\mathbb{C}\left[  y_{2}\right]  \mathbb{\oplus}y_{3}\cdot
\mathbb{C}\left[  y_{3}\right]  \mathbb{\oplus}y_{2}y_{3}\cdot\mathbb{C}%
\left[  y_{2},y_{3}\right] \\
y_{0}\cdot\mathbb{C}\left[  y_{0},y_{2}\right]   &  =y_{0}y_{2}\cdot
\mathbb{C}\left[  y_{0},y_{2}\right]  \mathbb{\oplus}y_{0}\cdot\mathbb{C}%
\left[  y_{0}\right] \\
y_{1}\cdot\mathbb{C}\left[  y_{1},y_{3}\right]   &  =y_{1}y_{3}\cdot
\mathbb{C}\left[  y_{1},y_{3}\right]  \mathbb{\oplus}y_{1}\cdot\mathbb{C}%
\left[  y_{1}\right]
\end{align*}

\end{example}

\paragraph{Stanley filtrations}

\begin{definition}
If $I\subset S$ is a monomial ideal, then a
\index{Stanley filtration|textbf}%
\textbf{Stanley filtration} is a Stanley decomposition with ordering of the
elements
\[
\mathcal{S}=\left\{  \left(  D_{1},\sigma_{1}\right)  ,...,\left(
D_{s},\sigma_{s}\right)  \right\}
\]
such that for all $j=1,...,s$%
\[
\mathcal{S}_{j}=\left\{  \left(  D_{1},\sigma_{1}\right)  ,...,\left(
D_{j},\sigma_{j}\right)  \right\}
\]
is a Stanley decomposition of
\[
S/\left(  I+\left\langle x^{D_{j+1}},...,x^{D_{s}}\right\rangle \right)
\]

\end{definition}

So a Stanley filtration gives Stanley decompositions%
\[%
\begin{tabular}
[c]{lll}%
$S/\left(  I+\left\langle x^{D_{2}},...,x^{D_{s}}\right\rangle \right)  $ &
$\cong$ & $S_{\sigma_{1}}\left(  -\left[  D_{1}\right]  \right)  $\\
$S/\left(  I+\left\langle x^{D_{3}},...,x^{D_{s}}\right\rangle \right)  $ &
$\cong$ & $S_{\sigma_{1}}\left(  -\left[  D_{1}\right]  \right)  \oplus
S_{\sigma_{2}}\left(  -\left[  D_{2}\right]  \right)  $\\
$\vdots$ &  & $\vdots$\\
$S/I$ & $\cong$ & $S_{\sigma_{1}}\left(  -\left[  D_{1}\right]  \right)
\oplus...\oplus S_{\sigma_{s}}\left(  -\left[  D_{s}\right]  \right)  $%
\end{tabular}
\]
Algorithm \ref{alg stanley decomposition} computes a Stanley filtration by
ordering the leaves of the generated tree by listing the $\left(
I+\left\langle y_{r}\right\rangle \right)  $ child prior to the $\left(
I:\left\langle y_{r}\right\rangle \right)  $ child. This algorithm is
implemented in the Macaulay2 library
\index{stanleyfiltration.m2|textbf}%
\textrm{stanleyfiltration.m2}.

\begin{example}
For $I=\left\langle y_{1}y_{2},y_{0}y_{3}\right\rangle \subset S=\mathbb{C}%
\left[  y_{0},...,y_{3}\right]  $ above algorithm computes the Stanley
filtration%
\begin{align*}
S/\left\langle y_{1},y_{0}\right\rangle  &  =1\cdot\mathbb{C}\left[
y_{2},y_{3}\right] \\
S/\left\langle y_{1}y_{2},y_{0}\right\rangle  &  =1\cdot\mathbb{C}\left[
y_{2},y_{3}\right]  \mathbb{\oplus}y_{1}\cdot\mathbb{C}\left[  y_{1}%
,y_{3}\right] \\
S/\left\langle y_{1}y_{2},y_{0}y_{3},y_{0}y_{1}\right\rangle  &
=1\cdot\mathbb{C}\left[  y_{2},y_{3}\right]  \mathbb{\oplus}y_{1}%
\cdot\mathbb{C}\left[  y_{1},y_{3}\right]  \mathbb{\oplus}y_{0}\cdot
\mathbb{C}\left[  y_{0},y_{2}\right] \\
S/\left\langle y_{1}y_{2},y_{0}y_{3}\right\rangle  &  =1\cdot\mathbb{C}\left[
y_{2},y_{3}\right]  \mathbb{\oplus}y_{1}\cdot\mathbb{C}\left[  y_{1}%
,y_{3}\right]  \mathbb{\oplus}y_{0}\cdot\mathbb{C}\left[  y_{0},y_{2}\right]
\mathbb{\oplus}y_{0}y_{1}\cdot\mathbb{C}\left[  y_{0},y_{1}\right]
\end{align*}
corresponding to%
\begin{align*}
\left\langle y_{1},y_{0}\right\rangle  &  =\left\langle y_{0},y_{1}%
\right\rangle \\
\left\langle y_{1}y_{2},y_{0}\right\rangle  &  =\left\langle y_{0}%
,y_{1}\right\rangle \cap\left\langle y_{0},y_{2}\right\rangle \\
\left\langle y_{1}y_{2},y_{0}y_{3},y_{0}y_{1}\right\rangle  &  =\left\langle
y_{0},y_{1}\right\rangle \cap\left\langle y_{0},y_{2}\right\rangle
\cap\left\langle y_{1},y_{3}\right\rangle \\
\left\langle y_{1}y_{2},y_{0}y_{3}\right\rangle  &  =\left\langle y_{0}%
,y_{1}\right\rangle \cap\left\langle y_{0},y_{2}\right\rangle \cap\left\langle
y_{1},y_{3}\right\rangle \cap\left\langle y_{2},y_{3}\right\rangle
\end{align*}

\end{example}

\subsubsection{Multigraded regularity\label{Sec Multigraded regularity}}

Let $C=\left\{  c_{1},...,c_{e}\right\}  \subset A_{n-1}\left(  Y\right)  $ be
a finite subset and $\mathbb{N}C\subset A_{n-1}\left(  Y\right)  $ the
semigroup generated by $C$. A subset $D\subset A_{n-1}\left(  Y\right)  $ is
called an $\mathbb{N}C$-module if $d+c\in D$ for all $d\in D$ and
$c\in\mathbb{N}C$. If $D$ is an $\mathbb{N}C$-module and $i\in\mathbb{Z}$,
then%
\[
D\left[  i\right]  =%
{\displaystyle\bigcup\limits_{\substack{\lambda_{1}+...+\lambda_{e}=\left\vert
i\right\vert \\\lambda_{j}\in\mathbb{Z}_{\geq0}}}}
\left(  \operatorname{sign}\left(  i\right)  \cdot%
{\textstyle\sum\nolimits_{j=1}^{e}}
\lambda_{j}c_{j}+D\right)  \subset A_{n-1}\left(  Y\right)
\]
is an $\mathbb{N}C$-module. For $m\in A_{n-1}\left(  Y\right)  $ it holds
$\left(  m+D\right)  \left[  i\right]  =m+D\left[  i\right]  $ and $D\left[
i+1\right]  \subset D\left[  i\right]  $.

\begin{definition}
Let $M$ be a finitely generated $A_{n-1}\left(  Y\right)  $-graded $S$-module
and let $m\in A_{n-1}\left(  Y\right)  $. Then $M$ is called
\index{regular|textbf}%
$m$\textbf{-regular} with respect to $C$ if for all $i\geq1$ and for all%
\[
a\in m+\mathbb{N}C\left[  1-i\right]
\]
the local cohomology
\newsym[$H_{B\left(  \Sigma\right)  }^{i}\left(  M\right)  $]{local cohomology}{}satisfies%
\[
H_{B\left(  \Sigma\right)  }^{i}\left(  M\right)  _{a}=0
\]
The regularity of $M$ with respect to $C$ is the
\newsym[$\operatorname{reg}_{C}\left(  M\right)  $]{regularity with respect to $C$}{}subset%
\[
\operatorname{reg}_{C}\left(  M\right)  =\left\{  m\in A_{n-1}\left(
Y\right)  \mid M\text{ is }m\text{-regular with respect to }C\right\}
\]

\end{definition}

\begin{remark}
With $C=\left\{  c_{1},...,c_{e}\right\}  $ the module $M$ is $m$-regular if
and only if%
\[
H_{B\left(  \Sigma\right)  }^{i}\left(  M\right)  _{a}=0
\]
holds for all $i\geq1$ and all
\[
a\in%
{\displaystyle\bigcup\limits_{\substack{\lambda_{1}+...+\lambda_{e}%
=i-1\\\lambda_{j}\in\mathbb{Z}_{\geq0}}}}
\left(  m-%
{\textstyle\sum\nolimits_{j=1}^{e}}
\lambda_{j}c_{j}+\mathbb{N}C\right)
\]
and for $i=0$ and all%
\[
a\in%
{\displaystyle\bigcup\limits_{j=1}^{e}}
\left(  m+c_{j}+\mathbb{N}C\right)
\]

\end{remark}

\begin{definition}
Let $m\in A_{n-1}\left(  Y\right)  $. A finitely generated $A_{n-1}\left(
Y\right)  $-graded $S$-module $M$ is called
\index{regular}%
$m$\textbf{-regular}, if it is $m$-regular with
\newsym[$\operatorname{reg}\left(  M\right)  $]{regularity}{}respect to the
unique minimal Hilbert basis $C=\left\{  c_{1},...,c_{e}\right\}  $ of
$\mathcal{K}$, giving any element of $\mathcal{K}$ as a $\mathbb{Z}_{\geq0}%
$-linear combination, i.e., with $\mathbb{N}C=\mathcal{K}$. The
\index{regularity|textbf}%
\textbf{regularity} of $M$ is%
\[
\operatorname{reg}\left(  M\right)  =\left\{  m\in A_{n-1}\left(  Y\right)
\mid M\text{ is }m\text{-regular}\right\}
\]

\end{definition}

The local cohomology groups may be computed in the following way:

Let $\Gamma$ be a finite regular
\index{cell complex|textbf}%
cell complex. \newsym[$\varepsilon$]{incidence function}{}A function
$\varepsilon:\Gamma\times\Gamma\rightarrow\left\{  -1,0,1\right\}  $ is called
an
\index{incidence function|textbf}%
\textbf{incidence function} if

\begin{itemize}
\item $\varepsilon\left(  F,G\right)  \neq0\Leftrightarrow G$ is a face of $F$.

\item $\varepsilon\left(  F,\varnothing\right)  =1$ for all $0$-cells $F\in$
$\Gamma^{0}$.

\item If $G\in\Gamma^{i-2}$ is a face of $F\in\Gamma^{i}$, then%
\[
\varepsilon\left(  F,H_{1}\right)  \varepsilon\left(  H_{1},G\right)
+\varepsilon\left(  F,H_{2}\right)  \varepsilon\left(  H_{2},G\right)  =0
\]
for the unique two faces $H_{1},H_{2}\in\Gamma^{i-1}$ of $F$ such that $G$ is
a face of $H_{1}$ and $H_{2}$.
\end{itemize}

\begin{lemma}
\cite[Sec. 6.3]{BrHe CohenMacaulay Rings} If $\Gamma$ is a finite regular cell
complex, then there is an incidence function on $\Gamma$ determined by an
orientation of the cells.
\end{lemma}

Associated to a cell complex $\Gamma$ of dimension $n$ together with an
\newsym[$\widetilde{\mathcal{C}}\left(  \Gamma\right)  $]{augumented chain complex}{}incidence
function $\varepsilon$ there is the
\index{augumented oriented chain complex|textbf}%
\textbf{augumented oriented chain complex}
\[
\widetilde{\mathcal{C}}\left(  \Gamma\right)  :0\rightarrow\mathcal{C}%
_{n-1}\overset{\delta}{\rightarrow}\mathcal{C}_{n-2}\rightarrow...\rightarrow
\mathcal{C}_{0}\overset{\delta}{\rightarrow}\mathcal{C}_{-1}\rightarrow0
\]
with coefficients in $R$ where%
\begin{align*}
\mathcal{C}_{i}  &  =%
{\textstyle\bigoplus_{F\in\Gamma^{i}}}
R\cdot F\\
\delta\left(  F\right)   &  =\sum_{G\in\Gamma^{i-1}}\varepsilon\left(
F,G\right)  G
\end{align*}
for $F\in\Gamma^{i}$ and extended linearily. For different incidence functions
the complexes $\widetilde{\mathcal{C}}\left(  \Gamma\right)  $ are isomorphic.
Denote by%
\[
\widetilde{H}_{i}\left(  \Gamma\right)  =H_{i}\left(  \widetilde{\mathcal{C}%
}\left(  \Gamma\right)  \right)
\]

\begin{theorem}
\cite[Sec. 6.3]{BrHe CohenMacaulay Rings} Let $\Gamma$ be a cell complex and
denote by $\left\vert \Gamma\right\vert $ the underlying topological space.
Then%
\[
\widetilde{H}_{i}\left(  \Gamma\right)  =\widetilde{H}_{i}\left(  \left\vert
\Gamma\right\vert \right)
\]
is the
\index{reduced singular homology}%
reduced singular homology of $\left\vert \Gamma\right\vert $.
\end{theorem}

An algorithm computing $\widetilde{H}_{i}\left(  \left\vert \Gamma\right\vert
\right)  $ via the augumented oriented chain complex is implemented in the
Macaulay2 library
\index{homology.m2|textbf}%
\textsf{homology.m2}.

Consider the cell complex $\Gamma$ given by the intersection of the fan
$\Sigma$ with a sphere of dimension $n-1$ together with the sphere as cell of
dimension $n-1$ and $\varnothing$ as $-1$-cell. Let $\varepsilon$ be an
incidence function given by an orientation. By abuse of notation identify
cones of $\Sigma$ and cells of $\Gamma$. For $\sigma\in\Gamma$ denote by
\[
S_{\left(  \sigma\right)  }=S_{x^{D_{\hat{\sigma}}}}%
\]
the localization of $S$ in the
\newsym[$S_{\left(  \sigma\right)  }$]{localization at the cone $\sigma$}{}multiplicatively
closed set generated by $x^{D_{\hat{\sigma}}}$. This relates to the Cox
quotient representation of $Y=X\left(  \Sigma\right)  $ as follows. We have
$S_{\left(  \sigma\right)  }=\mathbb{C}\left[  \check{\sigma}^{\prime}%
\cap\mathbb{Z}^{\Sigma\left(  1\right)  }\right]  $ where $\sigma^{\prime}%
\in\Sigma^{\prime}$ denotes the cone corresponding to $\sigma$, hence
$U_{Y^{\prime}}\left(  \sigma^{\prime}\right)  =\operatorname*{Spec}S_{\left(
\sigma\right)  }$ and $U_{Y}\left(  \sigma\right)  =U_{Y^{\prime}}\left(
\sigma\right)  /G\left(  \Sigma\right)  $. For the maximal cell $D_{\hat
{\sigma}}=0$ and $S_{\left(  \sigma\right)  }=S$.

Associate to $\Gamma$ the canonical \v{C}ech complex%
\[
C^{\ast}:0\rightarrow C^{0}\overset{\partial}{\rightarrow}C^{1}\rightarrow
...\overset{\partial}{\rightarrow}C^{n}\rightarrow0
\]
with%
\[
C^{i}=%
{\textstyle\bigoplus_{\sigma\in\Gamma^{n-i}}}
S_{\left(  \sigma\right)  }%
\]
and boundary map $\partial:C^{i-1}\rightarrow C^{i}$ given by the components%
\[
\partial:S_{\left(  \sigma_{1}\right)  }\rightarrow S_{\left(  \sigma
_{2}\right)  }%
\]
defined as $\varepsilon\left(  \sigma_{1},\sigma_{2}\right)  $ times the
natural map $S_{\left(  \sigma_{1}\right)  }\rightarrow S_{\left(  \sigma
_{2}\right)  }$ if $\sigma_{2}$ is a face of $\sigma_{1}$, and the $0$-map otherwise.

\begin{theorem}
\cite[Sec. 6.3]{BrHe CohenMacaulay Rings} If $M$ is a finitely generated
$A_{n-1}\left(  Y\right)  $-graded $S$-module, then%
\[
H_{B\left(  \Sigma\right)  }^{i}\left(  M\right)  \cong H^{i}\left(
M\otimes_{S}C^{\ast}\right)
\]

\end{theorem}

Consider the natural $\mathbb{Z}^{\Sigma\left(  1\right)  }$-grading refining
the $A_{n-1}\left(  Y\right)  $-grading. Recall that $\Sigma^{\prime\prime}$
is the fan over the simplex on $\Sigma\left(  1\right)  $ with $X\left(
\Sigma^{\prime\prime}\right)  =\mathbb{C}^{\Sigma\left(  1\right)  }$ and that
$\Sigma$ may be considered as a subfan of $\Sigma^{\prime\prime}$.

If $w\in\mathbb{Z}^{\Sigma\left(  1\right)  }$ define%
\[
G_{w}=\operatorname*{hull}\left\{  r\in\Sigma^{\prime\prime}\left(  1\right)
\mathbb{\mid}w_{r}<0\right\}
\]
Then%
\[
\left(  S_{\left(  \tau\right)  }\right)  _{w}=\left\{
\begin{tabular}
[c]{ll}%
$\mathbb{C}$ & if $G_{w}\subset\hat{\tau}$\\
$0$ & otherwise
\end{tabular}
\right\}
\]
For $G\in\Sigma^{\prime\prime}$ define%
\[
\Gamma_{G}=\left\{  F\in\Gamma\mid F\subset G\right\}
\]
Then%
\[
H^{i}\left(  C_{w}^{\ast}\right)  =\widetilde{H}_{n-i}\left(  \Gamma
_{\widehat{G_{w}}}\right)
\]
If $G_{w}\ $is the maximal cone of $\Sigma^{\prime\prime}$, then
$\Gamma_{G_{w}}=\Gamma$ is a sphere, hence%
\[
H^{i}\left(  C_{w}^{\ast}\right)  =\widetilde{H}_{n-i}\left(  \varnothing
\right)  =\left\{
\begin{tabular}
[c]{ll}%
$\mathbb{C}$ & $i=n+1$\\
$0$ & otherwise
\end{tabular}
\ \right\}  =\widetilde{H}^{i-2}\left(  \Gamma\right)
\]
If $G_{w}$ is the zero cone of $\Sigma^{\prime\prime}$, then
\[
H^{i}\left(  C_{w}^{\ast}\right)  =\widetilde{H}_{n-i}\left(  \Gamma\right)
=0=\widetilde{H}^{i-2}\left(  \varnothing\right)  \text{ for }i\neq1
\]
If $G_{w}$ lies between zero and maximal cone by Alexander duality%
\[
H^{i}\left(  C_{w}^{\ast}\right)  =\widetilde{H}_{n-i}\left(  \Gamma
_{\widehat{G_{w}}}\right)  \cong\widetilde{H}^{i-2}\left(  \Gamma
\backslash\Gamma_{\widehat{G_{w}}}\right)  \cong\widetilde{H}^{i-2}\left(
\Gamma_{G_{w}}\right)
\]
hence:

\begin{proposition}
For all $w\in\mathbb{Z}^{\Sigma\left(  1\right)  }$ and $i\neq1$
\[
\left(  H_{B\left(  \Sigma\right)  }^{i}\left(  S\right)  \right)  _{w}%
\cong\widetilde{H}^{i-2}\left(  \Gamma_{G_{w}}\right)
\]
and%
\[
\left(  H_{B\left(  \Sigma\right)  }^{1}\left(  S\right)  \right)  _{w}=0
\]

\end{proposition}

This allows computation of the local cohomology groups and regularity.

\begin{figure}
[h]
\begin{center}
\includegraphics[
height=1.465in,
width=1.4494in
]%
{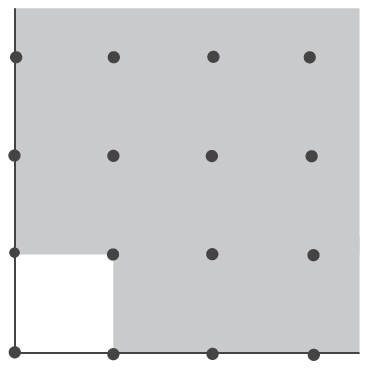}%
\caption{Regularity $\operatorname{reg}\left(  S\right)  $ for the Hirzebruch
surface $F_{2}$}%
\label{f2reg}%
\end{center}
\end{figure}

\begin{example}
Let
\[
Y=X\left(  \Sigma\right)  =F_{t}=\mathbb{P}\left(  \mathcal{O}_{\mathbb{P}%
^{1}}\oplus\mathcal{O}_{\mathbb{P}^{1}}\left(  t\right)  \right)
\]
be the \newsym[$F_{t}$]{Hirzebruch surface}{}Hirzebruch surface for $t\geq0$
given by the fan with the rays%
\[
\left(  1,0\right)  ,\left(  0,1\right)  ,\left(  -1,t\right)  ,\left(
0,-1\right)
\]
let%
\[
0\rightarrow M\overset{A}{\rightarrow}\mathbb{Z}^{4}\rightarrow A_{1}\left(
Y\right)  \rightarrow0
\]
with
\[
A=\left(
\begin{array}
[c]{cc}%
1 & 0\\
0 & 1\\
-1 & t\\
0 & -1
\end{array}
\right)
\]
be the presentation of the Chow group of $Y$ and, with respect to this
numbering of the rays, let%
\[
S=\mathbb{C}\left[  y_{1},y_{2},y_{3},y_{4}\right]
\]
be the Cox ring of $Y$, and
\[
B\left(  \Sigma\right)  =\left\langle y_{1},y_{3}\right\rangle \cap
\left\langle y_{2},y_{4}\right\rangle =\left\langle y_{1}y_{2},y_{2}%
y_{3},y_{3}y_{4}.y_{4}y_{1}\right\rangle
\]
the irrelevant ideal of $Y$. Fix an isomorphism%
\begin{gather*}
A_{1}\left(  Y\right)  \overset{B}{\rightarrow}\mathbb{Z}^{2}\\
B=\left(
\begin{array}
[c]{cccc}%
1 & -t & 1 & 0\\
0 & 1 & 0 & 1
\end{array}
\right)
\end{gather*}
Then%
\[
\operatorname{cpl}\left(  \Sigma\right)  =\operatorname*{hull}\left(  \left(
\begin{array}
[c]{c}%
1\\
0
\end{array}
\right)  ,\left(
\begin{array}
[c]{c}%
0\\
1
\end{array}
\right)  \right)  \subset\mathbb{Z}^{2}%
\]
and%
\[
\mathcal{K}=\operatorname{cpl}\left(  \Sigma\right)  \cap\mathbb{Z}%
^{2}=\mathbb{N}C
\]
with%
\[
C=\left\{  \left(
\begin{array}
[c]{c}%
1\\
0
\end{array}
\right)  ,\left(
\begin{array}
[c]{c}%
0\\
1
\end{array}
\right)  \right\}
\]

If $t=0,1$, then the regularity of $S$ is%
\[
\operatorname{reg}\left(  S\right)  =\mathcal{K}=\mathbb{Z}_{\geq0}^{2}%
\]
and for $t\geq2$%
\[
\operatorname{reg}\left(  S\right)  =\left(  \left(
\begin{array}
[c]{c}%
t-1\\
0
\end{array}
\right)  +\mathcal{K}\right)  \cup\left(  \left(
\begin{array}
[c]{c}%
0\\
1
\end{array}
\right)  +\mathcal{K}\right)
\]
shown in Figure \ref{f2reg} for $t=2$.
\end{example}

%

\begin{proposition}
\cite{MaSm Multigraded CastelnuovoMumford regularity} Let $M$ be a finitely
generated $A_{n-1}\left(  Y\right)  $-graded $S$-module $M$. Then
$\widetilde{M}$ is zero if and only if there is a $j>0$ such that
\[
\left(  B\left(  \Sigma\right)  ^{j}M\right)  _{a}=0\text{ }\forall
a\in\mathcal{K}%
\]

\end{proposition}

\begin{definition}
Let $\mathcal{F}$ be a coherent $\mathcal{O}_{Y}$-module and $m\in
A_{n-1}\left(  Y\right)  $. Then $\mathcal{F}$ is called
\index{regular sheaf|textbf}%
$m$\textbf{-regular} with respect to $C$ if
\[
H^{i}\left(  Y,\mathcal{F}\otimes\widetilde{S\left(  a\right)  }\right)  =0
\]
for all $i\geq1$ and for all $a\in m+\mathbb{N}C\left[  -i\right]  $.

The
\index{regularity}%
regularity of $\mathcal{F}$ with
\newsym[$\operatorname{reg}_{C}\left(  \mathcal{F}\right)  $]{regularity with respect to $C$}{}respect
to $C$ is%
\[
\operatorname{reg}_{C}\left(  \mathcal{F}\right)  =\left\{  m\in
A_{n-1}\left(  Y\right)  \mid\mathcal{F}\text{ is }m\text{-regular with
respect to }C\right\}
\]

$\mathcal{F}$ is called
\index{regular|textbf}%
$m$\textbf{-regular} if it is $m$-regular
\newsym[$\operatorname{reg}\left(  \mathcal{F}\right)  $]{regularity}{}with
respect to the minimal Hilbert basis of $\mathcal{K}$ and the regularity of
$\mathcal{F}$ is%
\[
\operatorname{reg}\left(  \mathcal{F}\right)  =\left\{  m\in A_{n-1}\left(
Y\right)  \mid\mathcal{F}\text{ is }m\text{-regular}\right\}
\]

\end{definition}

If $\mathcal{F}$ is $m$-regular with respect to $C$, then $\mathcal{F}$ is
$a$-regular with respect to $C$ for all $a\in m+\mathbb{N}C$. The regularity
of a module $M$ relates to the regularity of $\widetilde{M}$ as follows:

\begin{proposition}
\cite{MaSm Multigraded CastelnuovoMumford regularity} If $M$ is a finitely
generated $A_{n-1}\left(  Y\right)  $-graded $S$-module and $m\in
A_{n-1}\left(  Y\right)  $, then $M$ is $m$-regular if and only if
$\widetilde{M}$ is $m$-regular, the natural map
\[
M_{a}\rightarrow H^{0}\left(  Y,\mathcal{F}\otimes\widetilde{S\left(
a\right)  }\right)
\]
is surjective for all $a\in m+\mathbb{N}C$ and
\[
\left(  H_{B\left(  \Sigma\right)  }^{0}\left(  S\right)  \right)  _{a}=0
\]
for all
\[
a\in%
{\displaystyle\bigcup\limits_{j=1}^{e}}
\left(  m+c_{j}+\mathbb{N}C\right)
\]

\end{proposition}

Certain truncations do not change the sheafification:

\begin{lemma}
\label{Lem trunction multigraded module}\cite{MaSm Multigraded
CastelnuovoMumford regularity} Let $C\subset\mathcal{K}$ such that the cone
spanned by $C$ has full dimension, let $m\in\mathbb{Z}\mathcal{K}$ and let
$M^{\prime}$ be
\[
0\rightarrow M\mid_{\left(  m+\mathbb{N}C\right)  }\rightarrow M\rightarrow
M^{\prime}\rightarrow0
\]
Then there is $j>0$ such that
\[
\left(  B\left(  \Sigma\right)  ^{j}M^{\prime}\right)  _{a}=0\text{ }\forall
a\in\mathbb{Z}\mathcal{K}%
\]
so $\widetilde{M^{\prime}}=0$, hence%
\[
\widetilde{M}=\widetilde{M\mid_{\left(  m+\mathbb{N}C\right)  }}%
\]

\end{lemma}

The following proposition allows to pass to initial ideals:

\begin{proposition}
\cite{MaSm Multigraded CastelnuovoMumford regularity} If $>$ is a monomial
ordering on $S$ and $I\subset S$ is an ideal, then
\[
\operatorname{reg}\left(  S/in_{>}\left(  I\right)  \right)  \subset
\operatorname{reg}\left(  S/I\right)
\]
If $I$ is $B\left(  \Sigma\right)  $-saturated and $J=\left(  in_{>}\left(
I\right)  :B\left(  \Sigma\right)  ^{\infty}\right)  $, then%
\[
\operatorname{reg}\left(  S/J\right)  \subset\operatorname{reg}\left(
S/I\right)
\]

\end{proposition}

The Hilbert function of $S$ is $h_{S}\left(  t\right)  =\dim_{\mathbb{C}%
}\left(  S_{t}\right)  $ for $t\in\mathcal{K}$. Consider $\mathcal{K}$ as a
subset of $\mathbb{Z}^{a}\cong A_{n-1}\left(  Y\right)  $. The Hilbert
function of $S$ is given by a polynomial:

\begin{lemma}
\cite{MaSm Uniform bounds on multigraded regularity} There is a
\newsym[$P_{M}$]{Hilbert polynomial}{}polynomial $P_{S}\in$ $\mathbb{Q}\left[
t_{1},...,t_{a}\right]  $ such that $h_{S}\left(  t\right)  =P_{S}\left(
t\right)  $ for all $t\in\mathcal{K}$.
\end{lemma}

More generally if $M$ is a module, then the Hilbert function is given by a
polynomial for all $t\in\mathcal{K}$ sufficiently far from the boundary of
$\mathcal{K}$.

\begin{proposition}
\cite{MaSm Uniform bounds on multigraded regularity} Let $M$ be a finitely
generated graded $S$-module. There is a
\index{Hilbert polynomial}%
polynomial $P_{M}\in$ $\mathbb{Q}\left[  t_{1},...,t_{a}\right]  $ such that
$h_{M}\left(  t\right)  =P_{M}\left(  t\right)  $ for all $t$ in a finite
intersection of translates of $\mathcal{K}$.
\end{proposition}

Saturation does not change the Hilbert polynomial:

\begin{lemma}
\cite{MaSm Uniform bounds on multigraded regularity} Let $M$ be a finitely
generated graded $S$-module. Then%
\[
P_{M}=P_{M/H_{B\left(  \Sigma\right)  }^{0}\left(  M\right)  }%
\]

\end{lemma}

\begin{lemma}
\cite{MaSm Uniform bounds on multigraded regularity} Let $M$ be a finitely
generated graded $S$-module. For all $t\in A_{n-1}\left(  Y\right)  $%
\[
h_{M}\left(  t\right)  -P_{M}\left(  t\right)  =\sum_{i=0}^{n}\left(
-1\right)  ^{i}\left(  H_{B\left(  \Sigma\right)  }^{i}\left(  M\right)
\right)  _{t}%
\]

\end{lemma}

If $M$ is $m$-regular, then $\left(  H_{B\left(  \Sigma\right)  }^{i}\left(
M\right)  \right)  _{t}=0$ for all $i=0,...,n$ and all $t\in m+\mathcal{K}$
with $t\neq m$, hence on the $m$-translate of $\mathcal{K}$ the Hilbert
function of $M$ agrees with its Hilbert polynomial:

\begin{corollary}
\cite{MaSm Uniform bounds on multigraded regularity} Let $M$ be a finitely
generated graded $m$-regular $S$-module. Then%
\[
h_{M}\left(  t\right)  =P_{M}\left(  t\right)
\]
for all $t\in m+\mathcal{K}$ with $t\neq m$.
\end{corollary}

If $I\subset S$ is a monomial ideal, then a Stanley filtration of $S/I$ gives
a bound on the regularity of $I$:

\begin{theorem}
\cite{MaSm Uniform bounds on multigraded regularity} Let $I\subset S$ be a
monomial ideal and let
\[
\mathcal{S}=\left\{  \left(  D_{1},\sigma_{1}\right)  ,...,\left(
D_{s},\sigma_{s}\right)  \right\}
\]
with%
\[
S/\left(  I+\left\langle x^{D_{j+1}},...,x^{D_{s}}\right\rangle \right)  \cong
S_{\sigma_{1}}\left(  -\left[  D_{1}\right]  \right)  \oplus...\oplus
S_{\sigma_{j}}\left(  -\left[  D_{j}\right]  \right)
\]
for $j=1,...,s$ be a Stanley filtration of $S/I$. Then%
\[%
{\textstyle\bigcap\nolimits_{\left(  D,\sigma\right)  \in\mathcal{S}}}
\left(  \left[  D\right]  +\operatorname{reg}\left(  S_{\sigma}\right)
\right)  \subset\operatorname{reg}\left(  S/I\right)
\]

\end{theorem}

\begin{corollary}
Suppose $I$ is $B\left(  \Sigma\right)  $-saturated. Let
\[
\overline{\mathcal{S}}=\left\{  \left(  D,\sigma\right)  \in\mathcal{S}%
\mid\sigma\in\Sigma^{\prime}\right\}
\]
be the subset obtained by removing those Stanley pairs from $\mathcal{S}$,
which correspond to irrelevant ideals in the intersection%
\[
I=%
{\displaystyle\bigcap\limits_{\left(  D,\sigma\right)  \in\mathcal{S}}}
\left\langle y_{r}^{u_{r}+1}\mid r\in\Sigma\left(  1\right)  \text{, }%
r\subset\sigma\text{, }D=%
{\textstyle\sum\nolimits_{r\in\Sigma\left(  1\right)  }}
u_{r}D_{r}\right\rangle
\]
in the sense that they define the empty subset of $Y$. Then%
\[%
{\textstyle\bigcap\nolimits_{\left(  D,\sigma\right)  \in\overline
{\mathcal{S}}}}
\left(  \left[  D\right]  +\operatorname{reg}\left(  S_{\sigma}\right)
\right)  \subset\operatorname{reg}\left(  S/I\right)
\]

\end{corollary}

Let $>$ be a monomial ordering on $\mathbb{Q}\left[  t_{1},...,t_{a}\right]  $
refining the degree ordering with $\deg t_{i}=1$. By $>$ a partial ordering on
the fan $\Sigma$ is given via%
\[
\sigma_{1}>\sigma_{2}\Leftrightarrow\operatorname*{in}\nolimits_{>}\left(
P_{S_{\sigma_{1}}}\left(  t\right)  \right)  >\operatorname*{in}%
\nolimits_{>}\left(  P_{S_{\sigma_{2}}}\left(  t\right)  \right)
\]
refining the ordering of the cones of $\Sigma$ by inclusion, i.e.,
\[
V\left(  \sigma_{1}\right)  \subset V\left(  \sigma_{2}\right)
\Leftrightarrow\sigma_{1}\supset\sigma_{2}\Rightarrow\sigma_{1}<\sigma_{2}%
\]

\begin{algorithm}
\cite{MaSm Uniform bounds on multigraded regularity} Let $>$ be a total
ordering on $\Sigma$ refining above partial ordering. The following algorithm
computes a Stanley filtration%
\[
\mathcal{S}=\left(  \left(  D_{1},\sigma_{1}\right)  ,...,\left(  D_{s}%
,\sigma_{s}\right)  \right)
\]
of the monomial ideal $I\subset S$ such that if $\sigma_{i}\in\Sigma$ and
$D_{i}\neq0$ there is a $j<i$ with
\begin{align*}
\sigma_{j}  &  \in\Sigma\\
\sigma_{j}  &  <\sigma_{i}\\
D_{i}  &  =D_{j}+D_{r}\text{ with }r\in\Sigma\left(  1\right)  \text{ and
}r\subset\sigma_{j}%
\end{align*}

\begin{itemize}
\item If $I$ is a prime ideal and $I=\left\langle y_{r}\mid r\in\Sigma\left(
1\right)  \text{, }r\subset\sigma\right\rangle $ with $\sigma\in\Sigma$, then
return $\left(  \left(  0,\sigma\right)  \right)  $.

\item Otherwise:

\begin{itemize}
\item If $I\not \subset \left\langle y_{r}\mid r\in\Sigma\left(  1\right)
\text{, }r\subset\sigma\right\rangle $ for all $\sigma\in\Sigma$ then choose
$r\in\Sigma\left(  1\right)  $ such that $I\neq\left(  I:\left\langle
y_{r}\right\rangle \right)  \neq\left\langle 1\right\rangle $ .

\item Otherwise: Choose $\sigma\in\Sigma$ minimal with respect to $>$ with the
property%
\[
I\subsetneqq\left\langle y_{r}\mid r\in\Sigma\left(  1\right)  \text{,
}r\subset\sigma\right\rangle
\]
and choose $r\subset\sigma$ such that $I\neq\left(  I:\left\langle
y_{r}\right\rangle \right)  \neq\left\langle 1\right\rangle $
\end{itemize}

Compute Stanley decompositions $\mathcal{S}_{1}$ of $S/\left(  I+\left\langle
y_{r}\right\rangle \right)  $ and%
\[
\mathcal{S}_{2}=\left(  \left(  D_{1},\sigma_{1}\right)  ,...,\left(
D_{s},\sigma_{s}\right)  \right)
\]
of $S/\left(  I:\left\langle y_{r}\right\rangle \right)  $. Return
\[
\mathcal{S}=\mathcal{S}_{1}\operatorname*{join}\left(  \left(  D_{1}%
+D_{r},\sigma_{1}\right)  ,...,\left(  D_{s}+D_{r},\sigma_{s}\right)  \right)
\]

\end{itemize}
\end{algorithm}

With appropriate choice of the isomorphism $A_{n-1}\left(  Y\right)
\cong\mathbb{Z}^{a}$ we may assume that $\mathbb{R}_{\geq0}^{a}\subset
\operatorname*{cpl}\left(  \Sigma\right)  \subset\mathbb{R}^{a}$. From this it
follows that the lead coefficient of a Hilbert polynomial with respect to any
graded ordering of the monomials of $\mathbb{Q}\left[  t_{1},...,t_{a}\right]
$ is positive.

\begin{algorithm}
\label{Alg finitely many monomial ideals}\cite{MaSm Uniform bounds on
multigraded regularity} Let $>$ be a total ordering on $\Sigma$ induced by a
graded ordering $>$ on $\mathbb{Q}\left[  t_{1},...,t_{a}\right]  $, suppose
$\mathbb{R}_{\geq0}^{a}\subset\operatorname*{cpl}\left(  \Sigma\right)  $ and
let $P\left(  t\right)  \in\mathbb{Q}\left[  t_{1},...,t_{a}\right]  $. The
following algorithm returns all $B\left(  \Sigma\right)  $-saturated monomial
ideals with Hilbert polynomial $P\left(  t\right)  $.

\begin{itemize}
\item Let $\mathsf{finished}=\left\{  {}\right\}  $ and $\mathsf{todo}%
=\left\{  \left(  \varnothing,P\left(  t\right)  \right)  \right\}  $.

\item Let $\left(  \mathcal{S},Q\left(  t\right)  \right)  \in\mathsf{todo}$.

For all $\tau\in\Sigma$ and all $E\in\mathbb{Z}^{\Sigma\left(  1\right)  }$,
$E\geq0$ with the following properties

\begin{enumerate}
\item If $\mathcal{S}\neq\varnothing$ there is $\left(  D,\sigma\right)
\in\mathcal{S}$ with $\sigma\leq\tau$.

\item $\operatorname*{in}\nolimits_{>}\left(  Q\left(  t\right)  \right)
=\operatorname*{in}\nolimits_{>}\left(  P_{S_{\tau}}\left(  t\right)  \right)
$

\item $LC_{>}\left(  Q\left(  t\right)  -P_{S_{\tau}}\left(  t\right)
\right)  $ is positive.

\item If $\mathcal{S}=\varnothing$, then $E=0$.

\item If $\mathcal{S}\neq\varnothing$, then there is an $r\in\Sigma\left(
1\right)  $ with $r\subset\tau$ such that $E=D+D_{r}$.
\end{enumerate}

if $Q\left(  t\right)  =P_{S_{\tau}}\left(  t\right)  $, then
\[
\mathsf{finished}=\mathsf{finished}\cup\left\{  \mathcal{S}\cup\left\{
\left(  E,\tau\right)  \right\}  \right\}
\]
else%
\[
\mathsf{todo}=\mathsf{todo}\cup\left\{  \left(  \mathcal{S}\cup\left\{
\left(  E,\tau\right)  \right\}  ,Q\left(  t\right)  -P_{S_{\tau}}\left(
t\right)  \right)  \right\}
\]

\item Return all those monomial ideals%
\[%
{\displaystyle\bigcap\limits_{\left(  D,\sigma\right)  \in\mathcal{S}}}
\left\langle y_{r}^{u_{r}+1}\mid r\in\Sigma\left(  1\right)  \text{, }%
r\subset\sigma\text{, }D=%
{\textstyle\sum\nolimits_{r\in\Sigma\left(  1\right)  }}
u_{r}D_{r}\right\rangle
\]
for $\mathcal{S}\in\mathsf{finished}$, which have Hilbert polynomial $P\left(
t\right)  $.
\end{itemize}

The maximum of $\left\vert \mathcal{S}\right\vert $ for $\mathcal{S}%
\in\mathsf{finished}$ is called the
\index{Gotzmann number|textbf}%
\textbf{Gotzmann number} of $P\left(  t\right)  $.
\end{algorithm}

\begin{proposition}
For given $P\left(  t\right)  \in\mathbb{Q}\left[  t_{1},...,t_{a}\right]  $
there are only finitely many $B\left(  \Sigma\right)  $-saturated monomial
ideals with Hilbert polynomial $P\left(  t\right)  $.
\end{proposition}

Passing to the initial ideal we get:

\begin{theorem}
\label{Thm existence regularity}\cite{MaSm Uniform bounds on multigraded
regularity} Let $I\subset S$ be an $B\left(  \Sigma\right)  $-saturated ideal,
$m$ the Gotzmann number of $P_{S/I}\left(  t\right)  $ and $c\in%
{\textstyle\bigcap\nolimits_{r\in\Sigma\left(  1\right)  }}
\left(  \deg D_{r}+\mathcal{K}\right)  $, then%
\[%
{\textstyle\bigcap\nolimits_{\sigma\in\Sigma}}
\left(  \left(  m-1\right)  c+\operatorname{reg}\left(  S_{\sigma}\right)
\right)  \subset\operatorname{reg}\left(  S/I\right)
\]

\end{theorem}

\subsubsection{Multigraded Hilbert
schemes\label{Sec Multigraded Hilbert schemes}}

Consider the \newsym[$\mathbb{H}_{Y}^{P}$]{Hilbert functor}{}functor
$\mathbb{H}_{Y}^{P}$ with%
\[
\mathbb{H}_{Y}^{P}\left(  R\right)  =\left\{  \mathcal{J}\mid%
\begin{tabular}
[c]{l}%
$\mathcal{J}$ ideal sheaf of a family of subschemes $X\subset Y\times
_{\mathbb{C}}\operatorname*{Spec}R\rightarrow\operatorname*{Spec}R$\\
with Hilbert polynomial $P$%
\end{tabular}
\right\}
\]
for $R\in$\underline{$\mathbb{C}-Alg$} and fixed multigraded Hilbert
polynomial $P\in\mathbb{Q}\left[  t_{1},...,t_{s}\right]  $ with $s=\left\vert
\Sigma\left(  1\right)  \right\vert -n$. By Section
\ref{Sec Homogeneous coordinate representation of subvarieties and sheaves}
there is a one-to-one correspondence%
\[%
\begin{tabular}
[c]{lll}%
$\left\{  \text{ideal sheaves in }\mathbb{H}_{Y}^{P}\left(  R\right)
\right\}  $ & $\leftrightarrows$ & $\left\{  B\left(  \Sigma\right)
\text{-saturated ideals }I\subset S\otimes_{\mathbb{C}}R\right\}  $\\
\multicolumn{1}{c}{$\widetilde{I}$} & \multicolumn{1}{c}{$\leftarrow$} &
\multicolumn{1}{c}{$I$}\\
\multicolumn{1}{c}{$\mathcal{J}$} & \multicolumn{1}{c}{$\mapsto$} &
\multicolumn{1}{c}{$%
{\displaystyle\bigoplus\limits_{a\in A_{n-1}\left(  Y\right)  }}
H^{0}\left(  Y,\mathcal{J}\otimes_{\mathcal{O}_{Y}}\mathcal{O}_{Y}\left(
a\right)  \right)  $}%
\end{tabular}
\]

By Theorem \ref{Thm existence regularity} there is an $m\in\mathcal{K}$ such
that all $B\left(  \Sigma\right)  $-saturated ideals are $m$-regular. With%
\[
I\mid_{m+\mathcal{K}}=S\cdot\left(
{\displaystyle\bigoplus\limits_{a\in m+\mathcal{K}}}
I_{a}\right)
\]
by Lemma \ref{Lem trunction multigraded module} it holds%
\[
\widetilde{I\mid_{m+\mathcal{K}}}=\widetilde{I}%
\]

Define $h:A\cong\mathbb{Z}^{s}\rightarrow\mathbb{N}$ by $h\left(  a\right)
=P\left(  a\right)  $ and let $F$ be the multiplication by monomials on $S$.
Analogously to Corollary \ref{Cor existence finite set degrees} there is a
finite set $D\subset m+\mathcal{K}$ such that for all fields $K\in$%
\underline{$k-Alg$} and for all $L^{\prime}\in\mathbb{H}_{\left(  S_{D}%
,F_{D}\right)  }^{h}\left(  K\right)  $%
\[
\dim\left(  \left(  K\otimes S_{a}\right)  /L_{a}\right)  \leq h\left(
a\right)
\]
for all $a\in m+\mathcal{K}$, where $L\subset K\otimes S$ is the $F$-submodule
generated by $L^{\prime}$. Hence as $D$ is finite by Theorem
\ref{thm Hilbert scheme reduction to finite degree} the Hilbert functor
$\mathbb{H}_{\left(  S,F\right)  }^{h}$ is a subfunctor of $\mathbb{H}%
_{\left(  S_{D},F_{D}\right)  }^{h}$via the restriction map
\[
\mathbb{H}_{\left(  S,F\right)  }^{h}\rightarrow\mathbb{H}_{\left(
S_{D},F_{D}\right)  }^{h},\text{ }L\mapsto L_{D}%
\]
and $\mathbb{H}_{\left(  S,F\right)  }^{h}$ is represented by a closed
subscheme of the Hilbert scheme representing $\mathbb{H}_{\left(  S_{D}%
,F_{D}\right)  }^{h}$. As the Hilbert scheme representing $\mathbb{H}_{\left(
S_{D},F_{D}\right)  }^{h}$ is a closed subscheme of the Grassmann scheme
representing $\mathbb{G}_{S_{D}}^{h}$, the Hilbert scheme representing
$\mathbb{H}_{\left(  S,F\right)  }^{h}$ is projective.

\begin{theorem}
\cite{MaSm Uniform bounds on multigraded regularity} If $P\in\mathbb{Q}\left[
t_{1},...,t_{s}\right]  $ is a multigraded Hilbert polynomial, then
$\mathbb{H}_{Y}^{P}$ is represented by a projective scheme over $\mathbb{C}$.
\end{theorem}

\begin{algorithm}
\cite{MaSm Uniform bounds on multigraded regularity} The following algorithm
computes a subset $D\subset m+\mathcal{K}$ such that for all fields $K\in
$\underline{$k-Alg$} and all $L^{\prime}\in\mathbb{H}_{\left(  S_{D}%
,F_{D}\right)  }^{h}\left(  K\right)  $
\[
\dim\left(  \left(  K\otimes S_{a}\right)  /L_{a}\right)  \leq h\left(
a\right)
\]
for all $a\in m+\mathcal{K}$, where $L\subset K\otimes S$ is the $F$-submodule
generated by $L^{\prime}$.

\begin{enumerate}
\item $D:=\left\{  m\right\}  $

\item Compute by Algorithm \ref{Alg finitely many monomial ideals} the finite
set $M$ of all monomial ideals $I$ generated in degrees $D$ with
$h_{S/I}\left(  t\right)  =P\left(  t\right)  $ $\forall t\in D$.

\item Suppose $I\in M$ with $h_{S/I}\left(  t\right)  \neq P\left(  t\right)
$ for some $t\in m+\mathcal{K}$, then $D:=D\cup\left\{  t\right\}  $ and goto
$2.$ otherwise return $D$.
\end{enumerate}

If $I\in M$, then $J=\left(  I:B\left(  \Sigma\right)  ^{\infty}\right)  $ is
$m$-regular and has Hilbert polynomial $P_{S/J}=P$. Hence $J\mid
_{m+\mathcal{K}}$ is generated in degree $m$. As
\[
h_{S/I}\left(  m\right)  =h_{S/J}\left(  m\right)  =P\left(  m\right)
\]
it holds $J_{m}=I_{m}$, so by $I\subset J$ we have $I\mid_{m+\mathcal{K}%
}=J\mid_{m+\mathcal{K}}$, hence%
\[
h_{S/J}\left(  t\right)  =P\left(  t\right)  \text{ }\forall t\in
m+\mathcal{K}%
\]
So in any step $D$ satisfies the property required above.
\end{algorithm}

\begin{remark}
\label{Rem generalization state polytope}If $Y$ is just simplicial, then one
could replace $S$ by the ring%
\[%
{\displaystyle\bigoplus\limits_{a\in\operatorname*{Pic}\left(  Y\right)  }}
S_{a}%
\]
If $Y$ is a non-simplicial toric variety, then one has to introduce an
equivalence relation identifying different saturated ideals defining the same
subscheme of $Y$.
\end{remark}

\subsubsection{State polytope\label{Sec State polytope general toric}}

Let $I\subset S$ be a $B\left(  \Sigma\right)  $-saturated ideal with
multigraded\ Hilbert polynomial $P\left(  t\right)  $ and define
$h:A\cong\mathbb{Z}^{s}\rightarrow\mathbb{N}$ by $h\left(  a\right)  =P\left(
a\right)  $. Let $m\in\mathcal{K}$ such that all $B\left(  \Sigma\right)
$-saturated ideals are $m$-regular. Consider the finite set $D\subset
m+\mathcal{K}$ as constructed in Section \ref{Sec Multigraded Hilbert schemes}
such that $\mathbb{H}_{\left(  S,F\right)  }^{h}$ is a subfunctor of
$\mathbb{H}_{\left(  S_{D},F_{D}\right)  }^{h}$via the restriction map
\[
\mathbb{H}_{\left(  S,F\right)  }^{h}\rightarrow\mathbb{H}_{\left(
S_{D},F_{D}\right)  }^{h},\text{ }L\mapsto L_{D}%
\]
and $\mathbb{H}_{\left(  S,F\right)  }^{h}$ is represented by a closed
subscheme of the projective Hilbert scheme representing $\mathbb{H}_{\left(
S_{D},F_{D}\right)  }^{h}$.

By Section \ref{Sec Grassmann functor} the functor $\mathbb{G}_{S_{D}}^{h}$ is
a subfunctor of $\mathbb{G}_{S_{D}}^{r}$ with $r=\sum_{a\in D}h\left(
a\right)  $ and the corresponding morphism of schemes is a closed embedding.
Consider the Pl\"{u}cker embedding of $\mathbb{G}_{S_{D}}^{r}\rightarrow
\mathbb{P}\left(  W\right)  $ with $W=%
{\textstyle\bigwedge\nolimits^{\dim S_{D}-r}}
V$ and $V=S_{D}$. So we have closed embeddings%
\[
\mathbb{H}_{\left(  S,F\right)  }^{h}\rightarrow\mathbb{H}_{\left(
S_{D},F_{D}\right)  }^{h}\rightarrow\mathbb{G}_{S_{D}}^{h}\rightarrow
\mathbb{G}_{S_{D}}^{r}\rightarrow\mathbb{P}\left(  W\right)
\]

Denote by $T$ the torus of $Y$. With $\widehat{T}=M$ and $\widehat{T}^{\ast
}=N$ the bilinear pairing between characters and one parameter subgroups of
$T$%
\[%
\begin{tabular}
[c]{lll}%
$\widehat{T}\times\widehat{T}^{\ast}$ & $\rightarrow$ & $\widehat
{\mathbb{C}^{\ast}}=\mathbb{Z}$\\
$\left(  \chi,\lambda\right)  $ & $\mapsto$ & $\left\langle \chi
,\lambda\right\rangle =\chi\circ\lambda$%
\end{tabular}
\]
corresponds to the canonical bilinear pairing
\[
\left\langle -,-\right\rangle :M\times N\rightarrow\mathbb{Z}%
\]

Write the finite set $D=\left\{  \left[  D_{1}\right]  ,...,\left[
D_{p}\right]  \right\}  \subset A_{n-1}\left(  Y\right)  $ and $L_{i}%
=\mathcal{O}_{Y}\left(  D_{i}\right)  $. Choose linearizations of the torus
action $T\times Y\rightarrow Y$ on the $L_{i}$%
\[%
\begin{tabular}
[c]{rclll}%
$T$ & $\times$ & $L_{i}$ & $\overset{\overline{\sigma}}{\rightarrow}$ &
$L_{i}$\\
$id\times\pi$ & $\downarrow$ &  &  & $\downarrow\pi$\\
$T$ & $\times$ & $Y$ & $\overset{\sigma}{\rightarrow}$ & $Y$%
\end{tabular}
\]
which are unique up to translation in $\widehat{T}=M$. We will fix later
particular linearizations.

The action of $T$ on
\[
S_{D}=%
{\textstyle\bigoplus\nolimits_{i=1}^{p}}
H^{0}\left(  Y,L_{i}\right)  =%
{\textstyle\bigoplus\nolimits_{i=1}^{p}}
S_{\left[  D_{i}\right]  }%
\]
induces an action of $T$ on $W$ and the Pl\"{u}cker embedding $\mathfrak{p}%
:\mathbb{G}_{S_{D}}^{r}\mathbb{\rightarrow P}\left(  W\right)  $ is $T$-equivariant.

For $\chi\in\widehat{T}=M$ let%
\[
W_{\chi}=\left\{  v\in W\mid\Lambda v=\chi\left(  \Lambda\right)  v\text{
}\forall\Lambda\in T\right\}
\]
With%
\[
\operatorname*{State}\left(  W\right)  =\left\{  \chi\in M\mid W_{\chi}%
\neq\left\{  0\right\}  \right\}
\]
there is a decomposition%
\[
W=\bigoplus\nolimits_{\chi\in\operatorname*{State}\left(  W\right)  }W_{\chi}%
\]
Denote by $H\left(  I\right)  \in\mathbb{H}_{\left(  S,F\right)  }^{h}$ the
Hilbert point corresponding to $I$ and let $h^{\ast}\in W$ be a representative
of the image of $H\left(  I\right)  $ under the embedding $p:\mathbb{H}%
_{\left(  S,F\right)  }^{h}\mathbb{\rightarrow P}\left(  W\right)  $. Consider
the decomposition of $h^{\ast}$ corresponding to the decomposition of $W$%
\[
h^{\ast}=\sum\nolimits_{\chi\in\operatorname*{State}\left(  W\right)  }%
h_{\chi}%
\]
with $h_{\chi}\in W_{\chi}$. Define%
\[
\operatorname*{State}\left(  h\right)  =\left\{  \chi\in M\mid h_{\chi}%
\neq0\right\}
\]
and the state polytope of $I$ as the convex hull
\[
\operatorname*{State}\left(  I\right)  =\operatorname*{convexhull}\left(
\operatorname*{State}\left(  h\right)  \right)  \subset\widehat{T}%
\otimes_{\mathbb{Z}}\mathbb{R}=M_{\mathbb{R}}%
\]

Let $x_{0},...,x_{n}$ be a $T$-invariant basis of $V$ and%
\[
x_{B}=x_{b_{1}}\wedge...\wedge x_{b_{r}}%
\]
the corresponding $T$-invariant basis of $W$, compatible with the
decomposition of $W=\bigoplus\nolimits_{\chi\in\operatorname*{State}\left(
W\right)  }W_{\chi}$. With respect to the basis $\left(  x_{B}\right)  $ the
representation $\rho:T\rightarrow\operatorname*{GL}\left(  W\right)  $ given
by the action $T\times W\rightarrow W$ is of the form%
\[
\rho\left(  x\right)  =\operatorname*{diag}\left(  x^{m_{1}},...,x^{m_{\dim
W}}\right)
\]
with $m_{i}\in M$.

Let%
\begin{align*}
\lambda &  :\mathbb{C}^{\ast}\rightarrow T\\
\lambda\left(  t\right)   &  =\operatorname*{diag}\left(  t^{w_{1}%
},...,t^{w_{n}}\right)
\end{align*}
be a one parameter subgroup of $T$, then%
\[%
\begin{tabular}
[c]{llll}%
$\rho\circ\lambda:$ & $\mathbb{C}^{\ast}$ & $\rightarrow$ &
$\operatorname*{GL}\left(  W\right)  $\\
& \multicolumn{1}{c}{$t$} & $\mapsto$ & $\operatorname*{diag}\left(
t^{\left\langle w,m_{1}\right\rangle },...,t^{\left\langle w,m_{\dim
W}\right\rangle }\right)  $%
\end{tabular}
\]

With respect to the basis $\left(  x_{B}\right)  $%
\[
h^{\ast}=\left(  \alpha_{1},...,\alpha_{\dim W}\right)
\]
and%
\[
\lambda\left(  t\right)  \cdot h^{\ast}=\operatorname*{diag}\left(
t^{\left\langle w,m_{1}\right\rangle }\alpha_{1},...,t^{\left\langle w,m_{\dim
W}\right\rangle }\alpha_{\dim W}\right)
\]

Hence, with $p:\mathbb{H}_{\left(  S,F\right)  }^{h}\mathbb{\rightarrow
P}\left(  W\right)  $ and the line bundle%
\[
E=p^{\ast}\left(  \mathcal{O}_{\mathbb{P}\left(  W\right)  }\left(  1\right)
\right)
\]
we have%
\[
\mu^{E}\left(  h,\lambda\right)  =\min\left\{  \left\langle w,m_{i}%
\right\rangle \mid\alpha_{i}\neq0\right\}  =\min_{\chi\in\operatorname*{State}%
\left(  h\right)  }\left\langle \chi,\lambda\right\rangle
\]
so by Theorem \ref{Thm Hilbert-Mumford numerical stability} we obtain:

\begin{theorem}
\label{Thm stability}Suppose $Y=X\left(  \Sigma\right)  $ is a smooth toric
variety given by the fan $\Sigma\subset N_{\mathbb{R}}$ and let $S$ be the Cox
ring of $Y$ and $\mathcal{K}=\operatorname{cpl}\left(  \Sigma\right)  \cap
A_{n-1}\left(  Y\right)  $.

Let $I\subset S$ be a $B\left(  \Sigma\right)  $-saturated ideal with\ Hilbert
polynomial $P\left(  t\right)  $, $h$ the corresponding Hilbert function and
$D\subset m+\mathcal{K}$ such that the restriction map gives a closed
embedding $\mathbb{H}_{\left(  S,F\right)  }^{h}\rightarrow\mathbb{H}_{\left(
S_{D},F_{D}\right)  }^{h}$. Fix linearizations of the $T$-action on $Y$ on the
elements of $D$.

Then stability and semi-stability of the Hilbert point $H\left(  I\right)
\in\mathbb{H=H}_{\left(  S,F\right)  }^{h}$ are characterized as
\[%
\begin{tabular}
[c]{lll}%
$H\left(  I\right)  \in\mathbb{H}^{ss}$ & $\Leftrightarrow$ & $0\in
\operatorname*{State}\left(  I\right)  $\\
$H\left(  I\right)  \in\mathbb{H}^{s}$ & $\Leftrightarrow$ & $0\in
\operatorname*{int}\left(  \operatorname*{State}\left(  I\right)  \right)  $%
\end{tabular}
\]

\end{theorem}

\subsection{Toric homogeneous weight vectors and the Gr\"{o}bner
fan\label{1torichomogeneoussetting}}

In the same way as rational graded
\index{weight vector}%
weight vectors on the coordinate ring of $\mathbb{P}^{n}$ are up to multiples
parametrized by $\frac{\mathbb{Z}^{n+1}}{\mathbb{Z}\left(  1,...,1\right)  }$,
we want to parametrize
\index{weight vector}%
weight vectors, i.e., partial orderings given by
\index{weight vector}%
weight vectors for the variables, on the graded pieces of the
\index{Cox ring}%
Cox ring in the general toric setting.

Let $Y=X\left(  \Sigma\right)  $ be a complete toric variety, $v_{1}%
,...,v_{r}$ the
\index{minimal lattice generator}%
minimal lattice generators of the rays of $\Sigma$ forming the rows of the
presentation matrix $A$ of $A_{n-1}\left(  X\left(  \Sigma\right)  \right)  $
in%
\[
0\rightarrow M\overset{A}{\rightarrow}\mathbb{Z}^{\Sigma\left(  1\right)
}\rightarrow A_{n-1}\left(  X\left(  \Sigma\right)  \right)  \rightarrow
0\text{ }%
\]
Let $P=\operatorname*{convexhull}\left(  v_{1},...,v_{r}\right)  $ and
\index{Chow group}%
$S$ be the
\index{Cox ring}%
Cox ring of $Y$. Any rational
\index{weight vector}%
weight vector on $S$ is representable by an element $w\in\operatorname*{Hom}%
_{\mathbb{Z}}\left(  \mathbb{Z}^{\Sigma\left(  1\right)  },\mathbb{Z}\right)
$. Applying $\operatorname*{Hom}\nolimits_{\mathbb{Z}}\left(  \mathbb{-}%
,\mathbb{Z}\right)  $ to above sequence we get%
\[%
\begin{tabular}
[c]{lllrllll}
&  &  & $0$ & $\leftarrow$ & $\operatorname*{Ext}\nolimits_{\mathbb{Z}}%
^{1}\left(  A_{n-1}\left(  X\left(  \Sigma\right)  \right)  ,\mathbb{Z}%
\right)  $ & $\leftarrow$ & \\
$\leftarrow$ & $\overset{=N}{\operatorname*{Hom}\nolimits_{\mathbb{Z}}\left(
M,\mathbb{Z}\right)  }$ & $\overset{\_\circ A}{\leftarrow}$ &
\multicolumn{1}{l}{$\operatorname*{Hom}\nolimits_{\mathbb{Z}}\left(
\mathbb{Z}^{\Sigma\left(  1\right)  },\mathbb{Z}\right)  $} & $\leftarrow$ &
$\operatorname*{Hom}\nolimits_{\mathbb{Z}}\left(  A_{n-1}\left(  X\left(
\Sigma\right)  \right)  ,\mathbb{Z}\right)  $ & $\leftarrow$ & $0$%
\end{tabular}
\]
hence%
\[
\frac{\operatorname*{Hom}\nolimits_{\mathbb{Z}}\left(  \mathbb{Z}%
^{\Sigma\left(  1\right)  },\mathbb{Z}\right)  }{\operatorname*{Hom}%
\nolimits_{\mathbb{Z}}\left(  A_{n-1}\left(  X\left(  \Sigma\right)  \right)
,\mathbb{Z}\right)  }\cong\operatorname*{image}\left(  \_\circ A\right)
\subset N
\]
Now connect the left hand side to the
\index{weight vector}%
weight vectors on the graded pieces of the Cox ring $S$:

Note that scaled
\index{weight vector}%
weight vectors give the same ordering on the monomials. To take this into
account, define the following equivalence relation: For
\[
\overline{w_{1}},\overline{w_{2}}\in\frac{\operatorname*{Hom}%
\nolimits_{\mathbb{Z}}\left(  \mathbb{Z}^{\Sigma\left(  1\right)  }%
,\mathbb{Z}\right)  }{\operatorname*{Hom}\nolimits_{\mathbb{Z}}\left(
A_{n-1}\left(  X\right)  ,\mathbb{Z}\right)  }%
\]
let%
\[
\overline{w_{1}}\sim\overline{w_{2}}:\Leftrightarrow\exists\lambda_{1}%
,\lambda_{2}\in\mathbb{Z}_{>0}:\lambda_{1}\overline{w_{1}}=\lambda
_{2}\overline{w_{2}}%
\]
where $\lambda_{1}\overline{w_{1}}=\overline{\lambda_{1}w_{1}}$ is the induced
$\mathbb{Z}$-module structure inherited from $\operatorname*{Hom}%
\nolimits_{\mathbb{Z}}\left(  \mathbb{Z}^{\Sigma\left(  1\right)  }%
,\mathbb{Z}\right)  $. The map
\[%
\begin{tabular}
[c]{lll}%
$\frac{\operatorname*{Hom}\nolimits_{\mathbb{Z}}\left(  \mathbb{Z}%
^{\Sigma\left(  1\right)  },\mathbb{Z}\right)  }{\operatorname*{Hom}%
\nolimits_{\mathbb{Z}}\left(  A_{n-1}\left(  X\right)  ,\mathbb{Z}\right)  }$
& $\overset{\psi}{\rightarrow}$ & $\left\{  \text{graded wt. vec. on
}S\right\}  $\\
\multicolumn{1}{c}{$\overline{w}$} & $\mapsto$ & partial ordering given by $w$%
\end{tabular}
\
\]
is well defined, as
\begin{align*}
\overline{w_{1}}  &  =\overline{w_{2}}\\
&  \Leftrightarrow\left(  w_{1}-w_{2}\right)  \cdot\in\operatorname*{Hom}%
\nolimits_{\mathbb{Z}}\left(  A_{n-1}\left(  X\right)  ,\mathbb{Z}\right) \\
&  \Leftrightarrow\left(  w_{1}-w_{2}\right)  \cdot\in\ker\left(  \_\circ
A\right) \\
&  \Leftrightarrow\operatorname*{image}\left(  A\right)  \subset\ker\left(
\left(  w_{1}-w_{2}\right)  \cdot\right) \\
&  \Leftrightarrow w_{1}A=w_{2}A\\
&  \Leftrightarrow w_{1}\cdot=w_{2}\cdot\text{ on }\operatorname*{image}%
\left(  A\right) \\
&  \Rightarrow\left(  w_{1}a>0\Leftrightarrow w_{2}a>0\forall a\in
\operatorname*{image}\left(  A\right)  \right) \\
&  \Leftrightarrow\left(  y^{m_{1}}>_{w_{1}}y^{m_{2}}\Leftrightarrow y^{m_{1}%
}>_{w_{2}}y^{m_{2}}\right)  \forall\text{ Cox monomials }y^{m_{1}},y^{m_{2}%
}\text{ with }\deg y^{m_{1}}=\deg y^{m_{2}}\\
&  \Leftrightarrow>_{w_{1}}=>_{w_{2}}\text{ on }S_{\left[  D\right]  }\text{
}\forall\left[  D\right]  \in A_{n-1}\left(  X\right)
\end{align*}
Note that
\begin{align*}
\deg y^{m_{1}}  &  =\deg y^{m_{2}}\\
&  \Leftrightarrow\left[  \sum_{v\in\mathbb{Z}^{\Sigma\left(  1\right)  }%
}m_{1v}D_{v}\right]  =\left[  \sum_{v\in\mathbb{Z}^{\Sigma\left(  1\right)  }%
}m_{2v}D_{v}\right]  \in A_{n-1}\left(  X\right) \\
&  \Leftrightarrow m_{1}=m_{2}\operatorname{mod}\operatorname*{image}\left(
A\right) \\
&  \Leftrightarrow m_{1}-m_{2}\in\operatorname*{image}\left(  A\right)
\end{align*}
Surjectivity of $\psi$ is obvious, and%
\begin{align*}
&  >_{w_{1}}=\text{ }>_{w_{2}}\text{ on }S_{\left[  D\right]  }\forall\left[
D\right]  \in A_{n-1}\left(  X\right) \\
&  \Leftrightarrow\left(  w_{1}a>0\Leftrightarrow w_{2}a>0\ \forall
a\in\operatorname*{image}\left(  A\right)  \right) \\
&  \Leftrightarrow\exists\lambda_{1},\lambda_{2}\in\mathbb{Z}_{>0}:\lambda
_{1}w_{1}a=\lambda_{2}w_{2}a\ \forall a\in\operatorname*{image}\left(
A\right) \\
&  \Leftrightarrow\exists\lambda_{1},\lambda_{2}\in\mathbb{Z}_{>0}:\lambda
_{1}w_{1}\cdot=\lambda_{2}w_{2}\cdot\text{ on }\operatorname*{image}\left(
A\right) \\
&  \Leftrightarrow\left(  \lambda_{1}w_{1}-\lambda_{2}w_{2}\right)  A\\
&  \Leftrightarrow\operatorname*{image}\left(  A\right)  \subset\ker\left(
\left(  \lambda_{1}w_{1}-\lambda_{2}w_{2}\right)  \cdot\right) \\
&  \Leftrightarrow\left(  \lambda_{1}w_{1}-\lambda_{2}w_{2}\right)  \cdot
\in\ker\left(  \_\circ A\right) \\
&  \Leftrightarrow\left(  \lambda_{1}w_{1}-\lambda_{2}w_{2}\right)  \cdot
\in\operatorname*{Hom}\nolimits_{\mathbb{Z}}\left(  A_{n-1}\left(  X\right)
,\mathbb{Z}\right) \\
&  \Leftrightarrow\lambda_{1}\overline{w_{1}}=\lambda_{2}\overline{w_{2}}\\
&  \Leftrightarrow\overline{w_{1}}\sim\overline{w_{2}}%
\end{align*}
hence:

\begin{lemma}
The
\index{weight vector}%
map%
\[%
\begin{tabular}
[c]{lll}%
$\frac{\operatorname*{Hom}\nolimits_{\mathbb{Z}}\left(  \mathbb{Z}%
^{\Sigma\left(  1\right)  },\mathbb{Z}\right)  }{\operatorname*{Hom}%
\nolimits_{\mathbb{Z}}\left(  A_{n-1}\left(  X\right)  ,\mathbb{Z}\right)  }$
& $\overset{\psi}{\rightarrow}$ & $\left\{  \text{graded weight vectors on
}S\right\}  $\\
\multicolumn{1}{c}{$\overline{w}$} & $\mapsto$ & $>_{w}$%
\end{tabular}
\]
is
\index{Cox ring}%
well defined, surjective and%
\[
>_{w_{1}}=\ >_{w_{2}}\Leftrightarrow\quad\left(  \exists\lambda_{1}%
,\lambda_{2}\in\mathbb{Z}_{>0}:\lambda_{1}\overline{w_{1}}=\lambda
_{2}\overline{w_{2}}\right)  \quad\Leftrightarrow:\quad\overline{w_{1}}%
\sim\overline{w_{2}}%
\]

\end{lemma}

\begin{proposition}
After tensoring with $\mathbb{R}$, the map
\[%
\begin{tabular}
[c]{lllll}%
$N=\operatorname*{Hom}\nolimits_{\mathbb{Z}}\left(  M,\mathbb{Z}\right)  $ &
&  &  & \\
\multicolumn{1}{c}{$\cup$} &  &  &  & \\
\multicolumn{1}{r}{$\operatorname*{image}\left(  \_\circ A\right)  $} &
$\overset{\_\circ A}{\underset{\varphi}{\leftrightarrows}}$ & $\frac
{\operatorname*{Hom}\nolimits_{\mathbb{Z}}\left(  \mathbb{Z}^{\Sigma\left(
1\right)  },\mathbb{Z}\right)  }{\operatorname*{Hom}\nolimits_{\mathbb{Z}%
}\left(  A_{n-1}\left(  X\left(  \Sigma\right)  \right)  ,\mathbb{Z}\right)
}$ & $\rightarrow$ & $\left\{  \text{graded wt. vec. on }S\right\}  $\\
\multicolumn{1}{r}{$\sum_{i=1}^{r}w_{i}v_{i}\cdot$} & $\mapsto$ &
$\overline{\left(  w_{1},...,w_{r}\right)  }$ & $\mapsto$ & $>_{w}$%
\end{tabular}
\]
gives a one-to-one correspondence between half lines with origin $0$ in
$N_{\mathbb{R}}$ and the real
\index{weight vector}%
weight vectors on the graded parts
\index{Cox ring}%
of $S$.
\end{proposition}

\begin{lemma}
As $0\in\operatorname*{int}\left(  P\right)  $, there are $a_{i}>0$ such that
$\sum_{i=1}^{r}a_{i}v_{i}=0$, hence via translation by $\left(  a_{1}%
,...,a_{r}\right)  $ any
\index{weight ordering}%
weight ordering is equivalent to a
\index{global ordering}%
global one.
\end{lemma}

We extend the definition of the
\index{Gr\"{o}bner fan}%
Gr\"{o}bner fan to the general toric setting:

\begin{definition}
The
\index{Gr\"{o}bner fan}%
\textbf{Gr\"{o}bner fan }$GF\left(  J\right)  $ of a homogeneous ideal
$J\subset S$ is the
\index{complete fan}%
complete \newsym[$GF\left(  J\right)  $]{Gr\"{o}bner fan}{}polyhedral fan
formed by the cones $\varphi^{-1}\left(  \overline{C_{\varphi\left(  w\right)
}\left(  J\right)  }\right)  \subset N_{\mathbb{R}}$ for $w\in N_{\mathbb{R}}$.
\end{definition}

\begin{proposition}
If $Y$ is a smooth toric variety and $J\subset S$ is a homogeneous ideal, then
$GF\left(  J\right)  =\operatorname*{NF}\left(  \operatorname*{State}\left(
J\right)  \right)  $.
\end{proposition}

Note that the normal fan does not depend on translation of
$\operatorname*{State}\left(  J\right)  $ by choice of linearizations. Note
also, that the state polytope of $J$ and of its saturation have the same
normal fan.

\section{$\mathbb{Q}$-Gorenstein varieties and Fano
polytopes\label{Sec Q-Gorenstein varieties fano polytopes}}

\subsection{Singularities of toric
varieties\label{Singularities of toric varieties}}

Let $N\cong\mathbb{Z}^{n}$, $M=\operatorname*{Hom}\left(  N,\mathbb{Z}\right)
$, let $Y$ be an affine toric variety given by the rational polyhedral
$n$-dimensional cone $\sigma\subset N_{\mathbb{R}}$ and let $v_{1}%
,...,v_{s}\in N$ be the minimal lattice generators of $\sigma$.

\begin{lemma}
\cite{Dais Resolving 3dimensional toric singularities} The affine toric
variety $Y$ is
\index{Q-Gorenstein}%
$\mathbb{Q}$-Gorenstein if and only if there is an $m\in M_{\mathbb{Q}}$ with
$\left\langle m,v_{i}\right\rangle =-1$ $\forall i=1,...,s$.
\end{lemma}

\begin{definition}
The minimal $r\in\mathbb{Z}_{>0}$ such that there is an $m\in M$ with
$\left\langle m,v_{i}\right\rangle =-r$ $\forall i=1,...,s$ is called the
\index{index of a singularity|textbf}%
\textbf{index} of the singularity of $Y$. So $Y$ is
\index{Gorenstein}%
Gorenstein if and only if it has index $1$.
\end{definition}

\begin{lemma}
\cite{Dais Resolving 3dimensional toric singularities} Suppose $Y$ is
\index{Q-Gorenstein}%
$\mathbb{Q}$-Gorenstein and $m\in M_{\mathbb{Q}}$ with $\left\langle
m,v_{i}\right\rangle =-1$ $\forall i=1,...,s$. Then $Y$ is
\index{terminal singularities|textbf}%
terminal if and only if
\[
\sigma\cap\left\{  w\in N\mid\left\langle m,w\right\rangle \geq-1\right\}
=\left\{  0,v_{1},...,v_{s}\right\}
\]
and $Y$ is canonical if and only if%
\[
\sigma\cap\left\{  w\in N\mid\left\langle m,w\right\rangle >-1\right\}
=\left\{  0\right\}
\]

\end{lemma}

\begin{proposition}
\cite{Dais Resolving 3dimensional toric singularities}
\label{thm singularities of toric varieties}If $Y$ is
\index{Q-Gorenstein}%
$\mathbb{Q}$-Gorenstein, then it is log-terminal. If $Y$ is
\index{Gorenstein}%
Gorenstein, then it is
\index{canonical singularities}%
canonical.
\end{proposition}

Let $Y$ be a normal $\mathbb{Q}$-Gorenstein toric variety of dimension $n$,
given by the rational polyhedral fan $\Sigma\subset N_{\mathbb{R}}$. As $Y$ is
$\mathbb{Q}$-Gorenstein, there is a continuous function $\varphi_{K_{Y}%
}:N_{\mathbb{R}}\rightarrow\mathbb{R}_{\geq0}$ such that $\varphi_{K_{Y}}$ is
piecewise linear on the fan $\Sigma$ and $\varphi_{K_{Y}}\left(  \hat
{r}\right)  =1$ for the minimal lattice generators $\hat{r}$ of all rays
$r\in\Sigma\left(  1\right)  $.

\begin{proposition}
\cite{KMM Introduction to the Minimal Model Program} Suppose $\Sigma^{\prime}$
is a refinement of $\Sigma$ inducing a
\index{resolution of singularities}%
resolution of singularities by the birational morphism $f:X\left(
\Sigma^{\prime}\right)  \rightarrow X\left(  \Sigma\right)  $ and denote by
$D_{1},...,D_{r}$ the irreducible components of the exceptional divisor of
$f$. Then $D_{1},...,D_{r}$ have only
\index{normal crossings}%
normal crossings, $D_{1},...,D_{r}$ correspond to the rays of $\Sigma^{\prime
}$ not in $\Sigma$, and%
\[
K_{X\left(  \Sigma^{\prime}\right)  }=f^{\ast}K_{X\left(  \Sigma\right)
}+\sum_{r\in\Sigma^{\prime}\left(  1\right)  \backslash\Sigma\left(  1\right)
}a_{r}D_{r}%
\]
with%
\[
a_{r}=\varphi_{K_{Y}}\left(  \hat{r}\right)  -1
\]

\end{proposition}

In particular, $f$ is crepant if and only if $\varphi_{K_{Y}}\left(  \hat
{r}\right)  =1$ for all $r\in\Sigma^{\prime}\left(  1\right)  \backslash
\Sigma\left(  1\right)  $.

\subsection{Fano polytopes\label{Sec Fano polytopes}}

Section
\ref{Mirror symmetry for singular Calabi-Yau varieties and stringy Hodge numbers}
suggests to consider $\mathbb{Q}$-Gorenstein Fano varieties, so we generalize
Definition \ref{Def Gorenstein Fano} to the following:

\begin{definition}
A normal variety $Y$ is
\index{toric Fano|textbf}%
called
\index{Fano|textbf}%
\textbf{Fano} if some multiple of $-K_{Y}$ is an ample Cartier divisor.
\end{definition}

By $K_{Y}=-\sum_{v\in\Sigma\left(  1\right)  }D_{v}$ a toric variety $Y$ is
$\mathbb{Q}$-Gorenstein if and only if some multiple of $\sum_{v\in
\Sigma\left(  1\right)  }D_{v}$ is
\index{Cartier}%
Cartier.

\begin{lemma}
If $Y$ is a
\index{complete toric variety}%
complete toric variety then it is Fano if and only if some multiple of
$\sum_{v\in\Sigma\left(  1\right)  }D_{v}$
\index{T-Cartier divisor}%
is Cartier and
\index{ample}%
ample if and only if $Y\cong X\left(  \operatorname*{NF}\left(  \Delta
_{-K_{Y}}\right)  \right)  $.
\end{lemma}

\begin{definition}
A polytope $P\subset N_{\mathbb{R}}\cong\mathbb{R}^{n}$ of dimension $n$ is
called a
\index{Fano polytope|textbf}%
\textbf{Fano polytope} if $P$ is
\index{integral polytope}%
integral and $0$ is the unique lattice point in the interior of $P$.
\end{definition}

If $P\subset N_{\mathbb{R}}$ is a Fano polytope, then $P^{\ast}$ is cut out by
the equations $\left\langle m,w_{i}\right\rangle \geq-1$ for the vertices
$w_{i}\in N$, so if $m\in P^{\ast}\cap M$ is a lattice point in the interior
of $P^{\ast}$, then $\left\langle m,w_{i}\right\rangle \in\mathbb{Z}$ and
$\left\langle m,w_{i}\right\rangle >-1$ for all $i$, hence:

\begin{lemma}
\label{lemma Fano polytope 0 unique interior lattice point}If $P\subset
N_{\mathbb{R}}$ is a Fano polytope, then $0$ is the unique interior lattice
point of $P^{\ast}$.
\end{lemma}

\begin{definition}
Denote \newsym[$\Sigma\left(  P\right)  $]{fan over the faces}{}by
$\Sigma\left(  P\right)  $ the
\index{fan over the faces|textbf}%
\textbf{fan over the faces} of $P$.
\end{definition}

By the characterization of ample Cartier divisors in Section
\ref{Divisors on toric varieties}, we obtain:

\begin{proposition}
If $P$ is a Fano polytope, then $X\left(  \Sigma\left(  P\right)  \right)  $
is a toric Fano variety. It is $\mathbb{Q}$-Gorenstein, hence it has
\index{log terminal singularities|textbf}%
log terminal singularities by Proposition
\ref{thm singularities of toric varieties}.
\end{proposition}

Note that the vertices of $P$ are the minimal lattice generators of the rays
of $X\left(  \Sigma\right)  $. From Section
\ref{Singularities of toric varieties} we also get:

\begin{proposition}
If $P\subset N_{\mathbb{R}}\cong\mathbb{R}^{n}$ is a Fano polytope, then it holds:

\begin{enumerate}
\item If $P\cap N=\operatorname*{vert}\left(  P\right)  \cup\left\{
0\right\}  $, i.e., all lattice points of $\partial P$ are vertices, then
$X\left(  \Sigma\left(  P\right)  \right)  $ is terminal.

\item If all facets of $P$ are of the form $P\cap\left\{  w\in N_{\mathbb{R}%
}\mid\left\langle m,w\right\rangle =-1\right\}  $ with integral $m\in M$, then
$X\left(  \Sigma\left(  P\right)  \right)  $ is Gorenstein.
\end{enumerate}
\end{proposition}

So the second condition is equivalent to $P^{\ast}=\Delta_{-K_{Y}}$ being
integral: Writing all facets $F$ of $P$ as $F=P\cap\left\{  w\in
N_{\mathbb{R}}\mid\left\langle m_{F},w\right\rangle =-1\right\}  $ with
$m_{F}\in M$, the $m_{F}$ are the vertices of $P^{\ast}%
=\operatorname*{convexhull}\left\{  m_{F}\mid F\text{ facet of }P\right\}  $.

\begin{proposition}
A Fano polytope $P$ is reflexive if and only if $P^{\ast}$ is integral. Then
$X\left(  \Sigma\left(  P\right)  \right)  $ is a Gorenstein toric Fano
variety, hence it has
\index{canonical singularities}%
canonical singularities by Proposition
\ref{thm singularities of toric varieties}.
\end{proposition}

\begin{proposition}
Suppose $P\subset N_{\mathbb{R}}$ is a Fano polytope, $\Sigma=\Sigma\left(
P\right)  \subset N_{\mathbb{R}}$ is the fan over the faces of $P$ and
$Y=X\left(  \Sigma\right)  $. As $X$ is $\mathbb{Q}$-Gorenstein, there is a
continuous function $\varphi_{K_{Y}}:N_{\mathbb{R}}\rightarrow\mathbb{R}%
_{\geq0}$ such that $\varphi_{K_{Y}}$ is piecewise linear on the fan $\Sigma$
and $\varphi_{K_{Y}}\left(  v\right)  =1$ for all vertices of $P$, i.e.,
$P=\left\{  w\in N_{\mathbb{R}}\mid\varphi_{K_{Y}}\left(  w\right)
\leq1\right\}  $. If $\Sigma^{\prime}$ is a refinement of $\Sigma$ inducing a
\index{resolution of singularities}%
resolution of singularities via the birational morphism $f:X\left(
\Sigma^{\prime}\right)  \rightarrow X\left(  \Sigma\right)  $, then%
\[
K_{X\left(  \Sigma^{\prime}\right)  }=f^{\ast}K_{X\left(  \Sigma\right)
}+\sum_{r\in\Sigma^{\prime}\left(  1\right)  \backslash\Sigma\left(  1\right)
}\left(  \varphi_{K_{Y}}\left(  \hat{r}\right)  -1\right)  D_{r}%
\]
Hence $f$ is crepant if and only if the introduced rays $\Sigma^{\prime
}\left(  1\right)  \backslash\Sigma\left(  1\right)  $ have minimal lattice
generators on the boundary of $P$.
\end{proposition}

\section{The tropical
\index{mirror construction}%
mirror construction for complete intersections in toric
varieties\label{Sec tropical mirror construction for complete intersections}}

In the following, we give a tropical mirror construction for complete
intersections in toric varieties as defined in Section
\ref{Sec Batyrev and Borisov mirror construction}, and we show that the result
coincides with the Batyrev-Borisov mirror.

\subsection{The degeneration for toric complete
intersections\label{degenerationtoriccompleteintersections}}

Consider the setup from Section
\ref{Sec Batyrev and Borisov mirror construction}, i.e., let $Y=\mathbb{P}%
\left(  \Delta\right)  $ be
\index{mirror construction}%
a Gorenstein toric
\index{Fano}%
Fano
\index{toric Fano}%
variety, represented by the
\index{reflexive}%
reflexive polytope $\Delta\subset M_{\mathbb{R}}$, with
\index{normal fan}%
normal fan $\Sigma\subset N_{\mathbb{R}}$ and Cox ring $S$, and
\index{ray}%
let $\Sigma\left(  1\right)  =I_{1}\cup...\cup I_{c}$ be a nef partition, so
$E_{j}=\sum_{v\in I_{j}}D_{v}$ are Cartier,
\index{spanned by global sections}%
spanned by global sections and $\sum_{j=1}^{c}E_{j}=-K_{Y}$. Define
$\Delta_{j}=\Delta_{E_{j}}$ as the polytope of sections of $E_{j}$ and
\begin{align*}
\nabla_{j}  &  =\operatorname*{convexhull}\left\{  \left\{  0\right\}  \cup
I_{j}\right\} \\
\nabla_{BB}^{\ast}  &  =\operatorname*{convexhull}\left(  \Delta_{1}%
\cup...\cup\Delta_{c}\right)
\end{align*}
so%
\begin{align*}
\Delta &  =\Delta_{1}+...+\Delta_{c}\\
\nabla_{BB}  &  =\nabla_{1}+...+\nabla_{c}%
\end{align*}

Consider the monomial degeneration $\mathfrak{X}$ as defined in Section
\ref{Degenerationcompleteintersection}%
\begin{align*}
m_{j}  &  =\prod_{v\in I_{j}}y_{v}\text{ for }j=1,...,c\\
I_{0}  &  =\left\langle m_{j}\mid j=1,...,c\right\rangle \\
I  &  =\left\langle f_{j}=t\cdot g_{j}+m_{j}\mid j=1,...,c\right\rangle
\subset\mathbb{C}\left[  t\right]  \otimes S\\
g_{j}  &  \in S_{\left[  E_{j}\right]  }\text{, }j=1,...,c\text{ general,
reduced with respect to }I_{0}%
\end{align*}
of the complete intersection given by general sections of the Cartier divisors
$E_{1},...,E_{c}$ to the monomial ideal $I_{0}$.

The resolution of $I_{0}$ is given by the Koszul complex $K_{\bullet}$ on
$m=\left(  m_{1},...,m_{c}\right)  $, i.e., the complex of the simplex on
$m_{1},...,m_{c}$,%
\[
0\rightarrow K_{c}\overset{\partial}{\rightarrow}...\overset{\partial
}{\rightarrow}K_{1}\overset{\partial}{\rightarrow}K_{0}%
\]
with%
\begin{align*}
\mathcal{E}  &  =\mathcal{O}_{Y}\left(  E_{1}\right)  \oplus...\oplus
\mathcal{O}_{Y}\left(  E_{c}\right) \\
K_{0}  &  =\mathcal{O}_{Y}\\
K_{p}  &  =\bigwedge\nolimits^{p}\mathcal{E}^{\ast}\text{ for }p=1,...,c
\end{align*}
and the maps $\partial$ are given by contraction with the section $m$ of
$\mathcal{E}$. With respect to the standard frame $e_{i_{1}...i_{p}}=e_{i_{1}%
}\wedge...\wedge e_{i_{p}}$ for $1\leq i_{1}<...<i_{p}\leq c$ of $K_{p}$ we
can write more explicitly%
\begin{align*}
0  &  \rightarrow\mathcal{O}_{Y}\left(  -E_{1}-...-E_{c}\right)
\rightarrow...\rightarrow\bigoplus_{1\leq i_{1}<...<i_{p}\leq c}%
\mathcal{O}_{Y}\left(  -E_{i_{1}}-...-E_{i_{p}}\right)  \rightarrow...\\
...  &  \rightarrow\bigoplus_{i=1}^{c}\mathcal{O}_{Y}\left(  -E_{i}\right)
\rightarrow\mathcal{O}_{Y}%
\end{align*}
and%
\begin{gather*}
\partial:K_{p}\rightarrow K_{p-1}\\
\partial\left(  e_{i_{1}...i_{p}}\right)  =\sum_{j=1}^{p}\left(  -1\right)
^{j-1}m_{i_{j}}e_{i_{1}...i_{j-1}i_{j+1}...i_{p}}%
\end{gather*}

Denote by $\pi_{1}:Y\times\operatorname*{Spec}\mathbb{C}\left[  \left[
t\right]  \right]  \rightarrow Y$ the projection on the first component. So
for above family $\mathfrak{X}$ defined by $I=\left\langle m_{1}%
+tg_{1},...,m_{c}+tg_{c}\right\rangle $, the Koszul complex on $\left(
m_{1}+tg_{1},...,m_{c}+tg_{c}\right)  $ considered as a section of $\pi
_{1}^{\ast}\mathcal{E}$ gives a lift of all syzygies of $\left(
m_{1},...,m_{c}\right)  $, so $\mathfrak{X}$ is flat.

In the same way any first order deformation over $\operatorname*{Spec}\left(
\mathbb{C}\left[  t\right]  /\left\langle t^{2}\right\rangle \right)  $ gives
a deformation over $\operatorname*{Spec}\left(  \mathbb{C}\left[  \left[
t\right]  \right]  \right)  $, hence:

\begin{proposition}
The family $\mathfrak{X}\subset Y\times\operatorname*{Spec}\left(
\mathbb{C}\left[  \left[  t\right]  \right]  \right)  $ defined by $I$ is a
flat
\index{degeneration}%
degeneration with fibers polarized in $Y=\mathbb{P}\left(  \Delta\right)  $
and
\index{monomial degeneration}%
monomial
\index{special fiber}%
special fiber $X_{0}$ given by $I_{0}$. The fiber over the generic point of
$\operatorname*{Spec}\left(  \mathbb{C}\left[  \left[  t\right]  \right]
\right)  $ is a Calabi-Yau complete intersection of codimension $c$ in $Y$
given by general sections of $\mathcal{O}\left(  E_{1}\right)
,...,\mathcal{O}\left(  E_{c}\right)  $.

The
\index{deformation}%
deformations of $I_{0}$ are unobstructed and the base space is smooth.

Let $v_{1},...,v_{p}\in\operatorname*{Hom}\left(  I_{0},S/I_{0}\right)  _{0}$
be a
\index{genericy condition}%
basis of the tangent space of the Hilbert scheme of $X_{0}$. The degeneration
$\mathfrak{X}$ is general in the sense that if $v$ is the tangent vector of
$\mathfrak{X}$ and $v=\sum_{i=1}^{p}\lambda_{i}v_{i}$, then we have
$\lambda_{i}\neq0$ $\forall i$.
\end{proposition}

The ideals of the maximal strata of $X_{0}$ are the ideals $\left\langle
y_{j_{1}},...,y_{j_{c}}\right\rangle \subset S$ for $j_{1}\in I_{1}%
,...,j_{c}\in I_{c}$, hence are given by $c=\operatorname*{codim}\left(
X_{t}\right)  $ equations.

Note that in a toric variety the ideal of a stratum of codimension $c$ in the
Cox ring $S$ can have more than $c$ generators if the toric variety is not
simplicial, and hence the face of $\Delta$ corresponding to the stratum may be
contained in more than $c$ facets, see Example \ref{2example22B}.

\subsection{The Gr\"{o}bner cone associated to the special fiber and the
polytope $\nabla$\label{Groebnerconeassocitatedtothespecialfiber ci}}

Fix a tie break ordering $>$ on $\mathbb{C}\left[  t\right]  \otimes S$ with
$t$ local and respecting the Chow
\index{grading}%
grading on $S$, so $L_{>}\left(  f_{j}\right)  =m_{j}$. Denote by $\varphi$
the map from $N_{\mathbb{R}}$ to the graded weight vectors on $S$ as defined
in Section \ref{1torichomogeneoussetting}. The special fiber
\index{Gr\"{o}bner cone}%
Gr\"{o}bner cone
\[
C_{I_{0}}\left(  I\right)  =\left\{  -\left(  w_{t},w_{y}\right)
\in\mathbb{R}\oplus N_{\mathbb{R}}\mid L_{>_{\left(  w_{t},\varphi\left(
w_{y}\right)  \right)  }}\left(  I\right)  =I_{0}\right\}
\]
is given by%
\[
C_{I_{0}}\left(  I\right)  =\left\{  -\left(  w_{t},w_{y}\right)
\in\mathbb{R}\oplus N_{\mathbb{R}}\mid\operatorname*{trop}\left(
g_{j}\right)  \left(  \varphi\left(  w_{y}\right)  \right)  -w_{t}%
\leq\operatorname*{trop}\left(  m_{j}\right)  \left(  \varphi\left(
w_{y}\right)  \right)  \text{ }\forall j\right\}
\]
Note that the equalities of the lead terms respectively the tropical
inequalities are well defined by homogeneity.

As $\operatorname*{convexhull}\left(  \Delta_{1}\cup...\cup\Delta_{c}\right)
=\nabla_{BB}^{\ast}$ contains $0$ in its interior, for all $w\in
N_{\mathbb{R}}$ there is $j$ and a vertex $0\neq\tilde{m}$ of $\Delta_{j}$
such that $\left\langle \tilde{m},w\right\rangle >0$. Then $m=m_{j}\cdot
A\tilde{m}$ is a monomial of some $g_{j}$with $\varphi\left(  w_{y}\right)
\left(  \frac{m}{m_{j}}\right)  >0$, i.e.,%
\[
\operatorname*{trop}\left(  m\right)  \left(  \varphi\left(  w_{y}\right)
\right)  >\operatorname*{trop}m_{j}\left(  \varphi\left(  w_{y}\right)
\right)
\]
hence:

\begin{lemma}
\label{Lem cone in positive half space}The special fiber Gr\"{o}bner cone
satisfies%
\[
C_{I_{0}}\left(  I\right)  \cap\left\{  w_{t}=0\right\}  =\left\{  0\right\}
\]

\end{lemma}

So $I_{0}$ cannot appear as lead ideal of the general fiber ideal $I_{gen}$.

Intersecting $C_{I_{0}}\left(  I\right)  $
\index{Gr\"{o}bner cone}%
with the hyperplane $\left\{  w_{t}=1\right\}  $ we obtain the convex
polytope
\[
\nabla=C_{I_{0}}\left(  I\right)  \cap\left\{  w_{t}=1\right\}  \subset
N_{\mathbb{R}}%
\]
with%
\[
\nabla=\left\{  -w_{y}\in N_{\mathbb{R}}\mid\operatorname*{trop}\left(
g_{j}\right)  \left(  \varphi\left(  w_{y}\right)  \right)  -1\leq
\operatorname*{trop}\left(  m_{j}\right)  \left(  \varphi\left(  w_{y}\right)
\right)  \text{ }\forall j=1,...,c\right\}
\]
and $0\in\operatorname*{int}\nabla$.

As $\nabla$ is given by integral linear equations corresponding to the
deformations of $I_{0}$ appearing in $I$, the polytope $\nabla^{\ast}$
\index{lattice polytope}%
is
\index{integral polytope}%
integral.

Rewriting the tropical equations we have
\begin{align*}
\nabla &  =\left\{  -w_{y}\in N_{\mathbb{R}}\mid\operatorname*{trop}\left(
m\right)  \left(  \varphi\left(  w_{y}\right)  \right)  -1\leq
\operatorname*{trop}\left(  m_{j}\right)  \left(  \varphi\left(  w_{y}\right)
\right)  \text{ }\forall\text{ monomials }m\text{ of }g_{j}\text{ }\forall
j\right\} \\
&  =\left\{  w_{y}\in N_{\mathbb{R}}\mid\varphi\left(  w_{y}\right)  \left(
\frac{m}{m_{j}}\right)  \geq-1\text{ }\forall\text{ monomials }m\text{ of
}g_{j}\text{ and }\forall j\right\}
\end{align*}

The linear conditions defining $\nabla$ do not change if we do not require the
$f_{j}$ to be
\index{reduced}%
reduced with respect to $I_{0}$. To see this, let $\frac{m}{m_{j}}$ be a
degree $0$ Cox
\index{Laurent monomial}%
Laurent monomial, i.e., $\frac{m}{m_{j}}\in\operatorname*{image}\left(
A\right)  $. Any $\varphi\left(  w_{y}\right)  $ has a positive representative
in $\mathbb{R}^{\Sigma\left(  1\right)  }$, hence if $m_{j}\mid m$, then
\[
\varphi\left(  w_{y}\right)  \left(  \frac{m}{m_{j}}\right)  \geq0\geq-1
\]
If $m\in S_{\left[  E_{j}\right]  }$ is divisible by some
\[
m_{i}=\prod_{v\in I_{i}}y_{v}%
\]
then $A^{-1}\left(  \frac{m}{m_{j}}\right)  $ is an interior point of a face
$F$ of $\Delta_{j}$ of dimension $\dim\left(  F\right)  \geq1$, hence, the
defining inequality of $\nabla$%
\[
\left\langle A^{-1}\left(  \frac{m}{m_{j}}\right)  ,w_{y}\right\rangle \geq-1
\]
given by $m$ is redundant, so we get:

\begin{proposition}
The polytope $\nabla$ is given by
\[
\nabla=\left\{  w_{y}\in N_{\mathbb{R}}\mid\left\langle A^{-1}\left(  \frac
{m}{m_{i}}\right)  ,w_{y}\right\rangle \geq-1\text{ }\forall\text{ monomials
}m\in S_{\left[  E_{j}\right]  }\text{ }\forall j=1,...,c\right\}
\]
and $0\in\nabla^{\ast}$.
\end{proposition}

Reformulating this in terms of lattice monomials%
\begin{align*}
\nabla &  =\left\{  w_{y}\in N_{\mathbb{R}}\mid\left\langle \widetilde
{m},w_{y}\right\rangle \geq-1\text{ }\forall\widetilde{m}\in\Delta_{j}\text{
}\forall j\right\} \\
&  =\left\{  w_{y}\in N_{\mathbb{R}}\mid\left\langle \widetilde{m}%
,w_{y}\right\rangle \geq-1\text{ }\forall\widetilde{m}\in\nabla_{BB}^{\ast
}\right\} \\
&  =\nabla_{BB}%
\end{align*}
by $\nabla_{BB}^{\ast}=\operatorname*{convexhull}\left(  \Delta_{1}\cup
...\cup\Delta_{c}\right)  $, hence:

\begin{theorem}
The polytope $\nabla=C_{I_{0}}\left(  I\right)  \cap\left\{  w_{t}=1\right\}
$ coincides with $\nabla_{BB}$ given in the
\index{mirror construction}%
mirror construction by Batyrev and Borisov%
\[
\nabla=\nabla_{BB}%
\]

\end{theorem}

\begin{corollary}
$\nabla$ is
\index{reflexive}%
reflexive, so it defines a Gorenstein
\index{toric Fano}%
toric
\index{Fano}%
Fano variety $Y^{\circ}=\mathbb{P}\left(  \nabla\right)  $.
\end{corollary}

Denote by $S^{\circ}=\mathbb{C}\left[  z_{r}\mid r\in\Sigma^{\circ}\left(
1\right)  \right]  $ with $\Sigma^{\circ}=\operatorname*{NF}\left(
\nabla\right)  $ the Cox ring of $Y^{\circ}$.

\begin{figure}
[h]
\begin{center}
\includegraphics[
height=2.3592in,
width=2.0833in
]%
{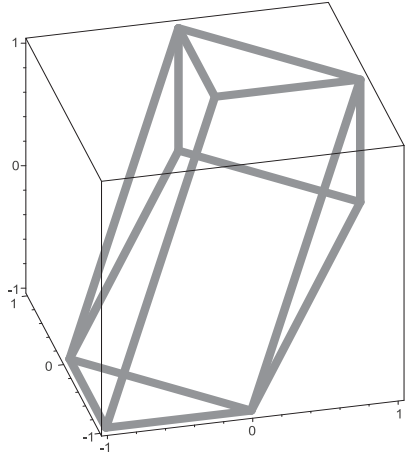}%
\caption{Reflexive supporting polyhedron of the special fiber Gr\"{o}bner cone
for the monomial degeneration of the complete intersection of two general
quadrics in $\mathbb{P}^{3}$}%
\label{Fig Nabla 22}%
\end{center}
\end{figure}

\begin{example}
\label{22Nabla}Consider the monomial degeneration $\mathfrak{X}\subset
\mathbb{P}^{3}\times\operatorname*{Spec}\mathbb{C}\left[  \left[  t\right]
\right]  $ over $\operatorname*{Spec}\mathbb{C}\left[  \left[  t\right]
\right]  $ of an
\index{elliptic curve!complete intersection}%
elliptic curve given as the complete intersection of two general quadrics in
$\mathbb{P}^{3}$ as given in Example \ref{22degen}, i.e., by the ideal
\[
I=\left\langle t\cdot g_{1}+x_{1}x_{2},\ t\cdot g_{2}+x_{0}x_{3}\right\rangle
\subset\mathbb{C}\left[  t\right]  \otimes S
\]
where $g_{1},g_{2}\in\mathbb{C}\left[  x_{0},...,x_{3}\right]  _{2}$ are
general, not involving monomials in $I_{0}=\left\langle x_{1}x_{2},x_{0}%
x_{3}\right\rangle $. Here $S=\mathbb{C}\left[  x_{0},...,x_{3}\right]  $
denotes the
\index{Cox ring}%
Cox ring of $\mathbb{P}\left(  \Delta\right)  =\mathbb{P}^{3}$ with variables
$x_{0},...,x_{3}$ corresponding to the vertices of $\Delta^{\ast}$ and
$\Delta$ denotes the degree $4$ Veronese polytope of $\mathbb{P}^{3}$.

For this example the reflexive polytope $\nabla=C_{I_{0}}\left(  I\right)
\cap\left\{  w_{t}=1\right\}  $ is depicted in Figure \ref{Fig Nabla 22}. As
shown above it agrees with $\nabla_{BB}$ given in Example \ref{22example1}.
\end{example}

%

\subsection{The initial ideals of the faces of $\nabla$%
\label{Sec initial ideals of the faces ci}}

In the following we explicitly give the correspondence of lead ideals of $I$
and faces of $\nabla^{\ast}$.

Consider the notation from the last section. Note that all $w$ in the interior
of a face $F$ of $\nabla$ lead to the same
\index{initial ideal}%
initial ideal of $I$ denoted by $\operatorname*{in}_{F}\left(  I\right)  $.
Let $F$ be a face of $\nabla$ and $m$ a monomial of $g_{j}$. Then $m$ is a
monomial of $\operatorname*{in}_{F}\left(  g_{j}\right)  $ if and only if
$m\in S_{\left[  E_{j}\right]  }$ and%
\[
\varphi\left(  w_{y}\right)  \left(  m\right)  +1=\varphi\left(  w_{y}\right)
\left(  m_{j}\right)  \text{ }\forall w_{y}\in F
\]
if and only if $A^{-1}\left(  \frac{m}{m_{i}}\right)  \in\Delta_{j}\cap M$
and
\[
\left\langle A^{-1}\left(  \frac{m}{m_{i}}\right)  ,w_{y}\right\rangle
=-1\text{ }\forall w_{y}\in F
\]
if and only if there is an $\widetilde{m}\in\Delta_{j}\cap M$ with
$m=m_{j}\cdot A\left(  \widetilde{m}\right)  $ and%
\[
\left\langle \widetilde{m},w_{y}\right\rangle =-1\text{ }\forall w_{y}\in F
\]
hence:

\begin{lemma}
\label{inF}The monomials appearing
\index{initial form}%
in $\operatorname*{in}_{F}\left(  g_{j}\right)  $ are%
\[
\left\{  m_{j}\right\}  \cup\left\{  m_{j}\cdot A\left(  \widetilde{m}\right)
\mid\widetilde{m}\in\Delta_{j}\cap M\text{ with }\left\langle \widetilde
{m},w_{y}\right\rangle =-1\text{ }\forall w_{y}\in F\text{ and }m_{j}\cdot
A\left(  \widetilde{m}\right)  \notin I_{0}\right\}
\]

\end{lemma}

\subsection{The
\index{dual complex}%
dual complex of $\nabla$\label{Sec dual complex ci}}

If $F$ is a face of $\nabla$
\index{initial form}%
write%
\[
\operatorname*{in}\nolimits_{F}\left(  f_{j}\right)  =t\sum_{m\in G_{j}\left(
F\right)  }c_{m}m+m_{j}%
\]
for $j=1,...,c$.

\begin{definition}
If $F$ is a face of $\nabla$, then define the
\index{dual|textbf}%
\textbf{dual face} of $F$ as the convex hull of all first order
\index{initial form}%
deformations appearing the initial ideal of $I$ with respect to $F$
\[
\operatorname*{dual}\left(  F\right)  =\operatorname*{convexhull}\left(
A^{-1}\left(  \frac{m}{m_{j}}\right)  \mid m\in G_{j}\left(  F\right)  ,\text{
}j=1,...,c\right)  \subset M_{\mathbb{R}}%
\]

\end{definition}

Then we have
\[
\operatorname*{dual}\left(  F\right)  =\operatorname*{convexhull}\left(
\bigcup_{j=1}^{c}\left\{  A^{-1}\left(  \frac{m}{m_{j}}\right)  \mid m\text{ a
monomial of }\operatorname*{in}\nolimits_{F}\left(  g_{j}\right)  \right\}
\right)
\]
so by Lemma \ref{inF}%
\begin{align*}
\operatorname*{dual}\left(  F\right)   &  =\operatorname*{convexhull}\left(
\bigcup_{j=1}^{c}\left\{  \widetilde{m}\in\Delta_{j}\cap M\mid\left\langle
\widetilde{m},w_{y}\right\rangle =-1\text{ }\forall w_{y}\in F\right\}
\right) \\
&  =\operatorname*{convexhull}\left(  \left\{  \widetilde{m}\in\bigcup
_{j=1}^{c}\Delta_{j}\cap M\mid\left\langle \widetilde{m},w_{y}\right\rangle
=-1\text{ }\forall w_{y}\in F\right\}  \right) \\
&  =\left\{  \widetilde{m}\in\nabla^{\ast}\mid\left\langle \widetilde{m}%
,w_{y}\right\rangle =-1\text{ }\forall w_{y}\in F\right\} \\
&  =F^{\ast}%
\end{align*}
\label{dualFFdual}hence:

\begin{proposition}
If $F$ is a face of $\nabla$, then%
\[
\operatorname*{dual}\left(  F\right)  =F^{\ast}\subset\nabla^{\ast}%
\]
in particular $\operatorname*{dual}\left(  F\right)  $ is a face of
$\nabla^{\ast}$, so%
\[
\operatorname*{dual}:\operatorname*{Poset}\left(  \nabla\right)
\rightarrow\operatorname*{Poset}\left(  \nabla^{\ast}\right)
\]
is the inclusion reversing map from the face poset of $\nabla$ to the face
poset of $\nabla^{\ast}$ given by dualization of the face.
\end{proposition}

\begin{figure}
[h]
\begin{center}
\includegraphics[
height=2.1949in,
width=5.4155in
]%
{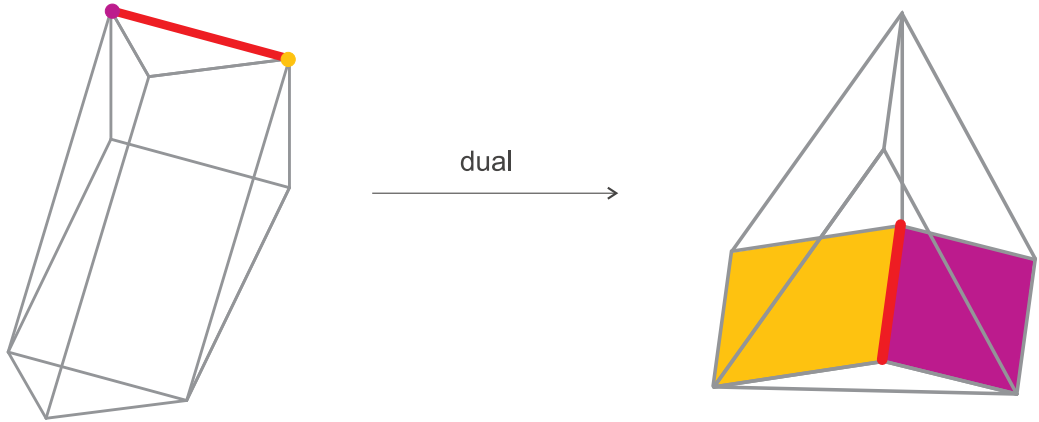}%
\caption{Faces of initial ideals for the monomial degeneration of the complete
intersection of two general quadrics in $\mathbb{P}^{3}$}%
\label{Fig dual Nabla 22}%
\end{center}
\end{figure}

\begin{example}
The complex of initial ideals $\operatorname*{dual}\left(  \nabla\right)  $
for above Example \ref{22Nabla} is visualized in Figure
\ref{Fig dual Nabla 22}. Some faces of $\nabla$ and their corresponding
\index{elliptic curve!complete intersection}%
images under $\operatorname*{dual}$ are highlighted.
\end{example}

%

\subsection{The Bergman subcomplex of $\nabla$%
\label{Sec Bergman subcomplex ci}}

Intersecting
\index{tropical subcomplex}%
the Bergman complex
\index{Bergman subcomplex|textbf}%
with $\nabla$, we obtain the following subcomplex of dimension $d$ of the
boundary complex of $\nabla$%
\[
B\left(  I\right)  =BC_{I_{0}}\left(  I\right)  =\left(  BF\left(  I\right)
\cap\operatorname*{Poset}\left(  C_{I_{0}}\left(  I\right)  \right)  \right)
\cap\left\{  w_{t}=1\right\}
\]
the Bergman subcomplex
\newsym[$BC_{I_{0}}\left(  I\right)$]{special fiber Bergman complex}{}or
\newsym[$B\left(  I\right)$]{special fiber Bergman complex}{}tropical
subcomplex of $\nabla$. Here $\operatorname*{Poset}\left(  C_{I_{0}}\left(
I\right)  \right)  $ is the fan generated by the cone $C_{I_{0}}\left(
I\right)  $. The intersection of the fan $BF\left(  I\right)  \cap
\operatorname*{Poset}\left(  C_{I_{0}}\left(  I\right)  \right)  $ with the
hyperplane $\left\{  w_{t}=1\right\}  $ is defined as the complex whose faces
are the intersections of the cones of the fan with $\left\{  w_{t}=1\right\}
$.

\begin{figure}
[h]
\begin{center}
\includegraphics[
height=2.8409in,
width=2.4587in
]%
{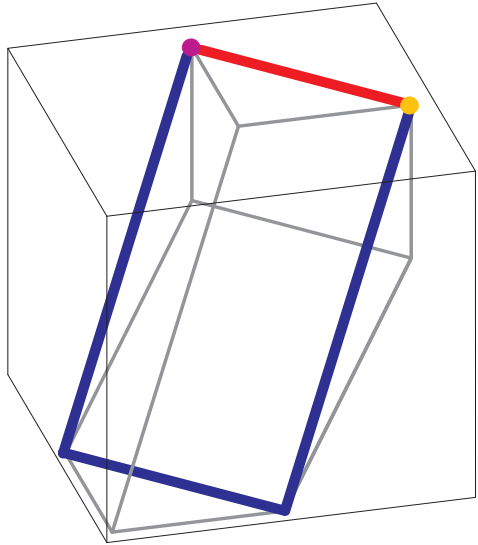}%
\caption{Bergman subcomplex for the monomial degeneration of the complete
intersection of two general quadrics in $\mathbb{P}^{3}$}%
\label{Fig Bergman subcomplex 22}%
\end{center}
\end{figure}

\begin{example}
\label{2example22B}For above
\index{elliptic curve!complete intersection}%
Example \ref{22Nabla} the tropical subcomplex $B\left(  I\right)
\subset\operatorname*{Poset}\left(  \nabla\right)  $ is shown in Figure
\ref{Fig Bergman subcomplex 22}.
\end{example}

%

\subsection{The mirror complex\label{Sec mirror complex complete intersection}%
}

If $F$ is a face
\index{mirror complex}%
of $B\left(  I\right)  $
\index{initial form}%
write%
\[
\operatorname*{in}\nolimits_{F}\left(  f_{j}\right)  =t\sum_{m\in G_{j}\left(
F\right)  }m+m_{j}%
\]
for $j=1,...,c$, then $G_{j}\left(  F\right)  \neq\emptyset$ $\forall j$. Then
we \newsym[$\mu\left(  F\right)  $]{mirror face}{}define the map
\[
\mu:B\left(  I\right)  \rightarrow\operatorname*{Poset}\left(  \Delta\right)
\]
mapping a face of $B\left(  I\right)  $ to
\index{mirror}%
the
\index{Minkowski sum}%
Minkowski sum of the initial forms%
\begin{align*}
\mu\left(  F\right)   &  =\sum_{j=1}^{c}\operatorname*{convexhull}\left(
A^{-1}\left(  \frac{m}{m_{j}}\right)  \mid m\in G_{j}\left(  F\right)  \right)
\\
&  =\sum_{j=1}^{c}\operatorname*{convexhull}\left(  \tilde{m}\mid\tilde{m}%
\in\Delta_{j},\left\langle w,\tilde{m}\right\rangle =-1\forall w\in F\right)
\end{align*}

\begin{proposition}
$\mu\left(  B\left(  I\right)  \right)  $ is a
\index{mirror complex}%
subcomplex of $\operatorname*{Poset}\left(  \Delta\right)  $
\index{mirror}%
and $\mu$ induces an isomorphism of complexes%
\[
B\left(  I\right)  ^{\vee}\rightarrow\mu\left(  B\left(  I\right)  \right)
\subset\operatorname*{Poset}\left(  \Delta\right)
\]
If $F$ is a facet of $B\left(  I\right)  $, then%
\[
\dim\left(  \mu\left(  F\right)  \right)  =n-c-\dim\left(  F\right)
=d-\dim\left(  F\right)
\]

\end{proposition}

\begin{figure}
[h]
\begin{center}
\includegraphics[
height=2.5918in,
width=5.8703in
]%
{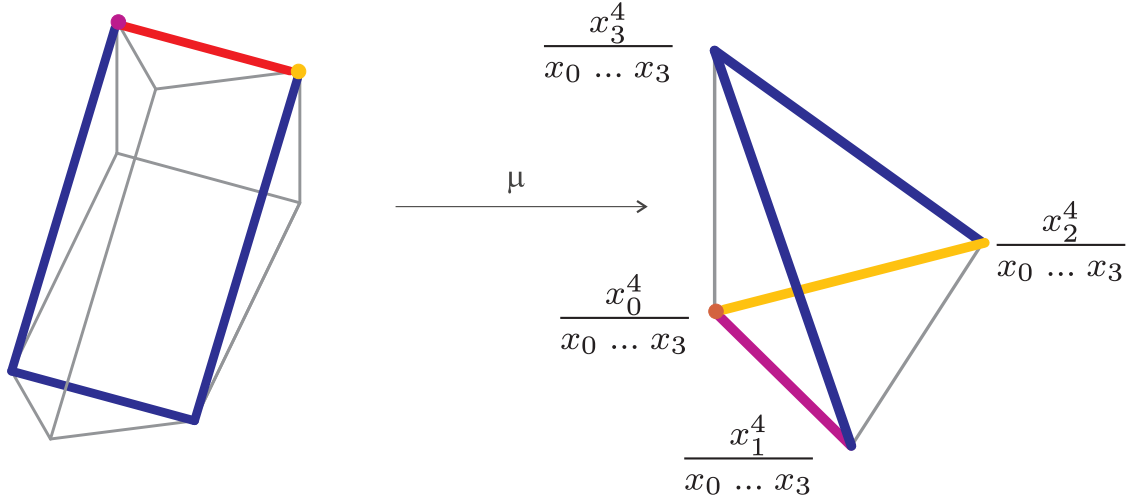}%
\caption{Mirror dual complex of the Bergman subcomplex for the monomial
degeneration of the complete intersection of two general quadrics in
$\mathbb{P}^{3}$}%
\label{Fig mirror complex 22}%
\end{center}
\end{figure}

\begin{example}
For above
\index{elliptic curve!complete intersection}%
Example \ref{22Nabla}, the complexes $B\left(  I\right)  \subset\nabla$ and
$\mu\left(  B\left(  I\right)  \right)  \subset\Delta$
\index{mirror complex}%
are shown in Figure \ref{Fig mirror complex 22}.
\end{example}

%

\subsection{The dual complex of $B\left(  I\right)  $%
\label{Sec dual complex BI ci}}

Now consider the image of $B\left(  I\right)  $ under the map
$\operatorname*{dual}$. By Lemma \ref{dualFFdual} the complex
$\operatorname*{dual}\left(  B\left(  I\right)  \right)  $
\index{dual complex}%
is a
\index{dual}%
subcomplex of $\nabla^{\ast}$. If $F$ is a face of $B\left(  I\right)  $,
then
\[
\dim\left(  \operatorname*{dual}\left(  F\right)  \right)  =n-1-\dim F
\]

Intersecting $\operatorname*{dual}\left(  F\right)  $
\index{dual}%
with $\Delta_{j}$, we can recover the
\index{initial form}%
initial monomials of the individual equations, as
\begin{align*}
\operatorname*{dual}\left(  F\right)  \cap\Delta_{j}  &
=\operatorname*{convexhull}\left\{  \widetilde{m}\in\nabla^{\ast}\cap
M\cap\Delta_{j}\mid\left\langle \widetilde{m},w_{y}\right\rangle =-1\text{
}\forall w_{y}\in F\right\} \\
&  =\operatorname*{convexhull}\left(  \widetilde{m}\in\Delta_{j}\cap
M\mid\left\langle \widetilde{m},w_{y}\right\rangle =-1\text{ }\forall w_{y}\in
F\right) \\
&  =\operatorname*{convexhull}\left(  A^{-1}\left(  \frac{m}{m_{j}}\right)
\mid m\text{ a monomial of }\operatorname*{in}\nolimits_{F}\left(
g_{j}\right)  \right)
\end{align*}
by Lemma \ref{inF}, hence:

\begin{lemma}
\label{dualFdeltaj}If $F$ is a face of the special fiber Bergman complex
$B\left(  I\right)  $, then the intersection of its dual face
$\operatorname*{dual}\left(  F\right)  \subset\nabla$ with $\Delta_{j}%
\subset\nabla$ is the face of $\Delta_{j}$ given as the convex hull of the
deformations appearing in $\operatorname*{in}\nolimits_{F}\left(
f_{j}\right)  $, i.e.,%
\[
\operatorname*{dual}\left(  F\right)  \cap\Delta_{j}%
=\operatorname*{convexhull}\left(  A^{-1}\left(  \frac{m}{m_{j}}\right)  \mid
m\text{ a monomial of }\operatorname*{in}\nolimits_{F}\left(  g_{j}\right)
\right)
\]

\end{lemma}

So we
\index{dual}%
obtain:

\begin{proposition}
If $F$ is a face
\index{dual}%
of $B\left(  I\right)  $, then $\mu\left(  F\right)  $ is the Minkowski sum%
\[
\mu\left(  F\right)  =\sum_{j=1}^{c}\operatorname*{dual}\left(  F\right)
\cap\Delta_{j}%
\]

\end{proposition}

\begin{figure}
[h]
\begin{center}
\includegraphics[
height=4.4884in,
width=3.9418in
]%
{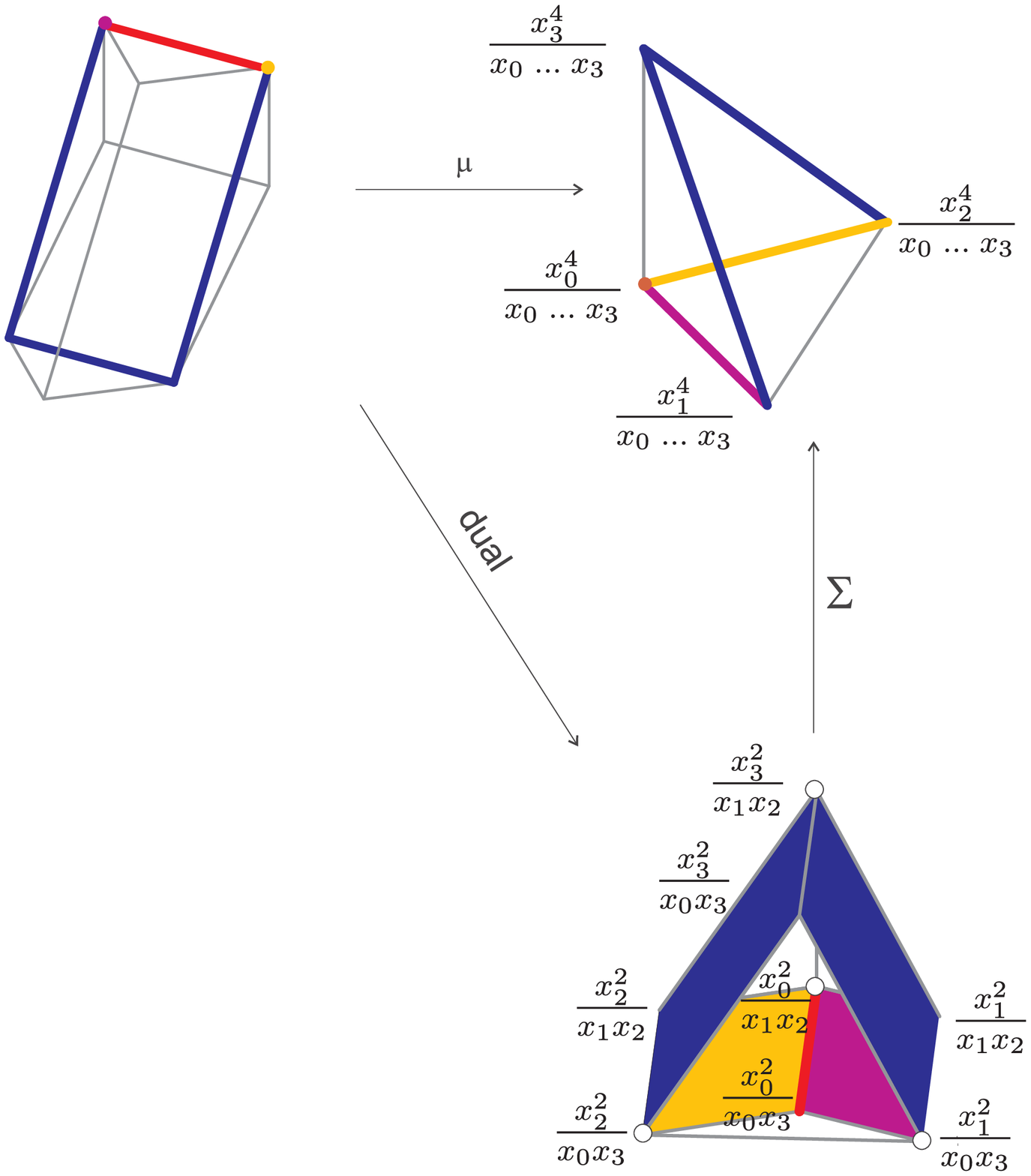}%
\caption{Complex of initial ideals for the monomial degeneration of the
complete intersection of two general quadrics in $\mathbb{P}^{3}$}%
\label{Fig dual complex 22}%
\end{center}
\end{figure}

\begin{example}
\label{Ex dual B for 22 complete intersection}For
\index{Bergman subcomplex}%
above
\index{elliptic curve!complete intersection}%
Example \ref{22Nabla} the complex
\index{dual}%
$\operatorname*{dual}\left(  B\left(  I\right)  \right)  \subset\nabla^{\ast
}$
\index{dual complex}
is shown in Figure \ref{Fig dual complex 22}.
\end{example}

%

\subsection{$B\left(  I\right)  $ and the complex of strata of $X_{0}%
$\label{Sec intersection complex ci}}

Consider the
\index{lim|textbf}%
map $\lim$ given by%
\[%
\begin{tabular}
[c]{llll}%
$\lim:$ & $B\left(  I\right)  $ & $\rightarrow$ & $\operatorname*{Strata}%
\left(  Y\right)  $\\
& \multicolumn{1}{c}{$F$} & $\mapsto$ & $\left\{  \lim_{t\rightarrow0}a\left(
t\right)  \mid a\in\operatorname*{val}\nolimits^{-1}\left(
\operatorname*{int}\left(  F\right)  \right)  \right\}  $%
\end{tabular}
\]
where $\operatorname*{Strata}\left(  Y\right)  $ denotes
\index{Strata}%
the
\index{poset}%
poset of closures of the
\index{toric strata}%
toric strata of the toric variety $Y=\frac{\mathbb{C}^{\Sigma\left(  1\right)
}-V\left(  B\left(  \Sigma\right)  \right)  }{G\left(  \Sigma\right)  }$.

\begin{proposition}
If $F$ is a face of $B\left(  I\right)  $, then $\lim\left(  F\right)
=V\left(  \left(  \mu\left(  F\right)  \right)  ^{\ast}\right)  $ and the
complexes $\lim\left(  B\left(  I\right)  \right)  $ and $\mu\left(  B\left(
I\right)  \right)  $ are isomorphic.
\end{proposition}

Note that the $k$-dimensional
\index{torus orbit closure}%
orbit closures $V\left(  \sigma\right)  $ correspond to the cones $\sigma$ of
$\Sigma$ of dimension $n-k$ (i.e., faces of $\Delta^{\ast}$ of dimension
$n-k-1$) .

\begin{example}
In the above
\index{elliptic curve!complete intersection}%
Example \ref{22Nabla} for $w=\left(  1,0,1\right)  =\left(  1,0,0\right)
+\left(  0,0,1\right)  $ we have $\varphi\left(  w\right)  =\left(
0,1,0,1\right)  +\mathbb{Z}\left(  1,1,1,1\right)  $ and%
\[
\lim\left(  \left\{  w\right\}  \right)  =V\left(  x_{1},x_{3}\right)
\subset\frac{\mathbb{C}^{4}-\left\{  0\right\}  }{\mathbb{C}^{\ast}}%
\]

\end{example}

Denote by $\operatorname*{Strata}\nolimits_{\Delta}\left(  I_{0}\right)  $ the
complex of faces of $\Delta$ corresponding to the strata in $Y$ of the reduced
monomial ideal $I_{0}$.

\begin{proposition}
The map%
\[%
\begin{tabular}
[c]{ccc}%
$B\left(  I\right)  ^{\vee}$ & $\rightarrow$ & $\operatorname*{Strata}%
\nolimits_{\Delta}\left(  I_{0}\right)  $\\
$F^{\vee}$ & $\mapsto$ & \multicolumn{1}{l}{$\lim\left(  F\right)  $}%
\end{tabular}
\]
is an isomorphism of complexes and%
\[
\dim\left(  \lim\left(  F\right)  \right)  =n-1-\dim\left(  \mu\left(
F\right)  \right)  ^{\ast}=\dim\left(  \mu\left(  F\right)  \right)
=d-\dim\left(  F\right)
\]

\end{proposition}

\begin{figure}
[h]
\begin{center}
\includegraphics[
height=4.4391in,
width=4.1814in
]%
{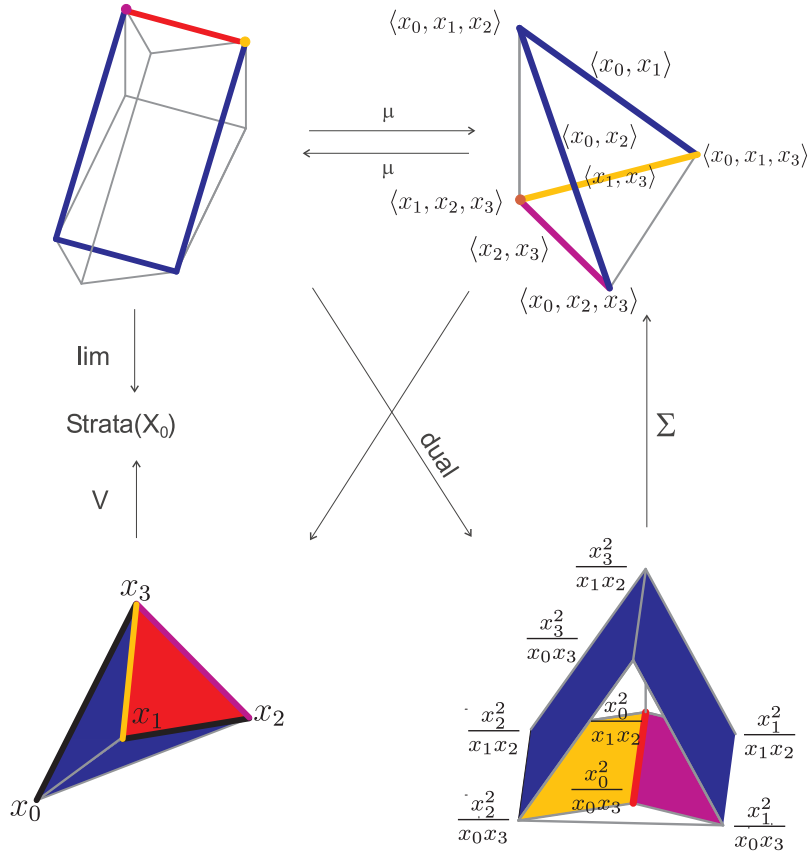}%
\caption{The image of $\lim$ for the monomial degeneration of the complete
intersection of two general quadrics in $\mathbb{P}^{3}$}%
\label{Fig uebersicht lim 22}%
\end{center}
\end{figure}

\begin{example}
For above
\index{elliptic curve!complete intersection}%
Example \ref{22Nabla} the relation between the maps $\lim$ and $\mu$ is shown
in Figure \ref{Fig uebersicht lim 22}.
\end{example}

%

\subsection{Remark on the topology of $B\left(  I\right)  $%
\label{Sec Remark on the toplogy of BI}}

In the case of a toric hypersurface $B\left(  I\right)  $ and $\lim\left(
B\left(  I\right)  \right)  \cong\mu\left(  B\left(  I\right)  \right)  $ are
the boundaries of the polytopes $\nabla$ and $\Delta$, respectively, hence are
homeomorphic to spheres.

\begin{remark}
Consider a
\index{degeneration}%
degeneration obtained from a general
\index{complete intersection}%
complete intersection inside projective space $\mathbb{P}^{n}=\mathbb{P}%
\left(  \Delta\right)  $ given by the partition $\Sigma\left(  1\right)
=I_{1}\cup...\cup I_{c}$ and denote by $S\left(  I_{j}\right)  $ the simplex
on $I_{j}$.
\index{mirror}%
Then $\lim\left(  B\left(  I\right)  \right)  \cong\mu\left(  B\left(
I\right)  \right)  \cong\operatorname*{Strata}\left(  X_{0}\right)  $ is
isomorphic to the join
\[
S\left(  I_{1}\right)  \ast...\ast S\left(  I_{c}\right)  =S^{\left\vert
I_{1}\right\vert -2}\ast...\ast S^{\left\vert I_{c}\right\vert -2}\cong
S^{n-2c+c}=S^{d}%
\]
via the complementary numbering%
\[
\mu\left(  F\right)  \mapsto\left\{  r\in\Sigma\left(  1\right)  \mid
r\not \subset \operatorname*{hull}\left(  \left(  \mu\left(  F\right)
\right)  ^{\ast}\right)  \right\}
\]

\end{remark}

Recall that the
\index{join|textbf}%
join of complexes $C_{1},C_{2}$ is the complex%
\[
C_{1}\ast C_{2}=\left\{  f\vee g\mid f\in C_{1},g\in C_{2}\right\}
\]
where $\vee$ denotes the disjoint union.

The complex $\operatorname*{Strata}\left(  X_{0}\right)  $ is homoemorphic to
a sphere also in the general complete intersection setup given by a nef
partition of $\Delta$. So, as the dual cell complex of $\lim\left(  B\left(
I\right)  \right)  $, also $B\left(  I\right)  $ is homeomorphic to a sphere.
Note that, as we will see below, the complex $B\left(  I\right)  $ corresponds
to the nef partition of $\nabla$ dual to the nef partition of $\Delta$.

\subsection{Covering of $B\left(  I\right)  ^{\vee}$ and reconstruction of $I$
from the tropical data\label{Sec covering complete intersection}}

The
\index{dual complex}%
map%
\[%
\begin{tabular}
[c]{lll}%
$\operatorname*{dual}\left(  B\left(  I\right)  \right)  $ & $\longrightarrow$
& $\mu\left(  B\left(  I\right)  \right)  $\\
\multicolumn{1}{c}{$F^{\ast}$} & $\mapsto$ & $\sum_{j=1}^{c}F^{\ast}\cap
\Delta_{j}=\mu\left(  F\right)  $%
\end{tabular}
\]
induces a $c:1$
\index{covering}%
covering of complexes:

\begin{proposition}
\label{covering in dual for complete intersection}The complex
$\operatorname*{dual}\left(  B\left(  I\right)  \right)  $
\index{dual complex}%
contains
\index{dual}%
a
\index{trivial covering}%
trivial $c:1$
\index{covering}%
covering%
\[
\bigcup_{j=1}^{c}\operatorname*{dual}\left(  B\left(  I\right)  \right)
\cap\Delta_{j}\overset{\pi}{\rightarrow}B\left(  I\right)  ^{\vee}%
\]
with
\index{covering}%
sheets
\index{sheets}%
$\operatorname*{dual}\left(  B\left(  I\right)  \right)  \cap\Delta_{j}$. Here
the intersection of the polytope $\Delta_{j}$ with the complex
$\operatorname*{dual}\left(  B\left(  I\right)  \right)  $ is defined as
intersection of each face of $\operatorname*{dual}\left(  B\left(  I\right)
\right)  $ with $\Delta_{j}$.

If $\operatorname*{dual}\left(  F\right)  $ is a minimal face of
$\operatorname*{dual}\left(  B\left(  I\right)  \right)  $, i.e., has
dimension $n-1-d=c-1$, then $F$ has precisely $c$ vertices (indeed precisely
$c$ lattice points).
\end{proposition}

The
\index{dual}%
union
\index{dual complex}%
of the
\index{sheets}%
sheets should be related to the
\index{tropical subcomplex}%
tropical subcomplex of infinity of the mirror
\index{degeneration}%
degeneration, i.e.,%
\[
\bigcup_{j=1}^{c}\operatorname*{dual}\left(  B\left(  I\right)  \right)
\cap\Delta_{j}=BF\left(  I^{\ast}\right)  \cap\left\{  w_{t}=0\right\}
\subset S^{n}\subset\mathbb{R}^{n+1}%
\]

Any face has $c$ disjoint non empty, but possibly degenerate, preimage faces.
There is an algorithm, computing the above
\index{covering}%
covering inductively
\index{dual}%
from
\index{dual complex}%
$\operatorname*{dual}\left(  B\left(  I\right)  \right)  $ without using the
polytopes $\Delta_{j}$. It starts with associating to any face
$\operatorname*{dual}\left(  F\right)  $ of lowest dimension $n-1-d$ the set
of its $c$ vertices. Inductively for growing dimension of
$\operatorname*{dual}\left(  F\right)  $, associate to it the set of those of
its faces, which intersect each previously computed set of sheet faces at most once:

\begin{algorithm}
The following algorithm computes the above $c:1$ covering $\pi$ of $B\left(
I\right)  ^{\vee}$:

\begin{itemize}
\item If $F$ is a face of $B\left(  I\right)  $ of $\dim\left(  F\right)  =d$
and $p_{1},...,p_{c}$ are the vertices of $\operatorname*{dual}\left(
F\right)  $ then set
\[
\pi\left(  p_{j}\right)  =F^{\vee}%
\]
for $j=1,...,c$.

\item If $l>0$ and $F$ is a face of $B\left(  I\right)  $ of $\dim\left(
F\right)  =d-l$ then the faces of the covering $\pi$ over $F^{\vee}$ are those
faces $H$ of $\operatorname*{dual}\left(  F\right)  $ with

\begin{itemize}
\item $H$ intersects at most one of the elements of $\pi^{-1}\left(  Q^{\vee
}\right)  $ for every face $Q^{\vee}\subsetneqq F^{\vee}$, i.e., for all faces
$Q$ of $B\left(  I\right)  $ with $F\subsetneqq Q$, and

\item $H\notin\pi^{-1}\left(  Q^{\vee}\right)  $ for all faces $Q^{\vee
}\subsetneqq F^{\vee}$.
\end{itemize}
\end{itemize}
\end{algorithm}

\begin{figure}
[h]
\begin{center}
\includegraphics[
height=2.4561in,
width=2.7268in
]%
{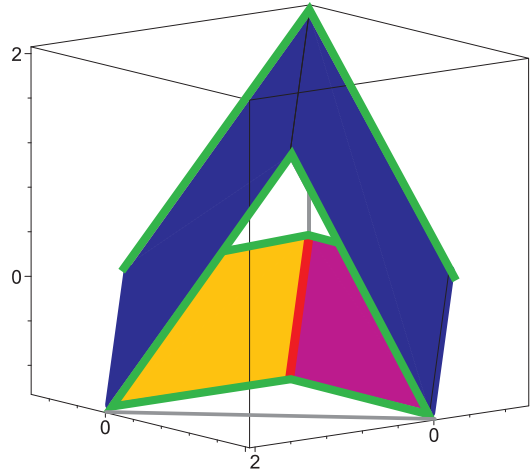}%
\caption{Covering of $B\left(  I\right)  ^{\vee}$ given by the initial terms
for the monomial degeneration of the complete intersection of two general
quadrics in $\mathbb{P}^{3}$}%
\label{Fig sheets 22}%
\end{center}
\end{figure}

\begin{example}
In the case of the degeneration of the complete intersection of two general
quadrics to the monomial ideal $\left\langle x_{1}x_{2},x_{0}x_{3}%
\right\rangle $, as defined in Example \ref{22Nabla}, the two
\index{sheets}%
sheets of the
\index{covering}%
covering inside $\operatorname*{dual}\left(  B\left(  I\right)  \right)  $ are
shown in Figure \ref{Fig sheets 22}. The sheets are formed by the
\index{lead monomial}%
initial terms of the defining equations $f_{1}=x_{1}x_{2}+tg_{1}$ and
$f_{2}=x_{0}x_{3}+tg_{2}$ at the faces of the Bergman complex.
\end{example}

%

\begin{corollary}
Above
\index{covering}%
covering allows to reconstruct the
\index{reduced standard basis}%
reduced
\index{Gr\"{o}bner basis}%
Gr\"{o}bner basis
\index{dual}%
equations
\index{dual complex}%
by clearing the denominators from%
\[
f_{j}=t\cdot\sum_{\tilde{m}\in\operatorname*{dual}\left(  B\left(  I\right)
\right)  \cap\Delta_{j}}A\tilde{m}+1
\]

\end{corollary}

\subsection{Covering of $\left(  \mu\left(  B\left(  I\right)  \right)
\right)  ^{\ast}$, construction of $I^{\circ}$ from the tropical data and
equivalence to the Batyrev-Borisov mirror\label{Sec mirror covering ci}}

We now apply a similar procedure to construct the mirror family $I^{\circ}$.
The inclusion
\index{dual}%
corresponding
\index{dual complex}%
to $\operatorname*{dual}\left(  B\left(  I\right)  \right)  =B\left(
I\right)  ^{\ast}\subset\nabla^{\ast}$ on the
\index{mirror}%
mirror
\index{mirror complex}%
side should be $\left(  \mu\left(  B\left(  I\right)  \right)  \right)
^{\ast}\subset\Delta^{\ast}$, indeed:

\begin{lemma}
\label{1mirrorcovering}For any face $F$ of $B\left(  I\right)  $ we have
$F=\sum_{j=1}^{c}\left(  \mu\left(  F\right)  \right)  ^{\ast}\cap\nabla_{j}$.

Applying the above algorithm yields a
\index{covering}%
covering
\[
\bigcup_{j=1}^{c}\left(  \mu\left(  B\left(  I\right)  \right)  \right)
^{\ast}\cap\nabla_{j}\rightarrow\left(  \mu\left(  B\left(  I\right)  \right)
\right)  ^{\vee}%
\]

\end{lemma}

The union of the
\index{sheets}%
sheets should be related to the
\index{tropical subcomplex}%
tropical subcomplex of infinity:%
\[
\bigcup_{j=1}^{c}\left(  \mu\left(  B\left(  I\right)  \right)  \right)
^{\ast}\cap\nabla_{j}=BF\left(  I\right)  \cap\left\{  w_{t}=0\right\}
\subset S^{n}\subset\mathbb{R}^{n+1}%
\]

\begin{figure}
[h]
\begin{center}
\includegraphics[
height=2.4206in,
width=2.2935in
]%
{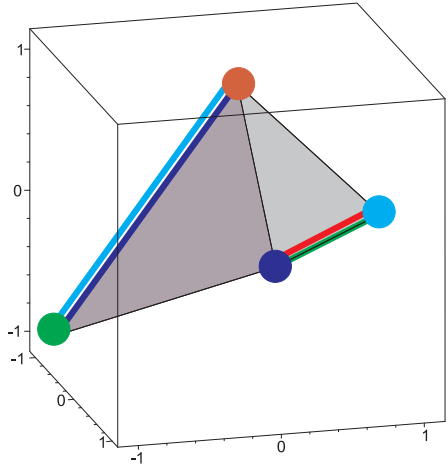}%
\caption{Covering of $\mu\left(  B\left(  I\right)  \right)  ^{\vee}$ inside
$\Delta^{\ast}$ for the monomial degeneration of the complete intersection of
two general quadrics in $\mathbb{P}^{3}$}%
\label{Fig covering mirror 22}%
\end{center}
\end{figure}

\begin{example}
In above Example \ref{22Nabla} the
\index{covering}%
faces of the covering inside the complex of mirror initial ideals are shown in
Figure \ref{Fig covering mirror 22}.
\end{example}

%

As a Corollary to Lemma \ref{1mirrorcovering} the sheets of the covering
correspond to the equations defining the complete intersection special fiber
$X_{0}\subset Y$.

\begin{corollary}
Denoting by $B_{i}=\left(  \mu\left(  B\left(  I\right)  \right)  \right)
^{\ast}\cap\nabla_{j}$ the sheets of this covering, we have%
\[
I_{0}=\left\langle
{\textstyle\prod\nolimits_{\substack{r\in\Sigma\left(  1\right)  \\\hat{r}%
\in\operatorname*{supp}\left(  B_{i}\right)  }}}
y_{r}\mid i=1,...,c\right\rangle \subset S
\]

\end{corollary}

In terms of the complex $\mu\left(  B\left(  I\right)  \right)  =\lim\left(
B\left(  I\right)  \right)  =\operatorname*{Strata}\nolimits_{\Delta}\left(
I_{0}\right)  $ we can define the ideal%
\begin{align*}
I_{0}^{\Sigma}  &  =%
{\textstyle\bigcap\nolimits_{F\in\operatorname*{Strata}\nolimits_{\Delta
}\left(  I_{0}\right)  _{d}}}
\left\langle y_{G^{\ast}}\mid G\text{ a facet of }\Delta\text{ with }F\subset
G\right\rangle \\
&  =\left\langle
{\displaystyle\prod\limits_{v\in J}}
y_{v}\mid J\subset\Sigma\left(  1\right)  \text{ with }\operatorname*{supp}%
\left(  \mu\left(  B\left(  I\right)  \right)  \right)  \subset%
{\displaystyle\bigcup\limits_{v\in J}}
F_{v}\right\rangle \subset S
\end{align*}
Passing from $I_{0}$ to $I_{0}^{\Sigma}$ is for reduced monomial ideals the
non-simplicial Cox ring analogue of saturation in the irrelevant ideal.

\begin{lemma}
The ideals $I_{0}^{\Sigma}$ and $I_{0}$ in $S$ both define the same subvariety
$X_{0}$ of $Y$.
\end{lemma}

Denote the map%
\[
0\rightarrow M^{\circ}\overset{A^{\circ}}{\rightarrow}\mathbb{Z}^{\Sigma
^{\ast}\left(  1\right)  }\rightarrow A_{n-1}\left(  X\left(  \Sigma^{\circ
}\right)  \right)  \rightarrow0
\]
by $A^{\circ}$. We consider the ideal generated by equations corresponding to
the
\index{sheets}%
sheets of the
\index{covering}%
covering given in Proposition \ref{1mirrorcovering}. The
\index{degeneration}%
degeneration $\mathfrak{X}^{\circ}\subset\mathbb{P}\left(  \nabla\right)
\times\operatorname*{Spec}\mathbb{C}\left[  t\right]  $ defined by%
\[
I^{\circ}=\left\langle t\cdot\sum_{\delta\in\left(  \mu\left(  B\left(
I\right)  \right)  \right)  ^{\ast}\cap\nabla_{j}}c_{\delta}\cdot A^{\circ
}\left(  \delta\right)  +1\mid j=1,...,c\right\rangle
\]
with generic coefficients $c_{\delta}$ coincides with the
\index{degeneration}%
degeneration associated to the Batyrev-Borisov mirror, i.e., after clearing
the denominators, the generators coincide with the reduced Gr\"{o}bner basis
of the ideal of the degeneration associated to the Batyrev-Borisov mirror.

\begin{theorem}
The
\index{mirror construction}%
mirror obtained from the tropical construction coincides with the
Batyrev-Borisov mirror.
\end{theorem}

\begin{corollary}
Denoting by $B_{i}^{\circ}=\operatorname*{dual}\left(  B\left(  I\right)
\right)  \cap\Delta_{j}$ the sheets of the covering given in Proposition
\ref{covering in dual for complete intersection}, the special fiber of
$\mathfrak{X}^{\circ}$ is given by%
\[
I_{0}^{\circ}=\left\langle
{\textstyle\prod\nolimits_{\substack{r\in\Sigma\left(  1\right)  \\\hat{r}%
\in\operatorname*{supp}\left(  B_{i}^{\circ}\right)  }}}
z_{r}\mid i=1,...,c\right\rangle \subset S^{\circ}=\mathbb{C}\left[  z_{r}\mid
r\in\Sigma^{\circ}\left(  1\right)  \right]
\]

\end{corollary}

In terms of the complex $B\left(  I\right)  =\operatorname*{Strata}%
\nolimits_{\nabla}\left(  I_{0}^{\circ}\right)  $ we have the ideal%
\begin{align*}
\left(  I_{0}^{\circ}\right)  ^{\Sigma^{\circ}}  &  =%
{\textstyle\bigcap\nolimits_{F\in B\left(  I\right)  _{d}}}
\left\langle z_{G^{\ast}}\mid G\text{ a facet of }\nabla\text{ with }F\subset
G\right\rangle \\
&  =\left\langle
{\displaystyle\prod\limits_{v\in J}}
z_{v}\mid J\subset\Sigma^{\circ}\left(  1\right)  \text{ with }%
\operatorname*{supp}\left(  B\left(  I\right)  \right)  \subset%
{\displaystyle\bigcup\limits_{v\in J}}
F_{v}\right\rangle \subset S^{\circ}%
\end{align*}

\begin{lemma}
The ideals $\left(  I_{0}^{\circ}\right)  ^{\Sigma^{\circ}}$ and $I_{0}%
^{\circ}$ in $S^{\circ}$ both define the same subvariety $X_{0}^{\circ}$ of
$Y^{\circ}$.
\end{lemma}

Indeed from the point of view of saturation in the sense of removing the
irrelevant components, we should associate to the special fiber $X_{0}^{\circ
}$ of $\mathfrak{X}^{\circ}$ the ideal $\left(  I_{0}^{\circ}\right)
^{\Sigma}$, and to the degeneration $\mathfrak{X}^{\circ}$ the ideal%
\[
\left\langle t\cdot\sum_{\substack{\delta\in\left(  \mu\left(  B\left(
I\right)  \right)  \right)  ^{\ast}\cap\nabla_{j}\\A^{\circ}\delta\cdot
m_{0}\in S^{\circ}}}c_{\delta}\cdot A^{\circ}\left(  \delta\right)  \cdot
m_{0}+m_{0}\mid m_{0}\text{ a minimal generator of }\left(  I_{0}^{\circ
}\right)  ^{\Sigma^{\circ}}\right\rangle
\]
in $S^{\circ}\otimes\mathbb{C}\left[  t\right]  $ with generic coefficients
$c_{\delta}$. The same holds true of course for $\mathfrak{X}$.

Note that passing to the saturated discription does not change the objects
involved in the tropical mirror construction, as the special fiber complex and
the set of first order deformations does not change.

\begin{figure}
[h]
\begin{center}
\includegraphics[
height=4.5057in,
width=4.4218in
]%
{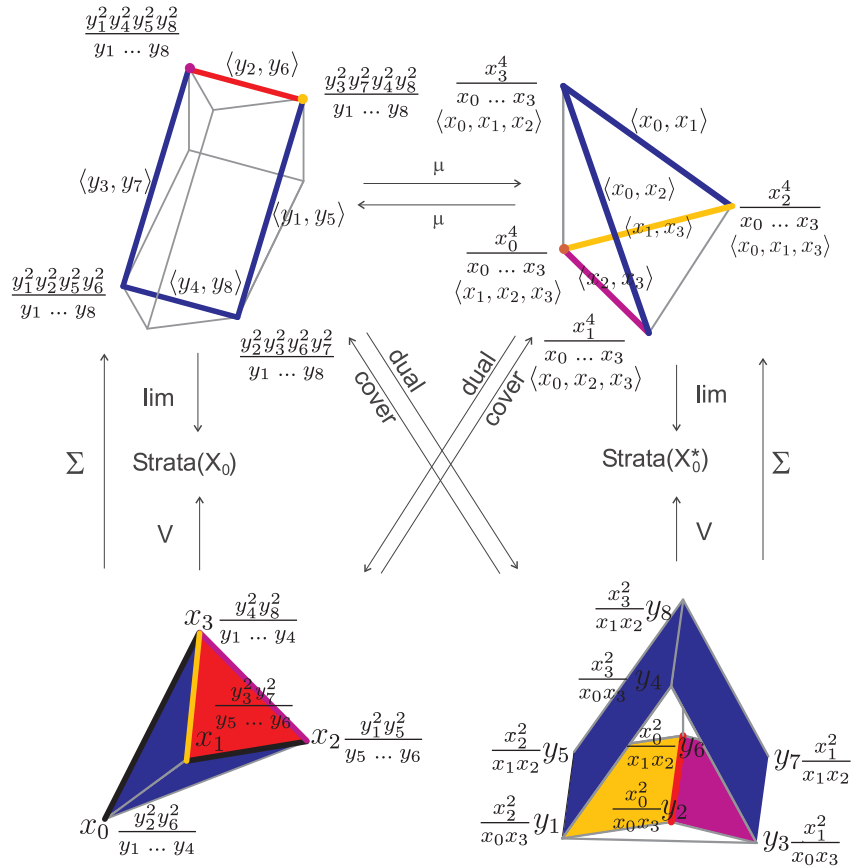}%
\caption{Summary of the tropical mirror construction for above monomial
degeneration of the complete intersection of two general quadrics in
$\mathbb{P}^{3}$}%
\label{Fig summary 22}%
\end{center}
\end{figure}

\begin{example}
Figure \ref{Fig summary 22} gives a
\index{Bergman subcomplex}%
summary
\index{dual}%
of
\index{dual complex}%
the
\index{mirror}%
tropical
\index{mirror complex}%
mirror
\index{tropical subcomplex}%
construction for above monomial degeneration of the complete intersection of
two general quadrics in $\mathbb{P}^{3}$.
\end{example}

%

\subsection{Examples\label{Sec Examples ci}}

In the following we give the explicit computations for some simple examples.
The text is computer generated from the output given by the Maple package
\textsf{tropicalmirror}, which implements the tropical mirror construction.
See also Section \ref{Sec tropicalmirror implementation} for a short
description of the \textsf{tropicalmirror} package.

\subsubsection{The elliptic curve given as the complete intersection of two
generic quadrics in $\mathbb{P}^{3}$%
\index{elliptic curve!complete intersection}%
}

\paragraph{Setup}

Let $Y=\mathbb{P}^{
3
}=X\left(  \Sigma\right)  $, $\Sigma
=\operatorname*{Fan}\left(  P\right)  =NF\left(  \Delta\right)  \subset
N_{\mathbb{R}}$ with the Fano polytope $P=\Delta^{\ast}$ given by%
\[
\Delta=\operatorname*{convexhull}\left(  
\begin{tabular}
[c]{lll}
$\left(3,-1,-1\right)$ &$\left(-1,3,-1\right)$ &$\left(-1,-1,3\right)$\\
$\left(-1,-1,-1\right)$ & &
\end{tabular}
\right)  \subset
M_{\mathbb{R}}%
\]
and let%
\[
S=\mathbb{C}
[x_0, x_1, x_2, x_3]
\]
be the Cox ring of $Y$ with the variables%
\[
\begin{tabular}
[c]{llll}
$x_{1} = x_{\left(1,0,0\right)}$ &$x_{2} = x_{\left(0,1,0\right)}$ &$x_{3} = x_{\left(0,0,1\right)}$ &$x_{0} = x_{\left(-1,-1,-1\right)}$
\end{tabular}
\]
associated to the rays of $\Sigma$.

Consider the degeneration $\mathfrak{X}\subset Y\times\operatorname*{Spec}%
\mathbb{C}\left[  t\right]  $ of 
complete intersection
elliptic curves
 of type $\left(
2,2
\right)  $ with monomial special fiber
\[
I_{0}=\left\langle 
\begin{tabular}
[c]{ll}
$x_{1}\,x_{2}$ & $x_{0}\,x_{3}$
\end{tabular}
\right\rangle
\]
The degeneration is given by the ideal $I\subset S\otimes\mathbb{C}\left[
t\right]  $ with $I_{0}$-reduced generators%
\[
\left\{  
\begin{tabular}
[c]{l}
$x_{1}\,x_{2}+t (s_{1}\,x_{1}^{2}+s_{2}\,x_{1}\,x_{3}+s_{3}\,x_{1}\,x_{0}+s_{4}\,x_{2}^{2}+s_{5}\,x_{2}\,x_{3}+s_{6}\,x_{2}\,x_{0}+s_{7}\,x_{3}^{2}+s_{8}\,x_{0}^{2} ),\medskip$
\\
$x_{0}\,x_{3}+t (s_{9}\,x_{1}^{2}+s_{10}\,x_{1}\,x_{3}+s_{11}\,x_{1}\,x_{0}+s_{12}\,x_{2}^{2}+s_{13}\,x_{2}\,x_{3}+s_{14}\,x_{2}\,x_{0}+s_{15}\,x_{3}^{2}+s_{16}\,x_{0}^{2} )$
\\
\end{tabular}
\right\}
\]

\paragraph{Special fiber Gr\"{o}bner cone}

The space of first order deformations of $\mathfrak{X}$ has dimension
$
16
$ and the deformations represented by the Cox Laurent monomials

\begin{center}
\begin{longtable}
[c]{lllllllllllll}
${\frac {x_{0}^{2}}{x_{1}\,x_{2}}}$ & ${\frac {x_{3}^{2}}{x_{1}\,x_{2}}}$ & ${\frac {x_{1}^{2}}{x_{0}\,x_{3}}}$ & ${\frac {x_{2}^{2}}{x_{0}\,x_{3}}}$ & ${\frac {x_{3}}{x_{1}}}$ & ${\frac {x_{2}}{x_{1}}}$ & ${\frac {x_{0}}{x_{3}}}$ & ${\frac {x_{1}}{x_{3}}}$ & ${\frac {x_{2}}{x_{3}}}$ & ${\frac {x_{1}}{x_{0}}}$ & ${\frac {x_{2}}{x_{0}}}$ & ${\frac {x_{3}}{x_{0}}}$ & ${\frac {x_{0}}{x_{1}}}$\\
${\frac {x_{0}}{x_{2}}}$ & ${\frac {x_{3}}{x_{2}}}$ & ${\frac {x_{1}}{x_{2}}}$ &  & & & & & & & & &
\end{longtable}
\end{center}

\noindent form a torus invariant basis. These deformations give linear
inequalities defining the special fiber Gr\"{o}bner cone, i.e.,%
\[
C_{I_{0}}\left(  I\right)  =\left\{  \left(  w_{t},w_{y}\right)  \in
\mathbb{R}\oplus N_{\mathbb{R}}\mid\left\langle w,A^{-1}\left(  m\right)
\right\rangle \geq-w_{t}\forall m\right\}
\]
for above monomials $m$ and the presentation matrix
\[
A=
\left [\begin {array}{ccc} 1&0&0\\\noalign{\medskip}0&1&0\\\noalign{\medskip}0&0&1\\\noalign{\medskip}-1&-1&-1\end {array}\right ]
\]
of $A_{
2
}\left(  Y\right)  $.

The vertices of special fiber polytope
\[
\nabla=C_{I_{0}}\left(  I\right)  \cap\left\{  w_{t}=1\right\}
\]
and the number of faces, each vertex is contained in, are given by the table

\begin{center}
\begin{longtable}
[c]{rl}
$3$ & \begin{tabular}[c]{llll}$\left(0,0,1\right)$ &$\left(-1,-1,-1\right)$ &$\left(0,1,0\right)$ &$\left(1,0,0\right)$\end{tabular}\medskip\\
$4$ & \begin{tabular}[c]{llll}$\left(1,0,1\right)$ &$\left(0,1,1\right)$ &$\left(0,-1,-1\right)$ &$\left(-1,0,-1\right)$\end{tabular}
\end{longtable}
\end{center}

\noindent The number of faces of $\nabla$ and their $F$-vectors are

\begin{center}
\begin{longtable}
[c]{rrll}
Dimension &Number of faces &F-vector &\\
$-1$ & $1$ & $\left(0,0,0,0,0\right)$ &\\
$0$ & $8$ & $\left(1,1,0,0,0\right)$ &point\\
$1$ & $14$ & $\left(1,2,1,0,0\right)$ &edge\\
$2$ & $4$ & $\left(1,4,4,1,0\right)$ &quadrangle\\
$2$ & $4$ & $\left(1,3,3,1,0\right)$ &triangle\\
$3$ & $1$ & $\left(1,8,14,8,1\right)$ &
\end{longtable}
\end{center}

The dual $P^{\circ}=\nabla^{\ast}$ of $\nabla$ is a Fano polytope with vertices

\begin{center}
\begin{longtable}
[c]{llllll}
$\left(-1,-1,0\right)$ &$\left(-1,1,0\right)$ &$\left(1,-1,0\right)$ &$\left(-1,-1,2\right)$ &$\left(0,2,-1\right)$ &$\left(0,0,-1\right)$\\
$\left(2,0,-1\right)$ &$\left(0,0,1\right)$ & & & &
\end{longtable}
\end{center}

\noindent and the $F$-vectors of the faces of $P^{\circ}$ are

\begin{center}
\begin{longtable}
[c]{rrll}
Dimension &Number of faces &F-vector &\\
$-1$ & $1$ & $\left(0,0,0,0,0\right)$ &\\
$0$ & $8$ & $\left(1,1,0,0,0\right)$ &point\\
$1$ & $14$ & $\left(1,2,1,0,0\right)$ &edge\\
$2$ & $4$ & $\left(1,3,3,1,0\right)$ &triangle\\
$2$ & $4$ & $\left(1,4,4,1,0\right)$ &quadrangle\\
$3$ & $1$ & $\left(1,8,14,8,1\right)$ &
\end{longtable}
\end{center}

The polytopes $\nabla$ and $P^{\circ}$ have $0$ as the unique interior lattice
point, and the Fano polytope $P^{\circ}$ defines the embedding toric Fano
variety $Y^{\circ}=X\left(  \operatorname*{Fan}\left(  P^{\circ}\right)
\right)  $ of the fibers of the mirror degeneration. The dimension of the
automorphism group of $Y^{\circ}$ is
\[
\dim\left(  \operatorname*{Aut}\left(  Y^{\circ}\right)  \right)
=
3
\]
Let%
\[
S^{\circ}=\mathbb{C}
[y_{1},y_{2},...,y_{8}]
\]
be the Cox ring of $Y^{\circ}$ with variables

\begin{center}
\begin{longtable}
[c]{lll}
$y_{1} = y_{\left(-1,-1,0\right)}={\frac {x_{0}^{2}}{x_{1}\,x_{2}}}$ &$y_{2} = y_{\left(-1,1,0\right)}={\frac {x_{2}}{x_{1}}}$ &$y_{3} = y_{\left(1,-1,0\right)}={\frac {x_{1}}{x_{2}}}$\\
$y_{4} = y_{\left(-1,-1,2\right)}={\frac {x_{3}^{2}}{x_{1}\,x_{2}}}$ &$y_{5} = y_{\left(0,2,-1\right)}={\frac {x_{2}^{2}}{x_{0}\,x_{3}}}$ &$y_{6} = y_{\left(0,0,-1\right)}={\frac {x_{0}}{x_{3}}}$\\
$y_{7} = y_{\left(2,0,-1\right)}={\frac {x_{1}^{2}}{x_{0}\,x_{3}}}$ &$y_{8} = y_{\left(0,0,1\right)}={\frac {x_{3}}{x_{0}}}$ &
\end{longtable}
\end{center}

\paragraph{Bergman subcomplex}

Intersecting the tropical variety of $I$ with the special fiber Gr\"{o}bner
cone $C_{I_{0}}\left(  I\right)  $ we obtain the special fiber Bergman
subcomplex
\[
B\left(  I\right)  =C_{I_{0}}\left(  I\right)  \cap BF\left(  I\right)
\cap\left\{  w_{t}=1\right\}
\]
The following table chooses an indexing of the vertices of $\nabla$ involved
in $B\left(  I\right)  $ and gives for each vertex the numbers $n_{\nabla}$
and $n_{B}$ of faces of $\nabla$ and faces of $B\left(  I\right)  $ it is
contained in

\begin{center}
\begin{longtable}
[c]{rrl}
$n_{\nabla}$ &$n_{B(I)}$ & \medskip\\
$4$ & $2$ & \begin{tabular}[c]{lll}$1 = \left(1,0,1\right)$ &$2 = \left(0,1,1\right)$ &$3 = \left(0,-1,-1\right)$\\$4 = \left(-1,0,-1\right)$ & &\end{tabular}
\end{longtable}
\end{center}

\noindent With this indexing the Bergman subcomplex $B\left(  I\right)  $ of
$\operatorname*{Poset}\left(  \nabla\right)  $ associated to the degeneration
$\mathfrak{X}$ is

\begin{center}
\begin{longtable}
[c]{l}
$[],\medskip$
\\
$[[1],[2],[3],[4]],\medskip$
\\
$[[3, 4],[2, 4],[1, 3],[1, 2]],\medskip$
\\
$[],\medskip$
\\
$[]$
\\
\end{longtable}
\end{center}

$B\left(  I\right)  $ has the $F$-vector%
\[
\begin{tabular}
[c]{rrrrrr}
Dimension &$-1$ & $0$ & $1$ & $2$ & $3$\\
Number of faces &$0$ & $4$ & $4$ & $0$ & $0$
\end{tabular}
\]
and the $F$-vectors of its faces are

\begin{center}
\begin{longtable}
[c]{rrll}
Dimension &Number of faces &F-vector &\\
$0$ & $4$ & $\left(1,1,0,0,0\right)$ &point\\
$1$ & $4$ & $\left(1,2,1,0,0\right)$ &edge
\end{longtable}
\end{center}

\paragraph{Dual complex}

The dual complex $\operatorname*{dual}\left(  B\left(  I\right)  \right)
=\left(  B\left(  I\right)  \right)  ^{\ast}$ of deformations associated to
$B\left(  I\right)  $ via initial ideals is given by

\begin{center}
\begin{longtable}
[c]{l}
$[],\medskip$
\\
$[],\medskip$
\\
$[[3, 4]^{*}=\left\langle {\frac {x_{3}^{2}}{x_{1}\,x_{2}}},{\frac {x_{3}}{x_{0}}}\right\rangle,[2, 4]^{*}=\left\langle {\frac {x_{1}^{2}}{x_{0}\,x_{3}}},{\frac {x_{1}}{x_{2}}}\right\rangle,[1, 3]^{*}=\left\langle {\frac {x_{2}^{2}}{x_{0}\,x_{3}}},{\frac {x_{2}}{x_{1}}}\right\rangle,[1, 2]^{*}=\left\langle {\frac {x_{0}^{2}}{x_{1}\,x_{2}}},{\frac {x_{0}}{x_{3}}}\right\rangle],\medskip$
\\
$[[1]^{*}=\left\langle {\frac {x_{0}^{2}}{x_{1}\,x_{2}}},{\frac {x_{2}^{2}}{x_{0}\,x_{3}}},{\frac {x_{0}}{x_{3}}},{\frac {x_{2}}{x_{1}}}\right\rangle,[2]^{*}=\left\langle {\frac {x_{0}^{2}}{x_{1}\,x_{2}}},{\frac {x_{1}^{2}}{x_{0}\,x_{3}}},{\frac {x_{0}}{x_{3}}},{\frac {x_{1}}{x_{2}}}\right\rangle,[3]^{*}=\left\langle {\frac {x_{3}^{2}}{x_{1}\,x_{2}}},{\frac {x_{2}^{2}}{x_{0}\,x_{3}}},{\frac {x_{3}}{x_{0}}},{\frac {x_{2}}{x_{1}}}\right\rangle,$
\\
$[4]^{*}=\left\langle {\frac {x_{3}^{2}}{x_{1}\,x_{2}}},{\frac {x_{1}^{2}}{x_{0}\,x_{3}}},{\frac {x_{3}}{x_{0}}},{\frac {x_{1}}{x_{2}}}\right\rangle],\medskip$
\\
$[]$
\\
\end{longtable}
\end{center}

\noindent when writing the vertices of the faces as deformations of $X_{0}$.
Note that the $T$-invariant basis of deformations associated to a face is
given by all lattice points of the corresponding polytope in $M_{\mathbb{R}}$.

When numbering the vertices of the faces of $\operatorname*{dual}\left(
B\left(  I\right)  \right)  $ by the Cox variables of the mirror toric Fano
variety $Y^{\circ}$ the complex $\operatorname*{dual}\left(  B\left(
I\right)  \right)  $ is

\begin{center}
\begin{longtable}
[c]{l}
$[],\medskip$
\\
$[],\medskip$
\\
$[[3, 4]^{*}=\left\langle y_{4},y_{8}\right\rangle,[2, 4]^{*}=\left\langle y_{7},y_{3}\right\rangle,[1, 3]^{*}=\left\langle y_{5},y_{2}\right\rangle,[1, 2]^{*}=\left\langle y_{1},y_{6}\right\rangle],\medskip$
\\
$[[1]^{*}=\left\langle y_{1},y_{5},y_{6},y_{2}\right\rangle,[2]^{*}=\left\langle y_{1},y_{7},y_{6},y_{3}\right\rangle,[3]^{*}=\left\langle y_{4},y_{5},y_{8},y_{2}\right\rangle,$
\\
$[4]^{*}=\left\langle y_{4},y_{7},y_{8},y_{3}\right\rangle],\medskip$
\\
$[]$
\\
\end{longtable}
\end{center}

\noindent The dual complex has the $F$-vector%
\[
\begin{tabular}
[c]{rrrrrr}
Dimension &$-1$ & $0$ & $1$ & $2$ & $3$\\
Number of faces &$0$ & $0$ & $4$ & $4$ & $0$
\end{tabular}
\]
and the $F$-vectors of the faces of $\operatorname*{dual}\left(  B\left(
I\right)  \right)  $ are

\begin{center}
\begin{longtable}
[c]{rrll}
Dimension &Number of faces &F-vector &\\
$1$ & $4$ & $\left(1,2,1,0,0\right)$ &edge\\
$2$ & $4$ & $\left(1,4,4,1,0\right)$ &quadrangle
\end{longtable}
\end{center}

Recall that in this example the toric variety $Y$ is projective space. The
number of lattice points of the support of $\operatorname*{dual}\left(
B\left(  I\right)  \right)  $ relates to the dimension $h^{1,
0
}%
\left(  X\right)  $ of the complex moduli space of the generic fiber $X$ of
$\mathfrak{X}$ and to the dimension $h^{1,1}\left(  \bar{X}^{\circ}\right)  $
of the K\"{a}hler moduli space of the $MPCR$-blowup $\bar{X}^{\circ}$ of the
generic fiber $X^{\circ}$ of the mirror degeneration
\begin{align*}
\left\vert \operatorname*{supp}\left(  \operatorname*{dual}\left(  B\left(
I\right)  \right)  \right)  \cap M\right\vert  &
=
16
=
15
+
1
=\dim\left(  \operatorname*{Aut}%
\left(  Y\right)  \right)  +h^{1,
0
}\left(  X\right) \\
&  =
12
+
3
+
1
\\
&  =\left\vert \operatorname*{Roots}\left(  Y\right)  \right\vert +\dim\left(
T_{Y}\right)  +h^{1,1}\left(  \bar{X}^{\circ}\right)
\end{align*}
There are%
\[
h^{1,
0
}\left(  X\right)  +\dim\left(  T_{Y^{\circ}}\right)
=
1
+
3
\]
non-trivial toric polynomial deformations of $X_{0}$

\begin{center}
\begin{longtable}
[c]{llll}
${\frac {x_{3}^{2}}{x_{1}\,x_{2}}}$ & ${\frac {x_{1}^{2}}{x_{0}\,x_{3}}}$ & ${\frac {x_{2}^{2}}{x_{0}\,x_{3}}}$ & ${\frac {x_{0}^{2}}{x_{1}\,x_{2}}}$
\end{longtable}
\end{center}

\noindent They correspond to the toric divisors

\begin{center}
\begin{longtable}
[c]{llll}
$D_{\left(-1,-1,2\right)}$ & $D_{\left(2,0,-1\right)}$ & $D_{\left(0,2,-1\right)}$ & $D_{\left(-1,-1,0\right)}$
\end{longtable}
\end{center}

\noindent on a MPCR-blowup of $Y^{\circ}$ inducing $
1
$ non-zero
toric divisor classes on the mirror $X^{\circ}$. The following $
12
$
toric divisors of $Y$ induce the trivial divisor class on $X^{\circ}$

\begin{center}
\begin{longtable}
[c]{lllll}
$D_{\left(0,0,1\right)}$ & $D_{\left(1,-1,0\right)}$ & $D_{\left(-1,1,0\right)}$ & $D_{\left(0,0,-1\right)}$ & $D_{\left(0,1,-1\right)}$\\
$D_{\left(-1,0,0\right)}$ & $D_{\left(1,0,-1\right)}$ & $D_{\left(0,-1,0\right)}$ & $D_{\left(0,1,0\right)}$ & $D_{\left(-1,0,1\right)}$\\
$D_{\left(1,0,0\right)}$ & $D_{\left(0,-1,1\right)}$ &  & &
\end{longtable}
\end{center}

\paragraph{Mirror special fiber}

The ideal $I_{0}^{\circ}$ of the monomial special fiber $X_{0}^{\circ}$ of the
mirror degeneration $\mathfrak{X}^{\circ}$ is generated by the following set
of monomials in $S^{\circ}$%
\[
\left\{  
\begin{tabular}
[c]{lllllll}
$y_{1}\,y_{2}\,y_{3}\,y_{8}$ & $y_{1}\,y_{2}\,y_{4}\,y_{7}$ & $y_{1}\,y_{5}\,y_{7}\,y_{8}$ & $y_{3}\,y_{5}\,y_{6}\,y_{8}$ & $y_{2}\,y_{3}\,y_{4}\,y_{6}$ & $y_{4}\,y_{5}\,y_{6}\,y_{7}$ & $y_{1}\,y_{3}\,y_{5}\,y_{8}$\\
$y_{1}\,y_{4}\,y_{5}\,y_{7}$ & $y_{3}\,y_{4}\,y_{5}\,y_{6}$ & $y_{2}\,y_{3}\,y_{6}\,y_{8}$ & $y_{2}\,y_{6}\,y_{7}\,y_{8}$ & $y_{1}\,y_{3}\,y_{4}\,y_{5}$ & $y_{2}\,y_{4}\,y_{6}\,y_{7}$ & $y_{1}\,y_{2}\,y_{7}\,y_{8}$\\
$y_{1}\,y_{2}\,y_{3}\,y_{4}$ & $y_{5}\,y_{6}\,y_{7}\,y_{8}$ &  & & & &
\end{tabular}
\right\}
\]
Indeed already the ideal%
\[
J_{0}^{\circ}=\left\langle 
\begin{tabular}
[c]{ll}
$y_{5}\,y_{6}\,y_{7}\,y_{8}$ & $y_{1}\,y_{2}\,y_{3}\,y_{4}$
\end{tabular}
\right\rangle
\]
defines the same subvariety of the toric variety $Y^{\circ}$, and
$J_{0}^{\circ\Sigma}=I_{0}^{\circ}$. Recall that passing from $J_{0}^{\circ}$
to $J_{0}^{\circ\Sigma}$ is the non-simplicial toric analogue of saturation.

The complex $B\left(  I\right)  ^{\ast}$ labeled by the variables of the Cox
ring $S^{\circ}$ of $Y^{\ast}$, as written in the last section, is the complex
$SP\left(  I_{0}^{\circ}\right)  $ of prime ideals of the toric strata of the
special fiber $X_{0}^{\circ}$ of the mirror degeneration $\mathfrak{X}^{\circ
}$, i.e., the primary decomposition of $I_{0}^{\circ}$ is%
\[
I_{0}^{\circ}=
\begin{tabular}
[c]{l}
$\left\langle y_{5},y_{2}\right\rangle \cap\left\langle y_{7},y_{3}\right\rangle \cap\left\langle y_{4},y_{8}\right\rangle \cap\left\langle y_{1},y_{6}\right\rangle$
\\
\end{tabular}
\]
Each facet $F\in B\left(  I\right)  $ corresponds to one of these ideals and this
ideal is generated by the products of facets of $\nabla$ containing $F$.

\paragraph{Covering structure in the deformation complex of the
degeneration $\mathfrak{X}$}

According to the local reduced Gr\"{o}bner basis each face of the complex of
deformations $\operatorname*{dual}\left(  B\left(  I\right)  \right)  $
decomposes into 
$2$
polytopes forming a $
2
:1$
trivial
covering of $B\left(  I\right)  $

\begin{center}
\begin{longtable}
[c]{l}
$[],\medskip$
\\
$[],\medskip$
\\
$[[\left\langle y_{4}\right\rangle,\left\langle y_{8}\right\rangle]\mapsto\left\langle y_{4},y_{8}\right\rangle^{\ast\vee}=[3, 4]^{\vee},[\left\langle y_{3}\right\rangle,\left\langle y_{7}\right\rangle]\mapsto\left\langle y_{7},y_{3}\right\rangle^{\ast\vee}=[2, 4]^{\vee},$
\\
$[\left\langle y_{2}\right\rangle,\left\langle y_{5}\right\rangle]\mapsto\left\langle y_{5},y_{2}\right\rangle^{\ast\vee}=[1, 3]^{\vee},[\left\langle y_{1}\right\rangle,\left\langle y_{6}\right\rangle]\mapsto\left\langle y_{1},y_{6}\right\rangle^{\ast\vee}=[1, 2]^{\vee}],\medskip$
\\
$[[\left\langle y_{1},y_{2}\right\rangle,\left\langle y_{5},y_{6}\right\rangle]\mapsto\left\langle y_{1},y_{5},y_{6},y_{2}\right\rangle^{\ast\vee}=[1]^{\vee},$
\\
$[\left\langle y_{1},y_{3}\right\rangle,\left\langle y_{7},y_{6}\right\rangle]\mapsto\left\langle y_{1},y_{7},y_{6},y_{3}\right\rangle^{\ast\vee}=[2]^{\vee},$
\\
$[\left\langle y_{4},y_{2}\right\rangle,\left\langle y_{5},y_{8}\right\rangle]\mapsto\left\langle y_{4},y_{5},y_{8},y_{2}\right\rangle^{\ast\vee}=[3]^{\vee},$
\\
$[\left\langle y_{4},y_{3}\right\rangle,\left\langle y_{7},y_{8}\right\rangle]\mapsto\left\langle y_{4},y_{7},y_{8},y_{3}\right\rangle^{\ast\vee}=[4]^{\vee}],\medskip$
\\
$[]$
\\
\end{longtable}
\end{center}

\noindent Here the faces $F\in B\left(  I\right)  $ are specified both via the
vertices of $F^{\ast}$ labeled by the variables of $S^{\circ}$ and by the
numbering of the vertices of $B\left(  I\right)  $ chosen above.

\noindent This covering has 
2 sheets forming the complexes
 
\[
\begin{tabular}
[c]{l}
$[],\medskip$
\\
$[],\medskip$
\\
$[\left\langle y_{8}\right\rangle,\left\langle y_{7}\right\rangle,\left\langle y_{5}\right\rangle,\left\langle y_{6}\right\rangle],\medskip$
\\
$[\left\langle y_{5},y_{6}\right\rangle,\left\langle y_{7},y_{6}\right\rangle,\left\langle y_{5},y_{8}\right\rangle,\left\langle y_{7},y_{8}\right\rangle],\medskip$
\\
$[]$
\\
\end{tabular}
\]
\[
\begin{tabular}
[c]{l}
$[],\medskip$
\\
$[],\medskip$
\\
$[\left\langle y_{4}\right\rangle,\left\langle y_{3}\right\rangle,\left\langle y_{2}\right\rangle,\left\langle y_{1}\right\rangle],\medskip$
\\
$[\left\langle y_{1},y_{2}\right\rangle,\left\langle y_{1},y_{3}\right\rangle,\left\langle y_{4},y_{2}\right\rangle,\left\langle y_{4},y_{3}\right\rangle],\medskip$
\\
$[]$
\\
\end{tabular}
\]
with $F$-vector

\begin{center}
\begin{longtable}
[c]{rrll}
Dimension &Number of faces &F-vector &\\
$0$ & $4$ & $\left(1,1,0,0,0\right)$ &point\\
$1$ & $4$ & $\left(1,2,1,0,0\right)$ &edge
\end{longtable}
\end{center}

\noindent Writing the vertices of the faces as deformations the covering is
given by

\begin{center}
\begin{longtable}
[c]{l}
$[],\medskip$
\\
$[],\medskip$
\\
$[[\left\langle {\frac {x_{3}^{2}}{x_{1}\,x_{2}}}\right\rangle,\left\langle {\frac {x_{3}}{x_{0}}}\right\rangle]\mapsto[3, 4]^{\vee},[\left\langle {\frac {x_{1}}{x_{2}}}\right\rangle,\left\langle {\frac {x_{1}^{2}}{x_{0}\,x_{3}}}\right\rangle]\mapsto[2, 4]^{\vee},$
\\
$[\left\langle {\frac {x_{2}}{x_{1}}}\right\rangle,\left\langle {\frac {x_{2}^{2}}{x_{0}\,x_{3}}}\right\rangle]\mapsto[1, 3]^{\vee},[\left\langle {\frac {x_{0}^{2}}{x_{1}\,x_{2}}}\right\rangle,\left\langle {\frac {x_{0}}{x_{3}}}\right\rangle]\mapsto[1, 2]^{\vee}],\medskip$
\\
$[[\left\langle {\frac {x_{0}^{2}}{x_{1}\,x_{2}}},{\frac {x_{2}}{x_{1}}}\right\rangle,\left\langle {\frac {x_{2}^{2}}{x_{0}\,x_{3}}},{\frac {x_{0}}{x_{3}}}\right\rangle]\mapsto[1]^{\vee},$
\\
$[\left\langle {\frac {x_{0}^{2}}{x_{1}\,x_{2}}},{\frac {x_{1}}{x_{2}}}\right\rangle,\left\langle {\frac {x_{1}^{2}}{x_{0}\,x_{3}}},{\frac {x_{0}}{x_{3}}}\right\rangle]\mapsto[2]^{\vee},$
\\
$[\left\langle {\frac {x_{3}^{2}}{x_{1}\,x_{2}}},{\frac {x_{2}}{x_{1}}}\right\rangle,\left\langle {\frac {x_{2}^{2}}{x_{0}\,x_{3}}},{\frac {x_{3}}{x_{0}}}\right\rangle]\mapsto[3]^{\vee},$
\\
$[\left\langle {\frac {x_{3}^{2}}{x_{1}\,x_{2}}},{\frac {x_{1}}{x_{2}}}\right\rangle,\left\langle {\frac {x_{1}^{2}}{x_{0}\,x_{3}}},{\frac {x_{3}}{x_{0}}}\right\rangle]\mapsto[4]^{\vee}],\medskip$
\\
$[]$
\\
\end{longtable}
\end{center}

\noindent with 
2 sheets forming the complexes

\begin{center}
\[
\begin{tabular}
[c]{l}
$[],\medskip$
\\
$[],\medskip$
\\
$[\left\langle {\frac {x_{3}}{x_{0}}}\right\rangle,\left\langle {\frac {x_{1}^{2}}{x_{0}\,x_{3}}}\right\rangle,\left\langle {\frac {x_{2}^{2}}{x_{0}\,x_{3}}}\right\rangle,\left\langle {\frac {x_{0}}{x_{3}}}\right\rangle],\medskip$
\\
$[\left\langle {\frac {x_{2}^{2}}{x_{0}\,x_{3}}},{\frac {x_{0}}{x_{3}}}\right\rangle,\left\langle {\frac {x_{1}^{2}}{x_{0}\,x_{3}}},{\frac {x_{0}}{x_{3}}}\right\rangle,\left\langle {\frac {x_{2}^{2}}{x_{0}\,x_{3}}},{\frac {x_{3}}{x_{0}}}\right\rangle,\left\langle {\frac {x_{1}^{2}}{x_{0}\,x_{3}}},{\frac {x_{3}}{x_{0}}}\right\rangle],\medskip$
\\
$[]$
\\
\end{tabular}
\]
\[
\begin{tabular}
[c]{l}
$[],\medskip$
\\
$[],\medskip$
\\
$[\left\langle {\frac {x_{3}^{2}}{x_{1}\,x_{2}}}\right\rangle,\left\langle {\frac {x_{1}}{x_{2}}}\right\rangle,\left\langle {\frac {x_{2}}{x_{1}}}\right\rangle,\left\langle {\frac {x_{0}^{2}}{x_{1}\,x_{2}}}\right\rangle],\medskip$
\\
$[\left\langle {\frac {x_{0}^{2}}{x_{1}\,x_{2}}},{\frac {x_{2}}{x_{1}}}\right\rangle,\left\langle {\frac {x_{0}^{2}}{x_{1}\,x_{2}}},{\frac {x_{1}}{x_{2}}}\right\rangle,\left\langle {\frac {x_{3}^{2}}{x_{1}\,x_{2}}},{\frac {x_{2}}{x_{1}}}\right\rangle,\left\langle {\frac {x_{3}^{2}}{x_{1}\,x_{2}}},{\frac {x_{1}}{x_{2}}}\right\rangle],\medskip$
\\
$[]$
\\
\end{tabular}
\]
\end{center}

\noindent Note that the torus invariant basis of deformations corresponding to
a Bergman face is given by the set of all lattice points of the polytope
specified above.

\paragraph{Limit map}

The limit map $\lim:B\left(  I\right)  \rightarrow\operatorname*{Poset}\left(
\Delta\right)  $ associates to a face $F$ of $B\left(  I\right)  $ the face of
$\Delta$ formed by the limit points of arcs lying over the weight vectors
$w\in F$, i.e. with lowest order term $t^{w}$.

Labeling the faces of the Bergman complex $B\left(  I\right)  \subset
\operatorname*{Poset}\left(  \nabla\right)  $ and the faces of
$\operatorname*{Poset}\left(  \Delta\right)  $ by the corresponding dual faces
of $\nabla^{\ast}$ and $\Delta^{\ast}$, hence considering the limit map
$\lim:B\left(  I\right)  \rightarrow\operatorname*{Poset}\left(
\Delta\right)  $ as a map $B\left(  I\right)  ^{\ast}\rightarrow
\operatorname*{Poset}\left(  \Delta^{\ast}\right)  $, the limit correspondence
is given by

\begin{center}
\begin{longtable}
[c]{l}
$[],\medskip$
\\
$[\left\langle y_{1},y_{5},y_{6},y_{2}\right\rangle\mapsto \left\langle x_{1},x_{3}\right\rangle,\left\langle y_{1},y_{7},y_{6},y_{3}\right\rangle\mapsto \left\langle x_{2},x_{3}\right\rangle,$
\\
$\left\langle y_{4},y_{5},y_{8},y_{2}\right\rangle\mapsto \left\langle x_{1},x_{0}\right\rangle,\left\langle y_{4},y_{7},y_{8},y_{3}\right\rangle\mapsto \left\langle x_{2},x_{0}\right\rangle],\medskip$
\\
$[\left\langle y_{4},y_{8}\right\rangle\mapsto \left\langle x_{1},x_{2},x_{0}\right\rangle,\left\langle y_{7},y_{3}\right\rangle\mapsto \left\langle x_{2},x_{3},x_{0}\right\rangle,$
\\
$\left\langle y_{5},y_{2}\right\rangle\mapsto \left\langle x_{1},x_{3},x_{0}\right\rangle,\left\langle y_{1},y_{6}\right\rangle\mapsto \left\langle x_{1},x_{2},x_{3}\right\rangle]$
\\
\end{longtable}
\end{center}

The image of the limit map coincides with the image of $\mu$ and with the
Bergman complex of the mirror, i.e. $\lim\left(  B\left(  I\right)  \right)
=\mu\left(  B\left(  I\right)  \right)  =B\left(  I^{\ast}\right)  $.

\paragraph{Mirror complex}

Numbering the vertices of the mirror complex $\mu\left(  B\left(  I\right)
\right)  $ as%
\[
\begin{tabular}
[c]{lll}
$1 = \left(3,-1,-1\right)$ &$2 = \left(-1,3,-1\right)$ &$3 = \left(-1,-1,3\right)$\\
$4 = \left(-1,-1,-1\right)$ & &
\end{tabular}
\]
$\mu\left(  B\left(  I\right)  \right)  \subset\operatorname*{Poset}\left(
\Delta\right)  $ is

\begin{center}
\begin{longtable}
[c]{l}
$[],\medskip$
\\
$[[3],[1],[2],[4]],\medskip$
\\
$[[2, 4],[1, 4],[2, 3],[1, 3]],\medskip$
\\
$[],\medskip$
\\
$[]$
\\
\end{longtable}
\end{center}

\noindent The ordering of the faces of $\mu\left(  B\left(  I\right)  \right)
$ is compatible with above ordering of the faces of $B\left(  I\right)  $. The
$F$-vectors of the faces of $\mu\left(  B\left(  I\right)  \right)  $ are%
\[
\begin{tabular}
[c]{rrrr}
Dimension &Number of faces &F-vector &\\
$1$ & $4$ & $\left(1,2,1,0,0\right)$ &edge\\
$2$ & $4$ & $\left(1,3,3,1,0\right)$ &triangle
\end{tabular}
\]
The first order deformations of the mirror special fiber $X_{0}^{\circ}$
correspond to the lattice points of the dual complex $\left(  \mu\left(
B\left(  I\right)  \right)  \right)  ^{\ast}\subset\operatorname*{Poset}%
\left(  \Delta^{\ast}\right)  $ of the mirror complex. We label the first
order deformations of $X_{0}^{\circ}$ corresponding to vertices of
$\Delta^{\ast}$ by the homogeneous coordinates of $Y$%
\[
\begin{tabular}
[c]{ll}
$x_{1} = x_{\left(1,0,0\right)}={\frac {y_{3}\,{y_{7}}^{2}}{y_{1}\,y_{2}\,y_{4}}}$ &$x_{2} = x_{\left(0,1,0\right)}={\frac {y_{2}\,{y_{5}}^{2}}{y_{1}\,y_{3}\,y_{4}}}$\\
$x_{3} = x_{\left(0,0,1\right)}={\frac {y_{4}^{2}y_{8}}{y_{5}\,y_{6}\,y_{7}}}$ &$x_{0} = x_{\left(-1,-1,-1\right)}={\frac {y_{1}^{2}y_{6}}{y_{5}\,y_{7}\,y_{8}}}$
\end{tabular}
\]
So writing the vertices of the faces of $\left(  \mu\left(  B\left(  I\right)
\right)  \right)  ^{\ast}$ as homogeneous coordinates of $Y$, the complex
$\left(  \mu\left(  B\left(  I\right)  \right)  \right)  ^{\ast}$ is given by

\begin{center}
\begin{longtable}
[c]{l}
$[],\medskip$
\\
$[],\medskip$
\\
$[\left\langle x_{1},x_{3}\right\rangle,\left\langle x_{2},x_{3}\right\rangle,\left\langle x_{1},x_{0}\right\rangle,\left\langle x_{2},x_{0}\right\rangle],\medskip$
\\
$[\left\langle x_{1},x_{2},x_{0}\right\rangle,\left\langle x_{2},x_{3},x_{0}\right\rangle,\left\langle x_{1},x_{3},x_{0}\right\rangle,\left\langle x_{1},x_{2},x_{3}\right\rangle],\medskip$
\\
$[]$
\\
\end{longtable}
\end{center}

\noindent The complex $\mu\left(  B\left(  I\right)  \right)  ^{\ast}$ labeled
by the variables of the Cox ring $S$ gives the ideals of the toric strata of
the special fiber $X_{0}$ of $\mathfrak{X}$, i.e. the complex $SP\left(
I_{0}\right)  $, so in particular the primary decomposition of $I_{0}$ is%
\[
I_{0}=
\begin{tabular}
[c]{l}
$\left\langle x_{1},x_{3}\right\rangle \cap\left\langle x_{1},x_{0}\right\rangle \cap\left\langle x_{2},x_{3}\right\rangle \cap\left\langle x_{2},x_{0}\right\rangle$
\\
\end{tabular}
\]

Labeling the vertices of the faces by the corresponding deformations the
complex is given by

\begin{center}
\begin{longtable}
[c]{l}
$[],\medskip$
\\
$[],\medskip$
\\
$[\left\langle {\frac {y_{3}\,{y_{7}}^{2}}{y_{1}\,y_{2}\,y_{4}}},{\frac {{y_{4}}^{2}y_{8}}{y_{5}\,y_{6}\,y_{7}}}\right\rangle,\left\langle {\frac {y_{2}\,{y_{5}}^{2}}{y_{1}\,y_{3}\,y_{4}}},{\frac {{y_{4}}^{2}y_{8}}{y_{5}\,y_{6}\,y_{7}}}\right\rangle,\left\langle {\frac {y_{1}^{2}y_{6}}{y_{5}\,y_{7}\,y_{8}}},{\frac {y_{3}\,{y_{7}}^{2}}{y_{1}\,y_{2}\,y_{4}}}\right\rangle,$
\\
$\left\langle {\frac {y_{1}^{2}y_{6}}{y_{5}\,y_{7}\,y_{8}}},{\frac {y_{2}\,{y_{5}}^{2}}{y_{1}\,y_{3}\,y_{4}}}\right\rangle],\medskip$
\\
$[\left\langle {\frac {y_{1}^{2}y_{6}}{y_{5}\,y_{7}\,y_{8}}},{\frac {y_{3}\,{y_{7}}^{2}}{y_{1}\,y_{2}\,y_{4}}},{\frac {y_{2}\,{y_{5}}^{2}}{y_{1}\,y_{3}\,y_{4}}}\right\rangle,\left\langle {\frac {y_{1}^{2}y_{6}}{y_{5}\,y_{7}\,y_{8}}},{\frac {y_{2}\,{y_{5}}^{2}}{y_{1}\,y_{3}\,y_{4}}},{\frac {{y_{4}}^{2}y_{8}}{y_{5}\,y_{6}\,y_{7}}}\right\rangle,$
\\
$\left\langle {\frac {y_{1}^{2}y_{6}}{y_{5}\,y_{7}\,y_{8}}},{\frac {y_{3}\,{y_{7}}^{2}}{y_{1}\,y_{2}\,y_{4}}},{\frac {{y_{4}}^{2}y_{8}}{y_{5}\,y_{6}\,y_{7}}}\right\rangle,\left\langle {\frac {y_{3}\,{y_{7}}^{2}}{y_{1}\,y_{2}\,y_{4}}},{\frac {y_{2}\,{y_{5}}^{2}}{y_{1}\,y_{3}\,y_{4}}},{\frac {{y_{4}}^{2}y_{8}}{y_{5}\,y_{6}\,y_{7}}}\right\rangle],\medskip$
\\
$[]$
\\
\end{longtable}
\end{center}

\paragraph{Covering structure in the deformation complex of the mirror
degeneration}

Each face of the complex of deformations $\left(  \mu\left(  B\left(
I\right)  \right)  \right)  ^{\ast}$ of the mirror special fiber $X_{0}%
^{\circ}$ decomposes as the convex hull of 
$2$
polytopes forming
a $
2
:1$ 
trivial
covering of $\mu\left(  B\left(
I\right)  \right)  ^{\vee}$

\begin{center}
\begin{longtable}
[c]{l}
$[],\medskip$
\\
$[],\medskip$
\\
$[[\left\langle x_{3}\right\rangle,\left\langle x_{1}\right\rangle]\mapsto\left\langle x_{1},x_{3}\right\rangle^{\ast\vee}=[2, 4]^{\vee},[\left\langle x_{3}\right\rangle,\left\langle x_{2}\right\rangle]\mapsto\left\langle x_{2},x_{3}\right\rangle^{\ast\vee}=[1, 4]^{\vee},$
\\
$[\left\langle x_{0}\right\rangle,\left\langle x_{1}\right\rangle]\mapsto\left\langle x_{1},x_{0}\right\rangle^{\ast\vee}=[2, 3]^{\vee},[\left\langle x_{0}\right\rangle,\left\langle x_{2}\right\rangle]\mapsto\left\langle x_{2},x_{0}\right\rangle^{\ast\vee}=[1, 3]^{\vee}],\medskip$
\\
$[[\left\langle x_{0}\right\rangle,\left\langle x_{1},x_{2}\right\rangle]\mapsto\left\langle x_{1},x_{2},x_{0}\right\rangle^{\ast\vee}=[3]^{\vee},[\left\langle x_{3},x_{0}\right\rangle,\left\langle x_{2}\right\rangle]\mapsto\left\langle x_{2},x_{3},x_{0}\right\rangle^{\ast\vee}=[1]^{\vee},$
\\
$[\left\langle x_{3},x_{0}\right\rangle,\left\langle x_{1}\right\rangle]\mapsto\left\langle x_{1},x_{3},x_{0}\right\rangle^{\ast\vee}=[2]^{\vee},[\left\langle x_{3}\right\rangle,\left\langle x_{1},x_{2}\right\rangle]\mapsto\left\langle x_{1},x_{2},x_{3}\right\rangle^{\ast\vee}=[4]^{\vee}],\medskip$
\\
$[]$
\\
\end{longtable}
\end{center}

\noindent 
Due to the singularities of
$Y^{\circ}$ this covering involves degenerate faces, i.e. faces $G\mapsto
F^{\vee}$ with $\dim\left(  G\right)  <\dim\left(  F^{\vee}\right)  $.

Writing the vertices of the faces as deformations of $X_{0}^{\circ}$ the
covering is given by

\begin{center}
\begin{longtable}
[c]{l}
$[],\medskip$
\\
$[],\medskip$
\\
$[[\left\langle {\frac {y_{4}^{2}y_{8}}{y_{5}\,y_{6}\,y_{7}}}\right\rangle,\left\langle {\frac {y_{3}\,{y_{7}}^{2}}{y_{1}\,y_{2}\,y_{4}}}\right\rangle]\mapsto[2, 4]^{\vee},[\left\langle {\frac {y_{4}^{2}y_{8}}{y_{5}\,y_{6}\,y_{7}}}\right\rangle,\left\langle {\frac {y_{2}\,{y_{5}}^{2}}{y_{1}\,y_{3}\,y_{4}}}\right\rangle]\mapsto[1, 4]^{\vee},$
\\
$[\left\langle {\frac {y_{1}^{2}y_{6}}{y_{5}\,y_{7}\,y_{8}}}\right\rangle,\left\langle {\frac {y_{3}\,{y_{7}}^{2}}{y_{1}\,y_{2}\,y_{4}}}\right\rangle]\mapsto[2, 3]^{\vee},[\left\langle {\frac {y_{1}^{2}y_{6}}{y_{5}\,y_{7}\,y_{8}}}\right\rangle,\left\langle {\frac {y_{2}\,{y_{5}}^{2}}{y_{1}\,y_{3}\,y_{4}}}\right\rangle]\mapsto[1, 3]^{\vee}],\medskip$
\\
$[[\left\langle {\frac {y_{1}^{2}y_{6}}{y_{5}\,y_{7}\,y_{8}}}\right\rangle,\left\langle {\frac {y_{3}\,{y_{7}}^{2}}{y_{1}\,y_{2}\,y_{4}}},{\frac {y_{2}\,{y_{5}}^{2}}{y_{1}\,y_{3}\,y_{4}}}\right\rangle]\mapsto[3]^{\vee},[\left\langle {\frac {y_{4}^{2}y_{8}}{y_{5}\,y_{6}\,y_{7}}},{\frac {{y_{1}}^{2}y_{6}}{y_{5}\,y_{7}\,y_{8}}}\right\rangle,\left\langle {\frac {y_{2}\,{y_{5}}^{2}}{y_{1}\,y_{3}\,y_{4}}}\right\rangle]\mapsto[1]^{\vee},$
\\
$[\left\langle {\frac {y_{4}^{2}y_{8}}{y_{5}\,y_{6}\,y_{7}}},{\frac {{y_{1}}^{2}y_{6}}{y_{5}\,y_{7}\,y_{8}}}\right\rangle,\left\langle {\frac {y_{3}\,{y_{7}}^{2}}{y_{1}\,y_{2}\,y_{4}}}\right\rangle]\mapsto[2]^{\vee},[\left\langle {\frac {y_{4}^{2}y_{8}}{y_{5}\,y_{6}\,y_{7}}}\right\rangle,\left\langle {\frac {y_{3}\,{y_{7}}^{2}}{y_{1}\,y_{2}\,y_{4}}},{\frac {y_{2}\,{y_{5}}^{2}}{y_{1}\,y_{3}\,y_{4}}}\right\rangle]\mapsto[4]^{\vee}],\medskip$
\\
$[]$
\\
\end{longtable}
\end{center}

\paragraph{Mirror degeneration}

The space of first order deformations of $X_{0}^{\circ}$ in the mirror
degeneration $\mathfrak{X}^{\circ}$ has dimension $
4
$ and the
deformations represented by the monomials%
\[
\left\{  
\begin{tabular}
[c]{llll}
${\frac {y_{4}^{2}y_{8}}{y_{5}\,y_{6}\,y_{7}}}$ & ${\frac {y_{3}\,{y_{7}}^{2}}{y_{1}\,y_{2}\,y_{4}}}$ & ${\frac {y_{1}^{2}y_{6}}{y_{5}\,y_{7}\,y_{8}}}$ & ${\frac {y_{2}\,{y_{5}}^{2}}{y_{1}\,y_{3}\,y_{4}}}$
\end{tabular}
\right\}
\]
form a torus invariant
basis.
The number of lattice points of the dual of the mirror complex of $I$ relates
to the dimension $h^{1,
0
}\left(  X^{\circ}\right)  $ of complex
moduli space of the generic fiber $X^{\circ}$ of $\mathfrak{X}^{\circ}$ and to
the dimension $h^{1,1}\left(  X\right)  $ of the K\"{a}hler moduli space of
the generic fiber $X$ of $\mathfrak{X}$ via%
\begin{align*}
\left\vert \operatorname*{supp}\left(  \left(  \mu\left(  B\left(  I\right)
\right)  \right)  ^{\ast}\right)  \cap N\right\vert  &
=
4
=
3
+
1
\\
&  =\dim\left(  \operatorname*{Aut}\left(  Y^{\circ}\right)  \right)
+h^{1,
0
}\left(  X^{\circ}\right)  =\dim\left(  T\right)
+h^{1,1}\left(  X\right)
\end{align*}

The 
mirror degeneration 
$\mathfrak{X}^{\circ}\subset Y^{\circ}\times\operatorname*{Spec}%
\mathbb{C}\left[  t\right]  $
 of $\mathfrak{X}$ is
given by the ideal 
$I^{\circ}\subset S^{\circ}\otimes\mathbb{C}\left[  t\right]  $
generated by%
\[

\]
Indeed already the ideal $J^{\circ}$ generated by%
\[
\left\{  
\begin{tabular}
[c]{l}
$t (s_{3}\,y_{1}^{2}y_{6}^{2}+s_{1}\,y_{4}^{2}y_{8}^{2} )+y_{5}\,y_{6}\,y_{7}\,y_{8},\medskip$
\\
$t (s_{4}\,y_{2}^{2}y_{5}^{2}+s_{2}\,y_{3}^{2}y_{7}^{2} )+y_{1}\,y_{2}\,y_{3}\,y_{4}$
\\
\end{tabular}
\right\}
\]
defines $\mathfrak{X}^{\circ}$.

\paragraph{Contraction of the mirror degeneration}

In the following we give a birational map relating the degeneration
$\mathfrak{X}^{\circ}$ to a Greene-Plesser type orbifolding mirror family by
contracting divisors on $Y^{\circ}$.
See also Section \ref{sec orbifolding mirror families} below.

In order to contract the divisors

\begin{center}
\begin{longtable}
[c]{ll}
$y_{2} = y_{\left(-1,1,0\right)}={\frac {x_{2}}{x_{1}}}$ &$y_{3} = y_{\left(1,-1,0\right)}={\frac {x_{1}}{x_{2}}}$\\
$y_{6} = y_{\left(0,0,-1\right)}={\frac {x_{0}}{x_{3}}}$ &$y_{8} = y_{\left(0,0,1\right)}={\frac {x_{3}}{x_{0}}}$
\end{longtable}
\end{center}

\noindent consider the $\mathbb{Q}$-factorial toric Fano variety $\hat
{Y}^{\circ}=X\left(  \hat{\Sigma}^{\circ}\right)  $, where $\hat{\Sigma
}^{\circ}=\operatorname*{Fan}\left(  \hat{P}^{\circ}\right)  \subset
M_{\mathbb{R}}$ is the fan over the 
Fano polytope $\hat
{P}^{\circ}\subset M_{\mathbb{R}}$ given as the convex hull of the remaining
vertices of $P^{\circ}=\nabla^{\ast}$ corresponding to the Cox variables%
\[
\begin{tabular}
[c]{ll}
$y_{4} = y_{\left(-1,-1,2\right)}={\frac {x_{3}^{2}}{x_{1}\,x_{2}}}$ &$y_{7} = y_{\left(2,0,-1\right)}={\frac {x_{1}^{2}}{x_{0}\,x_{3}}}$\\
$y_{5} = y_{\left(0,2,-1\right)}={\frac {x_{2}^{2}}{x_{0}\,x_{3}}}$ &$y_{1} = y_{\left(-1,-1,0\right)}={\frac {x_{0}^{2}}{x_{1}\,x_{2}}}$
\end{tabular}
\]
of the toric variety $\hat{Y}^{\circ}$ with Cox ring
\[
\hat{S}^{\circ}=\mathbb{C}
[y_{4},y_{7},y_{5},y_{1}]
\]
The Cox variables of $\hat{Y}^{\circ}$ correspond to the set of Fermat deformations of $\mathfrak{X}$.

Let%
\[
Y^{\circ}=X\left(  \Sigma^{\circ}\right)  \rightarrow X\left(  \hat{\Sigma
}^{\circ}\right)  =\hat{Y}^{\circ}%
\]
be a birational map from $Y^{\circ}$ to a minimal birational model $\hat
{Y}^{\circ}$, which contracts the divisors of the rays $\Sigma^{\circ}\left(
1\right)  -\hat{\Sigma}^{\circ}\left(  1\right)  $ corresponding to the Cox
variables
\[
\begin{tabular}
[c]{llll}
$y_{2}$ &$y_{3}$ &$y_{6}$ &$y_{8}$
\end{tabular}
\]

Representing $\hat{Y}^{\circ}$ as a quotient we have
\[
\hat{Y}^{\circ}=\left(  \mathbb{C}^{
4
}-V\left(  B\left(
\hat{\Sigma}^{\circ}\right)  \right)  \right)  //\hat{G}^{\circ}%
\]
with%
\[
\hat{G}^{\circ}=
\mathbb{Z}_{4}\times\left(  \mathbb{C}^{\ast}\right)  ^{1}
\]
acting via%
\[
\xi y=
\left( \,u_{1}^{2}\,v_{1} \cdot y_{4},\,u_{1}\,v_{1} \cdot y_{7},\,u_{1}^{3}\,v_{1} \cdot y_{5},\,v_{1} \cdot y_{1} \right)
\]
for $\xi=
\left(u_1,v_1\right)
\in\hat{G}^{\circ}$ and $y\in\mathbb{C}%
^{
4
}-V\left(  B\left(  \hat{\Sigma}^{\circ}\right)  \right)  $.

Hence with the group%
\[
\hat{H}^{\circ}=
\mathbb{Z}_{4}
\]
of order 
4
 the toric variety $\hat{Y}^{\circ}$ is the quotient%
\[
\hat{Y}^{\circ}=\mathbb{P}^{
3
}/\hat{H}^{\circ}%
\]
of projective space $\mathbb{P}^{
3
}$.

The first order mirror degeneration $\mathfrak{X}^{1\circ}$ induces via
$Y\rightarrow\hat{Y}^{\circ}$ a degeneration $\mathfrak{\hat{X}}^{1\circ
}\subset\hat{Y}^{\circ}\times\operatorname*{Spec}\mathbb{C}\left[  t\right]
/\left\langle t^{2}\right\rangle $ given by the ideal $\hat{I}^{1\circ}%
\subset\hat{S}^{\circ}\otimes\mathbb{C}\left[  t\right]  /\left\langle
t^{2}\right\rangle $ generated by the
 Fermat-type equations
\[
\left\{  
\begin{tabular}
[c]{l}
$y_{1}\,y_{4}+t (s_{1}\,y_{5}^{2}+s_{2}\,y_{7}^{2} ),\medskip$
\\
$y_{5}\,y_{7}+t (s_{4}\,y_{1}^{2}+s_{3}\,y_{4}^{2} )$
\\
\end{tabular}
\right\}
\]

Note, that%
\[
\left\vert \operatorname*{supp}\left(  \left(  \mu\left(  B\left(  I\right)
\right)  \right)  ^{\ast}\right)  \cap N-\operatorname*{Roots}\left(  \hat
{Y}^{\circ}\right)  \right\vert -\dim\left(  T_{\hat{Y}^{\circ}}\right)
=
4
-
3
=
1
\]
so this family has one independent parameter.

The special fiber $\hat{X}_{0}^{\circ}\subset\hat{Y}^{\circ}$ of
$\mathfrak{\hat{X}}^{1\circ}$ is cut out by the monomial ideal $\hat{I}%
_{0}^{\circ}\subset\hat{S}^{\circ}$ generated by%
\[
\left\{  
\begin{tabular}
[c]{ll}
$y_{1}\,y_{4}$ & $y_{5}\,y_{7}$
\end{tabular}
\right\}
\]

The complex $SP\left(  \hat{I}_{0}^{\circ}\right)  $ of prime ideals of the
toric strata of $\hat{X}_{0}^{\circ}$ is

\begin{center}
\begin{longtable}
[c]{l}
$[],\medskip$
\\
$[],\medskip$
\\
$[\left\langle y_{4},y_{5}\right\rangle,\left\langle y_{4},y_{7}\right\rangle,\left\langle y_{5},y_{1}\right\rangle,\left\langle y_{7},y_{1}\right\rangle],\medskip$
\\
$[\left\langle y_{7},y_{5},y_{1}\right\rangle,\left\langle y_{4},y_{5},y_{1}\right\rangle,\left\langle y_{4},y_{7},y_{1}\right\rangle,\left\langle y_{4},y_{7},y_{5}\right\rangle],\medskip$
\\
$[]$
\\
\end{longtable}
\end{center}

\noindent so $\hat{I}_{0}^{\circ}$ has the primary decomposition%
\[
\hat{I}_{0}^{\circ}=
\begin{tabular}
[c]{l}
$\left\langle y_{5},y_{1}\right\rangle \cap\left\langle y_{7},y_{1}\right\rangle \cap\left\langle y_{4},y_{5}\right\rangle \cap\left\langle y_{4},y_{7}\right\rangle$
\\
\end{tabular}
\]
The Bergman subcomplex $B\left(  I\right)  \subset\operatorname*{Poset}\left(
\nabla\right)  $ induces a subcomplex of $\operatorname*{Poset}\left(
\hat{\nabla}\right)  $, $\hat{\nabla}=\left(  \hat{P}^{\circ}\right)  ^{\ast}$
corresponding to $\hat{I}_{0}^{\circ}$. Indexing of the 
vertices
of $\hat{\nabla}$ by%
\[
\begin{tabular}
[c]{ll}
$1 = \left(\frac{1}{2},\frac{1}{2},2\right)$ &$2 = \left(\frac{3}{2},-\frac{1}{2},0\right)$\\
$3 = \left(-\frac{1}{2},\frac{3}{2},0\right)$ &$4 = \left(-\frac{3}{2},-\frac{3}{2},-2\right)$
\end{tabular}
\]
this complex is given by
\[
\left\{  
\begin{tabular}
[c]{l}
$[],\medskip$
\\
$[[1],[2],[3],[4]],\medskip$
\\
$[[2, 4],[3, 4],[1, 2],[1, 3]],\medskip$
\\
$[],\medskip$
\\
$[]$
\\
\end{tabular}
\right\}
\]

\subsubsection{The $K3$ surface given as the complete intersection of a
generic quadric and a generic cubic in $\mathbb{P}^{4}$}

\paragraph{Setup}

Let $Y=\mathbb{P}^{
4
}=X\left(  \Sigma\right)  $, $\Sigma
=\operatorname*{Fan}\left(  P\right)  =NF\left(  \Delta\right)  \subset
N_{\mathbb{R}}$ with the Fano polytope $P=\Delta^{\ast}$ given by%
\[
\Delta=\operatorname*{convexhull}\left(  
\begin{tabular}
[c]{lll}
$\left(4,-1,-1,-1\right)$ &$\left(-1,4,-1,-1\right)$ &$\left(-1,-1,4,-1\right)$\\
$\left(-1,-1,-1,4\right)$ &$\left(-1,-1,-1,-1\right)$ &
\end{tabular}
\right)  \subset
M_{\mathbb{R}}%
\]
and let%
\[
S=\mathbb{C}
[x_0, x_1, x_2, x_3, x_4]
\]
be the Cox ring of $Y$ with the variables%
\[
\begin{tabular}
[c]{ll}
$x_{1} = x_{\left(1,0,0,0\right)}$ &$x_{2} = x_{\left(0,1,0,0\right)}$\\
$x_{3} = x_{\left(0,0,1,0\right)}$ &$x_{4} = x_{\left(0,0,0,1\right)}$\\
$x_{0} = x_{\left(-1,-1,-1,-1\right)}$ &
\end{tabular}
\]
associated to the rays of $\Sigma$.

Consider the degeneration $\mathfrak{X}\subset Y\times\operatorname*{Spec}%
\mathbb{C}\left[  t\right]  $ of 
complete intersection
K3 surfaces
 of type $\left(
2,3
\right)  $ with monomial special fiber
\[
I_{0}=\left\langle 
\begin{tabular}
[c]{ll}
$x_{1}\,x_{2}$ & $x_{0}\,x_{3}\,x_{4}$
\end{tabular}
\right\rangle
\]
The degeneration is given by the ideal $I\subset S\otimes\mathbb{C}\left[
t\right]  $ with $I_{0}$-reduced generators%
\[
\left\{  
\begin{tabular}
[c]{l}
$x_{1}\,x_{2}+t(s_{1}\,x_{1}^{2}+...+s_{5}\,x_{2}^{2}+...+s_{9}\,x_{3}^{2}+...+s_{12}\,x_{4}^{2}+...+s_{14}\,x_{0}^{2}),\medskip$
\\
$x_{0}\,x_{3}\,x_{4}+t(s_{15}\,x_{1}^{3}+...+s_{25}\,x_{2}^{3}+...+s_{35}\,x_{3}^{3}+...+s_{40}\,x_{4}^{3}+...+s_{43}\,x_{0}^{3})$
\\
\end{tabular}
\right\}
\]

\paragraph{Special fiber Gr\"{o}bner cone}

The space of first order deformations of $\mathfrak{X}$ has dimension
$
42
$ and the deformations represented by the Cox Laurent monomials

\begin{center}
\begin{longtable}
[c]{llllllllll}
${\frac {x_{2}^{3}}{x_{0}\,x_{3}\,x_{4}}}$ & ${\frac {x_{1}^{3}}{x_{0}\,x_{3}\,x_{4}}}$ & ${\frac {x_{0}^{2}}{x_{1}\,x_{2}}}$ & ${\frac {x_{3}^{2}}{x_{1}\,x_{2}}}$ & ${\frac {x_{1}^{2}}{x_{0}\,x_{3}}}$ & ${\frac {x_{2}^{2}}{x_{0}\,x_{3}}}$ & ${\frac {x_{4}^{2}}{x_{1}\,x_{2}}}$ & ${\frac {x_{4}\,x_{0}}{x_{1}\,x_{2}}}$ & ${\frac {x_{3}\,x_{4}}{x_{1}\,x_{2}}}$ & ${\frac {x_{0}^{2}}{x_{3}\,x_{4}}}$\\
${\frac {x_{4}^{2}}{x_{0}\,x_{3}}}$ & ${\frac {x_{3}^{2}}{x_{0}\,x_{4}}}$ & ${\frac {x_{2}^{2}}{x_{0}\,x_{4}}}$ & ${\frac {x_{2}\,x_{0}}{x_{3}\,x_{4}}}$ & ${\frac {x_{1}\,x_{0}}{x_{3}\,x_{4}}}$ & ${\frac {x_{1}^{2}}{x_{3}\,x_{4}}}$ & ${\frac {x_{2}^{2}}{x_{3}\,x_{4}}}$ & ${\frac {x_{3}\,x_{1}}{x_{0}\,x_{4}}}$ & ${\frac {x_{4}\,x_{1}}{x_{0}\,x_{3}}}$ & ${\frac {x_{4}\,x_{2}}{x_{0}\,x_{3}}}$\\
${\frac {x_{3}\,x_{2}}{x_{0}\,x_{4}}}$ & ${\frac {x_{1}^{2}}{x_{0}\,x_{4}}}$ & ${\frac {x_{1}}{x_{4}}}$ & ${\frac {x_{3}}{x_{1}}}$ & ${\frac {x_{2}}{x_{1}}}$ & ${\frac {x_{0}}{x_{3}}}$ & ${\frac {x_{1}}{x_{3}}}$ & ${\frac {x_{2}}{x_{3}}}$ & ${\frac {x_{1}}{x_{0}}}$ & ${\frac {x_{2}}{x_{0}}}$\\
${\frac {x_{3}}{x_{0}}}$ & ${\frac {x_{0}}{x_{1}}}$ & ${\frac {x_{0}}{x_{2}}}$ & ${\frac {x_{3}}{x_{2}}}$ & ${\frac {x_{1}}{x_{2}}}$ & ${\frac {x_{4}}{x_{2}}}$ & ${\frac {x_{4}}{x_{1}}}$ & ${\frac {x_{3}}{x_{4}}}$ & ${\frac {x_{0}}{x_{4}}}$ & ${\frac {x_{4}}{x_{3}}}$\\
${\frac {x_{4}}{x_{0}}}$ & ${\frac {x_{2}}{x_{4}}}$ &  & & & & & & &
\end{longtable}
\end{center}

\noindent form a torus invariant basis. These deformations give linear
inequalities defining the special fiber Gr\"{o}bner cone, i.e.,%
\[
C_{I_{0}}\left(  I\right)  =\left\{  \left(  w_{t},w_{y}\right)  \in
\mathbb{R}\oplus N_{\mathbb{R}}\mid\left\langle w,A^{-1}\left(  m\right)
\right\rangle \geq-w_{t}\forall m\right\}
\]
for above monomials $m$ and the presentation matrix
\[
A=
\left [\begin {array}{cccc} 1&0&0&0\\\noalign{\medskip}0&1&0&0\\\noalign{\medskip}0&0&1&0\\\noalign{\medskip}0&0&0&1\\\noalign{\medskip}-1&-1&-1&-1\end {array}\right ]
\]
of $A_{
3
}\left(  Y\right)  $.

The vertices of special fiber polytope
\[
\nabla=C_{I_{0}}\left(  I\right)  \cap\left\{  w_{t}=1\right\}
\]
and the number of faces, each vertex is contained in, are given by the table

\begin{center}
\begin{longtable}
[c]{rl}
$6$ & \begin{tabular}[c]{llll}$\left(0,0,1,0\right)$ &$\left(0,0,0,1\right)$ &$\left(-1,-1,-1,-1\right)$ &$\left(0,1,0,0\right)$\\$\left(1,0,0,0\right)$ & & &\end{tabular}\medskip\\
$9$ & \begin{tabular}[c]{llll}$\left(1,0,1,0\right)$ &$\left(0,1,1,0\right)$ &$\left(1,0,0,1\right)$ &$\left(0,1,0,1\right)$\\$\left(0,-1,-1,-1\right)$ &$\left(-1,0,-1,-1\right)$ & &\end{tabular}
\end{longtable}
\end{center}

\noindent The number of faces of $\nabla$ and their $F$-vectors are

\begin{center}
\begin{longtable}
[c]{rrll}
Dimension &Number of faces &F-vector &\\
$-1$ & $1$ & $\left(0,0,0,0,0,0\right)$ &\\
$0$ & $11$ & $\left(1,1,0,0,0,0\right)$ &point\\
$1$ & $25$ & $\left(1,2,1,0,0,0\right)$ &edge\\
$2$ & $12$ & $\left(1,3,3,1,0,0\right)$ &triangle\\
$2$ & $12$ & $\left(1,4,4,1,0,0\right)$ &quadrangle\\
$3$ & $2$ & $\left(1,4,6,4,1,0\right)$ &tetrahedron\\
$3$ & $8$ & $\left(1,6,9,5,1,0\right)$ &prism\\
$4$ & $1$ & $\left(1,11,25,24,10,1\right)$ &
\end{longtable}
\end{center}

The dual $P^{\circ}=\nabla^{\ast}$ of $\nabla$ is a Fano polytope with vertices

\begin{center}
\begin{longtable}
[c]{lllll}
$\left(-1,-1,0,2\right)$ &$\left(-1,-1,0,0\right)$ &$\left(-1,1,0,0\right)$ &$\left(1,-1,0,0\right)$ &$\left(-1,-1,2,0\right)$\\
$\left(0,3,-1,-1\right)$ &$\left(0,0,-1,2\right)$ &$\left(0,0,-1,-1\right)$ &$\left(3,0,-1,-1\right)$ &$\left(0,0,2,-1\right)$
\end{longtable}
\end{center}

\noindent and the $F$-vectors of the faces of $P^{\circ}$ are

\begin{center}
\begin{longtable}
[c]{rrll}
Dimension &Number of faces &F-vector &\\
$-1$ & $1$ & $\left(0,0,0,0,0,0\right)$ &\\
$0$ & $10$ & $\left(1,1,0,0,0,0\right)$ &point\\
$1$ & $24$ & $\left(1,2,1,0,0,0\right)$ &edge\\
$2$ & $9$ & $\left(1,4,4,1,0,0\right)$ &quadrangle\\
$2$ & $16$ & $\left(1,3,3,1,0,0\right)$ &triangle\\
$3$ & $6$ & $\left(1,6,9,5,1,0\right)$ &prism\\
$3$ & $5$ & $\left(1,4,6,4,1,0\right)$ &tetrahedron\\
$4$ & $1$ & $\left(1,10,24,25,11,1\right)$ &
\end{longtable}
\end{center}

The polytopes $\nabla$ and $P^{\circ}$ have $0$ as the unique interior lattice
point, and the Fano polytope $P^{\circ}$ defines the embedding toric Fano
variety $Y^{\circ}=X\left(  \operatorname*{Fan}\left(  P^{\circ}\right)
\right)  $ of the fibers of the mirror degeneration. The dimension of the
automorphism group of $Y^{\circ}$ is
\[
\dim\left(  \operatorname*{Aut}\left(  Y^{\circ}\right)  \right)
=
4
\]
Let%
\[
S^{\circ}=\mathbb{C}
[y_{1},y_{2},...,y_{10}]
\]
be the Cox ring of $Y^{\circ}$ with variables

\begin{center}
\begin{longtable}
[c]{lll}
$y_{1} = y_{\left(-1,-1,0,2\right)}={\frac {x_{4}^{2}}{x_{1}\,x_{2}}}$ &$y_{2} = y_{\left(-1,-1,0,0\right)}={\frac {x_{0}^{2}}{x_{1}\,x_{2}}}$ &$y_{3} = y_{\left(-1,1,0,0\right)}={\frac {x_{2}}{x_{1}}}$\\
$y_{4} = y_{\left(1,-1,0,0\right)}={\frac {x_{1}}{x_{2}}}$ &$y_{5} = y_{\left(-1,-1,2,0\right)}={\frac {x_{3}^{2}}{x_{1}\,x_{2}}}$ &$y_{6} = y_{\left(0,3,-1,-1\right)}={\frac {x_{2}^{3}}{x_{0}\,x_{3}\,x_{4}}}$\\
$y_{7} = y_{\left(0,0,-1,2\right)}={\frac {x_{4}^{2}}{x_{0}\,x_{3}}}$ &$y_{8} = y_{\left(0,0,-1,-1\right)}={\frac {x_{0}^{2}}{x_{3}\,x_{4}}}$ &$y_{9} = y_{\left(3,0,-1,-1\right)}={\frac {x_{1}^{3}}{x_{0}\,x_{3}\,x_{4}}}$\\
$y_{10} = y_{\left(0,0,2,-1\right)}={\frac {x_{3}^{2}}{x_{0}\,x_{4}}}$ & &
\end{longtable}
\end{center}

\paragraph{Bergman subcomplex}

Intersecting the tropical variety of $I$ with the special fiber Gr\"{o}bner
cone $C_{I_{0}}\left(  I\right)  $ we obtain the special fiber Bergman
subcomplex
\[
B\left(  I\right)  =C_{I_{0}}\left(  I\right)  \cap BF\left(  I\right)
\cap\left\{  w_{t}=1\right\}
\]
The following table chooses an indexing of the vertices of $\nabla$ involved
in $B\left(  I\right)  $ and gives for each vertex the numbers $n_{\nabla}$
and $n_{B}$ of faces of $\nabla$ and faces of $B\left(  I\right)  $ it is
contained in

\begin{center}
\begin{longtable}
[c]{rrl}
$n_{\nabla}$ &$n_{B(I)}$ & \medskip\\
$9$ & $3$ & \begin{tabular}[c]{lll}$1 = \left(1,0,1,0\right)$ &$2 = \left(0,1,1,0\right)$ &$3 = \left(1,0,0,1\right)$\\$4 = \left(0,1,0,1\right)$ &$5 = \left(0,-1,-1,-1\right)$ &$6 = \left(-1,0,-1,-1\right)$\end{tabular}
\end{longtable}
\end{center}

\noindent With this indexing the Bergman subcomplex $B\left(  I\right)  $ of
$\operatorname*{Poset}\left(  \nabla\right)  $ associated to the degeneration
$\mathfrak{X}$ is

\begin{center}
\begin{longtable}
[c]{l}
$[],\medskip$
\\
$[[1],[2],[3],[4],[5],[6]],\medskip$
\\
$[[3, 5],[3, 4],[4, 6],[1, 5],[1, 3],[1, 2],[5, 6],[2, 6],[2, 4]],\medskip$
\\
$[[3, 4, 5, 6],[1, 3, 5],[2, 4, 6],[1, 2, 5, 6],[1, 2, 3, 4]],\medskip$
\\
$[],\medskip$
\\
$[]$
\\
\end{longtable}
\end{center}

$B\left(  I\right)  $ has the $F$-vector%
\[
\begin{tabular}
[c]{rrrrrrr}
Dimension &$-1$ & $0$ & $1$ & $2$ & $3$ & $4$\\
Number of faces &$0$ & $6$ & $9$ & $5$ & $0$ & $0$
\end{tabular}
\]
and the $F$-vectors of its faces are

\begin{center}
\begin{longtable}
[c]{rrll}
Dimension &Number of faces &F-vector &\\
$0$ & $6$ & $\left(1,1,0,0,0,0\right)$ &point\\
$1$ & $9$ & $\left(1,2,1,0,0,0\right)$ &edge\\
$2$ & $3$ & $\left(1,4,4,1,0,0\right)$ &quadrangle\\
$2$ & $2$ & $\left(1,3,3,1,0,0\right)$ &triangle
\end{longtable}
\end{center}

\paragraph{Dual complex}

The dual complex $\operatorname*{dual}\left(  B\left(  I\right)  \right)
=\left(  B\left(  I\right)  \right)  ^{\ast}$ of deformations associated to
$B\left(  I\right)  $ via initial ideals is given by

\begin{center}
\begin{longtable}
[c]{l}
$[],\medskip$
\\
$[],\medskip$
\\
$[[3, 4, 5, 6]^{*}=\left\langle {\frac {x_{3}^{2}}{x_{1}\,x_{2}}},{\frac {x_{3}^{2}}{x_{0}\,x_{4}}}\right\rangle,[1, 3, 5]^{*}=\left\langle {\frac {x_{2}^{3}}{x_{0}\,x_{3}\,x_{4}}},{\frac {x_{2}}{x_{1}}}\right\rangle,[2, 4, 6]^{*}=\left\langle {\frac {x_{1}^{3}}{x_{0}\,x_{3}\,x_{4}}},{\frac {x_{1}}{x_{2}}}\right\rangle,$
\\
$[1, 2, 5, 6]^{*}=\left\langle {\frac {x_{4}^{2}}{x_{1}\,x_{2}}},{\frac {x_{4}^{2}}{x_{0}\,x_{3}}}\right\rangle,[1, 2, 3, 4]^{*}=\left\langle {\frac {x_{0}^{2}}{x_{1}\,x_{2}}},{\frac {x_{0}^{2}}{x_{3}\,x_{4}}}\right\rangle],\medskip$
\\
$[[3, 5]^{*}=\left\langle {\frac {x_{3}^{2}}{x_{1}\,x_{2}}},{\frac {x_{2}^{3}}{x_{0}\,x_{3}\,x_{4}}},{\frac {x_{3}^{2}}{x_{0}\,x_{4}}},{\frac {x_{2}}{x_{1}}}\right\rangle,[3, 4]^{*}=\left\langle {\frac {x_{3}^{2}}{x_{1}\,x_{2}}},{\frac {x_{0}^{2}}{x_{1}\,x_{2}}},{\frac {x_{3}^{2}}{x_{0}\,x_{4}}},{\frac {x_{0}^{2}}{x_{3}\,x_{4}}}\right\rangle,$
\\
$[4, 6]^{*}=\left\langle {\frac {x_{3}^{2}}{x_{1}\,x_{2}}},{\frac {x_{1}^{3}}{x_{0}\,x_{3}\,x_{4}}},{\frac {x_{3}^{2}}{x_{0}\,x_{4}}},{\frac {x_{1}}{x_{2}}}\right\rangle,[1, 5]^{*}=\left\langle {\frac {x_{4}^{2}}{x_{1}\,x_{2}}},{\frac {x_{2}^{3}}{x_{0}\,x_{3}\,x_{4}}},{\frac {x_{4}^{2}}{x_{0}\,x_{3}}},{\frac {x_{2}}{x_{1}}}\right\rangle,$
\\
$[1, 3]^{*}=\left\langle {\frac {x_{0}^{2}}{x_{1}\,x_{2}}},{\frac {x_{2}^{3}}{x_{0}\,x_{3}\,x_{4}}},{\frac {x_{0}^{2}}{x_{3}\,x_{4}}},{\frac {x_{2}}{x_{1}}}\right\rangle,[1, 2]^{*}=\left\langle {\frac {x_{4}^{2}}{x_{1}\,x_{2}}},{\frac {x_{0}^{2}}{x_{1}\,x_{2}}},{\frac {x_{4}^{2}}{x_{0}\,x_{3}}},{\frac {x_{0}^{2}}{x_{3}\,x_{4}}}\right\rangle,$
\\
$[5, 6]^{*}=\left\langle {\frac {x_{4}^{2}}{x_{1}\,x_{2}}},{\frac {x_{3}^{2}}{x_{1}\,x_{2}}},{\frac {x_{4}^{2}}{x_{0}\,x_{3}}},{\frac {x_{3}^{2}}{x_{0}\,x_{4}}}\right\rangle,[2, 6]^{*}=\left\langle {\frac {x_{4}^{2}}{x_{1}\,x_{2}}},{\frac {x_{1}^{3}}{x_{0}\,x_{3}\,x_{4}}},{\frac {x_{4}^{2}}{x_{0}\,x_{3}}},{\frac {x_{1}}{x_{2}}}\right\rangle,$
\\
\end{longtable}
\end{center}

\noindent when writing the vertices of the faces as deformations of $X_{0}$.
Note that the $T$-invariant basis of deformations associated to a face is
given by all lattice points of the corresponding polytope in $M_{\mathbb{R}}$.

When numbering the vertices of the faces of $\operatorname*{dual}\left(
B\left(  I\right)  \right)  $ by the Cox variables of the mirror toric Fano
variety $Y^{\circ}$ the complex $\operatorname*{dual}\left(  B\left(
I\right)  \right)  $ is

\begin{center}
\begin{longtable}
[c]{l}
$[],\medskip$
\\
$[],\medskip$
\\
$[[3, 4, 5, 6]^{*}=\left\langle y_{5},y_{10}\right\rangle,[1, 3, 5]^{*}=\left\langle y_{6},y_{3}\right\rangle,[2, 4, 6]^{*}=\left\langle y_{9},y_{4}\right\rangle,[1, 2, 5, 6]^{*}=\left\langle y_{1},y_{7}\right\rangle,$
\\
$[1, 2, 3, 4]^{*}=\left\langle y_{2},y_{8}\right\rangle],\medskip$
\\
$[[3, 5]^{*}=\left\langle y_{5},y_{6},y_{10},y_{3}\right\rangle,[3, 4]^{*}=\left\langle y_{5},y_{2},y_{10},y_{8}\right\rangle,[4, 6]^{*}=\left\langle y_{5},y_{9},y_{10},y_{4}\right\rangle,$
\\
$[1, 5]^{*}=\left\langle y_{1},y_{6},y_{7},y_{3}\right\rangle,[1, 3]^{*}=\left\langle y_{2},y_{6},y_{8},y_{3}\right\rangle,[1, 2]^{*}=\left\langle y_{1},y_{2},y_{7},y_{8}\right\rangle,$
\\
$[5, 6]^{*}=\left\langle y_{1},y_{5},y_{7},y_{10}\right\rangle,[2, 6]^{*}=\left\langle y_{1},y_{9},y_{7},y_{4}\right\rangle,[2, 4]^{*}=\left\langle y_{2},y_{9},y_{8},y_{4}\right\rangle],\medskip$
\\
$[[1]^{*}=\left\langle y_{1},y_{2},y_{6},y_{7},y_{8},y_{3}\right\rangle,[2]^{*}=\left\langle y_{1},y_{2},y_{9},y_{7},y_{8},y_{4}\right\rangle,[3]^{*}=\left\langle y_{5},y_{2},y_{6},y_{10},y_{8},y_{3}\right\rangle,$
\\
$[4]^{*}=\left\langle y_{5},y_{2},y_{9},y_{10},y_{8},y_{4}\right\rangle,[5]^{*}=\left\langle y_{1},y_{5},y_{6},y_{7},y_{10},y_{3}\right\rangle,[6]^{*}=\left\langle y_{1},y_{5},y_{9},y_{7},y_{10},y_{4}\right\rangle],\medskip$
\\
$[]$
\\
\end{longtable}
\end{center}

\noindent The dual complex has the $F$-vector%
\[
\begin{tabular}
[c]{rrrrrrr}
Dimension &$-1$ & $0$ & $1$ & $2$ & $3$ & $4$\\
Number of faces &$0$ & $0$ & $5$ & $9$ & $6$ & $0$
\end{tabular}
\]
and the $F$-vectors of the faces of $\operatorname*{dual}\left(  B\left(
I\right)  \right)  $ are

\begin{center}
\begin{longtable}
[c]{rrll}
Dimension &Number of faces &F-vector &\\
$1$ & $5$ & $\left(1,2,1,0,0,0\right)$ &edge\\
$2$ & $9$ & $\left(1,4,4,1,0,0\right)$ &quadrangle\\
$3$ & $6$ & $\left(1,6,9,5,1,0\right)$ &prism
\end{longtable}
\end{center}

Recall that in this example the toric variety $Y$ is projective space. The
number of lattice points of the support of $\operatorname*{dual}\left(
B\left(  I\right)  \right)  $ relates to the dimension $h^{1,
1
}%
\left(  X\right)  $ of the complex moduli space of the generic fiber $X$ of
$\mathfrak{X}$ and to the dimension $h^{1,1}\left(  \bar{X}^{\circ}\right)  $
of the K\"{a}hler moduli space of the $MPCR$-blowup $\bar{X}^{\circ}$ of the
generic fiber $X^{\circ}$ of the mirror degeneration
\begin{align*}
\left\vert \operatorname*{supp}\left(  \operatorname*{dual}\left(  B\left(
I\right)  \right)  \right)  \cap M\right\vert  &
=
42
=
24
+
18
=\dim\left(  \operatorname*{Aut}%
\left(  Y\right)  \right)  +h^{1,
1
}\left(  X\right) \\
&  =
20
+
4
+
18
\\
&  =\left\vert \operatorname*{Roots}\left(  Y\right)  \right\vert +\dim\left(
T_{Y}\right)  +h^{1,1}\left(  \bar{X}^{\circ}\right)
\end{align*}
There are%
\[
h^{1,
1
}\left(  X\right)  +\dim\left(  T_{Y^{\circ}}\right)
=
18
+
4
\]
non-trivial toric polynomial deformations of $X_{0}$

\begin{center}
\begin{longtable}
[c]{llllllllll}
${\frac {x_{3}^{2}}{x_{0}\,x_{4}}}$ & ${\frac {x_{3}^{2}}{x_{1}\,x_{2}}}$ & ${\frac {x_{2}^{3}}{x_{0}\,x_{3}\,x_{4}}}$ & ${\frac {x_{1}^{3}}{x_{0}\,x_{3}\,x_{4}}}$ & ${\frac {x_{4}^{2}}{x_{0}\,x_{3}}}$ & ${\frac {x_{4}^{2}}{x_{1}\,x_{2}}}$ & ${\frac {x_{0}^{2}}{x_{3}\,x_{4}}}$ & ${\frac {x_{0}^{2}}{x_{1}\,x_{2}}}$ & ${\frac {x_{2}^{2}}{x_{0}\,x_{4}}}$ & ${\frac {x_{3}\,x_{2}}{x_{0}\,x_{4}}}$\\
${\frac {x_{1}^{2}}{x_{0}\,x_{4}}}$ & ${\frac {x_{3}\,x_{1}}{x_{0}\,x_{4}}}$ & ${\frac {x_{2}^{2}}{x_{0}\,x_{3}}}$ & ${\frac {x_{4}\,x_{2}}{x_{0}\,x_{3}}}$ & ${\frac {x_{2}\,x_{0}}{x_{3}\,x_{4}}}$ & ${\frac {x_{2}^{2}}{x_{3}\,x_{4}}}$ & ${\frac {x_{0}\,x_{4}}{x_{1}\,x_{2}}}$ & ${\frac {x_{3}\,x_{4}}{x_{1}\,x_{2}}}$ & ${\frac {x_{1}^{2}}{x_{0}\,x_{3}}}$ & ${\frac {x_{4}\,x_{1}}{x_{0}\,x_{3}}}$\\
${\frac {x_{1}\,x_{0}}{x_{3}\,x_{4}}}$ & ${\frac {x_{1}^{2}}{x_{3}\,x_{4}}}$ &  & & & & & & &
\end{longtable}
\end{center}

\noindent They correspond to the toric divisors

\begin{center}
\begin{longtable}
[c]{lllll}
$D_{\left(0,0,2,-1\right)}$ & $D_{\left(-1,-1,2,0\right)}$ & $D_{\left(0,3,-1,-1\right)}$ & $D_{\left(3,0,-1,-1\right)}$ & $D_{\left(0,0,-1,2\right)}$\\
$D_{\left(-1,-1,0,2\right)}$ & $D_{\left(0,0,-1,-1\right)}$ & $D_{\left(-1,-1,0,0\right)}$ & $D_{\left(0,2,0,-1\right)}$ & $D_{\left(0,1,1,-1\right)}$\\
$D_{\left(2,0,0,-1\right)}$ & $D_{\left(1,0,1,-1\right)}$ & $D_{\left(0,2,-1,0\right)}$ & $D_{\left(0,1,-1,1\right)}$ & $D_{\left(0,1,-1,-1\right)}$\\
$D_{\left(0,2,-1,-1\right)}$ & $D_{\left(-1,-1,0,1\right)}$ & $D_{\left(-1,-1,1,1\right)}$ & $D_{\left(2,0,-1,0\right)}$ & $D_{\left(1,0,-1,1\right)}$\\
$D_{\left(1,0,-1,-1\right)}$ & $D_{\left(2,0,-1,-1\right)}$ &  & &
\end{longtable}
\end{center}

\noindent on a MPCR-blowup of $Y^{\circ}$ inducing $
18
$ non-zero
toric divisor classes on the mirror $X^{\circ}$. The following $
20
$
toric divisors of $Y$ induce the trivial divisor class on $X^{\circ}$

\begin{center}
\begin{longtable}
[c]{lllll}
$D_{\left(-1,1,0,0\right)}$ & $D_{\left(1,-1,0,0\right)}$ & $D_{\left(-1,0,1,0\right)}$ & $D_{\left(0,0,0,-1\right)}$ & $D_{\left(0,0,1,-1\right)}$\\
$D_{\left(0,-1,1,0\right)}$ & $D_{\left(-1,0,0,1\right)}$ & $D_{\left(-1,0,0,0\right)}$ & $D_{\left(0,0,-1,0\right)}$ & $D_{\left(0,0,-1,1\right)}$\\
$D_{\left(0,0,1,0\right)}$ & $D_{\left(0,0,0,1\right)}$ & $D_{\left(0,-1,0,1\right)}$ & $D_{\left(0,-1,0,0\right)}$ & $D_{\left(0,1,-1,0\right)}$\\
$D_{\left(1,0,-1,0\right)}$ & $D_{\left(0,1,0,-1\right)}$ & $D_{\left(1,0,0,-1\right)}$ & $D_{\left(0,1,0,0\right)}$ & $D_{\left(1,0,0,0\right)}$
\end{longtable}
\end{center}

\paragraph{Mirror special fiber}

The ideal $I_{0}^{\circ}$ of the monomial special fiber $X_{0}^{\circ}$ of the
mirror degeneration $\mathfrak{X}^{\circ}$ is generated by the following set
of monomials in $S^{\circ}$%
\[
\left\{  
\begin{tabular}
[c]{lllll}
$y_{4}\,y_{5}\,y_{6}\,y_{7}\,y_{8}$ & $y_{1}\,y_{2}\,y_{5}\,y_{6}\,y_{9}$ & $y_{3}\,y_{4}\,y_{5}\,y_{7}\,y_{8}$ & $y_{10}\,y_{4}\,y_{6}\,y_{7}\,y_{8}$ & $y_{1}\,y_{3}\,y_{4}\,y_{5}\,y_{8}$\\
$y_{3}\,y_{5}\,y_{7}\,y_{8}\,y_{9}$ & $y_{10}\,y_{1}\,y_{4}\,y_{6}\,y_{8}$ & $y_{10}\,y_{1}\,y_{6}\,y_{8}\,y_{9}$ & $y_{10}\,y_{6}\,y_{7}\,y_{8}\,y_{9}$ & $y_{1}\,y_{5}\,y_{6}\,y_{8}\,y_{9}$\\
$y_{10}\,y_{1}\,y_{3}\,y_{4}\,y_{8}$ & $y_{1}\,y_{2}\,y_{3}\,y_{4}\,y_{5}$ & $y_{10}\,y_{1}\,y_{2}\,y_{3}\,y_{4}$ & $y_{1}\,y_{4}\,y_{5}\,y_{6}\,y_{8}$ & $y_{1}\,y_{3}\,y_{5}\,y_{8}\,y_{9}$\\
$y_{10}\,y_{3}\,y_{7}\,y_{8}\,y_{9}$ & $y_{1}\,y_{2}\,y_{4}\,y_{5}\,y_{6}$ & $y_{2}\,y_{3}\,y_{4}\,y_{5}\,y_{7}$ & $y_{5}\,y_{6}\,y_{7}\,y_{8}\,y_{9}$ & $y_{10}\,y_{2}\,y_{3}\,y_{4}\,y_{7}$\\
$y_{10}\,y_{1}\,y_{2}\,y_{6}\,y_{9}$ & $y_{10}\,y_{2}\,y_{6}\,y_{7}\,y_{9}$ & $y_{10}\,y_{1}\,y_{3}\,y_{8}\,y_{9}$ & $y_{10}\,y_{1}\,y_{2}\,y_{3}\,y_{9}$ & $y_{2}\,y_{5}\,y_{6}\,y_{7}\,y_{9}$\\
$y_{2}\,y_{3}\,y_{5}\,y_{7}\,y_{9}$ & $y_{10}\,y_{2}\,y_{3}\,y_{7}\,y_{9}$ & $y_{2}\,y_{4}\,y_{5}\,y_{6}\,y_{7}$ & $y_{1}\,y_{2}\,y_{3}\,y_{5}\,y_{9}$ & $y_{10}\,y_{3}\,y_{4}\,y_{7}\,y_{8}$\\
$y_{10}\,y_{2}\,y_{4}\,y_{6}\,y_{7}$ & $y_{10}\,y_{1}\,y_{2}\,y_{4}\,y_{6}$ &  & &
\end{tabular}
\right\}
\]
Indeed already the ideal%
\[
J_{0}^{\circ}=\left\langle 
\begin{tabular}
[c]{ll}
$y_{1}\,y_{2}\,y_{3}\,y_{4}\,y_{5}$ & $y_{10}\,y_{6}\,y_{7}\,y_{8}\,y_{9}$
\end{tabular}
\right\rangle
\]
defines the same subvariety of the toric variety $Y^{\circ}$, and
$J_{0}^{\circ\Sigma}=I_{0}^{\circ}$. Recall that passing from $J_{0}^{\circ}$
to $J_{0}^{\circ\Sigma}$ is the non-simplicial toric analogue of saturation.

The complex $B\left(  I\right)  ^{\ast}$ labeled by the variables of the Cox
ring $S^{\circ}$ of $Y^{\ast}$, as written in the last section, is the complex
$SP\left(  I_{0}^{\circ}\right)  $ of prime ideals of the toric strata of the
special fiber $X_{0}^{\circ}$ of the mirror degeneration $\mathfrak{X}^{\circ
}$, i.e., the primary decomposition of $I_{0}^{\circ}$ is%
\[
I_{0}^{\circ}=
\begin{tabular}
[c]{l}
$\left\langle y_{1},y_{7}\right\rangle \cap\left\langle y_{6},y_{3}\right\rangle \cap\left\langle y_{5},y_{10}\right\rangle \cap\left\langle y_{2},y_{8}\right\rangle \cap\left\langle y_{9},y_{4}\right\rangle$
\\
\end{tabular}
\]
Each facet $F\in B\left(  I\right)  $ corresponds to one of these ideals and this
ideal is generated by the products of facets of $\nabla$ containing $F$.

\paragraph{Covering structure in the deformation complex of the
degeneration $\mathfrak{X}$}

According to the local reduced Gr\"{o}bner basis each face of the complex of
deformations $\operatorname*{dual}\left(  B\left(  I\right)  \right)  $
decomposes into 
$2$
polytopes forming a $
2
:1$
trivial
covering of $B\left(  I\right)  $

\begin{center}
\begin{longtable}
[c]{l}
$[],\medskip$
\\
$[],\medskip$
\\
$[[\left\langle y_{5}\right\rangle,\left\langle y_{10}\right\rangle]\mapsto\left\langle y_{5},y_{10}\right\rangle^{\ast\vee}=[3, 4, 5, 6]^{\vee},[\left\langle y_{3}\right\rangle,\left\langle y_{6}\right\rangle]\mapsto\left\langle y_{6},y_{3}\right\rangle^{\ast\vee}=[1, 3, 5]^{\vee},$
\\
$[\left\langle y_{4}\right\rangle,\left\langle y_{9}\right\rangle]\mapsto\left\langle y_{9},y_{4}\right\rangle^{\ast\vee}=[2, 4, 6]^{\vee},[\left\langle y_{1}\right\rangle,\left\langle y_{7}\right\rangle]\mapsto\left\langle y_{1},y_{7}\right\rangle^{\ast\vee}=[1, 2, 5, 6]^{\vee},$
\\
$[\left\langle y_{2}\right\rangle,\left\langle y_{8}\right\rangle]\mapsto\left\langle y_{2},y_{8}\right\rangle^{\ast\vee}=[1, 2, 3, 4]^{\vee}],\medskip$
\\
$[[\left\langle y_{5},y_{3}\right\rangle,\left\langle y_{6},y_{10}\right\rangle]\mapsto\left\langle y_{5},y_{6},y_{10},y_{3}\right\rangle^{\ast\vee}=[3, 5]^{\vee},$
\\
$[\left\langle y_{5},y_{2}\right\rangle,\left\langle y_{10},y_{8}\right\rangle]\mapsto\left\langle y_{5},y_{2},y_{10},y_{8}\right\rangle^{\ast\vee}=[3, 4]^{\vee},$
\\
$[\left\langle y_{5},y_{4}\right\rangle,\left\langle y_{9},y_{10}\right\rangle]\mapsto\left\langle y_{5},y_{9},y_{10},y_{4}\right\rangle^{\ast\vee}=[4, 6]^{\vee},$
\\
\end{longtable}
\end{center}

\noindent Here the faces $F\in B\left(  I\right)  $ are specified both via the
vertices of $F^{\ast}$ labeled by the variables of $S^{\circ}$ and by the
numbering of the vertices of $B\left(  I\right)  $ chosen above.

\noindent This covering has 
2 sheets forming the complexes
 
\[
\begin{tabular}
[c]{l}
$[],\medskip$
\\
$[],\medskip$
\\
$[\left\langle y_{5}\right\rangle,\left\langle y_{3}\right\rangle,\left\langle y_{4}\right\rangle,\left\langle y_{1}\right\rangle,\left\langle y_{2}\right\rangle],\medskip$
\\
$[\left\langle y_{5},y_{3}\right\rangle,\left\langle y_{5},y_{2}\right\rangle,\left\langle y_{5},y_{4}\right\rangle,\left\langle y_{1},y_{3}\right\rangle,\left\langle y_{2},y_{3}\right\rangle,\left\langle y_{1},y_{2}\right\rangle,\left\langle y_{1},y_{5}\right\rangle,\left\langle y_{1},y_{4}\right\rangle,\left\langle y_{2},y_{4}\right\rangle],\medskip$
\\
$[\left\langle y_{1},y_{2},y_{3}\right\rangle,\left\langle y_{1},y_{2},y_{4}\right\rangle,\left\langle y_{5},y_{2},y_{3}\right\rangle,\left\langle y_{5},y_{2},y_{4}\right\rangle,\left\langle y_{1},y_{5},y_{3}\right\rangle,\left\langle y_{1},y_{5},y_{4}\right\rangle],\medskip$
\\
$[]$
\\
\end{tabular}
\]
\[
\begin{tabular}
[c]{l}
$[],\medskip$
\\
$[],\medskip$
\\
$[\left\langle y_{10}\right\rangle,\left\langle y_{6}\right\rangle,\left\langle y_{9}\right\rangle,\left\langle y_{7}\right\rangle,\left\langle y_{8}\right\rangle],\medskip$
\\
$[\left\langle y_{6},y_{10}\right\rangle,\left\langle y_{10},y_{8}\right\rangle,\left\langle y_{9},y_{10}\right\rangle,\left\langle y_{6},y_{7}\right\rangle,\left\langle y_{6},y_{8}\right\rangle,\left\langle y_{7},y_{8}\right\rangle,\left\langle y_{7},y_{10}\right\rangle,\left\langle y_{9},y_{7}\right\rangle,\left\langle y_{9},y_{8}\right\rangle],\medskip$
\\
$[\left\langle y_{6},y_{7},y_{8}\right\rangle,\left\langle y_{9},y_{7},y_{8}\right\rangle,\left\langle y_{6},y_{10},y_{8}\right\rangle,\left\langle y_{9},y_{10},y_{8}\right\rangle,\left\langle y_{6},y_{7},y_{10}\right\rangle,\left\langle y_{9},y_{7},y_{10}\right\rangle],\medskip$
\\
$[]$
\\
\end{tabular}
\]
with $F$-vector

\begin{center}
\begin{longtable}
[c]{rrll}
Dimension &Number of faces &F-vector &\\
$0$ & $5$ & $\left(1,1,0,0,0,0\right)$ &point\\
$1$ & $9$ & $\left(1,2,1,0,0,0\right)$ &edge\\
$2$ & $6$ & $\left(1,3,3,1,0,0\right)$ &triangle
\end{longtable}
\end{center}

\noindent Writing the vertices of the faces as deformations the covering is
given by

\begin{center}
\begin{longtable}
[c]{l}
$[],\medskip$
\\
$[],\medskip$
\\
$[[\left\langle {\frac {x_{3}^{2}}{x_{1}\,x_{2}}}\right\rangle,\left\langle {\frac {x_{3}^{2}}{x_{0}\,x_{4}}}\right\rangle]\mapsto[3, 4, 5, 6]^{\vee},[\left\langle {\frac {x_{2}}{x_{1}}}\right\rangle,\left\langle {\frac {x_{2}^{3}}{x_{0}\,x_{3}\,x_{4}}}\right\rangle]\mapsto[1, 3, 5]^{\vee},$
\\
$[\left\langle {\frac {x_{1}}{x_{2}}}\right\rangle,\left\langle {\frac {x_{1}^{3}}{x_{0}\,x_{3}\,x_{4}}}\right\rangle]\mapsto[2, 4, 6]^{\vee},[\left\langle {\frac {x_{4}^{2}}{x_{1}\,x_{2}}}\right\rangle,\left\langle {\frac {x_{4}^{2}}{x_{0}\,x_{3}}}\right\rangle]\mapsto[1, 2, 5, 6]^{\vee},$
\\
$[\left\langle {\frac {x_{0}^{2}}{x_{1}\,x_{2}}}\right\rangle,\left\langle {\frac {x_{0}^{2}}{x_{3}\,x_{4}}}\right\rangle]\mapsto[1, 2, 3, 4]^{\vee}],\medskip$
\\
$[[\left\langle {\frac {x_{3}^{2}}{x_{1}\,x_{2}}},{\frac {x_{2}}{x_{1}}}\right\rangle,\left\langle {\frac {x_{2}^{3}}{x_{0}\,x_{3}\,x_{4}}},{\frac {{x_{3}}^{2}}{x_{0}\,x_{4}}}\right\rangle]\mapsto[3, 5]^{\vee},[\left\langle {\frac {x_{3}^{2}}{x_{1}\,x_{2}}},{\frac {{x_{0}}^{2}}{x_{1}\,x_{2}}}\right\rangle,\left\langle {\frac {x_{3}^{2}}{x_{0}\,x_{4}}},{\frac {{x_{0}}^{2}}{x_{3}\,x_{4}}}\right\rangle]\mapsto[3, 4]^{\vee},$
\\
$[\left\langle {\frac {x_{3}^{2}}{x_{1}\,x_{2}}},{\frac {x_{1}}{x_{2}}}\right\rangle,\left\langle {\frac {x_{1}^{3}}{x_{0}\,x_{3}\,x_{4}}},{\frac {{x_{3}}^{2}}{x_{0}\,x_{4}}}\right\rangle]\mapsto[4, 6]^{\vee},[\left\langle {\frac {x_{4}^{2}}{x_{1}\,x_{2}}},{\frac {x_{2}}{x_{1}}}\right\rangle,\left\langle {\frac {x_{2}^{3}}{x_{0}\,x_{3}\,x_{4}}},{\frac {{x_{4}}^{2}}{x_{0}\,x_{3}}}\right\rangle]\mapsto[1, 5]^{\vee},$
\\
\end{longtable}
\end{center}

\noindent with 
2 sheets forming the complexes

\begin{center}
\[
\begin{tabular}
[c]{l}
$[],\medskip$
\\
$[],\medskip$
\\
$[\left\langle {\frac {x_{3}^{2}}{x_{1}\,x_{2}}}\right\rangle,\left\langle {\frac {x_{2}}{x_{1}}}\right\rangle,\left\langle {\frac {x_{1}}{x_{2}}}\right\rangle,\left\langle {\frac {x_{4}^{2}}{x_{1}\,x_{2}}}\right\rangle,\left\langle {\frac {x_{0}^{2}}{x_{1}\,x_{2}}}\right\rangle],\medskip$
\\
$[\left\langle {\frac {x_{3}^{2}}{x_{1}\,x_{2}}},{\frac {x_{2}}{x_{1}}}\right\rangle,\left\langle {\frac {x_{3}^{2}}{x_{1}\,x_{2}}},{\frac {{x_{0}}^{2}}{x_{1}\,x_{2}}}\right\rangle,\left\langle {\frac {x_{3}^{2}}{x_{1}\,x_{2}}},{\frac {x_{1}}{x_{2}}}\right\rangle,\left\langle {\frac {x_{4}^{2}}{x_{1}\,x_{2}}},{\frac {x_{2}}{x_{1}}}\right\rangle,\left\langle {\frac {x_{0}^{2}}{x_{1}\,x_{2}}},{\frac {x_{2}}{x_{1}}}\right\rangle,\left\langle {\frac {x_{4}^{2}}{x_{1}\,x_{2}}},{\frac {{x_{0}}^{2}}{x_{1}\,x_{2}}}\right\rangle,$
\\
$\left\langle {\frac {x_{4}^{2}}{x_{1}\,x_{2}}},{\frac {{x_{3}}^{2}}{x_{1}\,x_{2}}}\right\rangle,\left\langle {\frac {x_{4}^{2}}{x_{1}\,x_{2}}},{\frac {x_{1}}{x_{2}}}\right\rangle,\left\langle {\frac {x_{0}^{2}}{x_{1}\,x_{2}}},{\frac {x_{1}}{x_{2}}}\right\rangle],\medskip$
\\
$[\left\langle {\frac {x_{4}^{2}}{x_{1}\,x_{2}}},{\frac {{x_{0}}^{2}}{x_{1}\,x_{2}}},{\frac {x_{2}}{x_{1}}}\right\rangle,\left\langle {\frac {x_{4}^{2}}{x_{1}\,x_{2}}},{\frac {{x_{0}}^{2}}{x_{1}\,x_{2}}},{\frac {x_{1}}{x_{2}}}\right\rangle,\left\langle {\frac {x_{3}^{2}}{x_{1}\,x_{2}}},{\frac {{x_{0}}^{2}}{x_{1}\,x_{2}}},{\frac {x_{2}}{x_{1}}}\right\rangle,\left\langle {\frac {x_{3}^{2}}{x_{1}\,x_{2}}},{\frac {{x_{0}}^{2}}{x_{1}\,x_{2}}},{\frac {x_{1}}{x_{2}}}\right\rangle,$
\\
$\left\langle {\frac {x_{4}^{2}}{x_{1}\,x_{2}}},{\frac {{x_{3}}^{2}}{x_{1}\,x_{2}}},{\frac {x_{2}}{x_{1}}}\right\rangle,\left\langle {\frac {x_{4}^{2}}{x_{1}\,x_{2}}},{\frac {{x_{3}}^{2}}{x_{1}\,x_{2}}},{\frac {x_{1}}{x_{2}}}\right\rangle],\medskip$
\\
$[]$
\\
\end{tabular}
\]
\[
\begin{tabular}
[c]{l}
$[],\medskip$
\\
$[],\medskip$
\\
$[\left\langle {\frac {x_{3}^{2}}{x_{0}\,x_{4}}}\right\rangle,\left\langle {\frac {x_{2}^{3}}{x_{0}\,x_{3}\,x_{4}}}\right\rangle,\left\langle {\frac {x_{1}^{3}}{x_{0}\,x_{3}\,x_{4}}}\right\rangle,\left\langle {\frac {x_{4}^{2}}{x_{0}\,x_{3}}}\right\rangle,\left\langle {\frac {x_{0}^{2}}{x_{3}\,x_{4}}}\right\rangle],\medskip$
\\
$[\left\langle {\frac {x_{2}^{3}}{x_{0}\,x_{3}\,x_{4}}},{\frac {{x_{3}}^{2}}{x_{0}\,x_{4}}}\right\rangle,\left\langle {\frac {x_{3}^{2}}{x_{0}\,x_{4}}},{\frac {{x_{0}}^{2}}{x_{3}\,x_{4}}}\right\rangle,\left\langle {\frac {x_{1}^{3}}{x_{0}\,x_{3}\,x_{4}}},{\frac {{x_{3}}^{2}}{x_{0}\,x_{4}}}\right\rangle,\left\langle {\frac {x_{2}^{3}}{x_{0}\,x_{3}\,x_{4}}},{\frac {{x_{4}}^{2}}{x_{0}\,x_{3}}}\right\rangle,\left\langle {\frac {x_{2}^{3}}{x_{0}\,x_{3}\,x_{4}}},{\frac {{x_{0}}^{2}}{x_{3}\,x_{4}}}\right\rangle,$
\\
$\left\langle {\frac {x_{4}^{2}}{x_{0}\,x_{3}}},{\frac {{x_{0}}^{2}}{x_{3}\,x_{4}}}\right\rangle,\left\langle {\frac {x_{4}^{2}}{x_{0}\,x_{3}}},{\frac {{x_{3}}^{2}}{x_{0}\,x_{4}}}\right\rangle,\left\langle {\frac {x_{1}^{3}}{x_{0}\,x_{3}\,x_{4}}},{\frac {{x_{4}}^{2}}{x_{0}\,x_{3}}}\right\rangle,\left\langle {\frac {x_{1}^{3}}{x_{0}\,x_{3}\,x_{4}}},{\frac {{x_{0}}^{2}}{x_{3}\,x_{4}}}\right\rangle],\medskip$
\\
$[\left\langle {\frac {x_{2}^{3}}{x_{0}\,x_{3}\,x_{4}}},{\frac {{x_{4}}^{2}}{x_{0}\,x_{3}}},{\frac {{x_{0}}^{2}}{x_{3}\,x_{4}}}\right\rangle,\left\langle {\frac {x_{1}^{3}}{x_{0}\,x_{3}\,x_{4}}},{\frac {{x_{4}}^{2}}{x_{0}\,x_{3}}},{\frac {{x_{0}}^{2}}{x_{3}\,x_{4}}}\right\rangle,\left\langle {\frac {x_{2}^{3}}{x_{0}\,x_{3}\,x_{4}}},{\frac {{x_{3}}^{2}}{x_{0}\,x_{4}}},{\frac {{x_{0}}^{2}}{x_{3}\,x_{4}}}\right\rangle,$
\\
\end{tabular}
\]
\end{center}

\noindent Note that the torus invariant basis of deformations corresponding to
a Bergman face is given by the set of all lattice points of the polytope
specified above.

\paragraph{Limit map}

The limit map $\lim:B\left(  I\right)  \rightarrow\operatorname*{Poset}\left(
\Delta\right)  $ associates to a face $F$ of $B\left(  I\right)  $ the face of
$\Delta$ formed by the limit points of arcs lying over the weight vectors
$w\in F$, i.e. with lowest order term $t^{w}$.

Labeling the faces of the Bergman complex $B\left(  I\right)  \subset
\operatorname*{Poset}\left(  \nabla\right)  $ and the faces of
$\operatorname*{Poset}\left(  \Delta\right)  $ by the corresponding dual faces
of $\nabla^{\ast}$ and $\Delta^{\ast}$, hence considering the limit map
$\lim:B\left(  I\right)  \rightarrow\operatorname*{Poset}\left(
\Delta\right)  $ as a map $B\left(  I\right)  ^{\ast}\rightarrow
\operatorname*{Poset}\left(  \Delta^{\ast}\right)  $, the limit correspondence
is given by

\begin{center}
\begin{longtable}
[c]{l}
$[],\medskip$
\\
$[\left\langle y_{1},y_{2},y_{6},y_{7},y_{8},y_{3}\right\rangle\mapsto \left\langle x_{1},x_{3}\right\rangle,\left\langle y_{1},y_{2},y_{9},y_{7},y_{8},y_{4}\right\rangle\mapsto \left\langle x_{2},x_{3}\right\rangle,$
\\
$\left\langle y_{5},y_{2},y_{6},y_{10},y_{8},y_{3}\right\rangle\mapsto \left\langle x_{1},x_{4}\right\rangle,\left\langle y_{5},y_{2},y_{9},y_{10},y_{8},y_{4}\right\rangle\mapsto \left\langle x_{2},x_{4}\right\rangle,$
\\
$\left\langle y_{1},y_{5},y_{6},y_{7},y_{10},y_{3}\right\rangle\mapsto \left\langle x_{1},x_{0}\right\rangle,\left\langle y_{1},y_{5},y_{9},y_{7},y_{10},y_{4}\right\rangle\mapsto \left\langle x_{2},x_{0}\right\rangle],\medskip$
\\
$[\left\langle y_{5},y_{6},y_{10},y_{3}\right\rangle\mapsto \left\langle x_{1},x_{4},x_{0}\right\rangle,\left\langle y_{5},y_{2},y_{10},y_{8}\right\rangle\mapsto \left\langle x_{1},x_{2},x_{4}\right\rangle,$
\\
$\left\langle y_{5},y_{9},y_{10},y_{4}\right\rangle\mapsto \left\langle x_{2},x_{4},x_{0}\right\rangle,\left\langle y_{1},y_{6},y_{7},y_{3}\right\rangle\mapsto \left\langle x_{1},x_{3},x_{0}\right\rangle,$
\\
$\left\langle y_{2},y_{6},y_{8},y_{3}\right\rangle\mapsto \left\langle x_{1},x_{3},x_{4}\right\rangle,\left\langle y_{1},y_{2},y_{7},y_{8}\right\rangle\mapsto \left\langle x_{1},x_{2},x_{3}\right\rangle,$
\\
$\left\langle y_{1},y_{5},y_{7},y_{10}\right\rangle\mapsto \left\langle x_{1},x_{2},x_{0}\right\rangle,\left\langle y_{1},y_{9},y_{7},y_{4}\right\rangle\mapsto \left\langle x_{2},x_{3},x_{0}\right\rangle,$
\\
\end{longtable}
\end{center}

The image of the limit map coincides with the image of $\mu$ and with the
Bergman complex of the mirror, i.e. $\lim\left(  B\left(  I\right)  \right)
=\mu\left(  B\left(  I\right)  \right)  =B\left(  I^{\ast}\right)  $.

\paragraph{Mirror complex}

Numbering the vertices of the mirror complex $\mu\left(  B\left(  I\right)
\right)  $ as%
\[
\begin{tabular}
[c]{ll}
$1 = \left(4,-1,-1,-1\right)$ &$2 = \left(-1,4,-1,-1\right)$\\
$3 = \left(-1,-1,4,-1\right)$ &$4 = \left(-1,-1,-1,4\right)$\\
$5 = \left(-1,-1,-1,-1\right)$ &
\end{tabular}
\]
$\mu\left(  B\left(  I\right)  \right)  \subset\operatorname*{Poset}\left(
\Delta\right)  $ is

\begin{center}
\begin{longtable}
[c]{l}
$[],\medskip$
\\
$[[3],[2],[1],[4],[5]],\medskip$
\\
$[[2, 3],[3, 5],[1, 3],[2, 4],[2, 5],[4, 5],[3, 4],[1, 4],[1, 5]],\medskip$
\\
$[[2, 4, 5],[1, 4, 5],[2, 3, 5],[1, 3, 5],[2, 3, 4],[1, 3, 4]],\medskip$
\\
$[],\medskip$
\\
$[]$
\\
\end{longtable}
\end{center}

\noindent The ordering of the faces of $\mu\left(  B\left(  I\right)  \right)
$ is compatible with above ordering of the faces of $B\left(  I\right)  $. The
$F$-vectors of the faces of $\mu\left(  B\left(  I\right)  \right)  $ are%
\[
\begin{tabular}
[c]{rrrr}
Dimension &Number of faces &F-vector &\\
$1$ & $6$ & $\left(1,2,1,0,0,0\right)$ &edge\\
$2$ & $9$ & $\left(1,3,3,1,0,0\right)$ &triangle\\
$3$ & $5$ & $\left(1,4,6,4,1,0\right)$ &tetrahedron
\end{tabular}
\]
The first order deformations of the mirror special fiber $X_{0}^{\circ}$
correspond to the lattice points of the dual complex $\left(  \mu\left(
B\left(  I\right)  \right)  \right)  ^{\ast}\subset\operatorname*{Poset}%
\left(  \Delta^{\ast}\right)  $ of the mirror complex. We label the first
order deformations of $X_{0}^{\circ}$ corresponding to vertices of
$\Delta^{\ast}$ by the homogeneous coordinates of $Y$%
\[
\begin{tabular}
[c]{ll}
$x_{1} = x_{\left(1,0,0,0\right)}={\frac {y_{4}\,{y_{9}}^{3}}{y_{1}\,y_{2}\,y_{3}\,y_{5}}}$ &$x_{2} = x_{\left(0,1,0,0\right)}={\frac {y_{3}\,{y_{6}}^{3}}{y_{1}\,y_{2}\,y_{4}\,y_{5}}}$\\
$x_{3} = x_{\left(0,0,1,0\right)}={\frac {y_{5}^{2}{y_{10}}^{2}}{y_{6}\,y_{7}\,y_{8}\,y_{9}}}$ &$x_{4} = x_{\left(0,0,0,1\right)}={\frac {y_{1}^{2}{y_{7}}^{2}}{y_{6}\,y_{8}\,y_{9}\,y_{10}}}$\\
$x_{0} = x_{\left(-1,-1,-1,-1\right)}={\frac {y_{2}^{2}{y_{8}}^{2}}{y_{6}\,y_{7}\,y_{9}\,y_{10}}}$ &
\end{tabular}
\]
So writing the vertices of the faces of $\left(  \mu\left(  B\left(  I\right)
\right)  \right)  ^{\ast}$ as homogeneous coordinates of $Y$, the complex
$\left(  \mu\left(  B\left(  I\right)  \right)  \right)  ^{\ast}$ is given by

\begin{center}
\begin{longtable}
[c]{l}
$[],\medskip$
\\
$[],\medskip$
\\
$[\left\langle x_{1},x_{3}\right\rangle,\left\langle x_{2},x_{3}\right\rangle,\left\langle x_{1},x_{4}\right\rangle,\left\langle x_{2},x_{4}\right\rangle,\left\langle x_{1},x_{0}\right\rangle,\left\langle x_{2},x_{0}\right\rangle],\medskip$
\\
$[\left\langle x_{1},x_{4},x_{0}\right\rangle,\left\langle x_{1},x_{2},x_{4}\right\rangle,\left\langle x_{2},x_{4},x_{0}\right\rangle,\left\langle x_{1},x_{3},x_{0}\right\rangle,\left\langle x_{1},x_{3},x_{4}\right\rangle,\left\langle x_{1},x_{2},x_{3}\right\rangle,$
\\
$\left\langle x_{1},x_{2},x_{0}\right\rangle,\left\langle x_{2},x_{3},x_{0}\right\rangle,\left\langle x_{2},x_{3},x_{4}\right\rangle],\medskip$
\\
$[\left\langle x_{1},x_{2},x_{4},x_{0}\right\rangle,\left\langle x_{1},x_{3},x_{4},x_{0}\right\rangle,\left\langle x_{2},x_{3},x_{4},x_{0}\right\rangle,\left\langle x_{1},x_{2},x_{3},x_{0}\right\rangle,\left\langle x_{1},x_{2},x_{3},x_{4}\right\rangle],\medskip$
\\
$[]$
\\
\end{longtable}
\end{center}

\noindent The complex $\mu\left(  B\left(  I\right)  \right)  ^{\ast}$ labeled
by the variables of the Cox ring $S$ gives the ideals of the toric strata of
the special fiber $X_{0}$ of $\mathfrak{X}$, i.e. the complex $SP\left(
I_{0}\right)  $, so in particular the primary decomposition of $I_{0}$ is%
\[
I_{0}=
\begin{tabular}
[c]{l}
$\left\langle x_{1},x_{3}\right\rangle \cap\left\langle x_{1},x_{0}\right\rangle \cap\left\langle x_{2},x_{3}\right\rangle \cap\left\langle x_{2},x_{0}\right\rangle \cap\left\langle x_{1},x_{4}\right\rangle \cap\left\langle x_{2},x_{4}\right\rangle$
\\
\end{tabular}
\]

Labeling the vertices of the faces by the corresponding deformations the
complex is given by

\begin{center}
%
\end{center}

\paragraph{Covering structure in the deformation complex of the mirror
degeneration}

Each face of the complex of deformations $\left(  \mu\left(  B\left(
I\right)  \right)  \right)  ^{\ast}$ of the mirror special fiber $X_{0}%
^{\circ}$ decomposes as the convex hull of 
$2$
polytopes forming
a $
2
:1$ 
trivial
covering of $\mu\left(  B\left(
I\right)  \right)  ^{\vee}$

\begin{center}
%
\end{center}

\noindent 
Due to the singularities of
$Y^{\circ}$ this covering involves degenerate faces, i.e. faces $G\mapsto
F^{\vee}$ with $\dim\left(  G\right)  <\dim\left(  F^{\vee}\right)  $.

Writing the vertices of the faces as deformations of $X_{0}^{\circ}$ the
covering is given by

\begin{center}
%
\end{center}

\paragraph{Mirror degeneration}

The space of first order deformations of $X_{0}^{\circ}$ in the mirror
degeneration $\mathfrak{X}^{\circ}$ has dimension $
5
$ and the
deformations represented by the monomials%
\[
\left\{  
\begin{tabular}
[c]{lllll}
${\frac {y_{5}^{2}{y_{10}}^{2}}{y_{6}\,y_{7}\,y_{8}\,y_{9}}}$ & ${\frac {y_{4}\,{y_{9}}^{3}}{y_{1}\,y_{2}\,y_{3}\,y_{5}}}$ & ${\frac {y_{3}\,{y_{6}}^{3}}{y_{1}\,y_{2}\,y_{4}\,y_{5}}}$ & ${\frac {y_{1}^{2}{y_{7}}^{2}}{y_{6}\,y_{8}\,y_{9}\,y_{10}}}$ & ${\frac {y_{2}^{2}{y_{8}}^{2}}{y_{6}\,y_{7}\,y_{9}\,y_{10}}}$
\end{tabular}
\right\}
\]
form a torus invariant
basis.

The 
mirror degeneration 
$\mathfrak{X}^{\circ}\subset Y^{\circ}\times\operatorname*{Spec}%
\mathbb{C}\left[  t\right]  $
 of $\mathfrak{X}$ is
given by the ideal 
$I^{\circ}\subset S^{\circ}\otimes\mathbb{C}\left[  t\right]  $
generated by%
\[
%
\]
Indeed already the ideal $J^{\circ}$ generated by%
\[
\left\{  
%
\right\}
\]
defines $\mathfrak{X}^{\circ}$.

\paragraph{Contraction of the mirror degeneration}

In the following we give a birational map relating the degeneration
$\mathfrak{X}^{\circ}$ to a Greene-Plesser type orbifolding mirror family by
contracting divisors on $Y^{\circ}$.
See also Section \ref{sec orbifolding mirror families} below.

In order to contract the divisors

\begin{center}
%
\end{center}

\noindent consider the $\mathbb{Q}$-factorial toric Fano variety $\hat
{Y}^{\circ}=X\left(  \hat{\Sigma}^{\circ}\right)  $, where $\hat{\Sigma
}^{\circ}=\operatorname*{Fan}\left(  \hat{P}^{\circ}\right)  \subset
M_{\mathbb{R}}$ is the fan over the 
Fano polytope $\hat
{P}^{\circ}\subset M_{\mathbb{R}}$ given as the convex hull of the remaining
vertices of $P^{\circ}=\nabla^{\ast}$ corresponding to the Cox variables%
\[
%
\]
of the toric variety $\hat{Y}^{\circ}$ with Cox ring
\[
\hat{S}^{\circ}=\mathbb{C}
[y_{5},y_{6},y_{9},y_{1},y_{2}]
\]
The Cox variables of $\hat{Y}^{\circ}$ correspond to the set of Fermat deformations of $\mathfrak{X}$.

Let%
\[
Y^{\circ}=X\left(  \Sigma^{\circ}\right)  \rightarrow X\left(  \hat{\Sigma
}^{\circ}\right)  =\hat{Y}^{\circ}%
\]
be a birational map from $Y^{\circ}$ to a minimal birational model $\hat
{Y}^{\circ}$, which contracts the divisors of the rays $\Sigma^{\circ}\left(
1\right)  -\hat{\Sigma}^{\circ}\left(  1\right)  $ corresponding to the Cox
variables
\[
\begin{tabular}
[c]{lllll}
$y_{3}$ &$y_{4}$ &$y_{7}$ &$y_{8}$ &$y_{10}$
\end{tabular}
\]

Representing $\hat{Y}^{\circ}$ as a quotient we have
\[
\hat{Y}^{\circ}=\left(  \mathbb{C}^{
5
}-V\left(  B\left(
\hat{\Sigma}^{\circ}\right)  \right)  \right)  //\hat{G}^{\circ}%
\]
with%
\[
\hat{G}^{\circ}=
\mathbb{Z}_{2}\times\mathbb{Z}_{6}\times\left(  \mathbb{C}^{\ast}\right)  ^{1}
\]
acting via%
\[
\xi y=
\left( \,u_{1}\,v_{1} \cdot y_{5},\,u_{1}\,u_{2}\,v_{1} \cdot y_{6},\,u_{1}\,u_{2}^{5}\,v_{1} \cdot y_{9},\,u_{2}^{3}\,v_{1} \cdot y_{1},\,v_{1} \cdot y_{2} \right)
\]
for $\xi=
\left(u_1,u_2,v_1\right)
\in\hat{G}^{\circ}$ and $y\in\mathbb{C}%
^{
5
}-V\left(  B\left(  \hat{\Sigma}^{\circ}\right)  \right)  $.

Hence with the group%
\[
\hat{H}^{\circ}=
\mathbb{Z}_{2}\times\mathbb{Z}_{6}
\]
of order 
12
 the toric variety $\hat{Y}^{\circ}$ is the quotient%
\[
\hat{Y}^{\circ}=\mathbb{P}^{
4
}/\hat{H}^{\circ}%
\]
of projective space $\mathbb{P}^{
4
}$.

The mirror degeneration $\mathfrak{X}^{\circ}$ induces via $Y\rightarrow
\hat{Y}^{\circ}$ a degeneration $\mathfrak{\hat{X}}^{\circ}\subset\hat
{Y}^{\circ}\times\operatorname*{Spec}\mathbb{C}\left[  t\right]  $ given
by the ideal $\hat{I}^{\circ}\subset\hat{S}^{\circ}%
\otimes\mathbb{C}\left[  t\right]  $ generated by the
 Fermat-type equations
\[
\left\{  
\begin{tabular}
[c]{l}
$y_{1}\,y_{2}\,y_{5}+t (s_{2}\,y_{6}^{3}+s_{1}\,y_{9}^{3} ),\medskip$
\\
$y_{6}\,y_{9}+t (s_{3}\,y_{1}^{2}+s_{5}\,y_{2}^{2}+s_{4}\,y_{5}^{2} )$
\\
\end{tabular}
\right\}
\]
The special fiber $\hat{X}_{0}^{\circ}\subset\hat{Y}^{\circ}$ of
$\mathfrak{\hat{X}}^{\circ}$ is cut out by the monomial ideal $\hat{I}%
_{0}^{\circ}\subset\hat{S}^{\circ}$ generated by%
\[
\left\{  
%
\right\}
\]

The complex $SP\left(  \hat{I}_{0}^{\circ}\right)  $ of prime ideals of the
toric strata of $\hat{X}_{0}^{\circ}$ is

\begin{center}
%
\end{center}

\noindent so $\hat{I}_{0}^{\circ}$ has the primary decomposition%
\[
\hat{I}_{0}^{\circ}=
%
\]
The Bergman subcomplex $B\left(  I\right)  \subset\operatorname*{Poset}\left(
\nabla\right)  $ induces a subcomplex of $\operatorname*{Poset}\left(
\hat{\nabla}\right)  $, $\hat{\nabla}=\left(  \hat{P}^{\circ}\right)  ^{\ast}$
corresponding to $\hat{I}_{0}^{\circ}$. Indexing of the 
vertices
of $\hat{\nabla}$ by%
\[
%
\]
this complex is given by
\[
\left\{  
%
\right\}
\]

\subsubsection{The Calabi-Yau threefold given as the complete intersection of
a~generic quadric and a generic quartic in $\mathbb{P}^{5}$}

\paragraph{Setup}

Let $Y=\mathbb{P}^{
5
}=X\left(  \Sigma\right)  $, $\Sigma
=\operatorname*{Fan}\left(  P\right)  =NF\left(  \Delta\right)  \subset
N_{\mathbb{R}}$ with the Fano polytope $P=\Delta^{\ast}$ given by%
\[
\Delta=\operatorname*{convexhull}\left(  
%
\right)  \subset
M_{\mathbb{R}}%
\]
and let%
\[
S=\mathbb{C}
[x_0, x_1, x_2, x_3, x_4, x_5]
\]
be the Cox ring of $Y$ with the variables%
\[
%
\]
associated to the rays of $\Sigma$.

Consider the degeneration $\mathfrak{X}\subset Y\times\operatorname*{Spec}%
\mathbb{C}\left[  t\right]  $ of 
complete intersection
Calabi-Yau 3-folds
 of type $\left(
2,4
\right)  $ with monomial special fiber
\[
I_{0}=\left\langle 
%
\right\rangle
\]
The degeneration is given by the ideal $I\subset S\otimes\mathbb{C}\left[
t\right]  $ with $I_{0}$-reduced generators%
\[
\left\{  
%
\right\}
\]

\paragraph{Special fiber Gr\"{o}bner cone}

The space of first order deformations of $\mathfrak{X}$ has dimension
$
124
$ and the deformations represented by the Cox Laurent monomials

\begin{center}
%
\end{center}

\noindent form a torus invariant basis. These deformations give linear
inequalities defining the special fiber Gr\"{o}bner cone, i.e.,%
\[
C_{I_{0}}\left(  I\right)  =\left\{  \left(  w_{t},w_{y}\right)  \in
\mathbb{R}\oplus N_{\mathbb{R}}\mid\left\langle w,A^{-1}\left(  m\right)
\right\rangle \geq-w_{t}\forall m\right\}
\]
for above monomials $m$ and the presentation matrix
\[
A=
\left [%
\right ]
\]
of $A_{
4
}\left(  Y\right)  $.
The vertices of special fiber polytope
\[
\nabla=C_{I_{0}}\left(  I\right)  \cap\left\{  w_{t}=1\right\}
\]
and the number of faces, each vertex is contained in, are given by the table
\begin{center}
%
\end{center}
\noindent The number of faces of $\nabla$ and their $F$-vectors are
\begin{center}
%
\end{center}
The dual $P^{\circ}=\nabla^{\ast}$ of $\nabla$ is a Fano polytope with vertices
\begin{center}
%
\end{center}
\noindent and the $F$-vectors of the faces of $P^{\circ}$ are

\begin{center}
%
\end{center}

The polytopes $\nabla$ and $P^{\circ}$ have $0$ as the unique interior lattice
point, and the Fano polytope $P^{\circ}$ defines the embedding toric Fano
variety $Y^{\circ}=X\left(  \operatorname*{Fan}\left(  P^{\circ}\right)
\right)  $ of the fibers of the mirror degeneration. The dimension of the
automorphism group of $Y^{\circ}$ is
\[
\dim\left(  \operatorname*{Aut}\left(  Y^{\circ}\right)  \right)
=
5
\]
Let%
\[
S^{\circ}=\mathbb{C}
[y_{1},y_{2},...,y_{12}]
\]
be the Cox ring of $Y^{\circ}$ with variables

\begin{center}
%
\end{center}

\paragraph{Bergman subcomplex}

Intersecting the tropical variety of $I$ with the special fiber Gr\"{o}bner
cone $C_{I_{0}}\left(  I\right)  $ we obtain the special fiber Bergman
subcomplex
\[
B\left(  I\right)  =C_{I_{0}}\left(  I\right)  \cap BF\left(  I\right)
\cap\left\{  w_{t}=1\right\}
\]
The following table chooses an indexing of the vertices of $\nabla$ involved
in $B\left(  I\right)  $ and gives for each vertex the numbers $n_{\nabla}$
and $n_{B}$ of faces of $\nabla$ and faces of $B\left(  I\right)  $ it is
contained in

\begin{center}
%
\end{center}

\noindent With this indexing the Bergman subcomplex $B\left(  I\right)  $ of
$\operatorname*{Poset}\left(  \nabla\right)  $ associated to the degeneration
$\mathfrak{X}$ is

\begin{center}
%
\]
and the $F$-vectors of its faces are

\begin{center}
%
\end{center}

\paragraph{Dual complex}

The dual complex $\operatorname*{dual}\left(  B\left(  I\right)  \right)
=\left(  B\left(  I\right)  \right)  ^{\ast}$ of deformations associated to
$B\left(  I\right)  $ via initial ideals is given by

\begin{center}
%
\end{center}

\noindent when writing the vertices of the faces as deformations of $X_{0}$.
Note that the $T$-invariant basis of deformations associated to a face is
given by all lattice points of the corresponding polytope in $M_{\mathbb{R}}$.

In order to compress the output we list one representative in any set of faces
$G$ with fixed $F$-vector of $G$ and $G^{\ast}$.

When numbering the vertices of the faces of $\operatorname*{dual}\left(
B\left(  I\right)  \right)  $ by the Cox variables of the mirror toric Fano
variety $Y^{\circ}$ the complex $\operatorname*{dual}\left(  B\left(
I\right)  \right)  $ is
\begin{center}
\begin{longtable}
[c]{l}
$[],\medskip$
\\
$[],\medskip$
\\
$[[2, 4, 6, 8]^{*}=\left\langle y_{11},y_{5}\right\rangle,[3, 4, 5, 6, 7, 8]^{*}=\left\langle y_{6},y_{12}\right\rangle,...],\medskip$
\\
$[[3, 4, 7, 8]^{*}=\left\langle y_{6},y_{2},y_{12},y_{9}\right\rangle,[1, 3, 7]^{*}=\left\langle y_{2},y_{7},y_{9},y_{4}\right\rangle,...],\medskip$
\\
$[[5, 7]^{*}=\left\langle y_{1},y_{6},y_{7},y_{8},y_{12},y_{4}\right\rangle,...],\medskip$
\\
$[[1]^{*}=\left\langle y_{1},y_{2},y_{3},y_{7},y_{8},y_{9},y_{10},y_{4}\right\rangle,...],\medskip$
\\
$[]$
\\
\end{longtable}
\end{center}
\noindent The dual complex has the $F$-vector%
\[
\begin{tabular}
[c]{rrrrrrrr}
Dimension &$-1$ & $0$ & $1$ & $2$ & $3$ & $4$ & $5$\\
Number of faces &$0$ & $0$ & $6$ & $14$ & $16$ & $8$ & $0$
\end{tabular}
\]
and the $F$-vectors of the faces of $\operatorname*{dual}\left(  B\left(
I\right)  \right)  $ are
\begin{center}
\begin{longtable}
[c]{rrll}
Dimension &Number of faces &F-vector &\\
$1$ & $6$ & $\left(1,2,1,0,0,0,0\right)$ &edge\\
$2$ & $14$ & $\left(1,4,4,1,0,0,0\right)$ &quadrangle\\
$3$ & $16$ & $\left(1,6,9,5,1,0,0\right)$ &prism\\
$4$ & $8$ & $\left(1,8,16,14,6,1,0\right)$ &
\end{longtable}
\end{center}

Recall that in this example the toric variety $Y$ is projective space. The
number of lattice points of the support of $\operatorname*{dual}\left(
B\left(  I\right)  \right)  $ relates to the dimension $h^{1,
2
}%
\left(  X\right)  $ of the complex moduli space of the generic fiber $X$ of
$\mathfrak{X}$ and to the dimension $h^{1,1}\left(  \bar{X}^{\circ}\right)  $
of the K\"{a}hler moduli space of the $MPCR$-blowup $\bar{X}^{\circ}$ of the
generic fiber $X^{\circ}$ of the mirror degeneration
\begin{align*}
\left\vert \operatorname*{supp}\left(  \operatorname*{dual}\left(  B\left(
I\right)  \right)  \right)  \cap M\right\vert  &
=
124
=
35
+
89
=\dim\left(  \operatorname*{Aut}%
\left(  Y\right)  \right)  +h^{1,
2
}\left(  X\right) \\
&  =
30
+
5
+
89
\\
&  =\left\vert \operatorname*{Roots}\left(  Y\right)  \right\vert +\dim\left(
T_{Y}\right)  +h^{1,1}\left(  \bar{X}^{\circ}\right)
\end{align*}
There are%
\[
h^{1,
2
}\left(  X\right)  +\dim\left(  T_{Y^{\circ}}\right)
=
89
+
5
\]
non-trivial toric polynomial deformations of $X_{0}$

\begin{center}

\end{center}

\noindent They correspond to the toric divisors

\begin{center}
%
\end{center}

\noindent on a MPCR-blowup of $Y^{\circ}$ inducing $
89
$ non-zero
toric divisor classes on the mirror $X^{\circ}$. The following $
30
$
toric divisors of $Y$ induce the trivial divisor class on $X^{\circ}$

\begin{center}
%
\end{center}

\paragraph{Mirror special fiber}

The ideal $I_{0}^{\circ}$ of the monomial special fiber $X_{0}^{\circ}$ of the
mirror degeneration $\mathfrak{X}^{\circ}$ is generated by the following set
of monomials in $S^{\circ}$%
\[
\left\{  
%
\right\}
\]
Indeed already the ideal%
\[
J_{0}^{\circ}=\left\langle 
%
\right\rangle
\]
defines the same subvariety of the toric variety $Y^{\circ}$, and
$J_{0}^{\circ\Sigma}=I_{0}^{\circ}$. Recall that passing from $J_{0}^{\circ}$
to $J_{0}^{\circ\Sigma}$ is the non-simplicial toric analogue of saturation.

The complex $B\left(  I\right)  ^{\ast}$ labeled by the variables of the Cox
ring $S^{\circ}$ of $Y^{\ast}$, as written in the last section, is the complex
$SP\left(  I_{0}^{\circ}\right)  $ of prime ideals of the toric strata of the
special fiber $X_{0}^{\circ}$ of the mirror degeneration $\mathfrak{X}^{\circ
}$, i.e., the primary decomposition of $I_{0}^{\circ}$ is%
\[
I_{0}^{\circ}=

\]
Each facet $F\in B\left(  I\right)  $ corresponds to one of these ideals and this
ideal is generated by the products of facets of $\nabla$ containing $F$.

\paragraph{Covering structure in the deformation complex of the
degeneration $\mathfrak{X}$}

According to the local reduced Gr\"{o}bner basis each face of the complex of
deformations $\operatorname*{dual}\left(  B\left(  I\right)  \right)  $
decomposes into 
$2$
polytopes forming a $
2
:1$
trivial
covering of $B\left(  I\right)  $

\begin{center}
%
\end{center}

\noindent Here the faces $F\in B\left(  I\right)  $ are specified both via the
vertices of $F^{\ast}$ labeled by the variables of $S^{\circ}$ and by the
numbering of the vertices of $B\left(  I\right)  $ chosen above.

\noindent This covering has 
2 sheets forming the complexes
 
\[
%
\end{center}

\noindent Writing the vertices of the faces as deformations the covering is
given by

\begin{center}
%
\end{center}

\noindent with 
2 sheets forming the complexes

\begin{center}
\[
%
\]
\end{center}

\noindent Note that the torus invariant basis of deformations corresponding to
a Bergman face is given by the set of all lattice points of the polytope
specified above.

\paragraph{Limit map}

The limit map $\lim:B\left(  I\right)  \rightarrow\operatorname*{Poset}\left(
\Delta\right)  $ associates to a face $F$ of $B\left(  I\right)  $ the face of
$\Delta$ formed by the limit points of arcs lying over the weight vectors
$w\in F$, i.e. with lowest order term $t^{w}$.

Labeling the faces of the Bergman complex $B\left(  I\right)  \subset
\operatorname*{Poset}\left(  \nabla\right)  $ and the faces of
$\operatorname*{Poset}\left(  \Delta\right)  $ by the corresponding dual faces
of $\nabla^{\ast}$ and $\Delta^{\ast}$, hence considering the limit map
$\lim:B\left(  I\right)  \rightarrow\operatorname*{Poset}\left(
\Delta\right)  $ as a map $B\left(  I\right)  ^{\ast}\rightarrow
\operatorname*{Poset}\left(  \Delta^{\ast}\right)  $, the limit correspondence
is given by

\begin{center}

\end{center}

The image of the limit map coincides with the image of $\mu$ and with the
Bergman complex of the mirror, i.e. $\lim\left(  B\left(  I\right)  \right)
=\mu\left(  B\left(  I\right)  \right)  =B\left(  I^{\ast}\right)  $.

\paragraph{Mirror complex}

Numbering the vertices of the mirror complex $\mu\left(  B\left(  I\right)
\right)  $ as%
\[
%
\end{center}

\noindent The ordering of the faces of $\mu\left(  B\left(  I\right)  \right)
$ is compatible with above ordering of the faces of $B\left(  I\right)  $. The
$F$-vectors of the faces of $\mu\left(  B\left(  I\right)  \right)  $ are%
\[
%
\]
The first order deformations of the mirror special fiber $X_{0}^{\circ}$
correspond to the lattice points of the dual complex $\left(  \mu\left(
B\left(  I\right)  \right)  \right)  ^{\ast}\subset\operatorname*{Poset}%
\left(  \Delta^{\ast}\right)  $ of the mirror complex. We label the first
order deformations of $X_{0}^{\circ}$ corresponding to vertices of
$\Delta^{\ast}$ by the homogeneous coordinates of $Y$%
\[
%
\]
So writing the vertices of the faces of $\left(  \mu\left(  B\left(  I\right)
\right)  \right)  ^{\ast}$ as homogeneous coordinates of $Y$, the complex
$\left(  \mu\left(  B\left(  I\right)  \right)  \right)  ^{\ast}$ is given by

\begin{center}
%
\end{center}

\noindent The complex $\mu\left(  B\left(  I\right)  \right)  ^{\ast}$ labeled
by the variables of the Cox ring $S$ gives the ideals of the toric strata of
the special fiber $X_{0}$ of $\mathfrak{X}$, i.e. the complex $SP\left(
I_{0}\right)  $, so in particular the primary decomposition of $I_{0}$ is%
\[
I_{0}=
%
\]

\paragraph{Covering structure in the deformation complex of the mirror
degeneration}

Each face of the complex of deformations $\left(  \mu\left(  B\left(
I\right)  \right)  \right)  ^{\ast}$ of the mirror special fiber $X_{0}%
^{\circ}$ decomposes as the convex hull of 
$2$
polytopes forming
a $
2
:1$ 
trivial
covering of $\mu\left(  B\left(
I\right)  \right)  ^{\vee}$

\begin{center}
%
\end{center}

\noindent 
Due to the singularities of
$Y^{\circ}$ this covering involves degenerate faces, i.e. faces $G\mapsto
F^{\vee}$ with $\dim\left(  G\right)  <\dim\left(  F^{\vee}\right)  $.

\paragraph{Mirror degeneration}

The space of first order deformations of $X_{0}^{\circ}$ in the mirror
degeneration $\mathfrak{X}^{\circ}$ has dimension $
6
$ and the
deformations represented by the monomials%
\[
\left\{  
%
\right\}
\]
form a torus invariant
basis.
The number of lattice points of the dual of the mirror complex of $I$ relates
to the dimension $h^{1,
2
}\left(  X^{\circ}\right)  $ of complex
moduli space of the generic fiber $X^{\circ}$ of $\mathfrak{X}^{\circ}$ and to
the dimension $h^{1,1}\left(  X\right)  $ of the K\"{a}hler moduli space of
the generic fiber $X$ of $\mathfrak{X}$ via%
\begin{align*}
\left\vert \operatorname*{supp}\left(  \left(  \mu\left(  B\left(  I\right)
\right)  \right)  ^{\ast}\right)  \cap N\right\vert  &
=
6
=
5
+
1
\\
&  =\dim\left(  \operatorname*{Aut}\left(  Y^{\circ}\right)  \right)
+h^{1,
2
}\left(  X^{\circ}\right)  =\dim\left(  T\right)
+h^{1,1}\left(  X\right)
\end{align*}

The 
mirror degeneration 
$\mathfrak{X}^{\circ}\subset Y^{\circ}\times\operatorname*{Spec}%
\mathbb{C}\left[  t\right]  $
 of $\mathfrak{X}$ is
given by the ideal 
$I^{\circ}\subset S^{\circ}\otimes\mathbb{C}\left[  t\right]  $
generated by%
\[

\]
Indeed already the ideal $J^{\circ}$ generated by%
\[
\left\{  
%
\right\}
\]
defines $\mathfrak{X}^{\circ}$.

\paragraph{Contraction of the mirror degeneration}

In the following we give a birational map relating the degeneration
$\mathfrak{X}^{\circ}$ to a Greene-Plesser type orbifolding mirror family by
contracting divisors on $Y^{\circ}$.
See also Section \ref{sec orbifolding mirror families} below.

In order to contract the divisors

\begin{center}
%
\end{center}

\noindent consider the $\mathbb{Q}$-factorial toric Fano variety $\hat
{Y}^{\circ}=X\left(  \hat{\Sigma}^{\circ}\right)  $, where $\hat{\Sigma
}^{\circ}=\operatorname*{Fan}\left(  \hat{P}^{\circ}\right)  \subset
M_{\mathbb{R}}$ is the fan over the 
Fano polytope $\hat
{P}^{\circ}\subset M_{\mathbb{R}}$ given as the convex hull of the remaining
vertices of $P^{\circ}=\nabla^{\ast}$ corresponding to the Cox variables%
\[
%
\]
of the toric variety $\hat{Y}^{\circ}$ with Cox ring
\[
\hat{S}^{\circ}=\mathbb{C}
[y_{6},y_{11},y_{2},y_{1},y_{7},y_{3}]
\]
The Cox variables of $\hat{Y}^{\circ}$ correspond to the set of Fermat deformations of $\mathfrak{X}$.

Let%
\[
Y^{\circ}=X\left(  \Sigma^{\circ}\right)  \rightarrow X\left(  \hat{\Sigma
}^{\circ}\right)  =\hat{Y}^{\circ}%
\]
be a birational map from $Y^{\circ}$ to a minimal birational model $\hat
{Y}^{\circ}$, which contracts the divisors of the rays $\Sigma^{\circ}\left(
1\right)  -\hat{\Sigma}^{\circ}\left(  1\right)  $ corresponding to the Cox
variables
\[
\begin{tabular}
[c]{llllll}
$y_{4}$ &$y_{5}$ &$y_{8}$ &$y_{9}$ &$y_{10}$ &$y_{12}$
\end{tabular}
\]

Representing $\hat{Y}^{\circ}$ as a quotient we have
\[
\hat{Y}^{\circ}=\left(  \mathbb{C}^{
6
}-V\left(  B\left(
\hat{\Sigma}^{\circ}\right)  \right)  \right)  //\hat{G}^{\circ}%
\]
with%
\[
\hat{G}^{\circ}=
\mathbb{Z}_{2}\times\mathbb{Z}_{2}\times\mathbb{Z}_{8}\times\left(  \mathbb{C}^{\ast}\right)  ^{1}
\]
acting via%
\[
\xi y=
\left( \,u_{1}\,v_{1} \cdot y_{6},\,u_{1}\,u_{3}^{7}\,v_{1} \cdot y_{11},\,u_{2}\,v_{1} \cdot y_{2},\,u_{1}\,u_{2}\,u_{3}^{4}\,v_{1} \cdot y_{1},\,u_{1}\,u_{3}\,v_{1} \cdot y_{7},\,v_{1} \cdot y_{3} \right)
\]
for $\xi=
\left(u_1,u_2,u_3,v_1\right)
\in\hat{G}^{\circ}$ and $y\in\mathbb{C}%
^{
6
}-V\left(  B\left(  \hat{\Sigma}^{\circ}\right)  \right)  $.

Hence with the group%
\[
\hat{H}^{\circ}=
\mathbb{Z}_{2}\times\mathbb{Z}_{2}\times\mathbb{Z}_{8}
\]
of order 
32
 the toric variety $\hat{Y}^{\circ}$ is the quotient%
\[
\hat{Y}^{\circ}=\mathbb{P}^{
5
}/\hat{H}^{\circ}%
\]
of projective space $\mathbb{P}^{
5
}$.

The first order mirror degeneration $\mathfrak{X}^{1\circ}$ induces via
$Y\rightarrow\hat{Y}^{\circ}$ a degeneration $\mathfrak{\hat{X}}^{1\circ
}\subset\hat{Y}^{\circ}\times\operatorname*{Spec}\mathbb{C}\left[  t\right]
/\left\langle t^{2}\right\rangle $ given by the ideal $\hat{I}^{1\circ}%
\subset\hat{S}^{\circ}\otimes\mathbb{C}\left[  t\right]  /\left\langle
t^{2}\right\rangle $ generated by the
 Fermat-type equations
\[
\left\{  

\right\}
\]

Note, that%
\[
\left\vert \operatorname*{supp}\left(  \left(  \mu\left(  B\left(  I\right)
\right)  \right)  ^{\ast}\right)  \cap N-\operatorname*{Roots}\left(  \hat
{Y}^{\circ}\right)  \right\vert -\dim\left(  T_{\hat{Y}^{\circ}}\right)
=
6
-
5
=
1
\]
so this family has one independent parameter.

The special fiber $\hat{X}_{0}^{\circ}\subset\hat{Y}^{\circ}$ of
$\mathfrak{\hat{X}}^{1\circ}$ is cut out by the monomial ideal $\hat{I}%
_{0}^{\circ}\subset\hat{S}^{\circ}$ generated by%
\[
\left\{  
%
\right\}
\]

The complex $SP\left(  \hat{I}_{0}^{\circ}\right)  $ of prime ideals of the
toric strata of $\hat{X}_{0}^{\circ}$ is

\begin{center}
%
\end{center}

\noindent so $\hat{I}_{0}^{\circ}$ has the primary decomposition%
\[
\hat{I}_{0}^{\circ}=
%
\]
The Bergman subcomplex $B\left(  I\right)  \subset\operatorname*{Poset}\left(
\nabla\right)  $ induces a subcomplex of $\operatorname*{Poset}\left(
\hat{\nabla}\right)  $, $\hat{\nabla}=\left(  \hat{P}^{\circ}\right)  ^{\ast}$
corresponding to $\hat{I}_{0}^{\circ}$. Indexing of the 
vertices
of $\hat{\nabla}$ by%
\[
%
\]
this complex is given by
\[
\left\{  
%
\right\}
\]

\subsubsection{The Calabi-Yau threefold given as the complete intersection of
two~generic cubics in $\mathbb{P}^{5}$\label{Sec Ex 33}}

\paragraph{Setup}

Let $Y=\mathbb{P}^{
5
}=X\left(  \Sigma\right)  $, $\Sigma
=\operatorname*{Fan}\left(  P\right)  =NF\left(  \Delta\right)  \subset
N_{\mathbb{R}}$ with the Fano polytope $P=\Delta^{\ast}$ given by%
\[
\Delta=\operatorname*{convexhull}\left(  
%
\right)  \subset
M_{\mathbb{R}}%
\]
and let%
\[
S=\mathbb{C}
[x_0, x_1, x_2, x_3, x_4, x_5]
\]
be the Cox ring of $Y$ with the variables%
\[
%
\]
associated to the rays of $\Sigma$.

Consider the degeneration $\mathfrak{X}\subset Y\times\operatorname*{Spec}%
\mathbb{C}\left[  t\right]  $ of 
complete intersection
Calabi-Yau 3-folds
 of type $\left(
3,3
\right)  $ with monomial special fiber
\[
I_{0}=\left\langle 
%
\right\rangle
\]
The degeneration is given by the ideal $I\subset S\otimes\mathbb{C}\left[
t\right]  $ with $I_{0}$-reduced generators%
\[
\left\{  
%
\right\}
\]

\paragraph{Special fiber Gr\"{o}bner cone}

The space of first order deformations of $\mathfrak{X}$ has dimension
$
108
$ and the deformations represented by the Cox Laurent monomials

\begin{center}
%
\end{center}

\noindent form a torus invariant basis. These deformations give linear
inequalities defining the special fiber Gr\"{o}bner cone, i.e.,%
\[
C_{I_{0}}\left(  I\right)  =\left\{  \left(  w_{t},w_{y}\right)  \in
\mathbb{R}\oplus N_{\mathbb{R}}\mid\left\langle w,A^{-1}\left(  m\right)
\right\rangle \geq-w_{t}\forall m\right\}
\]
for above monomials $m$ and the presentation matrix
\[
A=
\left [%
\right ]
\]
of $A_{
4
}\left(  Y\right)  $.

The vertices of special fiber polytope
\[
\nabla=C_{I_{0}}\left(  I\right)  \cap\left\{  w_{t}=1\right\}
\]
and the number of faces, each vertex is contained in, are given by the table

\begin{center}
%
\end{center}

\noindent The number of faces of $\nabla$ and their $F$-vectors are

\begin{center}
%
\end{center}

The dual $P^{\circ}=\nabla^{\ast}$ of $\nabla$ is a Fano polytope with vertices

\begin{center}
%
\end{center}

\noindent and the $F$-vectors of the faces of $P^{\circ}$ are

\begin{center}
%
\end{center}

The polytopes $\nabla$ and $P^{\circ}$ have $0$ as the unique interior lattice
point, and the Fano polytope $P^{\circ}$ defines the embedding toric Fano
variety $Y^{\circ}=X\left(  \operatorname*{Fan}\left(  P^{\circ}\right)
\right)  $ of the fibers of the mirror degeneration. The dimension of the
automorphism group of $Y^{\circ}$ is
\[
\dim\left(  \operatorname*{Aut}\left(  Y^{\circ}\right)  \right)
=
5
\]
Let%
\[
S^{\circ}=\mathbb{C}
[y_{1},y_{2},...,y_{12}]
\]
be the Cox ring of $Y^{\circ}$ with variables

\begin{center}
%
\end{center}

\paragraph{Bergman subcomplex}

Intersecting the tropical variety of $I$ with the special fiber Gr\"{o}bner
cone $C_{I_{0}}\left(  I\right)  $ we obtain the special fiber Bergman
subcomplex
\[
B\left(  I\right)  =C_{I_{0}}\left(  I\right)  \cap BF\left(  I\right)
\cap\left\{  w_{t}=1\right\}
\]
The following table chooses an indexing of the vertices of $\nabla$ involved
in $B\left(  I\right)  $ and gives for each vertex the numbers $n_{\nabla}$
and $n_{B}$ of faces of $\nabla$ and faces of $B\left(  I\right)  $ it is
contained in

\begin{center}
%
\end{center}

\noindent With this indexing the Bergman subcomplex $B\left(  I\right)  $ of
$\operatorname*{Poset}\left(  \nabla\right)  $ associated to the degeneration
$\mathfrak{X}$ is

\begin{center}
%
\]
and the $F$-vectors of its faces are

\begin{center}
%
\end{center}

\paragraph{Dual complex}

The dual complex $\operatorname*{dual}\left(  B\left(  I\right)  \right)
=\left(  B\left(  I\right)  \right)  ^{\ast}$ of deformations associated to
$B\left(  I\right)  $ via initial ideals is given by

\begin{center}
%
\end{center}

\noindent when writing the vertices of the faces as deformations of $X_{0}$.
Note that the $T$-invariant basis of deformations associated to a face is
given by all lattice points of the corresponding polytope in $M_{\mathbb{R}}$.

In order to compress the output we list one representative in any set of faces
$G$ with fixed $F$-vector of $G$ and $G^{\ast}$.

When numbering the vertices of the faces of $\operatorname*{dual}\left(
B\left(  I\right)  \right)  $ by the Cox variables of the mirror toric Fano
variety $Y^{\circ}$ the complex $\operatorname*{dual}\left(  B\left(
I\right)  \right)  $ is

\begin{center}
\begin{longtable}
[c]{l}
$[],\medskip$
\\
$[],\medskip$
\\
$[[2, 3, 5, 6, 8, 9]^{*}=\left\langle y_{5},y_{11}\right\rangle,...],\medskip$
\\
$[[1, 2, 4, 5]^{*}=\left\langle y_{2},y_{4},y_{8},y_{10}\right\rangle,[4, 5, 6]^{*}=\left\langle y_{6},y_{2},y_{12},y_{10}\right\rangle,...],\medskip$
\\
$[[3, 9]^{*}=\left\langle y_{1},y_{3},y_{5},y_{7},y_{11},y_{9}\right\rangle,...],\medskip$
\\
$[[1]^{*}=\left\langle y_{1},y_{2},y_{3},y_{4},y_{7},y_{8},y_{9},y_{10}\right\rangle,...],\medskip$
\\
$[]$
\\
\end{longtable}
\end{center}

\noindent The dual complex has the $F$-vector%
\[
\begin{tabular}
[c]{rrrrrrrr}
Dimension &$-1$ & $0$ & $1$ & $2$ & $3$ & $4$ & $5$\\
Number of faces &$0$ & $0$ & $6$ & $15$ & $18$ & $9$ & $0$
\end{tabular}
\]
and the $F$-vectors of the faces of $\operatorname*{dual}\left(  B\left(
I\right)  \right)  $ are

\begin{center}
\begin{longtable}
[c]{rrll}
Dimension &Number of faces &F-vector &\\
$1$ & $6$ & $\left(1,2,1,0,0,0,0\right)$ &edge\\
$2$ & $15$ & $\left(1,4,4,1,0,0,0\right)$ &quadrangle\\
$3$ & $18$ & $\left(1,6,9,5,1,0,0\right)$ &prism\\
$4$ & $9$ & $\left(1,8,16,14,6,1,0\right)$ &
\end{longtable}
\end{center}

Recall that in this example the toric variety $Y$ is projective space. The
number of lattice points of the support of $\operatorname*{dual}\left(
B\left(  I\right)  \right)  $ relates to the dimension $h^{1,
2
}%
\left(  X\right)  $ of the complex moduli space of the generic fiber $X$ of
$\mathfrak{X}$ and to the dimension $h^{1,1}\left(  \bar{X}^{\circ}\right)  $
of the K\"{a}hler moduli space of the $MPCR$-blowup $\bar{X}^{\circ}$ of the
generic fiber $X^{\circ}$ of the mirror degeneration
\begin{align*}
\left\vert \operatorname*{supp}\left(  \operatorname*{dual}\left(  B\left(
I\right)  \right)  \right)  \cap M\right\vert  &
=
108
=
35
+
73
=\dim\left(  \operatorname*{Aut}%
\left(  Y\right)  \right)  +h^{1,
2
}\left(  X\right) \\
&  =
30
+
5
+
73
\\
&  =\left\vert \operatorname*{Roots}\left(  Y\right)  \right\vert +\dim\left(
T_{Y}\right)  +h^{1,1}\left(  \bar{X}^{\circ}\right)
\end{align*}
There are%
\[
h^{1,
2
}\left(  X\right)  +\dim\left(  T_{Y^{\circ}}\right)
=
73
+
5
\]
non-trivial toric polynomial deformations of $X_{0}$

\begin{center}

\end{center}

\noindent They correspond to the toric divisors

\begin{center}
%
\end{center}

\noindent on a MPCR-blowup of $Y^{\circ}$ inducing $
73
$ non-zero
toric divisor classes on the mirror $X^{\circ}$. The following $
30
$
toric divisors of $Y$ induce the trivial divisor class on $X^{\circ}$

\begin{center}
%
\end{center}

\paragraph{Mirror special fiber}

The ideal $I_{0}^{\circ}$ of the monomial special fiber $X_{0}^{\circ}$ of the
mirror degeneration $\mathfrak{X}^{\circ}$ is generated by the following set
of monomials in $S^{\circ}$%
\[
\left\{  
%
\right\}
\]
Indeed already the ideal%
\[
J_{0}^{\circ}=\left\langle 
%
\right\rangle
\]
defines the same subvariety of the toric variety $Y^{\circ}$, and
$J_{0}^{\circ\Sigma}=I_{0}^{\circ}$. Recall that passing from $J_{0}^{\circ}$
to $J_{0}^{\circ\Sigma}$ is the non-simplicial toric analogue of saturation.

The complex $B\left(  I\right)  ^{\ast}$ labeled by the variables of the Cox
ring $S^{\circ}$ of $Y^{\ast}$, as written in the last section, is the complex
$SP\left(  I_{0}^{\circ}\right)  $ of prime ideals of the toric strata of the
special fiber $X_{0}^{\circ}$ of the mirror degeneration $\mathfrak{X}^{\circ
}$, i.e., the primary decomposition of $I_{0}^{\circ}$ is%
\[
I_{0}^{\circ}=

\]
Each facet $F\in B\left(  I\right)  $ corresponds to one of these ideals and this
ideal is generated by the products of facets of $\nabla$ containing $F$.

\paragraph{Covering structure in the deformation complex of the
degeneration $\mathfrak{X}$}

According to the local reduced Gr\"{o}bner basis each face of the complex of
deformations $\operatorname*{dual}\left(  B\left(  I\right)  \right)  $
decomposes into 
$2$
polytopes forming a $
2
:1$
trivial
covering of $B\left(  I\right)  $

\begin{center}
%
\end{center}

\noindent Here the faces $F\in B\left(  I\right)  $ are specified both via the
vertices of $F^{\ast}$ labeled by the variables of $S^{\circ}$ and by the
numbering of the vertices of $B\left(  I\right)  $ chosen above.

\noindent This covering has 
2 sheets forming the complexes
 
\[
%
\end{center}

\noindent Writing the vertices of the faces as deformations the covering is
given by

\begin{center}
%
\end{center}

\noindent with 
2 sheets forming the complexes

\begin{center}
\[
%
\]
\end{center}

\noindent Note that the torus invariant basis of deformations corresponding to
a Bergman face is given by the set of all lattice points of the polytope
specified above.

\paragraph{Limit map}

The limit map $\lim:B\left(  I\right)  \rightarrow\operatorname*{Poset}\left(
\Delta\right)  $ associates to a face $F$ of $B\left(  I\right)  $ the face of
$\Delta$ formed by the limit points of arcs lying over the weight vectors
$w\in F$, i.e. with lowest order term $t^{w}$.

Labeling the faces of the Bergman complex $B\left(  I\right)  \subset
\operatorname*{Poset}\left(  \nabla\right)  $ and the faces of
$\operatorname*{Poset}\left(  \Delta\right)  $ by the corresponding dual faces
of $\nabla^{\ast}$ and $\Delta^{\ast}$, hence considering the limit map
$\lim:B\left(  I\right)  \rightarrow\operatorname*{Poset}\left(
\Delta\right)  $ as a map $B\left(  I\right)  ^{\ast}\rightarrow
\operatorname*{Poset}\left(  \Delta^{\ast}\right)  $, the limit correspondence
is given by

\begin{center}

\end{center}

The image of the limit map coincides with the image of $\mu$ and with the
Bergman complex of the mirror, i.e. $\lim\left(  B\left(  I\right)  \right)
=\mu\left(  B\left(  I\right)  \right)  =B\left(  I^{\ast}\right)  $.

\paragraph{Mirror complex}

Numbering the vertices of the mirror complex $\mu\left(  B\left(  I\right)
\right)  $ as%
\[
%
\end{center}

\noindent The ordering of the faces of $\mu\left(  B\left(  I\right)  \right)
$ is compatible with above ordering of the faces of $B\left(  I\right)  $. The
$F$-vectors of the faces of $\mu\left(  B\left(  I\right)  \right)  $ are%
\[
%
\]
The first order deformations of the mirror special fiber $X_{0}^{\circ}$
correspond to the lattice points of the dual complex $\left(  \mu\left(
B\left(  I\right)  \right)  \right)  ^{\ast}\subset\operatorname*{Poset}%
\left(  \Delta^{\ast}\right)  $ of the mirror complex. We label the first
order deformations of $X_{0}^{\circ}$ corresponding to vertices of
$\Delta^{\ast}$ by the homogeneous coordinates of $Y$%
\[
%
\]
So writing the vertices of the faces of $\left(  \mu\left(  B\left(  I\right)
\right)  \right)  ^{\ast}$ as homogeneous coordinates of $Y$, the complex
$\left(  \mu\left(  B\left(  I\right)  \right)  \right)  ^{\ast}$ is given by

\begin{center}
%
\end{center}

\noindent The complex $\mu\left(  B\left(  I\right)  \right)  ^{\ast}$ labeled
by the variables of the Cox ring $S$ gives the ideals of the toric strata of
the special fiber $X_{0}$ of $\mathfrak{X}$, i.e. the complex $SP\left(
I_{0}\right)  $, so in particular the primary decomposition of $I_{0}$ is%
\[
I_{0}=
%
\]

\paragraph{Covering structure in the deformation complex of the mirror
degeneration}

Each face of the complex of deformations $\left(  \mu\left(  B\left(
I\right)  \right)  \right)  ^{\ast}$ of the mirror special fiber $X_{0}%
^{\circ}$ decomposes as the convex hull of 
$2$
polytopes forming
a $
2
:1$ 
trivial
covering of $\mu\left(  B\left(
I\right)  \right)  ^{\vee}$

\begin{center}
%
\end{center}

\noindent 
Due to the singularities of
$Y^{\circ}$ this covering involves degenerate faces, i.e. faces $G\mapsto
F^{\vee}$ with $\dim\left(  G\right)  <\dim\left(  F^{\vee}\right)  $.

\paragraph{Mirror degeneration}

The space of first order deformations of $X_{0}^{\circ}$ in the mirror
degeneration $\mathfrak{X}^{\circ}$ has dimension $
6
$ and the
deformations represented by the monomials%
\[
\left\{  
%
\right\}
\]
form a torus invariant
basis.
The number of lattice points of the dual of the mirror complex of $I$ relates
to the dimension $h^{1,
2
}\left(  X^{\circ}\right)  $ of complex
moduli space of the generic fiber $X^{\circ}$ of $\mathfrak{X}^{\circ}$ and to
the dimension $h^{1,1}\left(  X\right)  $ of the K\"{a}hler moduli space of
the generic fiber $X$ of $\mathfrak{X}$ via%
\begin{align*}
\left\vert \operatorname*{supp}\left(  \left(  \mu\left(  B\left(  I\right)
\right)  \right)  ^{\ast}\right)  \cap N\right\vert  &
=
6
=
5
+
1
\\
&  =\dim\left(  \operatorname*{Aut}\left(  Y^{\circ}\right)  \right)
+h^{1,
2
}\left(  X^{\circ}\right)  =\dim\left(  T\right)
+h^{1,1}\left(  X\right)
\end{align*}

The 
mirror degeneration 
$\mathfrak{X}^{\circ}\subset Y^{\circ}\times\operatorname*{Spec}%
\mathbb{C}\left[  t\right]  $
 of $\mathfrak{X}$ is
given by the ideal 
$I^{\circ}\subset S^{\circ}\otimes\mathbb{C}\left[  t\right]  $
generated by%
\[

\]
Indeed already the ideal $J^{\circ}$ generated by%
\[
\left\{  
%
\right\}
\]
defines $\mathfrak{X}^{\circ}$.

\paragraph{Contraction of the mirror degeneration}

In the following we give a birational map relating the degeneration
$\mathfrak{X}^{\circ}$ to a Greene-Plesser type orbifolding mirror family by
contracting divisors on $Y^{\circ}$.
See also Section \ref{sec orbifolding mirror families} below.

In order to contract the divisors

\begin{center}
%
\end{center}

\noindent consider the $\mathbb{Q}$-factorial toric Fano variety $\hat
{Y}^{\circ}=X\left(  \hat{\Sigma}^{\circ}\right)  $, where $\hat{\Sigma
}^{\circ}=\operatorname*{Fan}\left(  \hat{P}^{\circ}\right)  \subset
M_{\mathbb{R}}$ is the fan over the 
Fano polytope $\hat
{P}^{\circ}\subset M_{\mathbb{R}}$ given as the convex hull of the remaining
vertices of $P^{\circ}=\nabla^{\ast}$ corresponding to the Cox variables%
\[
%
\]
of the toric variety $\hat{Y}^{\circ}$ with Cox ring
\[
\hat{S}^{\circ}=\mathbb{C}
[y_{11},y_{2},y_{8},y_{6},y_{1},y_{7}]
\]
The Cox variables of $\hat{Y}^{\circ}$ correspond to the set of Fermat deformations of $\mathfrak{X}$.

Let%
\[
Y^{\circ}=X\left(  \Sigma^{\circ}\right)  \rightarrow X\left(  \hat{\Sigma
}^{\circ}\right)  =\hat{Y}^{\circ}%
\]
be a birational map from $Y^{\circ}$ to a minimal birational model $\hat
{Y}^{\circ}$, which contracts the divisors of the rays $\Sigma^{\circ}\left(
1\right)  -\hat{\Sigma}^{\circ}\left(  1\right)  $ corresponding to the Cox
variables
\[
\begin{tabular}
[c]{llllll}
$y_{3}$ &$y_{4}$ &$y_{5}$ &$y_{9}$ &$y_{10}$ &$y_{12}$
\end{tabular}
\]

Representing $\hat{Y}^{\circ}$ as a quotient we have
\[
\hat{Y}^{\circ}=\left(  \mathbb{C}^{
6
}-V\left(  B\left(
\hat{\Sigma}^{\circ}\right)  \right)  \right)  //\hat{G}^{\circ}%
\]
with%
\[
\hat{G}^{\circ}=
\mathbb{Z}_{3}\times\mathbb{Z}_{3}\times\mathbb{Z}_{9}\times\left(  \mathbb{C}^{\ast}\right)  ^{1}
\]
acting via%
\[
\xi y=
\left( \,u_{2}^{2}\,v_{1} \cdot y_{11},\,u_{1}^{2}\,u_{2}\,u_{3}^{4}\,v_{1} \cdot y_{2},\,u_{2}\,u_{3}^{3}\,v_{1} \cdot y_{8},\,u_{2}^{2}\,u_{3}\,v_{1} \cdot y_{6},\,u_{1}\,u_{3}^{4}\,v_{1} \cdot y_{1},\,v_{1} \cdot y_{7} \right)
\]
for $\xi=
\left(u_1,u_2,u_3,v_1\right)
\in\hat{G}^{\circ}$ and $y\in\mathbb{C}%
^{
6
}-V\left(  B\left(  \hat{\Sigma}^{\circ}\right)  \right)  $.

Hence with the group%
\[
\hat{H}^{\circ}=
\mathbb{Z}_{3}\times\mathbb{Z}_{3}\times\mathbb{Z}_{9}
\]
of order 
81
 the toric variety $\hat{Y}^{\circ}$ is the quotient%
\[
\hat{Y}^{\circ}=\mathbb{P}^{
5
}/\hat{H}^{\circ}%
\]
of projective space $\mathbb{P}^{
5
}$.

The first order mirror degeneration $\mathfrak{X}^{1\circ}$ induces via
$Y\rightarrow\hat{Y}^{\circ}$ a degeneration $\mathfrak{\hat{X}}^{1\circ
}\subset\hat{Y}^{\circ}\times\operatorname*{Spec}\mathbb{C}\left[  t\right]
/\left\langle t^{2}\right\rangle $ given by the ideal $\hat{I}^{1\circ}%
\subset\hat{S}^{\circ}\otimes\mathbb{C}\left[  t\right]  /\left\langle
t^{2}\right\rangle $ generated by the
 Fermat-type equations
\[
\left\{  

\right\}
\]

Note, that%
\[
\left\vert \operatorname*{supp}\left(  \left(  \mu\left(  B\left(  I\right)
\right)  \right)  ^{\ast}\right)  \cap N-\operatorname*{Roots}\left(  \hat
{Y}^{\circ}\right)  \right\vert -\dim\left(  T_{\hat{Y}^{\circ}}\right)
=
6
-
5
=
1
\]
so this family has one independent parameter.

The special fiber $\hat{X}_{0}^{\circ}\subset\hat{Y}^{\circ}$ of
$\mathfrak{\hat{X}}^{1\circ}$ is cut out by the monomial ideal $\hat{I}%
_{0}^{\circ}\subset\hat{S}^{\circ}$ generated by%
\[
\left\{  
%
\right\}
\]

The complex $SP\left(  \hat{I}_{0}^{\circ}\right)  $ of prime ideals of the
toric strata of $\hat{X}_{0}^{\circ}$ is

\begin{center}
%
\end{center}

\noindent so $\hat{I}_{0}^{\circ}$ has the primary decomposition%
\[
\hat{I}_{0}^{\circ}=
%
\]
The Bergman subcomplex $B\left(  I\right)  \subset\operatorname*{Poset}\left(
\nabla\right)  $ induces a subcomplex of $\operatorname*{Poset}\left(
\hat{\nabla}\right)  $, $\hat{\nabla}=\left(  \hat{P}^{\circ}\right)  ^{\ast}$
corresponding to $\hat{I}_{0}^{\circ}$. Indexing of the 
vertices
of $\hat{\nabla}$ by%
\[
%
\]
this complex is given by
\[
\left\{  
%
\right\}
\]

This is the one parameter Greene-Plesser orbifolding mirror family of the
generic complete intersection of two cubics in $\mathbb{P}^{5}$, given in
\cite{LT Lines on CalabiYau complete intersections mirror symmetry and
PicardFuchs equations}.

\subsubsection{The Calabi-Yau threefold given as the Pfaffian complete
intersection of two generic quadrics and a generic cubic in $\mathbb{P}^{6}$}

\paragraph{Setup}

Let $Y=\mathbb{P}^{
6
}=X\left(  \Sigma\right)  $, $\Sigma
=\operatorname*{Fan}\left(  P\right)  =NF\left(  \Delta\right)  \subset
N_{\mathbb{R}}$ with the Fano polytope $P=\Delta^{\ast}$ given by%
\[
\Delta=\operatorname*{convexhull}\left(  
\begin{tabular}
[c]{ll}
$\left(6,-1,-1,-1,-1,-1\right)$ &$\left(-1,6,-1,-1,-1,-1\right)$\\
$\left(-1,-1,6,-1,-1,-1\right)$ &$\left(-1,-1,-1,6,-1,-1\right)$\\
$\left(-1,-1,-1,-1,6,-1\right)$ &$\left(-1,-1,-1,-1,-1,6\right)$\\
$\left(-1,-1,-1,-1,-1,-1\right)$ &
\end{tabular}
\right)  \subset
M_{\mathbb{R}}%
\]
and let%
\[
S=\mathbb{C}
[x_0, x_1, x_2, x_3, x_4, x_5, x_6]
\]
be the Cox ring of $Y$ with the variables%
\[
\begin{tabular}
[c]{ll}
$x_{1} = x_{\left(1,0,0,0,0,0\right)}$ &$x_{2} = x_{\left(0,1,0,0,0,0\right)}$\\
$x_{3} = x_{\left(0,0,1,0,0,0\right)}$ &$x_{4} = x_{\left(0,0,0,1,0,0\right)}$\\
$x_{5} = x_{\left(0,0,0,0,1,0\right)}$ &$x_{6} = x_{\left(0,0,0,0,0,1\right)}$\\
$x_{0} = x_{\left(-1,-1,-1,-1,-1,-1\right)}$ &
\end{tabular}
\]
associated to the rays of $\Sigma$.

Consider the degeneration $\mathfrak{X}\subset Y\times\operatorname*{Spec}%
\mathbb{C}\left[  t\right]  $ of 
Pfaffian complete intersection
Calabi-Yau 3-folds
 with Buchsbaum-Eisenbud
resolution%
\[
0\rightarrow\mathcal{O}_{Y}\left(  -
7
\right)  \rightarrow
\mathcal{E}^{\ast}\left(  -3\right)  \overset{A_{t}}{\rightarrow}%
\mathcal{E}\left(  -2\right)  \rightarrow\mathcal{O}_{Y}\rightarrow
\mathcal{O}_{X_{t}}\rightarrow0
\]
where%
\[
\mathcal{E}=\text{
$2\mathcal{O}\left(  1\right)  \oplus\mathcal{O}$
}%
\]%
\[
A_{t}=A_{0}+t\cdot A
\]%
\[
A_{0}=
\left [
\right ]
\]
the monomial special fiber of $\mathfrak{X}$ is given by%
\[
I_{0}=\left\langle 
%
\right\rangle
\]
and generic $A\in%
{\textstyle\bigwedge\nolimits^{2}}
\mathcal{E}\left(  1\right)  $.

The degeneration $\mathfrak{X}$ is given by the ideal $I\subset S\otimes
\mathbb{C}\left[  t\right]  $ with $I_{0}$-reduced generators of degrees
$
2,2,3
$%
\[
\left\{  
%
\right\}
\]

\paragraph{Special fiber Gr\"{o}bner cone}

The space of first order deformations of $\mathfrak{X}$ has dimension
$
121
$ and the deformations represented by the Cox Laurent monomials

\begin{center}
%
\end{center}

\noindent form a torus invariant basis. These deformations give linear
inequalities defining the special fiber Gr\"{o}bner cone, i.e.,%
\[
C_{I_{0}}\left(  I\right)  =\left\{  \left(  w_{t},w_{y}\right)  \in
\mathbb{R}\oplus N_{\mathbb{R}}\mid\left\langle w,A^{-1}\left(  m\right)
\right\rangle \geq-w_{t}\forall m\right\}
\]
for above monomials $m$ and the presentation matrix
\[
A=
\left [%
\right ]
\]
of $A_{
5
}\left(  Y\right)  $.

The vertices of special fiber polytope
\[
\nabla=C_{I_{0}}\left(  I\right)  \cap\left\{  w_{t}=1\right\}
\]
and the number of faces, each vertex is contained in, are given by the table

\begin{center}
%
\end{center}
\noindent The number of faces of $\nabla$ and their $F$-vectors are
\begin{center}
%
\end{center}
The dual $P^{\circ}=\nabla^{\ast}$ of $\nabla$ is a Fano polytope with vertices
\begin{center}
%
\end{center}
\noindent and the $F$-vectors of the faces of $P^{\circ}$ are
\begin{center}
%
\end{center}
The polytopes $\nabla$ and $P^{\circ}$ have $0$ as the unique interior lattice
point, and the Fano polytope $P^{\circ}$ defines the embedding toric Fano
variety $Y^{\circ}=X\left(  \operatorname*{Fan}\left(  P^{\circ}\right)
\right)  $ of the fibers of the mirror degeneration. The dimension of the
automorphism group of $Y^{\circ}$ is
\[
\dim\left(  \operatorname*{Aut}\left(  Y^{\circ}\right)  \right)
=
6
\]
Let%
\[
S^{\circ}=\mathbb{C}
[y_{1},y_{2},...,y_{21}]
\]
be the Cox ring of $Y^{\circ}$ with variables

\begin{center}
%
\end{center}

\paragraph{Bergman subcomplex}

Intersecting the tropical variety of $I$ with the special fiber Gr\"{o}bner
cone $C_{I_{0}}\left(  I\right)  $ we obtain the special fiber Bergman
subcomplex
\[
B\left(  I\right)  =C_{I_{0}}\left(  I\right)  \cap BF\left(  I\right)
\cap\left\{  w_{t}=1\right\}
\]
The following table chooses an indexing of the vertices of $\nabla$ involved
in $B\left(  I\right)  $ and gives for each vertex the numbers $n_{\nabla}$
and $n_{B}$ of faces of $\nabla$ and faces of $B\left(  I\right)  $ it is
contained in

\begin{center}
%
\end{center}

\noindent With this indexing the Bergman subcomplex $B\left(  I\right)  $ of
$\operatorname*{Poset}\left(  \nabla\right)  $ associated to the degeneration
$\mathfrak{X}$ is

\begin{center}
%
\]
and the $F$-vectors of its faces are

\begin{center}
%
\end{center}

\paragraph{Dual complex}

The dual complex $\operatorname*{dual}\left(  B\left(  I\right)  \right)
=\left(  B\left(  I\right)  \right)  ^{\ast}$ of deformations associated to
$B\left(  I\right)  $ via initial ideals is given by

\begin{center}
%
\end{center}

\noindent when writing the vertices of the faces as deformations of $X_{0}$.
Note that the $T$-invariant basis of deformations associated to a face is
given by all lattice points of the corresponding polytope in $M_{\mathbb{R}}$.

In order to compress the output we list one representative in any set of faces
$G$ with fixed $F$-vector of $G$ and $G^{\ast}$.

When numbering the vertices of the faces of $\operatorname*{dual}\left(
B\left(  I\right)  \right)  $ by the Cox variables of the mirror toric Fano
variety $Y^{\circ}$ the complex $\operatorname*{dual}\left(  B\left(
I\right)  \right)  $ is
\begin{center}
\begin{longtable}
[c]{l}
$[],\medskip$
\\
$[],\medskip$
\\
$[],\medskip$
\\
$[[1, 2, 3, 4, 9, 10, 11, 12]^{*}=\left\langle y_{1},y_{8},y_{17}\right\rangle,[1, 3, 5, 7, 9, 11]^{*}=\left\langle y_{4},y_{10},y_{15}\right\rangle,$
\\
$...],\medskip$
\\
$[[1, 5, 9]^{*}=\left\langle y_{3},y_{4},y_{10},y_{11},y_{15},y_{16}\right\rangle,[1, 3, 9, 11]^{*}=\left\langle y_{1},y_{4},y_{8},y_{10},y_{15},y_{17}\right\rangle,$
\\
$...],\medskip$
\\
$[[9, 11]^{*}=\left\langle y_{7},y_{1},y_{4},y_{14},y_{8},y_{10},y_{15},y_{21},y_{17}\right\rangle,$
\\
$...],\medskip$
\\
$[[1]^{*}=\left\langle y_{1},y_{2},y_{3},y_{4},y_{8},y_{9},y_{10},y_{11},y_{15},y_{16},y_{17},y_{18}\right\rangle,$
\\
$...],\medskip$
\\
$[]$
\\
\end{longtable}
\end{center}
\noindent The dual complex has the $F$-vector%
\[
\begin{tabular}
[c]{rrrrrrrrr}
Dimension &$-1$ & $0$ & $1$ & $2$ & $3$ & $4$ & $5$ & $6$\\
Number of faces &$0$ & $0$ & $0$ & $7$ & $19$ & $24$ & $12$ & $0$
\end{tabular}
\]
and the $F$-vectors of the faces of $\operatorname*{dual}\left(  B\left(
I\right)  \right)  $ are
\begin{center}
\begin{longtable}
[c]{rrll}
Dimension &Number of faces &F-vector &\\
$2$ & $7$ & $\left(1,3,3,1,0,0,0,0\right)$ &triangle\\
$3$ & $19$ & $\left(1,6,9,5,1,0,0,0\right)$ &prism\\
$4$ & $24$ & $\left(1,9,18,15,6,1,0,0\right)$ &\\
$5$ & $12$ & $\left(1,12,30,34,21,7,1,0\right)$ &
\end{longtable}
\end{center}
Recall that in this example the toric variety $Y$ is projective space. The
number of lattice points of the support of $\operatorname*{dual}\left(
B\left(  I\right)  \right)  $ relates to the dimension $h^{1,
2
}%
\left(  X\right)  $ of the complex moduli space of the generic fiber $X$ of
$\mathfrak{X}$ and to the dimension $h^{1,1}\left(  \bar{X}^{\circ}\right)  $
of the K\"{a}hler moduli space of the $MPCR$-blowup $\bar{X}^{\circ}$ of the
generic fiber $X^{\circ}$ of the mirror degeneration
\begin{align*}
\left\vert \operatorname*{supp}\left(  \operatorname*{dual}\left(  B\left(
I\right)  \right)  \right)  \cap M\right\vert  &
=
121
=
48
+
73
=\dim\left(  \operatorname*{Aut}%
\left(  Y\right)  \right)  +h^{1,
2
}\left(  X\right) \\
&  =
42
+
6
+
73
\\
&  =\left\vert \operatorname*{Roots}\left(  Y\right)  \right\vert +\dim\left(
T_{Y}\right)  +h^{1,1}\left(  \bar{X}^{\circ}\right)
\end{align*}
There are%
\[
h^{1,
2
}\left(  X\right)  +\dim\left(  T_{Y^{\circ}}\right)
=
73
+
6
\]
non-trivial toric polynomial deformations of $X_{0}$

\begin{center}

\end{center}

\noindent They correspond to the toric divisors

\begin{center}
%
\end{center}

\noindent on a MPCR-blowup of $Y^{\circ}$ inducing $
73
$ non-zero
toric divisor classes on the mirror $X^{\circ}$. The following $
42
$
toric divisors of $Y$ induce the trivial divisor class on $X^{\circ}$

\begin{center}
%
\end{center}

\paragraph{Mirror special fiber}

The ideal $I_{0}^{\circ}$ of the monomial special fiber $X_{0}^{\circ}$ of the
mirror degeneration $\mathfrak{X}^{\circ}$ is generated by the following set
of monomials in $S^{\circ}$%
\[
\left\{  
%
\right\}
\]
Indeed already the ideal%
\[
J_{0}^{\circ}=\left\langle 
%
\right\rangle
\]
defines the same subvariety of the toric variety $Y^{\circ}$, and
$J_{0}^{\circ\Sigma}=I_{0}^{\circ}$. Recall that passing from $J_{0}^{\circ}$
to $J_{0}^{\circ\Sigma}$ is the non-simplicial toric analogue of saturation.

The complex $B\left(  I\right)  ^{\ast}$ labeled by the variables of the Cox
ring $S^{\circ}$ of $Y^{\ast}$, as written in the last section, is the complex
$SP\left(  I_{0}^{\circ}\right)  $ of prime ideals of the toric strata of the
special fiber $X_{0}^{\circ}$ of the mirror degeneration $\mathfrak{X}^{\circ
}$, i.e., the primary decomposition of $I_{0}^{\circ}$ is%
\[
I_{0}^{\circ}=

\]
Each facet $F\in B\left(  I\right)  $ corresponds to one of these ideals and this
ideal is generated by the products of facets of $\nabla$ containing $F$.

\paragraph{Covering structure in the deformation complex of the
degeneration $\mathfrak{X}$}

According to the local reduced Gr\"{o}bner basis each face of the complex of
deformations $\operatorname*{dual}\left(  B\left(  I\right)  \right)  $
decomposes into 
$3$
polytopes forming a $
3
:1$
trivial
covering of $B\left(  I\right)  $

\begin{center}
%
\end{center}

\noindent Here the faces $F\in B\left(  I\right)  $ are specified both via the
vertices of $F^{\ast}$ labeled by the variables of $S^{\circ}$ and by the
numbering of the vertices of $B\left(  I\right)  $ chosen above.

\noindent This covering has 
3 sheets forming the complexes
 
\[
%
\end{center}

\noindent Writing the vertices of the faces as deformations the covering is
given by

\begin{center}
%
\end{center}

\noindent with 
3 sheets forming the complexes

\begin{center}
\[
%
\]
\end{center}

\noindent Note that the torus invariant basis of deformations corresponding to
a Bergman face is given by the set of all lattice points of the polytope
specified above.

\paragraph{Limit map}

The limit map $\lim:B\left(  I\right)  \rightarrow\operatorname*{Poset}\left(
\Delta\right)  $ associates to a face $F$ of $B\left(  I\right)  $ the face of
$\Delta$ formed by the limit points of arcs lying over the weight vectors
$w\in F$, i.e. with lowest order term $t^{w}$.

Labeling the faces of the Bergman complex $B\left(  I\right)  \subset
\operatorname*{Poset}\left(  \nabla\right)  $ and the faces of
$\operatorname*{Poset}\left(  \Delta\right)  $ by the corresponding dual faces
of $\nabla^{\ast}$ and $\Delta^{\ast}$, hence considering the limit map
$\lim:B\left(  I\right)  \rightarrow\operatorname*{Poset}\left(
\Delta\right)  $ as a map $B\left(  I\right)  ^{\ast}\rightarrow
\operatorname*{Poset}\left(  \Delta^{\ast}\right)  $, the limit correspondence
is given by

\begin{center}

\end{center}

The image of the limit map coincides with the image of $\mu$ and with the
Bergman complex of the mirror, i.e. $\lim\left(  B\left(  I\right)  \right)
=\mu\left(  B\left(  I\right)  \right)  =B\left(  I^{\ast}\right)  $.

\paragraph{Mirror complex}

Numbering the vertices of the mirror complex $\mu\left(  B\left(  I\right)
\right)  $ as%
\[
%
\end{center}

\noindent The ordering of the faces of $\mu\left(  B\left(  I\right)  \right)
$ is compatible with above ordering of the faces of $B\left(  I\right)  $. The
$F$-vectors of the faces of $\mu\left(  B\left(  I\right)  \right)  $ are%
\[
%
\]
The first order deformations of the mirror special fiber $X_{0}^{\circ}$
correspond to the lattice points of the dual complex $\left(  \mu\left(
B\left(  I\right)  \right)  \right)  ^{\ast}\subset\operatorname*{Poset}%
\left(  \Delta^{\ast}\right)  $ of the mirror complex. We label the first
order deformations of $X_{0}^{\circ}$ corresponding to vertices of
$\Delta^{\ast}$ by the homogeneous coordinates of $Y$%
\[
%
\]
So writing the vertices of the faces of $\left(  \mu\left(  B\left(  I\right)
\right)  \right)  ^{\ast}$ as homogeneous coordinates of $Y$, the complex
$\left(  \mu\left(  B\left(  I\right)  \right)  \right)  ^{\ast}$ is given by

\begin{center}
%
\end{center}

\noindent The complex $\mu\left(  B\left(  I\right)  \right)  ^{\ast}$ labeled
by the variables of the Cox ring $S$ gives the ideals of the toric strata of
the special fiber $X_{0}$ of $\mathfrak{X}$, i.e. the complex $SP\left(
I_{0}\right)  $, so in particular the primary decomposition of $I_{0}$ is%
\[
I_{0}=
%
\]

\paragraph{Covering structure in the deformation complex of the mirror
degeneration}

Each face of the complex of deformations $\left(  \mu\left(  B\left(
I\right)  \right)  \right)  ^{\ast}$ of the mirror special fiber $X_{0}%
^{\circ}$ decomposes as the convex hull of 
$3$
polytopes forming
a $
3
:1$ 
trivial
covering of $\mu\left(  B\left(
I\right)  \right)  ^{\vee}$

\begin{center}
%
\end{center}

\noindent 
Due to the singularities of
$Y^{\circ}$ this covering involves degenerate faces, i.e. faces $G\mapsto
F^{\vee}$ with $\dim\left(  G\right)  <\dim\left(  F^{\vee}\right)  $.

\paragraph{Mirror degeneration}

The space of first order deformations of $X_{0}^{\circ}$ in the mirror
degeneration $\mathfrak{X}^{\circ}$ has dimension $
7
$ and the
deformations represented by the monomials%
\[
\left\{  
%
\right\}
\]
form a torus invariant
basis.
The number of lattice points of the dual of the mirror complex of $I$ relates
to the dimension $h^{1,
2
}\left(  X^{\circ}\right)  $ of complex
moduli space of the generic fiber $X^{\circ}$ of $\mathfrak{X}^{\circ}$ and to
the dimension $h^{1,1}\left(  X\right)  $ of the K\"{a}hler moduli space of
the generic fiber $X$ of $\mathfrak{X}$ via%
\begin{align*}
\left\vert \operatorname*{supp}\left(  \left(  \mu\left(  B\left(  I\right)
\right)  \right)  ^{\ast}\right)  \cap N\right\vert  &
=
7
=
6
+
1
\\
&  =\dim\left(  \operatorname*{Aut}\left(  Y^{\circ}\right)  \right)
+h^{1,
2
}\left(  X^{\circ}\right)  =\dim\left(  T\right)
+h^{1,1}\left(  X\right)
\end{align*}

The 
conjectural first order
mirror degeneration 
$\mathfrak{X}^{\circ}\subset Y^{\circ}\times\operatorname*{Spec}%
\mathbb{C}\left[  t\right]  $
 of $\mathfrak{X}$ is
given by the ideal 
$I^{\circ}\subset S^{\circ}\otimes\mathbb{C}\left[  t\right]  $
generated by%
\[

\]
Indeed already the ideal $J^{\circ}$ generated by%
\[
\left\{  
%
\right\}
\]
defines $\mathfrak{X}^{\circ}$.

\paragraph{Contraction of the mirror degeneration}

In the following we give a birational map relating the degeneration
$\mathfrak{X}^{\circ}$ to a Greene-Plesser type orbifolding mirror family by
contracting divisors on $Y^{\circ}$.
See also Section \ref{sec orbifolding mirror families} below.

In order to contract the divisors

\begin{center}
%
\end{center}

\noindent consider the $\mathbb{Q}$-factorial toric Fano variety $\hat
{Y}^{\circ}=X\left(  \hat{\Sigma}^{\circ}\right)  $, where $\hat{\Sigma
}^{\circ}=\operatorname*{Fan}\left(  \hat{P}^{\circ}\right)  \subset
M_{\mathbb{R}}$ is the fan over the 
Fano polytope $\hat
{P}^{\circ}\subset M_{\mathbb{R}}$ given as the convex hull of the remaining
vertices of $P^{\circ}=\nabla^{\ast}$ corresponding to the Cox variables%
\[
%
\]
of the toric variety $\hat{Y}^{\circ}$ with Cox ring
\[
\hat{S}^{\circ}=\mathbb{C}
[y_{15},y_{8},y_{9},y_{3},y_{7},y_{19},y_{6}]
\]
The Cox variables of $\hat{Y}^{\circ}$ correspond a set of Fermat deformations of $\mathfrak{X}$. There are 6 sets of Fermat deformations with $\mathbb{C}^{\ast}$ acting with weight $1$ on all Cox variables 
\[
%
\]

Let%
\[
Y^{\circ}=X\left(  \Sigma^{\circ}\right)  \rightarrow X\left(  \hat{\Sigma
}^{\circ}\right)  =\hat{Y}^{\circ}%
\]
be a birational map from $Y^{\circ}$ to a minimal birational model $\hat
{Y}^{\circ}$, which contracts the divisors of the rays $\Sigma^{\circ}\left(
1\right)  -\hat{\Sigma}^{\circ}\left(  1\right)  $ corresponding to the Cox
variables
\[
%
\]

Representing $\hat{Y}^{\circ}$ as a quotient we have
\[
\hat{Y}^{\circ}=\left(  \mathbb{C}^{
7
}-V\left(  B\left(
\hat{\Sigma}^{\circ}\right)  \right)  \right)  //\hat{G}^{\circ}%
\]
with%
\[
\hat{G}^{\circ}=
\mathbb{Z}_{2}\times\mathbb{Z}_{2}\times\mathbb{Z}_{12}\times\left(  \mathbb{C}^{\ast}\right)  ^{1}
\]
acting via%
\[
\xi y=
\left( \,u_{1}\,u_{3}^{11}\,v_{1} \cdot y_{15},\,u_{2}\,u_{3}^{9}\,v_{1} \cdot y_{8},\,u_{2}\,u_{3}^{3}\,v_{1} \cdot y_{9},\,u_{2}\,u_{3}^{6}\,v_{1} \cdot y_{3},\,u_{1}\,u_{2}\,u_{3}^{3}\,v_{1} \cdot y_{7},\,u_{1}\,u_{3}^{7}\,v_{1} \cdot y_{19},\,v_{1} \cdot y_{6} \right)
\]
for $\xi=
\left(u_1,u_2,u_3,v_1\right)
\in\hat{G}^{\circ}$ and $y\in\mathbb{C}%
^{
7
}-V\left(  B\left(  \hat{\Sigma}^{\circ}\right)  \right)  $.

Hence with the group%
\[
\hat{H}^{\circ}=
\mathbb{Z}_{2}\times\mathbb{Z}_{2}\times\mathbb{Z}_{12}
\]
of order 
48
 the toric variety $\hat{Y}^{\circ}$ is the quotient%
\[
\hat{Y}^{\circ}=\mathbb{P}^{
6
}/\hat{H}^{\circ}%
\]
of projective space $\mathbb{P}^{
6
}$.

The first order mirror degeneration $\mathfrak{X}^{1\circ}$ induces via
$Y\rightarrow\hat{Y}^{\circ}$ a degeneration $\mathfrak{\hat{X}}^{1\circ
}\subset\hat{Y}^{\circ}\times\operatorname*{Spec}\mathbb{C}\left[  t\right]
/\left\langle t^{2}\right\rangle $ given by the ideal $\hat{I}^{1\circ}%
\subset\hat{S}^{\circ}\otimes\mathbb{C}\left[  t\right]  /\left\langle
t^{2}\right\rangle $ generated by the
 Fermat-type equations
\[
\left\{  

\right\}
\]

Note, that%
\[
\left\vert \operatorname*{supp}\left(  \left(  \mu\left(  B\left(  I\right)
\right)  \right)  ^{\ast}\right)  \cap N-\operatorname*{Roots}\left(  \hat
{Y}^{\circ}\right)  \right\vert -\dim\left(  T_{\hat{Y}^{\circ}}\right)
=
7
-
6
=
1
\]
so this family has one independent parameter.

The special fiber $\hat{X}_{0}^{\circ}\subset\hat{Y}^{\circ}$ of
$\mathfrak{\hat{X}}^{1\circ}$ is cut out by the monomial ideal $\hat{I}%
_{0}^{\circ}\subset\hat{S}^{\circ}$ generated by%
\[
\left\{  
%
\right\}
\]

The complex $SP\left(  \hat{I}_{0}^{\circ}\right)  $ of prime ideals of the
toric strata of $\hat{X}_{0}^{\circ}$ is

\begin{center}
%
\end{center}

\noindent so $\hat{I}_{0}^{\circ}$ has the primary decomposition%
\[
\hat{I}_{0}^{\circ}=
%
\]
The Bergman subcomplex $B\left(  I\right)  \subset\operatorname*{Poset}\left(
\nabla\right)  $ induces a subcomplex of $\operatorname*{Poset}\left(
\hat{\nabla}\right)  $, $\hat{\nabla}=\left(  \hat{P}^{\circ}\right)  ^{\ast}$
corresponding to $\hat{I}_{0}^{\circ}$. Indexing of the 
vertices
of $\hat{\nabla}$ by%
\[
%
\]
this complex is given by
\[
\left\{  
%
\right\}
\]

\section{The tropical
\index{mirror construction}%
mirror construction\label{Sec tropical mirror construction}}

\subsection{Concept of the tropical mirror
construction\label{Sec concept of the tropical mirror construction}}

In the following, the concepts involved in the tropical mirror construction
are summarized, omitting detailed conditions and technicalities.

Let $N\cong\mathbb{Z}^{n}$ be a lattice, $M=\operatorname*{Hom}\left(
N,\mathbb{Z}\right)  $ the dual lattice and $Y=X\left(  \Sigma\right)  $ a
toric Fano variety of dimension $n$ given by a Fano polytope $P\subset
N_{\mathbb{R}}$, i.e., $\Sigma=\Sigma\left(  P\right)  $ is the fan over $P$.
Denote by $\Delta=\Delta_{-K_{Y}}=P^{\ast}$ the dual polytope of $P$ (which is
not necessarily integral) and by $S$ the Cox ring of $Y$.

Let $\mathfrak{X}\subset Y\times\operatorname*{Spec}\left(  \mathbb{C}\left[
\left[  t\right]  \right]  \right)  $ be a flat family of Calabi-Yau varieties
of dimension $d$ given by the ideal $I\subset\mathbb{C}\left[  t\right]
\otimes S$. Suppose that the special fiber of $\mathfrak{X}$ over the zero
point $\operatorname*{Spec}\left(  \mathbb{C}\right)  $ is given by the
reduced monomial ideal $I_{0}$. We require that the tangent vector of
$\mathfrak{X}$ at $X_{0}$ is sufficiently general in the tangent space of the
component of moduli space of $X_{0}$ containing the family $\mathfrak{X}$.

The goal is to associate to $\mathfrak{X}$ a degeneration $\mathfrak{X}%
^{\circ}$ of Calabi-Yau varieties with fibers in a toric Fano variety such
that the general fibers of $\mathfrak{X}$ and $\mathfrak{X}^{\circ}$ form a
mirror pair.

The presentation of the Chow group of $Y$%
\[
0\rightarrow M\overset{A}{\rightarrow}\mathbb{Z}^{\Sigma\left(  1\right)
}\rightarrow A_{n-1}\left(  Y\right)  \rightarrow0\text{ }%
\]
induces a correspondence of weight vectors on the Cox ring $S$ and the
elements of%
\[
\frac{\operatorname*{Hom}\nolimits_{\mathbb{R}}\left(  \mathbb{R}%
^{\Sigma\left(  1\right)  },\mathbb{R}\right)  }{\operatorname*{Hom}%
\nolimits_{\mathbb{R}}\left(  A_{n-1}\left(  Y\right)  \otimes_{\mathbb{Z}%
}\mathbb{R},\mathbb{R}\right)  }\cong N_{\mathbb{R}}%
\]
This vector space naturally contains the lattice
\[
\operatorname*{image}\left(  \_\circ A\right)  \cong\mathbb{Z}^{n}%
\]

We associate to $\mathfrak{X}$ the special fiber Gr\"{o}bner cone%
\[
C_{0}=C_{I_{0}}\left(  I\right)  \subset\mathbb{R\oplus}N_{\mathbb{R}}%
\]
defined as the closure of the set of weight vectors on $\mathbb{C}\left[
t\right]  \otimes S$ which select $I_{0}$ as initial ideal of $I$. It is a
closed strongly convex rational polyhedral cone.

The cone $C_{0}$ intersects the hyperplane of $t$ -weight $w_{t}=1$, which
contains via stereographic projection the Bergman complex of $I$, in the
polytope%
\[
\nabla=C_{0}\cap\left\{  w_{t}=1\right\}  \subset\left\{  w_{t}=1\right\}
=N_{\mathbb{R}}%
\]
The dual polytope $\nabla^{\ast}$ is integral and contains $0$ as unique
interior point, i.e., $\nabla^{\ast}$ is a Fano polytope, defining a toric
Fano variety $Y^{\circ}=\mathbb{P}\left(  \Sigma^{\circ}\right)  $ by the fan
$\Sigma^{\circ}=\Sigma\left(  \nabla^{\ast}\right)  $ over the faces of
$\nabla^{\ast}$.

The intersection of the Bergman fan with the special fiber Gr\"{o}bner cone%
\[
B\left(  I\right)  =BF\left(  I\right)  \cap C_{0}\cap\left\{  w_{t}%
=1\right\}  \subset\operatorname*{Poset}\left(  \nabla\right)
\]
is a subcomplex of dimension $d$ of the boundary of $\nabla$.

As $B\left(  I\right)  $ is a subset of $\operatorname*{val}\left(
V_{K}\left(  I\right)  \right)  $ for the metric completion $K$ of the field
of Puisseux series, we have a map of complexes%
\[%
\begin{tabular}
[c]{llll}%
$\lim:$ & $B\left(  I\right)  $ & $\rightarrow$ & $\operatorname*{Strata}%
\left(  Y\right)  \cong\operatorname*{Poset}\left(  \Delta\right)  $\\
& \multicolumn{1}{c}{$F$} & $\mapsto$ & $\left\{  \lim_{t\rightarrow0}a\left(
t\right)  \mid a\in\operatorname*{val}^{-1}\left(  \operatorname*{int}\left(
F\right)  \right)  \right\}  $%
\end{tabular}
\
\]
Here $\operatorname*{val}$ is the valuation map defined in Section
\ref{1NonArchimedianAmoebas} and $\operatorname*{Strata}\left(  Y\right)  $ is
the complex of all closures of toric strata of $Y$. Note that
$\operatorname*{Strata}\left(  Y\right)  $ is isomorphic to the complex of
faces $\operatorname*{Poset}\left(  \Delta\right)  $ of $\Delta\subset
M_{\mathbb{R}}$, considered as a complex. The image of the map $\lim$ is the
complex $\operatorname*{Strata}\nolimits_{\Delta}\left(  I_{0}\right)  $ of
strata of $X_{0}$ considered as a subcomplex of $\operatorname*{Poset}\left(
\Delta\right)  $, i.e.,%
\[
\lim\left(  B\left(  I\right)  \right)  =\operatorname*{Strata}%
\nolimits_{\Delta}\left(  I_{0}\right)  \subset\operatorname*{Poset}\left(
\Delta\right)
\]
As a consequence we expect that $B\left(  I\right)  \subset
\operatorname*{Poset}\left(  \nabla\right)  $ is the complex of strata of the
special fiber $X_{0}^{\circ}$ of the mirror degeneration $\mathfrak{X}^{\circ
}$, i.e.,%
\[
B\left(  I\right)  =\operatorname*{Strata}\nolimits_{\Delta}\left(
I_{0}^{\circ}\right)
\]
and its ideal $I_{0}^{\circ}$ is obtained as follows:

Any ray $v$ of the normal fan $\Sigma^{\circ}=\operatorname*{NF}\left(
\nabla\right)  $ corresponds to a facet $F_{v}$ of $\nabla$. Write $S^{\circ
}=\mathbb{C}\left[  z_{v}\mid v\in\Sigma^{\circ}\left(  1\right)  \right]  $
for the Cox ring of $Y^{\circ}$. The subcomplex $B\left(  I\right)
\subset\nabla$ defines a monomial ideal
\[
I_{0}^{\circ}=\left\langle
{\displaystyle\prod\limits_{v\in J}}
z_{v}\mid J\subset\Sigma^{\circ}\left(  1\right)  \text{ with }%
\operatorname*{supp}\left(  B\left(  I\right)  \right)  \subset%
{\displaystyle\bigcup\limits_{v\in J}}
F_{v}\right\rangle \subset S^{\circ}%
\]
generated by the products of variables of $S^{\circ}$ such that the
corresponding union of facets contains $\operatorname*{supp}\left(  B\left(
I\right)  \right)  $ as a subset. Here $\operatorname*{supp}\left(  B\left(
I\right)  \right)  $ denotes the underlying set of the subcomplex $B\left(
I\right)  \subset\nabla$. The special fiber of the mirror degeneration is
expected to be given by
\[
X_{0}^{\circ}=V\left(  I_{0}^{\circ}\right)  \subset Y^{\circ}%
\]

Interpreting $N$ as the lattice of monomials of $Y^{\circ}$, a general first
order polynomial deformation of $I_{0}^{\circ}$, including trivial
deformations, is given as a general linear combination of the lattice points
of
\[
\left(  \lim\left(  B\left(  I\right)  \right)  \right)  ^{\ast}%
\subset\operatorname*{Poset}\left(  \Delta^{\ast}\right)  =P\subset
N_{\mathbb{R}}%
\]
which is the complex of faces $F^{\ast}\subset\Delta^{\ast}$ dual to the faces
$F$ of $\lim\left(  B\left(  I\right)  \right)  $. These lattice points map to
Cox Laurent monomials via the presentation of the Chow group of $Y^{\circ}$%
\[
0\rightarrow N\overset{A^{\circ}}{\rightarrow}\mathbb{Z}^{\Sigma^{\circ
}\left(  1\right)  }\rightarrow A_{n-1}\left(  Y^{\circ}\right)
\rightarrow0\text{ }%
\]
Denote their image by%
\[
\Xi^{\circ}=A^{\circ}\left(  N\cap\operatorname*{supp}\left(  \lim\left(
B\left(  I\right)  \right)  \right)  ^{\ast}\right)
\]

Define the first order deformation of $X_{0}^{\circ}$%
\[
\mathfrak{X}^{1\circ}\subset Y^{\circ}\times\operatorname*{Spec}%
\mathbb{C}\left[  t\right]  /\left\langle t^{2}\right\rangle
\]
by the ideal%
\[
I^{1\circ}=\left\langle u+t\cdot\sum_{\alpha\in\Xi^{\circ}}a_{\alpha}%
\cdot\alpha\left(  u\right)  \mid u\in I_{0}^{\circ}\right\rangle
\subset\mathbb{C}\left[  t\right]  /\left\langle t^{2}\right\rangle \otimes
S^{\circ}%
\]
with general coefficients $a_{\alpha}$. The family $\mathfrak{X}^{1\circ}$ is
expected to be up to first order the mirror degeneration of $\mathfrak{X}%
$.\medskip

This construction is motivated by the following structure on the first order
deformations of $X_{0}$: For any face $F$ of $B\left(  I\right)  $ denote the
associated initial ideal by $\operatorname*{in}_{F}\left(  I\right)  $. For
any tie break ordering $>$ inside $C_{0}$ we have $L_{>}\left(
\operatorname*{in}_{F}\left(  I\right)  \right)  =I_{0}$. Associated to $F$
there is a first order deformation%
\[
\mathfrak{X}_{F}\subset Y\times\operatorname*{Spec}\left(  \mathbb{C}\left[
t\right]  /\left\langle t^{2}\right\rangle \right)
\]
defined by the image of $\operatorname*{in}_{F}\left(  I\right)  $ under%
\[%
\begin{tabular}
[c]{ccc}%
$\mathbb{C}\left[  t\right]  \otimes S$ & $\rightarrow$ & $\mathbb{C}\left[
t\right]  /\left\langle t^{2}\right\rangle \otimes S$\\
$\cup$ &  & $\cup$\\
$\operatorname*{in}_{F}\left(  I\right)  $ & $\rightarrow$ & $\left\langle
m_{i0}+t\sum a_{ij}m_{ij}\mid i\right\rangle $%
\end{tabular}
\
\]
where $m_{10},...,m_{r0}$ are minimal generators of $I_{0}$. By homogeneity
the Cox Laurent monomials $\frac{m_{ij}}{m_{i0}}$ are in the image of $A$, and
the image of the map%
\[%
\begin{tabular}
[c]{llll}%
$\operatorname*{dual}:$ & $B\left(  I\right)  $ & $\rightarrow$ &
$\operatorname*{Poset}\left(  \nabla^{\ast}\right)  $\\
& $F$ & $\mapsto$ & $\operatorname*{convexhull}\left(  \left\{  A^{-1}\left(
\frac{m_{ij}}{m_{i0}}\right)  \mid i,j\right\}  \right)  $%
\end{tabular}
\
\]
associating to each face the convex hull of the first order deformations
appearing in its initial ideal, carries the structure of a complex, indeed
$\operatorname*{dual}\left(  F\right)  =F^{\ast}\subset\nabla^{\ast}$. The
lattice points%
\[
M\cap\operatorname*{supp}\left(  \operatorname*{dual}\left(  B\left(
I\right)  \right)  \right)
\]
of the image form a torus invariant basis of the space of those first order
polynomial deformations of $X_{0}$, which map modulo trivial deformations in
the tangent space direction of the component of the moduli space of $X_{0}$,
containing the family $\mathfrak{X}$.

Above construction depends on the degeneration $\mathfrak{X}$ only up to first
order. Applying the construction to the first order mirror family
$\mathfrak{X}^{\circ}$ recovers the original degeneration $\mathfrak{X}$ up to
first order.\medskip

We summarize the key identifications, made in above construction, denoting the
identifications by the symbols $\updownarrow$ and $\leftrightarrow$. To make a
clear distinction between the two mirror partners and their embedding toric
Fano varieties, denote the lattice of monomials of $Y^{\circ}=X\left(
\Sigma\left(  \nabla^{\ast}\right)  \right)  $ by $M^{\circ}$ and its dual
lattice by $N^{\circ}=\operatorname*{Hom}\left(  M,\mathbb{Z}\right)  $.
Denote by $X$ the general fiber of $\mathfrak{X}$.%
\[%

$}%
\end{tabular}
\]
and the analogous mirrored diagram. The key connections between both diagramms
are made by the maps $\lim$ relating the complexes of strata of $I_{0}$ and
$I_{0}^{\circ}$, and by the map $\operatorname*{dual}$, i.e., the
correspondence between weights and initial ideals.

As discussed in \cite{AGM The MonomialDivisor Mirror Map} in the case of
hypersurfaces, the identification of the lattice $N$ and the lattice of
monomials $M^{\circ}$ of $X\left(  \Sigma^{\circ}\right)  $ gives rise to a
mirror map between complex and K\"{a}hler moduli.

In an analogous way, above tropical mirror construction allows interpretation
of the vertices of the faces (or, via MPCP-blowup, of the lattice points) of
$\operatorname*{dual}\left(  B\left(  I^{\circ}\right)  \right)  $ as first
order polynomial deformations of $X_{0}^{\circ}$ or as toric K\"{a}hler
classes on $X$, which should induce a mirror correspondence between complex
moduli and K\"{a}hler moduli.

\subsection{First order deformations
\index{deformation}%
and degree $0$ Cox Laurent
monomials\label{first order deformations and degree 0 Cox Laurent monomials}}

Consider a toric variety $X\left(  \Sigma\right)  $ of dimension $n$ with
\index{Cox ring}%
Cox ring $S$, and recall that the map $\deg$ in
\[
0\rightarrow M\overset{A}{\rightarrow}\mathbb{Z}^{\Sigma\left(  1\right)
}\overset{\deg}{\rightarrow}A_{n-1}\left(  X\left(  \Sigma\right)  \right)
\rightarrow0
\]
can be considered as the map associating to a Cox
\index{Laurent monomial}%
Laurent monomial its degree in the
\index{Chow group}%
Chow group of divisors $A_{n-1}\left(  X\left(  \Sigma\right)  \right)  $.
Hence $\operatorname*{image}\left(  A\right)  =\ker\left(  \deg\right)  $ is
precisely the set of degree $0$ Cox
\index{Laurent monomial}%
Laurent monomials. So there is an isomorphism%
\[
M\overset{A}{\rightleftarrows}\operatorname*{image}\left(  A\right)
\subset\mathbb{Z}^{\Sigma\left(  1\right)  }%
\]
of $M$ and the degree $0$ Cox Laurent monomials, and $M_{\mathbb{R}}%
\subset\mathbb{R}^{\Sigma\left(  1\right)  }$ is a sub vector space containing
the lattice $M\cong\mathbb{Z}^{n}$.

The characters of the big torus $\left(  \mathbb{C}^{\ast}\right)
^{\Sigma\left(  1\right)  }\cong\operatorname*{Hom}\nolimits_{\mathbb{Z}%
}\left(  \mathbb{Z}^{\Sigma\left(  1\right)  },\mathbb{C}^{\ast}\right)  $ are
the Cox Laurent monomials, i.e., the elements of $\mathbb{Z}^{\Sigma\left(
1\right)  }$.

Let $I_{0}\subset S$ be a monomial ideal. As $I_{0}$ is generated by finitely
many elements and the space of elements of $S$ of this degree is finite
dimensional, the degree $0$ homomorphisms in $\operatorname*{Hom}\left(
I_{0},S/I_{0}\right)  $ form a finite dimensional vector space denoted by
$\operatorname*{Hom}\left(  I_{0},S/I_{0}\right)  _{0}$. The big torus
$\left(  \mathbb{C}^{\ast}\right)  ^{\Sigma\left(  1\right)  }$ acts by%
\[%
\begin{tabular}
[c]{lll}%
$\operatorname*{Hom}\nolimits_{\mathbb{Z}}\left(  \mathbb{Z}^{\Sigma\left(
1\right)  },\mathbb{C}^{\ast}\right)  \times\mathbb{C}\left[  \mathbb{Z}%
^{\Sigma\left(  1\right)  }\right]  $ & $\rightarrow$ & $\mathbb{C}\left[
\mathbb{Z}^{\Sigma\left(  1\right)  }\right]  $\\
\multicolumn{1}{c}{$\left(  \lambda,m\right)  $} & $\mapsto$ & $\lambda\left(
m\right)  \cdot m$%
\end{tabular}
\]
on $\mathbb{C}\left[  \mathbb{Z}^{\Sigma\left(  1\right)  }\right]  $ and on
$S$. The induced action of the abelian group $\operatorname*{Hom}%
\nolimits_{\mathbb{Z}}\left(  \mathbb{Z}^{\Sigma\left(  1\right)  }%
,\mathbb{C}^{\ast}\right)  $ on the vector space $\operatorname*{Hom}\left(
I_{0},S/I_{0}\right)  _{0}$ gives a representation
\[
\operatorname*{Hom}\nolimits_{\mathbb{Z}}\left(  \mathbb{Z}^{\Sigma\left(
1\right)  },\mathbb{C}^{\ast}\right)  \rightarrow\operatorname*{GL}\left(
\operatorname*{Hom}\left(  I_{0},S/I_{0}\right)  _{0}\right)
\]
which decomposes into characters, as any irreducible representation of an
abelian group over an algebraically closed field is $1$-dimensional.

So, denoting first order deformations which are characters
\index{T-deformation|textbf}%
as $\left(  \mathbb{C}^{\ast}\right)  ^{\Sigma\left(  1\right)  }%
$\textbf{-deformations}, the vector space $\operatorname*{Hom}\left(
I_{0},S/I_{0}\right)  _{0}$ has a basis of $\left(  \mathbb{C}^{\ast}\right)
^{\Sigma\left(  1\right)  }$-deformations . Any such homomorphism $\delta$ is
represented by a degree $0$ Cox Laurent monomial, i.e., by a character of
$\left(  \mathbb{C}^{\ast}\right)  ^{\Sigma\left(  1\right)  }$. There are
relatively prime monomials $q_{0},q_{1}\in S$ with $\frac{q_{1}}{q_{0}}%
\in\operatorname*{image}\left(  A\right)  $ such that for all minimal
generators $m\in I_{0}$ with $\delta\left(  m\right)  \neq0$ we have
$\frac{\delta\left(  m\right)  }{m}=\frac{q_{1}}{q_{0}}$. So $\delta$ is the
degree $0$ homomorphism $I_{0}\rightarrow S/I_{0}$ defined by%
\[
\delta\left(  m\right)  =\left\{
\begin{tabular}
[c]{ll}%
$\frac{q_{1}}{q_{0}}\cdot m$ & if $q_{0}\mid m$\\
$0$ & otherwise
\end{tabular}
\ \right\}
\]
for minimal generators $m\in I_{0}$.

\begin{lemma}
If $I_{0}$ is a
\index{monomial ideal}%
monomial ideal, then $\operatorname*{Hom}\left(  I_{0},S/I_{0}\right)  _{0}$
has a basis
\index{T-invariant}%
of $\left(  \mathbb{C}^{\ast}\right)  ^{\Sigma\left(  1\right)  }%
$-deformations represented by elements of $\operatorname*{image}\left(
A\right)  \cong M$.
\end{lemma}

With respect to
\index{weight vector}%
weights on $I$ recall from Section \ref{1torichomogeneoussetting} that there
is a bijection%

\[%
\begin{tabular}
[c]{lllll}%
$N\supset\operatorname*{image}\left(  \_\circ A\right)  $ & $\overset{\_\circ
A}{\underset{\varphi}{\leftrightarrows}}$ & $\frac{\operatorname*{Hom}%
\nolimits_{\mathbb{Z}}\left(  \mathbb{Z}^{\Sigma\left(  1\right)  }%
,\mathbb{Z}\right)  }{\operatorname*{Hom}\nolimits_{\mathbb{Z}}\left(
A_{n-1}\left(  X\left(  \Sigma\right)  \right)  ,\mathbb{Z}\right)  }$ &
$\rightarrow$ & $\left\{  \text{graded wt. vec. on }S\right\}  $%
\end{tabular}
\
\]
inducing an isomorphism of vector spaces%
\[
N_{\mathbb{R}}\underset{\varphi_{\mathbb{R}}}{\leftrightarrows}\frac
{\operatorname*{Hom}\nolimits_{\mathbb{R}}\left(  \mathbb{R}^{\Sigma\left(
1\right)  },\mathbb{R}\right)  }{\operatorname*{Hom}\nolimits_{\mathbb{R}%
}\left(  A_{n-1}\left(  X\left(  \Sigma\right)  \right)  \otimes_{\mathbb{Z}%
}\mathbb{R},\mathbb{R}\right)  }%
\]
i.e., $N_{\mathbb{R}}$ is a quotient of $\operatorname*{Hom}%
\nolimits_{\mathbb{R}}\left(  \mathbb{R}^{\Sigma\left(  1\right)  }%
,\mathbb{R}\right)  $.

The mirror correspondence between Calabi-Yau degenerations with fibers
polarized in toric Fano varieties $X\left(  \Sigma\right)  $ with lattices $N$
and $M$ respectively $X\left(  \Sigma^{\circ}\right)  $ with lattices
$N^{\circ}$ and $M^{\circ}$ will be induced by the identification of%
\[%
\begin{tabular}
[c]{ccccccc}%
$\frac{\operatorname*{Hom}\nolimits_{\mathbb{R}}\left(  \mathbb{R}%
^{\Sigma\left(  1\right)  },\mathbb{R}\right)  }{\operatorname*{Hom}%
\nolimits_{\mathbb{R}}\left(  A_{n-1}\left(  X\left(  \Sigma\right)  \right)
\otimes_{\mathbb{Z}}\mathbb{R},\mathbb{R}\right)  }$ & $\cong$ &
$N_{\mathbb{R}}$ & with & $M_{\mathbb{R}}^{\circ}$ & $\subset$ &
$\mathbb{R}^{\Sigma^{\circ}\left(  1\right)  }$%
\end{tabular}
\]
and of%
\[%
\begin{tabular}
[c]{ccccccl}%
$\mathbb{R}^{\Sigma\left(  1\right)  }$ & $\supset$ & $M_{\mathbb{R}}$ &
with & $N_{\mathbb{R}}^{\circ}$ & $\cong$ & $\frac{\operatorname*{Hom}%
\nolimits_{\mathbb{R}}\left(  \mathbb{R}^{\Sigma^{\circ}\left(  1\right)
},\mathbb{R}\right)  }{\operatorname*{Hom}\nolimits_{\mathbb{R}}\left(
A_{n-1}\left(  X\left(  \Sigma^{\circ}\right)  \right)  \otimes_{\mathbb{Z}%
}\mathbb{R},\mathbb{R}\right)  }$%
\end{tabular}
\]
and of the corresponding lattices.

\subsection{Monomial ideals in the Cox ring and the stratified toric primary
decomposition\label{Monomial ideals in the Cox ring and the stratified toric primary decomposition}%
}

Let $N\cong\mathbb{Z}^{n}$, let $M=\operatorname*{Hom}\left(  N,\mathbb{Z}%
\right)  $ be the dual lattice, $\Sigma\subset N_{\mathbb{R}}$ a complete fan,
$Y=X\left(  \Sigma\right)  $ the corresponding toric variety and
$S=\mathbb{C}\left[  y_{r}\mid r\in\Sigma\left(  1\right)  \right]  $ the Cox
ring of $Y$.

\begin{definition}
Let $I_{0}\subset S$ be a reduced
\newsym[$\operatorname*{rays}\nolimits_{m}\left(  \Sigma\right)  $]{rays of a monomial}{}monomial
ideal. If $m\in I_{0}$ is a monomial, define%
\[
\operatorname*{rays}\nolimits_{m}\left(  \Sigma\right)  =\left\{  r\in
\Sigma\left(  1\right)  \mid y_{r}\text{ divides }m\right\}
\]
The
\index{stratified toric primary decomposition}%
\textbf{stratified toric primary decomposition }$SP\left(  I_{0}\right)  $ of
$I_{0}$ is the complex with
\newsym[$SP$]{stratified toric primary decomposition}{}faces of dimension $s$
given by%
\[
SP\left(  I_{0}\right)  _{s}=\left\{  \left\langle y_{r}\mid r\in
\operatorname*{rays}\left(  \sigma\right)  \right\rangle \mid%
\begin{tabular}
[c]{l}%
$\sigma\in\Sigma\left(  n-s\right)  \text{ with }$\\
$\operatorname*{rays}\left(  \sigma\right)  \cap\operatorname*{rays}%
\nolimits_{m}\left(  \Sigma\right)  \neq\varnothing\text{ }$\\
for all monomials $m\in I_{0}$%
\end{tabular}
\right\}
\]

\end{definition}

\begin{remark}
Suppose that all maximal faces $SP\left(  I_{0}\right)  $ have the same
\newsym[$IS$]{intersection complex}{}dimension, i.e., the vanishing locus of
$I_{0}$ is equidimensional. The
\index{intersection complex|textbf}%
\textbf{intersection complex }$IS\left(  I_{0}\right)  $\textbf{ }of $I_{0}$
is the subcomplex of the simplex on the maximal faces $SP\left(  I_{0}\right)
$, containing the face $F$ if the ideal%
\[
\sum\limits_{J\in F}J\in SP\left(  I_{0}\right)
\]
i.e., if the ideal $\sum\nolimits_{J\in F}J$ is again a face of $SP\left(
I_{0}\right)  $. The complexes $SP\left(  I_{0}\right)  $ and $IS\left(
I_{0}\right)  $ are dual to each other.
\end{remark}

Suppose $D$ is a divisor on $Y=X\left(  \Sigma\right)  $ such that some
multiple of $D$ is ample Cartier, then $\Delta=\Delta_{D}$ is not necessarily
integral, but combinatorially dual to $\Sigma$, i.e., $\Sigma
=\operatorname*{NF}\left(  \Delta_{D}\right)  $.

For example we could consider a Fano polytope $P\subset N_{\mathbb{R}}$,
$\Sigma=\Sigma\left(  P\right)  $ the fan over $P$ and $Y=X\left(
\Sigma\right)  $ the corresponding toric Fano variety and $\Delta
=\Delta_{-K_{Y}}=P^{\ast}\subset M_{\mathbb{R}}$.

We can reformulate above notations in terms of a subcomplex of the polytope
$\Delta$ and the dimensions of the faces are the geometric dimension of the
corresponding faces of $\Delta$:

If $F$ is a face of $\Delta$, define%
\[
\operatorname*{facets}\nolimits_{F}\left(  \Delta\right)  =\left\{  G\mid
G\text{ facet of }\Delta\text{ with }F\subset G\right\}
\]
as the set of facets of $\Delta$ containing $F$. If $m\in I_{0}$ is a
monomial, define%
\[
\operatorname*{facets}\nolimits_{m}\left(  \Delta\right)  =\left\{  G\mid
G\text{ facet of }\Delta\text{ with }y_{G^{\ast}}\mid m\right\}
\]
as the set of those facets of $\Delta$, which appear, considered as Cox
variable, as a factor of $m$.

\begin{definition}
The
\index{complex of strata|textbf}%
\textbf{complex of strata} of $I_{0}$ is the subcomplex
$\operatorname*{Strata}_{\Delta}\left(  I_{0}\right)  $ of the associated
complex $\operatorname*{Poset}\left(  \Delta\right)  $ of $\Delta$ with
\newsym[$\operatorname*{Strata}$]{strata of a monomial ideal}{}faces of
$\operatorname*{Strata}_{\Delta}\left(  I_{0}\right)  $ of dimension $s$ given
by%
\[
\operatorname*{Strata}\nolimits_{\Delta}\left(  I_{0}\right)  _{s}=\left\{
F\mid%
\begin{tabular}
[c]{l}%
$F\text{ a face of }\Delta\text{ of }\dim\left(  F\right)  =s\text{ with}$\\
$\operatorname*{facets}\nolimits_{F}\left(  \Delta\right)  \cap
\operatorname*{facets}\nolimits_{m}\left(  \Delta\right)  \neq\varnothing$\\
for all monomials $m\in I_{0}$%
\end{tabular}
\right\}
\]

\end{definition}

\begin{lemma}
The faces of dimension $s$ of the stratified toric primary decomposition
$SP\left(  I_{0}\right)  $ of $I_{0}$ are given by
\[
SP\left(  I_{0}\right)  _{s}=\left\{  \left\langle y_{G^{\ast}}\mid G\text{ a
facet of }\Delta\text{ with }F\subset G\right\rangle \mid F\in
\operatorname*{Strata}\nolimits_{\Delta}\left(  I_{0}\right)  _{s}\right\}
\]
and $\operatorname*{Strata}\nolimits_{\Delta}\left(  I_{0}\right)  \cong
SP\left(  I_{0}\right)  $.
\end{lemma}

These definitions may be generalized to the case of non-reduced monomial
ideals, though this is not used in the following.

\begin{proposition}
Let $I_{0}\subset S$ be a reduced monomial ideal such that
$\operatorname*{Strata}\nolimits_{\Delta}\left(  I_{0}\right)  $ is
equidimensional of dimension $d$. Then there is a unique monomial ideal
$I_{0}^{\Sigma}\subset S$ maximal with respect to inclusion such that
$\operatorname*{Strata}\nolimits_{\Delta}\left(  I_{0}\right)
=\operatorname*{Strata}\nolimits_{\Delta}\left(  I_{0}^{\Sigma}\right)  $. It
holds%
\begin{align*}
I_{0}^{\Sigma}  &  =%
{\textstyle\bigcap\nolimits_{F\in\operatorname*{Strata}\nolimits_{\Delta
}\left(  I_{0}\right)  _{d}}}
\left\langle y_{G^{\ast}}\mid G\text{ a facet of }\Delta\text{ with }F\subset
G\right\rangle \\
&  =\left\langle
{\displaystyle\prod\limits_{v\in J}}
y_{v}\mid J\subset\Sigma\left(  1\right)  \text{ with }\operatorname*{supp}%
\left(  \operatorname*{Strata}\nolimits_{\Delta}\left(  I_{0}\right)  \right)
\subset%
{\displaystyle\bigcup\limits_{v\in J}}
F_{v}\right\rangle \subset S
\end{align*}

\end{proposition}

\begin{definition}
We denote $I_{0}^{\Sigma}$ as the $\Sigma$\textbf{-saturation} of $I_{0}$.
\end{definition}

\begin{remark}
If $\Sigma$ is simplicial, then%
\[
I_{0}^{\Sigma}=\left(  I_{0}:B\left(  \Sigma\right)  ^{\infty}\right)
\]

\end{remark}

In the special case of $Y=\mathbb{P}^{n}$ the complex $\operatorname*{Strata}%
\nolimits_{\Delta}\left(  I_{0}\right)  $ is related to the representation of
Stanley-Reisner ideals by the following remark (see also Section
\ref{Sec deformations and obstructions of a toric generalization of stanley reisner rings}%
):

\begin{remark}
Suppose $Y=\mathbb{P}\left(  \Delta\right)  \cong\mathbb{P}^{n}$ where
$\Delta$ is the degree $n+1$ Veronese polytope and let $S$ be the homogeneous
coordinate ring of $Y$. So $\Delta^{\ast}$ (and of course also $\Delta$) is a
simplex and the faces of $\Delta^{\ast}$ correspond to the subsets of the set
of vertices of $\Delta^{\ast}$. The vertices of $\Delta^{\ast}$ generate the
rays of $\Sigma=\operatorname*{NF}\left(  \Delta\right)  $, the cones of
$\Sigma$ correspond to the subsets of $\Sigma\left(  1\right)  $. The rays of
$\Sigma$ correspond to the variables of $S=\mathbb{C}\left[  y_{r}\mid
r\in\Sigma\left(  1\right)  \right]  $.

Let $Z$ be a simplicial subcomplex of $\operatorname*{Poset}\left(
\Sigma\right)  \cong\operatorname*{Poset}\left(  \Delta^{\ast}\right)  $,
where each face of $Z$ is considered as a set of rays of $\Sigma$, and%
\[
I_{0}=\left\langle
{\textstyle\prod\nolimits_{r\in M}}
y_{r}\mid M\subset\Sigma\left(  1\right)  \text{ a non-face of }Z\right\rangle
\subset S
\]
the corresponding
\index{Stanley-Reisner ideal|textbf}%
\textbf{Stanley-Reisner ideal}.

If $F\in\operatorname*{Strata}\nolimits_{\Delta}\left(  I_{0}\right)  $ is a
face and $F^{\ast}\subset\Delta^{\ast}$ the dual face of $F$ then denote the
set of all rays of $\Sigma$ in the complement of $\operatorname*{hull}\left(
F^{\ast}\right)  \in\Sigma$ by $\operatorname*{comp}\left(  F\right)  $, so,
e.g., if $F$ is a vertex of $\Delta$ then the complement of $F^{\ast}%
\subset\Delta^{\ast}$ contains precisely one vertex of $\Delta^{\ast}$. The
map%
\[%
\begin{tabular}
[c]{cccccc}
& $\operatorname*{Poset}\left(  \Delta\right)  $ &  & $\operatorname*{Poset}%
\left(  \Sigma\right)  $ & $\cong$ & $\operatorname*{Poset}\left(
\Delta^{\ast}\right)  $\\
& $\cup$ &  & $\cup$ &  & \\
$\operatorname*{comp}:$ & $\operatorname*{Strata}\nolimits_{\Delta}\left(
I_{0}\right)  $ & $\overset{\cong}{\rightarrow}$ & $Z$ &  & \\
& $F$ & $\mapsto$ & $\operatorname*{comp}\left(  F\right)  $ &  &
\end{tabular}
\]
is an isomorphism of complexes and%
\[
I_{0}=\left\langle
{\displaystyle\prod\limits_{v\in J}}
y_{v}\mid J\subset\Sigma\left(  1\right)  \text{ with }\operatorname*{supp}%
\left(  \operatorname*{Strata}\nolimits_{\Delta}\left(  I_{0}\right)  \right)
\subset%
{\displaystyle\bigcup\limits_{v\in J}}
F_{v}\right\rangle
\]

\end{remark}

\begin{figure}
[h]
\begin{center}
\includegraphics[
height=2.3826in,
width=3.0727in
]%
{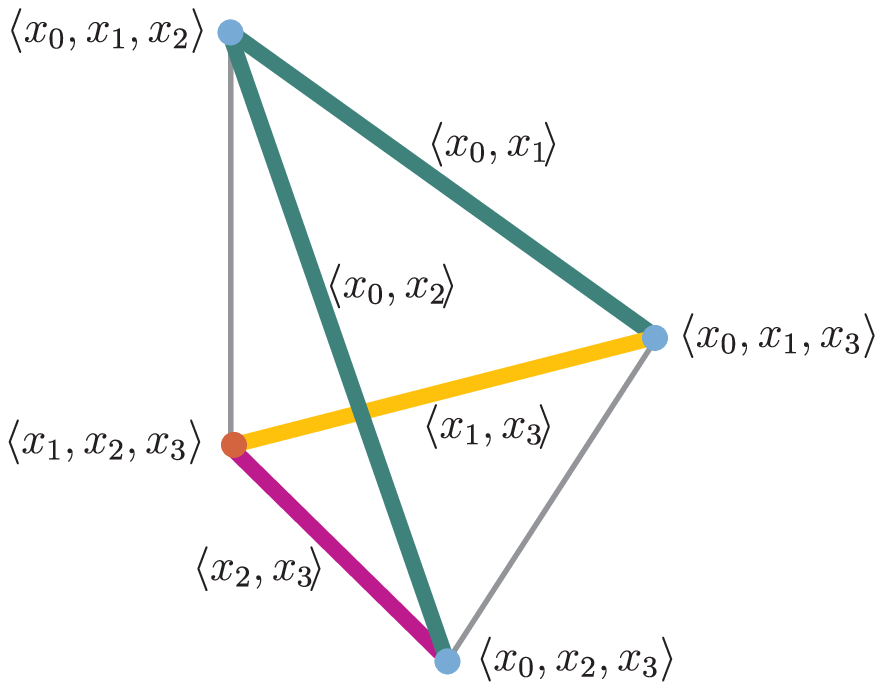}%
\caption{The complex $\operatorname*{Strata}\nolimits_{\Delta}\left(
I_{0}\right)  \subset\operatorname*{Poset}\left(  \Delta\right)  $ for the
ideal $I_{0}=\left\langle x_{0}x_{3},x_{1}x_{2}\right\rangle $}%
\label{Fig strata delta 22}%
\end{center}
\end{figure}

\begin{figure}
[hh]
\begin{center}
\includegraphics[
height=1.6821in,
width=1.9095in
]%
{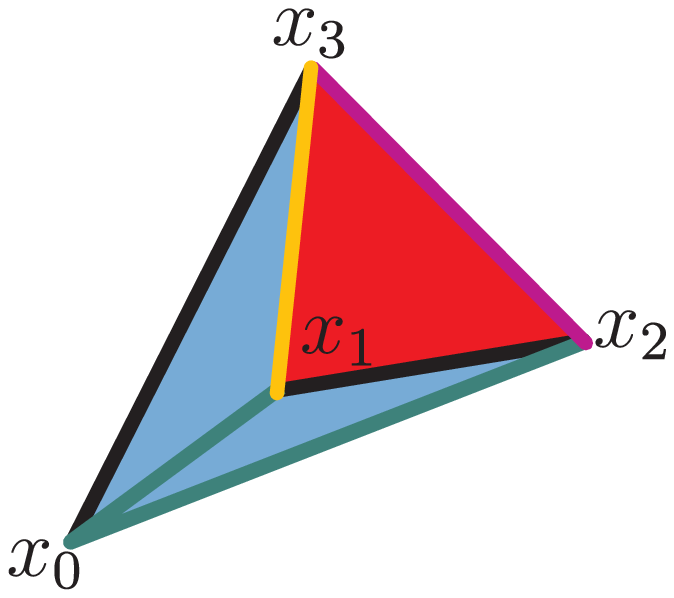}%
\caption{The complex $\left(  \operatorname*{Strata}\nolimits_{\Delta}\left(
I_{0}\right)  \right)  ^{\ast}\subset\operatorname*{Poset}\left(  \Delta
^{\ast}\right)  $ for the ideal $I_{0}=\left\langle x_{0}x_{3},x_{1}%
x_{2}\right\rangle $}%
\label{Fig dual strata delta 22}%
\end{center}
\end{figure}

\begin{figure}
[hhh]
\begin{center}
\includegraphics[
height=1.7288in,
width=1.951in
]%
{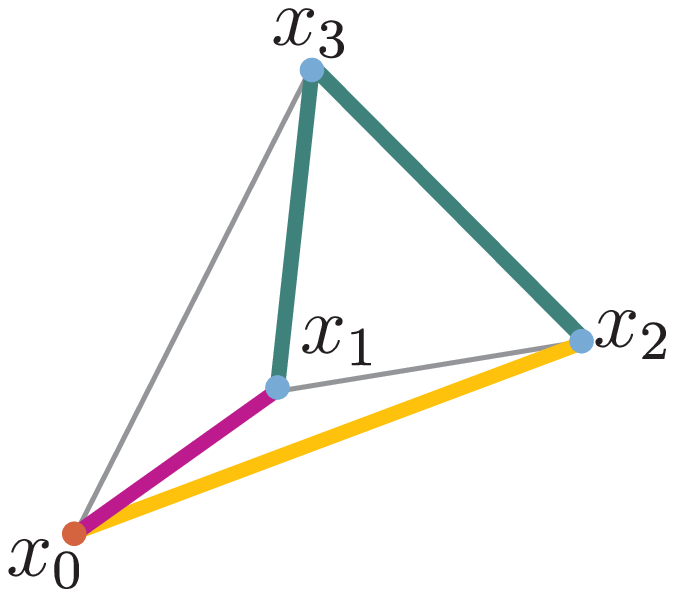}%
\caption{The subcomplex of $\operatorname*{Poset}\left(  \Delta^{\ast}\right)
$ defining the Stanley-Reisner ideal $I_{0}=\left\langle x_{0}x_{3},x_{1}%
x_{2}\right\rangle $}%
\label{Fig SR complex 22}%
\end{center}
\end{figure}

\begin{example}
Let $Y=\mathbb{P}\left(  \Delta\right)  \cong\mathbb{P}^{3}$ with%
\[
\Delta=\operatorname*{convexhull}\left\{  \left(  -1,-1,-1\right)  ,\left(
3,-1,-1\right)  ,\left(  -1,3,-1\right)  ,\left(  -1,-1,3\right)  \right\}
\]
so%
\[
\Delta^{\ast}=\operatorname*{convexhull}\left\{  \left(  -1,-1,-1\right)
,\left(  1,0,0\right)  ,\left(  0,1,0\right)  ,\left(  0,0,1\right)  \right\}
\]
and write%
\[%

\]

Consider the complex%
\[
Z=\left(
%
\right)
\]
The Stanley-Reisner ideal associated to $Z$ is the monomial ideal%
\begin{align*}
I_{0}  &  =\left\langle x_{0}x_{3},x_{1}x_{2},x_{0}x_{1}x_{2},x_{0}x_{1}%
x_{3},x_{0}x_{2}x_{3},x_{1}x_{2}x_{3},x_{0}x_{1}x_{2}x_{3}\right\rangle \\
&  =\left\langle x_{0}x_{3},x_{1}x_{2}\right\rangle \subset S=\mathbb{C}%
\left[  x_{0},x_{1},x_{2},x_{3}\right]
\end{align*}

The complex of strata associated to $I_{0}$ is%
\[
\operatorname*{Strata}\nolimits_{\Delta}\left(  I_{0}\right)
=\operatorname*{convexhull}\left(
%
\right)
\]
This notation is short for applying $\operatorname*{convexhull}$ to each face
of the complex in the argument. The dual of $\operatorname*{Strata}%
\nolimits_{\Delta}\left(  I_{0}\right)  $ is%
\[
\left(  \operatorname*{Strata}\nolimits_{\Delta}\left(  I_{0}\right)  \right)
^{\ast}=\operatorname*{convexhull}\left(
%
\right\}  \smallskip\\
\left\{  {}\right\}  \smallskip\\
\left\{  {}\right\}
\end{array}
\right) \\
&  =Z
\end{align*}

Figure \ref{Fig strata delta 22} shows the complexes $\operatorname*{Strata}%
\nolimits_{\Delta}\left(  I_{0}\right)  \subset\operatorname*{Poset}\left(
\Delta\right)  $, Figure \ref{Fig dual strata delta 22} the complex $\left(
\operatorname*{Strata}\nolimits_{\Delta}\left(  I_{0}\right)  \right)  ^{\ast
}\subset\operatorname*{Poset}\left(  \Delta^{\ast}\right)  $ and Figure
\ref{Fig SR complex 22} the corresponding Stanley-Reisner complex $Z$
considered as a subcomplex of $\operatorname*{Poset}\left(  \Delta^{\ast
}\right)  $.
\end{example}

%


\subsection{Locally relevant
deformations\label{Sec locally relevant deformations}}

Let $N\cong\mathbb{Z}^{n}$ and $M=\operatorname*{Hom}\left(  N,\mathbb{Z}%
\right)  $, let $\Sigma\subset N_{\mathbb{R}}$ be a complete fan and
$Y=X\left(  \Sigma\right)  $ the associated toric variety.

\begin{definition}
Let $X_{0}\subset Y$ be a union of equidimensional strata, let $X_{i}%
\in\operatorname*{Strata}\left(  X_{0}\right)  $ be a torus stratum of $X_{0}$
and consider a first order deformation $\mathfrak{X}\subset Y\times
\operatorname*{Spec}\left(  \mathbb{C}\left[  t\right]  /\left\langle
t^{2}\right\rangle \right)  $ of $X_{0}$. Then $\mathfrak{X}$ is called
\index{locally irrelevant|textbf}%
\textbf{locally irrelevant} at the stratum $X_{i}$ if there is a formal
analytic open neighborhood $\tilde{U}\subset Y$ of $X_{i}$ such that for the
open neighborhood $U=\tilde{U}\cap X_{0}$ of $X_{i}$ in $X_{0}$ there is an
isomorphism%
\[
U\times\operatorname*{Spec}\left(  \mathbb{C}\left[  t\right]  /\left\langle
t^{2}\right\rangle \right)  \cong\mathfrak{X}\cap\left(  \tilde{U}%
\times\operatorname*{Spec}\left(  \mathbb{C}\left[  t\right]  /\left\langle
t^{2}\right\rangle \right)  \right)
\]
which extends
\[
X_{i}\times\operatorname*{Spec}\left(  \mathbb{C}\left[  t\right]
/\left\langle t^{2}\right\rangle \right)  \subset\mathfrak{X}%
\]
Otherwise, $\mathfrak{X}$ is called
\index{locally relevant|textbf}%
\textbf{locally relevant} at $X_{i}$. The deformation $\mathfrak{X}$ is
called
\index{strongly locally relevant}%
\textbf{strongly locally relevant} at $X_{i}$ if $\mathfrak{X}$ is locally
relevant at $X_{i}$ and locally irrelevant for all strata $X_{j}%
\in\operatorname*{Strata}\left(  X_{0}\right)  $ with $X_{i}\cap
X_{j}=\varnothing$.
\end{definition}

\begin{figure}
[h]
\begin{center}
\includegraphics[
height=2.5702in,
width=2.8746in
]%
{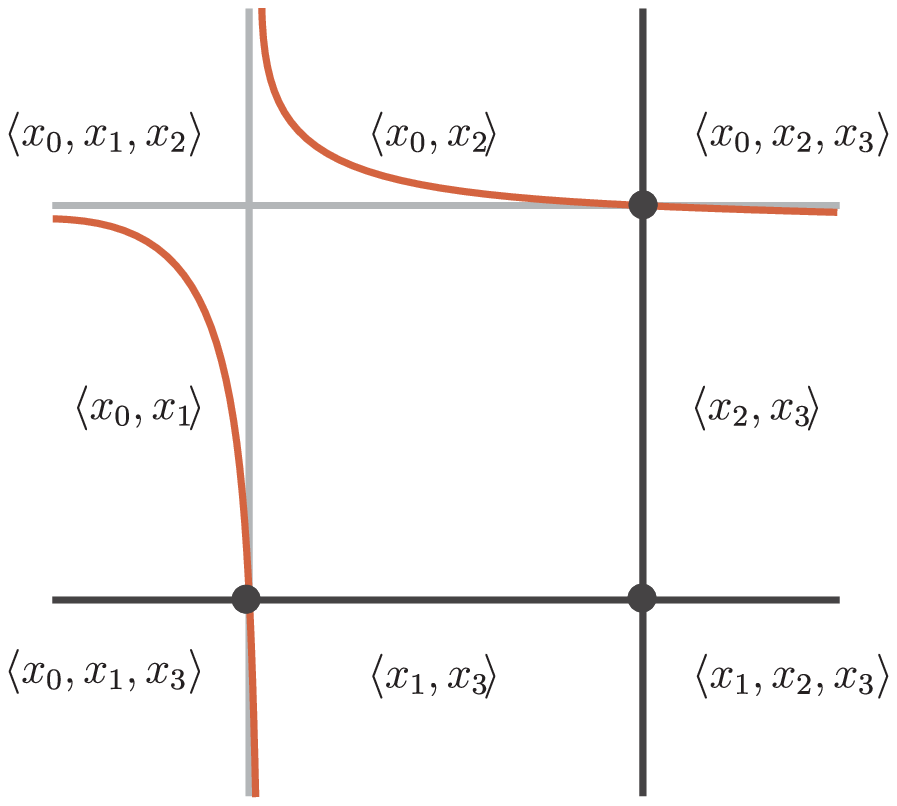}%
\caption{Visualization of the deformation $\frac{x_{3}^{2}}{x_{1}x_{2}}$ of
$I_{0}=\left\langle x_{0}x_{3},x_{1}x_{2}\right\rangle $}%
\label{Fig def x32x1x2}%
\end{center}
\end{figure}

\begin{figure}
[hh]
\begin{center}
\includegraphics[
height=2.4033in,
width=2.8746in
]%
{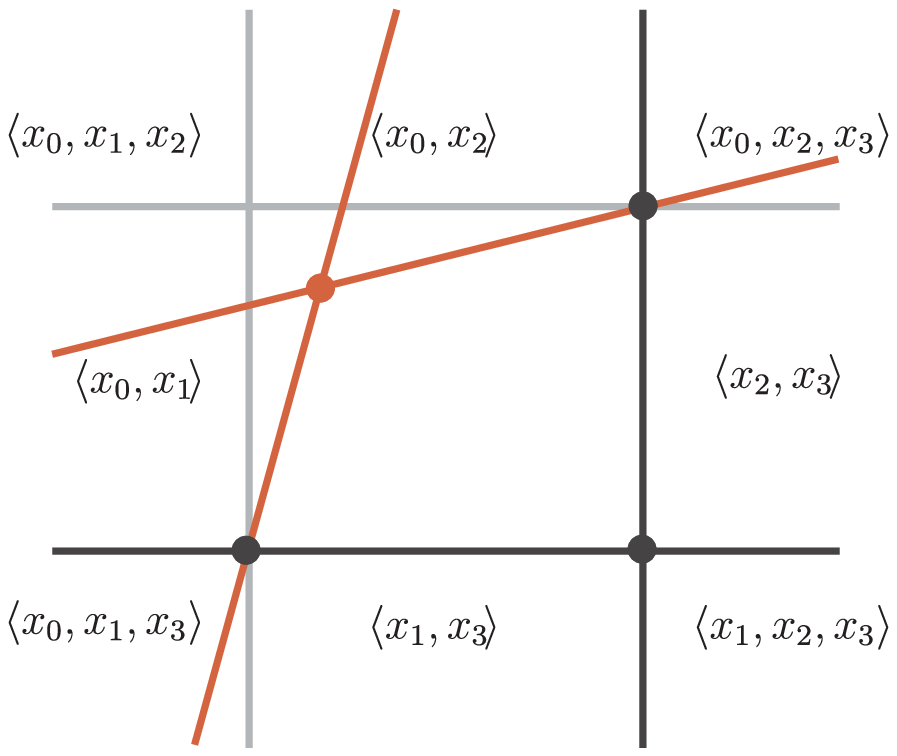}%
\caption{Visualization of the deformation $\frac{x_{3}}{x_{0}}$ of
$I_{0}=\left\langle x_{0}x_{3},x_{1}x_{2}\right\rangle $}%
\label{Fig def x3x0}%
\end{center}
\end{figure}

\begin{figure}
[hhh]
\begin{center}
\includegraphics[
height=2.3903in,
width=2.8746in
]%
{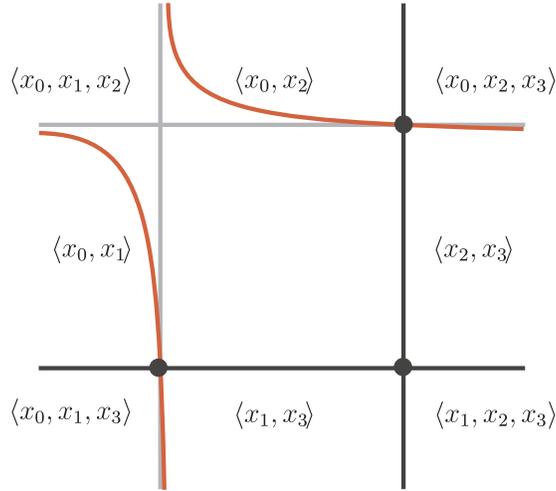}%
\caption{Visualization of the deformation $\frac{x_{1}}{x_{0}}$ of
$I_{0}=\left\langle x_{0}x_{3},x_{1}x_{2}\right\rangle $}%
\label{Fig x1x0}%
\end{center}
\end{figure}

\begin{example}
Consider $X_{0}\subset Y=\mathbb{P}^{3}$ given by the monomial ideal
$I_{0}=\left\langle x_{0}x_{3},x_{1}x_{2}\right\rangle \subset S=\mathbb{C}%
\left[  x_{0},...,x_{3}\right]  $. The ideals of the strata of $X_{0}$ are
shown in Figure \ref{Fig strata delta 22}. Consider the following torus
invariant deformations given by Cox Laurent monomials:

\begin{itemize}
\item For $\frac{x_{3}^{2}}{x_{1}x_{2}}$ the deformation $\mathfrak{X}$ is
given by%
\[
\left\langle x_{0}x_{3},x_{1}x_{2}+t\cdot x_{3}^{2}\right\rangle =\left\langle
x_{3},x_{1}\right\rangle \cap\left\langle x_{3},x_{2}\right\rangle
\cap\left\langle x_{0},x_{1}x_{2}+t\cdot x_{3}^{2}\right\rangle
\]
and is locally relevant at $\left\langle x_{0},x_{2}\right\rangle
,\left\langle x_{0},x_{1}\right\rangle ,\left\langle x_{0},x_{1}%
,x_{2}\right\rangle $. It is strongly locally relevant at $\left\langle
x_{0},x_{1},x_{2}\right\rangle $.

\item For $\frac{x_{3}}{x_{0}}$ the deformation $\mathfrak{X}$ is given by%
\[
\left\langle x_{0}x_{3}+t\cdot x_{3}^{2},x_{1}x_{2}\right\rangle =\left\langle
x_{3},x_{1}\right\rangle \cap\left\langle x_{3},x_{2}\right\rangle
\cap\left\langle x_{1},x_{0}+t\cdot x_{3}\right\rangle \cap\left\langle
x_{2},x_{0}+t\cdot x_{3}\right\rangle
\]
and is locally relevant at $\left\langle x_{0},x_{2}\right\rangle
,\left\langle x_{0},x_{1}\right\rangle ,\left\langle x_{0},x_{1}%
,x_{2}\right\rangle $. It is strongly locally relevant at $\left\langle
x_{0},x_{1},x_{2}\right\rangle $.

\item For $\frac{x_{1}}{x_{0}}$ the deformation $\mathfrak{X}$ is given by%
\[
\left\langle x_{0}x_{3}+t\cdot x_{1}x_{3},x_{1}x_{2}\right\rangle
=\left\langle x_{3},x_{1}\right\rangle \cap\left\langle x_{3},x_{2}%
\right\rangle \cap\left\langle x_{1},x_{0}\right\rangle \cap\left\langle
x_{2},x_{0}+t\cdot x_{1}\right\rangle
\]
and is locally relevant at $\left\langle x_{0},x_{1},x_{2}\right\rangle
,\left\langle x_{0},x_{2},x_{3}\right\rangle ,\left\langle x_{0}%
,x_{2}\right\rangle ,\left\langle x_{1},x_{0}\right\rangle ,\left\langle
x_{3},x_{2}\right\rangle $. It is strongly locally relevant at $\left\langle
x_{0},x_{2}\right\rangle $.
\end{itemize}

The Figures \ref{Fig def x32x1x2}, \ref{Fig def x3x0} and \ref{Fig x1x0}
visualize the strata of $X_{0}$ where these deformations are locally relevant
or irrelevant.
\end{example}

%

See also Example \ref{Ex locally relevant} below.

\subsection{Setup for the tropical mirror construction for monomial
degenerations of Calabi-Yau varieties polarized in toric Fano
varieties\label{genericy condition}}

Consider the following setup for the tropical mirror construction. The
conditions on the degeneration may be subject to generalization and redundancy.

We begin with the following setup:

\begin{itemize}
\item Let $N\cong\mathbb{Z}^{n}$, let $M=\operatorname*{Hom}\left(
N,\mathbb{Z}\right)  $ be the dual lattice, $P$ a Fano polytope,
$\Sigma=\Sigma\left(  P\right)  $ and $Y=X\left(  \Sigma\right)  $ the
corresponding toric Fano variety. Denote by $S=\mathbb{C}\left[  y_{r}\mid
r\in\Sigma\left(  1\right)  \right]  $ the Cox ring of $Y$ and by%
\[
0\rightarrow M\overset{A}{\rightarrow}\mathbb{Z}^{\Sigma\left(  1\right)
}\rightarrow A_{n-1}\left(  Y\right)  \rightarrow0\text{ }%
\]
the presentation of the Chow group of divisors of $Y$.

\item Let $I_{0}\subset S$ be an equidimensional reduced monomial ideal with
$I_{0}=I_{0}^{\Sigma}$ and $\mathfrak{X}\subset Y\times\operatorname*{Spec}%
\mathbb{C}\left[  \left[  t\right]  \right]  $ an irreducible flat family of
Calabi-Yau varieties of dimension $d$ given by the ideal $I\subset
\mathbb{C}\left[  t\right]  \otimes S$ and suppose that the special fiber
$X_{0}$ over $\operatorname*{Spec}\mathbb{C}$ is given by $I_{0}$.

Assume that the underlying topological space of the
\index{cell complex}%
cell complex $\operatorname*{Strata}_{\Delta}\left(  I_{0}\right)  $ is
homeomorphic to a sphere.
\end{itemize}

Define the following:

\begin{itemize}
\item Let $>$ be a monomial ordering on $\mathbb{C}\left[  t\right]  \otimes
S$, which is respecting the Chow
\index{grading}%
grading on $S$ and is local in $t$. Let%
\[
C_{I_{0}}\left(  I\right)  =\left\{  -\left(  w_{t},w_{y}\right)
\in\mathbb{R}\oplus N_{\mathbb{R}}\mid L_{>_{\left(  w_{t},\varphi\left(
w_{y}\right)  \right)  }}\left(  I\right)  =I_{0}\right\}
\]
be Gr\"{o}bner cone corresponding to the lead ideal $I_{0}$.

\item Let%
\[
BF_{I_{0}}\left(  I\right)  =BF\left(  I\right)  \cap\operatorname*{Poset}%
\left(  C_{I_{0}}\left(  I\right)  \right)
\]

\end{itemize}

\noindent We require $\mathfrak{X}$ to satisfy the following conditions:

\begin{enumerate}
\item $C_{I_{0}}\left(  I\right)  \cap\left\{  w_{t}=0\right\}  =\left\{
0\right\}  $.

\item $C_{I_{0}}\left(  I\right)  $ is the cone defined by the half-space
equations corresponding to the torus invariant first order deformations
appearing in the reduced standard basis of $I$ in $S\times\mathbb{C}\left[
t\right]  /\left\langle t^{2}\right\rangle $ with respect to a monomial
ordering in the interior of $C_{I_{0}}\left(  I\right)  $.

All lattice points of $F^{\ast}$ appear as deformations in $I$.

\item Any first order deformation appearing in $I$ is also a non-zero
deformation of the anticanonical Calabi-Yau hypersurface in $Y$.

\item Any facet of $\operatorname*{Strata}\nolimits_{\Delta}\left(
I_{0}\right)  $ is contained in precisely $c$ facets of $\Delta=P^{\ast}$.

\item Any facet of $BF_{I_{0}}\left(  I\right)  $ is contained in precisely
$c$ facets of $\operatorname*{Poset}\left(  C_{I_{0}}\left(  I\right)
\right)  $.
\end{enumerate}

\noindent In the following we give a geometric interpretation of these conditions.

\begin{enumerate}
\item We can satisfy requirement $1.$ via a condition on the position of the
Hilbert point of $I_{0}$ with respect to the state polytope of the general fiber:

Let $\mathcal{K}=\operatorname{cpl}\left(  \Sigma\right)  \cap
\operatorname*{Pic}\left(  Y\right)  $ and $I_{gen}\subset S$ be the saturated
ideal of the general fiber of $\mathfrak{X}$. Let $P\left(  t\right)  $ be
the\ Hilbert polynomial of $I_{gen}$, $h$ the corresponding Hilbert function,
$D\subset m+\mathcal{K}$ such that the restriction map gives a closed
embedding $\mathbb{H}_{\left(  S,F\right)  }^{h}\rightarrow\mathbb{H}_{\left(
S_{D},F_{D}\right)  }^{h}$. Fix linearizations of the action of $T$ on the
elements of $D$. We require that the Hilbert point
\[
H\left(  I_{0}\right)  \in\operatorname*{int}\left(  \operatorname*{State}%
\left(  I_{gen}\right)  \right)  \subset M_{\mathbb{R}}%
\]
If we fix the linearizations such that $H\left(  I_{0}\right)  $ corresponds
to $0\in M_{\mathbb{R}}$, then by Theorem \ref{Thm stability} this condition
is equivalent to $H\left(  I_{gen}\right)  \in\mathbb{H}^{s}\left(  E\right)
$ with $E=p^{\ast}\left(  \mathcal{O}_{\mathbb{P}\left(  W\right)  }\left(
1\right)  \right)  $, i.e., that the Hilbert point of $I_{gen}$ is in the
stable locus of the Hilbert scheme.

The construction of the Hilbert scheme in Section
\ref{Sec Computing the Bergman fan} assumes $Y$ to be a smooth toric variety.
With an appropriate definition of $\operatorname*{State}\left(  I_{gen}%
\right)  $ as discussed in Remark \ref{Rem generalization state polytope},
this condition is expected to be stated in the same form for a general
simplicial or even non simplicial toric variety $Y$.

If the Gr\"{o}bner cone $C_{I_{0}}\left(  I\right)  $ corresponding to $I_{0}$
intersects $BF\left(  I\right)  \cap\left\{  w_{t}=0\right\}  $, the Hilbert
point $H\left(  I_{0}\right)  $ of $I_{0}$ would lie on the boundary of
$\operatorname*{State}\left(  I_{gen}\right)  \subset M_{\mathbb{R}}$
contradicting $H\left(  I_{0}\right)  \in\operatorname*{int}\left(
\operatorname*{State}\left(  I_{gen}\right)  \right)  $, hence%
\[
C_{I_{0}}\left(  I\right)  \cap\left\{  w_{t}=0\right\}  =\left\{  0\right\}
\]

\item We can satisfy condition $2$. via a genericity condition on the tangent
vector with respect to the tangent space of the component of the Hilbert
scheme containing $\mathfrak{X}$:

Assume that $\mathfrak{X}$ lies in a smooth component of the complex moduli
space $\mathcal{M}$ of $X_{0}$ (for example normal crossing at $X_{0}$).

Let $v_{1},...,v_{p}\in\operatorname*{Hom}\left(  I_{0},S/I_{0}\right)  _{0}$
be a basis of the tangent space of that component of the Hilbert scheme at
$X_{0}$, which contains the tangent vector $v$ of $\mathfrak{X}$.

Assume that $\mathfrak{X}$ is maximal in its component of the Hilbert scheme,
i.e., writing $v=\sum_{i=1}^{p}\lambda_{i}v_{i}$ we have $\lambda_{i}\neq0$
$\forall i$. Consider the reduced standard basis of $I$ with respect to a
monomial ordering in the interior of $C_{I_{0}}\left(  I\right)  $. Then
already the first order deformations appearing in this standard basis, i.e.,
the Cox Laurent monomials corresponding to $t$-linear non special fiber terms,
give the linear half-space equations defining $C_{I_{0}}\left(  I\right)  $.

\item We can satisfy condition $3.$ via a condition of the resolution of
$\mathcal{O}_{X_{0}}$:

Assume that $\mathcal{O}_{X_{0}}$ has a resolution%
\[
0\rightarrow\mathcal{O}_{Y}\left(  -K_{Y}\right)  \rightarrow...\rightarrow
\mathcal{F}_{1}\rightarrow\mathcal{O}_{Y}\rightarrow\mathcal{O}_{X_{0}%
}\rightarrow0
\]
with direct sums $\mathcal{F}_{j}=%
{\textstyle\bigoplus\nolimits_{i}}
\mathcal{O}_{Y}\left(  D_{ji}\right)  $ with divisors $D_{ji}$.

Consider the reduced standard basis of $I$ in $S\times\mathbb{C}\left[
t\right]  /\left\langle t^{2}\right\rangle $ with respect to a monomial
ordering in the interior of $C_{I_{0}}\left(  I\right)  $. Then any first
order $\left(  \mathbb{C}^{\ast}\right)  ^{\Sigma\left(  1\right)  }%
$-deformation $\delta$ appearing in the standard basis is represented by a Cox
Laurent monomial such that the denominator divides $%
{\textstyle\prod\nolimits_{r\in\Sigma\left(  1\right)  }}
y_{r}$. Hence $\delta$ is also a deformation of the anticanonical Calabi-Yau
hypersurface in $Y$ defined by $\left\langle
{\textstyle\prod\nolimits_{r\in\Sigma\left(  1\right)  }}
y_{r}\right\rangle $.

\item Denote by $c$ the codimension of $X_{0}\subset Y$. We interpret
condition $4.$ as the condition that the $\Sigma$-saturated ideals defining
the components (i.e., strata of maximal dimension $d$) of $X_{0}$ are
generated by $c$ variables of $S$, i.e., are of the form $\left\langle
y_{r_{1}},...,y_{r_{c}}\right\rangle \subset S$.

\item We can satisfy condition $5.$ via a condition on the locally relevant
deformations of $\mathfrak{X}$ at the zero dimensional strata of $X_{0}$:

Let $p$ be a zero dimensional stratum of $X_{0}$ and $\mathfrak{X}_{p}$ the
flat family given by
\[
I_{p}=\left\langle m+t\cdot c_{j}\cdot\delta\left(  m\right)  \mid\delta\text{
strongly locally relevant at }p\text{, }m\in I_{0}\right\rangle
\]
with general coefficients $c_{j}$. For all zero dimensional strata $p$ of
$X_{0}$ we require: All initial ideals $\operatorname*{in}_{w}I_{i}$ for $w\in
C_{I_{0}}\left(  I\right)  $ which do not contain a monomial and are minimal
with respect to the set of contributing deformations involve precisely $c$
first order deformations.
\end{enumerate}

\begin{remark}
Note that the condition $H\left(  I_{0}\right)  \in\operatorname*{int}\left(
\operatorname*{State}\left(  I\right)  \right)  $ is independent of
rescalation of $\operatorname*{State}\left(  I\right)  $ by changing $D$, and
independent of translation of $\operatorname*{State}\left(  I\right)  $ by
changing the linearizations.

For hypersurfaces the requirement $H\left(  I_{0}\right)  \in
\operatorname*{int}\left(  \operatorname*{State}\left(  I\right)  \right)  $
is equivalent to the condition that the special fiber of $\mathfrak{X}$
corresponds to the unique interior lattice point of the Batyrev polytope
$\Delta=P^{\ast}$. We may fix a linearization of the torus action on
$\mathcal{O}_{Y}\left(  -K_{Y}\right)  $ such that $0\in M$ corresponds to the
unique interior lattice point, i.e., we fix the element
\[
V\left(
{\textstyle\prod\nolimits_{r\in\Sigma\left(  1\right)  }}
y_{r}\right)  \in\left\vert -K_{Y}\right\vert
\]
of the linear system $\left\vert -K_{Y}\right\vert $.

The condition that all facets of $\operatorname*{Strata}\nolimits_{\Delta
}\left(  I_{0}\right)  $ are contained in precisely $c$ facets of
$\Delta=P^{\ast}$ says, by flatness of the family $\mathfrak{X}$, that the
total space of $\mathfrak{X}$ is a local complete intersection at the generic
points of the strata of maximal dimension $d$ of $X_{0}$. So if $\left\langle
y_{r_{1}},...,y_{r_{c}}\right\rangle $ is a stratum of maximal dimension, then
$I$ is given by $c$ equations in the localization $S_{\left\langle y_{r_{1}%
},...,y_{r_{c}}\right\rangle }\otimes\mathbb{C}\left[  t\right]  $ at the
prime ideal $\left\langle y_{r_{1}},...,y_{r_{c}}\right\rangle $. Note that%
\[
S_{\left\langle y_{r_{1}},...,y_{r_{c}}\right\rangle }=\mathbb{C}\left(
y_{r}\mid r\notin\left\{  r_{1},...,r_{c}\right\}  \right)  \left[  y_{r_{1}%
},...,y_{r_{c}}\right]  _{>}%
\]
for any local ordering $>$ on the monomials in the variables $y_{r_{1}%
},...,y_{r_{c}}$.

The condition on the locally relevant deformations at the $0$-dimensional
strata of $X_{0}$, is a condition on the singularities of $\mathfrak{X}$ at
these strata. But note that this condition is far away from requiring the
total space of $\mathfrak{X}$ to be a local complete intersection there.
\end{remark}

\subsection{The Gr\"{o}bner cone associated to the
\index{special fiber}%
special fiber\label{Sec the groebner cone associated to the special fiber}}

Consider the setup given in Section \ref{genericy condition}. Let
$m_{1},...,m_{r}$ be minimal generators of the monomial ideal $I_{0}\subset
S$, let the
\index{flat family}%
flat family of Calabi-Yau varieties $\mathfrak{X}\subset Y\times
\operatorname*{Spec}\left(  \mathbb{C}\left[  \left[  t\right]  \right]
\right)  $ be given by the ideal
\[
I=\left\langle f_{j}=m_{j}+tg_{j}\mid j=1,...,r\right\rangle \subset
\mathbb{C}\left[  t\right]  \otimes S
\]
and suppose that the $f_{j}$ are reduced with respect to $I_{0}$, i.e., no
term of $g_{j}$ is in $I_{0}$ considered as an ideal in $\mathbb{C}\left[
t\right]  \otimes S$.

Fix a tie break ordering $>$ on $\mathbb{C}\left[  t\right]  \otimes S$, which
is respecting the Chow
\index{grading}%
grading on $S$ and is local in $t$, so $L_{>}\left(  f_{j}\right)  =m_{j}$.

\begin{definition}
Let $C_{I_{0}}\left(  I\right)  $ be the
\index{special fiber Gr\"{o}bner cone|textbf}%
cone of
\index{weight vector}%
weight vectors
\index{Gr\"{o}bner cone}%
selecting $I_{0}$ as
\index{lead ideal}%
lead ideal%
\[
C_{I_{0}}\left(  I\right)  =\left\{  -\left(  w_{t},w_{y}\right)
\in\mathbb{R}\oplus N_{\mathbb{R}}\mid L_{>_{\left(  w_{t},\varphi\left(
w_{y}\right)  \right)  }}\left(  I\right)  =I_{0}\right\}
\]

\end{definition}

Consider $w=-\left(  w_{t},w_{y}\right)  \in C_{I_{0}}\left(  I\right)  $ and
the
\index{weight ordering}%
weight ordering $>_{\left(  w_{t},\varphi\left(  w_{y}\right)  \right)  }$ on
$\mathbb{C}\left[  t\right]  \otimes S$ with tie break ordering $>$, so
\[
L_{>_{\left(  w_{t},\varphi\left(  w_{y}\right)  \right)  }}\left(
f_{j}\right)  =m_{j}%
\]
As $\mathfrak{X}$ is flat, for every
\index{syzygy}%
syzygy $s\in S^{r}$ of $m_{1},...,m_{r}$, i.e., with%
\[
\left(  m_{1},...,m_{r}\right)  \cdot s=0
\]
there is an $l\in\left(  \mathbb{C}\left[  \left[  t\right]  \right]  \otimes
S\right)  ^{c}$ such that%
\[
\left(  f_{1},...,f_{r}\right)  \cdot\left(  s-t\cdot l\right)  =0
\]
so%
\[
\frac{1}{t}\left(  f_{1},...,f_{r}\right)  \cdot s=\left(  g_{1}%
,...,g_{r}\right)  \cdot s=\left(  f_{1},...,f_{r}\right)  \cdot l
\]
i.e.,%
\[
\frac{1}{t}\left(  f_{1},...,f_{r}\right)  \cdot s\in\left\langle
f_{1},...,f_{r}\right\rangle
\]
and the Buchberger normal form in $\mathbb{C}\left[  \left[  t\right]
\right]  \otimes S$ yields%
\[
NF_{>_{\left(  w_{t},\varphi\left(  w_{y}\right)  \right)  }}\left(  \left(
f_{1},...,f_{r}\right)  \cdot s,I\right)  =0
\]
so $f_{1},...,f_{r}$ form a minimal
\index{Gr\"{o}bner basis}%
Gr\"{o}bner basis of $I$ with respect to $>_{\left(  w_{t},\varphi\left(
w_{y}\right)  \right)  }$.

As we have $f_{1},...,f_{r}$ assumed to be
\index{reduced standard basis}%
reduced, they form the reduced
\index{Gr\"{o}bner basis}%
Gr\"{o}bner basis of $I$ with respect to $>_{\left(  w_{t},\varphi\left(
w_{y}\right)  \right)  }$ and hence the condition
\[
L_{>_{\left(  w_{t},\varphi\left(  w_{y}\right)  \right)  }}\left(  I\right)
=I_{0}%
\]
is equivalent to%
\[
\operatorname*{trop}\left(  f_{j}-m_{j}\right)  \left(  w_{t},\varphi\left(
w_{y}\right)  \right)  \leq\operatorname*{trop}\left(  m_{j}\right)  \left(
\varphi\left(  w_{y}\right)  \right)  \text{ }\forall j
\]
Here $\operatorname*{trop}\left(  f_{j}-m_{j}\right)  $ denotes the
corresponding piecewise linear function of $f_{j}-m_{j}\in\mathbb{C}\left[
t\right]  \otimes S$ and $\operatorname*{trop}\left(  m_{j}\right)  $ the
piecewise linear function of $m_{j}\in S$.

\begin{lemma}
\label{1specialfibertropicalequations}With the notation from
\index{Gr\"{o}bner cone}%
above%
\[
C_{I_{0}}\left(  I\right)  =\left\{  -\left(  w_{t},w_{y}\right)
\in\mathbb{R}\oplus N_{\mathbb{R}}\mid\operatorname*{trop}\left(  f_{j}%
-m_{j}\right)  \left(  w_{t},\varphi\left(  w_{y}\right)  \right)
\leq\operatorname*{trop}\left(  m_{j}\right)  \left(  \varphi\left(
w_{y}\right)  \right)  \text{ }\forall j\right\}
\]
It is a closed polyhedral cone with $\left(  1,0,...,0\right)  \in C_{I_{0}%
}\left(  I\right)  $.
\end{lemma}

The defining equations of $C_{I_{0}}\left(  I\right)  $ are the deformations
appearing in $I$, represented as degree $0$ Cox Laurent monomials, which again
correspond to lattice monomials in $M$, i.e., if $t^{a}m\neq m_{j}$ is a
monomial of $f_{j}$, then%
\[
\left\langle w,A^{-1}\left(  \frac{m}{m_{j}}\right)  \right\rangle \geq-a\cdot
w_{t}%
\]
is a defining equation of $C_{I_{0}}\left(  I\right)  $, hence:

\begin{lemma}
The dual cone of $C_{I_{0}}\left(  I\right)  $ is spanned by the deformations
appearing in $I$, considered as degree $0$ Cox Laurent monomials, i.e.,%
\[
C_{I_{0}}\left(  I\right)  ^{\ast}=\operatorname*{hull}\left(  \left\{
\left(  \tilde{m}_{t},\tilde{m}\right)  \in\mathbb{R}\oplus M_{\mathbb{R}}\mid%
\begin{tabular}
[c]{l}%
$\exists j\text{ such that }t^{\tilde{m}_{t}}\cdot A\left(  \tilde{m}\right)
\cdot m_{j}\in\mathbb{C}\left[  t\right]  \otimes S$\\
$\text{and is a monomial of }f_{j}-m_{j}$%
\end{tabular}
\ \right\}  \right)
\]
and $\left(  1,0,...,0\right)  \in C_{I_{0}}\left(  I\right)  ^{\ast}$.
\end{lemma}

By assumption
\[
C_{I_{0}}\left(  I\right)  \cap\left\{  w_{t}=0\right\}  =\left\{  0\right\}
\]
hence:

\begin{lemma}
The cone $C_{I_{0}}\left(  I\right)  $ minus the zero point is contained in
the half-space $\left\{  w_{t}>0\right\}  $.
\end{lemma}

If $\left(  1,0,...,0\right)  \in C_{I_{0}}\left(  I\right)  ^{\ast}$ would
lie on the boundary of $C_{I_{0}}\left(  I\right)  ^{\ast}$, then $C_{I_{0}%
}\left(  I\right)  $ would contain a ray in $\left\{  w_{t}=0\right\}  $, hence:

\begin{lemma}
The monomial $\left(  1,0,...,0\right)  $ lies in the interior of the dual
cone $C_{I_{0}}\left(  I\right)  ^{\ast}$, i.e.,%
\[
\left(  1,0,...,0\right)  \in\operatorname*{int}\left(  C_{I_{0}}\left(
I\right)  ^{\ast}\right)
\]

\end{lemma}

The flat family $\mathfrak{X}\subset Y\times\operatorname*{Spec}\left(
\mathbb{C}\left[  \left[  t\right]  \right]  \right)  $ induces a first order
flat family
\[
\mathfrak{X}^{1}\subset Y\times\operatorname*{Spec}\left(  \mathbb{C}\left[
t\right]  /\left\langle t^{2}\right\rangle \right)
\]
given by
\[
I^{1}=\left\langle f_{j}^{1}=m_{j}+tg_{j}^{1}\mid j=1,...,r\right\rangle
=\mathbb{C}\left[  t\right]  /\left\langle t^{2}\right\rangle \otimes S
\]
with $g_{j}^{1}\in S$.

By above assumption the defining equations of $C_{I_{0}}\left(  I\right)  $
are given by first order deformations appearing in the reduced standard basis
of $I$, so%
\[
C_{I_{0}}\left(  I\right)  =\left\{  -\left(  w_{t},w_{y}\right)
\in\mathbb{R}\oplus N_{\mathbb{R}}\mid\operatorname*{trop}\left(  g_{j}%
^{1}\right)  \left(  w_{t},\varphi\left(  w_{y}\right)  \right)  +w_{t}%
\leq\operatorname*{trop}\left(  m_{j}\right)  \left(  \varphi\left(
w_{y}\right)  \right)  \text{ }\forall j\right\}
\]

\begin{corollary}
Intersecting $C_{I_{0}}\left(  I\right)  $
\index{Gr\"{o}bner cone}%
with the hyperplane $\left\{  w_{t}=1\right\}  $ we obtain the convex
polytope
\[
\nabla=C_{I_{0}}\left(  I\right)  \cap\left\{  w_{t}=1\right\}  \subset
N_{\mathbb{R}}%
\]
with%
\[
\nabla=\left\{  -w_{y}\in N_{\mathbb{R}}\mid\operatorname*{trop}\left(
g_{j}\right)  \left(  \varphi\left(  w_{y}\right)  \right)  -1\leq
\operatorname*{trop}\left(  m_{j}\right)  \left(  \varphi\left(  w_{y}\right)
\right)  \text{ }\forall j\right\}
\]
and $0$ in the interior of $\nabla$.
\end{corollary}

Rewriting the tropical equations%

\begin{align*}
\nabla &  =\left\{  -w_{y}\in N_{\mathbb{R}}\mid\operatorname*{trop}\left(
m\right)  \left(  \varphi\left(  w_{y}\right)  \right)  -1\leq
\operatorname*{trop}\left(  m_{j}\right)  \left(  \varphi\left(  w_{y}\right)
\right)  \text{ }\forall\text{ monomials }m\text{ of }g_{j}^{1}\text{ }\forall
j\right\} \\
&  =\left\{  w_{y}\in N_{\mathbb{R}}\mid\varphi\left(  w_{y}\right)  \left(
\frac{m}{m_{j}}\right)  \geq-1\text{ }\forall\text{ monomials }m\text{ of
}g_{j}^{1}\text{ and }\forall j\right\} \\
&  =\left\{  w_{y}\in N_{\mathbb{R}}\mid\left\langle A^{-1}\left(  \frac
{m}{m_{j}}\right)  ,w_{y}\right\rangle \geq-1\text{ }\forall\text{ monomials
}m\text{ of }g_{j}^{1}\text{ and }\forall j\right\}
\end{align*}
hence%
\[
\nabla^{\ast}=\operatorname*{convexhull}\left\{  A^{-1}\left(  \frac{m}{m_{j}%
}\right)  \in M_{\mathbb{R}}\mid\exists j\text{ such that }m\text{ is a
monomial of }g_{j}^{1}\right\}
\]
so it follows:

\begin{lemma}
$\nabla^{\ast}$
\index{lattice polytope}%
is an
\index{integral polytope}%
integral polytope.
\end{lemma}

Any first order deformation appearing in $g_{j}^{1}$ represented by a Cox
Laurent monomial $\frac{m}{m_{j}}$ is also a deformation of the anticanonical
Calabi-Yau hypersurface in $Y$ defined by $\left\langle
{\textstyle\prod\nolimits_{r\in\Sigma\left(  1\right)  }}
y_{r}\right\rangle $, hence $A^{-1}\left(  \frac{m}{m_{i}}\right)  \in
\Delta=\Delta_{-K_{Y}}$, i.e., $\nabla^{\ast}\subset\Delta$. As $\Delta$ is
dual to a Fano polytope, it has $0$ as unique interior lattice point by Lemma
\ref{lemma Fano polytope 0 unique interior lattice point}. Hence also
$\nabla^{\ast}$ has no interior lattice point besides $0$. As $\left(
1,0,...,0\right)  \in\operatorname*{int}\left(  C_{I_{0}}\left(  I\right)
^{\ast}\right)  $, the polytope $\nabla^{\ast}$ contains $0$ in its interior.

\begin{lemma}
$\nabla^{\ast}$ contains $0$ as unique interior lattice point.
\end{lemma}

\begin{theorem}
\label{thm fano polytope}$\nabla^{\ast}$ is a Fano polytope, hence the fan
$\Sigma^{\circ}=\Sigma\left(  \nabla^{\ast}\right)  $ over the faces of
$\nabla^{\ast}$ defines a
\index{Q-Gorenstein}%
$\mathbb{Q}$-Gorenstein
\index{toric Fano}%
toric
\index{Fano}%
Fano variety $Y^{\circ}=X\left(  \Sigma^{\circ}\right)  $.
\end{theorem}

\subsection{The dual complex of initial
ideals\label{Sec dual complex general setting}}

\begin{definition}
If $F$ is a face of $\nabla$, there is
\index{initial form}%
an associated
\index{initial ideal of face|textbf}%
\textbf{initial ideal of }$I$\textbf{ with respect to the face} $F$: For all
$w_{1},w_{2}$ in the relative interior $\operatorname*{int}\left(  F\right)  $
of $F$ we have%
\[
\operatorname*{in}\nolimits_{_{\left(  1,\varphi\left(  w_{1}\right)  \right)
}}\left(  I\right)  =\operatorname*{in}\nolimits_{_{\left(  1,\varphi\left(
w_{2}\right)  \right)  }}\left(  I\right)
\]
Denote this ideal by $\operatorname*{in}\nolimits_{F}\left(  I\right)  $. For
all $w_{1},w_{2}\in\operatorname*{int}\left(  F\right)  $ and $f\in I$%
\[
\operatorname*{in}\nolimits_{_{\left(  1,\varphi\left(  w_{1}\right)  \right)
}}\left(  f\right)  =\operatorname*{in}\nolimits_{_{\left(  1,\varphi\left(
w_{2}\right)  \right)  }}\left(  f\right)
\]
denote this initial term of $f$ by $\operatorname*{in}\nolimits_{F}\left(
f\right)  $.
\end{definition}

If $F$ is a face of $\nabla$, then
\[
\operatorname*{in}\nolimits_{F}\left(  I\right)  =\left\langle
\operatorname*{in}\nolimits_{F}\left(  f_{j}\right)  \mid
j=1,...,r\right\rangle
\]
as $f_{1},...,f_{r}$ form a Gr\"{o}bner basis of $I$ with respect to any
weight vector in $C_{I_{0}}\left(  I\right)  $.

Recall that we wrote
\[
I^{1}=\left\langle f_{j}^{1}=m_{j}+tg_{j}^{1}\mid j=1,...,r\right\rangle
\]
with $g_{j}^{1}\in S$ for the ideal of the first order deformation
$\mathfrak{X}^{1}$ associated to $\mathfrak{X}$.

For $j=1,...,r$ define $G_{j}\left(  F\right)  $ as%
\[
\operatorname*{in}\nolimits_{F}\left(  f_{j}^{1}\right)  =t\sum_{m\in
G_{j}\left(  F\right)  }c_{m}\cdot m+m_{j}%
\]

\begin{definition}
If $F$ is a face of $\nabla$, then define the
\index{dual|textbf}%
\textbf{dual face} of $F$ as
\[
\operatorname*{dual}\left(  F\right)  =\operatorname*{convexhull}\left(
A^{-1}\left(  \frac{m}{m_{j}}\right)  \mid m\in G_{j}\left(  F\right)  ,\text{
}j=1,...,r\right)  \subset M_{\mathbb{R}}%
\]
the convex hull of the first order deformations appearing in the initial ideal
with respect to $F$. The dual face is a lattice polytope in $M_{\mathbb{R}}$.
\end{definition}

By the genericity condition on the tangent vector of $\mathfrak{X}$ we have:

\begin{lemma}%
\[
\operatorname*{dual}\left(  F\right)  \cap M=\left\{  A^{-1}\left(  \frac
{m}{m_{j}}\right)  \mid m\in G_{j}\left(  F\right)  ,\text{ }%
j=1,...,r\right\}
\]

\end{lemma}

The dual of $\nabla$ is the convex hull of the first order deformations
appearing in $I$%
\[
\nabla^{\ast}=\operatorname*{convexhull}\left(  \left\{  A^{-1}\left(
\frac{m}{m_{j}}\right)  \in M_{\mathbb{R}}\mid\exists j\text{ such that
}m\text{ is a monomial of }g_{j}^{1}\right\}  \right)
\]
so for the dual face of $F$ we have%
\begin{align*}
\operatorname*{dual}\left(  F\right)   &  =\operatorname*{convexhull}\left(
\left\{  A^{-1}\left(  \frac{m}{m_{j}}\right)  \mid m\in G_{j}\left(
F\right)  ,\text{ }j=1,...,r\right\}  \right) \\
&  =\operatorname*{convexhull}\left(  \bigcup_{j=1}^{r}\left\{  A^{-1}\left(
\frac{m}{m_{j}}\right)  \mid m\in G_{j}\left(  F\right)  \right\}  \right) \\
&  =\operatorname*{convexhull}\left(  \bigcup_{j=1}^{r}\left\{  A^{-1}\left(
\frac{m}{m_{j}}\right)  \mid%
\begin{tabular}
[c]{l}%
$m\in g_{j}^{1}\text{ with }$\\
$\left\langle A^{-1}\frac{m}{m_{j}},w_{y}\right\rangle =-1\text{ }\forall
w_{y}\in F$%
\end{tabular}
\ \ \right\}  \right) \\
&  =\operatorname*{convexhull}\left(  \left\{  \widetilde{m}\in\nabla^{\ast
}\cap M\mid\left\langle \widetilde{m},w_{y}\right\rangle =-1\text{ }\forall
w_{y}\in F\right\}  \right) \\
&  =\left\{  \widetilde{m}\in\nabla^{\ast}\mid\left\langle \widetilde{m}%
,w_{y}\right\rangle =-1\text{ }\forall w_{y}\in F\right\}
\end{align*}
hence:

\begin{proposition}
\label{Prop dual face combinatorial dualization general setting}If $F$ is a
face of $\nabla$, then
\[
\operatorname*{dual}\left(  F\right)  =F^{\ast}%
\]
in particular $\operatorname*{dual}\left(  F\right)  $ is
\newsym[$\operatorname*{dual}\left(  F\right)  $]{dual face}{}a face of
$\nabla^{\ast}$,
\[%
\begin{tabular}
[c]{cccc}%
$\operatorname*{dual}:$ & $\operatorname*{Poset}\left(  \nabla\right)  $ &
$\rightarrow$ & $\operatorname*{Poset}\left(  \nabla^{\ast}\right)  $\\
& $F$ & $\mapsto$ & $\operatorname*{dual}\left(  F\right)  $%
\end{tabular}
\ \ \
\]
is an inclusion reversing map of complexes and
\[
\dim\left(  \operatorname*{dual}\left(  F\right)  \right)  =n-1-\dim\left(
F\right)
\]

\end{proposition}

The non-special fiber terms of $\operatorname*{in}\nolimits_{F}\left(
f_{j}^{1}\right)  $, $j=1,...,r$, i.e., the elements of $G_{j}\left(
F\right)  $, $j=1,...,r$ split into characters of the big torus $\left(
\mathbb{C}^{\ast}\right)  ^{\Sigma\left(  1\right)  }$. These characters are
the Cox Laurent monomials%
\[
\delta_{F}\left(  I^{1}\right)  =\left\{  \frac{m}{m_{j}}\mid m\in
G_{j}\left(  F\right)  ,\text{ }j=1,...,r\right\}
\]
and represent the first order $\left(  \mathbb{C}^{\ast}\right)
^{\Sigma\left(  1\right)  }$-deformations contributing to the degeneration
defined by $\operatorname*{in}\nolimits_{F}\left(  I^{1}\right)  $. Flat
families defined by initial ideals are also called
\index{Gr\"{o}bner deformation}%
Gr\"{o}bner deformations. Note that the syzygies of a monomial ideal are
binomial and a syzygy between $m_{i}$ and $m_{j}$ is represented by the
character of $\left(  \mathbb{C}^{\ast}\right)  ^{\Sigma\left(  1\right)  }$
given by $\operatorname*{lcm}\left(  m_{i},m_{j}\right)  $. Note also, that
monomials $m\in G_{i}\left(  F\right)  $ and $m^{\prime}\in G_{j}\left(
F\right)  $ with $\frac{m}{m_{i}}=\frac{m^{\prime}}{m_{j}}$ appear with the
same coefficient in the initial forms. On the other hand the elements of
$\delta_{F}\left(  I^{1}\right)  $ correspond via $A^{-1}$ to the lattice
points of $F^{\ast}$, so:

\begin{lemma}
The lattice points of $F^{\ast}$ are in one-to-one correspondence to the first
order deformations contributing to $\operatorname*{in}\nolimits_{F}\left(
I^{1}\right)  $. If $\delta\in F^{\ast}$ is a lattice point and $\frac{q_{1}%
}{q_{0}}=A\left(  \delta\right)  $ with relatively prime monomials
$q_{0},q_{1}\in S$ then%
\[
\delta\left(  m\right)  =\left\{
\begin{tabular}
[c]{ll}%
$\frac{q_{1}}{q_{0}}\cdot m$ & if $q_{0}\mid m$\\
$0$ & otherwise
\end{tabular}
\ \right\}
\]
for minimal generators $m\in I_{0}$ defines the corresponding $\left(
\mathbb{C}^{\ast}\right)  ^{\Sigma\left(  1\right)  }$-deformation in
$\operatorname*{Hom}\left(  I_{0},S/I_{0}\right)  _{0}$.
\end{lemma}

\subsection{Bergman subcomplex of $\nabla$%
\label{Sec Bergman subcomplex of Nabla general setting}}

\begin{definition}
The
\index{special fiber Bergman fan}%
\textbf{special fiber Bergman fan} is
\newsym[$\operatorname*{dual}\left(  F\right)  $]{dual face}{}defined as the
intersection of the fan $\operatorname*{Poset}\left(  C_{I_{0}}\left(
I\right)  \right)  $ over the special fiber Gr\"{o}bner cone $C_{I_{0}}\left(
I\right)  $ with the
\newsym[$BF_{I_{0}}\left( I\right)$]{special fiber Bergman fan}{}Bergman fan
$BF\left(  I\right)  $ of $I$%
\[
BF_{I_{0}}\left(  I\right)  =BF\left(  I\right)  \cap\operatorname*{Poset}%
\left(  C_{I_{0}}\left(  I\right)  \right)
\]
The
\index{special fiber Bergman complex}%
\textbf{special fiber Bergman complex}%
\[
B\left(  I\right)  =BC_{I_{0}}\left(  I\right)  =\left(  BF\left(  I\right)
\cap\operatorname*{Poset}\left(  C_{I_{0}}\left(  I\right)  \right)  \right)
\cap\left\{  w_{t}=1\right\}  \subset\operatorname*{Poset}\left(
\nabla\right)
\]
is defined as
\newsym[$BC_{I_{0}}\left(  I\right)$]{special fiber Bergman complex}{}the
\newsym[$B\left(  I\right)$]{special fiber Bergman complex}{}complex whose
faces are the intersections of the hyperplane $\left\{  w_{t}=1\right\}  $
with the faces of the
\index{Bergman complex}%
Bergman fan $BF\left(  I\right)  $ in $C_{I_{0}}\left(  I\right)  $.
\end{definition}

We also refer to $B\left(  I\right)  $ as the
\index{Bergman subcomplex|textbf}%
Bergman subcomplex or
\index{tropical subcomplex|textbf}%
tropical subcomplex of $\nabla$. By Theorem
\ref{thm tropical variety properties} we have:

\begin{remark}
The complex $B\left(  I\right)  $ consists of those faces $F$ of $\nabla$ such
that $\operatorname*{in}_{F}\left(  I\right)  $ does not contain a monomial.
\end{remark}

\begin{lemma}
The special fiber Bergman complex $B\left(  I\right)  $ is a polyhedral
\index{cell complex}%
cell complex, it is subcomplex of the boundary $\partial\nabla$ of $\nabla$.
\end{lemma}

\subsection{Remarks on the covering structure in $\operatorname*{dual}\left(
B\left(  I\right)  \right)  $\label{Sec covering structure in dual B}}

Interpreting the lattice points of the faces of $\operatorname*{dual}\left(
B\left(  I\right)  \right)  $ as deformations of $X_{0}$ and associating them
to the reduced standard basis equations $f_{i}$, $i=1,...r$ defining the total
space, we get:

If $F\in B\left(  I\right)  $ is a face, then denote by $G_{F}$ a minimal
standard basis of $I$ in%
\[
S_{I_{_{0}}\left(  F\right)  }\otimes\mathbb{C}\left[  t\right]  /\left\langle
t^{2}\right\rangle
\]
with the localization
\[
S_{I_{_{0}}\left(  F\right)  }=\mathbb{C}\left(  y_{j}\mid j\notin J\right)
\left[  y_{j}\mid j\in J\right]  _{>}%
\]
where the prime ideal $I_{_{0}}\left(  F\right)  =\left\langle y_{j}\mid j\in
J\right\rangle \subset S$ denotes the face of $SP\left(  I_{0}\right)  $
corresponding to $F$ and $>$ is a local ordering on the $y_{j}$. The standard
basis is computed using Mora normal form. Let $s$ be the maximum number of
elements of the $G_{F}$ over all faces $F\in B\left(  I\right)  $. Denote by
$\tilde{G}_{F}$ the standard basis reduced via Gr\"{o}bner normal form.

\begin{lemma}
If $F\in B\left(  I\right)  $ is a face, the lattice points of
$\operatorname*{dual}\left(  F\right)  $ are the first order deformations
appearing in the initial ideal of $\tilde{G}_{F}$ with respect to $F$.

The complex $\operatorname*{dual}\left(  B\left(  I\right)  \right)  $
contains an $s:1$ covering of faces: If $G$ is a face over $F^{\vee}\in
B\left(  I\right)  ^{\vee}$, then the lattice points of $G$ are the
deformations appearing in the initial form of one of the equations of the
reduced local standard basis $\tilde{G}_{F}$ of $I$ considered as an ideal in
$S_{I_{_{0}}\left(  F\right)  }\otimes\mathbb{C}\left[  t\right]
/\left\langle t^{2}\right\rangle $.

In general this covering is branched and the number of faces over $F^{\vee}\in
B\left(  I\right)  ^{\vee}$ is the number of elements of the reduced local
standard basis $\tilde{G}_{F}$ of $I$.
\end{lemma}

Note that this covering can have degenerate faces, i.e., faces $G$ over
$F^{\vee}\in B\left(  I\right)  ^{\vee}$ with $\dim\left(  G\right)
<\dim\left(  F^{\vee}\right)  $. It can be branched in the sense that if
$\mathfrak{X}$ is not a local complete intersection, the number of faces $G$
over a face of $B\left(  I\right)  ^{\vee}$ may be larger than the
codimension. Note that this number is bounded from below by the codimension.

If two first order deformations $\delta_{1}$ and $\delta_{2}$ lie in the same
face of the covering, then there is an element $f_{j}=m_{j}+tg_{j}$ of the
global reduced standard basis such that both $\delta_{1}$ and $\delta_{2}$
contribute in $f_{j}$, i.e., $g_{j}$ involves the monomials $\delta_{1}\left(
m_{j}\right)  $ and $\delta_{2}\left(  m_{j}\right)  $. If two deformations
contribute in the same element of the reduced global standard basis, they are
connected by a chain of faces of the covering.

The set of faces over $B\left(  I\right)  ^{\vee}$ can be totally
disconnected, e.g. if every element of the global reduced Gr\"{o}bner basis
involves at most one of the first order deformations, then all fibers of the
covering consist of points.

Removing all faces of the covering, which correspond to locally irrelevant
equations, removing multiple faces, which correspond to locally equivalent
equations, and keeping only faces, which involve only vertices of faces of the
covering of smaller dimension, we obtain a covering $\pi$ of $B\left(
I\right)  ^{\vee}$ denoted as the \textbf{reduced covering}.

\begin{remark}
If $I$ is a complete intersection the reduced covering $\pi$ is the $c:1$
covering given in Section \ref{Sec covering complete intersection}. If $I$ is
a local complete intersection then $\pi$ is also $c:1$.
\end{remark}

\begin{algorithm}
The following algorithm computes the reduced covering $\pi$:

\begin{itemize}
\item If $F$ is a face of $B\left(  I\right)  $ of $\dim\left(  F\right)  =d$
and $p_{1},...,p_{c}$ are the vertices of $\operatorname*{dual}\left(
F\right)  $ then set
\[
\pi\left(  p_{j}\right)  =F^{\vee}%
\]
for $j=1,...,c$.

\item If $l>0$ and $F$ is a face of $B\left(  I\right)  $ of $\dim\left(
F\right)  =d-l$ then the faces of the covering $\pi$ over $F^{\vee}$ are the
convex hulls $H$ of those subsets of the set of vertices of
$\operatorname*{dual}\left(  F\right)  $ with

\begin{itemize}
\item $H$ involves only vertices of faces $\pi^{-1}\left(  Q^{\vee}\right)  $
with $Q^{\vee}\in B\left(  I\right)  ^{\vee}$, $Q^{\vee}\subsetneqq F^{\vee}$,
i.e., of faces of the covering lying in some lower dimensional
$\operatorname*{dual}\left(  F\right)  $ for $F\in B\left(  I\right)  $.

\item $H$ intersects at most one of the elements of $\pi^{-1}\left(  Q^{\vee
}\right)  $ for all faces $Q^{\vee}\subsetneqq F^{\vee}$ of $B\left(
I\right)  ^{\vee}$, i.e., for all faces $Q$ of $B\left(  I\right)  $ with
$F\subsetneqq Q$,

\item $H\notin\pi^{-1}\left(  Q^{\vee}\right)  $ for all faces $Q^{\vee
}\subsetneqq F^{\vee}$.
\end{itemize}
\end{itemize}
\end{algorithm}

\subsection{Limit map\label{Sec limit map general setting}}

Recall that $N\cong\mathbb{Z}^{n}$, $M=\operatorname*{Hom}\left(
N,\mathbb{Z}\right)  $ is the dual lattice of $N$, $P$ is a Fano polytope,
$\Delta=P^{\ast}$, $\Sigma$ is the fan over $P$ and $Y=X\left(  \Sigma\right)
$ is the corresponding toric Fano variety with Cox ring $S=\mathbb{C}\left[
y_{r}\mid r\in\Sigma\left(  1\right)  \right]  $ and presentation%
\[%
\begin{tabular}
[c]{lllllll}%
$0\rightarrow$ & $M$ & $\overset{A}{\rightarrow}$ & $\operatorname*{WDiv}%
_{T}\left(  X\left(  \Sigma\right)  \right)  \cong\mathbb{Z}^{\Sigma\left(
1\right)  }$ & $\overset{\deg}{\rightarrow}$ & $A_{n-1}\left(  X\left(
\Sigma\right)  \right)  $ & $\rightarrow0$%
\end{tabular}
\ \ \ \ \
\]
of
\index{Chow group}%
$A_{n-1}\left(  X\left(  \Sigma\right)  \right)  $. By Section
\ref{1homogeneouscoordinate}%
\[
G\left(  \Sigma\right)  =\operatorname*{Hom}\nolimits_{\mathbb{Z}}\left(
A_{n-1}\left(  X\left(  \Sigma\right)  \right)  ,\mathbb{C}^{\ast}\right)
\]
acts on $\mathbb{C}^{\Sigma\left(  1\right)  }$ and with the
\index{irrelevant ideal}%
irrelevant ideal%
\[
B\left(  \Sigma\right)  =\left\langle \prod_{r\not \subset \sigma}y_{r}%
\mid\sigma\in\Sigma\right\rangle \subset S
\]
in the Cox ring, $\mathbb{C}^{\Sigma\left(  1\right)  }-V\left(  B\left(
\Sigma\right)  \right)  $ is invariant under $G\left(  \Sigma\right)  $ and we
have%
\[
Y=\left(  \mathbb{C}^{\Sigma\left(  1\right)  }-V\left(  B\left(
\Sigma\right)  \right)  \right)  //G\left(  \Sigma\right)
\]

Considering the setup from Section \ref{genericy condition}, recall that
$I_{0}\subset S$ is an equidimensional reduced monomial ideal and
$\mathfrak{X}\subset Y\times\operatorname*{Spec}\mathbb{C}\left[  \left[
t\right]  \right]  $ is a flat family of Calabi-Yau varieties of dimension $d$
given by the ideal $I\subset\mathbb{C}\left[  t\right]  \otimes S$ and special
fiber $X_{0}$, defined by $I_{0}$.

As defined in Section \ref{1NonArchimedianAmoebas}, we denote by $K$ the
metric completion of the field of Puisseux series $\overline{\mathbb{C}\left(
t\right)  }$ and by
\begin{align*}
\operatorname*{val}  &  :\left(  K^{\ast}\right)  ^{n}\rightarrow
\mathbb{R}^{n}\\
\left(  f_{1},...,f_{n}\right)   &  \mapsto\left(  val\left(  f_{1}\right)
,...,val\left(  f_{n}\right)  \right)
\end{align*}
the valuation map.

Applying $\operatorname*{Hom}\nolimits_{\mathbb{Z}}\left(  -,K^{\ast}\right)
$ to%
\[%

\ \ \ \
\]
we get an exact sequence%
\[%
%
\ \ \
\]
where $D_{r}$, $r\in\Sigma\left(  1\right)  $ denote the prime $T$-Weil
divisors, gives an isomorphism%
\[%
%
\ \ \
\]
and choosing a basis $e_{1},...,e_{n}$ of $N$, we have an isomorphism%
\[%
%
\ \ \
\]
Representing $A$ by $\left(  a_{rj}\right)  _{r\in\Sigma\left(  1\right)
,\text{ }j=1,...,n}$ with respect to these bases, we have%

\[%
%
\ \ \ \ \
\]
Then $V_{K}\left(  I\right)  \subset\left(  K^{\ast}\right)  ^{n}$ is the
image of the vanishing locus of $I\subset\mathbb{C}\left[  t\right]  \otimes
S$ in $\left(  K^{\ast}\right)  ^{\Sigma\left(  1\right)  }%
/\operatorname*{Hom}\nolimits_{\mathbb{Z}}\left(  A_{n-1}\left(  X\left(
\Sigma\right)  \right)  ,K^{\ast}\right)  $ under the isomorphism induced by
$\pi$%
\[
\left(  K^{\ast}\right)  ^{\Sigma\left(  1\right)  }/\operatorname*{Hom}%
\nolimits_{\mathbb{Z}}\left(  A_{n-1}\left(  X\left(  \Sigma\right)  \right)
,K^{\ast}\right)  \cong\left(  K^{\ast}\right)  ^{n}%
\]

If $F$ is a face of the special fiber Bergman complex%
\[
F\in B\left(  I\right)  \subset\operatorname*{val}\left(  V_{K}\left(
I\right)  \right)  =-\operatorname*{tropvar}\left(  I\right)
\]
then
\[
\operatorname*{val}\nolimits^{-1}\left(  \operatorname*{int}\left(  F\right)
\right)  \subset V_{K}\left(  I\right)  \subset\left(  K^{\ast}\right)  ^{n}%
\]
is the set of arc solutions of $I$ over the weight vectors in the relative
interior of $F$. Hence if $w\in\operatorname*{int}\left(  F\right)  $ there is
an arc
\[
a\left(  t\right)  =\left(  a_{i}t^{w_{i}}+\operatorname*{hot}\right)
_{i=1,...,n}\in V_{K}\left(  I\right)  \subset\left(  K^{\ast}\right)  ^{n}%
\]
with $a_{i}\in\mathbb{C}^{\ast}$. Using multi index notation write%
\[
a\left(  t\right)  =\left(  a_{i}t^{w_{i}}+\operatorname*{hot}\right)
_{i=1,...,n}=a_{w}\cdot t^{w}+\operatorname*{hot}\in\left(  K^{\ast}\right)
^{n}%
\]
with $a_{w}\in\left(  \mathbb{C}^{\ast}\right)  ^{n}$. In the following we
show that for all arcs $a\left(  t\right)  \in\operatorname*{val}%
\nolimits^{-1}\left(  \operatorname*{int}\left(  F\right)  \right)  $ the
limit point $\lim_{t\rightarrow0}a\left(  t\right)  $ lies in the same stratum
of the fiber $Y$ of $Y\times\operatorname*{Spec}\mathbb{C}\left[  \left[
t\right]  \right]  \rightarrow\operatorname*{Spec}\mathbb{C}\left[  \left[
t\right]  \right]  $ over $\operatorname*{Spec}\mathbb{C}$. We identify the
stratum and show that it is a stratum of $X_{0}$.

First suppose $a\left(  t\right)  =a_{w}\cdot t^{w}+\operatorname*{hot}%
\in\left(  K^{\ast}\right)  ^{n}$ is any element of $\left(  K^{\ast}\right)
^{n}$, then approximating a real vector $w\in N_{\mathbb{R}}\cong
\mathbb{R}^{n}$ by a sequence rational vectors $\left(  q_{j}\right)  $ with
$q_{j}\in\operatorname*{int}\left(  F\right)  $ and $\lim_{j\rightarrow\infty
}q_{j}=w$, we may assume that $w\in\mathbb{Q}^{n}$. The limit of a power
$\lim_{t\rightarrow0}a\left(  t\right)  ^{b}$ with $b\in\mathbb{Z}_{\geq1}$ of
the arc $a\left(  t\right)  $ exists if and only if $\lim_{t\rightarrow
0}a\left(  t\right)  $ exists and lies in the same stratum of $Y$. Taking the
power of the arc multiplies $w\in\mathbb{Q}^{n}$ with $b$, hence we may assume
that $w^{\prime}=bw\in N$.

Recall from Section \ref{Sec Toric varieties from fans} that there is a
one-to-one correspondence between lattice points of $N$ and $1$-parameter
subgroups
\index{torus}%
of $T=\operatorname*{Hom}\left(  M,\mathbb{C}^{\ast}\right)  $ given by%
\[%
\begin{tabular}
[c]{lll}%
$N$ & $\rightarrow$ & $\operatorname*{Hom}\left(  \mathbb{C}^{\ast},T\right)
$\\
$w$ & $\mapsto$ & $%
\begin{tabular}
[c]{llll}%
$\lambda_{w}:$ & $\mathbb{C}^{\ast}$ & $\rightarrow$ & $\operatorname*{Hom}%
\left(  M,\mathbb{C}^{\ast}\right)  $\\
& \multicolumn{1}{c}{$t$} & $\mapsto$ &
\begin{tabular}
[c]{llll}%
$\lambda_{w}\left(  t\right)  :$ & $M$ & $\rightarrow$ & $\mathbb{C}^{\ast}$\\
& $m$ & $\mapsto$ & $t^{\left\langle m,w\right\rangle }$%
\end{tabular}
\end{tabular}
\ \ \ \ $%
\end{tabular}
\ \ \ \
\]

So if $\tau$ is a cone of $\Sigma$ with $bw\in\operatorname*{int}\left(
\tau\right)  $ in the relative interior by Proposition
\ref{Prop Toric limit points}%
\[
\lim_{t\rightarrow0}\lambda_{bw}\left(  t\right)  =x_{\tau}%
\]
where $x_{\tau}$ is the
\index{distinguished point}%
distinguished point%
\[%
\begin{tabular}
[c]{llll}%
$x_{\tau}:$ & $\check{\tau}\cap M$ & $\rightarrow$ & $\mathbb{C}$\\
& $m$ & $\mapsto$ & $\left\{
\begin{tabular}
[c]{ll}%
$1$ & if $m\in\tau^{\perp}$\\
$0$ & otherwise
\end{tabular}
\ \ \right\}  $%
\end{tabular}
\ \
\]

As $\Sigma=\operatorname*{NF}\left(  \Delta\right)  $ is complete,
$\lim_{t\rightarrow0}a\left(  t\right)  $ exists in $Y$ and lies in the unique
stratum of $Y$ containing $x_{\tau}=\lim_{t\rightarrow0}\lambda_{bw}\left(
t\right)  $.

\begin{lemma}
If $a\left(  t\right)  =a_{w}\cdot t^{w}+\operatorname*{hot}\in\left(
K^{\ast}\right)  ^{n}\cong\operatorname*{Hom}\nolimits_{\mathbb{Z}}\left(
M,K^{\ast}\right)  $, then $\lim_{t\rightarrow0}a\left(  t\right)  $ exists in
$Y$ and lies in the unique stratum of $Y$ containing $x_{\tau}$ where $\tau$
is the cone of $\Sigma$ containing $w$ in its relative interior. This stratum
is $V\left(  \tau\right)  $.
\end{lemma}

Now suppose $F$ is a face of the special fiber Bergman complex $B\left(
I\right)  \subset\operatorname*{val}\left(  V_{K}\left(  I\right)  \right)  $
and $a\left(  t\right)  \in\operatorname*{val}\nolimits^{-1}\left(
\operatorname*{int}\left(  F\right)  \right)  $, so $a\left(  t\right)
=a_{w}\cdot t^{w}+\operatorname*{hot}\in\left(  K^{\ast}\right)  ^{n}$ with
$w\in\operatorname*{int}\left(  F\right)  $.

\begin{lemma}
If $F$ is a face of the special fiber Bergman complex $B\left(  I\right)  $
then there is a unique cone $\tau$ of $\Sigma$ such that $\operatorname*{int}%
\left(  F\right)  \subset\operatorname*{int}\left(  \tau\right)  $.
\end{lemma}

\begin{definition}
Hence we can define the map
\[%
\begin{tabular}
[c]{llll}%
$\mu:$ & $B\left(  I\right)  $ & $\rightarrow$ & $\operatorname*{Poset}\left(
\Delta\right)  $\\
& \multicolumn{1}{c}{$F$} & $\mapsto$ & $G$%
\end{tabular}
\]
where $G$ is the face of $\Delta$ with $\tau=\operatorname*{hull}\left(
G^{\ast}\right)  $, where $\tau$ the unique cone of $\Sigma$ such that
$\operatorname*{int}\left(  F\right)  \subset\operatorname*{int}\left(
\tau\right)  $.
\end{definition}

\begin{lemma}
Suppose $F$ is a face of the special fiber Bergman complex $B\left(  I\right)
$ and $a\left(  t\right)  \in\operatorname*{val}\nolimits^{-1}\left(
\operatorname*{int}\left(  F\right)  \right)  $. If $X_{i}$ is the stratum of
$Y$ containing $\lim_{t\rightarrow0}a\left(  t\right)  $ in its interior, then
$X_{i}$ is a stratum of $X_{0}$.
\end{lemma}

For any point $x_{0}\in X_{0}$, by taking a hyperplane section of
$\mathfrak{X}\subset Y\times\operatorname*{Spec}\mathbb{C}\left[  \left[
t\right]  \right]  $ through $x_{0}$, there is an arc $a\left(  t\right)  \in
V_{K}\left(  I\right)  $ such that $\lim_{t\rightarrow0}a\left(  t\right)
=x_{0}$, so $\operatorname*{val}\left(  a\left(  t\right)  \right)  \in
B\left(  I\right)  $, hence:

\begin{proposition}
If $F$ is a face of $B\left(  I\right)  $, then%
\[
\lim\left(  F\right)  =\left\{  \lim_{t\rightarrow0}a\left(  t\right)  \mid
a\in\operatorname*{val}\nolimits^{-1}\left(  \operatorname*{int}\left(
F\right)  \right)  \right\}
\]
is a closed stratum of $X_{0}$,
\newsym[$\lim\left(  F\right)  $]{limit stratum}{}called the
\index{limit stratum|textbf}%
\textbf{limit stratum of }$F$.

If $\tau$ is the unique cone of $\Sigma$ such that $\operatorname*{int}\left(
F\right)  \subset\operatorname*{int}\left(  \tau\right)  $, then
\[
\lim\left(  F\right)  =V\left(  \tau\right)  =V\left(  \operatorname*{hull}%
\left(  \left(  \mu\left(  F\right)  \right)  ^{\ast}\right)  \right)
\]

Associating to a face $F$ of the special fiber Bergman subcomplex $B\left(
I\right)  $ its limit stratum, we obtain an inclusion reversing map of
complexes
\[%
\begin{tabular}
[c]{llll}%
$\lim:$ & $B\left(  I\right)  $ & $\rightarrow$ & $\operatorname*{Strata}%
\left(  X_{0}\right)  \subset\operatorname*{Strata}\left(  Y\right)  $\\
& \multicolumn{1}{c}{$F$} & $\mapsto$ & $\left\{  \lim_{t\rightarrow0}a\left(
t\right)  \mid a\in\operatorname*{val}\nolimits^{-1}\left(
\operatorname*{int}\left(  F\right)  \right)  \right\}  $%
\end{tabular}
\
\]
where $\operatorname*{Strata}\left(  Y\right)  $ denotes the poset of closures
of toric strata of $Y$, and it holds%
\[
\lim\left(  B\left(  I\right)  \right)  \cong SP\left(  I_{0}\right)
\cong\operatorname*{Strata}\nolimits_{\Delta}\left(  I_{0}\right)
\]
We have the following correspondence%
\[%
\begin{tabular}
[c]{ccccccc}
&  & $\operatorname*{Strata}\left(  Y\right)  $ & $\rightleftarrows$ &
$\operatorname*{Poset}\left(  \Delta\right)  $ &  & \\
&  & $\cup$ &  & $\cup$ &  & \\
$B\left(  I\right)  $ & $\rightarrow$ & $\operatorname*{Strata}\left(
X_{0}\right)  $ & $\rightleftarrows$ & $\operatorname*{Strata}%
\nolimits_{\Delta}\left(  I_{0}\right)  $ & $\rightleftarrows$ & $SP\left(
I_{0}\right)  $\\
$F$ & $\mapsto$ & $\lim\left(  F\right)  =V\left(  \tau\right)  $ & $\mapsto$
& $H$ & $\mapsto$ & $\left\langle y_{G^{\ast}}\mid G\subset\Delta\text{ facet
with }H\subset G\right\rangle $\\
&  &  &  &  &  & $\Vert$\\
&  &  &  &  &  & $\left\langle y_{r}\mid r\in\Sigma\left(  1\right)  \text{,
}r\subset\tau\right\rangle $%
\end{tabular}
\
\]
where $\tau\in\Sigma=\Sigma\left(  \Delta^{\ast}\right)  $ is the cone with
$\operatorname*{int}\left(  F\right)  \subset\operatorname*{int}\left(
\tau\right)  $ and $H\subset\Delta$ is the face dual to $\tau$.
\end{proposition}

Elements of $\left(  K^{\ast}\right)  ^{n}\cong\left(  K^{\ast}\right)
^{\Sigma\left(  1\right)  }/\operatorname*{Hom}\nolimits_{\mathbb{Z}}\left(
A_{n-1}\left(  X\left(  \Sigma\right)  \right)  ,K^{\ast}\right)  $ are
represented by
\index{Cox arc|textbf}%
\textbf{Cox arcs} in $\left(  K^{\ast}\right)  ^{\Sigma\left(  1\right)  }$.
Using multi index notation we write $c\left(  t\right)  \in\left(  K^{\ast
}\right)  ^{\Sigma\left(  1\right)  }$ as%
\[
c\left(  t\right)  =\left(  c_{r}t^{j_{r}}+\operatorname*{hot}\right)
_{r\in\Sigma\left(  1\right)  }=c_{J}\cdot t^{J}+\operatorname*{hot}%
\]

\begin{remark}
Cox arcs $c_{1}\left(  t\right)  =c_{J_{1}}\cdot t^{J_{1}}+\operatorname*{hot}%
\in\left(  K^{\ast}\right)  ^{\Sigma\left(  1\right)  }$ and $c_{2}\left(
t\right)  =c_{J_{2}}\cdot t^{J_{2}}+\operatorname*{hot}\in\left(  K^{\ast
}\right)  ^{\Sigma\left(  1\right)  }$ represent the same arc in $\left(
K^{\ast}\right)  ^{n}$ if and only if
\[
q\left(  t\right)  =c_{1}\left(  t\right)  c_{2}\left(  t\right)  ^{-1}%
\]
satisfies%
\[
\left(  \prod\nolimits_{r\in\Sigma\left(  1\right)  }q_{r}^{a_{rj}}\right)
_{j=1,...,n}=\left(  1,...,1\right)
\]
in particular for the lowest order exponents we have $A^{t}\left(  J_{1}%
^{t}-J_{2}^{t}\right)  =0$.
\end{remark}

As $0$ is in the interior of%
\[
\Delta^{\ast}=\operatorname*{convexhull}\left\{  \hat{r}\mid r\in\Sigma\left(
1\right)  \right\}
\]
there are $\lambda_{r}\in\mathbb{R}_{>0}$ such that%
\[
\sum\nolimits_{r\in\Sigma\left(  1\right)  }\lambda_{r}\hat{r}=0
\]
So denoting by $\mathbb{R}_{>0}^{\Sigma\left(  1\right)  }\subset
\mathbb{R}^{\Sigma\left(  1\right)  }$ the positive orthant
\[
\left(  \lambda_{r}\right)  _{r\in\Sigma\left(  1\right)  }\in\ker\left(
A^{t}\right)  \cap\operatorname*{int}\left(  \mathbb{R}_{>0}^{\Sigma\left(
1\right)  }\right)
\]
is in the interior of the positive orthant, hence there is a basis
$w_{1},...,w_{s}\in\operatorname*{int}\left(  \mathbb{R}_{>0}^{\Sigma\left(
1\right)  }\right)  $ of $\ker\left(  A^{t}\right)  $. Denoting by
$\operatorname*{Poset}\left(  \mathbb{R}_{>0}^{\Sigma\left(  1\right)
}\right)  $ the simplex of faces of $\mathbb{R}_{>0}^{\Sigma\left(  1\right)
}$ we have:

\begin{lemma}
Suppose $a\left(  t\right)  \in\left(  K^{\ast}\right)  ^{n}$ and $c\left(
t\right)  \in\left(  K^{\ast}\right)  ^{\Sigma\left(  1\right)  }$ is a Cox
arc representing $a\left(  t\right)  =\pi\left(  c\left(  t\right)  \right)
$, and write $c\left(  t\right)  =c_{J}\cdot t^{J}+\operatorname*{hot}$ with
$c_{J}\in\left(  \mathbb{C}^{\ast}\right)  ^{\Sigma\left(  1\right)  }$. The
intersection of the affine space $J^{t}+\ker\left(  A^{t}\right)  $ with the
elements of $\operatorname*{Poset}\left(  \mathbb{R}_{>0}^{\Sigma\left(
1\right)  }\right)  $ is a poset. There is a minimal $0$-dimensional element
$\left(  J^{\prime}\right)  ^{t}$ and a Cox arc $c^{\prime}\left(  t\right)
=c_{J^{\prime}}\cdot t^{J^{\prime}}+\operatorname*{hot}$ with $a\left(
t\right)  =\pi\left(  c^{\prime}\left(  t\right)  \right)  $ such that%
\[
\lim_{t\rightarrow0}c^{\prime}\left(  t\right)  \in\mathbb{C}^{\Sigma\left(
1\right)  }-V\left(  B\left(  \Sigma\right)  \right)
\]
For any such minimal $J^{\prime}$ and any Cox arc $c^{\prime}\left(  t\right)
=c_{J^{\prime}}\cdot t^{J^{\prime}}+\operatorname*{hot}$ with $a\left(
t\right)  =\pi\left(  c^{\prime}\left(  t\right)  \right)  $ the limit point
$\lim_{t\rightarrow0}c^{\prime}\left(  t\right)  $ maps to%
\[
\lim_{t\rightarrow0}a\left(  t\right)  \in Y=\left(  \mathbb{C}^{\Sigma\left(
1\right)  }-V\left(  B\left(  \Sigma\right)  \right)  \right)  //G\left(
\Sigma\right)
\]
The limit point $\lim_{t\rightarrow0}a\left(  t\right)  $ lies in the interior
of the stratum of $Y$ given by the ideal%
\[
\left\langle y_{r}\mid r\in\Sigma\left(  1\right)  \text{ with }J_{r}^{\prime
}\neq0\right\rangle \subset S
\]

\end{lemma}

This allows to compute an ideal in the Cox ring defining the limit of a
Bergman face in terms of Cox arcs:

\begin{remark}
Let $F$ be a face of $B\left(  I\right)  $. Suppose $a\left(  t\right)
\in\operatorname*{val}\nolimits^{-1}\left(  \operatorname*{int}\left(
F\right)  \right)  $ and $c\left(  t\right)  \in\left(  K^{\ast}\right)
^{\Sigma\left(  1\right)  }$ is a Cox arc representing $a\left(  t\right)
=\pi\left(  c\left(  t\right)  \right)  $. Write $c\left(  t\right)
=c_{J}\cdot t^{J}+\operatorname*{hot}$ with $c_{J}\in\left(  \mathbb{C}^{\ast
}\right)  ^{\Sigma\left(  1\right)  }$ and let $J_{1}^{t},...,J_{q}^{t}$ be
those $0$-dimensional elements of the intersection of the affine space
$J^{t}+\ker\left(  A^{t}\right)  $ with the elements of $\operatorname*{Poset}%
\left(  \mathbb{R}_{>0}^{\Sigma\left(  1\right)  }\right)  $ such that%
\[
I_{F,i}=\left\langle y_{r}\mid r\in\Sigma\left(  1\right)  \text{ with
}J_{i,r}\neq0\right\rangle
\]
satisfies%
\[
\left(  I_{F,i}:B\left(  \Sigma\right)  ^{\infty}\right)  =I_{F,i}%
\]

Then $\lim\left(  F\right)  \subset Y$ is given by any of the ideals%
\[
\left\langle y_{r}\mid r\in\Sigma\left(  1\right)  \text{ with }J_{i,r}%
\neq0\right\rangle \subset S
\]
for $i=1,...,q$, hence $\lim\left(  F\right)  $ is also the vanishing locus of
the ideal%
\begin{align*}
I_{F} &  =I_{F,1}+...+I_{F,q}\\
&  =\left\langle y_{r}\mid r\in\Sigma\left(  1\right)  \text{ such that
}\exists i\text{ with }J_{i,r}^{\prime}\neq0\right\rangle
\end{align*}

\end{remark}

Representing $\lim\left(  F\right)  $ by the ideal $I_{F}$ has the advantage
that the intersection $\lim\left(  F_{1}\right)  \cap\lim\left(  F_{2}\right)
$ of two Bergman faces is given by the sum $I_{F_{1}}+I_{F_{2}}$ of the
corresponding ideals.

The ideal defining the stratum $\lim\left(  F\right)  $ is unique if $Y$ is simplicial.

Note that the $0$-dimensional elements of the intersection of the affine space
$J+\ker\left(  A^{t}\right)  $ with the elements of $\operatorname*{Poset}%
\left(  \mathbb{R}_{>0}^{\Sigma\left(  1\right)  }\right)  $ depend only on
the $1$-skeleton $\Sigma\left(  1\right)  $ of the fan $\Sigma$. The subset of
admissible limit strata represented by $J_{1},...,J_{q}$ are given via the
irrelevant ideal $B\left(  \Sigma\right)  $, i.e., by a subdivision of
$\Sigma\left(  1\right)  $ to build a fan $\Sigma$.

For an example see Section \ref{Sec tropicalmirror implementation}.

\begin{proposition}
By assumption the complex
\[
\lim\left(  B\left(  I\right)  \right)  \cong SP\left(  I_{0}\right)
\cong\operatorname*{Strata}\nolimits_{\Delta}\left(  I_{0}\right)
\]
is a polyhedral
\index{cell complex}%
cell complex homeomorphic to a sphere. By the map $\lim$ the complex $B\left(
I\right)  $ is a subdivision of the dual cell complex of $SP\left(
I_{0}\right)  $, hence $B\left(  I\right)  $ is homeomorphic to a sphere.

In particular $B\left(  I\right)  $ is equidimensional, connected in
codimension one and its dimension is the fiber dimension $d=\dim X_{t}$ of
$\mathfrak{X}$.
\end{proposition}

\begin{remark}
The primary decomposition of $I_{0}$ is given by
\[
I_{0}=%
{\textstyle\bigcap\nolimits_{P\in SP\left(  I_{0}\right)  _{d}}}
P=%
{\textstyle\bigcap\nolimits_{H\in\operatorname*{Strata}\nolimits_{\Delta
}\left(  I_{0}\right)  _{d}}}
\left\langle y_{r}\mid r\in\Sigma\left(  1\right)  \text{, }r\subset
\operatorname*{hull}\left(  H^{\ast}\right)  \right\rangle
\]

\end{remark}

\begin{remark}
We have the obvious representation of $I_{0}$ as the intersection of the prime
ideals
\[
I_{0}=%
{\textstyle\bigcap\nolimits_{P\in SP\left(  I_{0}\right)  _{d}}}
P=%
{\textstyle\bigcap\nolimits_{j=1}^{d}}
{\textstyle\bigcap\nolimits_{P\in SP\left(  I_{0}\right)  _{j}}}
P
\]
This intersection corresponds to a Stanley decomposition of $S/I_{0}$%
\[
S/I_{0}\cong%
{\textstyle\bigoplus\nolimits_{j=1}^{d}}
{\textstyle\bigoplus\nolimits_{P\in SP\left(  I_{0}\right)  _{j}}}
y^{D_{P}}\cdot\mathbb{C}\left[  y_{r}\mid y_{r}\notin P\right]
\]
with $D_{P}=\sum_{y_{r}\notin P}D_{r}$.
\end{remark}

The dual complex relates to the locally relevant deformations.

\begin{lemma}
The locally relevant deformations at the stratum $X_{i}$ of $X_{0}$ are the
lattice points%
\[
\left\{  m\in M\mid%
\begin{tabular}
[c]{l}%
$m\in\operatorname*{dual}\left(  G\right)  $ for some $G\in B\left(  I\right)
$ with $X_{i}\subset\lim\left(  G\right)  $\\
and $m\notin\operatorname*{dual}\left(  G^{\prime}\right)  $ for all $G\in
B\left(  I\right)  $ with $X_{i}\not \subset \lim\left(  G\right)  $%
\end{tabular}
\right\}
\]
i.e., the open star of the faces $\lim\left(  F\right)  \in
\operatorname*{dual}\left(  B\left(  I\right)  \right)  $ with $\lim\left(
F\right)  =X_{i}$.
\end{lemma}

Let $\delta_{1},...,\delta_{p}\in\operatorname*{Hom}\left(  I_{0}%
,S/I_{0}\right)  _{0}$ be a basis of the tangent space of the component of the
Hilbert scheme at $X_{0}$, which contains the tangent vector $v$ of
$\mathfrak{X}$. Let $\delta_{i}$ be a first order deformation contributing to
the tangent vector of the degeneration $\mathfrak{X}$ at $X_{0}$, i.e.,
writing $v=\sum_{i=1}^{p}\lambda_{i}\delta_{i}$ we have $\lambda_{i}\neq0$.
Then $\delta_{i}$ has to be locally relevant in at least one of the strata of
$X_{0}$, hence:

\begin{proposition}
All first order deformations contributing to the tangent space of
$\mathfrak{X}$ at $X_{0}$ are among the lattice points of
$\operatorname*{dual}\left(  B\left(  I\right)  \right)  $.
\end{proposition}

\begin{figure}
[h]
\begin{center}
\includegraphics[
height=4.0542in,
width=6.045in
]%
{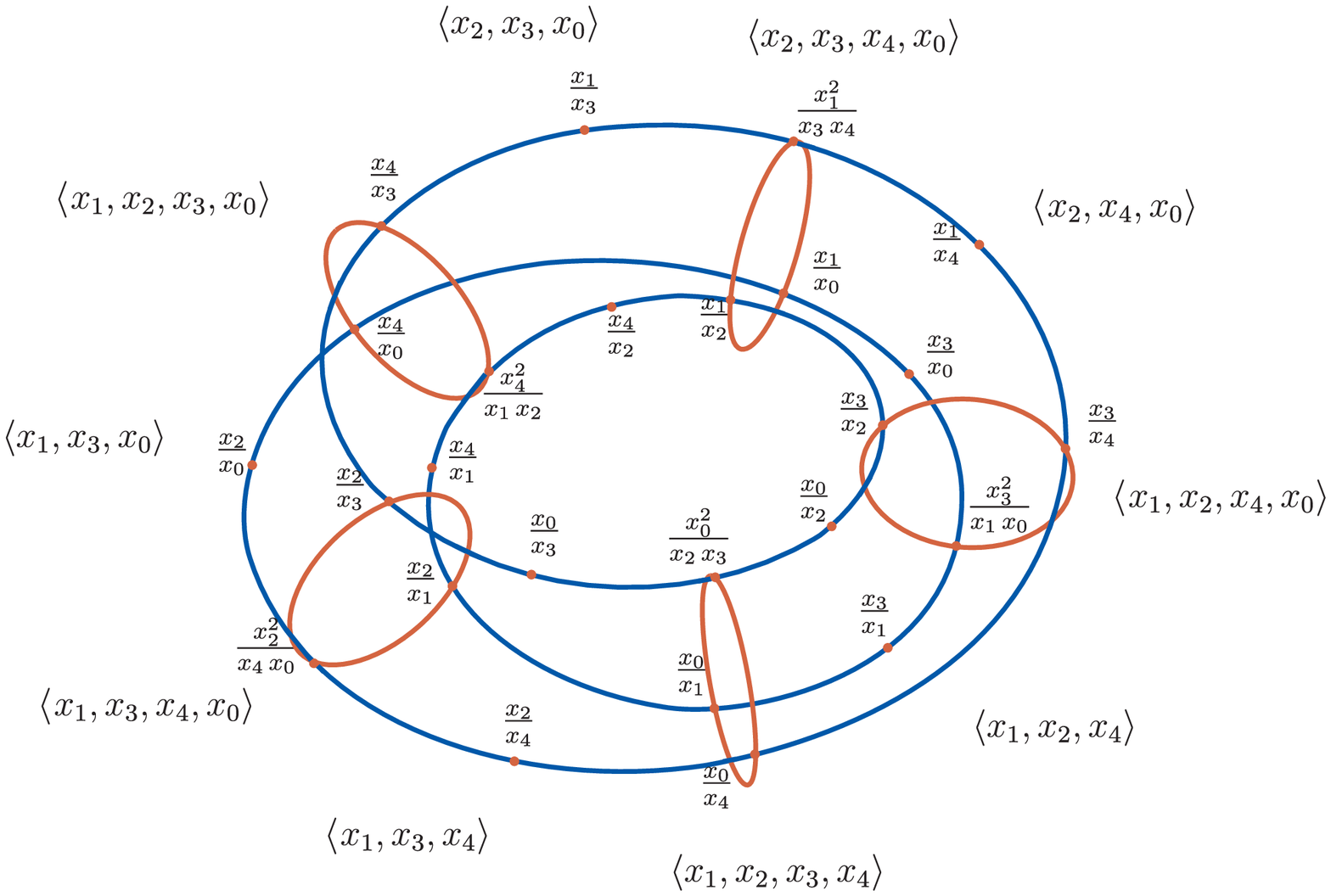}%
\caption{Complex of deformations for the Pfaffian elliptic curve}%
\label{Fig dual Pfaff ell}%
\end{center}
\end{figure}

\begin{example}
\label{Ex locally relevant}Consider degeneration $\mathfrak{X}$ of Pfaffian
elliptic curves given by the ideal defined in Example
\ref{Ex Degeneration general Pfaffian elliptic curve}. The first order
deformations of $X_{0}$ appearing in $\mathfrak{X}$ fit together in the
complex $\operatorname*{dual}\left(  B\left(  I\right)  \right)  $ consisting
of $5$ triangles and $5$ prisms. The triangles have $3$ lattice points forming
their vertices and the prisms have $7$ lattice points. Figure
\ref{Fig dual Pfaff ell} visualizes the complex $\operatorname*{dual}\left(
B\left(  I\right)  \right)  $. The faces of this complex are in one-to-one
correspondence to $SP\left(  I_{0}\right)  $, as also shown in Figure
\ref{Fig dual Pfaff ell}. For details see Section
\ref{Sec ex pfaffian elliptic curve}.

For example, the torus invariant locally relevant deformations of $X_{0}$ at
the stratum $\left(  0:0:0:0:1\right)  $ given by $\left\langle x_{0}%
,...x_{3}\right\rangle $ are%
\[%
\begin{tabular}
[c]{lllllll}%
$\frac{x_{2}}{x_{0}}$ & $\frac{x_{4}}{x_{1}}$ & $\frac{x_{4}}{x_{3}}$ &
$\frac{x_{4}}{x_{0}}$ & $\frac{x_{4}^{2}}{x_{1}x_{2}}$ & $\frac{x_{4}}{x_{2}}$
& $\frac{x_{1}}{x_{3}}$%
\end{tabular}
\]
and $\frac{x_{4}}{x_{3}},\frac{x_{4}}{x_{0}},\frac{x_{4}^{2}}{x_{1}x_{2}}$ are
the strongly locally relevant deformations.
\end{example}

%

\subsection{The special fiber $X_{0}^{\circ}$ of the mirror
degeneration\label{Sec special fiber mirror degeneration general setup}}

In the same way the spherical subcomplex $\lim\left(  B\left(  I\right)
\right)  \cong\operatorname*{Strata}\nolimits_{\Delta}\left(  I_{0}\right)
\subset\Delta$ corresponds to the special fiber monomial ideal of the
degeneration $\mathfrak{X}$, we expect the spherical subcomplex $B\left(
I\right)  \subset\nabla$ to correspond to the monomial special fiber of the
mirrror degeneration $\mathfrak{X}^{\circ}$.

By Theorem \ref{thm fano polytope} the polytope $\nabla^{\ast}$ is a Fano
polytope, so the fan $\Sigma^{\circ}=\Sigma\left(  \nabla^{\ast}\right)  $
over the faces of $\nabla^{\ast}$ defines a
\index{Q-Gorenstein}%
$\mathbb{Q}$-Gorenstein
\index{toric Fano}%
toric
\index{Fano}%
Fano variety $Y^{\circ}=X\left(  \Sigma^{\circ}\right)  $.

Denote by $S^{\circ}=\mathbb{C}\left[  z_{r}\mid r\in\Sigma^{\circ}\left(
1\right)  \right]  $ the Cox ring of $Y^{\circ}$, so the variables of
$S^{\circ}$ correspond to the vertices of the polytope $\nabla^{\ast}$ of
first order deformations appearing in $\mathfrak{X}$.

Define the monomial ideal%
\[
I_{0}^{\circ}=\left\langle
{\displaystyle\prod\limits_{v\in J}}
z_{v}\mid J\subset\Sigma^{\circ}\left(  1\right)  \text{ with }%
\operatorname*{supp}\left(  B\left(  I\right)  \right)  \subset%
{\displaystyle\bigcup\limits_{v\in J}}
F_{v}\right\rangle \subset S^{\circ}%
\]
where $F_{v}$ is the facet of $\nabla$ corresponding to the ray $v$ of the
normal fan $\Sigma^{\circ}=\operatorname*{NF}\left(  \nabla\right)  $ so:

\begin{proposition}
$I^{\circ}$ is a reduced monomial ideal. As $B\left(  I\right)  $ is a cell
complex homeomorphic to a sphere, $I^{\circ}$ defines a Calabi-Yau variety
$X_{0}^{\circ}\subset Y^{\circ}$, which is the union of toric strata of
$Y^{\circ}$%
\begin{align*}
I_{0}^{\circ}  &  =%
{\textstyle\bigcap\nolimits_{F\in B\left(  I\right)  _{d}}}
\left\langle z_{G^{\ast}}\mid G\text{ a facet of }\nabla\text{ with }F\subset
G\right\rangle \\
&  =%
{\textstyle\bigcap\nolimits_{H\in\left(  \operatorname*{dual}\left(  B\left(
I\right)  \right)  \right)  _{n-1-d}}}
\left\langle z_{r}\mid r\in\Sigma^{\circ}\left(  1\right)  ,\text{ }%
r\subset\operatorname*{hull}\left(  H\right)  \in\Sigma^{\circ}\right\rangle
\end{align*}

\end{proposition}

\subsection{First order mirror degeneration $\mathfrak{X}^{\circ}$ with
special fiber $X_{0}^{\circ}$\label{Sec first order mirror degeneration}}

In the same way as the lattice points of $\left(  B\left(  I\right)  \right)
^{\ast}\subset\nabla^{\ast}$ are the first order deformations of $I_{0}$
appearing in $I$, we consider the lattice points of $\left(  \lim\left(
B\left(  I\right)  \right)  \right)  ^{\ast}\subset\Delta^{\ast}$ as elements
in $\operatorname*{Hom}\left(  I_{0}^{\circ},S^{\circ}/I_{0}^{\circ}\right)
_{0}$, i.e., as first order deformations of $X_{0}^{\circ}$.

Note that the deformations of $X^{\circ}$ are represented independently of the
embedding of $X_{0}^{\circ}$ in $Y^{\circ}$. Indeed the deformations of
$X_{0}^{\circ}$ depend on the embedding of $X_{0}$.

A first order deformation of $X_{0}^{\circ}$%
\[
\mathfrak{X}^{1\circ}\subset Y^{\circ}\times\operatorname*{Spec}%
\mathbb{C}\left[  t\right]  /\left\langle t^{2}\right\rangle
\]
is defined by the ideal%
\[
I^{1\circ}=\left\langle m+t\cdot\sum_{\alpha\in\operatorname*{supp}\left(
\left(  \lim\left(  B\left(  I\right)  \right)  \right)  ^{\ast}\right)  \cap
N}a_{\alpha}\cdot\alpha\left(  m\right)  \mid m\in I_{0}^{\circ}\right\rangle
\subset\mathbb{C}\left[  t\right]  /\left\langle t^{2}\right\rangle \otimes
S^{\circ}%
\]
with general coefficients $a_{\alpha}$.

\begin{proposition}
The special fiber Gr\"{o}bner cone of $\mathfrak{X}^{1\circ}$ gives back
$\Delta\subset M_{\mathbb{R}}$, i.e.,%
\[
\Delta=C_{I_{0}^{\circ}}\left(  I^{1\circ}\right)  \cap\left\{  w_{t}%
=1\right\}
\]
Let $I_{gen}^{\circ}\subset S^{\circ}$ be the ideal of the general fiber of
$\mathfrak{X}^{1\circ}$. As $0\in N$ is the unique interior point of the Fano
polytope $\Delta^{\ast}$ the Hilbert point of $I_{0}^{\circ}$ lies in the
interior of the state polytope of $I_{gen}^{\circ}$, i.e.,%
\[
H\left(  I_{0}^{\circ}\right)  \in\operatorname*{int}\left(
\operatorname*{State}\left(  I_{gen}^{\circ}\right)  \right)  \subset
N_{\mathbb{R}}%
\]

\end{proposition}

\begin{conjecture}
If $\mathfrak{X}^{\circ}$ is the flat family with tangent direction in
$\mathfrak{X}^{1\circ}$, then the general fiber of $\mathfrak{X}^{\circ}$ and
the general fiber of $\mathfrak{X}$ form a mathematical mirror pair.
\end{conjecture}

\begin{remark}
Suppose $\mathfrak{X}$ and $\mathfrak{X}^{\prime}$ are degenerations as
defined above with fibers in $Y=X\left(  \Sigma\left(  P\right)  \right)  $
such that the special fiber ideals $I_{0}$ and $I_{0}^{\prime}$ define the
same subcomplex of $P^{\ast}$ and $\mathfrak{X}_{1}$ and $\mathfrak{X}_{2}$
involve the same first order deformations, then the tropical mirror
construction applied to $\mathfrak{X}$ respectively $\mathfrak{X}^{\prime}$
will lead to the same result. Hence, e.g., passing from $I_{0}$ to
$I_{0}^{\Sigma}$, i.e. the non simplicial analogue of saturation, and from
\[
I=\left\langle m+t\cdot%
{\textstyle\sum\limits_{\delta\in\operatorname*{supp}\left(
\operatorname*{dual}\left(  B\left(  I\right)  \right)  \right)  \cap M}}
c_{\delta}\cdot\delta\left(  m\right)  \mid m\in I_{0}\right\rangle
\]
to%
\[
I^{\prime}=\left\langle m+t\cdot%
{\textstyle\sum\limits_{\delta\in\operatorname*{supp}\left(
\operatorname*{dual}\left(  B\left(  I\right)  \right)  \right)  \cap M}}
c_{\delta}\cdot\delta\left(  m\right)  \mid m\in I_{0}^{\Sigma}\right\rangle
\]
does not change the geometry of the degeneration and the objects involved in
the tropical mirror construction.
\end{remark}

\subsection{Remarks on orbifolding mirror
families\label{sec orbifolding mirror families}}

Suppose $Y=X\left(  \Sigma\right)  $ is projective space $\mathbb{P}^{n}$
given by the fan $\Sigma\subset N_{\mathbb{R}}$ with the rays generated by
$\left(  1,0,...,0\right)  ,...,\left(  0,...,0,1\right)  ,\left(
-1,....,-1\right)  \in N$. Let $\left(  \delta_{i}\right)  _{i}$ be the torus
invariant basis of the space of first order deformations spanning the tangent
space of $\mathfrak{X}$ at $X_{0}$, so $\left\{  \delta_{i}\mid i\right\}
=\operatorname*{dual}\left(  B\left(  I\right)  \right)  \cap M$. Denote by
$c$ the codimension of the fibers in $Y$.

Representing each deformation $\delta_{i}$ as a Cox Laurent monomial
$\delta_{i}=\frac{a_{i}}{b_{i}}$ with relative prime Cox monomials $a_{i}$ and
$b_{i}$, there is a preordering on $\left\{  \delta_{i}\mid i\right\}  $ by
divisibility of the denominators $b_{i}$, i.e., $\frac{a_{i}}{b_{i}}%
=\delta_{i}\leq\delta_{j}=\frac{a_{j}}{b_{j}}\Leftrightarrow b_{i}\mid b_{j}$.

If $\delta=\frac{a}{b}$ is deformation with relative prime Cox monomials $a$
and $b$, then $\delta$ is called pure, if $a$ has the form $a=y_{s}^{d}$ for
some homogeneous variable $y_{s}$ of $Y$ and for some $d>0$.

For each $0$-dimensional stratum $p$ of $Y$ denote by%
\[
\mathfrak{D}_{p}=\left\{  m\in M\mid m\in F^{\ast}\text{ with }\lim
F=p\right\}
\]
the set of all torus invariant first order deformations of $X_{0}$ in
$\mathfrak{X}$ corresponding to $p$.

\begin{definition}
Let $\mathfrak{F}$ be a set of non-trivial deformations of $X_{0}$ in
$\mathfrak{X}$ corresponding to vertices of faces of $\operatorname*{dual}%
\left(  B\left(  I\right)  \right)  $ and denote by $\mathcal{R}$ the
corresponding set of rays of $\Sigma^{\circ}$. We call $\mathfrak{F}$ a set
\index{Fermat deformations|textbf}%
of \textbf{Fermat deformations} of $\mathfrak{X}$, if the following conditions
are satisfied:

\begin{itemize}
\item $\left\vert \mathcal{R}\cap\mathfrak{D}_{p}\right\vert =1$ for all
$0$-dimensional strata $p$ of $Y$ (in particular $\left\vert \mathcal{R}%
\right\vert =n+1$).

\item The convexhull of $\left\{  \hat{r}\mid r\in\mathcal{R}\right\}  $ is a
polytope of dimension $n=\dim\left(  N_{\mathbb{R}}\right)  $ containing $0$
in its interior, i.e.,$\mathcal{R}$ spans a projective fan $\hat{\Sigma
}^{\circ}$ (note that this fan is uniquely determined by $\mathcal{R}$).

\item The elements of $\mathfrak{F}$ are incomparable with respect to the
preordering $\leq$.
\end{itemize}
\end{definition}

Let $\hat{Y}^{\circ}=X\left(  \hat{\Sigma}^{\circ}\right)  $ be the toric Fano
variety defined by $\hat{\Sigma}^{\circ}$ with Cox ring $\hat{S}^{\circ
}=\mathbb{C}\left[  z_{r}\mid r\in\hat{\Sigma}^{\circ}\left(  1\right)
\right]  $ and%
\[
0\rightarrow M_{\mathbb{R}}\rightarrow\mathbb{Z}^{\mathcal{R}}\overset{\deg
}{\rightarrow}A_{n-1}\left(  \mathcal{R}\right)  \rightarrow0
\]
the corresponding presentation of $A_{n-1}\left(  \hat{Y}^{\circ}\right)
=A_{n-1}\left(  \mathcal{R}\right)  \cong H\oplus\mathbb{Z}$ with finite $H$
and let $\left(  h_{r},d_{r}\right)  =\deg\left(  z_{r}\right)  $. As $0$ is
in the interior of the convex hull of $\left\{  \hat{r}\mid r\in\hat{\Sigma
}^{\circ}\left(  1\right)  \right\}  $, we can assume that the $d_{r}$ are
positive integers. For all $w\in\left(  \lim F\right)  ^{\ast}\cap N$, $F\in
B\left(  I\right)  $ the Laurent monomials
\[%
{\textstyle\prod\nolimits_{r\in\hat{\Sigma}^{\circ}\left(  1\right)  }}
z_{r}^{\left\langle \hat{r},w\right\rangle }%
\]
are of degree $0$ with respect to the grading $\deg\left(  z_{r}\right)
=d_{r}$.

Note that for complete intersections the deformations corresponding to
vertices of faces of $\operatorname*{dual}\left(  B\left(  I\right)  \right)
$ are pure. The set $\mathfrak{F}$ is not unique in general, see for example
the complete intersection Calabi-Yau of degree $12$ in $\mathbb{P}^{6}$. If
$\mathfrak{X}$ is a degeneration of complete intersections of codimension
$c=2$, then there is a unique set of Fermat deformations $\mathfrak{F}$, which
is the set of maximal elements of $\operatorname*{dual}\left(  B\left(
I\right)  \right)  \cap M$ with respect to the preordering $\leq$ defined above.

\begin{remark}
With the notations of the preceding sections, let $\mathfrak{F}$ be a set of
Fermat deformations of $\mathfrak{X}$. Then
\[
\hat{P}^{\circ}=\operatorname*{convexhull}\left\{  A^{-1}\left(
\delta\right)  \mid\delta\in\mathfrak{F}\right\}  \subset P^{\circ}%
=\nabla^{\ast}%
\]
is a Fano polytope. $\hat{Y}^{\circ}$ is an orbifold $\hat{Y}^{\circ
}=\mathbb{P}\left(  d_{1},...,d_{n+1}\right)  /G$ with the $d_{i}$ defined as
above. Let $Y^{\circ}=X\left(  \Sigma^{\circ}\right)  \rightarrow X\left(
\hat{\Sigma}^{\circ}\right)  =\hat{Y}^{\circ}$ be a birational map contracting
all divisors of $Y^{\circ}$ corresponding to deformations not in
$\mathfrak{F}$. Then the first order flat family $\mathfrak{\hat{X}}^{1\circ
}\subset\hat{Y}^{\circ}\times\operatorname*{Spec}\mathbb{C}\left[  t\right]
/\left\langle t^{2}\right\rangle $ induced by $\mathfrak{X}^{1\circ}\subset
Y^{\circ}\times\operatorname*{Spec}\mathbb{C}\left[  t\right]  /\left\langle
t^{2}\right\rangle $ has special fiber given by%
\[
\hat{I}_{0}^{\circ}=\left\langle \hat{m}\mid\exists\text{ minimal generator
}m\text{ of }I_{0}^{\circ}\text{ divisible by }%
{\textstyle\prod\nolimits_{r\in\hat{\Sigma}^{\circ}\left(  1\right)  }}
z_{r}\right\rangle \subset\hat{S}^{\circ}%
\]
and involves the deformations
\[
\left\{
{\textstyle\prod\nolimits_{r\in\hat{\Sigma}^{\circ}\left(  1\right)  }}
z_{r}^{\left\langle \hat{r},w\right\rangle }\mid w\in\left(  \lim\left(
B\left(  I\right)  \right)  \right)  ^{\ast}\cap N\right\}
\]

\end{remark}

\section{Tropical
\index{mirror construction}%
mirror construction for the example of Pfaffian Calabi-Yau
varieties\label{Sec tropical mirror construction for Pfaffian Calabi-Yau varieties}%
}

\subsection{Pfaffian Calabi-Yau varieties\label{1PfaffianCalabiYauThreefolds}}

\begin{definition}
A subscheme $X$ of $\mathbb{P}_{K}^{n}$ of codimension $3$ is called
\index{Pfaffian|textbf}%
\textbf{Pfaffian subscheme} if there is

\begin{enumerate}
\item a vector bundle $E$ on $\mathbb{P}_{K}^{n}$ of rank $2k+1$ for some
$k\in\mathbb{Z}_{\geq0}$

\item and a
\index{skew symmetric}%
skew symmetric map $\varphi:\mathcal{E}\left(  -t\right)  \rightarrow
\mathcal{E}^{\ast}$, where $\mathcal{E}=\mathcal{O}_{\mathbb{P}_{K}^{n}%
}\left(  E\right)  $ such that

\begin{enumerate}
\item $\varphi$ is generically of rank $2k$

\item $\varphi$ degenerates to rank $2k-2$ in the expected codimension $3$
\end{enumerate}

\item $X$ is scheme theoretically the
\index{degeneracy locus}%
degeneracy locus of $\varphi$.
\end{enumerate}
\end{definition}

\begin{theorem}
[Buchsbaum-Eisenbud]\cite{BE Algebra structures for finite free resolutions
and some structure theorems for ideals of codimension 3},\newline\cite{Okonek
Notes on Varieties of Codimension 3 in mathbbP^N}, \cite{Walter Pfaffian
subschemes}\label{1BuchsbaumEisenbud} A Pfaffian subscheme $X$ of
$\mathbb{P}_{K}^{n}$ has a locally free resolution
\[
0\rightarrow\mathcal{O}_{\mathbb{P}_{K}^{n}}\left(  -t-2s\right)
\overset{\psi^{\ast}\left(  -t-2s\right)  }{\rightarrow}\mathcal{E}\left(
-t-s\right)  \overset{\varphi}{\rightarrow}\mathcal{E}^{\ast}\left(
-s\right)  \overset{\psi}{\rightarrow}\mathcal{O}_{\mathbb{P}_{K}^{n}%
}\rightarrow\mathcal{O}_{X}\rightarrow0
\]
where $s=c_{1}\left(  \mathcal{E}\right)  +kt$, and $\psi$ is locally given by
the
\index{Pfaffians}%
Pfaffians of order $2k$ of $\varphi$.
\end{theorem}

\begin{remark}
$X$ is locally
\index{Gorenstein}%
Gorenstein with
\[
\omega_{X}^{\circ}\cong\mathcal{O}_{X}\left(  t+2s-n-1\right)
\]
Thus $\omega_{X}^{\circ}\cong\mathcal{O}_{X}$ if and only if $t+2s=n+1$.
\end{remark}

\begin{theorem}
\label{thm walter}\cite{Walter Pfaffian subschemes} Let $K$ be a field with
$\operatorname*{char}\left(  K\right)  \neq2$. If $X\subset\mathbb{P}_{K}^{n}$
is an equidimensional, locally Gorenstein subscheme of dimension $n-3$, which
is subcanonical, i.e., $\omega_{X}^{\circ}\cong\mathcal{O}_{X}\left(
l\right)  $ for some integer $l$, then $X$ is Pfaffian if and only if the
following condition is satisfied%
\[
n\equiv0\operatorname{mod}4\text{ and }l=2s\text{ even }\Rightarrow\text{
}\chi\left(  \mathcal{O}_{X}\left(  s\right)  \right)  \text{ is even}%
\]

\end{theorem}

\begin{corollary}
A codimension $3$ subscheme of $\mathbb{P}^{6}$ is Pfaffian if and only if it
is locally Gorenstein and subcanonical.
\end{corollary}

\begin{example}
Using this construction, we get the following projectively Gorenstein
\index{Pfaffian}%
Pfaffian Calabi-Yau threefolds%
\[%
\begin{tabular}
[c]{llllll}%
$\mathcal{E}^{\ast}$ & $\operatorname*{rank}\left(  \mathcal{E}\right)  $ &
$\deg\left(  X\right)  $ & $h^{1,2}\left(  X\right)  $ & $h^{1,1}\left(
X\right)  $ & $\chi\left(  X\right)  $\\
&  &  &  &  & \\
$2\mathcal{O}\left(  1\right)  \oplus\mathcal{O}$ & $3$ & $12$ & $73$ & $1$ &
$-144$\\
$\mathcal{O}\left(  1\right)  \oplus4\mathcal{O}$ & $5$ & $13$ & $61$ & $1$ &
$-120$\\
$7\mathcal{O}$ & $7$ & $14$ & $50$ & $1$ & $-98$%
\end{tabular}
\]

\end{example}

\begin{example}
In \cite{Tonoli Canonical surfaces in mathbbP^5 and CalabiYau threefolds in
mathbbP^6} families of non
\index{projectively Cohen-Macaulay}%
projectively Cohen-Macaulay
\index{Pfaffian}%
Pfaffian Calabi-Yau threefolds with the following data were constructed and it
is shown that generic elements of each family are smooth:%
\[%
\begin{tabular}
[c]{llllll}%
$\mathcal{E}^{\ast}$ & $\operatorname*{rank}\left(  \mathcal{E}\right)  $ &
$\deg\left(  X\right)  $ & $h^{1,2}\left(  X\right)  $ & $h^{1,1}\left(
X\right)  $ & $\chi\left(  X\right)  $\\
&  &  &  &  & \\
$\Omega^{1}\left(  1\right)  \oplus3\mathcal{O}$ & $9$ & $15$ & $40$ & $1$ &
$-78$\\
$Syz^{1}\left(  M\right)  $ & $11$ & $16$ & $31$ & $1$ & $-60$\\
$Syz^{1}\left(  M^{\prime}\right)  $ & $13$ & $17$ & $23$ & $1$ & $-44$%
\end{tabular}
\
\]
where $M$ is a generic module of length $2$ generated in degree $-1$ with
Hilbert function $\left(  2,1,0,...\right)  $, and $M^{\prime}$ is a special
module of length $2$ generated in degree $-1$ with Hilbert function $\left(
3,5,0,...\right)  $ (the generic choice of $M^{\prime}$ gives a bundle
$\mathcal{E}=Syz^{1}\left(  M^{\prime}\right)  $, which does not admit any
alternating map $\mathcal{E}^{\ast}\left(  -1\right)  \rightarrow\mathcal{E}%
$). There are $3$ unirational families of smooth
\index{Pfaffian}%
Pfaffian Calabi-Yau threefolds of degree $17$ and all $3$ families have
$h^{1,2}\left(  X\right)  =23$. The
\index{Hodge numbers}%
Hodge numbers were obtained via computer algebra.
\end{example}

The Pfaffian given by $\mathcal{E=}2\mathcal{O}\left(  1\right)
\oplus\mathcal{O}$ is a complete intersection and the mirror construction is
given in Section
\ref{Sec tropical mirror construction for complete intersections}. In the
following, we will be concerned with the remaining two projectively Gorenstein examples.

\subsection{Deformations of Pfaffian
varieties\label{Sec deformations of Pfaffians}}

Let $Y=X\left(  \Sigma\right)  $ be a toric
\index{Fano}%
Fano
\index{toric Fano}%
variety given by the fan $\Sigma\subset N_{\mathbb{R}}$ over the Fano polytope
$P\subset N_{\mathbb{R}}$ and denote its Cox ring by $S$.

Suppose that $I_{0}=\left\langle m_{1},...,m_{r}\right\rangle \subset S$ is an
ideal generated by monomials $m_{i}\in H^{0}\left(  Y,\mathcal{O}_{Y}\left(
E_{i}\right)  \right)  $, $i=1,...,r$, which has a Pfaffian resolution%
\[
0\rightarrow\mathcal{O}_{Y}\left(  K_{Y}\right)  \rightarrow\mathcal{F}\left(
K_{Y}\right)  \overset{\varphi^{0}}{\rightarrow}\mathcal{F}^{\ast}\overset
{m}{\rightarrow}\mathcal{O}_{Y}%
\]
with%
\begin{align*}
\mathcal{F}  &  =\mathcal{O}_{Y}\left(  E_{1}\right)  \oplus...\oplus
\mathcal{O}_{Y}\left(  E_{r}\right) \\
m  &  =\left(  m_{1},...,m_{r}\right)
\end{align*}
and
\[
\varphi^{0}\in\bigwedge\nolimits^{2}\mathcal{F}^{\ast}\left(  -K_{Y}\right)
\]

Suppose $\mathfrak{X}^{1}\subset Y\times\operatorname*{Spec}\mathbb{C}\left[
t\right]  /\left\langle t^{2}\right\rangle $ is a first order deformation of
$I_{0}$ defined by%
\[
I^{1}=\left\langle f_{j}^{1}=t\cdot g_{j}+m_{j}\mid j=1,...,r\right\rangle
\subset\mathbb{C}\left[  t\right]  /\left\langle t^{2}\right\rangle \otimes S
\]
with $g_{j}\in S_{\left[  E_{j}\right]  }$.

Denote by $\overline{\pi}_{1}:Y\times\operatorname*{Spec}\mathbb{C}\left[
t\right]  /\left\langle t^{2}\right\rangle \rightarrow Y$ the projection on
the first component and by%
\[
K^{1}=K_{Y\times\operatorname*{Spec}\left(  \mathbb{C}\left[  t\right]
/\left\langle t^{2}\right\rangle \right)  /\operatorname*{Spec}\left(
\mathbb{C}\left[  t\right]  /\left\langle t^{2}\right\rangle \right)  }%
\]
the relative canonical sheaf. Then flatness of $\mathfrak{X}^{1}$ gives a lift
of the syzygies of $m$, hence a Pfaffian resolution of $I^{1}$ of the form%
\[
0\rightarrow\mathcal{O}_{Y\times\operatorname*{Spec}\mathbb{C}\left[
t\right]  /\left\langle t^{2}\right\rangle }\left(  K^{1}\right)
\rightarrow\mathcal{E}^{1}\left(  K^{1}\right)  \overset{\varphi^{1}%
}{\rightarrow}\left(  \mathcal{E}^{1}\right)  ^{\ast}\overset{f^{1}%
}{\rightarrow}\mathcal{O}_{Y\times\operatorname*{Spec}\mathbb{C}\left[
t\right]  /\left\langle t^{2}\right\rangle }%
\]
with $f^{1}=\left(  f_{1}^{1},...,f_{r}^{1}\right)  $ and
\[
\mathcal{E}^{1}=\overline{\pi}_{1}^{\ast}\mathcal{F}%
\]
and $\varphi^{1}$ is skew symmetric by the Theorem of Buchsbaum-Eisenbud ,
i.e.,%
\[
\varphi^{1}\in\bigwedge\nolimits^{2}\mathcal{E}^{1}\left(  -K^{1}\right)
\]

Denote by $\pi_{1}:Y\times\operatorname*{Spec}\mathbb{C}\left[  \left[
t\right]  \right]  \rightarrow Y$ the projection on the first component and by%
\[
K=K_{\left(  Y\times\operatorname*{Spec}\mathbb{C}\left[  \left[  t\right]
\right]  \right)  /\operatorname*{Spec}\mathbb{C}\left[  \left[  t\right]
\right]  }%
\]
the relative canonical sheaf. Let
\[
\mathcal{E}=\pi_{1}^{\ast}\mathcal{F}%
\]
and let $\varphi\in\bigwedge\nolimits^{2}\mathcal{E}\left(  -K\right)  $ be a
representative of $\varphi^{1}$ of $t$-degree $1$. Defining $f=\left(
f_{1},...,f_{r}\right)  $ as the Pfaffians of $\varphi$, one obtains a
Pfaffian resolution of the ideal generated by $f_{1},...,f_{r}$%
\[
0\rightarrow\mathcal{O}_{Y\times\operatorname*{Spec}\mathbb{C}\left[  \left[
t\right]  \right]  }\left(  K\right)  \rightarrow\mathcal{E}\left(  K\right)
\overset{\varphi}{\rightarrow}\mathcal{E}^{\ast}\overset{f}{\rightarrow
}\mathcal{O}_{Y\times\operatorname*{Spec}\mathbb{C}\left[  \left[  t\right]
\right]  }%
\]
hence a lift of $\mathfrak{X}^{1}$ to a flat family $\mathfrak{X}\subset
Y\times\operatorname*{Spec}\mathbb{C}\left[  \left[  t\right]  \right]  $, so:

\begin{proposition}
The
\index{deformation}%
deformations of $I_{0}$ are unobstructed and the base space is smooth.
\end{proposition}

By the same argument one obtains:

\begin{proposition}
Let $Y=X\left(  \Sigma\right)  $ be a toric
\index{Fano}%
Fano
\index{toric Fano}%
variety given by the fan $\Sigma\subset N_{\mathbb{R}}$ over the Fano polytope
$P\subset N_{\mathbb{R}}$. Denote the Cox ring of $Y$ by $S$.

Suppose that $X_{0}\subset Y$ is defined by an ideal $I_{0}=\left\langle
m_{1},...,m_{r}\right\rangle \subset S$, which is generated
\index{degeneration}%
by
\index{monomial degeneration}%
monomials $m_{i}\in H^{0}\left(  Y,\mathcal{O}_{Y}\left(  E_{i}\right)
\right)  $, $i=1,...,r$ and has a Pfaffian resolution%
\[
0\rightarrow\mathcal{O}_{Y}\left(  K_{Y}\right)  \rightarrow\mathcal{F}\left(
K_{Y}\right)  \overset{\varphi^{0}}{\rightarrow}\mathcal{F}^{\ast}\overset
{m}{\rightarrow}\mathcal{O}_{Y}%
\]
with $m=\left(  m_{1},...,m_{r}\right)  $, $\mathcal{F}=\mathcal{O}_{Y}\left(
E_{1}\right)  \oplus...\oplus\mathcal{O}_{Y}\left(  E_{r}\right)  $ and
$\varphi^{0}\in\bigwedge\nolimits^{2}\mathcal{F}^{\ast}\left(  -K_{Y}\right)
$.

Denote by $\pi_{1}:Y\times\operatorname*{Spec}\mathbb{C}\left[  \left[
t\right]  \right]  \rightarrow Y$ the projection on the first component.
Suppose that $\mathfrak{X}\subset Y\times\operatorname*{Spec}\mathbb{C}\left[
\left[  t\right]  \right]  $ is given by an ideal $I\subset\mathbb{C}\left[
\left[  t\right]  \right]  \otimes S$, which has a
\index{Pfaffian}%
Pfaffian resolution%
\[
0\rightarrow\mathcal{O}_{Y\times\operatorname*{Spec}\mathbb{C}\left[  \left[
t\right]  \right]  }\left(  K\right)  \rightarrow\mathcal{E}\left(  K\right)
\overset{\varphi}{\rightarrow}\mathcal{E}^{\ast}\rightarrow\mathcal{O}%
_{Y\times\operatorname*{Spec}\mathbb{C}\left[  \left[  t\right]  \right]  }%
\]
with $\mathcal{E}=\pi_{1}^{\ast}\mathcal{F}$ and $K=K_{\left(  Y\times
\operatorname*{Spec}\mathbb{C}\left[  \left[  t\right]  \right]  \right)
/\operatorname*{Spec}\mathbb{C}\left[  \left[  t\right]  \right]  }$, i.e.,
$I$ is generated by the Pfaffians of $\varphi\in\bigwedge\nolimits^{2}%
\mathcal{E}\left(  -K\right)  $. Suppose that $X_{0}\cong\mathfrak{X}%
\times_{k\left[  \left[  t\right]  \right]  }\operatorname*{Spec}k$.

Then $\mathfrak{X}$ is a flat degeneration of Pfaffian Calabi-Yau varieties
with fibers polarized in $Y$ and
\index{special fiber}%
special fiber $X_{0}$.
\end{proposition}

\begin{corollary}
The
\index{deformation}%
families given in Example
\ref{Ex Degeneration general Pfaffian elliptic curve} (monomials degenerations
of a general Pfaffian elliptic curve in $\mathbb{P}^{3}$), in Example
\ref{Ex degeneration 1 parameter Pfaffian degree 14} ($1$-parameter
degeneration of a Pfaffian Calabi-Yau threefold in an orbifold of
$\mathbb{P}^{6}$), in Example \ref{Ex Degeneration general Pfaffian degree 14}
(monomial degeneration of a general Pfaffian Calabi-Yau threefold of degree
$14$ in $\mathbb{P}^{6}$) and in Example
\ref{Ex Degeneration general Pfaffian degree 13} (monomial degeneration of a
general Pfaffian Calabi-Yau threefold of degree $13$ in $\mathbb{P}^{6}$) are
flat and satisfy the
\index{genericy condition}%
genericy condition on the tangent direction, given in Section
\ref{genericy condition}.
\end{corollary}

\subsection{Tropical mirror construction for the Pfaffian elliptic
curve\label{Sec ex pfaffian elliptic curve}%
\index{elliptic curve!complete intersection}%
}

\subsubsection{Setup}

Let $Y=\mathbb{P}^{
4
}=X\left(  \Sigma\right)  $, $\Sigma
=\operatorname*{Fan}\left(  P\right)  =NF\left(  \Delta\right)  \subset
N_{\mathbb{R}}$ with the Fano polytope $P=\Delta^{\ast}$ given by%
\[
\Delta=\operatorname*{convexhull}\left(  

\right)  \subset
M_{\mathbb{R}}%
\]
and let%
\[
S=\mathbb{C}
[x_0, x_1, x_2, x_3, x_4]
\]
be the Cox ring of $Y$ with the variables%
\[
%
\]
associated to the rays of $\Sigma$.

Consider the degeneration $\mathfrak{X}\subset Y\times\operatorname*{Spec}%
\mathbb{C}\left[  t\right]  $ of 
Pfaffian
elliptic curves
 with Buchsbaum-Eisenbud
resolution%
\[
0\rightarrow\mathcal{O}_{Y}\left(  -
5
\right)  \rightarrow
\mathcal{E}\left(  -3\right)  \overset{A_{t}}{\rightarrow}\mathcal{E}^{\ast
}\left(  -2\right)  \rightarrow\mathcal{O}_{Y}\rightarrow\mathcal{O}_{X_{t}%
}\rightarrow0
\]
where%
\[
\mathcal{E}=\text{
$5\mathcal{O}$
}%
\]%
\[
A_{t}=A_{0}+t\cdot A
\]%
\[
A_{0}=
\left [%
\right ]
\]
the monomial special fiber of $\mathfrak{X}$ is given by%
\[
I_{0}=\left\langle 
%
\]
The total space of the degeneration $\mathfrak{X}$ is a local complete intersection.

The induced first order degeneration%
\[
\mathfrak{X}^{1}\subset Y\times\operatorname*{Spec}\mathbb{C}\left[  t\right]
/\left\langle t^{2}\right\rangle
\]
is given by the ideal $I\subset S\otimes\mathbb{C}\left[  t\right]
/\left\langle t^{2}\right\rangle $ with $I_{0}$-reduced generators of degrees
$
2,2,2,2,2
$
\[
\left\{  
%
\]

The corresponding syzygy matrix is given by%
\[
%
\]

\subsubsection{Special fiber Gr\"{o}bner cone}

The space of first order deformations of $\mathfrak{X}$ has dimension
$
25
$ and the deformations represented by the Cox Laurent monomials
\begin{center}
%
\end{center}
\noindent form a torus invariant basis. These deformations give linear
inequalities defining the special fiber Gr\"{o}bner cone, i.e.,%
\[
C_{I_{0}}\left(  I\right)  =\left\{  \left(  w_{t},w_{y}\right)  \in
\mathbb{R}\oplus N_{\mathbb{R}}\mid\left\langle w,A^{-1}\left(  m\right)
\right\rangle \geq-w_{t}\forall m\right\}
\]
for above monomials $m$ and the presentation matrix
\[
A=
\left [%
\right ]
\]
of $A_{
3
}\left(  Y\right)  $.

The vertices of special fiber polytope
\[
\nabla=C_{I_{0}}\left(  I\right)  \cap\left\{  w_{t}=1\right\}
\]
and the number of faces, each vertex is contained in, are given by the table
\begin{center}
%
\end{center}
\noindent The number of faces of $\nabla$ and their $F$-vectors are
\begin{center}
%
\end{center}

The dual $P^{\circ}=\nabla^{\ast}$ of $\nabla$ is a Fano polytope with vertices

\begin{center}
%
\end{center}

\noindent and the $F$-vectors of the faces of $P^{\circ}$ are

\begin{center}
%
\end{center}

The polytopes $\nabla$ and $P^{\circ}$ have $0$ as the unique interior lattice
point, and the Fano polytope $P^{\circ}$ defines the embedding toric Fano
variety $Y^{\circ}=X\left(  \operatorname*{Fan}\left(  P^{\circ}\right)
\right)  $ of the fibers of the mirror degeneration. The dimension of the
automorphism group of $Y^{\circ}$ is
\[
\dim\left(  \operatorname*{Aut}\left(  Y^{\circ}\right)  \right)
=
4
\]
Let%
\[
S^{\circ}=\mathbb{C}
[y_{1},y_{2},...,y_{15}]
\]
be the Cox ring of $Y^{\circ}$ with variables

\begin{center}
%
\end{center}

\subsubsection{Bergman subcomplex}

Intersecting the tropical variety of $I$ with the special fiber Gr\"{o}bner
cone $C_{I_{0}}\left(  I\right)  $ we obtain the special fiber Bergman
subcomplex
\[
B\left(  I\right)  =C_{I_{0}}\left(  I\right)  \cap BF\left(  I\right)
\cap\left\{  w_{t}=1\right\}
\]
The following table chooses an indexing of the vertices of $\nabla$ involved
in $B\left(  I\right)  $ and gives for each vertex the numbers $n_{\nabla}$
and $n_{B}$ of faces of $\nabla$ and faces of $B\left(  I\right)  $ it is
contained in

\begin{center}
%
\end{center}

\noindent With this indexing the Bergman subcomplex $B\left(  I\right)  $ of
$\operatorname*{Poset}\left(  \nabla\right)  $ associated to the degeneration
$\mathfrak{X}$ is

\begin{center}
%
\]
and the $F$-vectors of its faces are

\begin{center}
%
\end{center}

\subsubsection{Dual complex}

The dual complex $\operatorname*{dual}\left(  B\left(  I\right)  \right)
=\left(  B\left(  I\right)  \right)  ^{\ast}$ of deformations associated to
$B\left(  I\right)  $ via initial ideals is given by

\begin{center}
%
\end{center}

\noindent when writing the vertices of the faces as deformations of $X_{0}$.
Note that the $T$-invariant basis of deformations associated to a face is
given by all lattice points of the corresponding polytope in $M_{\mathbb{R}}$.

When numbering the vertices of the faces of $\operatorname*{dual}\left(
B\left(  I\right)  \right)  $ by the Cox variables of the mirror toric Fano
variety $Y^{\circ}$ the complex $\operatorname*{dual}\left(  B\left(
I\right)  \right)  $ is

\begin{center}
%
\end{center}

\noindent The dual complex has the $F$-vector%
\[
%
\]
and the $F$-vectors of the faces of $\operatorname*{dual}\left(  B\left(
I\right)  \right)  $ are

\begin{center}
%
\end{center}

Recall that in this example the toric variety $Y$ is projective space. The
number of lattice points of the support of $\operatorname*{dual}\left(
B\left(  I\right)  \right)  $ relates to the dimension $h^{1,
0
}%
\left(  X\right)  $ of the complex moduli space of the generic fiber $X$ of
$\mathfrak{X}$ and to the dimension $h^{1,1}\left(  \bar{X}^{\circ}\right)  $
of the K\"{a}hler moduli space of the $MPCR$-blowup $\bar{X}^{\circ}$ of the
generic fiber $X^{\circ}$ of the mirror degeneration
\begin{align*}
\left\vert \operatorname*{supp}\left(  \operatorname*{dual}\left(  B\left(
I\right)  \right)  \right)  \cap M\right\vert  &
=
25
=
24
+
1
=\dim\left(  \operatorname*{Aut}%
\left(  Y\right)  \right)  +h^{1,
0
}\left(  X\right) \\
&  =
20
+
4
+
1
\\
&  =\left\vert \operatorname*{Roots}\left(  Y\right)  \right\vert +\dim\left(
T_{Y}\right)  +h^{1,1}\left(  \bar{X}^{\circ}\right)
\end{align*}
There are%
\[
h^{1,
0
}\left(  X\right)  +\dim\left(  T_{Y^{\circ}}\right)
=
1
+
4
\]
non-trivial toric polynomial deformations of $X_{0}$

\begin{center}

\end{center}

\noindent They correspond to the toric divisors

\begin{center}
%
\end{center}

\noindent on a MPCR-blowup of $Y^{\circ}$ inducing $
1
$ non-zero
toric divisor classes on the mirror $X^{\circ}$. The following $
20
$
toric divisors of $Y$ induce the trivial divisor class on $X^{\circ}$

\begin{center}
%
\end{center}

\subsubsection{Mirror special fiber}

The ideal $I_{0}^{\circ}$ of the monomial special fiber $X_{0}^{\circ}$ of the
mirror degeneration $\mathfrak{X}^{\circ}$ is generated by the following set
of monomials in $S^{\circ}$%
\[
\left\{  
%
\right\}
\]
Indeed already the ideal%
\[
J_{0}^{\circ}=\left\langle 
%
\right\rangle
\]
defines the same subvariety of the toric variety $Y^{\circ}$, and
$J_{0}^{\circ\Sigma}=I_{0}^{\circ}$. Recall that passing from $J_{0}^{\circ}$
to $J_{0}^{\circ\Sigma}$ is the non-simplicial toric analogue of saturation.

The complex $B\left(  I\right)  ^{\ast}$ labeled by the variables of the Cox
ring $S^{\circ}$ of $Y^{\ast}$, as written in the last section, is the complex
$SP\left(  I_{0}^{\circ}\right)  $ of prime ideals of the toric strata of the
special fiber $X_{0}^{\circ}$ of the mirror degeneration $\mathfrak{X}^{\circ
}$, i.e., the primary decomposition of $I_{0}^{\circ}$ is%
\[
I_{0}^{\circ}=

\]
Each facet $F\in B\left(  I\right)  $ corresponds to one of these ideals and this
ideal is generated by the products of facets of $\nabla$ containing $F$.

\subsubsection{Covering structure in the deformation complex of the
degeneration $\mathfrak{X}$}

According to the local reduced Gr\"{o}bner basis each face of the complex of
deformations $\operatorname*{dual}\left(  B\left(  I\right)  \right)  $
decomposes into 
$3$
polytopes forming a $
3
:1$
unramified
covering of $B\left(  I\right)  $

\begin{center}
%
\end{center}

\noindent Here the faces $F\in B\left(  I\right)  $ are specified both via the
vertices of $F^{\ast}$ labeled by the variables of $S^{\circ}$ and by the
numbering of the vertices of $B\left(  I\right)  $ chosen above.

The numbers of faces of the covering in each face of $\operatorname*{dual}%
\left(  B\left(  I\right)  \right)  $, i.e. over each face of $B\left(
I\right)  ^{\vee}$ are

\begin{center}
%
\end{center}

\noindent This covering has 
one sheet forming the complex
 
\[
%
\end{center}

\noindent Writing the vertices of the faces as deformations the covering is
given by

\begin{center}
%
\end{center}

\noindent with 
one sheet forming the complex

\begin{center}
%
\end{center}

\noindent Note that the torus invariant basis of deformations corresponding to
a Bergman face is given by the set of all lattice points of the polytope
specified above.
See also Figure \ref{Fig dual Pfaff ell} for a visualization of the dual
complex and the sheets of the covering.

In general we have for local complete intersections:

\begin{remark}
Let $\mathfrak{X}\subset Y\times\operatorname*{Spec}\mathbb{C}\left[
t\right]  $, $Y=X\left(  \Sigma\right)  =\mathbb{P}^{n}$ be a degeneration
satisfying the conditions for the tropical mirror construction with fibers of
codimension $c$. Let the special fiber be given by the monomial ideal
$I_{0}=I_{0}^{\Sigma}\subset S$ with minimal generators $m_{1},...,m_{r}$ and
the associated first order degeneration $\mathfrak{X}\subset Y\times
\operatorname*{Spec}\mathbb{C}\left[  t\right]  //\left\langle t^{2}%
\right\rangle $ be given by the ideal $I^{1}=\left\langle f_{i}=m_{i}%
+tg_{i}\right\rangle \subset S\otimes\mathbb{C}\left[  t\right]  /\left\langle
t^{2}\right\rangle $. If the total space $\mathfrak{X}$ is a local complete
intersection, then all first order deformations of $X_{0}$ contribute
precisely once in the local equations of $\mathfrak{X}$ at the strata of
$X_{0}$:

Suppose $F\subset SP\left(  I_{0}\right)  $ is the prime ideal of a stratum of
$X_{0}$ and%
\[
I_{F}^{1}\subset S_{F}\otimes\mathbb{C}\left[  t\right]  /\left\langle
t^{2}\right\rangle =\mathbb{C}\left(  y_{r}\mid y_{r}\notin F,\text{ }%
r\in\Sigma\left(  1\right)  \right)  \left[  y_{r}\mid y_{r}\in F,\text{ }%
r\in\Sigma\left(  1\right)  \right]  _{>}\otimes\mathbb{C}\left[  t\right]
/\left\langle t^{2}\right\rangle
\]
where $>$ is a local ordering on $y_{r}\in F,$ $r\in\Sigma\left(  1\right)  $,
is the localization of $I$ at $F$. Then for all deformations $\delta
\in\operatorname*{dual}\left(  B\right)  \cap M$ any monomial $\delta\left(
m_{i}\right)  $ appears at most once in the minimal reduced standard basis of
$I_{F}$.

The complex $\operatorname*{dual}\left(  B\right)  $ contains a $c:1$
unramified covering of $B^{\vee}$.
\end{remark}

\subsubsection{Limit map}

The limit map $\lim:B\left(  I\right)  \rightarrow\operatorname*{Poset}\left(
\Delta\right)  $ associates to a face $F$ of $B\left(  I\right)  $ the face of
$\Delta$ formed by the limit points of arcs lying over the weight vectors
$w\in F$, i.e., with lowest order term $t^{w}$.

Labeling the faces of the Bergman complex $B\left(  I\right)  \subset
\operatorname*{Poset}\left(  \nabla\right)  $ and the faces of
$\operatorname*{Poset}\left(  \Delta\right)  $ by the corresponding dual faces
of $\nabla^{\ast}$ and $\Delta^{\ast}$, hence considering the limit map
$\lim:B\left(  I\right)  \rightarrow\operatorname*{Poset}\left(
\Delta\right)  $ as a map $B\left(  I\right)  ^{\ast}\rightarrow
\operatorname*{Poset}\left(  \Delta^{\ast}\right)  $, the limit correspondence
is given by

\begin{center}

\end{center}

The image of the limit map coincides with the image of $\mu$ and with the
Bergman complex of the mirror, i.e. $\lim\left(  B\left(  I\right)  \right)
=\mu\left(  B\left(  I\right)  \right)  =B\left(  I^{\ast}\right)  $.

\subsubsection{Mirror complex}

Numbering the vertices of the mirror complex $\mu\left(  B\left(  I\right)
\right)  $ as%
\[
%
\end{center}

\noindent The ordering of the faces of $\mu\left(  B\left(  I\right)  \right)
$ is compatible with above ordering of the faces of $B\left(  I\right)  $. The
$F$-vectors of the faces of $\mu\left(  B\left(  I\right)  \right)  $ are%
\[
%
\]
The first order deformations of the mirror special fiber $X_{0}^{\circ}$
correspond to the lattice points of the dual complex $\left(  \mu\left(
B\left(  I\right)  \right)  \right)  ^{\ast}\subset\operatorname*{Poset}%
\left(  \Delta^{\ast}\right)  $ of the mirror complex. We label the first
order deformations of $X_{0}^{\circ}$ corresponding to vertices of
$\Delta^{\ast}$ by the homogeneous coordinates of $Y$%
\[
%
\]
So writing the vertices of the faces of $\left(  \mu\left(  B\left(  I\right)
\right)  \right)  ^{\ast}$ as homogeneous coordinates of $Y$, the complex
$\left(  \mu\left(  B\left(  I\right)  \right)  \right)  ^{\ast}$ is given by

\begin{center}
%
\end{center}

\noindent The complex $\mu\left(  B\left(  I\right)  \right)  ^{\ast}$ labeled
by the variables of the Cox ring $S$ gives the ideals of the toric strata of
the special fiber $X_{0}$ of $\mathfrak{X}$, i.e., the complex $SP\left(
I_{0}\right)  $, so in particular the primary decomposition of $I_{0}$ is%
\[
I_{0}=
%
\]

\subsubsection{Covering structure in the deformation complex of the mirror
degeneration}

Each face of the complex of deformations $\left(  \mu\left(  B\left(
I\right)  \right)  \right)  ^{\ast}$ of the mirror special fiber $X_{0}%
^{\circ}$ decomposes as the convex hull of 
$3$ respectively $4$
polytopes forming
a $
4
:1$ 
ramified
covering of $\mu\left(  B\left(
I\right)  \right)  ^{\vee}$

\begin{center}
%
\end{center}

\noindent 
Due to the singularities of
$Y^{\circ}$ this covering involves degenerate faces, i.e. faces $G\mapsto
F^{\vee}$ with $\dim\left(  G\right)  <\dim\left(  F^{\vee}\right)  $.

\subsubsection{Mirror degeneration}

The space of first order deformations of $X_{0}^{\circ}$ in the mirror
degeneration $\mathfrak{X}^{\circ}$ has dimension $
5
$ and the
deformations represented by the monomials%
\[
\left\{  
%
\right\}
\]
form a torus invariant
basis.
The number of lattice points of the dual of the mirror complex of $I$ relates
to the dimension $h^{1,
0
}\left(  X^{\circ}\right)  $ of complex
moduli space of the generic fiber $X^{\circ}$ of $\mathfrak{X}^{\circ}$ and to
the dimension $h^{1,1}\left(  X\right)  $ of the K\"{a}hler moduli space of
the generic fiber $X$ of $\mathfrak{X}$ via%
\begin{align*}
\left\vert \operatorname*{supp}\left(  \left(  \mu\left(  B\left(  I\right)
\right)  \right)  ^{\ast}\right)  \cap N\right\vert  &
=
5
=
4
+
1
\\
&  =\dim\left(  \operatorname*{Aut}\left(  Y^{\circ}\right)  \right)
+h^{1,
0
}\left(  X^{\circ}\right)  =\dim\left(  T\right)
+h^{1,1}\left(  X\right)
\end{align*}

The 
conjectural first order
mirror degeneration 
$\mathfrak{X}^{1\circ}\subset Y^{\circ}\times\operatorname*{Spec}%
\mathbb{C}\left[  t\right]  /\left\langle t^{2}\right\rangle $
 of $\mathfrak{X}$ is
given by the ideal 
$I^{1\circ}\subset S^{\circ}\otimes\mathbb{C}\left[  t\right]  /\left\langle
t^{2}\right\rangle $
generated by%
\[

\]
Indeed already the ideal $J^{\circ}$ generated by%
\[
\left\{  
%
\right\}
\]
defines $\mathfrak{X}^{1\circ}$.

\subsubsection{Contraction of the mirror degeneration}

In the following we give a birational map relating the degeneration
$\mathfrak{X}^{\circ}$ to a Greene-Plesser type orbifolding mirror family by
contracting divisors on $Y^{\circ}$.

In order to contract the divisors

\begin{center}
%
\end{center}

\noindent consider the $\mathbb{Q}$-factorial toric Fano variety $\hat
{Y}^{\circ}=X\left(  \hat{\Sigma}^{\circ}\right)  $, where $\hat{\Sigma
}^{\circ}=\operatorname*{Fan}\left(  \hat{P}^{\circ}\right)  \subset
M_{\mathbb{R}}$ is the fan over the 
reflexive
Fano polytope $\hat
{P}^{\circ}\subset M_{\mathbb{R}}$ given as the convex hull of the remaining
vertices of $P^{\circ}=\nabla^{\ast}$ corresponding to the Cox variables%
\[
%
\]
of the toric variety $\hat{Y}^{\circ}$ with Cox ring
\[
\hat{S}^{\circ}=\mathbb{C}
[y_{4},y_{1},y_{5},y_{2},y_{3}]
\]
The Cox variables of $\hat{Y}^{\circ}$ correspond to the set of Fermat deformations of $\mathfrak{X}$.

Let%
\[
Y^{\circ}=X\left(  \Sigma^{\circ}\right)  \rightarrow X\left(  \hat{\Sigma
}^{\circ}\right)  =\hat{Y}^{\circ}%
\]
be a birational map from $Y^{\circ}$ to a minimal birational model $\hat
{Y}^{\circ}$, which contracts the divisors of the rays $\Sigma^{\circ}\left(
1\right)  -\hat{\Sigma}^{\circ}\left(  1\right)  $ corresponding to the Cox
variables
\[
\begin{tabular}
[c]{llllllllll}
$y_{6}$ &$y_{7}$ &$y_{8}$ &$y_{9}$ &$y_{10}$ &$y_{11}$ &$y_{12}$ &$y_{13}$ &$y_{14}$ &$y_{15}$
\end{tabular}
\]

Representing $\hat{Y}^{\circ}$ as a quotient we have
\[
\hat{Y}^{\circ}=\left(  \mathbb{C}^{
5
}-V\left(  B\left(
\hat{\Sigma}^{\circ}\right)  \right)  \right)  //\hat{G}^{\circ}%
\]
with%
\[
\hat{G}^{\circ}=
\mathbb{Z}_{5}\times\left(  \mathbb{C}^{\ast}\right)  ^{1}
\]
acting via%
\[
\xi y=
\left( \,u_{1}^{4}\,v_{1} \cdot y_{4},\,u_{1}\,v_{1} \cdot y_{1},\,v_{1} \cdot y_{5},\,u_{1}^{3}\,v_{1} \cdot y_{2},\,u_{1}^{2}\,v_{1} \cdot y_{3} \right)
\]
for $\xi=
\left(u_1,v_1\right)
\in\hat{G}^{\circ}$ and $y\in\mathbb{C}%
^{
5
}-V\left(  B\left(  \hat{\Sigma}^{\circ}\right)  \right)  $.

Hence with the group%
\[
\hat{H}^{\circ}=
\mathbb{Z}_{5}
\]
of order 
5
 the toric variety $\hat{Y}^{\circ}$ is the quotient%
\[
\hat{Y}^{\circ}=\mathbb{P}^{
4
}/\hat{H}^{\circ}%
\]
of projective space $\mathbb{P}^{
4
}$.

The first order mirror degeneration $\mathfrak{X}^{1\circ}$ induces via
$Y\rightarrow\hat{Y}^{\circ}$ a degeneration $\mathfrak{\hat{X}}^{1\circ
}\subset\hat{Y}^{\circ}\times\operatorname*{Spec}\mathbb{C}\left[  t\right]
/\left\langle t^{2}\right\rangle $ given by the ideal $\hat{I}^{1\circ}%
\subset\hat{S}^{\circ}\otimes\mathbb{C}\left[  t\right]  /\left\langle
t^{2}\right\rangle $ generated by the
 Fermat-type equations
\[
\left\{  

\right\}
\]

Note, that%
\[
\left\vert \operatorname*{supp}\left(  \left(  \mu\left(  B\left(  I\right)
\right)  \right)  ^{\ast}\right)  \cap N-\operatorname*{Roots}\left(  \hat
{Y}^{\circ}\right)  \right\vert -\dim\left(  T_{\hat{Y}^{\circ}}\right)
=
5
-
4
=
1
\]
so this family has one independent parameter.

The special fiber $\hat{X}_{0}^{\circ}\subset\hat{Y}^{\circ}$ of
$\mathfrak{\hat{X}}^{1\circ}$ is cut out by the monomial ideal $\hat{I}%
_{0}^{\circ}\subset\hat{S}^{\circ}$ generated by%
\[
\left\{  
%
\right\}
\]

The complex $SP\left(  \hat{I}_{0}^{\circ}\right)  $ of prime ideals of the
toric strata of $\hat{X}_{0}^{\circ}$ is

\begin{center}
%
\end{center}

\noindent so $\hat{I}_{0}^{\circ}$ has the primary decomposition%
\[
\hat{I}_{0}^{\circ}=
%
\]
The Bergman subcomplex $B\left(  I\right)  \subset\operatorname*{Poset}\left(
\nabla\right)  $ induces a subcomplex of $\operatorname*{Poset}\left(
\hat{\nabla}\right)  $, $\hat{\nabla}=\left(  \hat{P}^{\circ}\right)  ^{\ast}$
corresponding to $\hat{I}_{0}^{\circ}$. Indexing of the 
vertices
of $\hat{\nabla}$ by%
\[
%
\]
this complex is given by
\[
\left\{  
%
\right\}
\]

The ideal $\hat{I}^{1\circ}\subset\hat{S}^{\circ}\otimes\mathbb{C}\left[
t\right]  /\left\langle t^{2}\right\rangle $ has a Pfaffian resolution%
\begin{gather*}
0\rightarrow\mathcal{O}_{\hat{Y}^{\circ}\times\operatorname*{Spec}%
\mathbb{C}\left[  t\right]  /\left\langle t^{2}\right\rangle }\left(
K^{1}\right)  \rightarrow\mathcal{E}^{1}\left(  K^{1}\right)  \overset
{\varphi^{1}}{\rightarrow}\left(  \mathcal{E}^{1}\right)  ^{\ast}%
\overset{f^{1}}{\rightarrow}\mathcal{O}_{\hat{Y}^{\circ}\times
\operatorname*{Spec}\mathbb{C}\left[  t\right]  /\left\langle t^{2}%
\right\rangle }\medskip\\
\text{where }\overline{\pi}_{1}:\hat{Y}^{\circ}\times\operatorname*{Spec}%
\mathbb{C}\left[  t\right]  /\left\langle t^{2}\right\rangle \rightarrow
\hat{Y}^{\circ}\text{ and }\mathcal{E}^{1}=\overline{\pi}_{1}^{\ast
}\mathcal{F}%
\end{gather*}
with%
\[
\mathcal{F}=
\begin{tabular}
[c]{l}
$\mathcal{O}_{\hat{Y}^{\circ}}\left(D_{ \left(2,0,-1,-1\right) }+D_{ \left(0,-1,-1,0\right) } \right) \oplus \mathcal{O}_{\hat{Y}^{\circ}}\left(D_{ \left(-1,0,2,0\right) }+D_{ \left(-1,-1,0,2\right) } \right) \oplus $
\\
$\mathcal{O}_{\hat{Y}^{\circ}}\left(D_{ \left(2,0,-1,-1\right) }+D_{ \left(0,2,0,-1\right) } \right) \oplus \mathcal{O}_{\hat{Y}^{\circ}}\left(D_{ \left(-1,0,2,0\right) }+D_{ \left(0,2,0,-1\right) } \right) \oplus $
\\
$\mathcal{O}_{\hat{Y}^{\circ}}\left(D_{ \left(-1,-1,0,2\right) }+D_{ \left(0,-1,-1,0\right) } \right)$
\\
\end{tabular}
\]
and $K^{1}=K_{\hat{Y}^{\circ}\times\operatorname*{Spec}\mathbb{C}\left[
t\right]  /\left\langle t^{2}\right\rangle /\operatorname*{Spec}%
\mathbb{C}\left[  t\right]  /\left\langle t^{2}\right\rangle }$ and
$\varphi^{1}\in\bigwedge\nolimits^{2}\mathcal{E}^{1}\left(  -K^{1}\right)  $
given by%
\[
\left [\begin {array}{ccccc} 0&ts_{3}\,y_{5}&y_{3}&-y_{4}&-ts_{2}\,y_{1}\\\noalign{\medskip}-ts_{3}\,y_{5}&0&ts_{4}\,y_{4}&y_{1}&-y_{2}\\\noalign{\medskip}-y_{3}&-ts_{4}\,y_{4}&0&ts_{5}\,y_{2}&y_{5}\\\noalign{\medskip}y_{4}&-y_{1}&-ts_{5}\,y_{2}&0&ts_{1}\,y_{3}\\\noalign{\medskip}ts_{2}\,y_{1}&y_{2}&-y_{5}&-ts_{1}\,y_{3}&0\end {array}\right ]
\]
Hence via the Pfaffians of $\varphi^{1}$ we obtain a resolution%
\begin{gather*}
0\rightarrow\mathcal{O}_{\hat{Y}^{\circ}\times\operatorname*{Spec}%
\mathbb{C}\left[  t\right]  }\left(  K\right)  \rightarrow\mathcal{E}\left(
K\right)  \rightarrow\mathcal{E}^{\ast}\rightarrow\mathcal{O}_{\hat{Y}^{\circ
}\times\operatorname*{Spec}\mathbb{C}\left[  t\right]  }\medskip\\
\text{where }\pi_{1}:Y\times\operatorname*{Spec}\mathbb{C}\left[  t\right]
\rightarrow Y\text{ , }\mathcal{E}=\pi_{1}^{\ast}\mathcal{F}\\
\text{and }K=K_{\hat{Y}^{\circ}\times\operatorname*{Spec}\mathbb{C}\left[
t\right]  /\operatorname*{Spec}\mathbb{C}\left[  t\right]  }%
\end{gather*}
of the ideal $\hat{I}^{\circ}\subset\hat{S}^{\circ}\otimes\mathbb{C}\left[
t\right]  $ generated by%
\[
\left\{  

\right\}
\]
which defines a flat family%
\[
\mathfrak{\hat{X}}^{\circ}\subset\hat{Y}^{\circ}\times\operatorname*{Spec}%
\mathbb{C}\left[  t\right]
\]

\subsection{Tropical mirror construction for the degree $14$ Pfaffian
Calabi-Yau threefold\label{Sec ex pfaffian degree 14}}

\subsubsection{Setup}

Let $Y=\mathbb{P}^{
6
}=X\left(  \Sigma\right)  $, $\Sigma
=\operatorname*{Fan}\left(  P\right)  =NF\left(  \Delta\right)  \subset
N_{\mathbb{R}}$ with the Fano polytope $P=\Delta^{\ast}$ given by%
\[
\Delta=\operatorname*{convexhull}\left(  
%
\right)  \subset
M_{\mathbb{R}}%
\]
and let%
\[
S=\mathbb{C}
[x_0, x_1, x_2, x_3, x_4, x_5, x_6]
\]
be the Cox ring of $Y$ with the variables%
\[
%
\]
associated to the rays of $\Sigma$.

Consider the degeneration $\mathfrak{X}\subset Y\times\operatorname*{Spec}%
\mathbb{C}\left[  t\right]  $ of 
Pfaffian
Calabi-Yau 3-folds
 with Buchsbaum-Eisenbud
resolution%
\[
0\rightarrow\mathcal{O}_{Y}\left(  -
7
\right)  \rightarrow
\mathcal{E}\left(  -
4
\right)  \overset{A_{t}}{\rightarrow
}\mathcal{E}^{\ast}\left(  -
3
\right)  \rightarrow\mathcal{O}%
_{Y}\rightarrow\mathcal{O}_{X_{t}}\rightarrow0
\]
where%
\[
\mathcal{E}=\text{
$7\mathcal{O}$
}%
\]%
\[
A_{t}=A_{0}+t\cdot A
\]%
\[
A_{0}=
\left [%
\right ]
\]
the monomial special fiber of $\mathfrak{X}$ is given by%
\[
I_{0}=\left\langle 
%
\]
The degeneration%
\[
\mathfrak{X}\subset Y\times\operatorname*{Spec}\mathbb{C}\left[  t\right]
\]
is given by the ideal $I\subset S\otimes\mathbb{C}\left[  t\right]  $
generated by the Pfaffians of $A_{0}+t\cdot A$ of degrees $
3,3,3,3,3,3,3
$.

\subsubsection{Special fiber Gr\"{o}bner cone}

The space of first order deformations of $\mathfrak{X}$ has dimension
$
98
$ and the deformations represented by the Cox Laurent monomials

\begin{center}
%
\end{center}

\noindent form a torus invariant basis. These deformations give linear
inequalities defining the special fiber Gr\"{o}bner cone, i.e.,%
\[
C_{I_{0}}\left(  I\right)  =\left\{  \left(  w_{t},w_{y}\right)  \in
\mathbb{R}\oplus N_{\mathbb{R}}\mid\left\langle w,A^{-1}\left(  m\right)
\right\rangle \geq-w_{t}\forall m\right\}
\]
for above monomials $m$ and the presentation matrix
\[
A=
\left [%
\right ]
\]
of $A_{
5
}\left(  Y\right)  $.

The vertices of special fiber polytope
\[
\nabla=C_{I_{0}}\left(  I\right)  \cap\left\{  w_{t}=1\right\}
\]
and the number of faces, each vertex is contained in, are given by the table

\begin{center}
%
\end{center}

\noindent The number of faces of $\nabla$ and their $F$-vectors are

\begin{center}
%
\end{center}

The dual $P^{\circ}=\nabla^{\ast}$ of $\nabla$ is a Fano polytope with vertices

\begin{center}
%
\end{center}

\noindent and the $F$-vectors of the faces of $P^{\circ}$ are

\begin{center}
%
\end{center}

The polytopes $\nabla$ and $P^{\circ}$ have $0$ as the unique interior lattice
point, and the Fano polytope $P^{\circ}$ defines the embedding toric Fano
variety $Y^{\circ}=X\left(  \operatorname*{Fan}\left(  P^{\circ}\right)
\right)  $ of the fibers of the mirror degeneration. The dimension of the
automorphism group of $Y^{\circ}$ is
\[
\dim\left(  \operatorname*{Aut}\left(  Y^{\circ}\right)  \right)
=
6
\]
Let%
\[
S^{\circ}=\mathbb{C}
[y_{1},y_{2},...,y_{70}]
\]
be the Cox ring of $Y^{\circ}$ with variables

\begin{center}
%
\end{center}

\subsubsection{Bergman subcomplex}

Intersecting the tropical variety of $I$ with the special fiber Gr\"{o}bner
cone $C_{I_{0}}\left(  I\right)  $ we obtain the special fiber Bergman
subcomplex
\[
B\left(  I\right)  =C_{I_{0}}\left(  I\right)  \cap BF\left(  I\right)
\cap\left\{  w_{t}=1\right\}
\]
The following table chooses an indexing of the vertices of $\nabla$ involved
in $B\left(  I\right)  $ and gives for each vertex the numbers $n_{\nabla}$
and $n_{B}$ of faces of $\nabla$ and faces of $B\left(  I\right)  $ it is
contained in

\begin{center}
%
\end{center}

\noindent With this indexing the Bergman subcomplex $B\left(  I\right)  $ of
$\operatorname*{Poset}\left(  \nabla\right)  $ associated to the degeneration
$\mathfrak{X}$ is

\begin{center}
%
\]
and the $F$-vectors of its faces are

\begin{center}
%
\end{center}

\subsubsection{Dual complex}

The dual complex $\operatorname*{dual}\left(  B\left(  I\right)  \right)
=\left(  B\left(  I\right)  \right)  ^{\ast}$ of deformations associated to
$B\left(  I\right)  $ via initial ideals is given by

\begin{center}
%
\end{center}

\noindent when writing the vertices of the faces as deformations of $X_{0}$.
Note that the $T$-invariant basis of deformations associated to a face is
given by all lattice points of the corresponding polytope in $M_{\mathbb{R}}$.

In order to compress the output we list one representative in any set of faces
$G$ with fixed $F$-vector of $G$ and $G^{\ast}$.

When numbering the vertices of the faces of $\operatorname*{dual}\left(
B\left(  I\right)  \right)  $ by the Cox variables of the mirror toric Fano
variety $Y^{\circ}$ the complex $\operatorname*{dual}\left(  B\left(
I\right)  \right)  $ is

\begin{center}

\end{center}

\noindent The dual complex has the $F$-vector%
\[
%
\]
and the $F$-vectors of the faces of $\operatorname*{dual}\left(  B\left(
I\right)  \right)  $ are

\begin{center}
%
\end{center}

Recall that in this example the toric variety $Y$ is projective space. The
number of lattice points of the support of $\operatorname*{dual}\left(
B\left(  I\right)  \right)  $ relates to the dimension $h^{1,
2
}%
\left(  X\right)  $ of the complex moduli space of the generic fiber $X$ of
$\mathfrak{X}$ and to the dimension $h^{1,1}\left(  \bar{X}^{\circ}\right)  $
of the K\"{a}hler moduli space of the $MPCR$-blowup $\bar{X}^{\circ}$ of the
generic fiber $X^{\circ}$ of the mirror degeneration
\begin{align*}
\left\vert \operatorname*{supp}\left(  \operatorname*{dual}\left(  B\left(
I\right)  \right)  \right)  \cap M\right\vert  &
=
98
=
48
+
50
=\dim\left(  \operatorname*{Aut}%
\left(  Y\right)  \right)  +h^{1,
2
}\left(  X\right) \\
&  =
42
+
6
+
50
\\
&  =\left\vert \operatorname*{Roots}\left(  Y\right)  \right\vert +\dim\left(
T_{Y}\right)  +h^{1,1}\left(  \bar{X}^{\circ}\right)
\end{align*}
There are%
\[
h^{1,
2
}\left(  X\right)  +\dim\left(  T_{Y^{\circ}}\right)
=
50
+
6
\]
non-trivial toric polynomial deformations of $X_{0}$

\begin{center}

\end{center}

\noindent They correspond to the toric divisors

\begin{center}
%
\end{center}

\noindent on a MPCR-blowup of $Y^{\circ}$ inducing $
50
$ non-zero
toric divisor classes on the mirror $X^{\circ}$. The following $
42
$
toric divisors of $Y$ induce the trivial divisor class on $X^{\circ}$

\begin{center}
%
\end{center}

\subsubsection{Mirror special fiber}

The complex $B\left(  I\right)  ^{\ast}$ labeled by the variables of the Cox
ring $S^{\circ}$ of $Y^{\ast}$, as written in the last section, is the complex
$SP\left(  I_{0}^{\circ}\right)  $ of prime ideals of the toric strata of the
monomial special fiber $X_{0}^{\circ}$ of the mirror degeneration
$\mathfrak{X}^{\circ}$, i.e. the primary decomposition of $I_{0}^{\circ}$ is%
\[
I_{0}^{\circ}=
%
\]
Each facet $F\in B\left(  I\right)  $ corresponds to one of these ideals and this
ideal is generated by the products of facets of $\nabla$ containing $F$.

\subsubsection{Covering structure in the deformation complex of the
degeneration $\mathfrak{X}$}

According to the local reduced Gr\"{o}bner basis each face of the complex of
deformations $\operatorname*{dual}\left(  B\left(  I\right)  \right)  $
decomposes into 
$3$ respectively $5$
polytopes forming a $
5
:1$
ramified
covering of $B\left(  I\right)  $

\begin{center}
%
\end{center}

\noindent Here the faces $F\in B\left(  I\right)  $ are specified both via the
vertices of $F^{\ast}$ labeled by the variables of $S^{\circ}$ and by the
numbering of the vertices of $B\left(  I\right)  $ chosen above.

The numbers of faces of the covering in each face of $\operatorname*{dual}%
\left(  B\left(  I\right)  \right)  $, i.e. over each face of $B\left(
I\right)  ^{\vee}$ are

\begin{center}
%
\end{center}

\noindent This covering has 
one sheet forming the complex
 
\[
%
\end{center}

\noindent Writing the vertices of the faces as deformations the covering is
given by

\begin{center}
%
\end{center}

\noindent with 
one sheet forming the complex

\begin{center}
%
\end{center}

\noindent Note that the torus invariant basis of deformations corresponding to
a Bergman face is given by the set of all lattice points of the polytope
specified above.

\subsubsection{Limit map}

The limit map $\lim:B\left(  I\right)  \rightarrow\operatorname*{Poset}\left(
\Delta\right)  $ associates to a face $F$ of $B\left(  I\right)  $ the face of
$\Delta$ formed by the limit points of arcs lying over the weight vectors
$w\in F$, i.e. with lowest order term $t^{w}$.

Labeling the faces of the Bergman complex $B\left(  I\right)  \subset
\operatorname*{Poset}\left(  \nabla\right)  $ and the faces of
$\operatorname*{Poset}\left(  \Delta\right)  $ by the corresponding dual faces
of $\nabla^{\ast}$ and $\Delta^{\ast}$, hence considering the limit map
$\lim:B\left(  I\right)  \rightarrow\operatorname*{Poset}\left(
\Delta\right)  $ as a map $B\left(  I\right)  ^{\ast}\rightarrow
\operatorname*{Poset}\left(  \Delta^{\ast}\right)  $, the limit correspondence
is given by

\begin{center}

\end{center}

The image of the limit map coincides with the image of $\mu$, i.e.
$\lim\left(  B\left(  I\right)  \right)  =\mu\left(  B\left(  I\right)
\right)  $.

\begin{figure}
[h]
\begin{center}
\includegraphics[
height=2.4509in,
width=1.695in
]%
{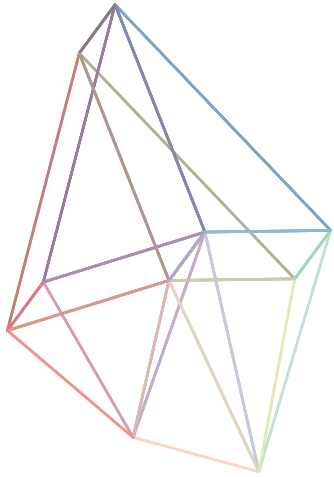}%
\caption{Projection of Bergman faces with same limit}%
\label{Fig dualcyclic}%
\end{center}
\end{figure}

Every zero dimensional stratum of $\mathbb{P}^{6}$ is the limit of $5$ Bergman
faces of dimension three, one tetrahedron, two pyramids and two prisms. Figure
\ref{Fig dualcyclic} shows a projection into $3$-space of the set of these
Bergman faces for one zero dimensional stratum. The union of the faces $F\in
B\left(  I\right)  $ which have as limit $\lim\left(  F\right)  =p$ the same
zero dimensional strata $p$ of $\mathbb{P}^{6}$ is not convex.%

\subsubsection{Mirror complex}

Numbering the vertices of the mirror complex $\mu\left(  B\left(  I\right)
\right)  $ as%
\[

\end{center}

\noindent The ordering of the faces of $\mu\left(  B\left(  I\right)  \right)
$ is compatible with above ordering of the faces of $B\left(  I\right)  $. The
$F$-vectors of the faces of $\mu\left(  B\left(  I\right)  \right)  $ are%
\[
%
\]
The first order deformations of the mirror special fiber $X_{0}^{\circ}$
correspond to the lattice points of the dual complex $\left(  \mu\left(
B\left(  I\right)  \right)  \right)  ^{\ast}\subset\operatorname*{Poset}%
\left(  \Delta^{\ast}\right)  $ of the mirror complex. We label the first
order deformations of $X_{0}^{\circ}$ corresponding to vertices of
$\Delta^{\ast}$ by the homogeneous coordinates of $Y$%
\[
%
\]
So writing the vertices of the faces of $\left(  \mu\left(  B\left(  I\right)
\right)  \right)  ^{\ast}$ as homogeneous coordinates of $Y$, the complex
$\left(  \mu\left(  B\left(  I\right)  \right)  \right)  ^{\ast}$ is given by

\begin{center}
%
\end{center}

\noindent The complex $\mu\left(  B\left(  I\right)  \right)  ^{\ast}$ labeled
by the variables of the Cox ring $S$ gives the ideals of the toric strata of
the special fiber $X_{0}$ of $\mathfrak{X}$, i.e. the complex $SP\left(
I_{0}\right)  $, so in particular the primary decomposition of $I_{0}$ is%
\[
I_{0}=
%
\]

\subsubsection{Covering structure in the deformation complex of the mirror
degeneration}

Each face of the complex of deformations $\left(  \mu\left(  B\left(
I\right)  \right)  \right)  ^{\ast}$ of the mirror special fiber $X_{0}%
^{\circ}$ decomposes as the convex hull of 
$3$,$4$,$5$,$6$ respectively $7$
polytopes forming
a $
7
:1$ 
ramified
covering of $\mu\left(  B\left(
I\right)  \right)  ^{\vee}$

\begin{center}
%
\end{center}

\noindent 
Due to the singularities of
$Y^{\circ}$ this covering involves degenerate faces, i.e. faces $G\mapsto
F^{\vee}$ with $\dim\left(  G\right)  <\dim\left(  F^{\vee}\right)  $.

\subsubsection{Mirror degeneration}

The space of first order deformations of $X_{0}^{\circ}$ in the mirror
degeneration $\mathfrak{X}^{\circ}$ has dimension $
7
$ and the
deformations represented by the monomials%
\[
\left\{  
%
\right\}
\]
form a torus invariant
basis $\mathfrak{B}^{\circ}$.
The number of lattice points of the dual of the mirror complex of $I$ relates
to the dimension $h^{1,
2
}\left(  X^{\circ}\right)  $ of complex
moduli space of the generic fiber $X^{\circ}$ of $\mathfrak{X}^{\circ}$ and to
the dimension $h^{1,1}\left(  X\right)  $ of the K\"{a}hler moduli space of
the generic fiber $X$ of $\mathfrak{X}$ via%
\begin{align*}
\left\vert \operatorname*{supp}\left(  \left(  \mu\left(  B\left(  I\right)
\right)  \right)  ^{\ast}\right)  \cap N\right\vert  &
=
7
=
6
+
1
\\
&  =\dim\left(  \operatorname*{Aut}\left(  Y^{\circ}\right)  \right)
+h^{1,
2
}\left(  X^{\circ}\right)  =\dim\left(  T\right)
+h^{1,1}\left(  X\right)
\end{align*}

The 
conjectural first order
mirror degeneration 
$\mathfrak{X}^{1\circ}\subset Y^{\circ}\times\operatorname*{Spec}%
\mathbb{C}\left[  t\right]  /\left\langle t^{2}\right\rangle $
 of $\mathfrak{X}$ is
given by the ideal 
$I^{1\circ}\subset S^{\circ}\otimes\mathbb{C}\left[  t\right]  /\left\langle
t^{2}\right\rangle $
generated by%
\[
\left\{  m+\sum_{\delta\in\mathfrak{B}^{\circ}}t\cdot s_{\delta}\cdot
\delta\left(  m\right)  \mid m\in I_{0}^{\circ}\right\}
\]

\subsubsection{Contraction of the mirror degeneration}

In the following we give a birational map relating the degeneration
$\mathfrak{X}^{\circ}$ to a Greene-Plesser type orbifolding mirror family by
contracting divisors on $Y^{\circ}$.

In order to contract the divisors

\begin{center}

\end{center}

\noindent consider the $\mathbb{Q}$-factorial toric Fano variety $\hat
{Y}^{\circ}=X\left(  \hat{\Sigma}^{\circ}\right)  $, where $\hat{\Sigma
}^{\circ}=\operatorname*{Fan}\left(  \hat{P}^{\circ}\right)  \subset
M_{\mathbb{R}}$ is the fan over the 
reflexive
Fano polytope $\hat
{P}^{\circ}\subset M_{\mathbb{R}}$ given as the convex hull of the remaining
vertices of $P^{\circ}=\nabla^{\ast}$ corresponding to the Cox variables%
\[
%
\]
of the toric variety $\hat{Y}^{\circ}$ with Cox ring
\[
\hat{S}^{\circ}=\mathbb{C}
[y_{22},y_{18},y_{30},y_{35},y_{16},y_{27},y_{38}]
\]
The Cox variables of $\hat{Y}^{\circ}$ correspond to the set of Fermat deformations of $\mathfrak{X}$.

Let%
\[
Y^{\circ}=X\left(  \Sigma^{\circ}\right)  \rightarrow X\left(  \hat{\Sigma
}^{\circ}\right)  =\hat{Y}^{\circ}%
\]
be a birational map from $Y^{\circ}$ to a minimal birational model $\hat
{Y}^{\circ}$, which contracts the divisors of the rays $\Sigma^{\circ}\left(
1\right)  -\hat{\Sigma}^{\circ}\left(  1\right)  $ corresponding to the Cox
variables
\[
\begin{tabular}
[c]{llllllllllllllll}
$y_{1}$ &$y_{2}$ &$y_{3}$ &$y_{4}$ &$y_{5}$ &$y_{6}$ &$y_{7}$ &$y_{8}$ &$y_{9}$ &$y_{10}$ &$y_{11}$ &$y_{12}$ &$y_{13}$ &$y_{14}$ &$y_{15}$ &$y_{17}$\\
$y_{19}$ &$y_{20}$ &$y_{21}$ &$y_{23}$ &$y_{24}$ &$y_{25}$ &$y_{26}$ &$y_{28}$ &$y_{29}$ &$y_{31}$ &$y_{32}$ &$y_{33}$ &$y_{34}$ &$y_{36}$ &$y_{37}$ &$y_{39}$\\
$y_{40}$ &$y_{41}$ &$y_{42}$ &$y_{43}$ &$y_{44}$ &$y_{45}$ &$y_{46}$ &$y_{47}$ &$y_{48}$ &$y_{49}$ &$y_{50}$ &$y_{51}$ &$y_{52}$ &$y_{53}$ &$y_{54}$ &$y_{55}$\\
$y_{56}$ &$y_{57}$ &$y_{58}$ &$y_{59}$ &$y_{60}$ &$y_{61}$ &$y_{62}$ &$y_{63}$ &$y_{64}$ &$y_{65}$ &$y_{66}$ &$y_{67}$ &$y_{68}$ &$y_{69}$ &$y_{70}$ &
\end{tabular}
\]

Representing $\hat{Y}^{\circ}$ as a quotient we have
\[
\hat{Y}^{\circ}=\left(  \mathbb{C}^{
7
}-V\left(  B\left(
\hat{\Sigma}^{\circ}\right)  \right)  \right)  //\hat{G}^{\circ}%
\]
with%
\[
\hat{G}^{\circ}=
\mathbb{Z}_{7}\times\left(  \mathbb{C}^{\ast}\right)  ^{1}
\]
acting via%
\[
\xi y=
\left( \,u_{1}^{5}\,v_{1} \cdot y_{22},\,u_{1}^{4}\,v_{1} \cdot y_{18},\,u_{1}^{6}\,v_{1} \cdot y_{30},\,u_{1}^{2}\,v_{1} \cdot y_{35},\,u_{1}^{3}\,v_{1} \cdot y_{16},\,u_{1}\,v_{1} \cdot y_{27},\,v_{1} \cdot y_{38} \right)
\]
for $\xi=
\left(u_1,v_1\right)
\in\hat{G}^{\circ}$ and $y\in\mathbb{C}%
^{
7
}-V\left(  B\left(  \hat{\Sigma}^{\circ}\right)  \right)  $.

Hence with the group%
\[
\hat{H}^{\circ}=
\mathbb{Z}_{7}
\]
of order 
7
 the toric variety $\hat{Y}^{\circ}$ is the quotient%
\[
\hat{Y}^{\circ}=\mathbb{P}^{
6
}/\hat{H}^{\circ}%
\]
of projective space $\mathbb{P}^{
6
}$.

The first order mirror degeneration $\mathfrak{X}^{1\circ}$ induces via
$Y\rightarrow\hat{Y}^{\circ}$ a degeneration $\mathfrak{\hat{X}}^{1\circ
}\subset\hat{Y}^{\circ}\times\operatorname*{Spec}\mathbb{C}\left[  t\right]
/\left\langle t^{2}\right\rangle $ given by the ideal $\hat{I}^{1\circ}%
\subset\hat{S}^{\circ}\otimes\mathbb{C}\left[  t\right]  /\left\langle
t^{2}\right\rangle $ generated by the
 Fermat-type equations
\[
\left\{  

\right\}
\]

Note, that%
\[
\left\vert \operatorname*{supp}\left(  \left(  \mu\left(  B\left(  I\right)
\right)  \right)  ^{\ast}\right)  \cap N-\operatorname*{Roots}\left(  \hat
{Y}^{\circ}\right)  \right\vert -\dim\left(  T_{\hat{Y}^{\circ}}\right)
=
7
-
6
=
1
\]
so this family has one independent parameter.

The special fiber $\hat{X}_{0}^{\circ}\subset\hat{Y}^{\circ}$ of
$\mathfrak{\hat{X}}^{1\circ}$ is cut out by the monomial ideal $\hat{I}%
_{0}^{\circ}\subset\hat{S}^{\circ}$ generated by%
\[
\left\{  
%
\right\}
\]

The complex $SP\left(  \hat{I}_{0}^{\circ}\right)  $ of prime ideals of the
toric strata of $\hat{X}_{0}^{\circ}$ is

\begin{center}
%
\end{center}

\noindent so $\hat{I}_{0}^{\circ}$ has the primary decomposition%
\[
\hat{I}_{0}^{\circ}=
%
\]
The Bergman subcomplex $B\left(  I\right)  \subset\operatorname*{Poset}\left(
\nabla\right)  $ induces a subcomplex of $\operatorname*{Poset}\left(
\hat{\nabla}\right)  $, $\hat{\nabla}=\left(  \hat{P}^{\circ}\right)  ^{\ast}$
corresponding to $\hat{I}_{0}^{\circ}$. Indexing of the 
vertices
of $\hat{\nabla}$ by%
\[
%
\]
this complex is given by
\[
\left\{  
%
\right\}
\]

The ideal $\hat{I}^{1\circ}\subset\hat{S}^{\circ}\otimes\mathbb{C}\left[
t\right]  /\left\langle t^{2}\right\rangle $ has a Pfaffian resolution%
\begin{gather*}
0\rightarrow\mathcal{O}_{\hat{Y}^{\circ}\times\operatorname*{Spec}%
\mathbb{C}\left[  t\right]  /\left\langle t^{2}\right\rangle }\left(
K^{1}\right)  \rightarrow\mathcal{E}^{1}\left(  K^{1}\right)  \overset
{\varphi^{1}}{\rightarrow}\left(  \mathcal{E}^{1}\right)  ^{\ast}%
\overset{f^{1}}{\rightarrow}\mathcal{O}_{\hat{Y}^{\circ}\times
\operatorname*{Spec}\mathbb{C}\left[  t\right]  /\left\langle t^{2}%
\right\rangle }\medskip\\
\text{where }\overline{\pi}_{1}:\hat{Y}^{\circ}\times\operatorname*{Spec}%
\mathbb{C}\left[  t\right]  /\left\langle t^{2}\right\rangle \rightarrow
\hat{Y}^{\circ}\text{ and }\mathcal{E}^{1}=\overline{\pi}_{1}^{\ast
}\mathcal{F}%
\end{gather*}
with%
\[
\mathcal{F}=
\begin{tabular}
[c]{l}
$\mathcal{O}_{\hat{Y}^{\circ}}\left(D_{ \left(0,2,0,0,-1,-1\right) }+D_{ \left(0,0,-1,-1,0,0\right) }+D_{ \left(2,0,0,-1,-1,0\right) } \right) \oplus $
\\
$\mathcal{O}_{\hat{Y}^{\circ}}\left(D_{ \left(0,-1,-1,0,0,2\right) }+D_{ \left(-1,-1,0,0,2,0\right) }+D_{ \left(0,0,-1,-1,0,0\right) } \right) \oplus $
\\
$\mathcal{O}_{\hat{Y}^{\circ}}\left(D_{ \left(0,2,0,0,-1,-1\right) }+D_{ \left(0,0,2,0,0,-1\right) }+D_{ \left(2,0,0,-1,-1,0\right) } \right) \oplus $
\\
$\mathcal{O}_{\hat{Y}^{\circ}}\left(D_{ \left(-1,0,0,2,0,0\right) }+D_{ \left(0,2,0,0,-1,-1\right) }+D_{ \left(0,0,2,0,0,-1\right) } \right) \oplus $
\\
$\mathcal{O}_{\hat{Y}^{\circ}}\left(D_{ \left(-1,0,0,2,0,0\right) }+D_{ \left(-1,-1,0,0,2,0\right) }+D_{ \left(0,0,2,0,0,-1\right) } \right) \oplus $
\\
$\mathcal{O}_{\hat{Y}^{\circ}}\left(D_{ \left(-1,0,0,2,0,0\right) }+D_{ \left(0,-1,-1,0,0,2\right) }+D_{ \left(-1,-1,0,0,2,0\right) } \right) \oplus $
\\
$\mathcal{O}_{\hat{Y}^{\circ}}\left(D_{ \left(0,-1,-1,0,0,2\right) }+D_{ \left(0,0,-1,-1,0,0\right) }+D_{ \left(2,0,0,-1,-1,0\right) } \right)$
\\
\end{tabular}
\]
and $K^{1}=K_{\hat{Y}^{\circ}\times\operatorname*{Spec}\mathbb{C}\left[
t\right]  /\left\langle t^{2}\right\rangle /\operatorname*{Spec}%
\mathbb{C}\left[  t\right]  /\left\langle t^{2}\right\rangle }$ and
$\varphi^{1}\in\bigwedge\nolimits^{2}\mathcal{E}^{1}\left(  -K^{1}\right)  $
given by%
\[
\left [\begin {array}{ccccccc} 0&-ty_{16}\,s_{6}&y_{30}&0&0&-y_{27}&ty_{18}\,s_{2}\\\noalign{\medskip}ty_{16}\,s_{6}&0&-ty_{35}\,s_{3}&y_{22}&0&0&-y_{38}\\\noalign{\medskip}-y_{30}&ty_{35}\,s_{3}&0&-ty_{27}\,s_{1}&y_{18}&0&0\\\noalign{\medskip}0&-y_{22}&ty_{27}\,s_{1}&0&-ty_{38}\,s_{7}&y_{16}&0\\\noalign{\medskip}0&0&-y_{18}&ty_{38}\,s_{7}&0&-ty_{30}\,s_{5}&y_{35}\\\noalign{\medskip}y_{27}&0&0&-y_{16}&ty_{30}\,s_{5}&0&-ty_{22}\,s_{4}\\\noalign{\medskip}-ty_{18}\,s_{2}&y_{38}&0&0&-y_{35}&ty_{22}\,s_{4}&0\end {array}\right ]
\]
Hence via the Pfaffians of $\varphi^{1}$ we obtain a resolution%
\begin{gather*}
0\rightarrow\mathcal{O}_{\hat{Y}^{\circ}\times\operatorname*{Spec}%
\mathbb{C}\left[  t\right]  }\left(  K\right)  \rightarrow\mathcal{E}\left(
K\right)  \rightarrow\mathcal{E}^{\ast}\rightarrow\mathcal{O}_{\hat{Y}^{\circ
}\times\operatorname*{Spec}\mathbb{C}\left[  t\right]  }\medskip\\
\text{where }\pi_{1}:Y\times\operatorname*{Spec}\mathbb{C}\left[  t\right]
\rightarrow Y\text{ , }\mathcal{E}=\pi_{1}^{\ast}\mathcal{F}\\
\text{and }K=K_{\hat{Y}^{\circ}\times\operatorname*{Spec}\mathbb{C}\left[
t\right]  /\operatorname*{Spec}\mathbb{C}\left[  t\right]  }%
\end{gather*}
of the ideal $\hat{I}^{\circ}\subset\hat{S}^{\circ}\otimes\mathbb{C}\left[
t\right]  $ generated by%
\[
\left\{  
\begin{tabular}
[c]{l}
$y_{18}\,y_{38}\,y_{16}+t\left( s_{3}\,y_{16}\,y_{35}^{2}+s_{4}\,y_{18}\,y_{22}^{2}\right)+{t}^{2}\left(-s_{1}\,s_{5}\,y_{27}\,y_{30}\,y_{38}\right)+{t}^{3}\left(-s_{3}\,s_{7}\,s_{4}\,y_{22}\,y_{35}\,y_{38}\right),\medskip$
\\
$y_{16}\,y_{30}\,y_{35}+t\left( s_{2}\,y_{16}\,y_{18}^{2}+s_{1}\,y_{27}^{2}y_{35}\right)+{t}^{2}\left(-s_{7}\,s_{4}\,y_{22}\,y_{30}\,y_{38}\right)+{t}^{3}\left(-s_{2}\,s_{1}\,s_{5}\,y_{18}\,y_{27}\,y_{30}\right),\medskip$
\\
$y_{22}\,y_{35}\,y_{27}+t\left( s_{6}\,y_{16}^{2}y_{35}+s_{7}\,y_{27}\,y_{38}^{2}\right)+{t}^{2}\left(-s_{2}\,s_{5}\,y_{18}\,y_{22}\,y_{30}\right)+{t}^{3}\left(-s_{6}\,s_{7}\,s_{4}\,y_{16}\,y_{22}\,y_{38}\right),\medskip$
\\
$y_{18}\,y_{27}\,y_{38}+t\left( s_{3}\,y_{27}\,y_{35}^{2}+s_{5}\,y_{30}^{2}y_{38}\right)+{t}^{2}\left(-s_{6}\,s_{4}\,y_{16}\,y_{18}\,y_{22}\right)+{t}^{3}\left(-s_{2}\,s_{3}\,s_{5}\,y_{18}\,y_{30}\,y_{35}\right),\medskip$
\\
$y_{16}\,y_{30}\,y_{38}+t\left( s_{4}\,y_{22}^{2}y_{30}+s_{1}\,y_{27}^{2}y_{38}\right)+{t}^{2}\left(-s_{2}\,s_{3}\,y_{16}\,y_{18}\,y_{35}\right)+{t}^{3}\left(-s_{6}\,s_{1}\,s_{4}\,y_{16}\,y_{22}\,y_{27}\right),\medskip$
\\
$y_{22}\,y_{35}\,y_{30}+t\left( s_{2}\,y_{18}^{2}y_{22}+s_{7}\,y_{30}\,y_{38}^{2}\right)+{t}^{2}\left(-s_{6}\,s_{1}\,y_{16}\,y_{27}\,y_{35}\right)+{t}^{3}\left(-s_{2}\,s_{3}\,s_{7}\,y_{18}\,y_{35}\,y_{38}\right),\medskip$
\\
$y_{18}\,y_{27}\,y_{22}+t\left( s_{6}\,y_{16}^{2}y_{18}+s_{5}\,y_{22}\,y_{30}^{2}\right)+{t}^{2}\left(-s_{3}\,s_{7}\,y_{27}\,y_{35}\,y_{38}\right)+{t}^{3}\left(-s_{6}\,s_{1}\,s_{5}\,y_{16}\,y_{27}\,y_{30}\right)$
\\
\end{tabular}
\right\}
\]
which defines a flat family%
\[
\mathfrak{\hat{X}}^{\circ}\subset\hat{Y}^{\circ}\times\operatorname*{Spec}%
\mathbb{C}\left[  t\right]
\]

This is the one parameter mirror family of the generic degree $14$ Pfaffian
Calabi-Yau threefold in $\mathbb{P}^{6}$, given in \cite{Ro dland The Pfaffian
CalabiYau its Mirror and their link to the Grassmannian mathbbG27}.

\subsection{Tropical mirror construction of the degree $13$ Pfaffian
Calabi-Yau\label{Sec Ex 23333} threefold\label{Sec Ex 3333333}}

\subsubsection{Hodge numbers}

The
\index{Hodge numbers}%
Hodge numbers of a general Pfaffian Calabi-Yau threefold $X$ of degree $13$ in
$\mathbb{P}^{6}$, which is smooth as observed above, can be determined as follows.

The
\index{Pfaffian}%
Pfaffian complex is%
\[
0\rightarrow\mathcal{O}_{\mathbb{P}^{6}}\left(  -7\right)  \rightarrow
\mathcal{O}_{\mathbb{P}^{6}}\left(  -5\right)  \oplus\mathcal{O}%
_{\mathbb{P}^{6}}\left(  -4\right)  ^{\oplus4}\rightarrow\mathcal{O}%
_{\mathbb{P}^{6}}\left(  -2\right)  \oplus\mathcal{O}_{\mathbb{P}^{6}}\left(
-3\right)  ^{\oplus4}\rightarrow\mathcal{O}_{\mathbb{P}^{6}}\rightarrow
\mathcal{O}_{X}%
\]
so note that%
\[
H^{i}\left(  X,\mathcal{O}_{X}\right)  \cong H^{i+3}\left(  \mathbb{P}%
^{6},\mathcal{O}_{\mathbb{P}^{6}}\left(  -7\right)  \right)  \cong\left\{
\begin{array}
[c]{cc}%
0 & i=1,2\\
\mathbb{C} & i=3
\end{array}
\right\}
\]

Decomposing the resolution of $\mathcal{J}_{X}^{2}$%
\begin{align*}
0  &  \rightarrow\mathcal{O}_{\mathbb{P}^{6}}\left(  -8\right)  ^{\oplus
6}\oplus\mathcal{O}_{\mathbb{P}^{6}}\left(  -9\right)  ^{\oplus4}%
\rightarrow\mathcal{O}_{\mathbb{P}^{6}}\left(  -6\right)  ^{\oplus4}%
\oplus\mathcal{O}_{\mathbb{P}^{6}}\left(  -7\right)  ^{\oplus16}%
\oplus\mathcal{O}_{\mathbb{P}^{6}}\left(  -8\right)  ^{\oplus4}\rightarrow\\
&  \rightarrow\mathcal{O}_{\mathbb{P}^{6}}\left(  -4\right)  \oplus
\mathcal{O}_{\mathbb{P}^{6}}\left(  -5\right)  ^{\oplus4}\oplus\mathcal{O}%
_{\mathbb{P}^{6}}\left(  -6\right)  ^{\oplus10}\rightarrow\mathcal{J}_{X}%
^{2}\rightarrow0
\end{align*}
into the short exact sequences%
\begin{align*}
0  &  \rightarrow\mathcal{O}_{\mathbb{P}^{6}}\left(  -8\right)  ^{\oplus
6}\oplus\mathcal{O}_{\mathbb{P}^{6}}\left(  -9\right)  ^{\oplus4}%
\rightarrow\mathcal{O}_{\mathbb{P}^{6}}\left(  -6\right)  ^{\oplus4}%
\oplus\mathcal{O}_{\mathbb{P}^{6}}\left(  -7\right)  ^{\oplus16}%
\oplus\mathcal{O}_{\mathbb{P}^{6}}\left(  -8\right)  ^{\oplus4}\rightarrow
\mathcal{K\rightarrow}0\\
0  &  \rightarrow\mathcal{K}\rightarrow\mathcal{O}_{\mathbb{P}^{6}}\left(
-4\right)  \oplus\mathcal{O}_{\mathbb{P}^{6}}\left(  -5\right)  ^{\oplus
4}\oplus\mathcal{O}_{\mathbb{P}^{6}}\left(  -6\right)  ^{\oplus10}%
\rightarrow\mathcal{J}_{X}^{2}\rightarrow0
\end{align*}
the long exact cohomology sequences give
\begin{align*}
...  &  \rightarrow H^{i}\left(  \mathbb{P}^{6},\mathcal{K}\right)
\rightarrow0\rightarrow H^{i}\left(  \mathbb{P}^{6},\mathcal{J}_{X}%
^{2}\right)  \rightarrow\\
&  \rightarrow H^{i+1}\left(  \mathbb{P}^{6},\mathcal{K}\right)
\rightarrow0\rightarrow H^{i+1}\left(  \mathbb{P}^{6},\mathcal{J}_{X}%
^{2}\right)  \rightarrow...
\end{align*}
for $i=0,...,6$, hence%
\[
H^{i}\left(  \mathbb{P}^{6},\mathcal{J}_{X}^{2}\right)  \cong H^{i+1}\left(
\mathbb{P}^{6},\mathcal{K}\right)  \text{ for }i=0,...,5
\]
and
\begin{align*}
0  &  \rightarrow0\rightarrow0\rightarrow H^{0}\left(  \mathbb{P}%
^{6},\mathcal{K}\right)  \rightarrow\\
&  \rightarrow0\rightarrow0\rightarrow H^{1}\left(  \mathbb{P}^{6}%
,\mathcal{K}\right)  \rightarrow\\
&  \rightarrow0\rightarrow0\rightarrow H^{2}\left(  \mathbb{P}^{6}%
,\mathcal{K}\right)  \rightarrow\\
&  \rightarrow0\rightarrow0\rightarrow H^{3}\left(  \mathbb{P}^{6}%
,\mathcal{K}\right)  \rightarrow\\
&  \rightarrow0\rightarrow0\rightarrow H^{4}\left(  \mathbb{P}^{6}%
,\mathcal{K}\right)  \rightarrow\\
&  \rightarrow0\rightarrow0\rightarrow H^{5}\left(  \mathbb{P}^{6}%
,\mathcal{K}\right)  \rightarrow\\
&  \rightarrow H^{6}\left(  \mathbb{P}^{6},\mathcal{O}_{\mathbb{P}^{6}}\left(
-8\right)  ^{\oplus6}\oplus\mathcal{O}_{\mathbb{P}^{6}}\left(  -9\right)
^{\oplus4}\right) \\
&  \rightarrow H^{6}\left(  \mathbb{P}^{6},\mathcal{O}_{\mathbb{P}^{6}}\left(
-6\right)  ^{\oplus4}\oplus\mathcal{O}_{\mathbb{P}^{6}}\left(  -7\right)
^{\oplus16}\oplus\mathcal{O}_{\mathbb{P}^{6}}\left(  -8\right)  ^{\oplus
4}\right) \\
&  \rightarrow H^{6}\left(  \mathbb{P}^{6},\mathcal{K}\right) \\
&  \rightarrow0
\end{align*}
hence%
\[
H^{i}\left(  \mathbb{P}^{6},\mathcal{K}\right)  =0\text{ for }i=0,...,4
\]
and%
\begin{gather*}
h^{5}\left(  \mathbb{P}^{6},\mathcal{K}\right)  -h^{6}\left(  \mathbb{P}%
^{6},\mathcal{O}_{\mathbb{P}^{6}}\left(  -8\right)  ^{\oplus6}\oplus
\mathcal{O}_{\mathbb{P}^{6}}\left(  -9\right)  ^{\oplus4}\right) \\
+h^{6}\left(  \mathbb{P}^{6},\mathcal{O}_{\mathbb{P}^{6}}\left(  -6\right)
^{\oplus4}\oplus\mathcal{O}_{\mathbb{P}^{6}}\left(  -7\right)  ^{\oplus
16}\oplus\mathcal{O}_{\mathbb{P}^{6}}\left(  -8\right)  ^{\oplus4}\right)
-h^{6}\left(  \mathbb{P}^{6},\mathcal{K}\right)  =0
\end{gather*}

Using%
\begin{align*}
h^{6}\left(  \mathbb{P}^{6},\mathcal{O}_{\mathbb{P}^{6}}\left(  -8\right)
^{\oplus6}\oplus\mathcal{O}_{\mathbb{P}^{6}}\left(  -9\right)  ^{\oplus
4}\right)   &  =6\cdot7+4\cdot28=154\\
h^{6}\left(  \mathbb{P}^{6},\mathcal{O}_{\mathbb{P}^{6}}\left(  -6\right)
^{\oplus4}\oplus\mathcal{O}_{\mathbb{P}^{6}}\left(  -7\right)  ^{\oplus
16}\oplus\mathcal{O}_{\mathbb{P}^{6}}\left(  -8\right)  ^{\oplus4}\right)   &
=16+4\cdot7=44
\end{align*}
one has%
\[
h^{5}\left(  \mathbb{P}^{6},\mathcal{K}\right)  -h^{6}\left(  \mathbb{P}%
^{6},\mathcal{K}\right)  =154-44=110
\]
hence:

\begin{lemma}%
\begin{align*}
H^{i}\left(  \mathbb{P}^{6},\mathcal{J}_{X}^{2}\right)   &  =0\text{ for
}i=0,...,3\\
h^{4}\left(  \mathbb{P}^{6},\mathcal{J}_{X}^{2}\right)  -h^{5}\left(
\mathbb{P}^{6},\mathcal{J}_{X}^{2}\right)   &  =110
\end{align*}

\end{lemma}

The long exact cohomology sequence of%
\[
0\rightarrow\mathcal{J}_{X}^{2}\rightarrow\mathcal{J}_{X}\rightarrow
\mathcal{N}_{X/\mathbb{P}^{6}}^{\vee}\rightarrow0
\]
reads%

\begin{align*}
H^{i}\left(  \mathbb{P}^{6},\mathcal{N}_{X/\mathbb{P}^{6}}^{\vee}\right)   &
=0\text{ for }i=0,1,2\\
H^{4}\left(  \mathbb{P}^{6},\mathcal{J}_{X}^{2}\right)  /H^{3}\left(
\mathbb{P}^{6},\mathcal{N}_{X/\mathbb{P}^{6}}^{\vee}\right)   &
\cong\mathbb{C}\\
H^{5}\left(  \mathbb{P}^{6},\mathcal{J}_{X}^{2}\right)   &  =0
\end{align*}
so $h^{4}\left(  \mathbb{P}^{6},\mathcal{J}_{X}^{2}\right)  =110$ and
$h^{3}\left(  \mathbb{P}^{6},\mathcal{N}_{X/\mathbb{P}^{6}}^{\vee}\right)
=109$, hence:

\begin{lemma}
The cohomology dimensions of $\mathcal{J}_{X}^{2}$ and $\mathcal{N}%
_{X/\mathbb{P}^{6}}^{\vee}$ are%
\begin{align*}
h^{i}\left(  \mathbb{P}^{6},\mathcal{J}_{X}^{2}\right)   &  =0\text{ for
}i=0,...,3\\
h^{4}\left(  \mathbb{P}^{6},\mathcal{J}_{X}^{2}\right)   &  =110\\
h^{5}\left(  \mathbb{P}^{6},\mathcal{J}_{X}^{2}\right)   &  =0
\end{align*}
and%
\begin{align*}
h^{i}\left(  \mathbb{P}^{6},\mathcal{N}_{X/\mathbb{P}^{6}}^{\vee}\right)   &
=0\text{ for }i=0,1,2\\
h^{3}\left(  \mathbb{P}^{6},\mathcal{N}_{X/\mathbb{P}^{6}}^{\vee}\right)   &
=109
\end{align*}

\end{lemma}

Using the Euler sequence and conormal sequence
\begin{align*}
0  &  \rightarrow\Omega_{\mathbb{P}^{6}}\mid_{X}\rightarrow\mathcal{O}%
_{X}\left(  -1\right)  ^{\oplus7}\rightarrow\mathcal{O}_{X}\rightarrow0\\
0  &  \rightarrow\mathcal{N}_{X}^{\vee}\rightarrow\Omega_{\mathbb{P}^{6}}%
\mid_{X}\rightarrow\Omega_{X}\rightarrow0
\end{align*}
as explained in Section \ref{Sec Remarks on the Hodge numbers of Calabi-Yau},
we get%
\begin{align*}
h^{1}\left(  X,\Omega_{X}\right)   &  =1\\
h^{2}\left(  X,\Omega_{X}\right)   &  =h^{3}\left(  X,\mathcal{N}%
_{X/\mathbb{P}^{6}}^{\vee}\right)  -h^{3}\left(  X,\Omega_{\mathbb{P}^{6}}%
\mid_{X}\right)  =109-48=61
\end{align*}

\begin{corollary}
The general degree $13$ Pfaffian Calabi-Yau threefold $X$ in $\mathbb{P}^{6}$
has%
\begin{align*}
h^{1,1}\left(  X\right)   &  =1\\
h^{1,2}\left(  X\right)   &  =61
\end{align*}

\end{corollary}

\subsubsection{Setup}

Let $Y=\mathbb{P}^{
6
}=X\left(  \Sigma\right)  $, $\Sigma
=\operatorname*{Fan}\left(  P\right)  =NF\left(  \Delta\right)  \subset
N_{\mathbb{R}}$ with the Fano polytope $P=\Delta^{\ast}$ given by%
\[
\Delta=\operatorname*{convexhull}\left(  
\begin{tabular}
[c]{ll}
$\left(6,-1,-1,-1,-1,-1\right)$ &$\left(-1,6,-1,-1,-1,-1\right)$\\
$\left(-1,-1,6,-1,-1,-1\right)$ &$\left(-1,-1,-1,6,-1,-1\right)$\\
$\left(-1,-1,-1,-1,6,-1\right)$ &$\left(-1,-1,-1,-1,-1,6\right)$\\
$\left(-1,-1,-1,-1,-1,-1\right)$ &
\end{tabular}
\right)  \subset
M_{\mathbb{R}}%
\]
and let%
\[
S=\mathbb{C}
[x_0, x_1, x_2, x_3, x_4, x_5, x_6]
\]
be the Cox ring of $Y$ with the variables%
\[
\begin{tabular}
[c]{ll}
$x_{1} = x_{\left(1,0,0,0,0,0\right)}$ &$x_{2} = x_{\left(0,1,0,0,0,0\right)}$\\
$x_{3} = x_{\left(0,0,1,0,0,0\right)}$ &$x_{4} = x_{\left(0,0,0,1,0,0\right)}$\\
$x_{5} = x_{\left(0,0,0,0,1,0\right)}$ &$x_{6} = x_{\left(0,0,0,0,0,1\right)}$\\
$x_{0} = x_{\left(-1,-1,-1,-1,-1,-1\right)}$ &
\end{tabular}
\]
associated to the rays of $\Sigma$.

Consider the degeneration $\mathfrak{X}\subset Y\times\operatorname*{Spec}%
\mathbb{C}\left[  t\right]  $ of 
Pfaffian
Calabi-Yau 3-folds
 with Buchsbaum-Eisenbud
resolution%
\[
0\rightarrow\mathcal{O}_{Y}\left(  -
7
\right)  \rightarrow
\mathcal{E}\left(  -
3
\right)  \overset{A_{t}}{\rightarrow
}\mathcal{E}^{\ast}\left(  -
2
\right)  \rightarrow\mathcal{O}%
_{Y}\rightarrow\mathcal{O}_{X_{t}}\rightarrow0
\]
where%
\[
\mathcal{E}=\text{
$\mathcal{O}\left(  1\right)  \oplus4\mathcal{O}$
}%
\]%
\[
A_{t}=A_{0}+t\cdot A
\]%
\[
A_{0}=
\left [
\right ]
\]
the monomial special fiber of $\mathfrak{X}$ is given by%
\[
I_{0}=\left\langle 
%
\]
The degeneration%
\[
\mathfrak{X}\subset Y\times\operatorname*{Spec}\mathbb{C}\left[  t\right]
\]
is given by the ideal $I\subset S\otimes\mathbb{C}\left[  t\right]  $
generated by the Pfaffians of $A_{0}+t\cdot A$ of degrees $
2,3,3,3,3
$.

\subsubsection{Special fiber Gr\"{o}bner cone}

The space of first order deformations of $\mathfrak{X}$ has dimension
$
109
$ and the deformations represented by the Cox Laurent monomials

\begin{center}
%
\end{center}

\noindent form a torus invariant basis. These deformations give linear
inequalities defining the special fiber Gr\"{o}bner cone, i.e.,%
\[
C_{I_{0}}\left(  I\right)  =\left\{  \left(  w_{t},w_{y}\right)  \in
\mathbb{R}\oplus N_{\mathbb{R}}\mid\left\langle w,A^{-1}\left(  m\right)
\right\rangle \geq-w_{t}\forall m\right\}
\]
for above monomials $m$ and the presentation matrix
\[
A=
\left [%
\right ]
\]
of $A_{
5
}\left(  Y\right)  $.

The vertices of special fiber polytope
\[
\nabla=C_{I_{0}}\left(  I\right)  \cap\left\{  w_{t}=1\right\}
\]
and the number of faces, each vertex is contained in, are given by the table

\begin{center}
%
\end{center}

\noindent The number of faces of $\nabla$ and their $F$-vectors are

\begin{center}
%
\end{center}

The dual $P^{\circ}=\nabla^{\ast}$ of $\nabla$ is a Fano polytope with vertices

\begin{center}
%
\end{center}

\noindent and the $F$-vectors of the faces of $P^{\circ}$ are

\begin{center}
%
\end{center}

The polytopes $\nabla$ and $P^{\circ}$ have $0$ as the unique interior lattice
point, and the Fano polytope $P^{\circ}$ defines the embedding toric Fano
variety $Y^{\circ}=X\left(  \operatorname*{Fan}\left(  P^{\circ}\right)
\right)  $ of the fibers of the mirror degeneration. The dimension of the
automorphism group of $Y^{\circ}$ is
\[
\dim\left(  \operatorname*{Aut}\left(  Y^{\circ}\right)  \right)
=
6
\]
Let%
\[
S^{\circ}=\mathbb{C}
[y_{1},y_{2},...,y_{49}]
\]
be the Cox ring of $Y^{\circ}$ with variables

\begin{center}
%
\end{center}

\subsubsection{Bergman subcomplex}

Intersecting the tropical variety of $I$ with the special fiber Gr\"{o}bner
cone $C_{I_{0}}\left(  I\right)  $ we obtain the special fiber Bergman
subcomplex
\[
B\left(  I\right)  =C_{I_{0}}\left(  I\right)  \cap BF\left(  I\right)
\cap\left\{  w_{t}=1\right\}
\]
The following table chooses an indexing of the vertices of $\nabla$ involved
in $B\left(  I\right)  $ and gives for each vertex the numbers $n_{\nabla}$
and $n_{B}$ of faces of $\nabla$ and faces of $B\left(  I\right)  $ it is
contained in

\begin{center}
%
\end{center}

\noindent With this indexing the Bergman subcomplex $B\left(  I\right)  $ of
$\operatorname*{Poset}\left(  \nabla\right)  $ associated to the degeneration
$\mathfrak{X}$ is

\begin{center}
%
\]
and the $F$-vectors of its faces are

\begin{center}
%
\end{center}

\subsubsection{Dual complex}

The dual complex $\operatorname*{dual}\left(  B\left(  I\right)  \right)
=\left(  B\left(  I\right)  \right)  ^{\ast}$ of deformations associated to
$B\left(  I\right)  $ via initial ideals is given by

\begin{center}
%
\end{center}

\noindent when writing the vertices of the faces as deformations of $X_{0}$.
Note that the $T$-invariant basis of deformations associated to a face is
given by all lattice points of the corresponding polytope in $M_{\mathbb{R}}$.

In order to compress the output we list one representative in any set of faces
$G$ with fixed $F$-vector of $G$ and $G^{\ast}$.

When numbering the vertices of the faces of $\operatorname*{dual}\left(
B\left(  I\right)  \right)  $ by the Cox variables of the mirror toric Fano
variety $Y^{\circ}$ the complex $\operatorname*{dual}\left(  B\left(
I\right)  \right)  $ is

\begin{center}

\end{center}

\noindent The dual complex has the $F$-vector%
\[
%
\]
and the $F$-vectors of the faces of $\operatorname*{dual}\left(  B\left(
I\right)  \right)  $ are

\begin{center}
%
\end{center}

Recall that in this example the toric variety $Y$ is projective space. The
number of lattice points of the support of $\operatorname*{dual}\left(
B\left(  I\right)  \right)  $ relates to the dimension $h^{1,
2
}%
\left(  X\right)  $ of the complex moduli space of the generic fiber $X$ of
$\mathfrak{X}$ and to the dimension $h^{1,1}\left(  \bar{X}^{\circ}\right)  $
of the K\"{a}hler moduli space of the $MPCR$-blowup $\bar{X}^{\circ}$ of the
generic fiber $X^{\circ}$ of the mirror degeneration
\begin{align*}
\left\vert \operatorname*{supp}\left(  \operatorname*{dual}\left(  B\left(
I\right)  \right)  \right)  \cap M\right\vert  &
=
109
=
48
+
61
=\dim\left(  \operatorname*{Aut}%
\left(  Y\right)  \right)  +h^{1,
2
}\left(  X\right) \\
&  =
42
+
6
+
61
\\
&  =\left\vert \operatorname*{Roots}\left(  Y\right)  \right\vert +\dim\left(
T_{Y}\right)  +h^{1,1}\left(  \bar{X}^{\circ}\right)
\end{align*}
There are%
\[
h^{1,
2
}\left(  X\right)  +\dim\left(  T_{Y^{\circ}}\right)
=
61
+
6
\]
non-trivial toric polynomial deformations of $X_{0}$

\begin{center}

\end{center}

\noindent They correspond to the toric divisors

\begin{center}
%
\end{center}

\noindent on a MPCR-blowup of $Y^{\circ}$ inducing $
61
$ non-zero
toric divisor classes on the mirror $X^{\circ}$. The following $
42
$
toric divisors of $Y$ induce the trivial divisor class on $X^{\circ}$

\begin{center}
%
\end{center}

\subsubsection{Mirror special fiber}

The complex $B\left(  I\right)  ^{\ast}$ labeled by the variables of the Cox
ring $S^{\circ}$ of $Y^{\ast}$, as written in the last section, is the complex
$SP\left(  I_{0}^{\circ}\right)  $ of prime ideals of the toric strata of the
monomial special fiber $X_{0}^{\circ}$ of the mirror degeneration
$\mathfrak{X}^{\circ}$, i.e. the primary decomposition of $I_{0}^{\circ}$ is%
\[
I_{0}^{\circ}=
%
\]
Each facet $F\in B\left(  I\right)  $ corresponds to one of these ideals and this
ideal is generated by the products of facets of $\nabla$ containing $F$.

\subsubsection{Covering structure in the deformation complex of the
degeneration $\mathfrak{X}$}

According to the local reduced Gr\"{o}bner basis each face of the complex of
deformations $\operatorname*{dual}\left(  B\left(  I\right)  \right)  $
decomposes into 
$3$ respectively $5$
polytopes forming a $
5
:1$
ramified
covering of $B\left(  I\right)  $

\begin{center}
%
\end{center}

\noindent Here the faces $F\in B\left(  I\right)  $ are specified both via the
vertices of $F^{\ast}$ labeled by the variables of $S^{\circ}$ and by the
numbering of the vertices of $B\left(  I\right)  $ chosen above.

The numbers of faces of the covering in each face of $\operatorname*{dual}%
\left(  B\left(  I\right)  \right)  $, i.e. over each face of $B\left(
I\right)  ^{\vee}$ are

\begin{center}
%
\end{center}

\noindent This covering has 
one sheet forming the complex
 
\[
%
\end{center}

\noindent Writing the vertices of the faces as deformations the covering is
given by

\begin{center}
%
\end{center}

\noindent with 
one sheet forming the complex

\begin{center}
%
\end{center}

\noindent Note that the torus invariant basis of deformations corresponding to
a Bergman face is given by the set of all lattice points of the polytope
specified above.

\subsubsection{Limit map}

The limit map $\lim:B\left(  I\right)  \rightarrow\operatorname*{Poset}\left(
\Delta\right)  $ associates to a face $F$ of $B\left(  I\right)  $ the face of
$\Delta$ formed by the limit points of arcs lying over the weight vectors
$w\in F$, i.e. with lowest order term $t^{w}$.

Labeling the faces of the Bergman complex $B\left(  I\right)  \subset
\operatorname*{Poset}\left(  \nabla\right)  $ and the faces of
$\operatorname*{Poset}\left(  \Delta\right)  $ by the corresponding dual faces
of $\nabla^{\ast}$ and $\Delta^{\ast}$, hence considering the limit map
$\lim:B\left(  I\right)  \rightarrow\operatorname*{Poset}\left(
\Delta\right)  $ as a map $B\left(  I\right)  ^{\ast}\rightarrow
\operatorname*{Poset}\left(  \Delta^{\ast}\right)  $, the limit correspondence
is given by

\begin{center}

\end{center}

The image of the limit map coincides with the image of $\mu$, i.e.
$\lim\left(  B\left(  I\right)  \right)  =\mu\left(  B\left(  I\right)
\right)  $.

\subsubsection{Mirror complex}

Numbering the vertices of the mirror complex $\mu\left(  B\left(  I\right)
\right)  $ as%
\[
%
\end{center}

\noindent The ordering of the faces of $\mu\left(  B\left(  I\right)  \right)
$ is compatible with above ordering of the faces of $B\left(  I\right)  $. The
$F$-vectors of the faces of $\mu\left(  B\left(  I\right)  \right)  $ are%
\[
%
\]
The first order deformations of the mirror special fiber $X_{0}^{\circ}$
correspond to the lattice points of the dual complex $\left(  \mu\left(
B\left(  I\right)  \right)  \right)  ^{\ast}\subset\operatorname*{Poset}%
\left(  \Delta^{\ast}\right)  $ of the mirror complex. We label the first
order deformations of $X_{0}^{\circ}$ corresponding to vertices of
$\Delta^{\ast}$ by the homogeneous coordinates of $Y$%
\[
%
\]
So writing the vertices of the faces of $\left(  \mu\left(  B\left(  I\right)
\right)  \right)  ^{\ast}$ as homogeneous coordinates of $Y$, the complex
$\left(  \mu\left(  B\left(  I\right)  \right)  \right)  ^{\ast}$ is given by

\begin{center}
%
\end{center}

\noindent The complex $\mu\left(  B\left(  I\right)  \right)  ^{\ast}$ labeled
by the variables of the Cox ring $S$ gives the ideals of the toric strata of
the special fiber $X_{0}$ of $\mathfrak{X}$, i.e. the complex $SP\left(
I_{0}\right)  $, so in particular the primary decomposition of $I_{0}$ is%
\[
I_{0}=
%
\]

\subsubsection{Covering structure in the deformation complex of the mirror
degeneration}

Each face of the complex of deformations $\left(  \mu\left(  B\left(
I\right)  \right)  \right)  ^{\ast}$ of the mirror special fiber $X_{0}%
^{\circ}$ decomposes as the convex hull of 
$3$,$4$ respectively $5$
polytopes forming
a $
5
:1$ 
ramified
covering of $\mu\left(  B\left(
I\right)  \right)  ^{\vee}$

\begin{center}
%
\end{center}

\noindent 
Due to the singularities of
$Y^{\circ}$ this covering involves degenerate faces, i.e. faces $G\mapsto
F^{\vee}$ with $\dim\left(  G\right)  <\dim\left(  F^{\vee}\right)  $.

\subsubsection{Mirror degeneration}

The space of first order deformations of $X_{0}^{\circ}$ in the mirror
degeneration $\mathfrak{X}^{\circ}$ has dimension $
7
$ and the
deformations represented by the monomials%
\[
\left\{  
%
\right\}
\]
form a torus invariant
basis $\mathfrak{B}^{\circ}$.
The number of lattice points of the dual of the mirror complex of $I$ relates
to the dimension $h^{1,
2
}\left(  X^{\circ}\right)  $ of complex
moduli space of the generic fiber $X^{\circ}$ of $\mathfrak{X}^{\circ}$ and to
the dimension $h^{1,1}\left(  X\right)  $ of the K\"{a}hler moduli space of
the generic fiber $X$ of $\mathfrak{X}$ via%
\begin{align*}
\left\vert \operatorname*{supp}\left(  \left(  \mu\left(  B\left(  I\right)
\right)  \right)  ^{\ast}\right)  \cap N\right\vert  &
=
7
=
6
+
1
\\
&  =\dim\left(  \operatorname*{Aut}\left(  Y^{\circ}\right)  \right)
+h^{1,
2
}\left(  X^{\circ}\right)  =\dim\left(  T\right)
+h^{1,1}\left(  X\right)
\end{align*}

The 
conjectural first order
mirror degeneration 
$\mathfrak{X}^{1\circ}\subset Y^{\circ}\times\operatorname*{Spec}%
\mathbb{C}\left[  t\right]  /\left\langle t^{2}\right\rangle $
 of $\mathfrak{X}$ is
given by the ideal 
$I^{1\circ}\subset S^{\circ}\otimes\mathbb{C}\left[  t\right]  /\left\langle
t^{2}\right\rangle $
generated by%
\[
\left\{  m+\sum_{\delta\in\mathfrak{B}^{\circ}}t\cdot s_{\delta}\cdot
\delta\left(  m\right)  \mid m\in I_{0}^{\circ}\right\}
\]

\subsubsection{Contraction of the mirror degeneration}

In the following we give a birational map relating the degeneration
$\mathfrak{X}^{\circ}$ to a Greene-Plesser type orbifolding mirror family by
contracting divisors on $Y^{\circ}$.

In order to contract the divisors

\begin{center}

\end{center}

\noindent consider the $\mathbb{Q}$-factorial toric Fano variety $\hat
{Y}^{\circ}=X\left(  \hat{\Sigma}^{\circ}\right)  $, where $\hat{\Sigma
}^{\circ}=\operatorname*{Fan}\left(  \hat{P}^{\circ}\right)  \subset
M_{\mathbb{R}}$ is the fan over the 
Fano polytope $\hat
{P}^{\circ}\subset M_{\mathbb{R}}$ given as the convex hull of the remaining
vertices of $P^{\circ}=\nabla^{\ast}$ corresponding to the Cox variables%
\[
%
\]
of the toric variety $\hat{Y}^{\circ}$ with Cox ring
\[
\hat{S}^{\circ}=\mathbb{C}
[y_{9},y_{6},y_{5},y_{22},y_{8},y_{10},y_{7}]
\]
The Cox variables of $\hat{Y}^{\circ}$ correspond to the set of Fermat deformations of $\mathfrak{X}$.

Let%
\[
Y^{\circ}=X\left(  \Sigma^{\circ}\right)  \rightarrow X\left(  \hat{\Sigma
}^{\circ}\right)  =\hat{Y}^{\circ}%
\]
be a birational map from $Y^{\circ}$ to a minimal birational model $\hat
{Y}^{\circ}$, which contracts the divisors of the rays $\Sigma^{\circ}\left(
1\right)  -\hat{\Sigma}^{\circ}\left(  1\right)  $ corresponding to the Cox
variables
\[
\begin{tabular}
[c]{lllllllllllllll}
$y_{1}$ &$y_{2}$ &$y_{3}$ &$y_{4}$ &$y_{11}$ &$y_{12}$ &$y_{13}$ &$y_{14}$ &$y_{15}$ &$y_{16}$ &$y_{17}$ &$y_{18}$ &$y_{19}$ &$y_{20}$ &$y_{21}$\\
$y_{23}$ &$y_{24}$ &$y_{25}$ &$y_{26}$ &$y_{27}$ &$y_{28}$ &$y_{29}$ &$y_{30}$ &$y_{31}$ &$y_{32}$ &$y_{33}$ &$y_{34}$ &$y_{35}$ &$y_{36}$ &$y_{37}$\\
$y_{38}$ &$y_{39}$ &$y_{40}$ &$y_{41}$ &$y_{42}$ &$y_{43}$ &$y_{44}$ &$y_{45}$ &$y_{46}$ &$y_{47}$ &$y_{48}$ &$y_{49}$ & & &
\end{tabular}
\]

Representing $\hat{Y}^{\circ}$ as a quotient we have
\[
\hat{Y}^{\circ}=\left(  \mathbb{C}^{
7
}-V\left(  B\left(
\hat{\Sigma}^{\circ}\right)  \right)  \right)  //\hat{G}^{\circ}%
\]
with%
\[
\hat{G}^{\circ}=
\mathbb{Z}_{13}\times\left(  \mathbb{C}^{\ast}\right)  ^{1}
\]
acting via%
\[
\xi y=
\left( \,u_{1}^{11}\,v_{1} \cdot y_{9},\,u_{1}^{10}\,v_{1} \cdot y_{6},\,u_{1}^{10}\,v_{1} \cdot y_{5},\,u_{1}^{4}\,v_{1} \cdot y_{22},\,u_{1}^{8}\,v_{1} \cdot y_{8},\,u_{1}^{11}\,v_{1} \cdot y_{10},\,v_{1} \cdot y_{7} \right)
\]
for $\xi=
\left(u_1,v_1\right)
\in\hat{G}^{\circ}$ and $y\in\mathbb{C}%
^{
7
}-V\left(  B\left(  \hat{\Sigma}^{\circ}\right)  \right)  $.

Hence with the group%
\[
\hat{H}^{\circ}=
\mathbb{Z}_{13}
\]
of order 
13
 the toric variety $\hat{Y}^{\circ}$ is the quotient%
\[
\hat{Y}^{\circ}=\mathbb{P}^{
6
}/\hat{H}^{\circ}%
\]
of projective space $\mathbb{P}^{
6
}$.

The first order mirror degeneration $\mathfrak{X}^{1\circ}$ induces via
$Y\rightarrow\hat{Y}^{\circ}$ a degeneration $\mathfrak{\hat{X}}^{1\circ
}\subset\hat{Y}^{\circ}\times\operatorname*{Spec}\mathbb{C}\left[  t\right]
/\left\langle t^{2}\right\rangle $ given by the ideal $\hat{I}^{1\circ}%
\subset\hat{S}^{\circ}\otimes\mathbb{C}\left[  t\right]  /\left\langle
t^{2}\right\rangle $ generated by the
 Fermat-type equations
\[
\left\{  

\right\}
\]

Note, that%
\[
\left\vert \operatorname*{supp}\left(  \left(  \mu\left(  B\left(  I\right)
\right)  \right)  ^{\ast}\right)  \cap N-\operatorname*{Roots}\left(  \hat
{Y}^{\circ}\right)  \right\vert -\dim\left(  T_{\hat{Y}^{\circ}}\right)
=
7
-
6
=
1
\]
so this family has one independent parameter.

The special fiber $\hat{X}_{0}^{\circ}\subset\hat{Y}^{\circ}$ of
$\mathfrak{\hat{X}}^{1\circ}$ is cut out by the monomial ideal $\hat{I}%
_{0}^{\circ}\subset\hat{S}^{\circ}$ generated by%
\[
\left\{  
%
\right\}
\]

The complex $SP\left(  \hat{I}_{0}^{\circ}\right)  $ of prime ideals of the
toric strata of $\hat{X}_{0}^{\circ}$ is

\begin{center}
%
\end{center}

\noindent so $\hat{I}_{0}^{\circ}$ has the primary decomposition%
\[
\hat{I}_{0}^{\circ}=
%
\]
The Bergman subcomplex $B\left(  I\right)  \subset\operatorname*{Poset}\left(
\nabla\right)  $ induces a subcomplex of $\operatorname*{Poset}\left(
\hat{\nabla}\right)  $, $\hat{\nabla}=\left(  \hat{P}^{\circ}\right)  ^{\ast}$
corresponding to $\hat{I}_{0}^{\circ}$. Indexing of the 
vertices
of $\hat{\nabla}$ by%
\[
%
\]
this complex is given by
\[
\left\{  
%
\right\}
\]

The ideal $\hat{I}^{1\circ}\subset\hat{S}^{\circ}\otimes\mathbb{C}\left[
t\right]  /\left\langle t^{2}\right\rangle $ has a Pfaffian resolution%
\begin{gather*}
0\rightarrow\mathcal{O}_{\hat{Y}^{\circ}\times\operatorname*{Spec}%
\mathbb{C}\left[  t\right]  /\left\langle t^{2}\right\rangle }\left(
K^{1}\right)  \rightarrow\mathcal{E}^{1}\left(  K^{1}\right)  \overset
{\varphi^{1}}{\rightarrow}\left(  \mathcal{E}^{1}\right)  ^{\ast}%
\overset{f^{1}}{\rightarrow}\mathcal{O}_{\hat{Y}^{\circ}\times
\operatorname*{Spec}\mathbb{C}\left[  t\right]  /\left\langle t^{2}%
\right\rangle }\medskip\\
\text{where }\overline{\pi}_{1}:\hat{Y}^{\circ}\times\operatorname*{Spec}%
\mathbb{C}\left[  t\right]  /\left\langle t^{2}\right\rangle \rightarrow
\hat{Y}^{\circ}\text{ and }\mathcal{E}^{1}=\overline{\pi}_{1}^{\ast
}\mathcal{F}%
\end{gather*}
with%
\[
\mathcal{F}=
\begin{tabular}
[c]{l}
$\mathcal{O}_{\hat{Y}^{\circ}}\left(D_{ \left(-1,-1,0,0,0,-1\right) }+D_{ \left(0,0,-1,-1,3,-1\right) } \right) \oplus $
\\
$\mathcal{O}_{\hat{Y}^{\circ}}\left(D_{ \left(2,1,-1,-1,0,0\right) }+D_{ \left(1,2,-1,-1,0,0\right) }+D_{ \left(0,0,-1,-1,3,-1\right) } \right) \oplus $
\\
$\mathcal{O}_{\hat{Y}^{\circ}}\left(D_{ \left(-1,-1,1,2,-1,0\right) }+D_{ \left(0,0,0,0,-1,2\right) }+D_{ \left(-1,-1,2,1,-1,0\right) } \right) \oplus $
\\
$\mathcal{O}_{\hat{Y}^{\circ}}\left(D_{ \left(-1,-1,1,2,-1,0\right) }+D_{ \left(-1,-1,0,0,0,-1\right) }+D_{ \left(-1,-1,2,1,-1,0\right) } \right) \oplus $
\\
$\mathcal{O}_{\hat{Y}^{\circ}}\left(D_{ \left(2,1,-1,-1,0,0\right) }+D_{ \left(1,2,-1,-1,0,0\right) }+D_{ \left(0,0,0,0,-1,2\right) } \right)$
\\
\end{tabular}
\]
and $K^{1}=K_{\hat{Y}^{\circ}\times\operatorname*{Spec}\mathbb{C}\left[
t\right]  /\left\langle t^{2}\right\rangle /\operatorname*{Spec}%
\mathbb{C}\left[  t\right]  /\left\langle t^{2}\right\rangle }$ and
$\varphi^{1}\in\bigwedge\nolimits^{2}\mathcal{E}^{1}\left(  -K^{1}\right)  $
given by%
\[
\left [\begin {array}{ccccc} 0&ts_{6}\,{y_{7}}^{2}&y_{9}\,y_{10}&-y_{5}\,y_{6}&ts_{7}\,{y_{8}}^{2}\\\noalign{\medskip}-ts_{6}\,{y_{7}}^{2}&0&t\left (s_{3}\,y_{5}+s_{2}\,y_{6}\right )&y_{8}&y_{22}\\\noalign{\medskip}-y_{9}\,y_{10}&-t\left (s_{3}\,y_{5}+s_{2}\,y_{6}\right )&0&ts_{1}\,y_{22}&-y_{7}\\\noalign{\medskip}y_{5}\,y_{6}&-y_{8}&-ts_{1}\,y_{22}&0&t\left (-s_{5}\,y_{9}-s_{4}\,y_{10}\right )\\\noalign{\medskip}-ts_{7}\,{y_{8}}^{2}&-y_{22}&y_{7}&-t\left (-s_{5}\,y_{9}-s_{4}\,y_{10}\right )&0\end {array}\right ]
\]
Hence via the Pfaffians of $\varphi^{1}$ we obtain a resolution%
\begin{gather*}
0\rightarrow\mathcal{O}_{\hat{Y}^{\circ}\times\operatorname*{Spec}%
\mathbb{C}\left[  t\right]  }\left(  K\right)  \rightarrow\mathcal{E}\left(
K\right)  \rightarrow\mathcal{E}^{\ast}\rightarrow\mathcal{O}_{\hat{Y}^{\circ
}\times\operatorname*{Spec}\mathbb{C}\left[  t\right]  }\medskip\\
\text{where }\pi_{1}:Y\times\operatorname*{Spec}\mathbb{C}\left[  t\right]
\rightarrow Y\text{ , }\mathcal{E}=\pi_{1}^{\ast}\mathcal{F}\\
\text{and }K=K_{\hat{Y}^{\circ}\times\operatorname*{Spec}\mathbb{C}\left[
t\right]  /\operatorname*{Spec}\mathbb{C}\left[  t\right]  }%
\end{gather*}
of the ideal $\hat{I}^{\circ}\subset\hat{S}^{\circ}\otimes\mathbb{C}\left[
t\right]  $ generated by%
\[
\left\{  
\begin{tabular}
[c]{l}
$y_{7}\,y_{8}+t\left( s_{1}\,y_{22}^{2}\right)+{t}^{2}\left(-s_{3}\,y_{5}\,s_{5}\,y_{9}-s_{3}\,y_{5}\,s_{4}\,y_{10}-s_{2}\,y_{6}\,s_{5}\,y_{9}-s_{2}\,y_{6}\,s_{4}\,y_{10}\right),\medskip$
\\
$y_{5}\,y_{6}\,y_{7}+t\left( s_{5}\,y_{9}^{2}y_{10}+s_{4}\,y_{9}\,y_{10}^{2}\right)+{t}^{2}\left(-s_{7}\,s_{1}\,y_{8}^{2}y_{22}\right),\medskip$
\\
$y_{5}\,y_{6}\,y_{22}+t\left( s_{7}\,y_{8}^{3}\right)+{t}^{2}\left(-s_{6}\,y_{7}^{2}s_{5}\,y_{9}-s_{6}\,y_{7}^{2}s_{4}\,y_{10}\right),\medskip$
\\
$y_{9}\,y_{10}\,y_{22}+t\left( s_{6}\,y_{7}^{3}\right)+{t}^{2}\left(-s_{7}\,y_{8}^{2}s_{3}\,y_{5}-s_{7}\,y_{8}^{2}s_{2}\,y_{6}\right),\medskip$
\\
$-y_{8}\,y_{9}\,y_{10}+t\left( -s_{3}\,y_{5}^{2}y_{6}-s_{2}\,y_{5}\,y_{6}^{2}\right)+{t}^{2}\left(s_{6}\,s_{1}\,y_{7}^{2}y_{22}\right)$
\\
\end{tabular}
\right\}
\]
which defines a flat family%
\[
\mathfrak{\hat{X}}^{\circ}\subset\hat{Y}^{\circ}\times\operatorname*{Spec}%
\mathbb{C}\left[  t\right]
\]

\section{Remarks on a tropical computation of the stringy $E$%
-function\label{Sec tropical stringy E}}

Suppose we are given the setup of the tropical mirror construction via a
degeneration $\mathfrak{X}$ given by the ideal $I$. In the following we make
some remarks on the computation of Hodge numbers and the stringy $E$-function
of the general fiber from the tropical data, i.e., from the polytopes $\Delta$
and $\nabla$ and the complexes $B\left(  I\right)  \subset
\operatorname*{Poset}\left(  \nabla\right)  $ and $\lim\left(  B\left(
I\right)  \right)  \subset\operatorname*{Poset}\left(  \Delta\right)  $.

We recall in Sections
\ref{Sec String theoretic Hodge formula for hypersurfaces} and
\ref{Sec string theoretic Hodge formula for complete intersections} the
formulas by Batyrev and Borisov for the stringy $E$-function of a general
Calabi-Yau hypersurface inside a Gorenstein toric Fano variety and for
complete intersections given by nef partitions. These formulas give evidence
that it should be possible to compute the stringy $E$-function from the
tropical data via a formula analogous to those for hypersurfaces. Note also
that stringy $E$-functions and tropical geometry share the concept of formal
arcs. Furthermore the special fiber $X_{0}$ of $\mathfrak{X}$ is a union of
toric varieties and, as noted in Proposition
\ref{prop stringy E-function and stratifications} below, the stringy
$E$-function respects stratifications.

As this gives the general direction, we begin by recalling in Section
\ref{Sec Remarks on the Hodge numbers of Calabi-Yau} the relation of
$h^{d-1,1}\left(  X\right)  $, $h^{0}\left(  X,N_{X/\mathbb{P}^{n}}\right)  $
and $\operatorname*{Aut}\left(  \mathbb{P}^{n}\right)  $ for Calabi-Yau
manifolds of dimension $d$ in projective space $\mathbb{P}^{n}$.

\subsection{Hodge numbers of Calabi-Yau manifolds in $\mathbb{P}^{n}$ and the
relation between $h^{d-1,1}\left(  X\right)  $, $h^{0}\left(
X,N_{X/\mathbb{P}^{n}}\right)  $ and $\operatorname*{Aut}\left(
\mathbb{P}^{n}\right)  $\label{Sec Remarks on the Hodge numbers of Calabi-Yau}%
}

Let $X\subset\mathbb{P}^{n}$ be a
\index{Hodge numbers}%
Calabi-Yau $d$-fold for $d\geq3$.

\begin{itemize}
\item Note that for a Calabi-Yau $d$-fold
\[
T_{X}=\wedge^{1}\Omega_{X}^{1^{\ast}}\cong\wedge^{d}\Omega_{X}^{1^{\ast}%
}\otimes\wedge^{d-1}\Omega_{X}^{1}=\left(  \wedge^{d}\Omega_{X}^{^{1}}\right)
^{\ast}\otimes\wedge^{n-1}\Omega_{X}^{1}=\Omega_{X}^{d-1}%
\]

\item Tensoring the
\index{Euler sequence}%
Euler sequence with $\mathcal{O}_{X}$ gives%
\[
0\rightarrow\mathcal{O}_{X}\rightarrow\mathcal{O}_{X}\left(  1\right)
^{n+1}\rightarrow T_{\mathbb{P}^{n}}\mid_{X}\rightarrow0
\]
hence the long exact sequence%
\[%
\begin{tabular}
[c]{lllllll}%
$0\rightarrow$ & $H^{0}\left(  X,\mathcal{O}_{X}\right)  $ & $\rightarrow$ &
$H^{0}\left(  X,\mathcal{O}_{X}\left(  1\right)  ^{n+1}\right)  $ &
$\rightarrow$ & $H^{0}\left(  X,T_{\mathbb{P}^{n}}\mid_{X}\right)  $ &
$\rightarrow$\\
$\rightarrow$ & $H^{1}\left(  X,\mathcal{O}_{X}\right)  =0$ & $\rightarrow$ &
$H^{1}\left(  X,\mathcal{O}_{X}\left(  1\right)  ^{n+1}\right)  $ &
$\rightarrow$ & $H^{1}\left(  X,T_{\mathbb{P}^{n}}\mid_{X}\right)  $ &
$\rightarrow$\\
$\rightarrow$ & $H^{2}\left(  X,\mathcal{O}_{X}\right)  =0$ &  &  &  &  &
\end{tabular}
\
\]
so
\[
H^{0}\left(  X,T_{\mathbb{P}^{n}}\mid_{X}\right)  =\frac{H^{0}\left(
X,\mathcal{O}_{X}\left(  1\right)  ^{n+1}\right)  }{H^{0}\left(
X,\mathcal{O}_{X}\right)  }%
\]
and
\[
H^{1}\left(  X,T_{\mathbb{P}^{n}}\mid_{X}\right)  =H^{1}\left(  X,\mathcal{O}%
_{X}\left(  1\right)  ^{n+1}\right)  =H^{1}\left(  X,\mathcal{O}_{X}\left(
1\right)  \right)  ^{n+1}%
\]
By
\index{Kodaira vanishing}%
Kodaira vanishing, as $\mathcal{O}_{X}\left(  1\right)  $ is positive and
$\Omega_{X}^{3}=\mathcal{O}_{X}$ we get
\[
H^{i}\left(  X,\mathcal{O}_{X}\left(  1\right)  \right)  =H^{i}\left(
X,\mathcal{O}_{X}\left(  1\right)  \otimes\Omega_{X}^{3}\right)  =0\text{ for
}i>0
\]
hence
\[
H^{1}\left(  X,T_{\mathbb{P}^{n}}\mid_{X}\right)  =0
\]

\item The
\index{normal bundle sequence}%
normal bundle sequence%
\[
0\rightarrow T_{X}\rightarrow T_{\mathbb{P}^{n}}\mid_{X}\rightarrow
N_{X/\mathbb{P}^{n}}\rightarrow0
\]
gives the long exact sequence%
\[%

\]
hence%
\[
H^{1}\left(  X,\Omega_{X}^{d-1}\right)  \cong H^{1}\left(  X,T_{X}\right)
\cong\frac{H^{0}\left(  X,N_{X/\mathbb{P}^{n}}\right)  }{H^{0}\left(
X,T_{\mathbb{P}^{n}}\mid_{X}\right)  }%
\]

\item The long exact sequence for%
\[
0\rightarrow\mathcal{I}_{X}\rightarrow\mathcal{O}_{\mathbb{P}^{n}}%
\rightarrow\mathcal{O}_{X}\rightarrow0
\]
reads%
\[%
%
\]
so
\begin{align*}
H^{1}\left(  \mathbb{P}^{n},\mathcal{I}_{X}\right)   &  =0\\
H^{2}\left(  \mathbb{P}^{n},\mathcal{I}_{X}\right)   &  =0
\end{align*}

\item The long exact sequence for%
\[
0\rightarrow\mathcal{I}_{X}\left(  1\right)  \rightarrow\mathcal{O}%
_{\mathbb{P}^{n}}\left(  1\right)  \rightarrow\iota_{\ast}\mathcal{O}%
_{X}\left(  1\right)  \rightarrow0
\]
gives%
\[%
%
\]
hence $H^{1}\left(  \mathbb{P}^{n},\mathcal{I}_{X}\left(  1\right)  \right)
=0$ is equivalent to
\[
H^{0}\left(  \mathbb{P}^{n},\mathcal{O}_{\mathbb{P}^{n}}\left(  1\right)
\right)  \rightarrow H^{0}\left(  \mathbb{P}^{n},\iota_{\ast}\mathcal{O}%
_{X}\left(  1\right)  \right)
\]
being surjective, i.e., to $X$ being embedded by a complete linear system.

\item Tensoring the
\index{Euler sequence}%
Euler sequence with $\mathcal{I}_{X}$ gives the exact sequence%
\[
0\rightarrow\mathcal{I}_{X}\rightarrow\mathcal{I}_{X}\left(  1\right)
^{n+1}\rightarrow T_{\mathbb{P}^{n}}\otimes\mathcal{I}_{X}\rightarrow0
\]
hence the long exact sequence%
\[%
\begin{tabular}
[c]{lllllll}
& $H^{0}\left(  \mathbb{P}^{n},\mathcal{I}_{X}\right)  =0$ & $\rightarrow$ &
$H^{0}\left(  \mathbb{P}^{n},\mathcal{I}_{X}\left(  1\right)  ^{n+1}\right)  $
& $\rightarrow$ & $H^{0}\left(  \mathbb{P}^{n},T_{\mathbb{P}^{n}}%
\otimes\mathcal{I}_{X}\right)  $ & $\rightarrow$\\
$\rightarrow$ & $H^{1}\left(  \mathbb{P}^{n},\mathcal{I}_{X}\right)  =0$ &
$\rightarrow$ & $H^{1}\left(  \mathbb{P}^{n},\mathcal{I}_{X}\left(  1\right)
^{n+1}\right)  $ & $\rightarrow$ & $H^{1}\left(  \mathbb{P}^{n},T_{\mathbb{P}%
^{n}}\otimes\mathcal{I}_{X}\right)  $ & $\rightarrow$\\
$\rightarrow$ & $H^{2}\left(  \mathbb{P}^{n},\mathcal{I}_{X}\right)  =0$ &
$\rightarrow$ & $...$ &  &  &
\end{tabular}
\]
If $X$ does not lie in a hyperplane, $H^{0}\left(  \mathbb{P}^{n}%
,\mathcal{I}_{X}\left(  1\right)  \right)  =0$, so
\[
H^{0}\left(  \mathbb{P}^{n},T_{\mathbb{P}^{n}}\otimes\mathcal{I}_{X}\right)
=0
\]
If $H^{1}\left(  \mathbb{P}^{n},\mathcal{I}_{X}\left(  1\right)  \right)  =0$,
then%
\[
H^{1}\left(  \mathbb{P}^{n},T_{\mathbb{P}^{n}}\otimes\mathcal{I}_{X}\right)
=0
\]

\item Tensoring
\[
0\rightarrow\mathcal{I}_{X}\rightarrow\mathcal{O}_{\mathbb{P}^{n}}%
\rightarrow\mathcal{O}_{X}\rightarrow0
\]
with $T_{\mathbb{P}^{n}}$ gives the exact sequence%
\[
0\rightarrow\mathcal{I}_{X}\otimes T_{\mathbb{P}^{n}}\rightarrow
T_{\mathbb{P}^{n}}\rightarrow T_{\mathbb{P}^{n}}\mid_{X}\rightarrow0
\]
and the long exact sequence%
\[%
\begin{tabular}
[c]{cccccccc}%
$0$ & $\rightarrow$ & $H^{0}\left(  \mathbb{P}^{n},\mathcal{I}_{X}\otimes
T_{\mathbb{P}^{n}}\right)  $ & $\rightarrow$ & $H^{0}\left(  \mathbb{P}%
^{n},T_{\mathbb{P}^{n}}\right)  $ & $\rightarrow$ & $H^{0}\left(
\mathbb{P}^{n},T_{\mathbb{P}^{n}}\mid_{X}\right)  $ & $\rightarrow$\\
& $\rightarrow$ & $H^{1}\left(  \mathbb{P}^{n},\mathcal{I}_{X}\otimes
T_{\mathbb{P}^{n}}\right)  $ & $\rightarrow$ & $...$ &  &  &
\end{tabular}
\]
so%
\[
H^{0}\left(  \mathbb{P}^{n},T_{\mathbb{P}^{n}}\mid_{X}\right)  =H^{0}\left(
\mathbb{P}^{n},T_{\mathbb{P}^{n}}\right)  =\frac{H^{0}\left(  X,\mathcal{O}%
_{\mathbb{P}^{n}}\left(  1\right)  ^{n+1}\right)  }{H^{0}\left(
\mathbb{P}^{n},\mathcal{O}_{\mathbb{P}^{n}}\right)  }=\frac{H^{0}\left(
X,\mathcal{O}_{\mathbb{P}^{n}}\left(  1\right)  \right)  ^{n+1}}{H^{0}\left(
\mathbb{P}^{n},\mathcal{O}_{\mathbb{P}^{n}}\right)  }%
\]
hence $h^{0}\left(  \mathbb{P}^{n},T_{\mathbb{P}^{n}}\mid_{X}\right)  =\left(
n+1\right)  ^{2}-1$. Note that any element in $H^{0}\left(  \mathbb{P}%
^{n},T_{\mathbb{P}^{n}}\right)  $ can be considered as a generator of an
element in $\operatorname*{Aut}\left(  \mathbb{P}^{n}\right)  $, so
$h^{0}\left(  \mathbb{P}^{n},T_{\mathbb{P}^{n}}\right)  =\dim
\operatorname*{Aut}\left(  \mathbb{P}^{n}\right)  $.
\end{itemize}

Summarizing these observations:

\begin{proposition}
For any Calabi-Yau $d$-fold $X\subset\mathbb{P}^{n}$ with $d\geq3$ and not in
a hyperplane and with $H^{1}\left(  X,\mathcal{I}_{X}\left(  1\right)
\right)  =0$
\[
H^{1}\left(  X,\Omega_{X}^{d-1}\right)  \cong H^{1}\left(  X,T_{X}\right)
\cong\frac{H^{0}\left(  X,N_{X/\mathbb{P}^{n}}\right)  }{H^{0}\left(
X,T_{\mathbb{P}^{n}}\mid_{X}\right)  }%
\]
and%
\[
H^{0}\left(  X,T_{\mathbb{P}^{n}}\mid_{X}\right)  \cong H^{0}\left(
\mathbb{P}^{n},T_{\mathbb{P}^{n}}\right)
\]
in particular%
\[
h^{d-1,1}\left(  X\right)  =h^{0}\left(  X,N_{X/\mathbb{P}^{n}}\right)
-\dim\left(  \operatorname*{Aut}\left(  \mathbb{P}^{n}\right)  \right)
\]

\end{proposition}

\begin{remark}
Note that $H^{1}\left(  X,\mathcal{I}_{X}\left(  1\right)  \right)  =0$ if $X$
is
\index{projectively Cohen-Macaulay}%
projectively
\index{Cohen-Macaulay}%
Cohen-Macaulay. But $H^{1}\left(  X,\mathcal{I}_{X}\left(  1\right)  \right)
=0$ is also true for the
\index{Pfaffian}%
Pfaffian examples of degree $15,16$ and $17$ given in \cite{Tonoli Canonical
surfaces in mathbbP^5 and CalabiYau threefolds in mathbbP^6} (see Section
\ref{1PfaffianCalabiYauThreefolds}), which are
\index{Cohen-Macaulay}%
not
\index{projectively Cohen-Macaulay}%
projectively Cohen-Macaulay. $H^{1}\left(  X,\mathcal{I}_{X}\left(  1\right)
\right)  =0$ is equivalent to $X$ being embedded by a
\index{complete linear system}%
complete linear system.
\end{remark}

Although for $K3$ surfaces and elliptic curves we know that $h^{1,1}\left(
X\right)  =20$, respectively $h^{1,0}\left(  X\right)  =1$, it is interesting
to see how the calculation behaves:

\begin{remark}
Recall that there are no algebraic families of dimension more than $19$,
whereas all $K3$ form a $20=h^{1,1}\left(  X\right)  $-dimensional
differentiable family.

For $K3$ surfaces $T_{X}\cong\Omega_{X}^{1}$ hence $H^{0}\left(
X,T_{X}\right)  =H^{0}\left(  X,\Omega_{X}^{1}\right)  =0$, but $H^{2}\left(
X,\mathcal{O}_{X}\right)  =1$, so from the
\index{Euler sequence}%
Euler sequence tensored by $\mathcal{O}_{X}$
\[%
\begin{tabular}
[c]{ccccccc}
& $H^{1}\left(  X,\mathcal{O}_{X}\right)  =0$ & $\rightarrow$ & $H^{1}\left(
X,\mathcal{O}_{X}\left(  1\right)  ^{n+1}\right)  $ & $\rightarrow$ &
$H^{1}\left(  X,T_{\mathbb{P}^{n}}\mid_{X}\right)  $ & $\rightarrow$\\
$\rightarrow$ & $H^{2}\left(  X,\mathcal{O}_{X}\right)  $ & $\rightarrow$ &
$H^{2}\left(  X,\mathcal{O}_{X}\left(  1\right)  ^{n+1}\right)  =0$ &  &  &
\end{tabular}
\]
where $H^{2}\left(  X,\mathcal{O}_{X}\left(  1\right)  \right)  =0$ by
\index{Kodaira vanishing}%
Kodaira vanishing, so%
\[
H^{1}\left(  X,T_{\mathbb{P}^{n}}\mid_{X}\right)  =H^{2}\left(  X,\mathcal{O}%
_{X}\right)  \cong\mathbb{C}%
\]
From the
\index{normal bundle sequence}%
normal bundle sequence%
\[%
\begin{tabular}
[c]{ccccccc}
& $0=H^{0}\left(  X,T_{X}\right)  $ & $\rightarrow$ & $H^{0}\left(
X,T_{\mathbb{P}^{n}}\mid_{X}\right)  $ & $\rightarrow$ & $H^{0}\left(
X,N_{X/\mathbb{P}^{n}}\right)  $ & $\rightarrow$\\
$\rightarrow$ & $H^{1}\left(  X,T_{X}\right)  $ & $\rightarrow$ &
$H^{1}\left(  X,T_{\mathbb{P}^{n}}\mid_{X}\right)  $ & $\rightarrow$ &
$H^{1}\left(  X,N_{X/\mathbb{P}^{n}}\right)  $ & $\rightarrow$\\
$\rightarrow$ & $H^{2}\left(  X,T_{X}\right)  =0$ &  &  &  &  &
\end{tabular}
\]
and the fact that $h^{1}\left(  X,T_{X}\right)  =h^{1,1}\left(  X\right)
=20$, but the image in $H^{1}\left(  X,T_{X}\right)  $ is at most $19$
dimensional, we have $H^{1}\left(  X,N_{X/\mathbb{P}^{n}}\right)  =0$, hence
\[
h^{1,1}\left(  X\right)  =h^{1}\left(  X,T_{X}\right)  =h^{0}\left(
X,N_{X/\mathbb{P}^{n}}\right)  -h^{0}\left(  X,T_{\mathbb{P}^{n}}\mid
_{X}\right)  +1
\]
Furthermore%
\begin{align*}
H^{1}\left(  \mathbb{P}^{n},\mathcal{I}_{X}\right)   &  =0\\
H^{2}\left(  \mathbb{P}^{n},\mathcal{I}_{X}\right)   &  =0
\end{align*}
so if $H^{j}\left(  \mathbb{P}^{n},\mathcal{I}_{X}\left(  1\right)  \right)
=0$ for $j=0,1$, i.e., $X$ does not lie in a hyperplane and is embedded by a
\index{complete linear system}%
complete linear system, then also
\[
H^{j}\left(  X,T_{\mathbb{P}^{n}}\otimes\mathcal{I}_{X}\right)  =0
\]
for $j=0,1$, hence%
\[
H^{0}\left(  \mathbb{P}^{n},T_{\mathbb{P}^{n}}\mid_{X}\right)  =H^{0}\left(
\mathbb{P}^{n},T_{\mathbb{P}^{n}}\right)
\]
and we get%
\[
h^{1,1}\left(  X\right)  =h^{0}\left(  X,N_{X/\mathbb{P}^{n}}\right)  -\left(
\left(  n+1\right)  ^{2}-1\right)  +1
\]

\end{remark}

\begin{remark}
For elliptic curves $T_{X}\cong\mathcal{O}_{X}$, hence $H^{0}\left(
X,T_{X}\right)  \cong\mathbb{C}$. So from the
\index{Euler sequence}%
Euler sequence tensored by $\mathcal{O}_{X}$
\[%

\]
hence%
\[
H^{1}\left(  X,T_{\mathbb{P}^{n}}\mid_{X}\right)  =0
\]
From the
\index{normal bundle sequence}%
normal bundle sequence%
\[%
%
\
\]
one has%
\[
H^{2}\left(  \mathbb{P}^{n},\mathcal{I}_{X}\left(  1\right)  \right)  =0
\]
so if $X$ is not contained in a hyperplane and $H^{1}\left(  \mathbb{P}%
^{n},\mathcal{I}_{X}\left(  1\right)  \right)  =0$ (i.e $X$ embedded by a
\index{complete linear system}%
complete linear system), then from
\[%
%
\
\]
it follows%
\[
h^{0}\left(  \mathbb{P}^{n},T_{\mathbb{P}^{n}}\mid_{X}\right)  =h^{0}\left(
\mathbb{P}^{n},T_{\mathbb{P}^{n}}\right)  +1
\]
so%
\[
h^{1,0}\left(  X\right)  =h^{0}\left(  X,N_{X/\mathbb{P}^{n}}\right)  -\left(
\left(  n+1\right)  ^{2}-1\right)
\]

\end{remark}

\subsection{Batyrev%
\'{}%
s Hodge formula\label{Sec Batyrevs Hodge formula}}

Let $\Delta\subset M_{\mathbb{R}}$ be a reflexive polytope and $X$ a general
anticanonical hypersurface in $Y=\mathbb{P}\left(  \Delta\right)  $. To
\index{Batyrev}%
given an idea on the proof of the Equations \ref{4BatyrevHodgeformula}%
\begin{align*}
h^{d-1,1}\left(  \bar{X}\right)   &  =\left\vert \Delta\cap M\right\vert
-n-1-\sum_{Q\text{ facet of }\Delta}\left\vert \operatorname*{int}%
\nolimits_{M}\left(  Q\right)  \right\vert \\
&  +\sum_{\substack{Q\text{ face of }\Delta\\\operatorname*{codim}%
Q=2}}\left\vert \operatorname*{int}\nolimits_{M}\left(  Q\right)  \right\vert
\cdot\left\vert \operatorname*{int}\nolimits_{N}\left(  Q^{\ast}\right)
\right\vert \\
h^{1,1}\left(  \bar{X}\right)   &  =\left\vert \Delta^{\ast}\cap M\right\vert
-n-1-\sum_{Q^{\ast}\text{ facet of }\Delta^{\ast}}\left\vert
\operatorname*{int}\nolimits_{N}\left(  Q^{\ast}\right)  \right\vert \\
&  +\sum_{\substack{Q^{\ast}\text{ face of }\Delta^{\ast}%
\\\operatorname*{codim}Q^{\ast}=2}}\left\vert \operatorname*{int}%
\nolimits_{N}\left(  Q^{\ast}\right)  \right\vert \cdot\left\vert
\operatorname*{int}\nolimits_{M}\left(  Q\right)  \right\vert
\end{align*}
via MPCP desingularizations, suppose $\bar{\Sigma}$ is a
\index{maximal projective subdivision}%
maximal projective subdivision of the
\index{normal fan}%
normal fan $\operatorname*{NF}\left(  \Delta\right)  \subset N_{\mathbb{R}}$
of the reflexive polytope $\Delta\subset M_{\mathbb{R}}$, let%
\[
f:X\left(  \bar{\Sigma}\right)  \rightarrow\mathbb{P}\left(  \Delta\right)
\]
be the corresponding birational morphism inducing a
\index{crepant resolution}%
crepant morphism $\bar{X}\rightarrow X$, and write $D_{w}$ with $w\in
\bar{\Sigma}\left(  1\right)  $ for the prime $T$-Weil divisors on $X\left(
\bar{\Sigma}\right)  $.

\subsubsection{Toric divisor classes}

Restriction of divisors
\index{Chow group}%
from $X\left(  \bar{\Sigma}\right)  $ to $\bar{X}$ gives%
\begin{equation}%
\begin{tabular}
[c]{cccccccc}%
$0\rightarrow$ & $M$ & $\rightarrow$ & $\operatorname*{WDiv}_{T}\left(
X\left(  \bar{\Sigma}\right)  \right)  $ & $\rightarrow$ &  & $A_{n-1}\left(
X\left(  \bar{\Sigma}\right)  \right)  $ & $\rightarrow0$\\
& $\shortparallel$ &  & $\downarrow$ &  &  & $\downarrow$ & \\
$0\rightarrow$ & $M$ & $\rightarrow$ & $\operatorname*{WDiv}_{T}\left(
\bar{X}\right)  $ & $\rightarrow$ & $A_{d-1}\left(  \bar{X}\right)  _{toric}$
& $\subset A_{d-1}\left(  \bar{X}\right)  $ &
\end{tabular}
\label{4chowrestriction}%
\end{equation}
The image of the
\index{T-Weil divisor}%
toric Weil divisors $\operatorname*{WDiv}_{T}\left(  \bar{X}\right)  $ in
$A_{d-1}\left(  \bar{X}\right)  $ is not surjective in general, so denote the
image by $A_{d-1}\left(  \bar{X}\right)  _{toric}$ and its
\newsym[$H_{toric}^{1,1}$]{toric divisor classes}{}complexification, i.e., the
subspace of $H^{1,1}\left(  \bar{X}\right)  $ of
\index{toric divisor classes}%
toric divisor classes of $\bar{X}$, by%
\[
H_{toric}^{1,1}\left(  \bar{X}\right)  =A_{d-1}\left(  \bar{X}\right)
_{toric}\otimes\mathbb{C}%
\]

A divisor has trivial restriction if and only if its support is disjoint from
the general hypersurface $\bar{X}$. If $w\in\bar{\Sigma}\left(  1\right)  $ is
a lattice point in the relative interior of a facet of $\Delta^{\ast}\subset
N_{\mathbb{R}}$, i.e., if
\[
w\in\bigcup_{\operatorname*{codim}\left(  Q^{\ast}\right)  =1}%
\operatorname*{int}\left(  Q^{\ast}\right)
\]
then $f\left(  D_{w}\right)  $ is a point, so $D_{w}$ is disjoint from any
general element $\bar{X}$ of $\left\vert -K_{X\left(  \bar{\Sigma}\right)
}\right\vert $. If $w\in\bar{\Sigma}\left(  1\right)  $ is not in the relative
interior of a facet, then $\dim\left(  f\left(  D_{w}\right)  \right)  >0$ so
$f\left(  D_{w}\right)  $ meets $X$. Hence with%
\[
\Xi_{0}^{\ast}=\Delta^{\ast}\cap N-\bigcup_{\operatorname*{codim}Q^{\ast}%
\leq1}\operatorname*{int}\nolimits_{N}\left(  Q^{\ast}\right)
\]
we have%
\[
\operatorname*{WDiv}\nolimits_{T}\left(  \bar{X}\right)  \cong\mathbb{Z}%
^{\Xi_{0}^{\ast}}%
\]
and as cokernel of $M\rightarrow\operatorname*{WDiv}_{T}\left(  \bar
{X}\right)  $%
\[
A_{d-1}\left(  \bar{X}\right)  _{toric}\cong\mathbb{Z}^{\Xi_{0}^{\ast}}/M
\]
so%
\[
H_{toric}^{1,1}\left(  \bar{X}\right)  \cong\mathbb{Z}^{\Xi_{0}^{\ast}}/M
\]
with dimension%
\[
h_{toric}^{1,1}\left(  \bar{X}\right)  =\left\vert \Delta^{\ast}\cap
N\right\vert -1-\sum_{\operatorname*{codim}Q^{\ast}=1}\left\vert
\operatorname*{int}\left(  Q^{\ast}\right)  \right\vert -n
\]

\subsubsection{Polynomial deformations and complex moduli space}

Define the \newsym[$H_{poly}^{d-1,1}$]{polynomial deformations}{}subspace of
\index{polynomial deformation}%
polynomial first order
\index{deformation}%
deformations
\[
H_{poly}^{d-1,1}\left(  \bar{X}\right)  \subset H^{d-1,1}\left(  \bar
{X}\right)  \cong H^{1}\left(  \bar{X},T_{\bar{X}}\right)
\]
as the subspace determined by $\left\vert -K_{X\left(  \bar{\Sigma}\right)
}\right\vert $. Any element is given by a linear combination of the lattice
monomials $\Delta\cap M$. Multiplication of the equation by a constant does
not affect the zero set and the
\index{automorphism}%
automorphism group of $X\left(  \bar{\Sigma}\right)  $ has dimension%
\[
\dim\left(  \operatorname*{Aut}\left(  X\left(  \bar{\Sigma}\right)  \right)
\right)  =n+\sum_{\operatorname*{codim}Q=1}\left\vert \operatorname*{int}%
\nolimits_{M}\left(  Q\right)  \right\vert
\]
hence%
\[
h_{poly}^{d-1,1}\left(  \bar{X}\right)  =\left\vert \Delta\cap M\right\vert
-1-n-\sum_{\operatorname*{codim}Q=1}\left\vert \operatorname*{int}%
\nolimits_{M}\left(  Q\right)  \right\vert
\]

The tangent space to the space of polynomial deformations is%
\[
H_{poly}^{d-1,1}\left(  \bar{X}\right)  \cong\left(  \mathbb{Z}^{\Xi_{0}%
}/N\right)  \otimes\mathbb{C}%
\]
with%
\[
\Xi_{0}=\Delta\cap M-\bigcup_{\operatorname*{codim}Q\leq1}\operatorname*{int}%
\nolimits_{M}\left(  Q\right)
\]
and $\mathbb{Z}^{\Xi_{0}}/N$ given as the cokernel of the lower row in%
\[%
\begin{tabular}
[c]{cccccccc}%
$0\rightarrow$ & $N$ & $\rightarrow$ & $WDiv_{T}\left(  X\left(  \bar{\Sigma
}^{\ast}\right)  \right)  $ & $\rightarrow$ &  & $A_{n-1}\left(  X\left(
\bar{\Sigma}^{\ast}\right)  \right)  $ & $\rightarrow0$\\
& $\shortparallel$ &  & $\downarrow$ &  &  & $\downarrow$ & \\
$0\rightarrow$ & $N$ & $\rightarrow$ & $WDiv_{T}\left(  \bar{X}^{\circ
}\right)  =\mathbb{Z}^{\Xi_{0}}$ & $\rightarrow$ & $A_{d-1}\left(  \bar
{X}^{\circ}\right)  _{toric}$ & $\subset A_{d-1}\left(  \bar{X}^{\circ
}\right)  $ &
\end{tabular}
\
\]

For a description of the non-toric divisor classes and non-polynomial
\index{deformation}%
deformations see, e.g., \cite[Sec. 4.1]{CK Mirror Symmetry and Algebraic
Geometry}.

\subsection{First approximation of a tropical Hodge
formula\label{Sec First approximation of a tropical Hodge formula}}

Let $Y=\mathbb{P}\left(  \Delta\right)  =\mathbb{P}_{\mathbb{C}}^{n}$ for the
degree $n+1$ reflexive Veronese simplex $\Delta$ and denote by $S$ the
homogeneous coordinate ring of $Y$. Consider the setup of Section
\ref{Sec tropical mirror construction}: So let $\mathfrak{X}\subset
Y\times\operatorname*{Spec}\mathbb{C}\left[  t\right]  $ be a
\index{degeneration}%
degeneration of projective Calabi-Yau varieties defined by the ideal
$I\subset\mathbb{C}\left[  t\right]  \otimes S$ with monomial special fiber
given by $I_{0}\subset S$, general fiber $X$ and satisfying the conditions
given in Section \ref{genericy condition}. So $X$ has only unobstructed
polynomial deformations and as $Y$ is assumed to be projective space $I_{0}$
is a Stanley-Reisner ideal.

\begin{proposition}
A $T$-invariant basis of $H^{0}\left(  X,N_{X_{0}/\mathbb{P}^{n}}\right)  $ is
given
\index{dual}%
by
\[
A\left(  \operatorname*{supp}\left(  \operatorname*{dual}\left(  B\left(
I\right)  \right)  \right)  \cap M\right)
\]
and%
\[%
\begin{tabular}
[c]{clcll}%
$M$ &  & $\mathbb{Z}^{\Sigma\left(  1\right)  }$ &  & $\operatorname*{Hom}%
\left(  I_{0},S/I_{0}\right)  _{0}$\\
$\cup$ &  & $\cup$ &  & \multicolumn{1}{c}{$\cup$}\\
\multicolumn{1}{l}{$\operatorname*{supp}\left(  \operatorname*{dual}\left(
B\left(  I\right)  \right)  \right)  \cap M$} & $\underset{1:1}{\overset
{A}{\rightleftarrows}}$ & \multicolumn{1}{l}{$A\left(  \operatorname*{supp}%
\left(  \operatorname*{dual}\left(  B\left(  I\right)  \right)  \right)  \cap
M\right)  $} & $\subset$ & $H^{0}\left(  X_{0},N_{X_{0}/\mathbb{P}^{n}%
}\right)  $%
\end{tabular}
\]
in particular $h^{0}\left(  X_{0},N_{X_{0}/\mathbb{P}^{n}}\right)  =\left\vert
\operatorname*{supp}\left(  \operatorname*{dual}\left(  B\left(  I\right)
\right)  \right)  \cap M\right\vert $ is the number of lattice points of
$\operatorname*{dual}\left(  B\left(  I\right)  \right)  $.
\end{proposition}

\begin{corollary}
If $X$ is a Calabi-Yau manifold, then%
\[
h^{1,\dim X-1}\left(  X\right)  =\left\vert \operatorname*{supp}\left(
\operatorname*{dual}\left(  B\left(  I\right)  \right)  \right)  \cap
M-\operatorname*{Roots}\left(  \mathbb{P}\left(  \Delta\right)  \right)
\right\vert -\overset{=n}{\dim\left(  T_{Y}\right)  }+\overset{K3}{1}%
\]

\end{corollary}

\begin{figure}
[h]
\begin{center}
\includegraphics[
height=3.0865in,
width=2.9464in
]%
{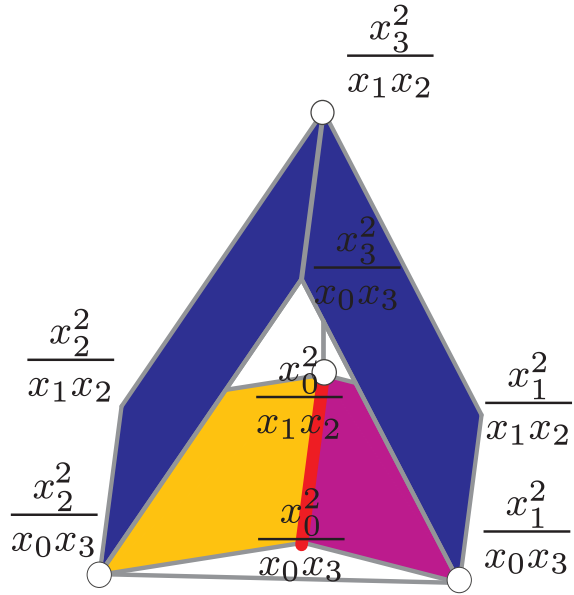}%
\caption{The dual complex and the monomials corresponding to the vertices of
$\nabla^{\ast}$ associated to the degeneration of the general complete
intersection elliptic curve in $\mathbb{P}^{3}$}%
\label{Fig dualT1}%
\end{center}
\end{figure}

\begin{example}
For the elliptic curve given as a complete intersection of two quadrics in
$\mathbb{P}^{3}$, as considered in Example \ref{22Nabla}, the
\index{dual}%
dual
\index{dual complex}%
complex $\operatorname*{dual}\left(  B\left(  I\right)  \right)  $ together
with the monomials corresponding to vertices of $\nabla^{\ast}$ is shown in
Figure \ref{Fig dualT1}. The $4$ lattice points of $\operatorname*{dual}%
\left(  B\left(  I\right)  \right)  $,
\index{dual complex}%
marked
\index{dual}%
with dots, form a
\index{first order deformation}%
basis
\index{deformation}%
of $T_{X_{0}}^{1}$, the remaining $12$ lattice points are
\index{root}%
roots, i.e., \newsym[$T_{X_{0}}^{1}$]{first order deformations}{}homomorphism
of the form $x_{i}\frac{\partial}{\partial x_{j}}\in\operatorname*{Hom}\left(
I_{0},S/I_{0}\right)  _{0}$ for $i\neq j$, of $\mathbb{P}\left(
\Delta\right)  =\mathbb{P}^{n}$. So with the
\index{torus}%
torus $T$ of $\mathbb{P}\left(  \Delta\right)  $ we have%
\begin{align*}
\dim\left(  T_{X_{0}}^{1}\right)  +\left\vert \operatorname*{Roots}\left(
\mathbb{P}\left(  \Delta\right)  \right)  \right\vert  &  =\\
h^{1,0}\left(  X\right)  +\dim\left(  T\right)  +\left\vert
\operatorname*{Roots}\left(  \mathbb{P}\left(  \Delta\right)  \right)
\right\vert  &  =\\
h^{1,0}\left(  X\right)  +\dim\left(  \operatorname*{Aut}\left(
\mathbb{P}\left(  \Delta\right)  \right)  \right)   &  =h^{0}\left(
X_{0},N_{X_{0}/\mathbb{P}^{n}}\right)
\end{align*}
The $h^{1,0}\left(  X\right)  =1$ dimensional tangent space to the moduli
space of $X$ is a quotient of the $4$-dimensional $T_{X_{0}}^{1}$ by the
$3$-dimensional torus $T$ of $Y$.
\end{example}

%

\begin{remark}
On the other hand lattice points of $\operatorname*{dual}\left(  B\left(
I\right)  \right)  $ correspond to rays of the MPCP-desingularization of the
toric variety $Y^{\circ}$ containing the Batyrev-Borisov mirror $X^{\circ}$ of
$X$, hence rays correspond to toric divisor classes of $X^{\circ}$, so we also
have an interpretation of the formula%
\[
h_{toric}^{1,1}\left(  X^{\circ}\right)  =\left\vert \operatorname*{supp}%
\left(  \operatorname*{dual}\left(  B\left(  I\right)  \right)  \right)  \cap
M-\operatorname*{Roots}\left(  \mathbb{P}\left(  \Delta\right)  \right)
\right\vert -\overset{=n}{\dim\left(  T_{Y}\right)  }%
\]

\end{remark}

Note that above formula agrees with the toric Batyrev formula for
\index{hypersurface}%
hypersurfaces%
\[
h^{1,\dim\left(  X\right)  -1}\left(  X\right)  =\left\vert \Delta\cap
M\right\vert -n-1-\sum_{\Gamma\text{ facet of }\Delta}\left\vert
\operatorname*{int}\nolimits_{M}\left(  \Gamma\right)  \right\vert
\]
as for any
\index{reflexive}%
reflexive polytope $\Delta$%
\[
\left\vert \partial\Delta\cap M\right\vert =\left\vert \Delta\cap M\right\vert
-1
\]
for any
\index{simplicial}%
simplicial polytope
\[
\dim\left(  \operatorname*{Aut}\left(  \mathbb{P}\left(  \Delta\right)
\right)  \right)  =n+\sum_{\Gamma\text{ facet of }\Delta}\left\vert
\operatorname*{int}\nolimits_{M}\left(  \Gamma\right)  \right\vert
\]
and the faces of the dual $\Delta^{\ast}$ of a Veronese polytope do not
contain any interior lattice points, hence%
\[
\sum_{\substack{Q\text{ face of }\Delta\\\operatorname*{codim}Q=2}}\left\vert
\operatorname*{int}\nolimits_{M}\left(  Q\right)  \right\vert \cdot\left\vert
\operatorname*{int}\nolimits_{N}\left(  Q^{\ast}\right)  \right\vert =0
\]

\begin{figure}
[h]
\begin{center}
\includegraphics[
height=2.6394in,
width=2.8279in
]%
{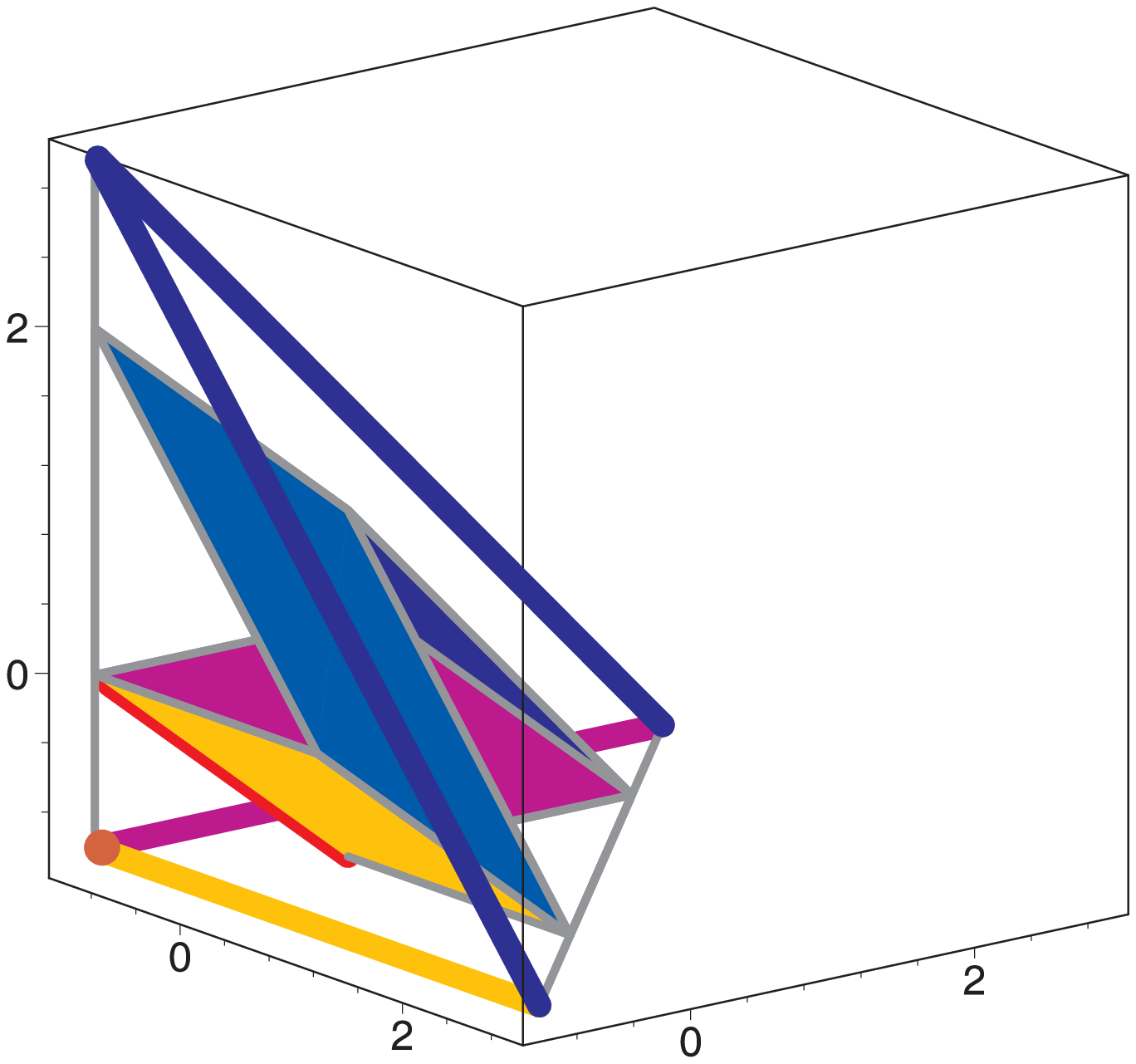}%
\caption{The complexes $\operatorname*{dual}\left(  B\left(  I\right)
\right)  \subset\operatorname*{dual}\left(  \operatorname*{Poset}\left(
\nabla\right)  \right)  $ and $\Delta$ for the degeneration of the complete
intersection of two general quadrics in $\mathbb{P}^{3}$}%
\label{Fig dual complex dual Nabla Delta for 22}%
\end{center}
\end{figure}

\begin{remark}
The lattice points of $\operatorname*{dual}\left(  B\left(  I\right)  \right)
$
\index{dual complex}%
corresponding
\index{dual}%
to
\index{root}%
roots of $\mathbb{P}\left(  \Delta\right)  $ are the lattice points of
$\operatorname*{supp}\left(  \operatorname*{dual}\left(  B\left(  I\right)
\right)  \right)  \subset\nabla^{\ast}\subset\Delta$ in the relative interior
of the facets of $\Delta$. The complex $\operatorname*{dual}\left(  B\left(
I\right)  \right)  \subset\operatorname*{dual}\left(  \operatorname*{Poset}%
\left(  \nabla\right)  \right)  $ and $\Delta$ are shown in Figure
\ref{Fig dual complex dual Nabla Delta for 22}.
\end{remark}

%

\subsection{String cohomology\label{Sec string cohomology}}

\subsubsection{Stringy $E$-function for toric
varieties\label{Sec Stringy Efunction of toric variety}}

Let $X$ be a normal $\mathbb{Q}$-Gorenstein toric variety of dimension $n$,
given by the rational polyhedral fan $\Sigma\subset N_{\mathbb{R}}$ and let
$\varphi_{K_{Y}}:N_{\mathbb{R}}\rightarrow\mathbb{R}_{\geq0}$ be the
continuous piecewise linear function with $\varphi_{K_{Y}}\left(  \hat
{r}\right)  =1$ for the minimal lattice generators $\hat{r}$ of all rays
$r\in\Sigma\left(  1\right)  $.

\begin{theorem}
\cite{Batyrev Stringy Hodge numbers of varieties with Gorenstein canonical
singularities. Integrable systems and algebraic geometry KobeKyoto 1997} The
stringy $E$-function of the normal $\mathbb{Q}$-Gorenstein toric variety $X$
of dimension $n$ is given by%
\[
E_{st}\left(  X,u,v\right)  =\left(  uv-1\right)  ^{n}\sum_{\sigma\in\Sigma
}\sum_{n\in N\cap\operatorname*{int}\left(  \sigma\right)  }\left(  uv\right)
^{-\varphi_{K_{Y}\left(  n\right)  }}%
\]

\end{theorem}

Recall that $\operatorname*{int}\left(  \sigma\right)  $ denotes the relative
interior of $\sigma$. For the $0$-cone we define $\operatorname*{int}\left(
0\right)  =\left\{  0\right\}  $.

\subsubsection{The combinatorics of
posets\label{Sec the combinatorics of posets}}

Recall that a
\index{poset|textbf}%
\textbf{poset} $P$ is a finite partially ordered set, i.e., a finite set $P$
with a
\index{reflexive}%
reflexive, antisymmetric ($x\leq y$ and $y\leq x$ implies $x=y$) and
transitive relation $\leq$.

\begin{lemma}
\cite[Sec. 2]{BB Mirror duality and stringtheoretic Hodge numbers} There is a
unique \newsym[$\mu_{P}$]{M\"{o}bius function}{}function%
\[
\mu_{P}:P\times P\rightarrow\mathbb{Z}%
\]
called
\index{M\"{o}bius function|textbf}%
\textbf{M\"{o}bius function}, such that for every function $f:P\rightarrow A$
to some abelian group $A$ and%
\[
g\left(  y\right)  =\sum_{x\leq y}f\left(  x\right)
\]
it holds%
\[
f\left(  y\right)  =\sum_{x\leq y}\mu_{P}\left(  x,y\right)  g\left(
x\right)
\]

\end{lemma}

\begin{definition}
Suppose that $P$ has a unique minimal element $\min\left(  P\right)  $ and
maximal element $\max\left(  P\right)  $ and that any maximal chain in $P$ has
the same length $d$. If $x\leq y$, then define
\[
\left[  x,y\right]  =\left\{  z\in P\mid x\leq z\leq y\right\}
\]

The
\index{rank function|textbf}%
\textbf{rank function} $\rho:P\rightarrow\left\{  0,...,d\right\}  $
associates to any $x\in P$ the length of any maximal chain in $\left[
\min\left(  P\right)  ,x\right]  $.
\end{definition}

\begin{definition}
A
\index{poset}%
poset $P$ with above properties is called
\index{Eulerian poset|textbf}%
\textbf{Eulerian} if its M\"{o}bius function satisfies%
\[
\mu_{P}\left(  x,y\right)  =\left(  -1\right)  ^{\rho\left(  y\right)
-\rho\left(  x\right)  }%
\]
for all $x\leq y$.
\end{definition}

\begin{lemma}
If $P$ is an Eulerian poset and $\left[  x,y\right]  \subset P$, then also
$\left[  x,y\right]  $ is in an Eulerian poset with rank function
\[%
\begin{tabular}
[c]{lll}%
$\left[  x,y\right]  $ & $\rightarrow$ & $\left\{  0,...,\rho\left(  y\right)
-\rho\left(  x\right)  \right\}  $\\
$z$ & $\mapsto$ & $\rho\left(  z\right)  -\rho\left(  x\right)  $%
\end{tabular}
\ \
\]

\end{lemma}

\begin{lemma}
Reversing the partial order, every Eulerian
\index{poset}%
poset $P$ has a
\index{dual poset|textbf}%
\textbf{dual poset} $P^{\ast}$, which is again Eulerian with rank function%
\[
\rho^{\ast}\left(  x\right)  =\rho\left(  P\right)  -\rho\left(  x\right)
\]

\end{lemma}

\begin{example}
For any $n$-dimensional
\index{strongly convex}%
strongly convex rational polyhedral cone $C\subset N_{\mathbb{R}}$, the set of
faces of $C$, together with inclusion, forms an Eulerian poset
$\operatorname*{Poset}\left(  C\right)  $ with rank function%
\[%
\begin{tabular}
[c]{llll}%
$\rho:$ & $\operatorname*{Poset}\left(  C\right)  $ & $\rightarrow$ &
$\left\{  0,...,\dim\left(  C\right)  \right\}  $\\
& $F$ & $\mapsto$ & $\dim\left(  F\right)  $%
\end{tabular}
\]
and with minimal respectively maximal element%
\begin{align*}
\min\left(  P\right)   &  =\left\{  0\right\} \\
\max\left(  P\right)   &  =C
\end{align*}

The dual poset of $P$ is the poset of the faces of the dual cone $\check
{C}\subset M_{\mathbb{R}}$.
\end{example}

Define \newsym[$\tau_{<s}$]{truncation operator}{}the
\index{truncation operator|textbf}%
\textbf{truncation operator} by%
\begin{align*}
\tau_{<s}  &  :\mathbb{Z}\left[  t\right]  \rightarrow\mathbb{Z}\left[
t\right] \\
\tau_{<s}\left(  \sum_{i=0}^{d}a_{i}t^{i}\right)   &  =\sum
_{\substack{i=0\\i<s}}^{d}a_{i}t^{i}%
\end{align*}

\begin{definition}
If $P$ is an Eulerian poset of rank $d$, then define the polynomials $G\left(
P,t\right)  ,H\left(  P,t\right)  \in\mathbb{Z}\left[  t\right]  $ recursively
by%
\begin{align*}
G\left(  P,t\right)   &  =1\\
H\left(  P,t\right)   &  =1
\end{align*}
for $d=0$ and%
\begin{align*}
H\left(  P,t\right)   &  =\sum_{\substack{x\in P\\x>\min\left(  P\right)
}}\left(  t-1\right)  ^{\rho\left(  x\right)  -1}G\left(  \left[  x,P\right]
,t\right) \\
G\left(  P,t\right)   &  =\tau_{<\frac{d}{2}}\left(  \left(  1-t\right)
H\left(  P,t\right)  \right)
\end{align*}
for $d>0$.
\end{definition}

\begin{example}
\label{Ex deg5 veronese}Suppose $P$ is the poset of the faces of a cone over
the
\index{reflexive}
degree $5$
\index{Veronese}%
Veronese simplex of $\mathbb{P}^{4}$, then%
\begin{align*}
H\left(  P,t\right)   &  =1+t+t^{2}+t^{3}+t^{4}\\
G\left(  P,t\right)   &  =1
\end{align*}
indeed for any boolean algebra $P$ of rank $n$, we have $H\left(  P,t\right)
=1+t+...+t^{n-1}$ and $G\left(  P,t\right)  =1$.
\end{example}

\subsubsection{String-theoretic Hodge formula for
hypersurfaces\label{Sec String theoretic Hodge formula for hypersurfaces}}

\begin{definition}
Suppose $N=\mathbb{Z}^{n}$ and $M=\operatorname*{Hom}\left(  N,\mathbb{Z}%
\right)  $. A cone $C$ of dimension $d\geq1$ in $M_{\mathbb{R}}$ is called
\index{Gorenstein cone|textbf}%
\textbf{Gorenstein cone} if there is a $w\in N$ with $\left\langle
m,w\right\rangle >0$ for all $0\neq m\in C$ and%
\[
\left\{  m\in C\mid\left\langle m,w\right\rangle =1\right\}
\]
is a rational convex polyhedron, called the \textbf{supporting polyhedron} of
$C$.
\end{definition}

\begin{remark}
Consider the setup of Section
\ref{Sec tropical mirror construction for complete intersections}, so let $Y$
be a toric Fano variety, $\mathfrak{X}$ a degeneration of Calabi-Yau
varieties, given by the ideal $I\subset\mathbb{C}\left[  t\right]  \otimes S$
and with monomial special fiber $I_{0}\subset S$. Applying the tropical mirror
construction, we obtain the strongly convex polyhedral cone
\[%
\begin{tabular}
[c]{lcc}%
$C_{I_{0}}\left(  I\right)  $ & $\subset$ & $N_{\mathbb{R}}\oplus\mathbb{R}$\\
&  & $\cup$\\
&  & $N\oplus\mathbb{Z}$%
\end{tabular}
\
\]
which is the closure of the set of weight vectors selecting $I_{0}$ as initial
ideal of $I$. Then the dual cone $C_{I_{0}}\left(  I\right)  ^{\vee}$ of
$C_{I_{0}}\left(  I\right)  $ is a Gorenstein cone.

If $\mathfrak{X}$ is a degeneration of complete intersections in a Gorenstein
toric Fano $Y=\mathbb{P}\left(  \Delta\right)  $, then also $C_{I_{0}}\left(
I\right)  $ is a Gorenstein cone with reflexive supporting polytope
$\nabla=C_{I_{0}}\left(  I\right)  \cap\left\{  w_{t}=1\right\}  $ and%
\[
Y^{\circ}=\mathbb{P}\left(  \nabla\right)  =\operatorname*{Proj}%
\mathbb{C}\left[  C_{I_{0}}\left(  I\right)  \right]
\]
with the natural grading on $\mathbb{C}\left[  C_{I_{0}}\left(  I\right)
\right]  $.
\end{remark}

\begin{example}
\label{Ex Gorenstein cone deg 5 veronese}The cone
\[
C=\left\{  \left(  \lambda,\lambda m\right)  \in\left(  \mathbb{Z}\oplus
M\right)  _{\mathbb{R}}\mid\lambda\in\mathbb{R}_{\geq0},\text{ }m\in
\Delta\right\}
\]
where $\Delta$ is the degree $5$ Veronese polytope is a Gorenstein cone over
the reflexive polyhedron $\Delta$.
\end{example}

\begin{definition}
Let $C$ be a
\newsym[$P_{\Delta}\left(  t\right)  $]{Erhard power series}{}Gorenstein cone
in $M_{\mathbb{R}}$ and $\Delta$ its supporting polyhedron. The
\index{Erhard power series|textbf}%
\textbf{Erhard power series} of $\Delta$ is%
\[
P_{\Delta}\left(  t\right)  =\sum_{k=0}^{\infty}\left\vert k\Delta\cap
M\right\vert \cdot t^{k}%
\]

\end{definition}

\begin{lemma}
\cite{Batyrev Dual polyhedra and mirror symmetry for CalabiYau hypersurfaces
in toric varieties} Let $C$ be a Gorenstein cone of dimension $d$ in
$M_{\mathbb{R}}$ and $\Delta$ its supporting polyhedron. Then there are
$\psi_{0},...,\psi_{d}\in\mathbb{Z}_{\geq0}$ such that%
\[
P_{\Delta}\left(  t\right)  =\frac{\psi_{0}+\psi_{1}\cdot t...+\psi_{d-1}\cdot
t^{d-1}}{\left(  1-t\right)  ^{d}}%
\]
Define
\[
S\left(  C,t\right)  =\psi_{0}+\psi_{1}\cdot t+...+\psi_{d-1}\cdot t^{d-1}%
\]

\end{lemma}

\begin{remark}
Note that%
\[
S\left(  C,t\right)  =\psi_{0}+\psi_{1}\cdot t+...+\psi_{d-1}\cdot
t^{d-1}=\left(  1-t\right)  ^{d}\cdot\sum_{k=0}^{\infty}\left\vert k\Delta\cap
M\right\vert \cdot t^{k}%
\]
depends only on the values $\left\vert k\Delta\cap M\right\vert $ for
$k=0,...,d-1$, because of the recursion relation $\left(  1-t\right)  ^{d}$.
\end{remark}

\begin{example}
For the Gorenstein cone $C$ over the degree $5$ Veronese polyhedron $\Delta$
as defined in Example \ref{Ex Gorenstein cone deg 5 veronese} we have%
\[%
\begin{tabular}
[c]{llllll}%
$k$ & $0$ & $1$ & $2$ & $3$ & $4$\\
$\left\vert k\Delta\cap M\right\vert $ & $1$ & $126$ & $1001$ & $3876$ &
$10626$%
\end{tabular}
\
\]
hence%
\[
S\left(  C,t\right)  =1+121t+381t^{2}+121t^{3}+t^{4}%
\]

\end{example}

\begin{definition}
If $C$ is a Gorenstein cone, define%
\[
\tilde{S}\left(  C,t\right)  =\sum_{C_{1}\text{ face of }C}S\left(
C_{1},t\right)  \left(  -1\right)  ^{\dim\left(  C\right)  -\dim\left(
C_{1}\right)  }G\left(  \left[  C_{1},C\right]  ,t\right)
\]

\end{definition}

\begin{example}
If $C$ is the Gorenstein cone over the Veronese polyhedron of degree $4$ of
$\mathbb{P}^{3}$, we have for the faces $C_{1}\subset C$%
\[%
\begin{tabular}
[c]{llll}%
$\dim\left(  C_{1}\right)  $ &
\begin{tabular}
[c]{l}%
number of faces of $C$\\
of this dimension
\end{tabular}
& $S\left(  C_{1},t\right)  $ & $G\left(  \left[  C_{1},C\right]  ,t\right)
$\\
$0$ & $1$ & $1$ & $1$\\
$1$ & $4$ & $1$ & $1$\\
$2$ & $6$ & $1+3t$ & $1$\\
$3$ & $4$ & $1+12t+3t^{2}$ & $1$\\
$4$ & $1$ & $1+31t+31t^{2}+t^{3}$ & $1$%
\end{tabular}
\]
hence%
\[
\tilde{S}\left(  C,t\right)  =t+19t^{2}+t^{3}%
\]

\end{example}

\begin{example}
If $C$ is the Gorenstein cone over the degree $5$ Veronese polyhedron from
Example \ref{Ex Gorenstein cone deg 5 veronese}, the $S\left(  C_{1},t\right)
$ for all faces $C_{1}\subset C$ is as follows%
\[%
\begin{tabular}
[c]{llll}%
$\dim\left(  C_{1}\right)  $ &
\begin{tabular}
[c]{l}%
number of faces of $C$\\
of this dimension
\end{tabular}
& $S\left(  C_{1},t\right)  $ & $G\left(  \left[  C_{1},C\right]  ,t\right)
$\\
$0$ & $1$ & $1$ & $1$\\
$1$ & $5$ & $1$ & $1$\\
$2$ & $10$ & $1+4t$ & $1$\\
$3$ & $10$ & $1+18t+6t^{2}$ & $1$\\
$4$ & $5$ & $1+52t+68t^{2}+4t^{3}$ & $1$\\
$5$ & $1$ & $1+121t+381t^{2}+121t^{3}+t^{4}$ & $1$%
\end{tabular}
\]
hence%
\[
\tilde{S}\left(  C,t\right)  =t+101t^{2}+101t^{3}+t^{4}%
\]

\end{example}

\begin{theorem}
\label{thm BB stringy E formula}\cite{BB Mirror duality and stringtheoretic
Hodge numbers} Let $C$ be a Gorenstein cone supported on a reflexive
polyhedron $\Delta$. If $X$ is an ample nondegenerate Calabi-Yau hypersurface
of dimension $d$ in $\mathbb{P}\left(  \Delta\right)  =\operatorname*{Proj}%
\mathbb{C}\left[  C\right]  $, then%
\begin{align*}
E_{st}\left(  X,u,v\right)   &  =\left(  uv\right)  ^{-1}\left(  -u\right)
^{\dim\left(  C\right)  }\tilde{S}\left(  C,u^{-1}v\right)  +\left(
uv\right)  ^{-1}\tilde{S}\left(  C^{\vee},uv\right) \\
&  +\left(  uv\right)  ^{-1}\sum_{0\subsetneq C_{1}\subsetneq C}\left(
-u\right)  ^{\dim\left(  C_{1}\right)  }\tilde{S}\left(  C_{1},u^{-1}v\right)
\tilde{S}\left(  C_{1}^{\vee},uv\right)
\end{align*}
where the sum goes over the faces $C_{1}$ of $C$.
\end{theorem}

We may write this formula as%
\[
E_{st}\left(  X,u,v\right)  =\left(  uv\right)  ^{-1}\sum_{C_{1}\subset
C}\left(  -u\right)  ^{\dim\left(  C_{1}\right)  }\tilde{S}\left(
C_{1},u^{-1}v\right)  \tilde{S}\left(  C_{1}^{\vee},uv\right)
\]

\begin{corollary}
Consider the setup of Theorem \ref{thm BB stringy E formula}. If $X^{\circ}$
is an ample nondegenerate Calabi-Yau hypersurface of dimension $d$ in
$\mathbb{P}\left(  \Delta^{\ast}\right)  =\operatorname*{Proj}\mathbb{C}%
\left[  K^{\ast}\right]  $, then the stringy $E$-functions of $X$ and
$X^{\circ}$ satisfy the mirror duality relation
\[
E_{st}\left(  X;u,v\right)  =\left(  -u\right)  ^{d}E_{st}\left(  X^{\circ
};u^{-1},v\right)
\]

\end{corollary}

\begin{example}
For the quadric $K3$ surface in $\mathbb{P}^{3}$ given by a general section in
the degree $4$ Veronese polytope, we obtain%
\[%
\begin{tabular}
[c]{llll}%
$\dim\left(  C_{1}\right)  $ &
\begin{tabular}
[c]{l}%
number of faces of $C$\\
of this dimension
\end{tabular}
& $\tilde{S}\left(  C_{1},t\right)  $ & $\tilde{S}\left(  C_{1}^{\ast
},t\right)  $\\
$0$ & $1$ & $1$ & $t+t^{2}+t^{3}$\\
$1$ & $4$ & $0$ & $0$\\
$2$ & $6$ & $3t$ & $0$\\
$3$ & $4$ & $3t+3t^{2}$ & $0$\\
$4$ & $1$ & $t+19t^{2}+t^{3}$ & $1$%
\end{tabular}
\]
hence%
\begin{align*}
E_{st}\left(  X,u,v\right)   &  =1+uv+\left(  uv\right)  ^{2}\\
&  +u^{2}+19uv+v^{2}\\
&  =1+\left(  u^{2}+20uv+v^{2}\right)  +\left(  uv\right)  ^{2}%
\end{align*}

\end{example}

\begin{example}
For the quintic Calabi-Yau threefold in $\mathbb{P}^{4}$ given by a general
section in the degree $5$ Veronese polytope%
\[%
\begin{tabular}
[c]{llll}%
$\dim\left(  C_{1}\right)  $ &
\begin{tabular}
[c]{l}%
number of faces of $C$\\
of this dimension
\end{tabular}
& $\tilde{S}\left(  C_{1},t\right)  $ & $\tilde{S}\left(  C_{1}^{\ast
},t\right)  $\\
$0$ & $1$ & $1$ & $t+t^{2}+t^{3}+t^{4}$\\
$1$ & $5$ & $0$ & $0$\\
$2$ & $10$ & $4t$ & $0$\\
$3$ & $10$ & $6t+6t^{2}$ & $0$\\
$4$ & $5$ & $4t+44t^{2}+4t^{3}$ & $0$\\
$5$ & $1$ & $t+101t^{2}+101t^{3}+t^{4}$ & $1$%
\end{tabular}
\]
hence%
\begin{align*}
E_{st}\left(  X,u,v\right)   &  =1+uv+\left(  uv\right)  ^{2}+\left(
uv\right)  ^{3}\\
&  -\left(  u^{3}+101u^{2}v+101uv^{2}+v^{3}\right) \\
&  =1+uv-\left(  u^{3}+101u^{2}v+101uv^{2}+v^{3}\right)  +\left(  uv\right)
^{2}+\left(  uv\right)  ^{3}%
\end{align*}

\end{example}

\subsubsection{String-theoretic Hodge formula for complete
intersections\label{Sec string theoretic Hodge formula for complete intersections}%
}

Let
\index{reflexive}%
$\Delta\subset M_{\mathbb{R}}$ be a reflexive polytope and $Y=\mathbb{P}%
\left(  \Delta\right)  $ the corresponding
\index{mirror construction}%
Gorenstein toric
\index{Fano}%
Fano
\index{toric Fano}%
variety, denote by $\Sigma\subset N_{\mathbb{R}}$ the
\index{normal fan}%
normal fan of $\Delta$, and let
\index{ray}%
$\Sigma\left(  1\right)  =I_{1}\cup...\cup I_{c}$ be a nef partition, i.e.,
$E_{j}=\sum_{v\in I_{j}}D_{v}$ are Cartier divisors,
\index{spanned by global sections}%
spanned by global sections and $\sum_{j=1}^{c}E_{j}=-K_{Y}$. Denote by
$\Delta_{j}=\Delta_{E_{j}}$ the polytope of sections of $E_{j}$ and by $X$ a
Calabi-Yau complete intersection given by general sections $s_{j}\in
H^{0}\left(  Y,\mathcal{O}_{Y}\left(  E_{j}\right)  \right)  $ for $j=1,...,c$.

Define%
\[
Z=\mathbb{P}\left(  \mathcal{O}_{Y}\left(  E_{1}\right)  \oplus...\oplus
\mathcal{O}_{Y}\left(  E_{c}\right)  \right)
\]
with canonical projection%
\[
\pi:Z\rightarrow Y
\]
Then $\pi_{\ast}\mathcal{O}_{Z}\left(  1\right)  =\mathcal{O}_{Y}\left(
D_{1}\right)  \oplus...\oplus\mathcal{O}_{Y}\left(  D_{c}\right)  $ and%
\[
H^{0}\left(  Z,\mathcal{O}_{Z}\left(  1\right)  \right)  \cong H^{0}\left(
Y,\mathcal{O}_{Y}\left(  E_{1}\right)  \right)  \oplus...\oplus H^{0}\left(
Y,\mathcal{O}_{Y}\left(  E_{c}\right)  \right)
\]
so $\left(  s_{1},...,s_{c}\right)  $ corresponds to a section $s\in
H^{0}\left(  Z,\mathcal{O}_{Z}\left(  1\right)  \right)  $. Let $\bar{X}$ be
the zero set of $s$.

As $X$ is transversal to the toric strata of $Y$%
\[
E_{st}\left(  X,u,v\right)  =E_{st}\left(  Y,u,v\right)  -E_{st}\left(
Y\backslash X,u,v\right)
\]

\begin{proposition}
\cite{BB Mirror duality and stringtheoretic Hodge numbers} $\pi\mid
_{Z\backslash\bar{X}}:Z\backslash\bar{X}\rightarrow Y\backslash X$ is in the
Zariski topology a locally trivial $\mathbb{C}^{c-1}$-bundle, hence%
\[
E_{st}\left(  Y\backslash X,u,v\right)  =\left(  uv\right)  ^{1-c}%
E_{st}\left(  Z\backslash\bar{X},u,v\right)
\]

\end{proposition}

As $Z$ is a $\mathbb{P}_{\mathbb{C}}^{c-1}$-bundle over $Y$%
\[
E_{st}\left(  Y,u,v\right)  =\left(  \left(  uv\right)  ^{c}-1\right)
^{-1}\left(  uv-1\right)  E_{st}\left(  Z,u,v\right)
\]

\begin{proposition}
\cite{BB Mirror duality and stringtheoretic Hodge numbers} The sheaf
$\mathcal{O}_{Z}\left(  1\right)  $ is Cartier, spanned by global sections,
the morphism%
\[
\alpha:Z\rightarrow W=\operatorname*{Proj}\bigoplus\nolimits_{k\geq0}%
H^{0}\left(  Z,\mathcal{O}_{Z}\left(  k\right)  \right)
\]
is crepant, $\mathcal{O}_{Z}\left(  c\right)  $ is the anticanonical sheaf on
$Z$ and $W$ is a Gorenstein toric Fano variety.

$\alpha\left(  \bar{X}\right)  $ is an ample hypersurface in $W$.
\end{proposition}

Note that
\[
W=\operatorname*{Proj}\mathbb{C}\left[  C\right]
\]
with the cone%
\[
C=\left\{  \left(  \lambda_{1},...,\lambda_{c},\sum_{i=1}^{c}\lambda_{i}%
m_{i}\right)  \in\left(  \mathbb{Z}^{r}\oplus M\right)  _{\mathbb{R}}%
\mid\lambda_{i}\in\mathbb{R}_{\geq0},\text{ }m_{i}\in\Delta_{i},\text{
}i=1,...c\right\}
\]
which is a Gorenstein cone with respect to $w\in N$ uniquely defined by%
\begin{align*}
\left\langle m,w\right\rangle  &  =0\text{ for all }m\in M_{\mathbb{R}}%
\subset\left(  \mathbb{Z}^{r}\oplus M\right)  _{\mathbb{R}}\\
\left\langle e_{i},w\right\rangle  &  =0\text{ for all }i=1,...,c
\end{align*}
and has reflexive supporting polyhedron.

Observing that%
\begin{align*}
E_{st}\left(  Z,u,v\right)   &  =E_{st}\left(  W,u,v\right) \\
E_{st}\left(  Z\backslash\bar{X},u,v\right)   &  =E_{st}\left(  W\backslash
\alpha\left(  \bar{X}\right)  ,u,v\right)
\end{align*}
we have%
\begin{align*}
E_{st}\left(  X,u,v\right)   &  =E_{st}\left(  Y,u,v\right)  -E_{st}\left(
Y\backslash X,u,v\right) \\
&  =\left(  \left(  uv\right)  ^{c}-1\right)  ^{-1}\left(  uv-1\right)
E_{st}\left(  Z,u,v\right)  -\left(  uv\right)  ^{1-c}E_{st}\left(
Z\backslash\bar{X},u,v\right) \\
&  =\left(  \left(  uv\right)  ^{c}-1\right)  ^{-1}\left(  uv-1\right)
E_{st}\left(  W,u,v\right)  -\left(  uv\right)  ^{1-c}E_{st}\left(
W\backslash\alpha\left(  \bar{X}\right)  ,u,v\right) \\
&  =\left(  \left(  uv\right)  ^{c}-1\right)  ^{-1}\left(  uv-1\right)
E_{st}\left(  W,u,v\right) \\
&  -\left(  uv\right)  ^{1-c}\left(  E_{st}\left(  W,u,v\right)
-E_{st}\left(  \alpha\left(  \bar{X}\right)  ,u,v\right)  \right)
\end{align*}
The stringy $E$-function $E_{st}\left(  W,u,v\right)  $ can be computed by the
following Proposition \ref{prop stringy E-function and stratifications}, which
shows equality of the stringy $E$-function and the original string-theoretic
$E$-function defined by Batyrev and Dais in \cite{BD Strong McKay
correspondence and stringtheoretic Hodge numbers}.

\begin{proposition}
\cite{BoMa String cohomology of CalabiYau hypersurfaces via mirror symmetry}
\label{prop stringy E-function and stratifications}Let $X=\bigcup
\nolimits_{i\in I}X_{i}$ be a stratified algebraic variety of dimension $n$
with the following properties (satisfied by $W$):

\begin{itemize}
\item $X$ has at most Gorenstein toroidal singularities such that for each
$i\in I$ the singularities of $X$ along the stratum $X_{i}$ of codimension
$c_{i}$ are given by some $c_{i}$-dimensional finite rational polyhedral cone
$\sigma_{i}$. This is equivalent to $X$ being locally isomorphic to
$\mathbb{C}^{n-c_{i}}\times U\left(  \sigma_{i}\right)  $ at all points $x\in
X_{i}$.

\item There is a desingularization $\pi:\overline{X}\rightarrow X$ such that
its restriction to the preimage if $X_{i}$ is a locally trivial fibration in
the Zariski topology.

\item For all points $x\in X_{i}$ the preimage of an analytic neighborhood of
$x$ under $\pi$ is analytically isomorphic the product of a complex disc and a
preimage of a neighborhood of $\left\{  0\right\}  $ in $U\left(  \sigma
_{i}\right)  $ under a resolution of singularities of $U\left(  \sigma
_{i}\right)  $ such that the isomorphism is compatible with the resolutions.
\end{itemize}

Then%
\[
E_{st}\left(  X,u,v\right)  =\sum_{i\in I}E\left(  X_{i},u,v\right)  \cdot
S\left(  \sigma_{i},uv\right)
\]

\end{proposition}

Hence if we denote by $P$ the Eulerian poset of the faces of the cone $C$ with
rank function%
\[%
\begin{tabular}
[c]{llll}%
$\rho:$ & $P$ & $\rightarrow$ & $\left\{  0,...,\dim\left(  C\right)
\right\}  $\\
& $F$ & $\mapsto$ & $\dim\left(  F\right)  $%
\end{tabular}
\]
then%
\[
E_{st}\left(  W,u,v\right)  =\sum_{\substack{x\in P\\x>\min\left(  P\right)
}}\left(  uv-1\right)  ^{\rho\left(  x\right)  -1}S\left(  x^{\ast},uv\right)
\]

In order to compute $E_{st}\left(  \alpha\left(  \bar{X}\right)  ,u,v\right)
$ we can apply Section
\ref{Sec String theoretic Hodge formula for hypersurfaces} to the Gorenstein
cone $C$.

\begin{theorem}
\cite{BB Mirror duality and stringtheoretic Hodge numbers} Let $X\subset
Y=\mathbb{P}\left(  \Delta\right)  $ and $X^{\circ}\subset Y^{\circ
}=\mathbb{P}\left(  \nabla\right)  $ be general complete intersections of
dimension $d$ defined by nef partitions, which are dual to each other with
respect to the construction by Batyrev and Borisov as given in Section
\ref{Sec Batyrev and Borisov mirror construction}. Then the stringy
$E$-functions of $X$ and $X^{\circ}$ satisfy the mirror duality relation%
\[
E_{st}\left(  X;u,v\right)  =\left(  -u\right)  ^{d}E_{st}\left(  X^{\circ
};u^{-1},v\right)
\]

\end{theorem}

\subsubsection{Remarks on a tropical computation of the stringy $E$%
-function\label{Sec first approximation of a tropical computation of the stringy E function}%
}

Consider the setup from Section \ref{Sec tropical mirror construction}. So
denote by $Y=X\left(  \Sigma\right)  $ the toric Fano variety given by the fan
$\Sigma$ over the Fano polytope $\Delta^{\ast}$, and denote by $S$ the Cox
ring of $Y$. Let $\mathfrak{X}\subset Y\times\operatorname*{Spec}%
\mathbb{C}\left[  \left[  t\right]  \right]  $ be the Calabi-Yau degeneration
given by the ideal $I\subset\mathbb{C}\left[  t\right]  \otimes S$ with
monomial special fiber given by $I_{0}\subset S$.

Recall that
\[
\nabla=C_{I_{0}}\left(  I\right)  \cap\left\{  w_{t}=1\right\}
\]
and $\nabla^{\ast}$ is a Fano polytope.

A first approximation of a tropical expression of the stringy $E$-function of
the general fiber $X^{\circ}$ of $\mathfrak{X}^{\circ}$ would be%
\[
E_{st}\left(  X^{\circ},u,v\right)  =\sum_{\substack{x\in BF_{I_{0}}\left(
I\right)  \\\dim\left(  x\right)  >0}}\left(  uv-1\right)  ^{\dim\left(
x\right)  -1}S\left(  x^{\vee},uv\right)
\]
where
\[
x^{\vee}\subset C_{I_{0}}\left(  I\right)  ^{\vee}%
\]
denotes the face dual to $x$ of the Gorenstein cone $C_{I_{0}}\left(
I\right)  ^{\vee}$ over the Fano polytope $\nabla^{\ast}$, i.e., $x^{\ast}$ is
the cone over
\[
\operatorname*{dual}\left(  x\cap\left\{  w_{t}=1\right\}  \right)
\subset\operatorname*{dual}\left(  B\left(  I\right)  \right)  \subset
\operatorname*{Poset}\left(  \nabla^{\ast}\right)
\]
Of course this will not work due to the nature of the singularities of the
reducible $X_{0}^{\circ}$.

One may ask for a formula for $E_{st}\left(  X,u,v\right)  $ in terms of the
data given by the Gorenstein cones $C_{I_{0}}\left(  I\right)  ^{\vee}$ and
$C_{I_{0}^{\circ}}\left(  I^{\circ}\right)  ^{\vee}$ and the subfans
$BF_{I_{0}}\left(  I\right)  ^{\vee}\subset\operatorname*{Poset}\left(
C_{I_{0}}\left(  I\right)  ^{\vee}\right)  $ and $BF_{I_{0}^{\circ}}\left(
I^{\circ}\right)  ^{\vee}\subset\operatorname*{Poset}\left(  C_{I_{0}^{\circ}%
}\left(  I^{\circ}\right)  ^{\vee}\right)  $. This formula should be mirror
symmetric with respect to the tropical mirror construction, i.e., should
satisfy%
\[
E_{st}\left(  X;u,v\right)  =\left(  -u\right)  ^{d}E_{st}\left(  X^{\circ
};u^{-1},v\right)
\]
when exchanging $BF_{I_{0}}\left(  I\right)  $ and $C_{I_{0}}\left(  I\right)
$ with $BF_{I_{0}^{\circ}}\left(  I^{\circ}\right)  $ and $C_{I_{0}^{\circ}%
}\left(  I^{\circ}\right)  $ and should specialize to the formula for
hypersurfaces from Section
\ref{Sec String theoretic Hodge formula for hypersurfaces}.

\section{Implementation of the tropical mirror
construction\label{Sec implementation tropical mirror construction}}

In order to implement the tropical mirror construction, the following packages
for the computer algebra systems Macaulay2 \cite{M2 A Computer Software System
Designed to Support Research in Commutative Algebra and Algebraic Geometry}
and Maple \cite{Maple Maple} have been written by the author:

\subsection{mora.m2\label{Sec moram2}}

The Macaulay2 library
\index{mora.m2}%
\textsf{mora.m2} provides an implementation of the standard basis algorithm.

Polynomials are represented as elements in the Macaulay2 type PolynomialRing
and ideals are represented via the type Ideal.

\begin{itemize}
\item Monomial orderings:

Denoting by $M$ the semigroup of monomials in a polynomial ring, they are
implemented as\textsf{ }functions $f:M\times M\rightarrow\left\{
\mathsf{true,}\text{ }\mathsf{false}\right\}  $ comparing two monomials, where
$f\left(  m_{1},m_{2}\right)  =\mathsf{true}$ if and only if $m_{1}>m_{2}$.

The following monomial orderings as defined in Section
\ref{Sec monomial orderings} are provided by \textsf{mora.m2}. They are
selected by the global variable $\mathsf{monord}$.

\begin{itemize}
\item lexicographical
\index{lexicographic ordering}%
$\mathsf{lp}$

\item reverse lexicographical
\index{reverse lexicographic ordering}%
$\mathsf{rp}$

\item degree reverse lexicographical
\index{degree reverse lexicographic ordering}%
$\mathsf{dp}$

\item negative lexicographical
\index{negative lexicographic ordering}%
$\mathsf{ls}$
\end{itemize}

The following weight orderings depend on a weight vector specified by the
global variable $\mathsf{ww}$ of Macaulay2 type List with rational entries,
whose length is the number of variables of the polynomial ring.

\begin{itemize}
\item weighted reverse lexicographical
\index{weighted reverse lexicographic ordering}%
$\mathsf{wp}$

\item weighted lexicographical
\index{weighted lexicographic ordering}%
$\mathsf{Wp}$

\item local weighted reverse lexicographical $\mathsf{ws}$

\item local weighted lexicographical $\mathsf{Ws}$
\end{itemize}

The matrix ordering
\index{matrix ordering}%
$\mathsf{Mat}$ depends on a matrix $\mathsf{mm}$ of Macaulay2 type Matrix with
rational entries. The number of columns of $\mathsf{Mat}$ is the number of
variables of the polynomial ring.

\item $\mathsf{L}\left(  f\right)  $

Computes the lead monomial of the polynomial $f$ with respect to the semigroup
ordering specified by $\mathsf{monord}$.

\item $\mathsf{SPolynomial}\left(  f,g\right)  $

Returns the s-polynomial%
\[
\operatorname*{SPolynomial}\left(  f,g\right)  =\frac{\operatorname{lcm}%
\left(  L\left(  f\right)  ,L\left(  g\right)  \right)  }{L\left(  f\right)
}f-\frac{LC\left(  f\right)  }{LC\left(  g\right)  }\frac{\operatorname{lcm}%
\left(  L\left(  f\right)  ,L\left(  g\right)  \right)  }{L\left(  g\right)
}g
\]
of the polynomials $f$ and $g$ in the given polynomial ring with semigroup
ordering $\mathsf{monord}$.

\item $\mathsf{NFG}\left(  f,G\right)  $

Computes the Gr\"{o}bner normal form of the polynomial $f$ with respect to the
finite Macaulay2 type list $G$ of polynomials and semigroup ordering
$\mathsf{monord}$ via Algorithm \ref{Alg Groebner normal form}.

\item $\mathsf{redNFG}\left(  f,G\right)  $

Returns the Gr\"{o}bner reduced normal form of the polynomial $f$ with respect
to the list $G$ and semigroup ordering $\mathsf{monord}$, using Algorithm
\ref{Alg Groebner reduced normal form}.

\item $\mathsf{NF}\left(  f,G\right)  $

Computes the Mora normal form of the polynomial $f$ with respect to the list
$G$ and semigroup ordering $\mathsf{monord}$ by Algorithm
\ref{Alg Mora normal form}.

\item $\mathsf{Std}\left(  G\right)  $

Implements Algorithm \ref{Alg standard basis} to compute a standard basis of
the ideal $\left\langle G\right\rangle $ for a list $G$ of elements in a
polynomial ring and semigroup ordering $\mathsf{monord}$.

\item $\mathsf{Minimize}\left(  G\right)  $

Given a list $G$ of polynomials computes an interreduced subset with respect
to the semigroup ordering $\mathsf{monord}$.

\item $\mathsf{MStd}\left(  G\right)  $

Returns a minimal standard basis of the ideal $\left\langle G\right\rangle $
for a list $G$ of elements in a polynomial ring and semigroup ordering
$\mathsf{monord}$.

\item $\mathsf{ReduceGb}\left(  G\right)  $

Given a minimal Gr\"{o}bner basis $G$ of the ideal $\left\langle
G\right\rangle $ with respect to $\mathsf{monord}$, returns a reduced
Gr\"{o}bner basis of $\left\langle G\right\rangle $ via Algorithm
\ref{Alg reduced std}.

\item $\mathsf{ReduceStd}\left(  G\right)  $

Applying Algorithm \ref{Alg reduced std mora} takes a minimal standard basis
$G$ of the ideal $\left\langle G\right\rangle $ with respect to
$\mathsf{monord}$ and computes a reduced standard basis of $\left\langle
G\right\rangle $ by the Gr\"{o}bner normal form. If the reduction does not
terminate, the procedure stops after a finite number of reductions of each
element of $G$ specified by the global variable $\mathsf{iterlimit}$.
\end{itemize}

The global variable $\mathsf{verbose}\in\left\{  0,1,2\right\}  $ controls the
output of intermediate results, e.g. of syzygies in Gr\"{o}bner computations.

\begin{example}
Load the package and create a polynomial ring:\textrm{\smallskip}%
\newline\texttt{load "mora.m2";}\newline\texttt{R=QQ[x,y,z];}\textrm{\medskip
}\newline Lead monomials with respect to various orderings:\textrm{\smallskip
}\newline\texttt{f=x\symbol{94}4+y\symbol{94}7+z\symbol{94}5+x\symbol{94}%
4*y*z+x\symbol{94}3*y\symbol{94}3;\newline monord=lp;\newline L(f)}%
\newline\textrm{x\symbol{94}4*y*z\medskip}\newline\texttt{monord=dp;}%
\newline\texttt{L(f)}\newline\textrm{y\symbol{94}7\medskip}\newline%
\texttt{monord=ls;\newline L(f)\newline}\textrm{z\symbol{94}5\medskip\newline%
}\texttt{monord=ds;\newline L(f)\newline}\textrm{x\symbol{94}4\medskip
\newline}\texttt{monord=Wp;\newline ww=\{2,1,-1\};\newline L(f)\newline%
}\textrm{x\symbol{94}3*y\symbol{94}3\medskip\newline}%
\texttt{monord=Wp;\newline ww=\{-3,-1,-2\}\newline L(f)\newline}%
\textrm{y\symbol{94}7\medskip\newline}\texttt{monord=Mat;\newline MM=matrix
\{\{-3,-1,-2\},\{1,0,0\},\{0,1,0\},\{0,0,1\}\};\newline L(f)\newline%
}\textrm{y\symbol{94}7\medskip}\newline Computing standard bases, division
with remainder:\textrm{\smallskip}\newline\texttt{monord=lp;\newline
G=\{x*y-1,y\symbol{94}2-1\};\newline std(G);\newline}\textrm{\{x*y-1,
y\symbol{94}2-1, x-y\}\medskip\newline}\texttt{GB=minimalstd(G);\newline%
}\textrm{\{y\symbol{94}2-1, x-y\}\medskip\newline}\texttt{f=x\symbol{94}%
2*y+x*y\symbol{94}2+y\symbol{94}2;\newline NFB(f,G)\newline}%
\textrm{x+y\symbol{94}2+y\medskip\newline}\texttt{NFB(f,GB)}\newline%
\textrm{2y+1\medskip}\newline\texttt{redNFB(f,GB)}\newline%
\textrm{y+1/2\medskip}\newline\texttt{G=\{x\symbol{94}2+y,x*y+x\};\newline
GB=minimalstd(G)\newline}\textrm{\{x\symbol{94}2+y, x*y+x, y\symbol{94}%
2+y\}}\texttt{\textrm{\medskip}\newline f=x\symbol{94}2-y\symbol{94}2;\newline
redNFB(f,G)\newline}\textrm{y\symbol{94}2+y}\texttt{\textrm{\medskip}\newline
redNFB(f,GB)}\newline\textrm{0}\texttt{\textrm{\medskip}}

Mora normal form and Gr\"{o}bner normal form for local
orderings:\textrm{\smallskip}\newline\texttt{monord=ls;\newline
NF(x,\{x-x\symbol{94}2\})\newline}\textrm{0}\texttt{\textrm{\medskip}\newline
iterlimit=50;\newline NFB(x,\{x-x\symbol{94}2\})}\newline\textrm{x\symbol{94}%
51}\texttt{\textrm{\medskip}}\newline\texttt{monord=ls;\newline f=z\symbol{94}%
2+y*z+y\symbol{94}2+x\symbol{94}2;\newline G=\{x,z\};\newline NF(f,G)\newline%
}\textrm{y\symbol{94}2+x\symbol{94}2}\texttt{\textrm{\medskip}}

Minimal standard bases of the ideal $\left\langle G\right\rangle $ for the
monomial orderings $dp$, $lp$, $ds$ and $ls$:\textrm{\smallskip}%
\newline\texttt{G=\{x\symbol{94}6+x\symbol{94}5*y\symbol{94}2,y\symbol{94}%
4-x\symbol{94}2*y\symbol{94}3\};\newline monord=dp;\newline
minimalstd(G)\newline}\textrm{\{x\symbol{94}6+x\symbol{94}5*y\symbol{94}2,
-x\symbol{94}2*y\symbol{94}3+y\symbol{94}4, x\symbol{94}6*y+x\symbol{94}%
3*y\symbol{94}4, -x\symbol{94}7+x*y\symbol{94}6, x\symbol{94}2*y\symbol{94}%
5+x*y\symbol{94}7, y\symbol{94}8+x\symbol{94}7\}}\texttt{\textrm{\medskip
}\newline monord=lp\newline minimalstd(G)\newline}\textrm{\{x\symbol{94}%
6+x\symbol{94}5*y\symbol{94}2, -x\symbol{94}2*y\symbol{94}3+y\symbol{94}4,
x*y\symbol{94}6+y\symbol{94}8, -y\symbol{94}9+y\symbol{94}6\}}%
\texttt{\textrm{\medskip}\newline monord=ds\newline minimalstd(G)\newline%
}\textrm{\{x\symbol{94}6+x\symbol{94}5*y\symbol{94}2,y\symbol{94}%
4-x\symbol{94}2*y\symbol{94}3\}}\texttt{\textrm{\medskip}\newline
monord=ls\newline GB=minimalstd(G)\newline}\textrm{\{x\symbol{94}%
5*y\symbol{94}2+x\symbol{94}6, y\symbol{94}4-x\symbol{94}2*y\symbol{94}3,
x\symbol{94}7*y\symbol{94}3-x\symbol{94}7\}}\texttt{\textrm{\medskip}\newline
iterlimit=10\symbol{94}6;\newline reducestd(GB)\newline}\textrm{\{x\symbol{94}%
5*y\symbol{94}2+x\symbol{94}6, y\symbol{94}4-x\symbol{94}2*y\symbol{94}3,
x\symbol{94}7\}}

\noindent Note that $y^{3}-1$ is a unit in $R_{>}$ for the negative
lexicographic ordering $>=ls$ and $x^{7}y^{3}-x^{7}=x^{7}\left(
y^{3}-1\right)  $.\texttt{\textrm{\medskip}}

The ideal of a line and a plane in the global setting and in the local ring
$\mathbb{Q}\left[  x,y,z\right]  _{\left\langle x,y,z\right\rangle }%
$:\textrm{\smallskip}\newline\texttt{G=\{x*y+y,x*z+z\};\newline
monord=dp;\newline minimalstd(G)}\newline\textrm{\{x*y+y,x*z+z\}}%
\texttt{\textrm{\medskip}}\newline\texttt{NF(x,GB)\newline}%
\textrm{x\texttt{\medskip}}\newline\texttt{NF(y,GB)\newline}%
\textrm{y\texttt{\medskip}}\newline\texttt{NF(z,GB)\newline}%
\textrm{z\texttt{\medskip}}\texttt{\newline monord=ls;\newline
GB=minimalstd(G);\newline}\textrm{\{x*y + y, x*z + z\}\texttt{\medskip
}\newline}\texttt{NF(x,GB)\newline}\textrm{x\texttt{\medskip}}\newline%
\texttt{NF(y,GB)\newline}\textrm{0\texttt{\medskip}}\newline%
\texttt{NF(z,GB)\newline}\textrm{0\texttt{\medskip}}\newline%
\texttt{reducestd(GB)\newline}\textrm{\{y, z\}}

\noindent Note that $1+x$ is a unit in $R_{>}$ for the negative lexicographic
ordering $>=ls$.
\end{example}

\subsection{homology.m2\label{Sec homologym2}}

The Macaulay2 library
\index{homology.m2}%
\textsf{homology.m2} provides the following functions:

Let $C$ be a cell complex given as a list, where the $d$-th element is a list
of the faces of dimension $d$, and each face is given as a list of vertices.

\begin{itemize}
\item $\mathsf{IsSimplicial}\left(  C\right)  $

Checks if $C$ is simplicial.

\item $\mathsf{AssociatedChainComplex}\left(  C,R\right)  $

Associates to the simplicial complex $C$ the associated chain complex with
coefficients in the Macaulay2 type ring $R$. The resulting chain complex is
represented as a Macaulay2 type chain complex. The orientation of the cells of
$C$ is represented by the ordering of the vertices in the lists representing
the faces of $C$ and the boundary maps are given by%
\[
\partial\left(  v_{i_{0}},...,v_{i_{d}}\right)  =\sum_{j=0}^{d}\left(
-1\right)  ^{j}\left(  v_{i_{0}},...,\widehat{v_{i_{j}}},...,v_{i_{d}}\right)
\]

\item $\mathsf{BoundaryMap}\left(  C,R,d\right)  $

Returns above boundary map $\partial:C_{d}\rightarrow C_{d-1}$.

\item $\mathsf{HomologyChainComplex}\left(  C,R\right)  $

Computes a list with the homology groups of $C$ with coefficients in $R$.
\end{itemize}

\begin{example}
Consider the following triangulation of the Klein bottle%
\[%
{\includegraphics[
height=2.5702in,
width=2.5227in
]%
{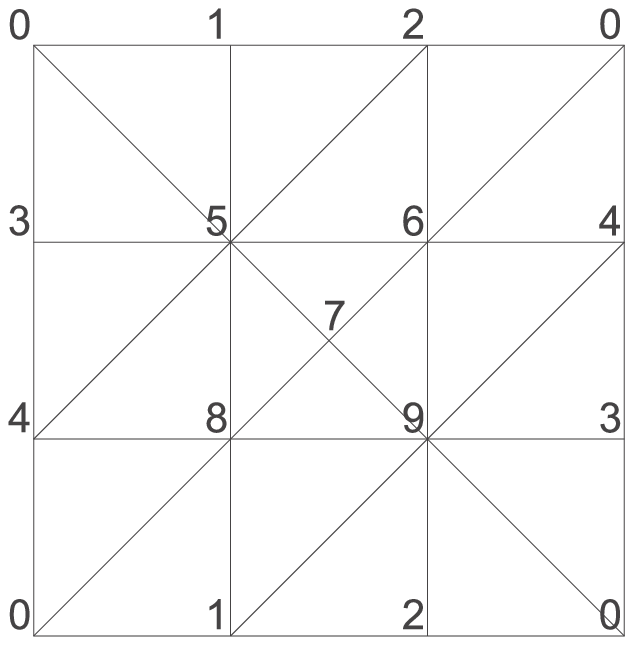}%
}%
\]
\texttt{load "homology.m2";}\newline\newline%
\texttt{C=\{\{\{0\},\{1\},\{2\},\{3\},\{4\},\{5\},\{6\},\{7\},\{8\},\{9\}\},}%
\newline%
\texttt{\{\{0,1\},\{1,2\},\{0,2\},\{0,8\},\{1,8\},\{1,9\},\{2,9\},\{0,9\},\{0,3\},\{4,8\},}%
\newline%
\texttt{\{8,9\},\{3,9\},\{3,4\},\{4,5\},\{5,8\},\{7,8\},\{7,9\},\{6,7\},\{5,7\},\{3,5\},}%
\newline%
\texttt{\{5,6\},\{6,9\},\{4,6\},\{4,9\},\{0,4\},\{0,5\},\{1,5\},\{2,5\},\{2,6\},\{0,6\}\},}%
\newline%
\texttt{\{\{0,1,8\},\{1,2,9\},\{2,0,9\},\{0,3,9\},\{1,9,8\},\{0,8,4\},\{3,4,5\},\{4,8,5\},}%
\newline%
\texttt{\{5,8,7\},\{7,8,9\},\{5,7,6\},\{6,7,9\},\{4,6,9\},\{3,4,9\},\{0,3,5\},\{0,5,1\},}%
\newline\texttt{\{1,5,2\},\{2,5,6\},\{0,2,6\},\{0,6,4\}\}\};}\newline%
\texttt{\newline cC=SimplicialChainComplex(C,ZZ)}\newline$0\leftarrow
\mathbb{Z}^{10}\leftarrow\mathbb{Z}^{30}\leftarrow\mathbb{Z}^{20}\leftarrow
0$\newline\texttt{HomologyChainComplex(cC)}\newline$\left(  \mathbb{Z}%
,\mathbb{Z\oplus Z}/2,0\right)  $\newline\newline%
\texttt{cC=SimplicialChainComplex(C,QQ);}\newline%
\texttt{HomologyChainComplex(cC)}\newline$\left(  \mathbb{Q},\mathbb{Q}%
,0\right)  $\newline\newline\texttt{cC=SimplicialChainComplex(C,ZZ/2);}%
\newline\texttt{HomologyChainComplex(cC)}\newline$\left(  \mathbb{Z}%
/2,\mathbb{Z}/2\mathbb{\oplus Z}/2,\mathbb{Z}/2\right)  $\newline
\end{example}

\subsection{stanleyfiltration.m2\label{Sec stanley filtration}}

The Macaulay2 library
\index{stanleyfiltration.m2}%
\textsf{stanleyfiltration.m2} provides the following functions:

\begin{itemize}
\item $\mathsf{StanleyDecomposition}\left(  I\right)  $

Implements Algorithm \ref{alg stanley decomposition} to compute a Stanley
decomposition%
\begin{align*}
S/I  &  \cong%
{\displaystyle\bigoplus\limits_{\left(  D,\sigma\right)  \in\mathcal{S}}}
S_{\sigma}\left(  -\left[  D\right]  \right) \\
\mathcal{S}  &  \subset\left\{  \left(  D,\sigma\right)  \mid D\in
\operatorname*{WDiv}\nolimits_{T}\left(  Y^{\prime\prime}\right)  \text{,
}D\text{ effective, }\sigma\in\Sigma^{\prime\prime}\right\}
\end{align*}
of a monomial ideal $I$ in the polynomial ring $S$, where $\Sigma
^{\prime\prime}$ is the fan over the simplex on the variables of $S$ and
$Y^{\prime\prime}=\mathbb{A}^{\Sigma^{\prime\prime}\left(  1\right)
}=\operatorname*{Spec}S$.

The ring $S$ is represented via the Macaulay2 type ring and $I$ via the
Macaulay2 type ideal.

The output is a set of tuples $\left(  m,P\right)  $ representing $\left(
D,\sigma\right)  \in\mathcal{S}$, where $m$ is a monomial in $S$ defining the
divisor $D$ and $P$ is a set of variables of $S$ generating the cone
$\sigma\in\Sigma^{\prime\prime}$.

\item $\mathsf{StanleyFiltration}\left(  I\right)  $

Returns a list with a Stanley filtration of the monomial ideal $I\subset S$.
The elements of the list are represented in the same way as for the output of
$\mathsf{StanleyDecomposition}$.

\item $\mathsf{MonomialIdealsFixedHilbertPolynomial}\left(  S,P,A,B\right)  $

Returns the set of monomial ideals in the multigraded polynomial ring
$S=\mathbb{Q}\left[  y_{1},...,y_{r}\right]  $ with multigraded Hilbert
polynomial $P\in\mathbb{Q}\left[  t_{1},...,t_{a}\right]  $, where
$A\in\mathbb{Z}^{d\times r}$ is the presentation matrix of the Chow group of a
smooth toric variety $Y$ and $B$ is the irrelevant ideal of $Y$.
\end{itemize}

\begin{example}
Consider the ideal%
\[
I=\left\langle y_{1}y_{2},y_{0}y_{3}\right\rangle \subset S=\mathbb{C}\left[
y_{0},...,y_{3}\right]
\]
\texttt{load "stanleyfiltration.m2";}\newline\newline%
\texttt{S=QQ[y\_0..y\_3];}\newline\newline%
\texttt{I=ideal(y\_1*y\_2,y\_0*y\_3);}\newline\newline%
\texttt{StanleyFiltration(I)}\newline\texttt{\{\{1, \{y\_0,y\_1\}\}, \{y1,
\{y\_0,y\_2\}\}, \{y0,\{y\_1,y\_3\}\}, \{y0*y1, \{y\_2,y\_3\}\}\}\newline%
}\newline This corresponds to the Stanley decomposition%
\begin{align*}
S/I  &  =1\cdot\mathbb{C}\left[  y_{2},y_{3}\right]  \mathbb{\oplus}y_{1}%
\cdot\mathbb{C}\left[  y_{1},y_{3}\right]  \mathbb{\oplus}y_{0}\cdot
\mathbb{C}\left[  y_{0},y_{2}\right]  \mathbb{\oplus}y_{0}y_{1}\cdot
\mathbb{C}\left[  y_{0},y_{1}\right] \\
I  &  =\left\langle y_{0},y_{1}\right\rangle \cap\left\langle y_{0}%
,y_{2}\right\rangle \cap\left\langle y_{1},y_{3}\right\rangle \cap\left\langle
y_{2},y_{3}\right\rangle
\end{align*}
and to the Stanley filtration given by the Stanley decompositions%
\begin{align*}
S/\left\langle y_{1},y_{0}\right\rangle  &  =1\cdot\mathbb{C}\left[
y_{2},y_{3}\right] \\
S/\left\langle y_{1}y_{2},y_{0}\right\rangle  &  =1\cdot\mathbb{C}\left[
y_{2},y_{3}\right]  \mathbb{\oplus}y_{1}\cdot\mathbb{C}\left[  y_{1}%
,y_{3}\right] \\
S/\left\langle y_{1}y_{2},y_{0}y_{3},y_{0}y_{1}\right\rangle  &
=1\cdot\mathbb{C}\left[  y_{2},y_{3}\right]  \mathbb{\oplus}y_{1}%
\cdot\mathbb{C}\left[  y_{1},y_{3}\right]  \mathbb{\oplus}y_{0}\cdot
\mathbb{C}\left[  y_{0},y_{2}\right] \\
S/\left\langle y_{1}y_{2},y_{0}y_{3}\right\rangle  &  =1\cdot\mathbb{C}\left[
y_{2},y_{3}\right]  \mathbb{\oplus}y_{1}\cdot\mathbb{C}\left[  y_{1}%
,y_{3}\right]  \mathbb{\oplus}y_{0}\cdot\mathbb{C}\left[  y_{0},y_{2}\right]
\mathbb{\oplus}y_{0}y_{1}\cdot\mathbb{C}\left[  y_{0},y_{1}\right]
\end{align*}
\newline
\end{example}

\subsection{tropicalmirror\label{Sec tropicalmirror implementation}}

In the Maple
\index{tropical mirror|textbf}%
package \textsf{tropicalmirror} we provide an implementation of the tropical
mirror construction given in the Sections
\ref{Sec tropical mirror construction} and
\ref{Sec tropical mirror construction for complete intersections}. It also
contains an implementation of the algorithms from Section
\ref{Sec Computing the Bergman fan} computing the Gr\"{o}bner and Bergman fan.

In addition to standard Maple packages, \textsf{tropicalmirror} assumes the
\textsf{convex }package for convex geometry to be present. For local
Gr\"{o}bner computations \textsf{tropicalmirror} allows to call:

\begin{itemize}
\item Macaulay2 with \textsf{mora.m2}.

\item Macaulay2 with Lazard ordering.

\item Singular with built in monomial orderings.

\item Singular with Lazard ordering.
\end{itemize}

The weight orderings can be represented in Macaulay2 and Singular as $Wp,$
$wp$ or by a matrix ordering. \textsf{tropicalmirror} assumes that the
following variables of type string are present:

\begin{itemize}
\item \texttt{runM2} with the command running Macaulay2 in the shell.

\item \texttt{runSingular} with the command running Singular in the shell.

\item \texttt{stdSystem} with value \texttt{M2} or \texttt{Singular} selecting
the computer algebra system for Gr\"{o}bner calculations.

\item \texttt{pathConvex} with the path to the convex package.\smallskip
\end{itemize}

Let $N=\mathbb{Z}^{n}$, $P\subset N_{\mathbb{R}}\mathbb{=}N\otimes\mathbb{R}$
be a Fano polytope, $\Sigma$ the fan over the faces of $P$ and $Y=X\left(
\Sigma\right)  $ the corresponding toric Fano variety of dimension $n$ as
defined in Section \ref{Sec Fano polytopes}. The polytope $P$ is represented
as type polytope and the fan $\Sigma$ as type fan in the convex package.
Choosing a numbering of rays of $\Sigma$, let $A$ be a Maple type matrix
presenting the Chow group of $Y$ via%
\[
0\rightarrow\mathbb{Z}^{n}\overset{A}{\rightarrow}\mathbb{Z}^{r}\rightarrow
A_{n-1}\left(  Y\right)  \rightarrow0
\]
as given in Section \ref{Divisors on toric varieties}. Let $v$ be a list of
names for the variables corresponding to the rays of $\Sigma$ in the rows of
$A$. Denote the $T$-Weil divisors of $Y$ corresponding to the $j$-th row of
$A$ by $D_{j}$.

The package \textsf{tropicalmirror} provides the following functions. They are
organized in a way to avoid multiple computations of the same result.

\begin{itemize}
\item $\mathsf{FanOverFaces}\left(  P\right)  $

Returns the fan over the faces of the polytope $P$ containing $0$ as defined
in Section \ref{Sec Fano polytopes}.

\item $\mathsf{RandomPolynomial}(A,v,a,c)$

Let $a$ be an element of $\mathbb{Z}^{r}$ corresponding the Weil divisor
$D=\sum_{r\in\Sigma\left(  1\right)  }a_{j}D_{j}$ representing the class
$\left[  D\right]  \in A_{n-1}\left(  Y\right)  $. The function
$\mathsf{RandomPolynomial}$ returns a Maple type polynomial $f\in S_{\left[
D\right]  }$ in the variables given by the Maple type list $v$ with
coefficients in $\left\{  1,...,c-1\right\}  $ such that all monomials in
$S_{\left[  D\right]  }$ appear in the polynomial. The Cox polynomial $f$ is
obtained as explained in Section \ref{Sec Global sections a cox monomials} and
corresponds to a generic linear combination of the lattice points of%
\[
\Delta_{D}=\left\{  m\in M_{\mathbb{R}}\mid\left\langle m,\hat{r}\right\rangle
\geq-a_{r}\forall r\in\Sigma\left(  1\right)  \right\}
\]
with $M=\operatorname*{Hom}\left(  N,\mathbb{Z}\right)  $, which form a
$T$-invariant basis of the space of
\index{global sections}%
global sections%
\[
H^{0}\left(  Y,\mathcal{O}_{Y}\left(  D\right)  \right)  \cong\bigoplus
_{m\in\Delta_{D}\cap M}\mathbb{C}x^{m}%
\]
of the
\index{reflexive sheaf}%
reflexive sheaf $\mathcal{O}_{Y}\left(  D\right)  $ as explained in Section
\ref{Divisors on toric varieties}.

\item $\mathsf{ReduceGenerators}\left(  v,t,gI\right)  $

Let $gI=\left[  f_{1},...,f_{r}\right]  $ be a list of Maple type polynomials
representing Cox homogeneous elements in $\mathbb{C}\left[  t\right]
\otimes_{\mathbb{C}}S$ such that for each polynomial $f_{j}$ the degree $0$
part with respect to the $t$-degree is a monomial $m_{j}$ in $S$. The function
$\mathsf{ReduceGenerators}$ removes all terms of $f_{j}-m_{j}$ which are
divisible by some $m_{i}$. Up to first order this amounts to Gr\"{o}bner
reduction of $gI$.

\item $\mathsf{AssociatedFirstOrderDegeneration}\left(  v,t,gI\right)  $

Deletes all terms of $t$-degree bigger than $1$ from the polynomials $f_{j}%
\in\mathbb{C}\left[  t\right]  \otimes_{\mathbb{C}}S$ in the list $gI=\left[
f_{1},...,f_{r}\right]  $.

\item $\mathsf{SpecialFiberGroebnerCone}\left(  A,v,t,gI\right)  $

Suppose $\mathfrak{X}\subset\mathbb{A}^{1}\times Y$ is a flat family of
Calabi-Yau varieties of dimension $d$ with monomial special fiber, given by
the ideal $I\subset\mathbb{C}\left[  t\right]  \otimes_{\mathbb{C}}S$
generated by the Cox homogeneous elements $f_{j}$ of the list $gI=\left[
f_{1},...,f_{r}\right]  $. Assume that the monomials of $f_{j}$ of $t$-degree
$0$ are minimal generators of the $B\left(  \Sigma\right)  $-saturated
monomial ideal $I_{0}$ of the special fiber of $\mathfrak{X}$. The function
$\mathsf{SpecialFiberGroebnerCone}$ returns the special fiber Gr\"{o}bner cone
$C_{I_{0}}\left(  I\right)  $ as defined in Section
\ref{Sec the groebner cone associated to the special fiber}. It is represented
by the convex type cone.

\item $\mathsf{GroebnerFan}\left(  \Sigma,A,v,t,s,gI\right)  $

Computes the ideal $J\subset\mathbb{C}\left[  t,s\right]  \otimes_{\mathbb{C}%
}S$ of the projective closure $\overline{\mathfrak{X}}\subset\mathbb{P}%
^{1}\times Y$ of the flat family $\mathfrak{X}\subset\mathbb{A}^{1}\times Y$
given by the ideal $I$ generated by the elements $f_{j}\in\mathbb{C}\left[
t\right]  \otimes_{\mathbb{C}}S$ of the list $gI=\left[  f_{1},...,f_{r}%
\right]  $. Returns the Gr\"{o}bner fan of $J$ as a subfan of $\mathbb{R\oplus
}N_{\mathbb{R}}$, computed as explained in Section
\ref{Sec computing the Groebner fan}. It is represented by the convex package
type fan.

\item $\mathsf{BergmanFan}\left(  \Sigma,A,v,t,s,gI,GF\right)  $

Computes the Bergman subfan of the Gr\"{o}bner fan $GF$ as explained in
Section \ref{Sec Computing the Bergman fan}, where $GF$ is the result returned
by $\mathsf{GroebnerFan}\left(  v,t,s,A,gI\right)  $. The result is
represented by the convex package type fan.

\item $\mathsf{AssociatedAnticanonicalSectionsPolytope}\left(  C\right)  $

Intersects the special fiber Gr\"{o}bner cone $C\subset\mathbb{R\oplus
}N_{\mathbb{R}}$ with the hyperplane $\left\{  w_{t}=1\right\}  $ and returns
the resulting polytope in $N_{\mathbb{R}}$.

\item $\mathsf{AssociatedFanoPolytope}\left(  C\right)  $

Computes $\nabla=C\cap\left\{  w_{t}=1\right\}  \subset N_{\mathbb{R}}$ and
returns the polytope $\nabla^{\ast}\subset M_{\mathbb{R}}$.

\item $\mathsf{FacePoset}\left(  \nabla\right)  $

Returns the complex of faces of a polytope $\nabla\subset N_{\mathbb{R}}$,
represented as a list of lists $L=\left[  L_{-1},L_{0},...,L_{n}\right]  $.
The list $L_{j}$ contains the faces of $\nabla$ of dimension $j$ and each face
is represented by the convex package type face of polyhedron.

\item $\mathsf{VertexRepresentation}\left(  B\right)  $

Given a complex $B$ represented as a list of lists, where each face is of typ
face of polyhedron, returns a list of lists, where each face is represented as
a list of its vertices.

\item $\mathsf{ChowGroup}\left(  A\right)  $

Returns a group isomorphic to the Chow group $A_{n-1}\left(  Y\right)  $ of
$Y=X\left(  \Sigma\right)  $, given as the cokernel of a diagonal matrix
$A^{\prime}$ of the same dimensions as $A$.

\item $\mathsf{ChowGroupAction}\left(  A\right)  $

Computes isomorphisms $W\in\operatorname*{GL}\left(  n,\mathbb{Z}\right)  $
and $U=\left(  u_{ij}\right)  \in\operatorname*{GL}\left(  r,\mathbb{Z}%
\right)  $%
\[%
\begin{tabular}
[c]{lllllll}%
$0\rightarrow$ & $\mathbb{Z}^{n}$ & $\overset{A}{\rightarrow}$ &
$\mathbb{Z}^{r}$ & $\rightarrow$ & $A_{n-1}\left(  Y\right)  $ &
$\rightarrow0$\\
& $\downarrow W$ &  & $\downarrow U$ &  & $\downarrow$ & \\
$0\rightarrow$ & $\mathbb{Z}^{n}$ & $\overset{A^{\prime}}{\rightarrow}$ &
$\mathbb{Z}^{r}$ & $\rightarrow$ & $H$ & $\rightarrow0$%
\end{tabular}
\]
such that $A^{\prime}$ is a matrix with non zero entries only on the diagonal.
As explained in Section \ref{Sec Cox quotient presentation}, the group%
\[
G\left(  \Sigma\right)  =\operatorname*{Hom}\nolimits_{\mathbb{Z}}\left(
A_{n-1}\left(  Y\right)  ,\mathbb{C}^{\ast}\right)
\]
acts on the affine Cox space $\operatorname*{Hom}\nolimits_{sg}\left(
\operatorname*{WDiv}\nolimits_{T}\left(  Y\right)  ,\mathbb{C}\right)
\cong\mathbb{C}^{r}$ by%
\begin{align*}
G\left(  \Sigma\right)  \times\operatorname*{Hom}\nolimits_{sg}\left(
\operatorname*{WDiv}\nolimits_{T}\left(  Y\right)  ,\mathbb{C}\right)   &
\rightarrow\operatorname*{Hom}\nolimits_{sg}\left(  \operatorname*{WDiv}%
\nolimits_{T}\left(  Y\right)  ,\mathbb{C}\right) \\
\left(  g,a\right)   &  \mapsto%
\begin{tabular}
[c]{llll}%
$ga:$ & $\operatorname*{WDiv}\nolimits_{T}\left(  Y\right)  $ & $\rightarrow$
& $\mathbb{C}$\\
& $D_{r}$ & $\mapsto$ & $g\left(  \left[  D_{r}\right]  \right)  a\left(
D_{r}\right)  $%
\end{tabular}
\end{align*}
hence
\[
G\left(  \Sigma\right)  ^{\prime}=\operatorname*{Hom}\nolimits_{\mathbb{Z}%
}\left(  H,\mathbb{C}^{\ast}\right)
\]
acts by%
\begin{align*}
G\left(  \Sigma\right)  ^{\prime}\times\mathbb{C}^{r}  &  \rightarrow
\mathbb{C}^{r}\\
\left(  \left(  t_{j}\right)  ,\left(  a_{j}\right)  \right)   &
\mapsto\left(
{\textstyle\prod\nolimits_{i=1}^{r}}
t_{i}^{u_{ij}}a_{j}\right)  _{j=1,...,r}%
\end{align*}
Defining
\[
T=\left(
{\textstyle\prod\nolimits_{i=1}^{r}}
t_{i}^{u_{ij}}\right)  _{j=1,...,r}%
\]

the function $\mathsf{ChowGroupAction}$ returns the list $\left[  A^{\prime
},T\right]  $.

\item $\mathsf{IrrelevantIdeal}\left(  \Sigma,A,v\right)  $

Returns the irrelevant ideal $B\left(  \Sigma\right)  \subset S$ as defined in
Section \ref{Sec Cox quotient presentation}. The variables in the list $v$
corresponds to the rows of $A$.

\item $\mathsf{BergmanSubfanOfGroebnerCone}\left(  \Sigma,A,v,t,gI,C\right)  $

Computes the Bergman subfan of the fan of faces of $C$ as defined in Section
\ref{Sec Bergman subcomplex of Nabla general setting} for the ideal $I$
generated by the elements of the list $gI$.

\item $\mathsf{BergmanSubcomplexOfSectionsPolytope}\left(  \Sigma
,A,v,t,gI,posetNabla\right)  $

Returns the Bergman subcomplex $B\left(  I\right)  $ as defined in Section
\ref{Sec Bergman subcomplex of Nabla general setting}, i.e., the intersection
of the output of $\mathsf{BergmanSubfanOfGroebnerCone}$ with the hyperplane
$\left\{  w_{t}=1\right\}  $. The result is a subcomplex of the complex of
faces $posetNabla=\mathsf{FacePoset}\left(  \nabla\right)  $ of $\nabla
=C\cap\left\{  w_{t}=1\right\}  \subset N_{\mathbb{R}}$ and is represented as
a list of faces of $\nabla$ of the form $\left[  ...,B_{j},...\right]  $ where
$B_{j}$ is a list of faces of dimension $j$. For practical reasons it is
useful to fix a numbering of the faces in each dimension, so we represent
$B_{j}$ as a list. Each face is represented by the convex package type face of
a polyhedron.

\item $\mathsf{SpecialFiberIdeal}\left(  AMirror,z,posetNabla,B\right)  $

Suppose $z$ is a list of names for the variables of the Cox ring of $Y^{\circ
}=X\left(  \operatorname*{NF}\left(  \nabla\right)  \right)  $ corresponding
the rays of the normal fan of $\nabla$, which are numbered by the rows in the
matrix $AMirror$. Suppose $posetNabla=\mathsf{FacePoset}\left(  \nabla\right)
$ and $B$ is the Bergman subcomplex of $posetNabla$ as given by the function
$\mathsf{BergmanSubcomplexOfSectionsPolytope}$. Then
$\mathsf{SpecialFiberIdeal}$ returns the ideal%
\[
I_{0}^{\circ}=%
{\textstyle\bigcap\nolimits_{F\in B_{d}}}
\left\langle z_{G^{\ast}}\mid G\text{ a facet of }\nabla\text{ with }F\subset
G\right\rangle
\]
which gives the subvariety $X_{0}^{\circ}\subset Y^{\circ}$ as defined in
Section \ref{Sec special fiber mirror degeneration general setup}. The ideal
$I_{0}^{\circ}$ is represented as a list of generators.

\item $\mathsf{ToricStrataDecomposition}\left(  posetDelta,I_{0}\right)  $

Returns the subcomplex $\operatorname*{Strata}_{\Delta}\left(  I_{0}\right)
\subset\operatorname*{Poset}\left(  \Delta\right)  $ as defined in Section
\ref{Monomial ideals in the Cox ring and the stratified toric primary decomposition}%
. It is represented as a list of lists and each face is of the convex package
type face of polyhedron.

\item $\mathsf{ComplexOfInitialIdeals}\left(  v,t,gI,C\right)  $

Gives the complex of initial ideals $\operatorname*{in}_{F}\left(  I\right)  $
for the faces $F$ of $B$. It is represented as a list of lists in the same way
as the Bergman subcomplex $B$. Each ideal $\operatorname*{in}_{F}\left(
I\right)  $ is represented by a list containing a standard basis with respect
to a monomial ordering in the interior of $C$.

\item $\mathsf{DualComplex}\left(  A,v,t,inI,B\right)  $

Computes the dual complex $\operatorname*{dual}\left(  B\left(  I\right)
\right)  \subset\operatorname*{Poset}\left(  \nabla^{\ast}\right)  \subset
M_{\mathbb{R}}$ as given by the map $\operatorname*{dual}$ defined in Section
\ref{Sec dual complex general setting}. Here $B$ denotes the Bergman
subcomplex $B\left(  I\right)  $ as returned by
$\mathsf{BergmanSubcomplexOfSectionsPolytope}$ and $inI$ is the complex of
initial ideals as returned by \linebreak$\mathsf{ComplexOfInitialIdeals}$. The
complex $\operatorname*{dual}\left(  B\left(  I\right)  \right)  $ is
represented as a list of lists in the same way as the Bergman subcomplex $B$.
The faces of $\operatorname*{dual}\left(  B\left(  I\right)  \right)  $ are
represented by the convex package type face of a polyhedron.

\item $\mathsf{CombinatorialDualization}\left(  posetNabla,B\right)  $

If $B$ is a subcomplex of the complex $posetNabla=\mathsf{FacePoset}\left(
\nabla\right)  $, then the poset of dual faces $F^{\ast}\subset\nabla^{\ast}$
is returned. It is represented as a list of lists in the same way as $B$ and
the faces are represented by the convex package type face of a polyhedron.

Suppose $B$ is the Bergman subcomplex $B\left(  I\right)  $ of
$\mathsf{FacePoset}\left(  \nabla\right)  $ as returned by
$\mathsf{BergmanSubcomplexOfSectionsPolytope}$, then by Proposition
\ref{Prop dual face combinatorial dualization general setting} we have
\[
\mathsf{CombinatorialDualization}\left(  posetNabla,B\right)
=\mathsf{DualComplex}\left(  A,v,t,inI,B\right)
\]
where $inI$ is the complex of initial ideals as returned by \linebreak%
$\mathsf{ComplexOfInitialIdeals}$.

\item $\mathsf{EqualityofFaceComplexes}\left(  B_{1},B_{2}\right)  $

Given two complexes $B_{1}$ and $B_{2}$ represented as a list of list of faces
of the same polyhedron, returns $true$ if $B_{1}=B_{2}$, i.e., if the lists in
each dimension agree up to permutation, otherwise returns $false$.

\item $\mathsf{MirrorComplex}\left(  P,A,v,t,inI,B\right)  $

Suppose the ideal $I$ defines a degeneration of complete intersections, $B$
denotes the Bergman subcomplex $B\left(  I\right)  $ of
$posetNabla=\mathsf{FacePoset}\left(  \nabla\right)  $ as returned by
$\mathsf{BergmanSubcomplexOfSectionsPolytope}$ and $inI$ is the complex of
initial ideals as given by $\mathsf{ComplexOfInitialIdeals}$. Then the
function $\mathsf{MirrorComplex}$ returns the complex $\mu\left(  B\left(
I\right)  \right)  \subset\mathsf{FacePoset}\left(  \Delta\right)  $ where
$\Delta=P^{\ast}$ as defined in Section
\ref{Sec mirror complex complete intersection}. It is represented as a list of
lists in the same way as $B$. The faces of $\mu\left(  B\left(  I\right)
\right)  $ are represented by the convex package type face of a polyhedron.

\item $\mathsf{LimitComplex}\left(  Iirr,A,v,t,B\right)  $

If $B$ denotes the Bergman subcomplex as returned by the function
$\mathsf{BergmanSubcomplexOfSectionsPolytope}$ and $Iirr$ the irrelevant ideal
of $X\left(  \Sigma\right)  $ as returned by the function
$\mathsf{IrrelevantIdeal}$, then the limit complex $\lim\left(  B\right)  $ is
computed as described in Section \ref{Sec limit map general setting}. The
limit complex is represented as a list of lists in the same way as $B$. The
strata corresponding to the face $F$ is represented by the ideal $I_{F}$ as
defined in Section \ref{Sec limit map general setting} and $I_{F}$ is given by
a list with minimal generators.

\begin{example}
\label{Ex limit Cox arc}Consider the monomial degeneration $\mathfrak{X}$ of
the complete intersection of two general quadrics in $\mathbb{P}^{3}$ defined
by the ideal $I\subset\mathbb{C}\left[  t\right]  \otimes\mathbb{C}\left[
x_{0},...,x_{3}\right]  $ as considered in the examples in Section
\ref{Sec tropical mirror construction for complete intersections}. To give a
computation of the ideal $I_{F}$ of $\lim\left(  F\right)  $ in a non
simplicial setting, we consider a face $F\in\mu\left(  B\left(  I\right)
\right)  \cong\operatorname*{Strata}\nolimits_{\Delta}\left(  I_{0}\right)  $
of the Bergman complex of the mirror. We have%
\[
\Delta=\operatorname*{convexhull}\left\{  \left(  3,-1,-1\right)  ,\left(
-1,3,-1\right)  ,\left(  -1,-1,3\right)  ,\left(  -1,-1,-1\right)  \right\}
\]
and
\[
\mu\left(  B\left(  I\right)  \right)  =\left(
\begin{array}
[c]{c}%
\left\{  \left\{  \left(  -1,3,-1\right)  \right\}  ,\left\{  \left(
-1,-1,-1\right)  \right\}  ,\left\{  \left(  -1,-1,3\right)  \right\}
,\left\{  \left(  3,-1,-1\right)  \right\}  \right\} \\
\left\{
\begin{array}
[c]{c}%
\left\{  \left(  -1,3,-1\right)  ,\left(  -1,-1,-1\right)  \right\}  ,\left\{
\left(  -1,-1,-1\right)  ,\left(  3,-1,-1\right)  \right\}  ,\\
\left\{  \left(  -1,3,-1\right)  ,\left(  -1,-1,3\right)  \right\}  ,\left\{
\left(  3,-1,-1\right)  ,\left(  -1,-1,3\right)  \right\}
\end{array}
\right\}
\end{array}
\right)
\]
omitting empty dimensions. The rays of $\Sigma^{\circ}$ are the rows of the
presentation matrix%
\[
A^{\circ}=\left(
\begin{array}
[c]{ccc}%
0 & 2 & -1\\
0 & 0 & -1\\
2 & 0 & -1\\
0 & 0 & 1\\
-1 & 1 & 0\\
-1 & -1 & 0\\
1 & -1 & 0\\
-1 & -1 & 2
\end{array}
\right)
\]
of $A_{2}\left(  Y^{\circ}\right)  $ and we denote the corresponding variables
of the Cox ring $S^{\circ}$ by $y_{1},...,y_{8}$. The irrelevant ideal of
$Y^{\circ}$ is%
\begin{align*}
B\left(  \Sigma^{\circ}\right)   &  =\left\langle y_{1},y_{7}\right\rangle
\cap\left\langle y_{3},y_{5}\right\rangle \cap\left\langle y_{2}%
,y_{8}\right\rangle \cap\left\langle y_{4},y_{6}\right\rangle \cap\left\langle
y_{2},y_{4}\right\rangle \cap\left\langle y_{5},y_{7}\right\rangle \\
&  \cap\left\langle y_{1},y_{3},y_{8}\right\rangle \cap\left\langle
y_{1},y_{3},y_{6}\right\rangle \cap\left\langle y_{1},y_{6},y_{8}\right\rangle
\cap\left\langle y_{3},y_{6},y_{8}\right\rangle
\end{align*}

Let $F$ be the face%
\begin{align*}
F  &  =\operatorname*{convexhull}\left\{  \left(  -1,3,-1\right)  ,\left(
-1,-1,3\right)  \right\} \\
w  &  =\left(  -1,3,-1\right)  +\left(  -1,-1,3\right)  =\left(
-2,2,2\right)  \in\operatorname*{int}\left(  F\right)
\end{align*}
and $a\left(  t\right)  \in\left(  K^{\ast}\right)  ^{3}$ with
$\operatorname*{val}\left(  a\left(  t\right)  \right)  =w$ and $c\left(
t\right)  =c_{J}\cdot t^{J}+\operatorname*{hot}$ a Cox arc representing
$a\left(  t\right)  $, i.e., $J^{t}+\ker\left(  A^{t}\right)  $ is the space
of solutions of the linear system of equations $A^{t}J^{t}=w^{t}$. The
intersection with this affine space with the positive orthant has the minimal
$0$-dimensional strata%
\[
\left(  2,0,0,0,0,0,0,2\right)  ,\left(  0,0,0,2,2,0,0,0\right)  ,\left(
0,0,2,0,4,0,0,2\right)  ,\left(  2,0,0,4,0,2,0,0\right)
\]
corresponding to the ideals%
\[
\left\langle y_{1},y_{8}\right\rangle ,\left\langle y_{4},y_{5}\right\rangle
,\left\langle y_{3},y_{5},y_{8}\right\rangle ,\left\langle y_{1},y_{4}%
,y_{6}\right\rangle \subset S^{\circ}%
\]
As $B\left(  \Sigma^{\circ}\right)  \subset\left\langle y_{3},y_{5}%
,y_{8}\right\rangle $ and $B\left(  \Sigma^{\circ}\right)  \subset\left\langle
y_{1},y_{4},y_{6}\right\rangle $ the limit face $\lim\left(  F\right)  $ is
given by any of the ideals%
\[
\left\langle y_{1},y_{8}\right\rangle ,\text{ }\left\langle y_{4}%
,y_{5}\right\rangle
\]
hence also by%
\[
\left\langle y_{1},y_{4},y_{5},y_{8}\right\rangle
\]
the ideal of all facets of $\nabla$ containing $\lim\left(  F\right)  $.
Figure \ref{Fig nabla22vars} shows the limit face $\lim\left(  F\right)
\subset B\left(  I\right)  \subset\nabla$ and the numbering of the facets of
$\nabla$ by the variables of $S^{\circ}$.
\end{example}

\item $\mathsf{DualLimitComplex}\left(  P,B\right)  $

Given the Bergman subcomplex $B\subset\nabla$, this function computes for all
faces $F$ the faces $H$ of $P$ such that $\dim\left(  F\cap H\right)
=\dim\left(  F\right)  $ and among those returns the set of minimal faces with
respect to inclusion. The output forms the dual limit complex of $B$ and is
represented as a list of lists compatible to $B$.

\item $\mathsf{LatticePoints}\left(  dB\right)  $

Returns the set of lattice points of the complex $dB\subset P\subset
\mathbb{Z}^{n}\otimes\mathbb{R}$, which is represented as a list of lists and
each face $F$ of $dB$ is a face of the polyhedron $P$.

\item $\mathsf{DeformationsFromCombinatorialData}\left(  A,v,dB\right)  $

Suppose $dB=\operatorname*{dual}\left(  B\left(  I\right)  \right)  $ as given
by the function $\mathsf{DualComplex}$. The function
$\mathsf{DeformationsFromCombinatorialData}$ returns the complex of first
order deformations of $X_{0}$, represented as a list of lists in the same way
as $\operatorname*{dual}\left(  B\left(  I\right)  \right)  $. Each face is a
set of degree $0$ Cox Laurent monomials corresponding to the lattice points of
the corresponding face of $\operatorname*{dual}\left(  B\left(  I\right)
\right)  $.

\item $\mathsf{FirstOrderDegenerationFromCombinatorialData}\left(
I0,defs,c\right)  $

Suppose $I_{0}$ is the ideal of the special fiber of $\mathfrak{X}$ given as a
list $I0$ of monomial generators and $defs$ is the complex of first order
deformations as returned by the complex
$\mathsf{DeformationsFromCombinatorialData}$ then the list%
\[
\left[  m+\sum_{\alpha\text{ in a face of }defs}t\cdot c_{\alpha}\cdot
\alpha\left(  m\right)  \mid m\in I0\right]
\]
is returned.

\item $\mathsf{ExtendFirstOrderPfaffian}\left(  gI\right)  $

If the ideal $I^{1}\subset\mathbb{C}\left[  t\right]  /\left\langle
t^{2}\right\rangle \times S$ generated by $gI$ is Pfaffian with syzygy matrix
$A$, a list with the Pfaffians of $A$ in $\mathbb{C}\left[  t\right]  \times
S$ is returned to extend the family defined by $I^{1}$ as explained in Section
\ref{Sec deformations of Pfaffians}.

\item $\mathsf{ModuliDimStanleyReisner}\left(  dB\right)  $

The function $\mathsf{ModuliDimStanleyReisner}$ computes the number of lattice
points of $dB=\operatorname*{dual}\left(  B\left(  I\right)  \right)  $ as
given by $\mathsf{DualComplex}$ and the number $r_{0}$ of roots of $Y$ among
them, and returns
\[
\left\vert \operatorname*{dual}\left(  B\left(  I\right)  \right)  \cap
M\right\vert -\dim Y-r_{0}%
\]
If $Y=\mathbb{P}\left(  \Delta\right)  $ where $\Delta$ is a Veronese polytope
of $\mathbb{P}^{n}$ this number is $h^{1,\dim\left(  X\right)  -1}\left(
X\right)  $ of the general fiber $X$ of $\mathfrak{X}$, as discussed in
Section \ref{Sec First approximation of a tropical Hodge formula}.

\item $\mathsf{StringyEFunctionOfGorensteinCone}\left(  C\right)  $

Given a Gorenstein cone $C$ returns%
\[
\left(  uv\right)  ^{-1}\sum_{C_{1}\subset C}\left(  -u\right)  ^{\dim\left(
C_{1}\right)  }\tilde{S}\left(  C_{1},u^{-1}v\right)  \tilde{S}\left(
C_{1}^{\vee},uv\right)
\]

\end{itemize}

%

\begin{figure}
[h]
\begin{center}
\includegraphics[
height=2.1006in,
width=2.0989in
]%
{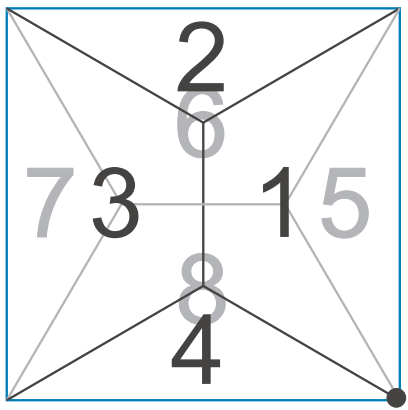}%
\caption{$\lim\left(  F\right)  \subset\nabla$ and the numbering of the facets
of $\nabla$ by the variables of the Cox ring for the complete intersection of
two general quadrics in $\mathbb{P}^{3}$}%
\label{Fig nabla22vars}%
\end{center}
\end{figure}

\begin{example}
Let $\mathfrak{X}\subset\mathbb{P}\left(  \Delta\right)  \times
\operatorname*{Spec}\mathbb{C}\left[  \left[  t\right]  \right]  $ be the
monomial degeneration of an
\index{elliptic curve!complete intersection}%
elliptic curve given as the complete intersection of two general quadrics in
$\mathbb{P}^{3}$, as considered in the examples in Section
\ref{Sec tropical mirror construction for complete intersections}. So let
$\mathfrak{X}$ be given by the ideal
\[
I=\left\langle t\cdot g_{1}+x_{1}x_{2},\ t\cdot g_{2}+x_{0}x_{3}\right\rangle
\subset\mathbb{C}\left[  t\right]  \otimes S
\]
where $g_{1},g_{2}\in S=\mathbb{C}\left[  x_{0},...,x_{3}\right]  $ are
general quadrics reduced with respect to $I_{0}=\left\langle x_{1}x_{2}%
,x_{0}x_{3}\right\rangle $.\textrm{\texttt{\medskip}}\newline%
\texttt{runM2:="M2":\smallskip\newline runSingular:="Singular":\smallskip
\newline stdSystem:="M2":\smallskip\newline
pathConvex:="/usr/local/convex":\smallskip}\newline%
\texttt{read("tropicalmirror"):}\textrm{\texttt{\medskip}}\newline%
\texttt{v:=[x0,x1,x2,x3]:\smallskip}\newline%
\texttt{A:=matrix([[-1,-1,-1],[1,0,0],[0,1,0],[0,0,1]]):\smallskip\newline
P:=convexhull(op(convert(A,listlist))):}%
\[
P:=polytope(3,3,4,4)
\]
\texttt{Sigma:=FanOverFaces(P):}%
\[
Sigma:=FAN(3,3,0,4,[0,0,4])
\]
\texttt{g1:=RandomPolynomial(A,v,matrix([[2],[0],[0],[0]]),13):\smallskip
\newline g2:=RandomPolynomial(A,v,matrix([[2],[0],[0],[0]]),13):\smallskip
\newline gI:=[t*g1+x1*x2,t*g2+x0*x3]:\smallskip\newline
gI:=ReduceGenerators(v,t,gI);}%
\begin{gather*}
gI:=[2tx0^{2}+tx1x0+8tx2x0+tx1^{2}+tx3x1+9tx3^{2}+12tx2x3+8tx2^{2}+x1x2,\\
15tx0^{2}+8tx1x0+11tx2x0+16tx1^{2}+16tx3x1+14tx2x3+15tx2^{2}+7tx3^{2}+x0x3]
\end{gather*}
\texttt{gI1:=AssociatedFirstOrderDegeneration(v,t,gI):\smallskip\newline
C:=SpecialFiberGroebnerCone(A,v,t,gI);}%
\[
C:=cone(4,4,0,8,8)
\]
\texttt{rays(C);}%
\begin{align*}
&  [[1,1,0,1],[1,0,1,1],[1,0,-1,-1],[1,-1,0,-1],\\
&  [1,1,0,0],[1,0,1,0],[1,0,0,1],[1,-1,-1,-1]]
\end{align*}
\texttt{Nabla:=AssociatedAnticanonicalSectionsPolytope(C);}%
\[
\nabla:=polytope(3,3,8,8)
\]
\texttt{vertices(Nabla);}%
\[
\lbrack\lbrack
1,0,1],[0,1,1],[0,-1,-1],[-1,0,-1],[0,0,1],[-1,-1,-1],[0,1,0],[1,0,0]]
\]
\texttt{PMirror:=AssociatedFanoPolytope(C):\smallskip\newline
vertices(PMirror);}%
\[
\lbrack\lbrack
-1,-1,0],[-1,-1,2],[0,2,-1],[2,0,-1],[0,0,-1],[0,0,1],[1,-1,0],[-1,1,0]]
\]
\texttt{SigmaMirror:=FanOverFaces(PMirror):\smallskip}\newline\texttt{\# By
fixing the presentation matrix AMirror of the Chow group of }\newline%
\texttt{\# YMirror=X(SigmaMirror), we choose a numbering of the Cox variables
}\newline\texttt{\# of YMirror compatible with Example \ref{Ex limit Cox arc}%
.\smallskip\newline AMirror:=matrix([[0,2,-1], [0,0,-1], [2,0,-1], [0,0,1],
[-1,1,0], [-1,-1,0], [1,-1,0], [-1,-1,2]]):\smallskip\newline
ChowGroup(AMirror);}%
\[
\left[  \left[

\]
\texttt{\# Numbers of faces of dimensions -1,0,1,2,3 are 1,8,14,8,1\smallskip
}\newline\texttt{VertexRepresentation}(\texttt{posetNabla});\newline%
\[%
%
\]
\texttt{\# Numbers of faces of dimensions -1,0,1,2,3 are 0,4,4,0,0\smallskip
}\newline\texttt{y:=[y1,y2,y3,y4,y5,y6,y7,y8];\smallskip}\newline%
\texttt{I0Mirror:=SpecialFiberIdeal(AMirror,y,posetNabla,B);}%
\begin{align*}
I0mirror  &  :=[y5y6y7y8,y1y6y7y8,y2y5y7y8,y1y2y7y8,y3y5y6y8,y1y3y6y8,\\
&  y2y3y5y8,y1y2y3y8,y4y5y6y7,y1y4y6y7,y2y4y5y7,y1y2y4y7,y3y4y5y6,\\
&  y1y3y4y6,y2y3y4y5,y1y2y3y4]
\end{align*}
\texttt{\newline inC:=ComplexOfInitialIdeals(v,t,gI,C):\smallskip\newline
dualB:=DualComplex(A,v,t,inI,B):\smallskip}\newline\texttt{\# Numbers of faces
of dimensions -1,0,1,2,3 are 0,0,4,4,0\smallskip}\newline%
\texttt{VertexRepresentation}(\texttt{dualB})\texttt{;}%
\[%

\]
\texttt{BMirror:=MirrorComplex(P,A,v,t,inI,B):\smallskip}\newline\texttt{\#
Numbers of faces of dimensions -1,0,1,2,3 are 0,4,4,0,0\smallskip}%
\newline\texttt{VertexRepresentation}(\texttt{BMirror})\texttt{;}%
\[%
%
\]
\texttt{cdualBMirror:=CombinatorialDualization(posetDelta,BMirror):\smallskip
}\newline\texttt{duallimB:=DualLimitComplex(P,B);\smallskip}\newline\texttt{\#
Numbers of faces of dimensions -1,0,1,2,3 are 0,0,4,4,0\smallskip}%
\newline\texttt{EqualityofFaceComplexes(duallimB,cdualBMirror);}%
\[
true
\]
\texttt{Covering(duallimB);}%
\[%
%
\]
\texttt{FirstOrderDegenerationFromCombinatorialData(I0Mirror,DefsMirror):\smallskip
}\newline\texttt{I0mirrorShort:=[y1*y2*y3*y4,y5*y6*y7*y8]:\smallskip}%
\newline%
\texttt{ts1:=ToricStrataDecomposition(posetNabla,I0mirrorShort):\smallskip
}\newline\texttt{ts2:=ToricStrataDecomposition(posetNabla,I0mirror):\smallskip
}\newline\texttt{EqualityofFaceComplexes(ts1,ts2);}%
\[
true
\]
\texttt{FirstOrderDegenerationFromCombinatorialData(I0mirrorShort,DefsMirror):}%
\[
\left[  t\cdot y2^{2}y6^{2}+t\cdot y4^{2}y8^{2}+y1y2y3y4,\,t\cdot y1^{2}%
y5^{2}+t\cdot y3^{2}y7^{2}+y5y6y7y8\right]
\]

\end{example}

\section{Perspectives\label{Sec perspectives}}

\subsection{Tropical computation of string
cohomology\label{Sec tropical computation of the stringy E function and string cohomology}%
}

Stringy $E$-functions and tropical geometry share many technical concepts, for
example formal arcs, which are used in motivic integration to prove key
theorems on stringy $E$-functions, as indicated in Section
\ref{Mirror symmetry for singular Calabi-Yau varieties and stringy Hodge numbers}%
. As explained in Section
\ref{Sec String theoretic Hodge formula for hypersurfaces}, Batyrev and
Borisov give in \cite{BB Mirror duality and stringtheoretic Hodge numbers} a
formula to compute the stringy $E$-function for a general anticanonical toric
hypersurface. The stringy $E$-function of a hypersurface in $\mathbb{P}\left(
\Delta\right)  $ is computed in terms of combinatorial data of the reflexive
polytope $\Delta$. As explained in Section
\ref{Sec string theoretic Hodge formula for complete intersections}, Batyrev
and Borisov represent complete intersections in $Y=\mathbb{P}\left(
\Delta\right)  $ given by general sections $s_{j}\in H^{0}\left(
Y,\mathcal{O}_{Y}\left(  E_{j}\right)  \right)  $ with $\sum_{j=1}^{c}%
E_{j}=-K_{Y}$ as the zero set of a section of%
\[
Z=\mathbb{P}\left(  \mathcal{O}_{Y}\left(  D_{1}\right)  \oplus...\oplus
\mathcal{O}_{Y}\left(  D_{c}\right)  \right)  \rightarrow Y
\]
so reducing the computation of the stringy $E$-function to the case of a
hypersurface and the bundle $Z$.

The tropical mirror construction, as given in Section
\ref{Sec tropical mirror construction}, provides the additional data given by
the complexes $B\left(  I\right)  \subset\operatorname*{Poset}\left(
\nabla\right)  $ and $\lim\left(  B\left(  I\right)  \right)  \subset
\operatorname*{Poset}\left(  \Delta\right)  $. So it is natural to ask whether
it is possible to give a direct formula for the stringy $E$-function of the
general complete intersection Calabi-Yau in terms of this data, as indicated
in Section
\ref{Sec first approximation of a tropical computation of the stringy E function}%
.

Further evidence to expect that the stringy $E$-function $E_{st}\left(
X\right)  $ of $X$ should be computable from the tropical data, is given the
fact that the special fiber $X_{0}$ of $\mathfrak{X}$ is a union of toric
varieties and, as noted in Section \ref{Sec tropical stringy E}, the stringy
$E$-function respects stratifications. See also Section
\ref{Sec Stringy Efunction of toric variety} for an explicit formula of the
stringy $E$-function of a toric variety.

So, in general one may ask if the stringy $E$-function of the generic element
of the degeneration $\mathfrak{X}$ is computable in terms of the tropical data.

\subsection{Hilbert schemes and moduli
spaces\label{Sec Hilbert schemes and moduli spaces}}

The multigraded Hilbert scheme described in Sections
\ref{Sec Grassmann functor}-\ref{Sec Multigraded Hilbert schemes} for smooth
toric varieties may be generalized to the non smooth and further to the non
simplicial setup by using the ideas of Section
\ref{Monomial ideals in the Cox ring and the stratified toric primary decomposition}%
. For reduced monomial ideals $I_{0}$ the saturation $\left(  I_{0}:B\left(
\Sigma\right)  ^{\infty}\right)  $ is generalized by the ideal $I_{0}^{\Sigma
}$ in the non-simplicial setup: Consider a toric variety $Y=X\left(
\Sigma\right)  $ given by the Fano polytope $P\subset N_{\mathbb{R}}$,
$\Sigma=\Sigma\left(  P\right)  $, let $\Delta=P^{\ast}\subset M_{\mathbb{R}}$
and $S=\mathbb{C}\left[  y_{r}\mid r\in\Sigma\left(  1\right)  \right]  $ be
the Cox ring of $Y$. Given a monomial ideal $I_{0}\subset S$ we associate to
$I_{0}$ the complex of strata $\operatorname*{Strata}\nolimits_{\Delta}\left(
I_{0}\right)  $ with faces of dimension $s$ given by%
\[
\operatorname*{Strata}\nolimits_{\Delta}\left(  I_{0}\right)  _{s}=\left\{
F\mid%
\begin{tabular}
[c]{l}%
$F\text{ a face of }\Delta\text{ of }\dim\left(  F\right)  =s\text{ with}$\\
$\operatorname*{facets}\nolimits_{F}\left(  \Delta\right)  \cap
\operatorname*{facets}\nolimits_{m}\left(  \Delta\right)  \neq\varnothing$\\
for all monomial $m\in I_{0}$%
\end{tabular}
\right\}
\]
where%
\begin{align*}
\operatorname*{facets}\nolimits_{F}\left(  \Delta\right)   &  =\left\{  G\mid
G\text{ facet of }\Delta\text{ with }F\subset G\right\} \\
\operatorname*{facets}\nolimits_{m}\left(  \Delta\right)   &  =\left\{  G\mid
G\text{ facet of }\Delta\text{ with }y_{G^{\ast}}\mid m\right\}
\end{align*}
To this complex we can associate the ideal%
\[
I_{0}^{\Sigma}=\left\langle
{\displaystyle\prod\limits_{v\in J}}
y_{v}\mid J\subset\Sigma\left(  1\right)  \text{ with }\operatorname*{supp}%
\left(  B\left(  I\right)  \right)  \subset%
{\displaystyle\bigcup\limits_{v\in J}}
F_{v}\right\rangle \subset S
\]
Replacing $I_{0}$ by $I_{0}^{\Sigma}$ also amounts to passing to a reduced
ideal, so without modification this can only work for the local Hilbert scheme
around a reduced point, but this is what is relevant for the tropical mirror construction.

Generalizing the multigraded Hilbert scheme, also the state polytope, as
defined in Section \ref{Sec State polytope general toric}, can be generalized
to the non simplicial setup.

Local moduli stacks may be computed by taking the quotient of the local
Hilbert scheme around a given ideal by a (in general non-reductive)
automorphism group, generalizing the ideas of Section
\ref{Sec First approximation of a tropical Hodge formula}. For geometric
invariant theory in the non-reductive setting see \cite{DoKi Towards
nonreductive geometric invariant theory}.

Let $\Delta\subset M_{\mathbb{R}}$ be a reflexive polytope and consider an
anticanonical toric hypersurface in $Y=\mathbb{P}\left(  \Delta\right)  $.
Generalizing the ideas from Batyrev%
\'{}%
s Hodge formulas for toric hypersurfaces, as outlined in Section
\ref{Sec Batyrevs Hodge formula}, the subset
\[
\Xi_{0}=\Delta\cap M-\left\{  0\right\}  -\bigcup_{\operatorname*{codim}%
Q=1}\operatorname*{int}\nolimits_{M}\left(  Q\right)
\]
of the set of lattice points of the polytope of sections $\Delta$ of $-K_{Y}$
plays the key role in the construction of the simplified complex moduli space
of anticanonical hypersurfaces. By mirror symmetry this set is also used in
the construction of the K\"{a}hler moduli space of the mirror. See, e.g.,
\cite{AGM The MonomialDivisor Mirror Map} and \cite[Section 6]{CK Mirror
Symmetry and Algebraic Geometry} for details. Now consider the setup of the
tropical mirror construction with degeneration $\mathfrak{X}$ given by $I$.
Then the set $\Xi_{0}$ generalizes to the set of lattice points of the support
of $\operatorname*{dual}\left(  B\left(  I\right)  \right)  \subset
\operatorname*{Poset}\left(  \nabla^{\ast}\right)  $ which do not correspond
to a trivial deformations of the special fiber. Note that for hypersurfaces
$\nabla^{\ast}=\Delta$. So one may ask if the complex $\operatorname*{dual}%
\left(  B\left(  I\right)  \right)  $ may be useful for the construction of
moduli spaces.

\subsection{Integrally affine
structures\label{Sec integrally affine structures}}

Gross and Siebert use in \cite{GrSi Affine manifolds log structures and mirror
symmetry}, \cite{GrSi Mirror symmetry via logarithmic degeneration data I} and
\cite{GrSi From real affine geometry to complex geometry} toric degenerations,
integrally affine structures and the discrete Legendre transform to give a
mirror construction.

They consider a degeneration $f:\mathfrak{X}\rightarrow S$ of Calabi-Yau
varieties whose total space is a complex analytic space its base is a complex
disc $S$ and which has a special fiber whose normalization is a disjoint union
of toric varieties. Outside a set of codimension $2$ any point $x$ in the
total space is assumed to have a neighborhood $U_{x}$ such that there is an
affine toric variety $Y_{x}$, a regular function $f_{x}$ given by a monomial
and a commutative diagram%
\[%
\begin{tabular}
[c]{lll}%
$U_{x}$ & $\rightarrow$ & $Y_{x}$\\
$\downarrow f\mid_{U_{x}}$ &  & $\downarrow f_{x}$\\
$S$ & $\rightarrow$ & $\mathbb{C}$%
\end{tabular}
\
\]
with open embeddings $U_{x}\rightarrow Y_{x}$ and $S\rightarrow\mathbb{C}$.

Given this data a manifold with an integral affine structure with
singularities is constructed, i.e., a topological manifold such that outside a
finite union of locally closed submanifolds of codimension at least $2$ there
are charts whose transition functions are integral affine transformations.
Furthermore a polyhedral decomposition of this manifold is constructed. This
process is reversible, i.e., from an integral affine manifold with
singularities and a polyhedral decomposition one can construct a degeneration.
Using the descrete Legendre transform one can obtain a mirror integral affine
manifold with polyhedral decomposition from which one obtains the mirror degeneration.

So one may ask if it is possible to obtain from the polytopes $\Delta$ and
$\nabla$ and the embedded subcomplexes $B\left(  I\right)  \subset
\operatorname*{Poset}\left(  \nabla\right)  $ and $\lim\left(  B\left(
I\right)  \right)  \subset\operatorname*{Poset}\left(  \Delta\right)  $
integral affine structures and how the tropical mirror construction relates to
the construction by Gross and Siebert.

\subsection{Torus fibrations\label{Sec torus fibrations}}

Consider the setup of the tropical mirror construction. So let $Y$ be a toric
Fano variety given by a Fano polytope $P\subset N_{\mathbb{R}}$ with Cox ring
$S$ and a monomial degeneration $\mathfrak{X}$ of Calabi-Yau varieties of
dimension $d$ which is given by the ideal $I\subset S\otimes\mathbb{C}\left[
t\right]  $.

The monomial special fiber $X_{0}$ has a degenerate torus fibration over the
sphere $\operatorname*{Strata}\left(  X_{0}\right)  \cong\lim\left(  B\left(
I\right)  \right)  \subset\operatorname*{Poset}\left(  \Delta\right)  $ with
$\Delta=P^{\ast}$: The strata of $X_{0}$ of dimension $s=0,...,d$ are complex
tori $\left(  \mathbb{C}^{\ast}\right)  ^{s}$ which contain $\left(
S^{1}\right)  ^{s}$.

In the same way also the mirror special fiber $X_{0}^{\circ}$ has a degenerate
torus fibration over $B\left(  I\right)  $. The dimensions of the tori and
their base faces are related via the map $\lim$ by $s\leftrightarrow d-s$.
This agrees with the large-small interchange of $T$-duality.

One may ask if the data provided by the spheres $B\left(  I\right)
\subset\operatorname*{Poset}\left(  \nabla\right)  $ and $\lim\left(  B\left(
I\right)  \right)  \subset\operatorname*{Poset}\left(  \Delta\right)  $ and
the complexes of deformations $\operatorname*{dual}\left(  B\left(  I\right)
\right)  $ and $\lim\left(  B\left(  I\right)  \right)  ^{\ast}$ in the
initial ideals provide sufficient information to obtain a SYZ-fibration of the
general fiber of $\mathfrak{X}$ as introduced in\linebreak\cite{SYZ Mirror
symmetry is Tduality}. Furthermore one may ask how these fibrations relate to
mirror symmetry via $T$-duality.

\subsection{Mirrors for further Stanley-Reisner Calabi-Yau
degenerations\label{1FurtherStanleyReisnerexamples}}

In Section
\ref{Monomial ideals in the Cox ring and the stratified toric primary decomposition}
we connected the combinatorial represenation of monomial ideals in the Cox
ring of a toric variety $Y$ in the special case of $Y=\mathbb{P}^{n}$ to the
Stanley-Reisner setup:

Let $Y=\mathbb{P}\left(  \Delta\right)  \cong\mathbb{P}^{n}$ where $\Delta$ is
the degree $n+1$ Veronese polytope and let $\Sigma$ be the normal fan of
$\Delta$. Let $\mathcal{R}=\Sigma\left(  1\right)  $ be the set of rays of
$\Sigma$ and $S=\mathbb{C}\left[  y_{r}\mid r\in\mathcal{R}\right]  $ the
homogeneous coordinate ring of $Y$. The faces of the simplex $\Delta^{\ast}$
correspond to the subsets of $\Sigma\left(  1\right)  $, i.e., the set of
variables of the homogeneous coordinate ring $S$.

Let $Z$ be a simplicial subcomplex of $\operatorname*{Poset}\left(
\Sigma\right)  $. In the following we represent the faces as subsets of
$\mathcal{R}$, so $\operatorname*{Poset}\left(  \Sigma\right)  $ is
represented as the complex $2^{\mathcal{R}}$ of all subsets of $\mathcal{R}$
and $Z\subset2^{\mathcal{R}}$. Any face $F$ of $\operatorname*{Poset}\left(
\Sigma\right)  \cong2^{\mathcal{R}}$ can be considered as a square free
monomial
\[
y_{F}=%
{\textstyle\prod\nolimits_{r\in F}}
y_{r}%
\]
Then the monomial ideal generated by the $y_{F}$ for the non-faces of $Z$
\[
I_{Z}=\left\langle
{\textstyle\prod\nolimits_{r\in F}}
y_{r}\mid F\in2^{\mathcal{R}}\text{ not a face of }Z\right\rangle \subset S
\]
is the
\index{Stanley-Reisner ideal}%
Stanley-Reisner ideal corresponding to $Z\subset2^{\mathcal{R}}$ and
$A_{Z}=S/I_{Z}$ is the
\index{Stanley-Reisner ring}%
Stanley-Reisner ring of $Z$. Note that for any monomial ideal $I_{Z}\subset S$
the ring $S/I_{0}$ is $\mathbb{Z}^{\mathcal{R}}$-graded. The complex $Z$
defines the affine scheme $\mathbb{A}_{Z}=\operatorname*{Spec}\left(
A_{Z}\right)  $ and the projective scheme $\mathbb{P}_{Z}=\operatorname*{Proj}%
\left(  A_{Z}\right)  $.

The complex $Z$ relates to the complex of strata $\operatorname*{Strata}%
\nolimits_{\Delta}\left(  I_{Z}\right)  $, which we defined in Section
\ref{Monomial ideals in the Cox ring and the stratified toric primary decomposition}%
, by the isomorphism of complexes%
\[%
\begin{tabular}
[c]{cccl}
& $\operatorname*{Poset}\left(  \Delta\right)  $ & $\overset{\cong
}{\rightarrow}$ & $2^{\mathcal{R}}$\\
& $\cup$ &  & $\cup$\\
$\operatorname*{comp}:$ & $\operatorname*{Strata}\nolimits_{\Delta}\left(
I_{Z}\right)  $ & $\overset{\cong}{\rightarrow}$ & $Z$\\
& $F$ & $\mapsto$ & \multicolumn{1}{c}{$\left\{  r\in\mathcal{R}\mid
r\not \subset \operatorname*{hull}\left(  F^{\ast}\right)  \right\}  $}%
\end{tabular}
\]

In \cite{AC Cotangent cohomology of StanleyReisner rings} and \cite{AC
Deforming StanleyReisner rings} the deformation theory of Stanley-Reisner
rings is addressed, computing the first order deformations and obstructions.
We give a short outline of the computation of the first order deformations of
$\mathbb{P}_{Z}\subset\mathbb{P}^{n}$ by Altmann and Christophersen.

Consider the following notation. Given a subcomplex $Z\subset2^{\mathcal{R}}$
we denote by
\[
\operatorname*{vert}\left(  Z\right)  =\left\{  r\in\mathcal{R}\mid\left\{
r\right\}  \in Z\right\}
\]
the set of vertices of $Z$. If $a\in2^{\mathcal{R}}$ is a face then we can
define complex of faces of $a$ as%
\[
\operatorname*{Poset}\left(  a\right)  =\left\{  b\in2^{\mathcal{R}}\mid
b\subset a\right\}
\]
the boundary of $a$ as%
\[
\partial a=\left\{  b\in2^{\mathcal{R}}\mid b\subsetneqq a\right\}
\]
and the link of $a$ in $Z$%

\[
\operatorname{lk}\left(  a,Z\right)  =\left\{  b\in Z\mid b\cap a=\emptyset
\text{, }b\cup a\in Z\right\}
\]

An element $c\in\mathbb{Z}^{\mathcal{R}}$ has a support $\operatorname*{supp}%
\left(  c\right)  \in2^{\mathcal{R}}$ defined as%
\[
\operatorname*{supp}\left(  c\right)  =\left\{  r\in\mathcal{R}\mid c_{r}%
\neq0\right\}
\]

Let $S$ be a polynomial $K$-algebra mapping onto $A$ with $A\cong S/I$ for
some ideal $I$ and%
\[
0\rightarrow R\rightarrow F\rightarrow S\rightarrow A\rightarrow0
\]
with free $F$ a presentation of $A$ as an $S$-module. If $M$ is an $A$-module
then define%
\[
T^{1}\left(  A/K,M\right)  =\operatorname*{coker}\left(  \operatorname*{Der}%
\nolimits_{K}\left(  S,M\right)  \rightarrow\operatorname*{Hom}\nolimits_{A}%
\left(  I/I^{2},A\right)  \right)
\]

In the Stanley-Reisner setup write $T_{A_{Z}}^{1}=T^{1}\left(  A_{Z}%
/\mathbb{C},A_{Z}\right)  $. The grading of $A_{Z}$ induces a grading on
$T_{A_{Z}}^{1}$.

For $c\in\mathbb{Z}^{\mathcal{R}}$ homomorphisms in $\operatorname*{Hom}%
\nolimits_{S}\left(  I_{Z},S/I_{Z}\right)  _{c}$ can be represented by Cox
Laurent monomials.

Computation of $T^{1}$ reduces to links of faces:

\begin{theorem}
\cite{AC Cotangent cohomology of StanleyReisner rings} Let $D\in
\mathbb{Z}^{\mathcal{R}}$ and write $D=D_{+}-D_{-}$ where $D_{+},D_{-}%
\in\mathbb{Z}_{\geq0}^{\mathcal{R}}$ with disjoint support. Denote by
$a=\operatorname*{supp}\left(  D_{+}\right)  $ and $b=\operatorname*{supp}%
\left(  D_{-}\right)  $.

Then $T_{A_{Z},D}^{1}=0$ unless $a\in Z$, $D_{-}\in\left\{  0,1\right\}
^{\mathcal{R}}$ and $b\neq\emptyset$. Suppose these conditions are satisfied.
$T_{A_{Z},D}^{1}$ depends only on $a$ and $b$, so write $T_{A_{Z},a-b}^{1}$
for $T_{A_{Z},D}^{1}$. Then%
\[
T_{A_{Z},a-b}^{1}=0
\]
unless $a\in Z$ and $b\subset\operatorname*{vert}\left(  \operatorname{lk}%
\left(  a,Z\right)  \right)  $ and if these conditions are satisfied%
\[
T_{A,a-b}^{1}\left(  Z\right)  \cong T_{A,\emptyset-b}^{1}\left(
\operatorname{lk}\left(  a,T\right)  \right)
\]

\end{theorem}

Suppose that $Z$ is a
\index{combinatorial manifold|textbf}%
\textbf{combinatorial manifold}, i.e., for all faces $a\in Z$ the link
$\operatorname{lk}\left(  a,Z\right)  $ is a sphere of dimension $\dim\left(
Z\right)  -\dim\left(  a\right)  -1$.

\begin{lemma}
\cite{AC Cotangent cohomology of StanleyReisner rings} For $b\in
2^{\Sigma\left(  1\right)  }$ with $\left\vert b\right\vert \geq2$ it is equivalent:

\begin{itemize}
\item $T_{A,\varnothing-b}^{1}\left(  Z\right)  \neq\emptyset$

\item $\dim\left(  T_{A,\varnothing-b}^{1}\left(  Z\right)  \right)  =1$

\item It holds%
\[
Z=\left\{
\begin{tabular}
[c]{ll}%
$L\ast\partial b$ & $\text{if }b\notin Z$\\
$L\ast\partial b\cup\partial L\ast\operatorname*{Poset}\left(  b\right)  $ &
$\text{if }b\in Z$%
\end{tabular}
\right\}
\]
where the geometric realization of $L$ is a $\dim\left(  Z\right)
+1-\left\vert b\right\vert $ sphere. In any case $Z$ is a sphere.
\end{itemize}
\end{lemma}

\begin{theorem}
\cite{AC Cotangent cohomology of StanleyReisner rings} If $Z$ is a manifold,
then
\[
T_{A_{Z}}^{1}=\sum_{\substack{D\in\mathbb{Z}^{\mathcal{R}}%
\\a=\operatorname*{supp}\left(  D\right)  \in Z}}T_{<0}^{1}\left(
\operatorname{lk}\left(  a,Z\right)  \right)
\]
with%
\[
T_{<0}^{1}\left(  \operatorname{lk}\left(  F,Z\right)  \right)  =\sum
T_{\emptyset-b}^{1}\left(  \operatorname{lk}\left(  F,Z\right)  \right)
\]
where the sum goes over all $b\subset\operatorname{lk}\left(  F,Z\right)  $
with $\left\vert b\right\vert \geq2$ and%
\[
\operatorname{lk}\left(  F,Z\right)  =\left\{
\begin{tabular}
[c]{ll}%
$L\ast\partial b$ & $\text{if }b\notin\operatorname{lk}\left(  F,Z\right)  $\\
$L\ast\partial b\cup\partial L\ast\operatorname*{Poset}\left(  b\right)  $ &
$\text{if }b\in\operatorname{lk}\left(  F,Z\right)  $%
\end{tabular}
\right\}
\]
where the geometric realization $L$ is a $\dim\left(  \operatorname{lk}\left(
F,Z\right)  \right)  +1-\left\vert b\right\vert $ sphere.

Note that all $T_{\emptyset-b}^{1}\left(  \operatorname{lk}\left(  F,Z\right)
\right)  $ are of dimension one.
\end{theorem}

From this one can compute for case of manifolds of dimension $\leq2$:

\begin{proposition}
\cite{AC Cotangent cohomology of StanleyReisner rings} If $Z$ is a manifold of
dimension $\leq2$, then $T_{<}^{1}\left(  Z\right)  $ is trivial or%
\[%
\begin{tabular}
[c]{llll}%
$\dim\left(  Z\right)  $ & $Z$ &  & $\dim\left(  T_{<}^{1}\left(  Z\right)
\right)  $\\
$0$ & $\partial\Delta_{1}$ & two points & $1$\\
$1$ & $E_{3}$ & triangle & $4$\\
$1$ & $E_{4}$ & quadrangle & $2$\\
$2$ & $\partial\Delta_{3}$ & tetrahedron & $11$\\
$2$ & $\Sigma\left(  E_{3}\right)  $ & suspension of a triangle & $5$\\
$2$ & $\Sigma\left(  E_{4}\right)  $ & octahedron & $3$\\
$2$ & $\Sigma\left(  E_{m}\right)  $ & suspension of an $m$-gon, $m\geq5$ &
$1$\\
$2$ & $C\left(  m,2\right)  $, $m\geq6$ & cyclic polytope & $1$%
\end{tabular}
\]
Here $\Delta_{m}$ denotes the $m$-simplex, $E_{m}$ the $m$-gon and
$\Sigma\left(  C\right)  $ the suspension of $C$, i.e., the double pyramid on
$C$.
\end{proposition}

This can be applied to compute $T_{<}^{1}\left(  Z\right)  $ for the links in
a threefold.

\begin{proposition}
\cite{AC Cotangent cohomology of StanleyReisner rings} Given a simplicial
complex $Z\subset2^{\Sigma\left(  1\right)  }$ and the corresponding
Stanley-Reisner ideal $I_{Z}$ we have%
\[
T_{\mathbb{P}_{Z}/\mathbb{P}^{n}}^{1}=H^{0}\left(  \mathbb{P}_{Z}%
,N_{\mathbb{P}_{Z}/\mathbb{P}^{n}}\right)  \cong\operatorname*{Hom}%
\nolimits_{S}\left(  I_{Z},S/I_{Z}\right)  _{0}%
\]

\end{proposition}

The kernel of $\operatorname*{Hom}\nolimits_{S}\left(  I_{Z},S/I_{Z}\right)
_{0}\rightarrow T_{A_{Z},0}^{1}$ is generated by the homomorphisms $x_{r_{1}%
}\frac{\partial}{\partial x_{r_{2}}}$ and%
\[
\dim\left(  \operatorname*{Hom}\nolimits_{S}\left(  I_{Z},S/I_{Z}\right)
_{0}\right)  =\dim\left(  T_{A_{Z},0}^{1}\right)  +\left(  n+1\right)  ^{2}%
\]

As outlined above for $T^{1}$, the methods given by Altmann and Christophersen
allow computation of the first order deformations and obstructions, hence
should provide the necessary data to apply the tropical mirror construction
given in Section \ref{Sec tropical mirror construction}.

Let $X_{0}\subset Y$ be defined by a Stanley-Reisner ideal. We may consider,
if existent, a component of the local Hilbert scheme of $X_{0}$ such that a
degeneration $\mathfrak{X}$ with general tangent vector in this component
given by an ideal $I$ satisfies $C_{I_{0}}\left(  I\right)  \cap\left\{
w_{t}=0\right\}  =\left\{  0\right\}  $. The first order deformations in the
tangent space of the $\mathfrak{X}$-component form the complex
$\operatorname*{dual}\left(  B\left(  I\right)  \right)  $ as defined in
Section \ref{Sec dual complex general setting} and span the Fano polytope
$P^{\circ}$, which gives the embedding toric Fano variety $Y^{\circ}=X\left(
\Sigma\right)  $ with $\Sigma=\Sigma\left(  P^{\circ}\right)  $ for the mirror fibers.

As an example, in \cite{GS An Enumeration of Simplicial 4Polytopes with 8
Vertices} an enumeration of all combinatorial types of
\index{simplicial}%
simplicial $4$-polytopes with $7$ and $8$ vertices is given. These correspond
to reduced monomial Calabi-Yau threefolds $X_{0}$ in $\mathbb{P}^{6}$ and
$\mathbb{P}^{7}$ via the Stanley-Reisner construction. For codimension $4$,
due to the lack of a structure theorem analogous to Theorem
\ref{1BuchsbaumEisenbud}, smoothing of $X_{0}$ has to be addressed by the
deformation theory of Stanley-Reisner rings.

\subsection{Deformations and obstructions of a non-simplicial generalization
of Stanley-Reisner
rings\label{Sec deformations and obstructions of a toric generalization of stanley reisner rings}%
}

Consider the setup of the previous Section
\ref{1FurtherStanleyReisnerexamples}. So let $Y=\mathbb{P}\left(
\Delta\right)  \cong\mathbb{P}^{n}$ with the degree $n+1$ Veronese polytope
$\Delta$, $\Sigma=\operatorname*{NF}\left(  \Delta\right)  $ and
$S=\mathbb{C}\left[  y_{r}\mid r\in\Sigma\left(  1\right)  \right]  $ the Cox
ring of $Y$. The faces of the simplex $\Delta^{\ast}$ correspond to the
subsets of $\Sigma\left(  1\right)  $, i.e., the set of variables of the
homogeneous coordinate ring $S$. Let $Z$ be a simplicial subcomplex of
$\operatorname*{Poset}\left(  \Sigma\right)  $ representing faces as sets of
rays and $I_{Z}\subset S$ the corresponding Stanley-Reisner ideal. As noted in
the previous Section \ref{1FurtherStanleyReisnerexamples} and Section
\ref{Monomial ideals in the Cox ring and the stratified toric primary decomposition}
the isomorphism%
\[%
\begin{tabular}
[c]{cccl}
& $\operatorname*{Poset}\left(  \Delta\right)  $ & $\overset{\cong
}{\rightarrow}$ & $2^{\Sigma\left(  1\right)  }$\\
& $\cup$ &  & $\cup$\\
$\operatorname*{comp}:$ & $\operatorname*{Strata}\nolimits_{\Delta}\left(
I_{Z}\right)  $ & $\overset{\cong}{\rightarrow}$ & $Z$\\
& $F$ & $\mapsto$ & \multicolumn{1}{c}{$\left\{  r\in\Sigma\left(  1\right)
\mid r\not \subset \operatorname*{hull}\left(  F^{\ast}\right)  \right\}  $}%
\end{tabular}
\]
transfers the combinatorial data to a subcomplex of $\operatorname*{Poset}%
\left(  \Delta\right)  $ and
\[
I_{Z}=\left\langle
{\displaystyle\prod\limits_{v\in J}}
y_{v}\mid J\subset\Sigma\left(  1\right)  \text{ with }\operatorname*{supp}%
\left(  \operatorname*{Strata}\nolimits_{\Delta}\left(  I_{Z}\right)  \right)
\subset%
{\displaystyle\bigcup\limits_{v\in J}}
F_{v}\right\rangle
\]
Note that this also works if $Y=X\left(  \Sigma\right)  $ is a toric variety
such that $\Sigma$ is the fan over the faces of a simplex $\Delta^{\ast}$. The
dual description of the ideal $I_{Z}$ via the subcomplex
$\operatorname*{Strata}\nolimits_{\Delta}\left(  I_{Z}\right)  \subset
\operatorname*{Poset}\left(  \Delta\right)  $ should allow for a reformulation
of the formulas for $T^{1}$ and $T^{2}$ by Altmann and Christophersen in terms
of the complex $\operatorname*{Strata}\nolimits_{\Delta}\left(  I_{Z}\right)
$. So one may ask if this allows for a non-simplicial generaliziation of these formulas.

\subsection{Mirrors of Calabi-Yau varieties given by ideals with Pfaffian
resolutions in the Cox rings of toric Fano
varieties\label{Sec mirror pfaffian resolution in Cox ring}}

Let $Y=X\left(  \Sigma\right)  $ be a $\mathbb{Q}$-Gorenstein toric variety of
dimension $n$ with Cox ring $S$. We call a subscheme $X\subset Y$ of
codimension $3$
\index{Pfaffian}%
Pfaffian, if

\begin{enumerate}
\item there is a vector bundle $\mathcal{F}$ on $Y$ of rank $2k+1$ for some
$k\in\mathbb{Z}_{\geq0}$

\item and a
\index{skew symmetric}%
skew symmetric map $\varphi:\mathcal{F}\left(  D\right)  \rightarrow
\mathcal{F}^{\ast}$ for some divisor $D$ such that $\varphi$ degenerates to
rank $\leq2k-2$ in codimension $3$

\item $X$ is scheme theoretically the
\index{degeneracy locus}%
degeneracy locus of $\varphi$.
\end{enumerate}

If $X$ is Pfaffian in $Y$ given by the skew symmetric map $\varphi
:\mathcal{F}\left(  D\right)  \rightarrow\mathcal{F}^{\ast}$ and $\det\left(
\mathcal{F}^{\ast}\right)  =\mathcal{O}_{X}\left(  E\right)  $ then the
resolution of $X$ is of the form%
\[
0\rightarrow\mathcal{O}_{Y}\left(  D-2E\right)  \rightarrow\mathcal{F}\left(
D-E\right)  \rightarrow\mathcal{F}^{\ast}\left(  -E\right)  \rightarrow
\mathcal{O}_{X}%
\]
$\omega_{X}^{\circ}\cong\mathcal{O}_{X}\left(  -D+2E+K_{Y}\right)  $ and $X$
is locally defined by the Pfaffians of $\varphi$.

One may ask if the following generalization of Walters theorem from Section
\ref{1PfaffianCalabiYauThreefolds} holds: If $X$ is an equidimensional,
locally Gorenstein subscheme $X\subset Y$ of dimension $n-3$, $\omega
_{X}^{\circ}\cong\mathcal{O}_{X}\left(  D\right)  $ for some divisor $D$ and
some parity condition similar to that in Theorem \ref{thm walter} is satisfied
then $X$ is Pfaffian.

With respect to these topics see also \cite{EPW Enriques surfaces and other
nonPfaffian subcanonical subschemes of codimension 3}.

We call $X$ globally defined by Pfaffians, if $X$ is Pfaffian where
$\mathcal{F}=\mathcal{O}_{Y}\left(  E_{1}\right)  \oplus...\oplus
\mathcal{O}_{Y}\left(  E_{r}\right)  $ is a direct sum with divisors $E_{j}$.
Then $X$ is defined by the Pfaffians of $\varphi$ in the Cox ring $S$.

Generalizing the work of Tonoli in \cite{Tonoli Canonical surfaces in
mathbbP^5 and CalabiYau threefolds in mathbbP^6} one may construct Pfaffian
Calabi-Yau varieties in toric Fano varieties.

Suppose $\mathfrak{X}\subset Y\times\operatorname*{Spec}\mathbb{C}\left[
t\right]  $ is a monomial degeneration of Calabi-Yau varieties with fibers in
the toric Fano variety $Y$ with Cox ring $S$ given by the ideal $I\subset
\mathbb{C}\left[  t\right]  \otimes S$ as in the setup of the tropical mirror
construction. Let $\mathfrak{X}^{\circ}\subset Y^{\circ}\times
\operatorname*{Spec}\mathbb{C}\left[  t\right]  $ be the mirror degeneration
with fibers in the toric Fano variety $Y^{\circ}$ with Cox ring $S^{\circ}$
given by the ideal $I^{\circ}\subset\mathbb{C}\left[  t\right]  \otimes
S^{\circ}$. For hypersurfaces and complete intersections we have:

Assume that $\mathfrak{X}$ is the degeneration associated to a general
anticanonical toric hypersurface in a Gorenstein toric Fano variety as given
in Section \ref{Degenerationcompleteintersection}. Then the mirror
degeneration is again a degeneration of toric hypersurfaces.

If $\mathfrak{X}\ $is the degeneration associated in Section
\ref{Degenerationcompleteintersection} to a general Calabi-Yau complete
intersection given by a nef partition in the Gorenstein toric Fano variety
$Y^{\circ}$, then the mirror degeneration $\mathfrak{X}^{\circ}\subset
Y^{\circ}\times\operatorname*{Spec}\mathbb{C}\left[  t\right]  $ can be
defined by an ideal $I^{\circ}\subset\mathbb{C}\left[  t\right]  \otimes
S^{\circ}$ with Koszul resolution. Note that the corresponding special fiber
ideal $I_{0}^{\circ}$ will not be $\Sigma^{\circ}$-saturated in general.

So if $\mathfrak{X}$ can be defined by an ideal $J\subset I$ with Pfaffian
resolution we may ask: Is there always a birational model $\mathfrak{\hat{X}%
}^{\circ}\subset\hat{Y}^{\circ}\times\operatorname*{Spec}\mathbb{C}\left[
t\right]  $ of the mirror degeneration $\mathfrak{X}^{\circ}\subset Y^{\circ
}\times\operatorname*{Spec}\mathbb{C}\left[  t\right]  $ which has fibers in a
toric Fano variety $\hat{Y}^{\circ}=X\left(  \Sigma^{\circ}\right)  $,
$\hat{\Sigma}^{\circ}\left(  1\right)  \subset\Sigma^{\circ}\left(  1\right)
$ with Cox ring $\hat{S}^{\circ}$ and is defined by an ideal $\hat{J}^{\circ
}\subset\mathbb{C}\left[  t\right]  \otimes\hat{S}^{\circ}$ with Pfaffian
resolution? Note that again the special fiber ideals in $S$ and $\hat
{S}^{\circ}$ corresponding to $J$ and $\hat{J}^{\circ}$ are not $\Sigma$-
respectively $\hat{\Sigma}^{\circ}$-saturated in general.

\subsection{Tropical geometry and mirror symmetry over finite
fields\label{Sec tropical geometry finite fields}}

Consider the Fermat one parameter family of quintics threefold hypersurfaces
given by%
\[
f\left(  x,t\right)  =\sum_{i=1}^{5}x_{i}^{5}+5t\cdot x_{1}x_{2}x_{3}%
x_{4}x_{5}%
\]
in projective space over $\mathbb{F}_{q}$ with $q=p^{r}$, $p\neq5$. Denote by
$N_{r}\left(  t\right)  $ the number of solutions of $f\left(  x,t\right)  $
in $\mathbb{P}_{\mathbb{F}_{q}}^{2}$. In \cite{COV CalabiYau Manifolds Over
Finite Fields I} these numbers are computed in terms of the periods and
in\linebreak\cite{COV CalabiYau Manifolds Over Finite Fields II} the structure
of the $\zeta$-function%
\[
\zeta\left(  s,t\right)  =\exp\left(
{\textstyle\sum\nolimits_{r=1}^{\infty}}
N_{r}\left(  t\right)  \frac{s^{r}}{r}\right)
\]
is discussed and related to mirror symmetry.

Non-Archimedian and $p$-adic geometry share many similarities. Also, as
observed by L. Tabera, tropical geometry behaves well with respect to finite
fields. One may ask how the $\zeta$-function relates to the tropical data
associated to monomial degenerations.

\subsection{Tropical curves and the $A$-model instanton
numbers\label{Sec instanton numbers}}

Mikhalkin gives in \cite{Mikhalkin Enumerative tropical algebraic geometry in
mathbbR^2} a formula enumerating curves of arbitriary genus in toric surfaces
via tropical geometry. He computes the finite number of curves of genus $g$
and degree $d$ passing through $3d-1+g$ points in general position, i.e., the
Gromov-Witten invariants of $\mathbb{P}^{2}$, by counting tropical curves via
lattice paths of length $3d-1+g$ in the degree $d$ Veronese polytope $\Delta$
of $\mathbb{P}^{2}$. This generalizes to other toric surfaces by replacing the
polytope $\Delta$.

In the context of Calabi-Yau varieties and mirror symmetry we are interested
in the $A$-model correlation functions defined via Gromov-Witten invariants.
The instanton numbers appearing in the Gromov-Witten invariants are related to
the number of rational curves of given degree on the Calabi-Yau variety.

Consider the setup of the tropical mirror construction. One may ask if it is
possible to compute instanton numbers, Gromov-Witten invariants and the
$A$-model correlation functions of the general fiber of the monomial
degeneration $\mathfrak{X}$ in terms of a tropical curve count using the
components of the special fiber of $\mathfrak{X}$.

\subsection{GKZ-hypergeometric differential equations and quantum cohomology
rings of Calabi-Yau varieties\label{Sec GKZ}}

Consider a Calabi-Yau hypersurface in a toric Fano variety $\mathbb{P}\left(
\Delta\right)  $ of dimension $n$ given by a reflexive polytope $\Delta\subset
M_{\mathbb{R}}$. For this setup the Gelfand-Kapranov-Zelevinski hypergeometric
systems are analysed in \cite{Hosono GKZ Systems Groebner Fans and Moduli
Spaces of CalabiYau Hypersurfaces} in the context of mirror symmetry. Via the
local Torelli theorem one can give a local coordinate on the moduli space in
terms of period integrals. For the hypersurface given by%
\[
f_{c}=\sum_{m\in\Delta^{\ast}\cap M}c_{m}x^{m}%
\]
we have one canonical period integral%
\[
\Pi\left(  c\right)  =\frac{1}{\left(  2\pi i\right)  ^{n}}\int_{C_{0}}%
\frac{1}{f_{c}\left(  x\right)  }%
{\textstyle\prod\nolimits_{i=1}^{n}}
\frac{dx_{i}}{x_{i}}%
\]
for the cycle
\[
C_{0}=\left\{  \left\vert x_{i}\right\vert =1\mid i=1,...,n\right\}  \subset
T=\operatorname*{Hom}\nolimits_{\mathbb{Z}}\left(  M,\mathbb{C}^{\ast}\right)
\]
As shown by Batyrev in \cite{Batyrev Variations of the mixed Hodge structure
of affine hypersurfaces in algebraic tori} the period integral satisfies the
following GKZ-hypergeometric system associated to $\mathcal{A}=\left\{
1\right\}  \times\left(  \Delta^{\ast}\cap N\right)  $ and the exponent
$\beta=1\times0$. With the lattice%
\[
L=\left\{  \left(  l_{\delta}\right)  \in\mathbb{Z}^{\mathcal{A}}\mid
\sum_{\delta\in\mathcal{A}}^{n}l_{\delta}\delta=0\right\}
\]
of relations on the elements of $\mathcal{A}$ (see also Section
\ref{Sec toric mori theory}) this system of differential equations is given by%
\begin{align*}
\left(
{\textstyle\prod\nolimits_{l_{\delta}>0}}
\left(  \frac{\partial}{\partial c_{\delta}}\right)  ^{l_{\delta}}-%
{\textstyle\prod\nolimits_{l_{\delta}<0}}
\left(  \frac{\partial}{\partial c_{\delta}}\right)  ^{-l_{\delta}}\right)
\Psi\left(  c\right)   &  =0\text{ for }l\in L\\
\left(  \sum_{\delta\in\mathcal{A}}\delta\cdot c_{\delta}\frac{\partial
}{\partial c_{\delta}}-\beta\right)  \Psi\left(  c\right)   &  =0
\end{align*}
The GKZ-hypergeometric system relates to the period integrals about the
maximal degeneration point of the hypersurface degeneration.

We may ask for a generalization of the hypersurface setup using the complexes
involved in the tropical mirror construction.

\newpage%

\printindex


\noindent\textsc{Department of Mathematics, Universitaet des Saarlandes,}%
\newline Campus E2 4, D-66123 Saarbr\"{u}cken

\noindent\textit{E-mail address: }boehm@math.uni-sb.de

\end{document}